\documentclass[a4paper,
               11pt,
               twoside]{scrreprt}
\usepackage[utf8]{inputenc} %avoid encoding problems
\usepackage[german,british]{babel}
\usepackage[a4paper,width=150mm,top=25mm,bottom=25mm,bindingoffset=6mm]{geometry} %marings

\usepackage[natbib,sortcites=true,sorting=nyt,maxbibnames=20,backend=biber,defernumbers=true]{biblatex}
\DefineBibliographyExtras{british}{}
\usepackage{microtype}

\usepackage{fancyhdr} 
\usepackage{amsmath}
\AtBeginDocument{\numberwithin{equation}{chapter}} % numbering of equations
\usepackage{amssymb}
\usepackage{amsthm}
\usepackage{mathtools} %mathclap
\usepackage{dsfont} %math font
\usepackage{mathrsfs}

\usepackage{graphicx} % Required for inserting images
\usepackage{svg}
\usepackage{float}

\usepackage{tikz}
\usetikzlibrary{backgrounds}
\usetikzlibrary{intersections}
\tikzstyle{block} = [draw, fill=white, rectangle, minimum height=3em, minimum width=6em]
\usepackage{tkz-euclide}

\usepackage{subcaption} %allow for subfigures
\usepackage{scalerel} 

\usepackage[shortlabels]{enumitem} %simplify enumerates

\usepackage[colorlinks=true,hypertexnames=false]{hyperref} %convert references to links
\usepackage{cleveref}

\usepackage{array}
\usepackage{longtable}

\usepackage{calc}    %measure length of symbols 

\graphicspath{{resources/}}
\newcommand{\rbl}{\left (}
\newcommand{\rbr}{\right )}
\newcommand{\sbl}{\left [}
\newcommand{\sbr}{\right ]}
\newcommand{\al}{\left \langle}
\newcommand{\ar}{\right \rangle}
\newcommand{\bl}{\left |}
\newcommand{\br}{\right |}
\newcommand{\nl}{\left\|}
\newcommand{\nr}{\right\|}
\newcommand{\cbl}{\left\lbrace }
\newcommand{\cbr}{\right\rbrace }
\newcommand{\Abs}[1]{\bl #1 \br}

\newcommand{\Norm}[2][ ]{\nl #2 \nr_{#1}}
\newcommand{\SNorm}[1]{\Norm[\infty]{#1}}
\newcommand{\LNorm}[2][2]{\Norm[L^{#1}]{#2}}

\newcommand{\N}{\mathds{N}}
\newcommand{\Z}{\mathds{Z}}
\newcommand{\R}{\mathds{R}}
\newcommand{\Rp}{\R_{\geq0}}
\newcommand{\Rpp}{\R_{>0}}
\newcommand{\C}{\mathds{C}}

\newcommand{\eps}{\varepsilon}
\renewcommand{\phi}{\varphi}

\newcommand{\cB}{\mathcal{B}}
\newcommand{\cC}{\mathcal{C}}
\newcommand{\cD}{\mathcal{D}}
\newcommand{\cE}{\mathcal{E}}
\newcommand{\cF}{\mathcal{F}}
\newcommand{\cG}{\mathcal{G}}

\newcommand{\cK}{\mathcal{K}}
\newcommand{\cL}{\mathcal{L}}
\newcommand{\cM}{\mathcal{M}}
\newcommand{\cN}{\mathcal{N}}
\newcommand{\cP}{\mathcal{P}}
\newcommand{\cR}{\mathcal{R}}

\newcommand{\cT}{\mathcal{T}}
\newcommand{\cU}{\mathcal{U}}
\newcommand{\cY}{\mathcal{Y}}

\newcommand{\impl}{\Rightarrow}
\newcommand{\Impl}{\Longrightarrow}
\newcommand{\To}{\longrightarrow}

\newcommand{\sub}{\subseteq}

\newcommand{\fa}{\forall \, }
\newcommand{\ex}{\exists \, }

\renewcommand{\d}{\mathrm{d}}
\newcommand{\me}{\mathrm{e}} % nice e function

\newcommand{\setdef}[2]{\cbl#1\left|\vphantom{#1} #2\right.\cbr}

\newcommand{\textover}[3][l]{%
 \makebox[\widthof{#3}][#1]{#2}%
 }

\newcommand{\dd}[2][ ]{\tfrac{\text{\normalfont d}#1}{\text{\normalfont d}#2}}
\newcommand{\Dd}[2][ ]{\frac{\text{\normalfont d}#1}{\text{\normalfont d}#2}}

\makeatletter
\newcommand{\Itemlabel}[2]{#2\def\@currentlabel{#2}\label{#1}}
\makeatother

\makeatletter
\newcommand\setcurrentname[1]{\def\@currentlabelname{#1}}
\makeatother

\DeclareMathOperator*{\esssup}{ess\,sup}
\newcommand{\GL}{\text{GL}}
\newcommand{\SGroup}{\text{S}}
\DeclareMathOperator*{\graph}{graph}
\DeclareMathOperator*{\Lip}{Lip}
\DeclareMathOperator*{\im}{im}
\DeclareMathOperator*{\spec}{spec}

\DeclareMathOperator*{\rf}{ref}
\DeclareMathOperator*{\loc}{loc}

\newcommand{\nocontentsline}[3]{}
\newcommand{\tocless}[2]{\bgroup\let\addcontentsline=\nocontentsline#1{#2}\egroup}
\newcommand{\toclesslab}[3]{\bgroup\let\addcontentsline=\nocontentsline#1{#2\label{#3}}\egroup} %with label

\counterwithout{equation}{section}
\counterwithout{figure}{section}
\newtheorem{definition}{Definition}[section]
\theoremstyle{definition}
\newtheorem{remark}[definition]{Remark}\Crefname{remark}{Remark}{Remarks}
\newtheorem{algo}[definition]{Algorithm}\Crefname{algo}{Algorithm}{Algorithms}
\newtheorem{example}[definition]{Example}\Crefname{example}{Example}{Examples}
\newtheorem{assn}[definition]{Assumption}\Crefname{assn}{Assumption}{Assumptions}
\theoremstyle{plain}
\newtheorem{prop}[definition]{Proposition}\Crefname{prop}{Proposition}{Propositions}
\newtheorem{corollary}[definition]{Corollary}\Crefname{corollary}{Corollary}{Corollaries}
\Crefname{assertion}{Assertion}{Assertions}
\newtheorem{theorem}[definition]{Theorem}\Crefname{theorem}{Theorem}{Theorems}
\newtheorem{lemma}[definition]{Lemma}\Crefname{lemma}{Lemma}{Lemmata}
\AtEndEnvironment{example}{\unskip\ \null\hfill\ensuremath{\diamond}}
\AtEndEnvironment{remark}{\unskip\ \null\hfill\ensuremath{\bullet}}
\AtEndEnvironment{algo}{\unskip\ \null\hfill\ensuremath{\blacktriangle}}
\newcommand{\oT}{\mathbf{T}}
\newcommand{\oTM}{\oT_{\mathrm{M}}}
\newcommand{\umax}{u_{\max}}
\newcommand{\xM}{x_{\mathrm{M}}}
\newcommand{\xMd}{\dot{x}_{\mathrm{M}}}
\newcommand{\xMh}{\hat{x}_{\mathrm{M}}}
\newcommand{\xMt}{\tilde{x}_{\mathrm{M}}}
\newcommand{\oTMh}{\hat{\oT}_{\mathrm{M}}}
\newcommand{\oTMt}{\tilde{\oT}_{\mathrm{M}}}

\newcommand{\InitValues}{\mathfrak{I}_{t_0,\tau}^{\Psi}}
\newcommand{\PropInitValues}{\mathfrak{PI}_{t_0,\tau}^{\Psi,\eps,\lambda}}
\newcommand{\InitState}{\hat{\mathfrak{X}}}
\newcommand{\InitStateK}{\mathfrak{X}}
\newcommand{\FunnelTrajectories}{\cY^{\Psi}}
\newcommand{\FCTrajectories}{\mathfrak{Y}^{\phi}}
\newcommand{\Controls}{\cU_{[\hat{t},\hat{t}+T]}}
\newcommand{\FunnelBoundaryFuncs}{\mathscr{G}}

\newcommand{\InitStrategy}{\kappa}

\newcommand{\Funnel}{\psi}
\newcommand{\UltMFunnel}{\psi_{r}}

\newcommand{\FunnelPenaltyFunc}{\nu_{\Funnel}}
\newcommand{\FunnelPenaltyFuncR}{\nu_{\Funnel_r}}
\newcommand{\FunnelPenaltyFuncWithout}{\nu}

\newcommand{\OrigFunnelPenaltyFunc}{\tilde{\ell}_{\psi}}

\newcommand{\OrigFunnelStageCost}{\ell_{\psi}}
\newcommand{\FunnelStageCost}{\ell_{\psi_{r}}}

\newcommand{\cDSet}{\cD^{\Psi}_{t}}

\newcommand{\cEFC}[1]{\cE_{r}^{#1}}

\newcommand{\Lpath}{\gamma}
\newcommand{\Lebesgue}{\lambda}
\newcommand{\Indic}{\mathds{1}}
\newcommand{\OpChi}[1][r]{\chi_{#1}}

\newcommand{\yM}{y_{\mathrm{M}}}

\newcommand{\yMd}{\dot{y}_{\mathrm{M}}}
\newcommand{\yMdk}{\dot{y}^k_{\mathrm{M}}}
\newcommand{\yMh}{\hat{y}_{\mathrm{M}}}
\newcommand{\hM}{h_{\mathrm{M}}}
\newcommand{\fM}{f_{\mathrm{M}}}
\newcommand{\fMmax}{f_{\mathrm{M}}^{\max}}
\newcommand{\fmax}{f^{\max}}
\newcommand{\FM}{F_{\mathrm{M}}}
\newcommand{\gM}{g_{\mathrm{M}}}

\newcommand{\gMmax}{g_{\mathrm{M}}^{\max}}
\newcommand{\gmax}{g^{\max}}
\newcommand{\gmin}{g^{\min}}
\newcommand{\gMmin}{g_{\mathrm{M}}^{\min}}
\newcommand{\gMInvmax}{g_{\mathrm{M}}^{-1\max}}
\newcommand{\GM}{G_{\mathrm{M}}}

\newcommand{\SampleTime}{\mathfrak{r}}

\newcommand{\ModSysClassDiscr}{\mathfrak{N}^{m,r}_{t_0}}

\newcommand{\eSTrack}{e_{\mathrm{S}}}
\newcommand{\eSTrackdot}{\dot{e}_{\mathrm{S}}}
\newcommand{\eMTrack}{e_{\mathrm{M}}}

\newcommand{\eM}{\xi}

\newcommand{\eS}{e}

\newcommand{\uFMPC}{u_{\mathrm{FMPC}}}
\newcommand{\uFMPCh}{\hat{u}_{\mathrm{FMPC}}}
\newcommand{\uFC}{u_{\mathrm{FC}}}
\newcommand{\uFCh}{\hat{u}_{\mathrm{FC}}}
\newcommand{\uZoH}{u_{\mathrm{ZoH}}}

\newcommand{\uFMPCk}[1][k]{u_{\mathrm{FMPC},{#1}}}
\newcommand{\uFCk}[1][k]{u_{\mathrm{FC},{#1}}}

\newcommand{\SetOfSignals}{\mathfrak{S}}

\newcommand{\LShift}{S_m}
\newcommand{\ActivFunc}{\mathfrak{a}}
\newcommand{\FCSurjec}{\mathcal{N}}
\newcommand{\FCBijec}{\gamma}
\newcommand{\HighGainFunc}{\mathfrak{h}}

\newcommand{\FunContr}{\gamma}
\newcommand{\FunDeriv}{\alpha}
\newcommand{\FunDiam}{\beta}

\newcommand{\FCDiscreteGain}{\nu}
\newcommand{\FCDiscreteThresh}{\iota}

\allowdisplaybreaks % allow for page break within align

\begin{document}

\setcounter{page}{1}
\pagenumbering{Roman}

\begin{titlepage}
\centering
\begin{tabular}{rl}
\renewcommand{\baselinestretch}{1.25}\normalsize
\includesvg{logo-upd.svg}
&
\begin{minipage}{2\linewidth}
\large{\textbf{Universität Paderborn}}\\
\large{Fakultät für Elektrotechnik, Informatik und Mathematik}\\
\large{Fachgebiet Systemtheorie}\\[4ex]
\end{minipage}
\end{tabular}
\vfill

{\Huge \textbf{Model Predictive Control for output tracking with prescribed performance}}\\[16ex]
{Dissertation zur Erlangung des akademischen Grades}\\[1ex]
{\large{\textbf{Doctor rerum naturalium (Dr. rer. nat.)}}}\\[4ex]

{von }\\[1ex]
{\Large{\textbf{Dario Rudolf Walter Dennstädt}}}\\[8ex]

{betreut durch}\\[1ex]
{\Large{\textbf{Jun. Prof.~Dr.~Thomas Berger}}}\\[8ex]
\vfill

\begin{tabular}{c}
Paderborn 2025 
\end{tabular}

\end{titlepage}

\counterwithin{figure}{chapter}
\newpage\null\thispagestyle{empty}
\newpage\null\thispagestyle{empty}
\newpage\null\thispagestyle{empty}

\begin{otherlanguage}{german}
\chapter*{Zusammenfassung}
Modellprädiktive Regelung (MPC) stellt einen Eckpfeiler der modernen
Regelungstheorie dar und erlaubt die simultane Berücksichtigung von
Nebenbedingungen sowie die multikriterielle Optimierung durch iterative
Vorhersage und Receding‑Horizon‑Optimierung. 
In der Praxis sehen sich MPC-Verfahren jedoch drei wesentlichen
Herausforderungen konfrontiert:
der Sicherstellung initialer und rekursiver Zulässigkeit (d.\,h. der 
dauerhaften Lösbarkeit des zugrundeliegenden Optimierungsproblems), ihrer Robustheit
gegenüber Modellabweichungen und unbekannten Störungen sowie
Restriktionen bei der Realisierung als Abtastsystem.

Diese Dissertation entwickelt ein innovatives MPC-Framework für nichtlineare,
zeitkontinuierliche Systeme, die mittels funktionaler Differentialgleichungen
beschrieben werden. Ziel ist die Ausgangsfolgeregelung glatter Referenzsignale
innerhalb vorgegebener Fehlertoleranzen zu gewährleisten und die genannten
Herausforderungen systematisch zu adressieren.

Im Mittelpunkt steht \emph{Funnel MPC} -- ein neuartiger Regelungsansatz, der
auf herkömm\-liche Endbedingungen und restriktiv lange Prädiktionshorizonte
verzichtet. Sein Fundament bilden sogenannte Funnel-Penalty-Funktionen:
Kostenfunktionen, die Abweichungen des Trackingfehlers von zeitvarianten
Toleranzschranken gezielt bestrafen. Angelehnt an Techniken der adaptiven
Funnel-Regelung garantiert dieser Ansatz sowohl initiale als auch rekursive Zulässigkeit
und gewährleistet zugleich strikte Einhaltung der Soll-Regelgüte.

Darauf aufbauend wird Funnel MPC mit der modellfreien Funnel-Regelung 
in eine hybride Zwei-Komponenten-Architektur
verschmolzen. Diese vereint modellbasierte Optimierung mit adaptiver Ausgangs\-rückführung,
um die konkurrierenden Ziele Optimalität und Robustheit auszubalancieren.
Ergebnis ist ein Regler, der die geforderte Regelgüte selbst bei strukturellen
Modellungenauigkeiten, unmodellierten Dynamiken und Störungen zuverlässig einhält.

Zur Steigerung der Prädiktionsgenauigkeit integrieren wir ein datengesteuertes
Lernverfahren, welches das Systemmodell fortlaufend basierend auf Online-Messungen
adaptiert. Diese Komponente reduziert Modell-System-Diskrepanzen kontinuierlich
und verbessert auf diese Weise langfristig die Regelgüte, ohne dabei
Robustheitsgarantien zu kompromittieren.

Schließlich überführen wir die zeitkontinuierlichen Regelgesetze in
eine Abtastimplementierung. Durch Herleitung expliziter Schranken für Abtastrate
und Stellaufwand garantieren wir Stabilität unter treppenförmigen Stellsignalen –
ein essenzieller Schritt für die praktische Umsetzung auf digitaler Hardware.

Durch systematische Verknüpfung von Zulässigkeit, Robustheit, Lernfähigkeit und
Abtastimplementierung entsteht ein ganzheitliches Framework zur Einhaltung vorgegebener 
Fehlertoleranzen bei der Ausgangsfolgeregelung für eine breite Klasse dynamischer
Systeme. Die vorgestellten Ergebnisse ebnen den Weg für zukünftige Entwicklungen im Bereich
des lernunterstützten und samplingbasierten, robusten MPC.
\end{otherlanguage}
\addcontentsline{toc}{chapter}{Abstracts}
\newpage\null\thispagestyle{empty}
\chapter*{Abstract}
Model Predictive Control (MPC) is a cornerstone of modern control theory,
offering a versatile framework for constraint handling and multi-objective
optimisation through iterative prediction and receding-horizon optimisation.
However, its practical application can face critical challenges: ensuring initial
and recursive feasibility (guaranteeing solvability of the underlying optimisation problem),
robustness against system-model mismatches
and unknown disturbances, and sampled-data implementation constraints.

This thesis develops a novel MPC framework for a class of non-linear continuous-time
systems governed by functional differential equations, targeting output tracking
of smooth reference signals within prescribed error bounds, while systematically
addressing the aforementioned challenges.

We first introduce \emph{funnel MPC}, a novel algorithm that eliminates reliance
on commonly used terminal conditions or restrictive long prediction horizons. At
its core are funnel penalty functions -- state costs that penalise deviations of
the tracking error from prescribed time-varying boundaries. Inspired by adaptive
funnel control principles, this framework ensures initial and recursive
feasibility while rigorously enforcing tracking performance guarantees.

Building on this foundation, we unify funnel MPC with model-free funnel feedback
into a two-component hybrid architecture. This structure synergises model-based
optimisation with adaptive feedback compensation, reconciling the competing
objectives of optimality and robustness. The resulting controller 
achieves prescribed tracking performance despite structural model-plant
mismatches, unmodelled dynamics, and disturbances.

To further enhance predictive accuracy, we introduce a data-driven learning
framework that iteratively refines the model using system  measurements.
This component enables the controller to mitigate model-plant discrepancies over
time, improving long-term performance without compromising robustness
guarantees. 
Bridging theory and practice, we finally formalise the transition from
continuous-time control laws to sampled-data implementations, deriving explicit
bounds on sampling rates and control effort to guarantee stability under
piecewise constant control signals -- 
a critical step toward deploying the algorithm on digital hardware.

By systematically addressing feasibility, robustness, learning integration,
and sampled-data implementation, this thesis establishes a cohesive
framework to ensure output tracking within prescribed error bounds for a large system class.
The results pave the way for future advances in learning-enhanced and sampled-data robust MPC.

\newpage\null\thispagestyle{empty}
\chapter*{Acknowledgement}
The completion of this dissertation marks the culmination of a significant
chapter in my life, and it would not have been possible without the unwavering
support, guidance, and encouragement of numerous individuals. It is with
profound gratitude that I acknowledge their invaluable contributions.\\

First and foremost, I extend my deepest appreciation to my supervisor,
Jun.\,Prof.\,Dr. Thomas Berger.
His expert guidance, insightful critiques, and constant
encouragement were indispensable throughout this research journey. His
intellectual rigour, patience, and unwavering belief in this project were
fundamental to the development and completion of this work. I am immensely
grateful for his dedication and mentorship.\\

I am also deeply indebted to Prof.\,Dr.\,Karl\,Worthmann for his extraordinary
generosity with his time and expertise. Throughout this process, he provided
invaluable advice, critical feedback, and insightful perspectives that
significantly shaped my research and enhanced the quality of this dissertation.
His willingness to engage deeply with my work and offer his support was
immensely appreciated and crucial to my progress.\\

A very special and heartfelt thank you goes to Prof.\,Dr.\,Achim\,Ilchmann for his
mentorship during my studies. It is not an exaggeration to say that without his
initial encouragement and unwavering belief in my potential, I would not have
returned to university to pursue this PhD. His inspiration and guidance
set me on this path, and for that, I am deeply grateful.\\

I sincerely thank the reviewers of this dissertation, Prof.\,Dr.-Ing.\,Timm\,Faulwasser
and Prof.\,Dr.\,Felix\,Schwenninger, for their time and careful consideration of my
work, and their constructive feedback. 
Their expertise and scrutiny are greatly appreciated.\\

My sincere thanks also go to my dear friend and colleague, Dr.\,Lukas\,Lanza, for his
collaboration, stimulating discussions, and numerous helpful suggestions during
the writing process. His camaraderie and intellectual input were a constant
source of support, inspiration, and motivation.\\

Finally, and most importantly, my boundless gratitude goes to my family and
friends. Their constant moral support, understanding, and unwavering belief in
me carried me through the most challenging times of this PhD journey. They
offered encouragement when my own hope faltered, provided perspective when I
needed it most, and celebrated every small victory along the way. This
achievement is deeply shared with them.
\newpage\null\thispagestyle{empty}

\tableofcontents 
\newpage\null\thispagestyle{empty}

\chapter{Introduction}
\pagenumbering{arabic}
\setcounter{page}{1}
Model Predictive Control (MPC) is a versatile, optimisation-based control
technique widely recognised for its effectiveness in managing linear and
non-linear multi-input multi-output systems, as discussed in
textbooks~\cite{GrunPann17,rawlings2017model}. We also refer to~\cite{Lee2011}
for an overview of the historical development of MPC. A hallmark of MPC is its
ability to explicitly incorporate both control and state constraints into the
optimisation framework, a feature that has propelled its adoption in diverse
applications, see e.g.~\cite{QinBadg03} and also~\cite{samad2020industry}. At
its core, MPC leverages a dynamic model of the system to iteratively forecast
its behaviour over a finite-time horizon. These predictions enable the
controller to solve a receding-horizon Optimal Control Problem (OCP), optimising
control inputs to balance competing objectives -- such as setpoint tracking or
energy efficiency -- against hard constraints like actuator saturation or
safety-critical state bounds. After applying the first part of this optimal
control to the system, the prediction horizon is shifted forward in time and the
model is re-initialised with measurement data from the system. Repeating this
process ad infinitum forms a closed-loop control system.

Despite its simplicity and conceptual elegance, practical implementation of MPC
demands rigorous attention to mathematical foundations. Foremost among these is
ensuring both \emph{initial feasibility} (existence of a valid solution of the
OCP at startup) and \emph{recursive feasibility}, which guarantees that
solvability of the optimal control problem at one time step automatically
implies solvability at the successor time instant. Providing these guarantees
becomes precarious if the utilised model deviates from the actual system, as MPC
relies heavily on model accuracy in order to predict the behaviour of the actual
system. Robustness to these discrepancies is a ubiquitous challenge as all
models are inherently approximate and real-world systems face unmeasurable
disturbances, parametric drift, or unmodelled dynamics. Consequently, designing
robust MPC algorithms capable of handling plant-model mismatches and external
disturbances remains an active research area. One branch of research focuses on
enhancing the MPC algorithm itself via robustification methods to harden the
controller against bounded uncertainties, see e.g. \cite{bemporad99} for an
overview of available techniques. Another branch of research explores adapting
the model to ensure robust constraint satisfaction of the actual system. The
latter has gained momentum with recent advances in machine learning and spurred
interest in integrating techniques like reinforcement learning (RL). On the
implementation front, practical limitations persist: sampled-data architectures
restrict controllers to discrete-time measurements, and hardware constraints
often necessitate piecewise constant control signals, introducing discretisation
errors that further complicate theoretical analyses.

This thesis progressively develops an MPC algorithm for continuous-time systems
to address output tracking of smooth reference signals with prescribed error
bounds. We systematically resolve the aforementioned challenges by first
outlining alternative methods and then proposing a novel approach. By
integrating principles from the adaptive control technique \emph{funnel
control}, we establish initial and recursive feasibility without relying on
terminal conditions or restrictive assumptions, such as demanding excessively
long prediction horizons. Building on this foundation, we unify the two control
strategies into a single, hybrid framework -- combining MPC's predictive
optimisation with funnel control’s adaptability -- to ensure robustness against
unknown disturbances and structural plant-model mismatches. Next, we investigate
the incorporation of a learning mechanism into the framework, enabling
data-driven adaptation of the underlying model to refine predictions using
system measurements. Finally, we derive sufficient conditions for the sampling
rate to guarantee stability when operating the controller in a sampled-data
setting with piecewise constant control signals, bridging the gap between
theoretical continuity and practical digital implementation.

\section{Problem formulation}
We consider non-linear multi-input multi-output control systems of order $r\in\N$ of the form 
\begin{equation}\label{eq:Sys}
    \begin{aligned}
    &\textover[r]{$y^{(r)}(t)$}{$\big(y(t_0),\ldots,y^{(r-1)}(t_0)\big)\ $}   = F(\oT(y,\dot{y},\ldots, y^{(r-1)})(t),u(t)),\\
    &\left.
        \begin{aligned}
       y|_{[0,t_0]}&= y^0  \in \cC^{r-1}([0,t_0],\R^m), && \mbox{if } t_0 >0,\\
       \big(y(t_0),\ldots,y^{(r-1)}(t_0)\big)&= y^0  \in\R^{rm}, && \mbox{if } t_0 = 0,
        \end{aligned} \right\}
    \end{aligned}
\end{equation}
with $t_0\geq 0$,  initial trajectory $y^0$, input $u\in L_{\loc}^{\infty}([t_0,\infty),\R^m)$,
and output $y(t)\in\R^m$ at time $t\geq t_0$.
Note that $u$ and $y$ have the same dimension~$m\in\N$.
The system consists of an \emph{unknown} continuous function $F \in \cC(\R^q\times \R^m,\R^m)$
satisfying the so-called \emph{perturbation high-gain property} introduced in~\Cref{Def:SystemClass}~\ref{Item:PerturbationHighGain},
and an \emph{unknown} operator $\oT$.
The operator~$\oT$ is causal, locally Lipschitz, and satisfies a bounded-input bounded-output property.
These properties will be introduced in detail in~\Cref{Def:OperatorClass} and 
the system under consideration is characterised in~\Cref{Def:SystemClass}.
Note that the system may also incorporate bounded disturbances $d\in L^\infty(\Rp,\R^p)$.  
They can be modelled as part of the unknown operator $\oT$, as we will discuss
in \Cref{Rem:SystemClass} \ref{Item:OperatorContainsDisturbences}.  For reasons
of simplicity, we however refrain from explicitly including them in
equation~\eqref{eq:Sys}.

\subsection{Control objective}\label{Sec:ControlObjective}
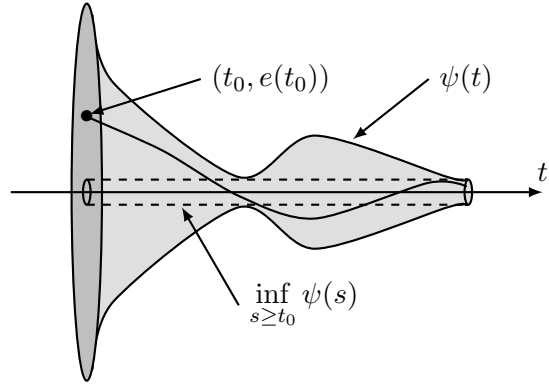
\begin{figure}[ht]
  \begin{center}
    \begin{tikzpicture}[scale=0.5]
\tikzset{>=latex}
  \filldraw[color=gray!25] plot[smooth] coordinates {(0.15,4.7)(0.7,2.9)(4,0.4)(6,1.5)(9.5,0.4)(10,0.333)(10.01,0.331)(10.041,0.3) (10.041,-0.3)(10.01,-0.331)(10,-0.333)(9.5,-0.4)(6,-1.5)(4,-0.4)(0.7,-2.9)(0.15,-4.7)};
  \draw[thick] plot[smooth] coordinates {(0.15,4.7)(0.7,2.9)(4,0.4)(6,1.5)(9.5,0.4)(10,0.333)(10.01,0.331)(10.041,0.3)};
  \draw[thick] plot[smooth] coordinates {(10.041,-0.3)(10.01,-0.331)(10,-0.333)(9.5,-0.4)(6,-1.5)(4,-0.4)(0.7,-2.9)(0.15,-4.7)};
  \draw[thick,fill=lightgray] (0,0) ellipse (0.4 and 5);
  \draw[thick] (0,0) ellipse (0.1 and 0.333);
  \draw[thick,fill=gray!25] (10.041,0) ellipse (0.1 and 0.333);
  \draw[thick] plot[smooth] coordinates {(0,2)(2,1.1)(4,-0.1)(6,-0.7)(9,0.25)(10,0.15)};
  \draw[thick,->] (-2,0)--(12,0) node[right,above]{\normalsize$t$};
  \draw[thick,dashed](0,0.333)--(10,0.333);
  \draw[thick,dashed](0,-0.333)--(10,-0.333);
  \node [black] at (0,2) {\textbullet};
  \draw[->,thick](4,-3)node[right]{\normalsize$\inf\limits_{s\geq t_0}\psi(s)$}--(2.5,-0.4);
  \draw[->,thick](3,3)node[right]{\normalsize$(t_0,e(t_0))$}--(0.07,2.07);
  \draw[->,thick](9,3)node[right]{\normalsize$\psi(t)$}--(7,1.4);
\end{tikzpicture}
\end{center}
 \vspace*{-2mm}
 \caption{Error evolution in a funnel $\mathcal F_{\psi}$ with boundary $\psi(t)$. The figure is based on~\cite[Fig. 1]{BergLe18a}, adapted to the current setting.}
 \label{Fig:funnel}
 \end{figure}

Our objective is to design a control strategy which allows tracking of a given reference
trajectory~$y_{\rf}\in W^{r,\infty}(\Rp,\R^{m})$ within pre-specified error bounds. To be more
precise, the tracking error ~$t\mapsto e(t)\coloneqq y(t)-y_{\rf}(t)$ shall evolve within the prescribed
performance funnel
\begin{equation}\label{eq:DefFunnel}
    \cF_{\Funnel}\coloneqq  \setdef{(t,e)\in \Rp\times\R^{m}}{\Norm{e} < \Funnel(t)}.
\end{equation}
This funnel is determined by the choice of the function~$\Funnel$
belonging  to the set
\begin{equation}\label{eq:DefSetOfFunnelFunctions}
    \cG\coloneqq \setdef
        {\Funnel\in W^{1,\infty}(\Rp,\R)}
        {
            \inf_{t\geq 0}\Funnel(t) > 0
        },
\end{equation}
see also Figure~\ref{Fig:funnel}.
Note that, for a function $\psi\in\cG$, there exists $\lambda>0$ such that
$\psi(t)\geq\lambda$ for all $t \ge 0$. Therefore, signals evolving in $\mathcal{F}_{\psi}$
are not forced to converge to $0$ asymptotically.
\begin{remark}
    In many practical applications perfect tracking is neither possible nor desired. Usually, the
    objective rather is to ensure the tracking error to be less than an (arbitrary small)
    a priori specified constant after a pre-specified period of time and to guarantee that the error does
    not exceed this bound at a later time.
    Tracking within a funnel, or in other words practical tracking, is advantageous since it allows
    tracking for system classes where asymptotic tracking is not possible or requires – when compared to
    asymptotic tracking – much less control effort.
    Note  that the function~$\psi$ is a design parameter, thus its choice is completely up to the
    designer.
    Moreover, arbitrary funnel functions -- and not restricted to
    constant or monotonous decreasing funnels -- give the user more flexibility in
    finding a suitable trade-off between tracking performance and control effort. Typically, the
    specific application dictates the constraints on the tracking error and thus indicates suitable
    choices for $\psi$. During safety critical system phases, the funnel will be small, while
    during non-critical phases the funnel can be widened again to reduce the control effort.
\end{remark}

\subsection{Funnel control}\label{Sec:FunnelControl}
In the context of output reference tracking within prescribed, possibly
time-varying, performance boundaries, \emph{funnel control} is an established
adaptive high-gain control methodology. The concept has been introduced in the
seminal work~\cite{IlchRyan02b} with the design goal of achieving the control
objective laid out in~\Cref{Sec:ControlObjective} for a broad class of
non-linear multi-input, multi-output systems. It has since received a lot
of research attention, see e.g.~\cite{HackHopf13,BergLe18a,BergIlch21}. For
a comprehensive literature overview, we also recommend the recent survey
paper~\cite{BergIlch23}. The funnel control approach solely invokes certain
structural assumptions about the system, namely stable internal dynamics and a
known globally defined relative degree with a globally pointwise sign-definite
high-frequency matrix. Under these conditions, the adaptive controller offers
robustness against disturbances and guarantees specified transient behaviour
without relying on explicit knowledge about the system to be controlled.
Since it ensures output tracking of reference signals within prescribed
performance bounds without having to resort to a model of the system, funnel
control proved useful for tracking problems in various applications such as
DC-link power flow control, see~\cite{SenfPaug14}, control of industrial
servo-systems, see~\cite{Hackl17}, and temperature control of chemical reactor
models, see~\cite{IlchTren04}.

\emph{Prescribed performance control} (PPC) is a methodology closely related to
funnel control. It was first introduced in \cite{BECHLIOULIS2008}. The core
idea of PPC involves transforming the original controlled system into a new
state-space representation using predefined performance functions that encode
desired transient and, potentially, steady-state behaviours. By ensuring the
uniform boundedness of the transformed system’s states via appropriate control
laws, the tracking problem for the original system is solved -- a result that is
both necessary and sufficient under this framework. While early PPC designs
relied on neural networks to approximate unknown non-linearities, later work in
\cite{BECHLIOULIS2014} developed an approximation-free scheme tailored for
systems in so-called pure feedback form. For a detailed and comprehensive overview 
of prescribed performance control, we also recommend
the survey paper~\cite{BuPPC23}. Although PPC and funnel control share a
similar objective -- enforcing error trajectories within predefined bounds --
they differ in their system classes and structural assumptions. Funnel control
applies to systems of the form~\eqref{eq:Sys}, whereas PPC addresses systems
structured as:
\begin{align*}
   \dot{x}_k(t)&=f_k(x_1(t),\ldots,x_{k+1}(t)),\qquad k=1,\ldots, r-1,\\
   \dot{x}_r(t)&=f_r(d(t),x_1(t),\ldots,x_{r}(t),z(t),u(t)),\\
   \dot{z}(t)&=g(d(t),x_1(t),\ldots,x_{r}(t),z(t)),\\
   y(t)&=x_1(t),
\end{align*}
with $f_k:\R^{km}\to\R^m$ for $k=1,\ldots, r-1$, $f_r:\R^{n+rm+m}\to\R^m$, 
$g:\R^{n+rm+q}\to\R^q$, and  $d\in L^\infty(\Rp,R^n)$ represents a bounded disturbance.
Crucially, both methods operate under minimal system knowledge, requiring only
generic structural assumptions rather than explicit functional details. 
While prescribed performance control presumes the partial derivatives
$\frac{\partial f_i}{\partial x_i}$  and $\frac{\partial f_r}{\partial u}$ to be
uniformly positive definite, see \cite{BECHLIOULIS2014}, 
funnel control assumes the system~\eqref{eq:Sys} to have bounded-input bounded-state stable
internal dynamics and the function $F$ to  satisfy the so-called \emph{high-gain property}, 
see \cite{BergIlch21}.
A key distinction between the two control techniques lies in their information
requirements: funnel control relies solely on the output $y$ and its derivatives
while PPC necessitates full state feedback. For the latter, this requirement was
softened in~\cite{Dimanidis20} via the incorporation of a high-gain observer and
\cite{bechlioulis2011} allows internal dynamics of a certain hierarchical
structure, so-called \emph{dynamical uncertainties}. Despite their conceptual
overlap, a rigorous comparative analysis of these approaches remains an open
research question.

Both funnel control and prescribed performance control face inherent limitations
due to their model-free nature. Since neither approach utilises a system model,
the controllers lack predictive capabilities, leaving the error evolution within
time-varying boundaries uncertain. For instance, the error trajectory may
approach the funnel boundary arbitrarily closely, triggering excessively large
feedback gains. This can lead to high-magnitude control inputs, peaking signals,
and -- from an implementation perspective -- significant sensitivity to measurement
noise. While theoretical guarantees ensure bounded control signals, their
precise upper bounds require knowledge of the system and remain a priori unknown. 
For practical implementation on digital devices, both schemes also demand high
sampling rates to maintain feasibility, imposing stringent hardware
requirements. The recent work~\cite{Berger24Internal} demonstrates that
incorporating an internal model into funnel control can markedly enhance
performance. This integration reduces noise sensitivity, mitigates extreme gain
behaviour, and can even achieve asymptotic tracking without a performance funnel
whose width converges to zero. Furthermore, numerical simulations
in~\cite{berger2019learningbased} reveal that an MPC strategy, blending funnel
control principles with predictive optimisation, outperforms pure funnel
control. This approach achieves smaller control actions and relaxed sampling
rate demands for zero-order-hold implementations, highlighting the value of a
model integration in the controller design.

Building upon the ideas from \cite{berger2019learningbased}, this thesis
develops an MPC scheme that integrates ideas from funnel control in order to
achieve the control objective of output tracking with prescribed performance for
systems of the form~\eqref{eq:Sys}. By combining the control methodologies,
we circumvent the shortcomings of both individual approaches. This enables us to
benefit from the best of both worlds: guaranteed feasibility  and robustness
(funnel control), and a superior control performance (MPC).

Before we focus on developing this MPC scheme, we would like to give an intuition 
of the funnel controller's functioning.
\begin{prop}\label{Intro:Prop:FC}
    Let $f:\R^m\to\R^m$ be a locally Lipschitz continuous function, $\Funnel\in\cG$,
    $y_{\rf}\in W^{1,\infty}(\Rp,\R^{m})$, and $y^0\in\R^m$ with
    $\Norm{y^0-y_{\rf}(0)}<\Funnel(0)$.
    Then, the application of the output feedback $u(t)\coloneqq \mu_{\mathrm{FC}}(t,y(t))$ with
    \begin{equation}\label{Intro:eq:FC}
        \mu_{\mathrm{FC}}(t,y) = -k(t,y) e(t,y),\qquad
        k(t,y)          = \frac{1}{\Funnel^2(t)- \Norm{e(t,y)}^2},\qquad
        e(t,y)          = y - y_{\rf}(t)
    \end{equation}
    to the system 
    \[
        \dot{y}(t)= f(y(t)) + u(t),\quad y(0)=y^0,
    \]
    leads to the closed-loop initial value problem
    \begin{align*}
        \dot{y}(t)  = f(y(t)) - \frac{y(t) - y_{\rf}(t)}{\Funnel^2(t)-\Norm{y(t) - y_{\rf}(t)}^2},\quad y(0)=y^0,
    \end{align*}
    which has a solution.  Moreover, every solution can be extended to a unique global
    solution ${y:[0,\infty)\to\R^m}$ and both $y$ and $u$ are bounded with
    essentially bounded weak derivatives. The tracking error evolves uniformly
    within the performance funnel, i.e.
    \[
        \ex \eps>0\ \fa t>0:\quad \Norm{e(t)} \le \Funnel(t)^{-1} - \eps.
    \]
\end{prop}
In \Cref{Prop:FC_dist}, we will prove that a more advanced funnel controller design 
achieves the specified control objective for systems of the form \eqref{eq:Sys}. 
The simpler controller \eqref{Intro:eq:FC}, however, is a special case of the
controller methodology proposed in \cite{IlchRyan02b} and provides an intuitive
demonstration of the underlying principles.
When the tracking error ${e=y-y_{\rf}}$ approaches zero,
the control gain $k(t,y(t))$ diminishes, effectively deactivating the controller. 
Conversely, as the tracking error norm nears the funnel boundary $\Funnel(t)$,
the gain $k(t,y(t))$ grows rapidly, producing a control input $u(t)=-k(t,y(t)) e(t,y(t))$
that aggressively steers the error away from the funnel boundary and towards the reference trajectory.
Consequently, if the initial tracking error $e(0)$ lies within the funnel boundaries, 
its evolution remains strictly confined within the funnel $\cF_{\Funnel}$ for all time.

\subsection{Model predictive control}\label{Sec:OriginalMPCScheme}
The idea of model predictive control (MPC) is, 
after measuring/obtaining the output $\hat{y}\in\cC^{r-1}([\hat{t}-\tau,\hat{t}],\R^m)$ 
of the system~\eqref{eq:Sys} over a short time window of length $\tau\geq0$ at the current time $\hat{t}$ 
with $\hat{t}-\tau\geq t_0$, to repeatedly calculate a control function $u^\star=u^\star(\cdot;\hat{t},\hat{y})$ 
minimising the integral of a \emph{state cost}~$\ell$ on the future time interval
$[\hat{t},\hat{t}+T]$ for $T>0$, called the \emph{prediction horizon}, and
to implement the computed optimal solution~$u^{\star}$ to system~\eqref{eq:Sys} over an interval of
length $\delta<T$, called the \emph{time shift}.
The prediction horizon $T$  determines how far ahead the controller plans, while
the time shift $\delta$ specifies the implementation period before re-optimisation.
To make predictions about the future system behaviour and its output and, based on them, to compute optimal control signals,
MPC uses a model of the form
\begin{equation}\label{eq:Intro:ModelEquation}
    \begin{aligned}
       \yM^{(r)}(t)   &= \FM(\oTM(\yM,\ldots, \yM^{(r-1)})(t),u(t)),\\
       \yM|_{[\hat{t}-\tau,\hat{t}]}&= \yMh,
    \end{aligned}
\end{equation}
where $\FM \in \cC(\R^q\times \R^m,\R^m)$ is a \emph{known} continuous function
and $\oTM$ a \emph{known} operator, as a surrogate for the unknown system~\eqref{eq:Sys}.
The initial history function $\yMh$ of the model~\eqref{eq:Intro:ModelEquation} at time $\hat{t}\geq t_0$
is an element of $\cC^{r-1}([\hat{t}-\tau,\hat{t}],\R^m)$ and selected based on the system measurement $\hat{y}$,
though not necessarily identical to them.
Note that we assume that both the output and input dimension $m\in\N$ 
as well as the order of the differential equation $r\in\N$ match between the system~\eqref{eq:Sys} and the model~\eqref{eq:Intro:ModelEquation}.

\begin{remark}
    The model~\eqref{eq:Intro:ModelEquation} is intentionally formulated in a
    quite general form using an abstract operator~$\oTM$ and initial values
    $\yMh$ given on a time interval of length $\tau\geq0$ in order to explain
    the general idea of model predictive control while avoiding the accidental
    exclusion of particular edge cases.
    However, it is probably most common to use a control affine multi-input multi-output model of the form 
    \begin{equation}\label{eq:Intro:ModelControlAffine}
    \begin{aligned}
        \xMd(t)&=\fM(\xM(t))+\gM(\xM(t))u(t),\qquad \xM(\hat{t})=\xMh,\\
        \yM(t)&=\hM(\xM(t)),
    \end{aligned}
    \end{equation}
    with initial data $\xMh\in\R^n$ at the current time $\hat{t}\geq t_0$,
    and functions $f:\R^n\to\R^n$, $g:\R^n\to \R^{n\times m}$, and $h:\R^n\to\R^m$. 
    Here, $\xM(t)\in\R^n$ is the state of the model, $\yM(t)\in\R^m$ the model's output, and  $u(t)\in\R^m$ the control input.
    In \eqref{eq:Intro:ModelControlAffine}, the states $\xM$ of the model
    are laid out in an explicit way contrary to the
    formulation~\eqref{eq:Intro:ModelEquation} where they are, in a certain
    sense, hidden within the operator~$\oTM$.
    Many works on MPC use even simpler models, namely linear time invariant (LTI) models, i.e. 
    they assume a model of the form
    \begin{align*}
        \xMd(t)&=A_{\mathrm{M}}\xM(t)+B_{\mathrm{M}}u(t),\qquad \xM(\hat{t})=\xMh\\
        \yM(t)&=C_{\mathrm{M}}\xM(t),
    \end{align*}
    with $A_{\mathrm{M}} \in \R^{n \times n}$ and $C_{\mathrm{M}}^\top,B_{\mathrm{M}}\in\R^{n\times m}$.
    In the later part of this thesis, we will restrict ourselves also to a, in comparison to~\eqref{eq:Intro:ModelEquation}, 
    simpler model class, which  we will formally define in~\Cref{Def:ModelClass}.
    It will, however, contain the class of control affine multi-input
    multi-output models of the form~\eqref{eq:Intro:ModelControlAffine}
    (under certain assumptions on $\fM$, $\gM$, and $\hM$) and LTI models as we will
    see in~\Cref{Ex:LTISystem,Ex:ControlAffineMod}.
\end{remark}

In addition to the model~\eqref{eq:Intro:ModelEquation}, another key component of the MPC algorithm is the \emph{stage cost function} $\ell$.
It formulates rewards for desired model behaviour and penalties for undesired behaviour, which are balanced out during the optimisation process.
When solving the problem of tracking a given reference signal $y_{\rf}$, a commonly used stage cost function is 
\begin{equation}\label{eq:stageCostClassicalMPC}
    \begin{aligned}
        \ell:\Rp\times\R^m\times\R^{m}&\to\R,\qquad
        (t,\yM,u) \mapsto              \Norm{\yM-y_{\rf}(t)}^2+\lambda_u \Norm{u}^2
    \end{aligned}
\end{equation}
with $\lambda_u> 0$. While the term $\Norm{\yM-y_{\rf}(t)}^2$ penalises the distance of
the model's output~$\yM$ to the reference signal $y_{\rf}$, the term $\Norm{u}^2$ penalises the control
effort. The parameter~$\lambda_u$ allows to adjust a suitable trade-off between tracking performance and
required control effort. 
Of course, if a reference input signal~$u_{\rf}$ is known, the second
summand may be replaced by $\| u - u_{\rf}(t) \|^2$.

With the concepts introduced so far at hand, a general MPC algorithm can be formulated as follows. 
\begin{algo}[MPC]\label{Algo:MPC}\ \\
    \textbf{Given:} System~\eqref{eq:Sys}, model~\eqref{eq:Intro:ModelEquation},
    reference signal $y_{\rf}\in W^{r,\infty}(\Rp,\R^{m})$,
    initial time $t^0\in\Rp$, 
    boundary $\umax\geq0$ on the control input, and a stage cost function~$\ell$ as in \eqref{eq:stageCostClassicalMPC}.\\
    \textbf{Set} the time shift $\delta >0$, 
                 the prediction horizon $T\geq\delta$, and index $k\coloneqq 0$.\\
    \textbf{Define} the time sequence~$(t_k)_{k\in\N_0} $ by $t_k \coloneqq  t_0+k\delta$.\\ 
    \textbf{Steps:}
    \begin{enumerate}[(a)]
    \item\label{agostep:MPCFirst} Obtain a measurement of the state $y$
    at current time~$t_k$ over the last time period $[t_k-\tau,t_k]$ and set $\hat{y}^k \coloneqq y|_{[t_k-\tau,t_k]}$.
    Select an initial value  $\yMh^k$ for the model~\eqref{eq:Intro:ModelEquation} based on the measurement data $\hat{y}^k$.
    \item Compute a solution $u^{\star}\in L^\infty([t_k,t_k +T],\R^{m})$ of the \emph{Optimal Control Problem (OCP)}
    \begin{equation}\label{eq:MpcOCP}
    \begin{alignedat}{2}
            &\!\mathop
            {\operatorname{minimise}}_{u\in L^{\infty}([t_k,t_k+T],\R^{m})}  &\quad&
            \int_{t_k}^{t_k+T}\ell(t,\yM(t),u(t))\d t \\
            &\text{subject to}     &\yM^{(r)}(t)   &= \FM(\oTM(\yM,\ldots, \yM^{(r-1)})(t),u(t)),\\
            &                      &\yM|_{[\hat{t}-\tau,\hat{t}]}&= \yMh^k,\\
            &                      &\Norm{u(t)}  &\leq \umax.\\
    \end{alignedat}
    \end{equation}
    \item Apply the feedback law
        \[
            \mu:[t_k,t_{k+1})\times\cC^{r}([t_k-\tau,t_k],\R^m)\to\R^m, \quad \mu(t,\hat y) =u^{\star}(t)
        \]
        to system~\eqref{eq:Sys}.
        Increment $k$ by $1$ and go to Step~\ref{agostep:MPCFirst}.
    \end{enumerate}
\end{algo}
To ensure a bounded control signal with a maximal predefined control
value~$\umax\geq0$, the constraint~$\SNorm{u}\leq \umax$ has been added as an
additional constraint to the OCP~\eqref{eq:MpcOCP}. A key reason for the success
and popularity of the MPC~\Cref{Algo:MPC} is its ability to directly take
additional constraints into account by either adding them as hard constraints to
the OCP~\eqref{eq:MpcOCP}, as done with $\SNorm{u}\leq \umax$, or by including
soft constraints via penalty terms in the cost
function~$\ell$~\eqref{eq:stageCostClassicalMPC}.
For example, one could add the additional constraint
\begin{equation}\label{eq:ConstraintClassicalMPC} 
    \fa t\in[\hat{t},\hat{t}+T]:\quad \Norm{\yM(t)-y_{\rf}(t)} <\Funnel(t) 
\end{equation}
to the optimisation problem in order to guarantee that the model's output $\yM$ tracks the reference $y_{\rf}$
with prescribed performance, cf.~\cite{berger2019learningbased}.

We like to emphasise that both the system~\eqref{eq:Sys} to be controlled and 
the model~\eqref{eq:Intro:ModelEquation} used in the MPC~\Cref{Algo:MPC} 
are continuous-time (functional) differential equations.
While some other MPC approaches focusing on tracking problems 
also consider continuous-time systems, see e.g.~\cite{Facchino23}, 
most address discrete-time systems~\cite{aydiner2016periodic, limon2018nonlinear, kohler2019nonlinear, kohler2022constrained}.

\subsubsection{Difficulties and Drawbacks}
Although utilising the stage cost~$\ell$ in~\eqref{eq:stageCostClassicalMPC} and
constraints~\eqref{eq:ConstraintClassicalMPC} in \Cref{Algo:MPC} might seem like a canonical
choice when solving the reference tracking problem with MPC, this approach has several drawbacks.
In particular, one has to guarantee initial and recursive feasibility of the MPC
Algorithm~\ref{Algo:MPC}. This means it is necessary to prove that the optimisation
problem~\eqref{eq:MpcOCP} has initially (i.e. at $t_k = t_0$) and recursively (i.e. at $t_k=t_0+\delta
k$ after $k$ steps of Algorithm~\ref{Algo:MPC}) a solution.
First of all, one has to show existence of an $L^\infty$-control $u$ bounded by $\umax\geq0$ which,
if applied to a model of type~\eqref{eq:Intro:ModelEquation}, guarantees that the model's tracking error~
\begin{equation}\label{eq:Intro:ModelTrackingerror}
    \eMTrack(t)\coloneqq \yM(t)-y_{\rf}(t)
\end{equation}
evolves within the performance funnel~$\cF_{\Funnel}$ given by $\Funnel\in\cG$, i.e. fulfilling~\eqref{eq:ConstraintClassicalMPC}.
Furthermore, one has to prove that there exists a solution~$u^{\star}$ of the
optimisation problem~\eqref{eq:MpcOCP} and that this solution
fulfils~\eqref{eq:ConstraintClassicalMPC}. To show recursive feasibility, it is
further necessary to ensure that after applying a solution~$u^{\star}$ of the
optimal control problem~\eqref{eq:MpcOCP} at time $t_k=t^0+\delta k$ to the
system~\eqref{eq:Sys} the optimisation problem is still well defined at the next
time instant $t_{k+1}=t^0+\delta (k+1)$ when re-initialised with $\yMh^k$ based
on new measurements from the system. We will discuss potential methods to guarantee 
initial and recursive feasibility and their advantages and disadvantages in more
detail in \Cref{Chapter:FunnelMPC}.

Both ensuring feasibility and achieving the control objective are already
challenging problems when one assumes the model to coincide with the system.
However, since every model, no matter how good, deviates from the actual system
and disturbances are omnipresent, we assume the
model~\eqref{eq:Intro:ModelEquation} to differ from the system~\eqref{eq:Sys}
rendering the problem even more demanding. One has to account for discrepancies
between the model predictions~$\yM(t)$ and the actual system output~$y(t)$, i.e.
the model-plant output mismatch
\begin{equation}\label{eq:Intro:ModelPlantMismatch}
        \eSTrack(t)  \coloneqq  y(t)-\yM(t).
\end{equation}
In the presented form, the MPC~\Cref{Algo:MPC} with
constraints~\eqref{eq:ConstraintClassicalMPC} can only achieve that the model's
tracking error $\eMTrack$ evolves within the funnel $\cF_{\Funnel}$. However,
additional robustification methods need to be utilised in order to compensate
for the model-plant output mismatch $\eSTrack$ and achieve the control objective
for the actual tracking error $e$ as laid out in~\Cref{Sec:ControlObjective}.
One aspect that merits particular attention is the selection of the initial value $\yMh^k$ 
for the model based on the system measurement $\hat{y}^k$ in Step~\ref{agostep:MPCFirst} of~\Cref{Algo:MPC}.
To reduce any occurring model-plant mismatch, one ideally would like to initialise the model~\eqref{eq:Intro:ModelEquation} 
directly with the system measurement $\hat{y}^k$.
This, however, might result in a constraint violation of applied restrictions like~\eqref{eq:ConstraintClassicalMPC} 
and thus render the optimal control problem~\eqref{eq:MpcOCP} unsolvable.
Although undesirable, it may be necessary to allow for a higher tolerance of the model-plant mismatch in order to mathematically guarantee 
feasibility and functioning of the control scheme.
Although it may seem evident, we would like to point out that initialising the model~\eqref{eq:Intro:ModelEquation}
with $\yMh^k$ based on measurement data $\hat{y}^k$ will result in the model's solution on the interval $[t_k,t_{k+1}]$ 
not being a continuous extension of the solution on the previous interval $[t_{k-1},t_{k}]$.
Although these individual solutions fulfil the differential equation~\eqref{eq:Intro:ModelEquation},
the entire trajectory is, in general, not a solution.
To refer to this trajectory of the concatenated solutions of the model's differential equation~\eqref{eq:Intro:ModelEquation},
we will use the term~\emph{concatenated solution} of the MPC algorithm.
This will be made mathematically precise in~\Cref{Def:SolutionClosedLoop}.
It is clear that the concatenated solution can jump at the time instants $t_k$ because of the initialisation with $\yMh^k$. 
Thus, it is not even a continuous function but merely a regulated function.
In order to avoid any potential ambiguities, we want to briefly recall the definition.
\begin{definition}[Regulated function]\label{Def:RegulatedFunction}
    On an interval $I$, we call a function $f:I\to\R^n$ a \emph{regulated function}, if 
    the left and right limits $f(t-)$ and $f(t+)$ exist for all interior points
    $t\in I$ and $f(a-)$ and $f(b+)$ exist whenever $a= \inf I \in I$ or $b=\sup I \in I$. 
    We denote the space of all regulated functions on~$I$ by~$\cR(I,\R^n)$.
\end{definition}

\section{Structure of this thesis and previously published results}
This thesis is subdivided into four parts. Every chapter focuses on one of the
laid-out difficulties related to solving the control objective with
model predictive control.

\Cref{Chapter:FunnelMPC} presents the theoretical foundations for solving the
output tracking problem as described in~\Cref{Sec:ControlObjective} using a
dedicated MPC algorithm. To this end, it is initially assumed that no
disturbances are present and that the model matches the system perfectly.
\Cref{Sec:FunnelStageCostFunctions} introduces the general concept of funnel
stage cost functions. Incorporating these into the MPC's optimal control problem
guarantees that the output tracking error evolves within prescribed performance
boundaries. This guarantee is rigorously proven
in~\Cref{Sec:MPCWithFunnelStageCost}, following the introduction of the model
class in~\Cref{Sec:ModelClass}. The resulting funnel MPC~\Cref{Algo:FunnelMPC}
is defined in~\Cref{Sec:FunnelMPC}, with the chapter's main result -- the proof
of its initial and recursive feasibility -- being established in~\Cref{Thm:FMPC}.
This proof leverages auxiliary error variables (introduced
in~\Cref{Sec:HighRelativeDegree}) to ensure satisfaction of the control
objective independently of the model order. It further relies on two key prior
results: \Cref{Th:ExUmax}, proving the existence of control functions confining
the error within funnel boundaries, and \Cref{Th:SolutionExists}, showing the
solvability of the optimal control problem with funnel stage costs.
While building upon prior works~\cite{BergDenn21,BergDenn22,BergDenn24}, this
chapter significantly extends them by formalising the general concept of funnel
stage cost functions and providing detailed proofs for all mathematical aspects,
addressing omissions due to page limitations in earlier publications.
Furthermore, the auxiliary error framework developed here enabled the proposal
of a low-complexity funnel controller for higher-order non-linear systems
in~\cite{Dennst25}.

\Cref{Chapter:RobustFunnelMPC} lifts the standing assumptions made in the
previous chapter, namely the absence of disturbances and the perfect
model-system alignment. Building upon the ideas from~\cite{BergDenn24b}, the
funnel MPC algorithm is robustified by incorporating the funnel controller as a
second controller component. After introducing the system class under
consideration in~\Cref{Sec:SystemClass}, the structure of the two component
control scheme is presented in~\Cref{Sec:RobustFMPCControllerStructure}. The
main result is~\Cref{Thm:RobustFMPC} showing that the robust funnel
MPC~\Cref{Algo:RobustFMPC} achieves the control objective in presence of
disturbances and even a structural system model mismatch.
\Cref{Chapter:RobustFunnelMPC} extends the results from~\cite{BergDenn24b} to
encompass systems and models of higher order with non-linear time delays and
potentially infinite-dimensional internal dynamics.

The control scheme is further extended by a (machine) learning component
in~\Cref{Chapter:LearningRobustFMPC}.
The structure of this three-component controller is laid out
in~\Cref{Sec:StructureLearningFMPC}. While the funnel controller component
mitigates model-plant mismatches, bounded disturbances, and uncertainties, the
machine learning component adapts the underlying model to the system data and,
thus, improves the contribution of the MPC component over time. 
\Cref{Def:LearningScheme} summarises the structural assumptions on a learning algorithm 
that ensure the successful interplay of the three components. 
The main result \Cref{Thm:learningFRMPC} shows the functioning of the
three-component controller and that this control~\Cref{Algo:LearningRFMPC}
achieves the control objective as laid out in \Cref{Sec:ControlObjective}.
To illustrate the abstract requirements and assumptions, a possible learning
approach is discussed in \Cref{Sec:LearningLinearModel}. The chapter builds upon
the work~\cite{LanzaDenn24b} and extends its results to the model and system
class of functional differential equations of arbitrary order which was
considered in the previous chapters.

In \Cref{Chapter:FunnelMPC,Chapter:RobustFunnelMPC,Chapter:LearningRobustFMPC},
it is assumed that the system output can be continuously measured and that an
arbitrary measurable control signal can be applied to the system.
\Cref{Chapter:DiscretFMPC} lifts this assumption and shows that  the robust funnel MPC~\Cref{Algo:RobustFMPC} 
from \Cref{Chapter:RobustFunnelMPC} can be modified to achieve the control objective with sampled-data control.
In line with the two-component structure of the controller, the chapter is divided into two parts,
each of which shows that the respective component can be designed in a sampled data manner.
\Cref{Sec:FCZoH} is based on the work~\cite{LanzaDenn24} and proposes a funnel controller that 
achieves the control objective while only measuring the system output at discrete time instants and 
only applying piecewise constant control signals.
Uniform bounds on the required sampling rate and the maximal applied control signal are derived.
These results are summed up in \Cref{Thm:DiscreteFC}.
As a small extension, \Cref{Sec:FCAsSaftyFilter} shows how this controller can also be used as 
a safety filter for other data-driven control approaches.
In \Cref{Sec:FMCPZoH}, the results are carried over to the funnel MPC scheme.
\Cref{Thm:DiscreteFMPC} shows that the sampled funnel MPC~\Cref{Algo:DiscrFunnelMPC},
which only applies piecewise constant control signals to the model, achieves the control objective.
This section builds upon the work~\cite{Dennst24} and extends previously published results 
to higher order models.

Within the course of this dissertation, the scientific article \cite{Oppeneiger24} was also published 
in addition to the mentioned works.
In that paper, a mathematical model for a magnetic levitation train is
developed. The funnel MPC algorithm from~\cite{BergDenn21} is then applied to
the system with the objective of guaranteeing the safe and dependable operation
by ensuring that the distance between the magnet and the reaction rail is kept
within a given range.
The control scheme is then compared to two different control approaches, one
linear state feedback controller and a model predictive control scheme with a
quadratic cost function, with respect to performance criteria such as
robustness, travel comfort, control effort, and computation time in an extensive
numerical simulation study. As this dissertation is mainly concerned with the
underlying mathematical theory of the funnel MPC algorithm and the work
\cite{Oppeneiger24} focuses on the modelling aspect of the magnetic levitation systems,
we will not present a separate evaluation of the results from~\cite{Oppeneiger24} in this thesis.

\chapter{Funnel Model Predictive Control}\label{Chapter:FunnelMPC}

Guaranteeing \emph{initial and recursive feasibility} is essential for the
successful application of MPC and it is one of the major challenges. This requires
guaranteeing that the optimal control problem is solvable at the initial time
step and that solvability at any subsequent step follows recursively. While
initial feasibility is often simply assumed, a common strategy to achieve
recursive feasibility involves augmenting the OCP with carefully designed
terminal conditions -- such as terminal costs and constraints -- as discussed in
~\cite{chen1998quasi} and~\cite{rawlings2017model}. However, these artificially
imposed terminal conditions introduce two key challenges: they raise the
computational complexity of solving the OCP and complicate the identification of
an initially feasible solution. Consequently, the domain of attraction of the
MPC controller may be significantly restricted, as noted in
\cite{chen2003terminal, gonzalez2009enlarging}. Furthermore, designing such
conditions becomes markedly more intricate under time-varying state constraints
\cite{manrique2014mpc}.
An alternative approach circumvents terminal conditions
by leveraging cost controllability principles
\cite{CoroGrun20,Tuna2006,worthmann2011stability} and employing a sufficiently
long prediction horizon, see \cite{boccia2014stability} and reference therein or
\cite{EsteWort21} for an extension to continuous-time systems. Notably, both
terminal-condition-based and horizon-based techniques face heightened complexity
when applied to systems with time-varying state or output constraints,
underscoring the need for tailored solutions in such cases. Additionally,
feasibility assurances depend critically on the considered system class:
non-linear or uncertain systems often necessitate more conservative designs,
such as tube-based methods \cite{MaynSero05} or adaptive mechanisms for model
refinement~\cite{HewingWaber20}. For systems with periodic constraints or
references (e.g. tracking periodic trajectories), feasibility analysis often
requires periodicity-aware terminal sets or horizon lengths
\cite{limon2014periodic}. Practical implementations must also contend with
computational limits, where inexact solvers or early termination can undermine
theoretical guarantees \cite{zanelli2021inexact}. For systems with persistent
infeasibility, slack variables or softened constraints may be introduced, albeit
at the expense of performance~\cite{kerrigan2000soft}. Lastly, in economic MPC,
where stability is not the primary objective, feasibility frameworks must
reconcile transient constraints with long-term economic goals~\cite{Amrit2009}.

To overcome the restrictions of mentioned methods to ensure initial and recursive feasibility,
\emph{funnel MPC (FMPC)} was originally proposed
in~\cite{berger2019learningbased} and then further developed in~\cite{BergDenn21,BergDenn22,BergDenn24}.
It is an MPC scheme that allows for output tracking of an
a priori given reference signal $y_{\rf}\in W^{r,\infty}(\Rp,\R^{m})$ such that the
tracking error evolves in a pre-specified, potentially time-varying, performance
funnel given by a function $\Funnel\in\cG$.
The core idea involves replacing $\ell$ from~\eqref{eq:stageCostClassicalMPC} 
in the MPC~\Cref{Algo:MPC} with a novel funnel stage cost function $\ell_{\Funnel}: \Rp
\times \R^m \times \R^m \to \R\cup{\infty}$. 
This function, parametrised by $\lambda_u \in \Rp$, is defined as:
\begin{align}\label{eq:stageCostFunnelMPC}
    \ell_{\Funnel}(t,y,u) =
    \begin{dcases}
        \frac {\Norm{y-y_{\rf}(t)}^2}{\Funnel(t)^2 - \Norm{y-y_{\rf}(t)}^2} + \lambda_u \Norm{u}^2,
            & \Norm{y-y_{\rf}(t)} \neq \Funnel(t),\\
        \infty,&\text{else}.
    \end{dcases}
\end{align}
The term~$\frac {\Norm{y-y_{\rf}(t)}^{2}}{\Funnel(t)^2 -\Norm{y-y_{\rf}(t)}^2}$
penalises the proximity of the tracking error to the funnel boundary, whereas the 
term $\Norm{u}$ serves as a penalisation of the control input as in~\eqref{eq:stageCostClassicalMPC}.  
The  parameter $\lambda_u$ can be used to adjust the balance between these two control objectives.
The cost function~$\ell_{\Funnel}$ is motivated by the results on funnel control which we briefly introduced in~\Cref{Sec:FunnelControl}.
Contrary to the cost function $\ell$ as in~\eqref{eq:stageCostClassicalMPC}, this ``funnel-like'' stage cost does not directly penalise 
the norm of the tracking error  but its proximity to the funnel boundary~$\Funnel$ 
and grows unbounded when the error norm approaches $\Funnel$.

While many reference tracking MPC approaches focus on ensuring the asymptotic
stability of the tracking error -- often through terminal constraints or costs
-- they generally do not enforce strict boundaries on the output signal. For
instance, \cite{aydiner2016periodic, koehler2020} guarantee stability of the
tracking error by designing terminal sets and costs around specific reference
trajectories, whereas \cite{kohler2019nonlinear} achieves this by relying on a
sufficiently long prediction horizon instead of terminal constraints. 
Similarly, \cite{kohler2022constrained} ensures constraint satisfaction for
discrete-time systems through stabilisability and detectability assumptions and
long horizons.
At first glance, tube-based MPC schemes, see e.g. \cite{Limon2008, Limon2010},
share similarities with funnel MPC, as they confine the tracking error to a
controllable and potentially time-varying range. However, their primary goal is
to compensate for model uncertainties and disturbances acting on the system by
constructing tubes around the reference trajectory. These tubes are often
offline-computed and are not user-definable as they must inherently account for
system uncertainties; for example \cite{Lopez2019} dynamically optimises
both tubes and reference trajectories  based on proximity to tube boundaries.

\emph{Control barrier functions (CBFs)} determine the control input and enforce
safety-critical constraints by solving a quadratic program (QP) that directly
regulates the derivative of the barrier function. This ensures that the system
remains within a safe set by design, guaranteeing positive invariance and
asymptotic stability without requiring predictive optimisation. Due to their
simplicity and modularity, CBFs are widely adopted in robotics,
see~\cite{Ames2017,ames2019control} for an overview. Building on this concept,
barrier function-based MPC integrates barrier functions into the MPC framework,
e.g. \cite{WILLS20041415, Marvi2019,Pfitz2021}. Unlike stand-alone CBFs, this approach
incorporates a barrier term directly into the MPC cost function, penalising
proximity to constraint boundaries over a prediction horizon. Alternatively,
control Lyapunov-barrier functions are incorporated in the MPC scheme to ensure
recursively feasibility and stabilisation of the closed-loop system, see e.g.
\cite{Wu2019}. However, similar to classical MPC, terminal costs or constraints
are often still required to ensure recursive feasibility and constraint
satisfaction. In contrast, funnel MPC eliminates the need for terminal
constraints or costs by employing a unique cost function that diverges as the system
output approaches the funnel boundary. It  circumvents the reliance on terminal
conditions and avoids the need for a sufficiently long prediction horizon
entirely.

In this chapter, we will analyse how the utilisation of stage cost functions
such as~\eqref{eq:stageCostFunnelMPC} in the MPC~\Cref{Algo:MPC} (which we 
will then call \emph{funnel MPC}) ensures fulfilment of the control objective: the
tracking of a given reference signal~$y_{\rf}$ with the model's output~$\yM$
within predefined funnel boundaries~$\Funnel$. For a large class of models,
including but not limited to models with non-linear time delays and potentially
infinite-dimensional internal dynamics, we will rigorously prove initial and
recursive feasibility of this MPC scheme. This will be achieved without
incorporating additional constraints in the optimal control
problem~\eqref{eq:MpcOCP}, without imposing additional terminal conditions, and
independent of the length of the prediction horizon $T>0$. For our analysis
in this chapter, we will assume that the system~\eqref{eq:Sys} and the surrogate
model~\eqref{eq:Intro:ModelEquation} coincide. In particular, we assume that
the model-plant mismatch~$\eSTrack$ defined
in~\eqref{eq:Intro:ModelPlantMismatch} is identically zero. 
Both assumptions will be relaxed in the later parts of this thesis.

\section{Funnel stage cost functions}\label{Sec:FunnelStageCostFunctions}
The key distinction between funnel MPC and the classical MPC~\Cref{Algo:MPC}
lies in how their respective stage costs penalise the tracking error
$e(t)=y(t)-y_{\rf}(t)$. While the classical stage cost function $\ell$
in~\eqref{eq:stageCostClassicalMPC} directly penalises the squared norm of the
tracking error $e$, the funnel MPC stage cost $\ell_{\Funnel}$
in~\eqref{eq:stageCostFunnelMPC} imposes a particular penalty tied to the
proximity of $e$ to the funnel boundary $\Funnel$.
This raises a natural question: how does this modified cost function ensure that
the tracking error~$e$ evolves within the funnel~$\cF_{\Funnel}$ defined by a
function $\psi\in\cG$ when a solution of the optimisation
problem~\eqref{eq:MpcOCP} is applied to the
model~\eqref{eq:Intro:ModelEquation}?  
If the initial error lies within the funnel, without explicit constraints in the
optimal control problem~\eqref{eq:MpcOCP}, then the error $e$ could theoretically
still touch or even exceed the boundary and evolve outside of the funnel
boundary after some time.
The previous work~\cite{berger2019learningbased} addressed this issue by
enforcing explicit hard state constraints of the
form~\eqref{eq:ConstraintClassicalMPC}.
In contrast, we will show that such constraints are unnecessary.
Instead, the specific structure of the cost function $\ell_{\Funnel}$
in~\eqref{eq:stageCostFunnelMPC} inherently enforces compliance with the
time-varying funnel boundaries. To analyse this mechanism of implicit
constraints, we will, in this section, examine the
function~$\OrigFunnelPenaltyFunc:\Rp\times\R^n$ defined as 
\begin{align}\label{eq:PenaltyOriginalFMPC}
    \OrigFunnelPenaltyFunc(t,e) =
    \begin{dcases}
        \frac {\Norm{e}^2}{\Funnel(t)^2-\Norm{e}^2}
            & \Norm{e} \neq \Funnel(t)\\
        \infty,&\text{else},
    \end{dcases}
\end{align}
which quantifies the ``cost'' of being close to the funnel boundary.
To isolate its essential features from any specific system dynamics or control
problem, we examine, for a given $\hat{t}\in\Rp$ and $T>0$, the integral 
\begin{equation}\label{eq:OriginalCostFunction}
    \int_{\hat{t}}^{\hat{t}+T}\OrigFunnelPenaltyFunc(s,\Lpath(s))\d{s}
\end{equation}
evaluated along an arbitrary \emph{Lipschitz path}
$\Lpath\in\Lip([\hat{t},\hat{t}+T],\R^n)$, i.e. a Lipschitz continuous
function~$\Lpath$ defined on the interval $[\hat{t},\hat{t}+T]$ with values in
$\R^n$. Identifying the essential aspects of the function
$\OrigFunnelPenaltyFunc$ will lead to the \Cref{Def:FunnelStageCostFunc} of
\emph{funnel penalty functions}. These penalty functions, when used as a stage
cost in the optimal control problem~\eqref{eq:MpcOCP} within the
MPC~\Cref{Algo:MPC}, ensure the tracking error $e$ evolves within the
funnel~$\cF_{\Funnel}$ given by a function $\Funnel\in\cG$.
Restricting our analysis to Lipschitz paths for now allows us 
to avoid technical complications and to 
focus on the essential characteristics of the penalty function
$\OrigFunnelPenaltyFunc$ and the associated cost function.
This simplification is justified because both the model output~$\yM$ in
\eqref{eq:Intro:ModelEquation} and the reference trajectory~$y_{\rf}\in
W^{r,\infty}(\Rp,\R^m)$ are differentiable functions; hence 
their restrictions to any compact interval $[\hat{t},\hat{t}+T]$ are Lipschitz
continuous (provided $\yMd$ is bounded). 
This assumption will later be validated in~\Cref{Prop:SolutionIsLPath,Prop:ErrorIsLPath}.

Before proceeding, note a subtle point regarding the interpretation of the integral~\eqref{eq:OriginalCostFunction}.
There is a distinction between a function that is Lebesgue integrable
(i.e. it belongs to~$L^1$) and a function for which the Lebesgue integral merely
exists but which does not necessarily have to be an element of $L^1$.
To make this difference clearer, we call a measurable function $\zeta : B \to \R$ on a Borel set
$B\subseteq\R$ \emph{quasi-integrable} if at least one of the Lebesgue integrals
\[
    \int_B \zeta^+(t) \d t\qquad \text{or} \qquad \int_B \zeta^-(t) \d t
\]
(with $\zeta^+ \coloneqq  \max\{\zeta,0\}$ and $\zeta^-\coloneqq \max\{-\zeta,0\}$) is finite. 
If both integrals diverge, the overall integral is defined to be infinity.
In particular, if $\OrigFunnelPenaltyFunc(\cdot,\Lpath(\cdot))$ is not quasi-integrable
for a given Lipschitz path $\Lpath\in\Lip([\hat{t},\hat{t}+T],\R^n)$,
then the Lebesgue integral in~\eqref{eq:OriginalCostFunction} is treated as infinity. 
Moreover, it may happen that $\Funnel(t) = \Norm{\Lpath(t)}$ for some
$t\in[\hat{t},\hat{t}+T]$ and with that
${\OrigFunnelPenaltyFunc(t,\Lpath(t))=\infty}$.
If the set of such points does not have Lebesgue measure zero, then the
integral~\eqref{eq:OriginalCostFunction} is infinity as well. In the following,
we will prove that if the integral~\eqref{eq:OriginalCostFunction} is finite,
i.e. $\OrigFunnelPenaltyFunc(\cdot,\Lpath(\cdot))$ is quasi-integrable over
$[\hat{t},\hat{t}+T]$, and the integral does not diverge, then the Lipschitz path $\Lpath$ 
must evolve entirely within the funnel~$\cF_{\Funnel}$. To
show this, an elementary lemma is proved first.

\begin{lemma}\label{Lem:PosLipschitzCont}
   Let $T>0$, $\hat{t}\in\Rp$ and $\Lpath\in\Lip([\hat{t},\hat{t}+T],\Rp)$ be a Lipschitz path.
   If~${\int_{\hat{t}}^{\hat{t}+T}\frac{1}{\Lpath(s)}\d s<\infty}$, then $\Lpath(s)>0$ for all $s\in[\hat{t},\hat{t}+T]$.
\end{lemma}
\begin{proof}
    Assume that there exists $t\in(\hat{t},\hat{t}+T)$ such that $\Lpath(t)=0$. Choose $\eps>0$ such that
    ${(t-\eps,t+\eps)\subset [\hat{t},\hat{t}+T]}$. Since $\Lpath$ is Lipschitz continuous, we have  that
    \[
      \ex C>0\ \fa s\in(t-\eps,t+\eps):\  \Lpath(s) =     \Abs{ \Lpath(s)-\Lpath(t)} \leq  C\Abs{ s-t      }.
    \]
    Therefore,
    \begin{align*}
        \infty>\int_{0}^{T}\frac{1}{\Lpath(s)}\d s
        \geq \int_{t-\eps}^{t+\eps}\frac{1}{\Lpath(s)}\d s
        \geq \int_{t-\eps}^{t+\eps}\frac{1}{C\Abs{s-t}}\d s
        =    \int_{-\eps}^{\eps}\frac{1}{C\Abs{s}}\d s
        =    \infty,
    \end{align*}
    a contradiction. A similar proof applies in the cases $t=\hat{t}$ and $t=\hat{t}+T$.
\end{proof}

\begin{remark}
    \Cref{Lem:PosLipschitzCont} is not true for all uniformly
    continuous functions in general. Consider the example:
    \[
        \int_{0}^1\frac{1}{\sqrt{x}}\d x
            = 2\sqrt{x}\Big|_0^1
            = 2.
    \]
\end{remark}
\begin{prop}\label{Prop:LpathGraph}
For $\Funnel\in\cG$, $\hat{t}\in\Rp$ and $T>0$, let 
$\Lpath\in\Lip([\hat{t},\hat{t}+T],\R^n)$ be a Lipschitz path with $(\hat{t},\Lpath(\hat{t}))\in\cF_{\Funnel}$.
Then,
\[
    \int_{\hat{t}}^{\hat{t}+T}\OrigFunnelPenaltyFunc(s,\Lpath(s))\d{s}<\infty\iff\graph(\Lpath)\subset\cF_{\Funnel}.
\]
\end{prop}
\begin{proof}
We show the two implications separately.

    \noindent 
    \emph{``$\impl$'':} We show that $\Norm{\Lpath(s)}<\Funnel(s)$ for all $s\in[\hat{t},\hat{t}+T]$.
    Since $(\hat{t},\Lpath(\hat{t}))\in\cF_{\Funnel}$, we have $\Norm{\Lpath(\hat{t})}<\Funnel(\hat{t})$. 
    Assume there exists $s\in [\hat{t},\hat{t}+T]$ with $\Norm{\Lpath(s)}\geq\Funnel(s)$.
    By continuity of $\Lpath$ and $\Funnel$, there exists 
    \[
        \hat{s}\coloneqq \min \setdef{t\in[\hat{t},\hat{t}+T]}{\Norm{\Lpath(t)}=\Funnel(t)}.
    \]
    Note that $\Norm{\Lpath(s)}< \Funnel(s)$ for all $s\in[\hat{t},\hat{s})$.
    Recalling the definition of the Lebesgue integral, see e.g.~\cite[Def 11.22]{Rudi76},
    $\int_{\hat{t}}^{\hat{t} + T}\OrigFunnelPenaltyFunc(s,\Lpath(s))\, \d{s}< \infty$ implies
    $\int_{\hat{t}}^{\hat{t} + T}\rbl \OrigFunnelPenaltyFunc(s,\Lpath(s))\rbr^+ \d{s} < \infty$. 
    Thus,
    \begin{align*}
       \int_{\hat{t}}^{\hat s}    \frac{1}{1-\frac{\Norm{\Lpath(s)}^2}{\Funnel(s)^2}}            \d{s}
       &= \int_{\hat{t}}^{\hat s}    \frac {\Norm{\Lpath(s)}^2}{\Funnel(s)^2  - \Norm{\Lpath(s)}^2} + 1           \d{s} \\
       & \leq  \int_{\hat{t}}^{\hat{t}+T}\rbl\frac {\Norm{\Lpath(s)}^2}{\Funnel(t)^2  - \Norm{\Lpath(s)}^2}\rbr^{+}    \d{s} + T\\
       &\leq      \int_{\hat{t}}^{\hat{t}+T}\rbl \OrigFunnelPenaltyFunc(s,\Lpath(s))\rbr^{+} \d{s} + T < \infty.
    \end{align*}
    Both the path $\Lpath$ and the funnel function $\Funnel$, being an element of $W^{1,\infty}(\Rp,\R),$ are Lip\-schitz continuous functions. 
    Since products and sums of Lipschitz continuous functions on a compact interval are again Lipschitz continuous,
    we may infer that $1-\frac{\Norm{\Lpath(\cdot)}^2}{\Funnel(\cdot)^2}$ is  Lipschitz continuous on $[\hat t, \hat s]$.
    By definition of $\hat{s}$, it moreover is non-negative. 
    Now,~\Cref{Lem:PosLipschitzCont} yields that it is strictly positive, i.e. $\Funnel(s)^2  > \Norm{\Lpath(s)}^2$
    for all $s\in[\hat t, \hat s]$, which contradicts the definition of~$\hat{s}$.
    
    \noindent 
    \emph{``$\Leftarrow$'':}
    Since $\graph(\Lpath)\subset\cF_{\Funnel}$, we have $\Norm{\Lpath(s)} < \Funnel(s)$ for all $s\in[\hat{t},\hat{t}+T]$.
    Due to continuity of the involved functions and the compactness of the interval $[\hat{t},\hat{t}+T]$,
    there exists $\eps\in(0,1)$ with $\Norm{\Lpath(s)} \le \Funnel(s) - \eps$ for all $s \in [\hat{t},\hat{t}+T]$.
    Then, $\OrigFunnelPenaltyFunc(s,\Lpath(s)) \ge 0$ for all $s \in [\hat{t},\hat{t}+T]$ and
    \begin{align*}
               \int_{\hat{t}}^{\hat{t}+T}        \Abs{ \OrigFunnelPenaltyFunc(s,\Lpath(s))}\d{s}
        =      \int_{\hat{t}}^{\hat{t}+T}\Abs{ \frac{\Norm{\Lpath(s)}^2}{\Funnel(s)^2-\Norm{\Lpath(s)}^2}}\d{s}
        \leq   \int_{\hat{t}}^{\hat{t}+T} \frac{\|\Funnel\|_\infty}{\eps} \d{s} 
        =   T\frac{\SNorm{\Funnel}}{\eps}  <      \infty.
    \end{align*}
    This completes the proof.
\end{proof}

\begin{remark}
    In contrast to funnel MPC, \emph{barrier function based MPC}, see
    e.g.~\cite{WILLS20041415, Marvi2019}, employs (relaxed) logarithmic barrier
    functions to penalise states near constraint boundaries. This might
    initially seem to be merely a subtle difference since both methods involve
    stage cost functions that grow unbounded as the state approaches the
    constraint boundaries. However, the distinction has significant theoretical
    implications. The results in \Cref{Lem:PosLipschitzCont} and, consequently,
    \Cref{Prop:LpathGraph}, arise from the non-integrability of $x\mapsto\tfrac{1}{x}$
    over the interval~$[0,1]$. Specifically, \Cref{Prop:LpathGraph} asserts that
    a finite value of the integral in~\eqref{eq:OriginalCostFunction} guarantees
    that any Lipschitz path starting within $\cF_{\Funnel}$ remains confined to
    the prescribed funnel boundaries. Consequently, when the stage cost function
    \eqref{eq:stageCostFunnelMPC} is used in the optimal control
    problem~\eqref{eq:MpcOCP} (within the MPC~\Cref{Algo:MPC}), the tracking
    error~$e\coloneqq y-y_{\text{ref}}$ is ensured to evolve within the funnel
    boundaries defined by $\Funnel\in\cG$.
    In contrast, a logarithmic barrier function is integrable over the
    interval~$[0,1]$:
    \[
        \int_{0}^{1}\ln(x^n)\text{d} x= 
        \ln(x^n)x\Big\vert_{0}^{1} -\int_{0}^{1}x\frac{n}{x}\text{d} x=0-n=-n.
    \]
    This integrability implies that logarithmic penalities alone can, in
    general, not guarantee that a Lipschitz path or, in the context of MPC, the model state
    always remain  within the desired region. 
    As a result, the usage of terminal conditions (costs and constraints)
    remains essential in the optimal control problem~\eqref{eq:MpcOCP}.
\end{remark}

\Cref{Prop:LpathGraph} establishes that the
integral~\eqref{eq:OriginalCostFunction} is finite if and only if the Lipschitz
path $\Lpath\in\Lip([\hat{t},\hat{t}+T],\R^n)$ evolves entirely within the
funnel~$\cF_{\Funnel}$, provided $\Lpath$ starts within the funnel.
Translating this result to the stage cost $\ell_{\Funnel}$ in~\eqref{eq:stageCostFunnelMPC} (to be used in the optimal control
problem~\eqref{eq:MpcOCP}), a finite stage cost guarantees that the model's
tracking error $e=\yM-y_{\rf}$ remains within the funnel boundaries defined
by~$\Funnel$.
However, the non-linearity and discontinuity of the function $\ell_{\Funnel}$
in~\eqref{eq:stageCostFunnelMPC} raise concerns about the solvability of the
optimal control problem~\eqref{eq:MpcOCP} when employing this stage cost
function. Moreover, even if a solution exists, additional analysis is required
to ensure that its solution also guarantees the evolution of the tracking error
within the funnel boundaries.
To address this in \Cref{Th:SolutionExists}, we will construct a sequence of control
functions converging to the infimum of the minimisation
problem~\eqref{eq:MpcOCP} and analyse the corresponding sequence of error
trajectories. We will then invoke the following~\Cref{Lem:OriginalCostFuncLimit}
to prove that if all trajectories in this sequence remain within~$\cF_{\Funnel}$, 
then their limit will also remain within the funnel boundaries.

\begin{lemma}\label{Lem:OriginalCostFuncLimit}
Let $\Funnel\in\cG$, $\hat{t}\in\Rp$, $T>0$ and $\Lpath^\star\in\Lip([\hat{t},\hat{t}+T],\R^n)$ be a Lipschitz path with
$(\hat{t},\Lpath^\star(\hat{t}))\in\cF_{\Funnel}$.
Further, let  $(\Lpath_n)\in\Lip([\hat{t},\hat{t}+T],\R^n)^\N$ be a sequence of
Lipschitz paths  with $(\hat{t},\Lpath_n(\hat{t}))\in\cF_{\Funnel}$ for all
$n\in\N$ that converges uniformly to $\Lpath^\star$.
If $\int_{\hat{t}}^{\hat{t}+T}\OrigFunnelPenaltyFunc(s,\Lpath_n(s))\d{s}$ is uniformly bounded by some constant $M\geq0$, then 
\[
    \int_{\hat{t}}^{\hat{t}+T}\OrigFunnelPenaltyFunc(s,\Lpath^\star(s))\d{s}<\infty.
\]
\end{lemma}
\begin{proof}
    It suffices to show that $\Norm{\Lpath^\star(s)}<\Funnel(s)$ for all ${s\in[\hat{t},\hat{t}+T]}$ according to \Cref{Prop:LpathGraph}.
    Since $(\hat{t},\Lpath^\star(\hat{t}))\in\cF_{\Funnel}$, we have $\Norm{\Lpath^\star(\hat{t})}<\Funnel(\hat{t})$. 
    Assume there exists $s\in [\hat{t},\hat{t}+T]$ with $\Norm{\Lpath^\star(s)}\geq\Funnel(s)$.
    By continuity of $\Lpath$ and $\Funnel$, there exists 
    \[
        \hat{s}\coloneqq \min \setdef{t\in[\hat{t},\hat{t}+T]}{\Norm{\Lpath^\star(t)}=\Funnel(t)}.
    \]
    We have $\graph(\Lpath_n)\subset\cF_{\Funnel}$ for all $n\in\N$ by assumption, cf. \Cref{Prop:LpathGraph}.
    Since, in addition, $\Lpath^\star$ is a bounded function,
    there exists a compact set $\cK$ such that $\im(\Lpath^\star)\subset\cK$ and $\im(\Lpath_n)\subset\cK$ for all $n\in\N$.
    Define the continuously differentiable function 
    \[
        \omega:[\hat{t},\hat{t}+T]\times\cK\to\R,\quad  (s,x)\mapsto 1-\frac{\Norm{x}^2}{\Funnel(s)^2}.
    \]
    Due to the compactness of $[\hat{t},\hat{t}+T]$ and $\cK$,
    the  function $\omega$ is Lipschitz continuous with Lipschitz constant $L_{\omega}>0$.
    We have $\omega(s,\Lpath_n(s))>0$ for all $n\in\N$ and all $s\in[\hat{t},\hat{t}+T]$ because
    $\graph(\Lpath_n)\subset\cF_{\Funnel}$ for all  $n\in\N$.
    Let $L^\star>0$ be the Lipschitz constant of~$\Lpath^\star$.
    Since $\omega(\hat{s},\Lpath^\star(\hat{s}))=0$, 
    we estimate the following for all $s\in[\hat{t},\hat{s}]$ and all $n\in\N$.
    \begin{align*}
        \omega(s,\Lpath_n(s))
        &=\Abs{\omega(s,\Lpath_n(s))}=\Abs{\omega(s,\Lpath_n(s))- \omega(\hat{s},\Lpath^\star(\hat{s}))}\\
        &\leq L_{\omega}\Norm{
            \begin{pmatrix}
                s-\hat{s}\\
                \Lpath_n(s)- \Lpath^\star(\hat{s})
        \end{pmatrix}}
        = L_{\omega}\Norm{
        \begin{pmatrix}
                s-\hat{s}\\
            \Lpath_n(s)-\Lpath^\star(s)+ \Lpath^\star(s)-\Lpath^\star(\hat{s})
        \end{pmatrix}}\\
        &\leq  L_{\omega}\Abs{s-\hat{s}}+
        L_{\omega}\Norm{\Lpath_n(s)-\Lpath^\star(s)}+L_{\omega}L^\star\Abs{s-\hat{s}}.
    \end{align*}
    Since  $\int_{\hat{t}}^{\hat{s}}({(L_{\omega}+L_{\omega}L^\star)\Abs{s-\hat{s}}})^{-1}\d
    s=\infty$, there exists $\eps>0$ with
    \[
    \int_{\hat{t}}^{\hat{s}}({(L_{\omega}+L_{\omega}L^\star)\Abs{s-\hat{s}}+L_{\omega}\eps})^{-1}\d
    s-\hat{s}>M.
    \]
    As a consequence of the uniform convergence of $\Lpath_n$ to $\Lpath^{\star}$, there
    exists $N\in\N$ such that $\Norm{\Lpath_n(s)-\Lpath^\star(s)}<\eps$
    for all $n\geq N$ and all $s\in[\hat{t},\hat{s}]$.
    Thus, we arrive at the following contradiction for $n\geq N$.
   \begin{align*}
        M &\geq\int_{\hat{t}}^{\hat{t}+T}\OrigFunnelPenaltyFunc(s,\Lpath_n(s))\d {s} 
        =\int_{\hat{t}}^{\hat{t}+T}
           \frac{\Norm{\Lpath_n(s)}^2}{\Funnel(s)^2-\Norm{\Lpath_n(s)}^2} 
           \d {s}
        = \int_{\hat{t}}^{\hat{t}+T}
            \frac {1}{\omega(s,\Lpath_n(s))}-1\d s\\
        &\geq \int_{\hat{t}}^{\hat{s}}
            \frac {1}{\omega(s,\Lpath_n(s))}-1\d s
        \geq \int_{\hat{t}}^{\hat{s}}
            \frac {1}{
        (L_{\omega}+L_{\omega}L^\star)\Abs{s-\hat{s}}+
        L_{\omega}\Norm{\Lpath_n(s)-\Lpath^\star(s)}}\d s -\hat{s}\\
        &> \int_{\hat{t}}^{\hat{s}} 
            \frac {1}{
        (L_{\omega}+L_{\omega}L^\star)\Abs{s-\hat{s}}+
        L_{\omega}\eps}\d s -\hat{s}>M.
    \end{align*}
    This completes the proof.
\end{proof}

\Cref{Prop:LpathGraph} and \Cref{Lem:OriginalCostFuncLimit} will be essential in proving that the optimal control problem~\eqref{eq:MpcOCP}
using the $\OrigFunnelPenaltyFunc$ from \eqref{eq:PenaltyOriginalFMPC} has a solution and that this solution ensures that 
the tracking error~$e$ evolves within the funnel~$\cF_{\Funnel}$.
To generalise these properties beyond the specific function $\OrigFunnelPenaltyFunc$, 
the following~\Cref{Def:FunnelPenaltyFunction} introduces the concept of \emph{funnel penalty functions}.
These are the functions complying with \Cref{Prop:LpathGraph} and \Cref{Lem:OriginalCostFuncLimit}.

\begin{definition}[Strict funnel penalty function]\label{Def:FunnelPenaltyFunction}
Given $\Funnel\in\cG$, consider a measurable function ${\FunnelPenaltyFunc:\Rp\times\R^{m}\to\R\cup\cbl\infty\cbr}$ whose 
restriction $\FunnelPenaltyFunc|_{\cF_{\Funnel}}$ is non-negative and continuous.
We call $\FunnelPenaltyFunc$ a \emph{strict funnel penalty function} for $\Funnel$, if, for all $\hat{t}\in\Rp$, $T>0$, and
every Lipschitz path $\Lpath\in\Lip([\hat{t},\hat{t}+T],\R^m)$ with $(\hat{t}, \Lpath(\hat{t}))\in\cF_{\Funnel}$ the following holds:
\begin{enumerate}
    \item[\Itemlabel{Item:FunnelPenaltyIntegral}{(F.1)}]
        ${\displaystyle \int_{\hat{t}}^{\hat{t}+T}\FunnelPenaltyFunc(s,\Lpath(s))\d{s}<\infty\iff\graph(\Lpath)\subset\cF_{\Funnel}}$.
    \item[\Itemlabel{Item:FunnelPenaltyConvergence}{(F.2)}] If a sequence $(\Lpath_n)\in\Lip([\hat{t},\hat{t}+T],\R^m)^\N$ with $(\hat{t}, \Lpath_n(\hat{t}))\in\cF_{\Funnel}$ for all $n\in\N$
    converges uniformly to $\Lpath$ and there exists $M\geq0$ such that 
    $\int_{\hat{t}}^{\hat{t}+T}\FunnelPenaltyFunc(s,\Lpath_n(s))\d{s}\leq M$ for all $n\in\N$, then 
    \[
        \int_{\hat{t}}^{\hat{t}+T}\FunnelPenaltyFunc(s,\Lpath(s))\d{s}<\infty.
    \]
\end{enumerate}
\end{definition}
\begin{example}
    Let $\OrigFunnelPenaltyFunc:\Rp\times\R^n$ be given as in~\eqref{eq:PenaltyOriginalFMPC} for $\Funnel\in\cG$.
    For every non-negative function $\FunnelPenaltyFuncWithout \in\cC(\Rp\times\R^m, \R)$,
    the function $\FunnelPenaltyFunc:\Rp\times \R^m\to \R\cup\{\infty\}$ defined by
    \[
        \FunnelPenaltyFunc(t,y)\coloneqq  \OrigFunnelPenaltyFunc(t,y) + \FunnelPenaltyFuncWithout(t,y)
    \]
    is a strict funnel penalty function. Since $\FunnelPenaltyFuncWithout$ is
    bounded on $\cF_{\Funnel}$, this is a direct result of
    \Cref{Prop:LpathGraph} and \Cref{Lem:OriginalCostFuncLimit}. 
    Consequently, one can model additional soft constraints via the function $\FunnelPenaltyFuncWithout$ 
    to be penalised in the MPC~\Cref{Algo:MPC}, without 
    losing the property of having a strict funnel penalty function.
\end{example}

In~\Cref{Def:FunnelPenaltyFunction}, $\FunnelPenaltyFunc$ is called \emph{strict}, because 
condition~\ref{Item:FunnelPenaltyIntegral} ensures that the Lipschitz path~$\Lpath$ evolves within the interior of the funnel~$\cF_{\Funnel}$, i.e.
$\Norm{\Lpath(s)}<\Funnel(s)$ for all $s\in[\hat{t},\hat{t}+T]$.
If, in addition, one also wants to allow for equality, i.e. only requires $\Norm{\Lpath(s)}\leq\Funnel(s)$ for all $s\in[\hat{t},\hat{t}+T]$, 
and replaces $\cF_{\Funnel}$ in~\Cref{Def:FunnelPenaltyFunction} by
\[
    \bar{\cF}_{\Funnel}\coloneqq  \setdef{(t,e)\in \Rp\times\R^{m}}{\Norm{e} \leq \Funnel(t)},
\]
then property~\ref{Item:FunnelPenaltyConvergence} can be omitted.
In this case, property~\ref{Item:FunnelPenaltyIntegral} and continuity of the function~$\FunnelPenaltyFunc$ on
the larger set~${\bar{\cF}_{\Funnel}}$ already imply condition~\ref{Item:FunnelPenaltyConvergence},
as the following~\Cref{Def:NonStrictFunnelPenaltyFunction} shows.
We therefore call a measurable function $\FunnelPenaltyFunc:\Rp\times\R^{m}\to\R\cup\cbl\infty\cbr$ 
that fulfils property~\ref{Item:FunnelPenaltyIntegral}
for the set $\bar{\cF}_{\Funnel}$ and whose restriction $\FunnelPenaltyFunc|_{\bar{\cF}_{\Funnel}}$ is non-negative and continuous
a \emph{non-strict funnel penalty function} for $\Funnel$.
\begin{lemma}\label{Def:NonStrictFunnelPenaltyFunction}
    Let $\Funnel\in\cG$ and $\FunnelPenaltyFunc:\Rp\times\R^{m}\to\R\cup\cbl\infty\cbr$ be a non-strict funnel penalty function.
    Then, $\FunnelPenaltyFunc$ fulfils condition~\ref{Item:FunnelPenaltyConvergence} of~\Cref{Def:FunnelPenaltyFunction} for the set $\bar{\cF}_{\Funnel}$.
\end{lemma}
\begin{proof}
    Let $\hat{t}\in\Rp$, $T>0$, and  $\Lpath\in\Lip([\hat{t},\hat{t}+T],\R^m)$ be a Lipschitz path with ${(\hat{t}, \Lpath(\hat{t}))\in\bar{\cF}_{\Funnel}}$.
    Further, let $M\geq0$ and $(\Lpath_n)\in\Lip([\hat{t},\hat{t}+T],\R^m)^\N$ be a to $\Lpath$ uniformly converging sequence  with 
    $(\hat{t}, \Lpath_n(\hat{t}))\in\cF_{\Funnel}$ and $\int_{\hat{t}}^{\hat{t}+T}\FunnelPenaltyFunc(s,\Lpath_n(s))\d{s}\leq M$ 
    for all $n\in\N$. Due to property~\ref{Item:FunnelPenaltyIntegral}, we have $\graph(\Lpath_n)\subset\bar{\cF}_{\Funnel}$.
    Recall that, for non-strict funnel penalty functions, \Cref{Def:FunnelPenaltyFunction} is formulated in terms of $\bar{\cF}_{\Funnel}$.
    We show that $\Norm{\Lpath(s)}\leq\Funnel(s)$ for all $s\in[\hat{t},\hat{t}+T]$.
    Assume there exists $\hat{s}\in(\hat{t},\hat{t}+T]$ with $\Norm{\Lpath(\hat{s})}>\Funnel(\hat{s})$.
    Then, there exists $\eps>0$ with $\Norm{\Lpath(\hat{s})}>\Funnel(\hat{s})+\eps$.
    Since the uniform convergence of $(\Lpath_n)$ towards $\Lpath$ implies pointwise convergence of $(\Lpath_n)$,
    there exists $K>0$ such that $\Norm{\Lpath(\hat{s})-\Lpath_k(\hat{s})}<\eps$ for all $k\geq K$.
    Furthermore, $\Norm{\Lpath_k(\hat{s})}\leq\Funnel(\hat{s})$ since $\graph(\Lpath_k)\subset\bar{\cF}_{\Funnel}$ for all $k\in\N$.
    This raises the following contradiction for $k\geq K$
    \[
    \Funnel(\hat{s})+\varepsilon<\Norm{\Lpath(\hat{s})}
    \leq\Norm{\Lpath(\hat{s})-\Lpath_k(\hat{s})}+\Norm{\Lpath_k(\hat{s})}
    \leq\varepsilon+\Funnel(\hat{s}).
    \]
    This completes the proof.
\end{proof}

\begin{remark}
    A strict funnel penalty function $\FunnelPenaltyFunc$ cannot be continuous on the whole set $\bar{\cF}_{\Funnel}$.
    Otherwise, for a Lipschitz path $\Lpath\in\Lip([\hat{t},\hat{t}+T],\R^m)$ with
    $\graph(\Lpath)\subset\bar{\cF}_{\Funnel}$ and $\Norm{\Lpath(s)}=\Funnel(s)$ for
    some $s\in [\hat{t},\hat{t}+T]$,
    the function $t\mapsto\FunnelPenaltyFunc(t,\Lpath(t))$ would be bounded and, thus, integrable.
    Hence, it is clear that the two concepts of a strict and non-strict funnel penalty function are mutually exclusive, meaning a single function cannot be both.
\end{remark}
We have already seen that $\OrigFunnelPenaltyFunc$ as in~\eqref{eq:PenaltyOriginalFMPC} is a strict funnel penalty function.
Now, we want to give an example for a non-strict funnel penalty function.

\begin{example}
    Let $\Funnel\in\cG$ and $\tilde\FunnelPenaltyFuncWithout_{\Funnel}\in\cC(\bar{\cF}_{\Funnel},\R)$ be a non-negative function.
    Then,
    \[
        \FunnelPenaltyFunc:\Rp\times\R^m\to \R\cup\{\infty\},  \quad
        \FunnelPenaltyFunc(t,e)=\begin{cases}
        \tilde\FunnelPenaltyFuncWithout_{\Funnel}(t,e),& (t,e)\in\bar{\cF}_{\Funnel}\\
        \infty, &\text{else}
        \end{cases}
    \]
     is a non-strict funnel penalty function. 
     To see this, let $\Lpath\in\Lip([\hat{t},\hat{t}+T],\R^m)$ be a Lipschitz path with $(\hat{t}, \Lpath(\hat{t}))\in\bar{\cF}_{\Funnel}$.
     We have to show that \ref{Item:FunnelPenaltyIntegral} from~\Cref{Def:FunnelPenaltyFunction} holds for $\bar{\cF}_{\Funnel}$.
     First, let $\graph(\Lpath)\subset\bar{\cF}_{\Funnel}$. Then,
     $\FunnelPenaltyFuncWithout_{\Funnel}(s,\Lpath(s))=\tilde\FunnelPenaltyFuncWithout_{\Funnel}(s,\Lpath(s))$
     for all $s\in[\hat{t},\hat{t}+T]$. Due to the continuity of the involved functions,
     $\tilde\FunnelPenaltyFuncWithout_{\Funnel}(s,\Lpath(s))$ is bounded on the compact interval~$[\hat{t},\hat{t}+T]$.
     Thus, the integral $\int_{\hat{t}}^{\hat{t}+T}\FunnelPenaltyFunc(s,\Lpath(s))\d{s}$ is finite. 
     To show the reverse implication, assume now that the integral is finite but that there exists 
     $\hat{s}\in(\hat{t},\hat{t}+T]$ with $\Norm{\Lpath(\hat{s})}>\Funnel(\hat{s})$.
     Then, there exists $\eps>0$ with $\Norm{\Lpath(s)}>\Funnel(s)$ for all $s \in [\hat{s}-\eps,\hat{s}]$ because $\Lpath$ and $\Funnel$ are continuous functions.
     Hence, the following contradiction arises. 
    \[
        \infty>\int_{\hat{t}}^{\hat{t}+T}\FunnelPenaltyFunc(s,\Lpath(s))\d{s}\geq \int_{\hat{s}-\eps}^{\hat{s}}\FunnelPenaltyFunc(s,\Lpath(s))\d{s}=\int_{\hat{s}-\eps}^{\hat{s}}\infty\ \d{s}=\infty.
    \]
\end{example}

We discussed the essential properties of the function~$\OrigFunnelPenaltyFunc$ and summed them up in~\Cref{Def:FunnelPenaltyFunction}
of funnel penalty functions. To use this concept in optimal control problems of the form~\eqref{eq:MpcOCP}, 
we additionally want to be able to penalise the necessary control effort in the cost function. To this end, 
the following~\Cref{Def:FunnelStageCostFunc} introduces \emph{funnel stage cost} functions. 
An example is the funnel MPC stage cost function~$\ell_{\Funnel}$ as in~\eqref{eq:stageCostFunnelMPC}.

\begin{definition}[Funnel stage cost function]\label{Def:FunnelStageCostFunc}
    Let $\Funnel\in\cG$ and $\FunnelPenaltyFunc:\Rp\times\R^{m}\to\R\cup\cbl\infty\cbr$ be a (non)-strict funnel penalty function. 
    For $\lambda_u\in\Rp$, we call a function 
    \begin{align*}
        \OrigFunnelStageCost:\Rp\times\R^{m}\times\R^{m}\to\R\cup\cbl\infty\cbr,\quad
        (t,z,u)\mapsto \FunnelPenaltyFunc(t,z)+\lambda_u\Norm{u}^2
    \end{align*}
    a \emph{(non)-strict funnel stage cost}.
\end{definition}
In \Cref{Def:FunnelStageCostFunc}, the penalisation term for the control input
$u$ consists of the squared norm of $u$ multiplied by the parameter~$\lambda_u$.
This parameter allows to adjust a suitable trade-off between tracking
performance and required control effort.
Note that $\lambda_u=0$ is explicitly allowed contrary to~\eqref{eq:stageCostClassicalMPC}.
Of course, utilising more sophisticated penalty terms is also possible. 
One option, for example, is the usage of norms induced by a positive definite matrix.
If a reference input signal~$u_{\rf}$ is known, the second
summand may also be replaced by $\| u - u_{\rf}(t) \|^2$.
In this work however, we will restrict ourselves to the presented case.

If not explicitly mentioned otherwise, we will use strict funnel stage cost functions in the rest of the present thesis and we will refer to them merely by the term \emph{funnel stage cost}.
The presented results generally also remain valid for non-strict funnel stage cost functions,
but then only with respect to the set~$\bar{\cF}_{\Funnel}$,
i.e. inequalities of the form $\Norm{e}<\Funnel(t)$ have to be replaced by $\Norm{e}\leq\Funnel(t)$. 
\begin{remark}\label{Rem:FunnelStageCostPositive}
    Note that every (non)-strict funnel stage cost $\OrigFunnelStageCost$ is non-negative for every element $(t,z,u)\in\cF_{\Funnel}\times\R^m$.
\end{remark}

\section{Model class}\label{Sec:ModelClass}

In the previous \Cref{Sec:FunnelStageCostFunctions}, the essential aspects of
the cost function $\ell_{\Funnel}$ from \eqref{eq:stageCostFunnelMPC} were
identified to introduce the more general concept of funnel stage cost
functions. This was done by considering Lipschitz paths in order to conduct this
analysis in isolation from any differential equation.
In this section, however, we will introduce the class of surrogate models for
the system~\eqref{eq:Sys} to be utilised in the MPC~\Cref{Algo:MPC}.
We consider non-linear control affine multi-input multi-output models of order $r\in\N$ of the form
\begin{equation} \label{eq:Model_r}
    \begin{aligned}
   &\textover[r]{$\yM^{(r)}(t)$}{$\big(\yM(t_0),\ldots,\yM^{(r-1)}(t_0)\big)\ $} = \fM \big(\oTM(\yM,\ldots,\yM^{(r-1)} )(t) \big) + \gM \big(\oTM(\yM,\ldots,\yM^{(r-1)} )(t) \big) u(t), \\
   &\left.
        \begin{aligned}
       \yM|_{[0,t_0]} &= \yM^0  \in \cC^{r-1}([0,t_0],\R^m), && \mbox{if } t_0 >0,\\
       \big(\yM(t_0),\ldots,\yM^{(r-1)}(t_0)\big) &= \yM^0  \in\R^{rm}, && \mbox{if } t_0 = 0,
        \end{aligned} \right\}
    \end{aligned}
\end{equation}
with $t_0\geq 0$,  initial trajectory $\yM^0$, control input $u\in L_{\loc}^{\infty}([t_0,\infty),\R^m)$,
and output ${\yM(t)\in\R^m}$ at time $t\geq t_0$.
Note that, like the system~\eqref{eq:Sys}, $u$ and $\yM$ have the same dimension~$m\in\N$.
The model consists of two locally Lipschitz continuous function $\fM  \in \Lip_{\loc}(\R^q,\R^m)$, $\gM\in \Lip_{\loc}(\R^q,\R^{m\times m})$, 
and an  operator $\oTM$.
To ensure that the control $u$ can always influence the dynamics, we assume that $\gM$ 
is everywhere point-wise invertible, i.e. $\gM$ satisfies $\gM(z)\in\GL_{m}(\R)$ for all $z\in\R^q$.
The operator~$\oTM$ is causal, locally Lipschitz, satisfies a bounded-input bounded-output 
and a limited memory property. It is characterised in detail in the following~\Cref{Def:OperatorClass}.

\begin{definition}[Operator class~$\cT_{t_0}^{n,q}$]\label{Def:OperatorClass} 
For $n,q\in\N$ and $t_0\in\Rp$, the set $\cT_{t_0}^{n,q}$ denotes the class of operators $\oT:
\cR(\Rp,\R^n) \to L^\infty_{\loc} ([t_0,\infty), \R^{q})$
for which the following properties hold:
\begin{enumerate}
    \item[\Itemlabel{Item:OperatorPropCasuality}{(T.1)}]\emph{Causality}:  $\fa y_1,y_2\in\cR(\Rp,\R^n)$  $\fa t\geq t_0$:
    \[
        y_1\vert_{[0,t]} = y_2\vert_{[0,t]}
        \ \Impl\ 
            \oT(y_1)\vert_{[t_0,t]}=\oT(y_2)\vert_{[t_0,t]}.
    \]
    \item[\Itemlabel{Item:OperatorPropLipschitz}{(T.2)}]\emph{Local Lipschitz}: 
    $\fa t \ge t_0 $ $\fa y \in \cR([0,t] ; \R^n)$ 
    $\ex \Delta, \delta, c > 0$ 
    $\fa y_1, y_2 \in \cR(\Rp,\R^n)$ with
    $y_1|_{[0,t]} = y_2|_{[0,t]} = y $ 
    and $\Norm{y_1(s) - y(t)} < \delta$,  $\Norm{y_2(s) - y(t)} < \delta $ for all $s \in [t,t+\Delta]$:
    \[
     \esssup_{\mathclap{s \in [t,t+\Delta]}}  \Norm{\oT(y_1)(s) - \oT(y_2)(s) }  
        \le c \ \sup_{\mathclap{s \in [t,t+\Delta]}}\ \Norm{y_1(s)- y_2(s)}.
    \] 
    \item[\Itemlabel{Item:OperatorPropBIBO}{(T.3)}]\emph{Bounded-input bounded-output (BIBO)}:
    $\fa c_0 > 0$ $\ex c_1>0$  $\fa y \in \cR(\Rp, \R^n)$:
    \[
    \sup_{t \in \Rp} \Norm{y(t)} \le c_0 \ 
    \Impl \ \sup_{t \in [t_0,\infty)} \Norm{\oT(y)(t)}  \le c_1.
    \]
    \item[\Itemlabel{Item:OperatorPropLimitMemory}{(T.4)}]\emph{Limited memory}: 
    $\ex \tau\geq0$ $\fa t\geq t_0$
    $\fa y_1, y_2 \in \cR(\Rp, \R^n)$ with 
    $y_1|_{I}= y_2|_{I}$ on the interval $I\coloneqq  [t-\tau,\infty)\cap\Rp$ and $\oT(y_1)|_{J} = \oT(y_2)|_{J}$ on the interval ${J\coloneqq [t-\tau,t]\cap[t_0,t]}$:
    \[
        \oT(y_1)\vert_{[t,\infty)}=\oT(y_2)\vert_{[t,\infty)}.
    \]
    The value $\tau$ in is called \emph{memory limit} of the operator~$\oT$.
\end{enumerate}
\end{definition}

Note that an operator $\oTM\in\cT_{t_0}^{n,q}$ can model \emph{non-linear time delays}, where $t_0$ corresponds to the initial delay,
and that it can even be the solution operator of an infinite-dimensional dynamical system, e.g. a partial differential equation.

We summarise our assumptions and define the general model class under consideration.
\begin{definition}[Model class~$\cM^{m,r}_{t_0}$]\label{Def:ModelClass}
    We say the model~\eqref{eq:Model_r} belongs to the \emph{model class}~$\cM^{m,r}_{t_0}$ for $m,r\in\N$, and $t_0\in\Rp$,
    written~$(\fM,\gM,\oTM)\in\cM^{m,r}_{t_0}$, if, for some $q\in\N$, the following holds:
    $\fM\in \Lip_{\loc}(\R^q ,\R^m)$,
    $\gM\in \Lip_{\loc}(\R^q ,\R^{m\times m})$ satisfies ${\gM(z)\in\GL_{m}(\R)}$ for all $z\in\R^q$,
    and ${\oTM\in\cT^{rm,q}_{t_0}}$.
\end{definition}

Due to their quite technical nature,
the properties of the~$\oTM\in\cT_{t_0}^{n,q}$ deserve some additional explanation. 
However, before commenting on them in~\Cref{Rem:PropertiesOperator},
we briefly discuss a simple example of a model belonging 
to the considered model class in order to have a familiar picture in mind.
In particular, this provides a simple candidate for an operator~$\oT$.
\begin{example}[Linear time-invariant model]\label{Ex:LTISystem}
    Consider a linear multi-input, multi-output differential equation of the form  
    \begin{equation}\label{eq:LTISystem}
        \begin{aligned}
        \dot x(t) &= A x(t) + B u(t),  \quad x(t_0) = x^0\\
        y(t) &= C x(t),
        \end{aligned}
    \end{equation}
    with $A \in \R^{n \times n}$ and $C^\top,B\in\R^{n\times m}$.
    Further, assume that this model has a \emph{strict relative degree} $r\in\N$, i.e.
    $CA^kB=0$ for all $k<r-1$ and $CA^{r-1}B\in\GL_m(\R)$.
    By~\cite[Lemma~3.5]{IlchRyan07}, there exists an invertible $U \in \R^{n \times n}$ such that with 
    $[z_1^\top,\ldots, z_r^\top,\eta^\top]^\top = U x$ the above system can be transformed into the \emph{(linear) Byrnes-Isidori form}
    \begin{align}\label{eq:LTIByrnesIsidori}
        \dot z_{i}(t) &= z_{i+1}, &z_i(t_0)&=z_i^0,\nonumber\\
        \dot z_{r}(t) &= \sum_{j=1}^{r}R_jz_j(t)+S\eta+\Gamma u(t),&z_r(t_0)&=z_r^0,\\
        \dot \eta(t)  &= Q \eta(t) + P z_1(t), &\eta(t_0)&=\eta^0,\nonumber 
        \shortintertext{with output}
        y(t) &= z_1(t),\nonumber
    \end{align}
    where $R_j \in \R^{m \times m}$ for all $j=1,\ldots,r$, $S, P^\top\in\R^{m\times (n-rm)}$, $Q\in\R^{(n-rm)\times (n-rm)}$,
    and $\Gamma = CA^{r-1}B$.
    Define the linear integral operator $\oT: \cR(\Rp,\R^{m})^r \to L^\infty_{\loc} ([t_0,\infty), \R^{q})$ by
    \[
        \oT(z_1,\ldots,z_r)(t)\coloneqq \sum_{j=1}^{r}R_jz_j(t)
        +S\rbl\me^{Q(t-t_0)}\eta_0+\int_{t_0}^t\me^{Q(t-s)}Pz_1(s)\d{s}\rbr.
    \]
    Utilising this operator~$\oT$ and the Byrnes-Isidori form~\eqref{eq:LTIByrnesIsidori}, 
    the differential equation~\eqref{eq:LTISystem} 
    can be put in the equivalent form~\eqref{eq:Model_r} with initial value 
    \[
    (y(t_0),\dot{y}(t_0),\ldots,y^{(r-1)}(t_0)) = (Cx^0,CAx^0,\ldots, CA^{r-1}x^0).
    \]
    In the following, we examine the properties~\ref{Item:OperatorPropCasuality}--\ref{Item:OperatorPropLimitMemory} of the operator~$\oT$.
    It is easy to see that it satisfies properties~\ref{Item:OperatorPropCasuality}
    and \ref{Item:OperatorPropLipschitz}.
    To also verify the limited memory property~\ref{Item:OperatorPropLimitMemory} for $\oT$ with $\tau=0$, 
    let $\hat{t}\geq t_0$  and $(y_1,\ldots,y_r),(\tilde{y}_1,\ldots,\tilde{y}_r)\in\cR(\Rp,\R^m)^r$ 
    with $\oT(y_1,\ldots,y_r)(\hat{t})=\oT(\tilde{y}_1,\ldots,\tilde{y}_r)(\hat{t})$ 
    and $(y_1,\ldots,y_r)(t)=(\tilde{y}_1,\ldots,\tilde{y}_r)(t)$ for all ${t\geq\hat{t}}$.
    Using the shorthand notation 
    \[
    V(z_1,\ldots,z_r)(t)\coloneqq \sum_{j=1}^{r}R_jz_j(t)+S\me^{Q(t-t_0)}\eta_0\  \text{ and }\ L(z_1)(t)\coloneqq S\int_{t_0}^t\me^{Q(t-s)}Pz_1(s)\d{s}, 
    \]
    we have $\oT(z_1,\ldots,z_r)(t)=V(z_1,\ldots,z_r)(t)+L(z_1)(t)$.
    It is clear that $V(y_1,\ldots,y_r)(t)$ and $V(\tilde{y}_1,\ldots,\tilde{y}_r)(t)$ are identical for all $t\geq\hat{t}$.
    To show the limited memory property~\ref{Item:OperatorPropLimitMemory} for the operator~$\oT$,
    it therefore is sufficient to show $L(y_1)(t)=L(\tilde{y}_1)(t)$ for all $t\geq\hat{t}$.
    For $t\geq \hat{t}$, we have
    \begin{align*}
     L(y_1)(t)&=S\int_{t_0}^t\me^{Q(t-s)}Py_1(s)\d{s}\\
                        &=L(y_1)(\hat t)+S\int_{\hat{t}}^t\me^{Q(t-s)}Py_1(s)\d{s}\\
                        &=\oT(y_1,\ldots,y_1)(\hat{t})-V(y_1,\ldots,y_r)(\hat{t})+S\int_{\hat{t}}^t\me^{Q(t-s)}Py_1(s)\d{s}\\
                        &=\oT(\tilde{y}_1,\ldots,\tilde{y}_r)(\hat{t})-V(\tilde{y}_1,\ldots,\tilde{y}_r)(\hat{t})+S\int_{\hat{t}}^t\me^{Q(t-s)}P\tilde{y}_1(s)\d{s}\\
                        &=L(\tilde{y}_1)(\hat t)+S\int_{\hat{t}}^t\me^{Q(t-s)}P\tilde{y}_1(s)\d{s}\\
                        &=S\int_{t_0}^t\me^{Q(t-s)}P\tilde{y}_1(s)\d{s}\\
                        &=L(\tilde{y}_1)(t).
    \end{align*}%
    Thus, $\oT$ has the property~\ref{Item:OperatorPropLimitMemory} with $\tau=0$.
    Additionally assume that~\eqref{eq:LTISystem} has~\emph{asymptotically stable zero dynamics}, i.e.
    \[
        \fa \lambda\in\C_{\geq 0}:\quad \det
        \begin{bmatrix}
           \lambda I_n - A & B\\ 
           C &0
        \end{bmatrix}\neq 0.
    \]
    This, also called \emph{minimum phase} property in literature, see e.g.~\cite{Isid95,IlchWirt13}, 
    implies that all eigenvalues of the matrix~$Q$ in~\eqref{eq:LTIByrnesIsidori} are in the open left plane,
    i.e. $\spec(Q)\subset\C_{<0}$, see~\cite[Lemma~3.5]{IlchRyan07}.
    Consequently, the operator~$\oT$ fulfils the BIBO property~\ref{Item:OperatorPropBIBO}, 
    see also the diagram in~\cite[Section 2.1.2]{BergIlch21} nicely illustrating the relationship 
    between the minimum phase property of model~\eqref{eq:LTISystem} and the BIBO property of~$\oT$.
\end{example}

\begin{example}[Non-linear model with state space representation]\label{Ex:ControlAffineMod}
    Consider a non-linear differential equation of the form 
    \begin{equation}\label{eq:NonLinSysStates}
    \begin{aligned}
        \dot{x}(t)  & = f(x(t)) + g(x(t)) u(t),\quad x(t_0)=x^0,\\
        y(t)        & = h(x(t)),
    \end{aligned}
    \end{equation}
    with~$t_0\in\R_{\ge 0}$, $x^0\in\R^n$, and non-linear functions~$f:\R^n\to \R^n$, $g:\R^n\to \R^{n\times m}$ and $h : \R^n \to \R^m$.
    We show in the following that~\eqref{eq:NonLinSysStates} can, under certain conditions, be put in the form~\eqref{eq:Model_r}
    and thus is an admissible candidate for a model.
    We recall the notion of relative degree for the differential equation~\eqref{eq:NonLinSysStates}, see e.g.~\cite[Sec. 5.1]{Isid95}.
    Assuming that $f,g,h$ are sufficiently smooth, the Lie derivative of~$h$ along~$f$ is defined by $\rbl L_f h\rbr(x) = h'(x) f(x)$
    and, successively, we define~$L_f^k h = L_f (L_{f}^{k-1} h)$ with~$L_f^0 h = h$.
    Furthermore, for the matrix-valued function~$g$, we have
    \[
        (L_gh)(x) = \sbl (L_{g_1}h)(x), \ldots, (L_{g_m}h)(x) \sbr,
    \]
    where~$g_i$ denotes the~$i$-th column of~$g$ for $i=1,\ldots, m$.
    Then, the differential equation~\eqref{eq:NonLinSysStates} is said to have \emph{strict (global) relative degree}~$r \in \mathbb{N}$, if
    \begin{align*}
        \fa k \in \{1,\ldots,r-1\}\ \fa x \in \R^n:\
            (L_g L_f^{k-1} h)(x)  &= 0  \\
          \text{and}\quad
            (L_g L_f^{r-1} h)(x)  &\in \GL_{m}(\R).
    \end{align*}
    If~\eqref{eq:NonLinSysStates} has relative degree~$r$, then, under the additional assumptions provided in~\cite[Cor.~5.6]{ByrnIsid91a},
    differential equation~\eqref{eq:NonLinSysStates} can be transformed into \emph{(non-linear) Byrnes-Isidori form} -- 
    a generalisation of \eqref{eq:LTIByrnesIsidori}.
    This means there exists a diffeomorphism~$\Phi:\R^n\to\R^n$ such that the coordinate transformation 
    $(y(t), \dot y(t),\ldots, y^{(r-1)}(t),\eta(t)) = \Phi(x(t))$ puts 
    the differential equation \eqref{eq:NonLinSysStates} into the form
    \begin{subequations}\label{eq:BIF}
        \begin{align}
            y^{(r)}(t)   &= p\big(y(t),\ldots,y^{(r-1)}(t),\eta(t)\big) + \Gamma\big(y(t),\ldots,y^{(r-1)}(t),\eta(t)\big)\,u(t),\\
            \dot \eta(t) &= q\big(y(t),\ldots,y^{(r-1)}(t),\eta(t)\big),\label{eq:zero_dyn}
        \end{align}
    \end{subequations}
    where $p:\R^n\to \R^{m}$, $q: \R^n\to \R^{n-rm}$, $\Gamma = L_g L_f^{r-1} h:\R^n\to\R^{m\times m}$ 
    are continuously differentiable and $(y(t_0),\dot y(t_0),\ldots, y^{(r-1)}(t_0),\eta(t_0)) =  \Phi(x_0)$.
    Note that, under these assumptions, the derivatives of the output~$y$ of~\eqref{eq:NonLinSysStates} are given by
    $y^{(i)}(t) = (L_f^i h)(x(t))$ for $i=0,\ldots,r-1$.
    In the following, we assume the existence of such diffeomorphism $\Phi:\R^n\to\R^n$
    but not necessarily the conditions stated in~\cite[Cor.~5.6]{ByrnIsid91a} as these are sufficient but not necessary for the existence of~$\Phi$.
    We further assume that internal dynamics~\eqref{eq:zero_dyn} satisfy the following \emph{bounded-input, bounded-state} (BIBS) condition:
    \begin{multline}\label{eq:BIBO-ID}
        \fa c_0 >0\ \ex c_1 >0\ \fa t_0\ge 0\ \fa  \eta^0\in\R^{n-rm} 
        \fa  \zeta\in L^\infty_{\loc}([t_0,\infty),\R^{rm}):\\ 
        \Norm{\eta^0}+ \SNorm{\zeta}  \leq c_0 \implies\ \SNorm{\eta (\cdot;t_0,\eta^0,\zeta)} \leq c_1,
    \end{multline}
    where $\eta (\cdot;t_0, \eta^0,\zeta):[t_0,\infty)\to\R^{n-rm}$ denotes
    the unique global solution of~\eqref{eq:zero_dyn} when $(y(t),\ldots,y^{(r-1)}(t))$ is substituted by~$\zeta$.
    Note that, in view of condition~\eqref{eq:BIBO-ID}, the maximal solution $\eta (\cdot;t_0, \eta^0,\zeta)$ can indeed be extended to a global solution, cf.~\cite[\S~10, Thm.~XX]{Walt98}.
    Utilising the unique global solution~$\eta$ of~\eqref{eq:zero_dyn},
    define operator
    \[
        \oT: \cR(\Rp,\R^{rm}) \to L^\infty_{\loc} ([t_0,\infty), \R^{n}),\quad \zeta\mapsto \oT(\zeta)\coloneqq (\zeta(\cdot),\eta(\cdot;t_0,\eta^0,\zeta)).
    \]
    It is easy to see that $\oT$ satisfies the causality property~\ref{Item:OperatorPropCasuality}.
    The bounded-input bounded-output property~\ref{Item:OperatorPropBIBO} 
    is a direct consequence of BIBS condition~\eqref{eq:BIBO-ID} on internal dynamics~\eqref{eq:zero_dyn}.
    Utilising the fact that the continuously differentiable function $q$ is local Lipschitz continuous
    in combination with the BIBS condition~\eqref{eq:BIBO-ID}, the 
    local Lipschitz property~\ref{Item:OperatorPropLipschitz} can be verified via straightforward calculations. 
    To verify that it also satisfies the limited memory property~\ref{Item:OperatorPropLimitMemory} for $\tau=0$, 
    let $\hat{t}\geq t_0$  and $\zeta_1,\zeta_2\in\cR(\Rp,\R^{rm})$ 
    with $\oT(\zeta_1)(\hat{t})=\oT(\zeta_2)(\hat{t})$ 
    and $\zeta_1(t)=\zeta_2$ for all $t\geq\hat{t}$.
    As $\eta$ is the maximal solution of~\eqref{eq:zero_dyn}, it can be represented for $\zeta\in\cR(\Rp,\R^{rm})$ as 
    \[
        \eta(t;t_0,\eta^0,\zeta)=\eta^0+\int_{t_0}^tq(\zeta(s),\eta(s;t_0,\eta^0,\zeta))\d{s}
    \]
    for all $t\geq t_0$. $\oT(\zeta_1)(\hat{t})=\oT(\zeta_2)(\hat{t})$ implies 
    \[
        \int_{t_0}^{\hat{t}}q(\zeta_{1}(s),\eta(s;t_0,\eta^0,\zeta_{1}))\d{s}
        =\int_{t_0}^{\hat{t}}q(\zeta_2(s),\eta(s;t_0,\eta^0,\zeta_2))\d{s}.
    \]
    Utilising $\zeta_1|_{[\hat{t},\infty)}=\zeta_2|_{[\hat{t},\infty)}$, we have
    \begin{align*} 
        \oT(\zeta_1)(t)&=(\zeta_1(t), \eta^0+\int_{t_0}^tq(\zeta_1(s),\eta(s;t_0,\eta^0,\zeta_1))\d{s})\\
        &=(0,\eta^0+\int_{t_0}^{\hat{t}}q(\zeta_1(s),\eta(s;t_0,\eta^0,\zeta_1))\d{s})
        +(\zeta_1(t), \int_{\hat{t}}^tq(\zeta_1(s),\eta(s;t_0,\eta^0,\zeta_1))\d{s})\\
        &=(0,\eta^0+\int_{t_0}^{\hat{t}}q(\zeta_2(s),\eta(s;t_0,\eta^0,\zeta_2))\d{s})
        +(\zeta_2(t),\int_{\hat{t}}^tq(\zeta_2(s),\eta(s;t_0,\eta^0,\zeta_2))\d{s})\\
        &=(\zeta_2(t), \eta^0+\int_{t_0}^tq(\zeta_2(s),\eta(s;t_0,\eta^0,\zeta_2))\d{s})=\oT(\zeta_2)(t)
    \end{align*}
    for $t\geq \hat{t}$.
    Thus, $\oT$ has the property~\ref{Item:OperatorPropLimitMemory} with $\tau=0$.
    The differential equation~\eqref{eq:NonLinSysStates} therefore can be put in the form~\eqref{eq:Model_r} 
    and it is an admissible model.
\end{example}

\begin{remark}\label{Rem:PropertiesOperator}
    We comment on several aspects of the operator class $\cT_{t_0}^{n,q}$ and its
    properties.
    \begin{enumerate}[(a)]
        \item\label{Item:Rem:PropertiesOperator:Causality}
        Let~$\oT\in\cT_{t_0}^{n,q}$ and $I$ be the interval~$I=[0,\hat{t}]$ or~$I=[0,\hat{t})$ for $\hat{t}\in\Rp$.
        For a given function~$\zeta\in\cR(I,\R^n)$, let $\zeta^{e}$ denote an arbitrary right extension of $\zeta$ 
        on the entire interval of non-negative real numbers, i.e. $\zeta^{e}\in\cR(\Rp,\R^n)$ with $\zeta^{e}|_{I}=\zeta$.
        By virtue of the causality property~\ref{Item:OperatorPropCasuality}, the restriction of $\oT(\zeta^{e})$ to the interval~$I$
        is uniquely determined by the function $\zeta$ in the sense that $\oT(\zeta^{e})|_{I}$ is independent of the chosen extension~$\zeta^{e}$.
        This observation made in~\cite[Remark 2 (iii)]{IlchRyan02b}
        allows us to apply the operator~$\oT$ in a certain sense to functions $\zeta\in\cR(I,\R^n)$
        by utilising an arbitrary right extension $\zeta^{e}\in\cR(\Rp,\R^n)$ instead.
        We will therefore use throughout this work the following notation.
        For $s\in I$, we write $\oT(\zeta)(s)$ in place of $\oT(\zeta^{e})(s)$.
        \item The Lipschitz property~\ref{Item:OperatorPropLipschitz} is a rather technical assumption to ensure
        the existence of a solution of the closed-loop initial value problems~\eqref{eq:Sys} and~\eqref{eq:Model_r}
        when a control law is applied, see~\Cref{Appendix:Th:SolutionExists}.
        We will summarise the corresponding results, tailored to the context,
        in~\Cref{Prop:SolutionExists}. For the technical details, however, we
        would also like to refer the reader to the
        works~\cite{ryan2001controlled,IlchRyan02a,IlchRyan02b} and also to the
        \nameref{Chapter:Appendix}
        where we recall the solution theory for the class of
        functional differential equations considered in this thesis.
        The Lipschitz property will moreover be used in~\Cref{Prop:SolutionUnique} 
        in order to prove uniqueness of the solution of the initial value problem~\eqref{eq:Model_r}.
        \item
        To motivate the BIBO property~\ref{Item:OperatorPropBIBO} of operator~$\oT$, we consider the example of a differential equation 
        of the form
        \begin{subequations}\label{eq:SusmannPeaking}
        \begin{align}
            \dot{y}(t)&=Ay(t)+Bu(t),\label{eq:SusmannPeakingLinearPart}\\
            \dot{\eta}(t)&=f(y(t),\eta(t)),\label{eq:SusmannPeakingNonLinearPart}
        \end{align}
        \end{subequations}
        with matrices $A\in\mathds{R}^{n\times n}$ and $B\times\mathds{R}^{n\times m}$
        and a continuously differentiable function ${f:\mathds{R}^n\times\mathds{R}^{\ell}\to\mathds{R}^{\ell}}$.
        As similarly shown in~\Cref{Ex:LTISystem}, the differential equation~\eqref{eq:SusmannPeaking} can be put in the from~\eqref{eq:Model_r},
        where $\oT$ is the solution operator of the non-linear equation~\eqref{eq:SusmannPeakingNonLinearPart}.
        If the matrices $(A,B)$ are controllable, then 
        a stabilising state feedback $u=Ky$ with $K\in\R^{m\times n}$ can be applied to~\eqref{eq:SusmannPeaking} and the
        linear part~\eqref{eq:SusmannPeakingLinearPart} can be estimated, for $t\geq0$ and $a,b>0$, by $\|x(t)\|\leq 
        b\me^{-at}\|x(0)\|$. Any prespecified~$a$ can be realised by the choice of~$K$.
        However, as stated by Sussmann and Kokotovic in~\cite{sussmann1991peaking}, \emph{one
        cannot, in general, choose $K$ so as to make the number $a$ large without making~$b$ large as
        well}. As first pointed out by Sussmann in~\cite{sussmann1990limitations}, the so called
        \emph{peaking-phenomenon} can cause the non-linear part~\eqref{eq:SusmannPeakingNonLinearPart} of
        the system to have finite escape time even if the system
        \[
            \dot{\eta}(t)=f(0,\eta(t))
        \]
        has $0$ as a global asymptotically stable equilibrium.
        The presumed BIBO property~\ref{Item:OperatorPropBIBO} of operator~$\oT$ not only avoids this problem but is even more
        essential since our control objective is to guarantee that the output~$y$ of the system~\eqref{eq:Sys} respectively 
        the output~$\yM$ of the model~\eqref{eq:Model_r} evolves within
        the funnel around the reference signal~$y_{\rf}$. Without this assumption and even with
        perfect tracking, the non-linear dynamics~\eqref{eq:SusmannPeakingNonLinearPart} might be unbounded and thus cause an
        unbounded control effort, or worse, its solution might even have finite escape time. 
        \item\label{Item:Rem:PropertiesOperator:MemoryLimit} Compared to previous works on funnel control, see e.g.~\cite{IlchRyan02b,BergIlch21,BergIlch23},
        the limited memory property~\ref{Item:OperatorPropLimitMemory} was newly introduced in~\cite{BergDenn24}
        and is essential in the context of MPC in order to ensure that in Step~\ref{agostep:MPCFirst}
        of each iteration of the MPC~\Cref{Algo:MPC} only the history 
        of the state of length up to the memory limit $\tau\ge 0$ is utilised, 
        instead of requiring the full signal history, which would be infeasible in practice.
        Let~$\oT\in\cT_{t_0}^{n,q}$ with memory limit $\tau\geq 0$ and $I$ be the interval~$I=[\hat{t}-\tau,\hat{t}+T]$ 
        for $\hat{t}\geq\tau$ and $T\in\Rp$.
        For given functions~$\zeta\in\cR(I,\R^n)$ and  $\hat{\oT}\in L^\infty_{\loc} ([\hat{t}-\tau,\hat{t}]\cap[t_0,\hat{t}], \R^{q})$, 
        let $\prescript{e}{}{\zeta}$ be a left extension of $\zeta$ on the interval $[0,\hat{t}]$ with $\oT(\prescript{e}{}{\zeta})|_{[\hat{t}-\tau,\hat{t}]\cap[t_0,\hat{t}]}=\hat{\oT}$.
        Similar to the observations in~\ref{Item:Rem:PropertiesOperator:Causality}, the restriction of $\oT(\prescript{e}{}{\zeta})$ to the interval~$I$
        is uniquely determined by the functions $\zeta$ and $\hat{\oT}$ in the sense that $\oT(\prescript{e}{}{\zeta})|_{I}$ 
        is independent of the chosen left extension~$\prescript{e}{}{\zeta}$ due to the limited memory property~\ref{Item:OperatorPropLimitMemory}.
        If the value $\hat{\oT}$ is fixed, this allows us to write 
        $\oT(\zeta)(s)$ in place of $\oT(\prescript{e}{}{\zeta})(s)$ for $s\in I$ 
        (assuming there exists a left extension $\prescript{e}{}{\zeta}$ with $\oT(\prescript{e}{}{\zeta})|_{[\hat{t}-\tau,\hat{t}]\cap[t_0,\hat{t}]}=\hat{\oT}$).
        \item In the literature on funnel control, the used operator $\oT$  belonging
        to~$\cT_{t_0}^{n,q}$ is usually defined on the space of continuous
        functions~$\cC(\Rp,\R^n)$, see e.g.~\cite{IlchRyan02b,BergIlch21}.
        In the context of MPC, this is too restrictive.
        Since the model is re-initialised with data from system measurements in Step~\ref{agostep:MPCFirst}
        at the beginning of every iteration of the MPC~\Cref{Algo:MPC}, one cannot assume
        continuity of the global solution trajectory of model~\eqref{eq:Model_r}.
        Measurement errors, disturbances, and a potential model-system mismatch will inevitably
        result in discontinuities at the points of model re-initialisation.
        To account for this, \Cref{Def:OperatorClass} generalises the operator's domain to the space
        of regulated functions~$\cR(\Rp,\R^n)$.
    \end{enumerate}
\end{remark}

We now define a concept of a solution for the initial value
problem~\eqref{eq:Model_r}. As we will use this differential equation as a model
for the MPC~\Cref{Algo:MPC}, a certain degree of care is required for this
definition. Since the initial value problem is solved at every iteration of the
algorithm at different time instants $\hat{t}\geq t_0$ with varying initial
values based on system measurements obtained in Step~\ref{agostep:MPCFirst} of
the algorithm, certain theoretical problems arise mainly caused by the domain of
operator~$\oTM$.

\begin{definition}[Model solution]\label{Def:ModSolution}
Given a control function $u\in L^\infty_{\loc}([\hat{t},\infty), \R^m)$,
a regulated function 
${\xM = (x_{{\mathrm{M}},1},\ldots,x_{{\mathrm{M}},r})}$ with $x_{\mathrm{M}, i}: [0,\omega) \to \R^{m}$, $\omega \in(\hat{t},\infty]$, $i=1,\ldots,r$,
is called a solution of the initial value problem~\eqref{eq:Model_r} 
at initial time $\hat{t} \geq t_0\geq0$ and with initial data $\xMh\in\cR([\hat{t}-\tau,\hat{t}]\cap[0,\hat{t}],\R^{rm})$ and
$\oTMh\in L^\infty_{\loc} ([\hat{t}-\tau,\hat{t}]\cap[t_0,\hat{t}], \R^{q})$ with $\tau\geq0$,  if 
\begin{equation} \label{eq:Model_InitialValue}
\begin{aligned}
    \xM|_{[\hat{t}-\tau,\hat{t}]\cap [0,\hat{t}]}&= \xMh,\\
    \oTM(\xM)|_{[\hat{t}-\tau,\hat{t}]\cap [t_0,\hat{t}]}&= \oTMh,
\end{aligned}
\end{equation}
and $\xM\vert_{[\hat{t},\omega)}$ is absolutely continuous such that, for almost all~$t\in[\hat{t},\omega)$, it fulfils 
\begin{equation}\label{eq:ModDiff}
    \begin{bmatrix}
         \dot{x}_{\mathrm{M},1}(t)\\
        \vdots\\
         \dot{x}_{\mathrm{M},r-1}(t)\\
         \dot{x}_{\mathrm{M},r}(t)
    \end{bmatrix}= 
    \begin{bmatrix}
         x_{\mathrm{M},2}(t)\\
        \vdots\\
         x_{\mathrm{M},r}(t)\\
         \fM \big(\oTM(\xM)(t) \big) 
    \end{bmatrix}
    \ +\ 
   \begin{bmatrix}
    0\\   
    \vdots\\
    0\\
    \gM \big(\oTM(\xM)(t)\big)
    \end{bmatrix}
    u(t).
    \\
\end{equation} 

A solution $\xM$ is said to be \emph{maximal} if it has no proper right extension that is also a solution.
A maximal solution is called a \emph{response} of the model associated with $u$ 
and we denote it by~$\xM(\cdot;\hat{t},\xMh,\oTMh,u)$.
Its first component $x_{\mathrm{M},1}$ is  denoted by~$\yM(\cdot;\hat{t},\xMh,\oTMh,u)$.
\end{definition}
Note that, in~\eqref{eq:Model_InitialValue} of \Cref{Def:ModSolution}, we
did not distinguish between the cases $t_0>0$ and $t_0=0$ as
in~\eqref{eq:Model_r}, since a larger variety of cases is possible here.
Essentially, one needs to distinguished whether the interval
$[\hat{t}-\tau,\hat{t}]\cap [0,\hat{t}]$ is a perfect interval or merely a
single point. In the latter case, we will implicitly assume a representation of
the initial condition as in~\eqref{eq:Model_r} and that $\xMh$ is an element of
$\R^{rm}$ and analogously that $\oTMh$ is an element of $\R^q$.
Or in other words, we interpret the restriction of the considered functions to an interval of the form
$I=[\hat{t}]$ as a function evaluation at $\hat{t}$ and identify the function spaces 
$\cR(I,\R^{rm})$ and $L^\infty_{\loc} (I, \R^{q})$ with the vector spaces $\R^{rm}$ and $\R^{q}$, respectively.
The parameter $\tau\geq0$ can be thought of as a memory length of past signal information 
used to initialise the differential equation.
\begin{remark}
    The above~\Cref{Def:ModSolution} contains certain peculiarities and differs from the conventional definition for the solution of an
    initial value problem. We comment on that.
    \begin{enumerate}[(a)]
    \item \Cref{Def:ModSolution} uses a solution concept in the sense of \emph{Carath\'{e}odory}, see e.g.~\cite[\S~10, Supplement II]{Walt98}.
    The solution~$\xM$ does not automatically possess a continuous first derivative and it fulfils differential equation~\eqref{eq:ModDiff}
    on $[\hat{t},\omega)$ with the exception of a set of Lebesgue measure zero.
    Equivalently, the solution concept can be formulated using an integral representation instead of equation~\eqref{eq:ModDiff}.
    Then, the function~$\xM$ is considered to be a solution to differential equation~\eqref{eq:Model_r}
    when satisfying 
    \begin{equation}\label{eq:ModIntegral}
        \xM(t)=\xMh(\hat{t})+\int_{\hat{t}}^t\FM(\xM(s),\oTM(\xM)(s))+\GM(\oTM(\xM)(s))u(s)\d{s},
    \end{equation}
    for all~$t\in[\hat{t},\omega)$, where
    $\FM(\xM,\oTM(\xM))\coloneqq \begin{bmatrix}
         x_{\mathrm{M},2},&
        \ldots,&
         x_{\mathrm{M},r},&
         \fM \big(\oTM(\xM) \big) \big) 
    \end{bmatrix}^\top$
    and 
    $\GM(\oTM(\xM))\coloneqq \begin{bmatrix}
    0, &  
    \hdots,&
    0,&
    \gM \big(\oTM(\xM)\big)
    \end{bmatrix}^\top$. At some instances, we will also use this representation. 
    \item Note that in~\Cref{Def:ModSolution} the solution~$\xM$ is defined on the whole interval $[0,\omega)$
    while the initial values are given at $\hat{t}\geq t_0\geq 0$ and while $\xM$ satisfies
    the differential equation merely on the interval $[\hat{t},\omega)$.
    The rational behind this is for the domain of the solution~$\xM$ to be in accordance 
    with the domain of the operator~$\oTM$ which is $\cR(\Rp,\R^{rm})$.
    By virtue of the causality property~\ref{Item:OperatorPropCasuality},
    one can, when solving the initial value problem~\eqref{eq:ModDiff} with~\eqref{eq:Model_InitialValue}, 
    evaluate $\oTM(\xM)(t)$ for $t\in[\hat{t},\omega)$, in a certain sense,
    implicitly utilising right extensions $\xM^{e}$ as discussed in \Cref{Rem:PropertiesOperator}~\ref{Item:Rem:PropertiesOperator:Causality}.
    However, it is necessary for $\xM$ to be also defined in the past, meaning on the interval~$[0,\hat{t}]$.
    The initial values $\xMh$ and $\oTMh$ determine, not necessarily uniquely,  $\xM$ on the interval~$[0,\hat{t}]$.
    \item
    It might seem, at first glance, like $\oTMh$ is already completely determined by $\xMh$. 
    The latter however is, in general,  not an element of the domain of~$\oTM$. 
    In order for $\xMh$ to be evaluable by the operator $\oTM$,
    one has to choose a left extension $\prescript{e}{}{\hat{x}}_{\mathrm{M}}$ on the interval~$[0,\hat{t}]$.
    As this left extension can initially be chosen arbitrarily, 
    the second condition on the initial value in~\eqref{eq:Model_InitialValue} 
    entwines the left extension $\prescript{e}{}{\hat{x}}_{\mathrm{M}}$ with the initial datum~$\oTMh$.
    As discussed in~\Cref{Rem:PropertiesOperator}~\ref{Item:Rem:PropertiesOperator:MemoryLimit}, 
    if $\tau\geq 0$ is greater than or equal to the memory limit of $\oTM$, then the left extension 
    $\prescript{e}{}{\hat{x}}_{\mathrm{M}}$, assuming its existence, is uniquely determined by 
    $\xMh$ and $\oTMh$ in the sense that $\oTM(\xM)(t)$ for $t\geq\hat{t}$ does not depend on the chosen extension.
    \item 
    To take the initial value $\yM^0$ from~\eqref{eq:Model_r} into account at $\hat{t}=t_0$,
    one can replace in \eqref{eq:Model_InitialValue} the initial value $\xMh$ with  $(\yM^0,\ldots,{\yM^0}^{(r-1)})$  
    and $\oTMh$ with $\oTM\big((\yM^0,\ldots,{\yM^0}^{(r-1)})\big)$.
    However, \Cref{Def:ModSolution} is intentionally formulated merely in terms of $\xMh$ and $\oTMh$
    independent of $\yM^0$ in order to avoid treating the initial time $\hat{t}=t_0$ as a special case.
    \end{enumerate}
\end{remark}

One of the difficulties we will have to address in the following is that the
initial values $\xMh$ and $\oTM$ in~\eqref{eq:Model_InitialValue} have to allow
for the existence of a solution of initial value problem~\eqref{eq:ModDiff}.
For that $\xMh$ is required to have an admissible left extension on the interval~$[0,\hat{t}]$, i.e. 
a left extension $\prescript{e}{}{\hat{x}}_{\mathrm{M}}$ with $\oTM(\prescript{e}{}{\hat{x}}_{\mathrm{M}})|_{[\hat{t}-\tau,\hat{t}]\cap[t_0,\hat{t}]}=\oTMh$.
Moreover, in the context of the MPC~\Cref{Algo:MPC}, the initial values have to be chosen in a way 
to entwine the initial value problem~\eqref{eq:ModDiff} to be solved at the current time instant~$\hat{t}$
with the history of $x_{\mathrm{M}}$ from the previous iterations and the initial trajectory $\yM^0$ given at time $t_0\geq 0$.
Before addressing these two questions, we want to consider initial value problem~\eqref{eq:ModDiff} individually, independent of the
MPC~\Cref{Algo:MPC}, and assume the existence of an admissible left extension for $\xMh$.
In this case, there exists a solution to the initial value problem as the following~\Cref{Prop:SolutionExists} shows.

\begin{prop}\label{Prop:SolutionExists}
    Consider model~\eqref{eq:Model_r} with $(\fM,\gM,\oTM)\in\cM^{m,r}_{t_0}$ and $t_0\geq0$.
    For $\hat{t}\geq t_0$ and a control function $u\in L^\infty_{\loc}([\hat{t},\infty), \R^m)$,
    let $\xMh\in\cR([\hat{t}-\tau,\hat{t}]\cap[0,\hat{t}],\R^{rm})$ and $\oTMh\in
    L^\infty_{\loc} ([\hat{t}-\tau,\hat{t}]\cap[t_0,\hat{t}], \R^{q})$ with $\tau\geq0$.
    If $\xMh$ has an admissible left extension $\prescript{e}{}{\hat{x}}_{\mathrm{M}}$  on the interval~$[0,\hat{t}]$, i.e. 
    $\prescript{e}{}{\hat{x}}_{\mathrm{M}}\in\cR([0,\hat{t}],\R^{rm})$ with 
    $\prescript{e}{}{\hat{x}}_{\mathrm{M}}|_{[\hat{t}-\tau,\hat{t}]\cap[0,\hat{t}]}=\xMh$ and
    $\oTM(\prescript{e}{}{\hat{x}}_{\mathrm{M}})|_{[\hat{t}-\tau,\hat{t}]\cap[t_0,\hat{t}]}=\oTMh$, then
\begin{enumerate}[label = (\roman{enumi}), ref=(\roman{enumi})]
    \item\label{Item:Theorem:ExistenceSolution} the initial value problem~\eqref{eq:ModDiff} with~\eqref{eq:Model_InitialValue}  has a solution $\xM:[0,\omega)\to \R^{rm}$ in the sense of~\Cref{Def:ModSolution},
    \item every solution can be extended to a maximal solution,
    \item\label{Item:Theorem:BoundedSolution} if $\xM:[0,\omega)\to \R^{rm}$ is a bounded maximal solution, then $\omega=\infty$.
\end{enumerate}
\end{prop}
\begin{proof}
    For given $\hat{t}\geq t_0$ and $u\in L^\infty_{\loc}([\hat{t},\infty), \R^m)$,
    let $\xMh\in\cC([\hat{t}-\tau,\hat{t}]\cap[0,\hat{t}],\R^{rm})$ and $\oTMh\in
    L^\infty_{\loc} ([\hat{t}-\tau,\hat{t}]\cap[t_0,\hat{t}], \R^{q})$ with $\tau\geq0$ be arbitrary but fixed.
    Further, let $\prescript{e}{}{\hat{x}}_{\mathrm{M}}$ be an admissible left extension of $\xMh$ on the interval~$[0,\hat{t}]$
    and let $\FM$ and $\GM$ be defined as in~\eqref{eq:ModIntegral}.
    As a consequence of~\Cref{Appendix:Cor:ModelSolutionExists} to be found in
    the \nameref{Chapter:Appendix},
    there exists a solution $x:[0,\omega)\to \R^{rm}$ with $\omega>\hat{t}$ of the initial value problem
    \begin{align}
        \dot{x}(t)&=\FM(x(t),\oTM(x)(t))+\GM(\oTM(x)(t))u(t),\label{eq:Prop:SolutionExists:Equation}\\
        x|_{[0,\hat{t}]}&=\prescript{e}{}{\hat{x}}_{\mathrm{M}}\nonumber,
    \end{align}
    where $x$ fulfils the differential equation for almost all $t\in [\hat{t},\omega)$. 
    Since $\prescript{e}{}{\hat{x}}_{\mathrm{M}}$ is an admissible left extension of $\xMh$, the function $x$ is also a solution of 
    initial value problem~\eqref{eq:ModDiff} with~\eqref{eq:Model_InitialValue} in the sense of~\Cref{Def:ModSolution}.
    This shows~\ref{Item:Theorem:ExistenceSolution}.

    Let $\tilde{x}:[0,\omega)\to \R^{rm}$ be an arbitrary solution of the initial value problem~\eqref{eq:ModDiff} in the sense of~\Cref{Def:ModSolution}
    with initial values $\xMh$ and $\oTMh$ as in~\eqref{eq:Model_InitialValue}.
    Due to the used solution concept,  the function $\tilde{x}$ is also a solution of the 
    initial value problem~\eqref{eq:Prop:SolutionExists:Equation} with initial value $x|_{[0,\hat{t}]}=\tilde{x}|_{[0,\hat{t}]}$.
    According to~\Cref{Appendix:Cor:ModelSolutionExists}, this function can be extended to a maximal solution.
    Moreover, if any maximal solution of this initial value problem is bounded, then $\omega=\infty$.
    Since $\tilde{x}|_{[0,\hat{t}]}$ is an admissible left extension of $\xMh$, both findings carry over to the 
    initial value problem~\eqref{eq:ModDiff} with initial values as in~\eqref{eq:Model_InitialValue} and 
    the solution definition in the sense of~\Cref{Def:ModSolution}.
\end{proof}

    While \Cref{Prop:SolutionExists} ensures the existence of a solution $\xM$ of initial value problem~\eqref{eq:ModDiff} with \eqref{eq:Model_InitialValue},
    it is in general not unique, in particular since $\xM$ does not need to comply with the differential equation~\eqref{eq:ModDiff} for $t<\hat{t}$.
    On the interval $[0,\hat{t}]$, it is merely a feasible left extension of $\xMh$ and, depending on $\oTM$, there might exist several feasible left extensions.
    However, if $\tau$ is greater than or equal to the memory limit of $\oTM$, then $\xM$ is, in a certain sense, unique as the following \Cref{Prop:SolutionUnique} shows.
    Namely, it is uniquely determined for $t\geq\hat{t}$.
    
\begin{prop}\label{Prop:SolutionUnique}
    Consider model~\eqref{eq:Model_r} with $(\fM,\gM,\oTM)\in\cM^{m,r}_{t_0}$ and $t_0\geq0$.
    For ${\hat{t}\geq t_0}$ and a control function $u\in L^\infty_{\loc}([\hat{t},\infty), \R^m)$,
    let  $\xMh\in\cR([\hat{t}-\tau,\hat{t}]\cap[0,\hat{t}],\R^{rm})$ and $\oTMh\in
    L^\infty_{\loc} ([\hat{t}-\tau,\hat{t}]\cap[t_0,\hat{t}], \R^{q})$ with $\tau\geq0$.
    Let $\xM^1:[0,\omega_1)\to \R^{rm}$ and ${\xM^2:[0,\omega_2)\to \R^{rm}}$ 
    be two solutions of the initial value problem~\eqref{eq:ModDiff} with~\eqref{eq:Model_InitialValue} in the sense of~\Cref{Def:ModSolution}.
    If the value $\tau$ is greater than or equal to the memory limit of $\oTM$, then 
    $\xM^1(t)=\xM^2(t)$ for all $t\in[\hat{t},\min\{\omega_1,\omega_2\})$.
\end{prop}

\begin{proof}
    \emph{Step 1}: We show that $\xM^1$ and $\xM^2$ coincide on the interval $[\hat{t},\hat{t}+\eps]$ for some~$\eps>0$.
    As $\xM^1$ and $\xM^2$ are solutions of the initial value problem~\eqref{eq:ModDiff} 
    with~\eqref{eq:Model_InitialValue} in the sense of~\Cref{Def:ModSolution},
    there exists a feasible left extension of $\xMh$, i.e. $\prescript{e}{}{\hat{x}}_{\mathrm{M}}\in\cR([0,\hat{t}],\R^{rm})$ with 
    $\prescript{e}{}{\hat{x}}_{\mathrm{M}}|_{[\hat{t}-\tau,\hat{t}]\cap[0,\hat{t}]}=\xMh$ and
    $\oTM(\prescript{e}{}{\hat{x}}_{\mathrm{M}})|_{[\hat{t}-\tau,\hat{t}]\cap[t_0,\hat{t}]}=\oTMh$.
    According to the local Lipschitz property~\ref{Item:OperatorPropLipschitz} of operator $\oTM$, there exists
    constants $\Delta, \delta, c > 0$  such that for all functions
    ${y_1, y_2 \in \cR(\Rp,\R^{rm})}$ with
    $y_1|_{[0,\hat{t}]} = y_2|_{[0,\hat{t}]} = \prescript{e}{}{\hat{x}}_{\mathrm{M}}$ 
    and $\Norm{y_1(s) - \xMh(\hat{t})} < \delta$,  $\Norm{y_2(s) - \xMh(\hat{t})} < \delta $ for all $s \in [t,t+\Delta]$:
    \[
     \esssup_{\mathclap{s \in [t,t+\Delta]}}  \Norm{\oTM(y_1)(s) - \oTM(y_2)(s) }  
        \le c \ \sup_{\mathclap{s \in [t,t+\Delta]}}\ \Norm{y_1(s)- y_2(s)}.
    \] 
    The functions $\xM^1$ and $\xM^2$ are continuous  on the interval $[\hat{t},\min\{\omega_1,\omega_2\})$ and satisfy ${\xM^1(\hat{t})=\xM^2(\hat{t})=\xMh(\hat{t})}$.
    Therefore, there exists some $\eps\in(0,\min\{\Delta,\omega_1,\omega_2\})$ such that 
    $\Norm{\xM^1(t)-\xMh(\hat{t})}<\delta$ and $\Norm{\xM^2(t)-\xMh(\hat{t})}<\delta$ for all $t\in[\hat{t},\hat{t}+\eps]$.
    
    We will estimate $\Norm{\oTM(\xM^{1})(t)-\oTM(\xM^{2})(t)}$ for $t\in[\hat{t},\hat{t}+\eps]$ in the following.
    To that end, let $\Indic_{I}$ denote the indicator function of an interval $I\subset\R$.
    Note that, as $\tau$ is greater than or equal to the memory limit  of $\oTM$, see property \ref{Item:OperatorPropLimitMemory}, 
    we have, for $t\geq\hat{t}$ and $i=1,2$,
    \[ 
        \oTM(\xM^{i})(t)=\oTM\rbl\Indic_{[0,\hat{t})}\prescript{e}{}{\hat{x}}+\Indic_{[\hat{t},t]}\xM^{i}\rbr(t)
        =\oTM\rbl\Indic_{[0,\hat{t})}\prescript{e}{}{\hat{x}}+\Indic_{[\hat{t},t]}\xM^{i}+\Indic_{(t,\hat{t}+\Delta]}\xMh(\hat{t})\rbr(t),
    \]
    where the causality property~\ref{Item:OperatorPropCasuality} was used in the second equation.
    Hence, for $t\in[\hat{t},\hat{t}+\eps]$, it follows
    \begin{small}
    \begin{align*}
        &\Norm{\oTM(\xM^1)(t)-\oTM(\xM^2)(t)}
        \leq\esssup_{\nu\in[\hat{t},t]}\Norm{\oTM(\xM^1)(\nu)-\oTM(\xM^2)(\nu)}\\
        &=\!\textover{$\esssup\limits_{\nu\in[\hat{t},t]}$}{$\esssup\limits_{\nu\in[\hat{t},\hat{t}+\Delta]}$}\!\!\!\!\!\Norm{
         \oTM\!\!\left(\Indic_{[0,\hat{t})}\prescript{e}{}{\hat{x}}\!+\!\Indic_{[\hat{t},t]}\xM^{1}\!+\!\Indic_{(t,\hat{t}+\Delta]}\xMh(\hat{t})\right)\!\!(\nu)
        \!-\!\oTM\!\!\left(\Indic_{[0,\hat{t})}\prescript{e}{}{\hat{x}}\!+\!\Indic_{[\hat{t},t]}\xM^{2}\!+\!\Indic_{(t,\hat{t}+\Delta]}\xMh(\hat{t})\right)\!\!(\nu)}\\
        &\!\leq\!\!\!\esssup_{\nu\in[\hat{t},\hat{t}+\Delta]}\!\Norm{
         \oTM\!\!\left(\Indic_{[0,\hat{t})}\prescript{e}{}{\hat{x}}\!+\!\Indic_{[\hat{t},t]}\xM^{1}\!+\!\Indic_{(t,\hat{t}+\Delta]}\xMh(\hat{t})\right)\!\!(\nu)
        \!-\!\oTM\!\!\left(\Indic_{[0,\hat{t})}\prescript{e}{}{\hat{x}}\!+\!\Indic_{[\hat{t},t]}\xM^{2}\!+\!\Indic_{(t,\hat{t}+\Delta]}\xMh(\hat{t})\right)\!\!(\nu)}\\
        &\leq c\!\!\!\!\sup_{\nu\in[\hat{t},\hat{t}+\Delta]}\Norm{
         \rbl\Indic_{[0,\hat{t})}\prescript{e}{}{\hat{x}}+\Indic_{[\hat{t},t]}\xM^{1}+\Indic_{(t,\hat{t}+\Delta]}\xMh(\hat{t})\rbr(\nu) 
        -\rbl\Indic_{[0,\hat{t})}\prescript{e}{}{\hat{x}}+\Indic_{[\hat{t},t]}\xM^{2}+\Indic_{(t,\hat{t}+\Delta]}\xMh(\hat{t})\rbr(\nu)
        } \\
        &= c\sup_{\nu\in[\hat{t},t]}\Norm{\Indic_{[\hat{t},t]}(\nu)\xM^{1}(\nu)-\Indic_{[\hat{t},t]}(\nu)\xM^{2}(\nu)}
         = c\sup_{\nu\in[\hat{t},t]}\Norm{\xM^{1}(\nu)-\xM^{2}(\nu)}.
    \end{align*}
    \end{small}
    We now estimate $\sup_{\nu\in[\hat{t},\hat{t}+\eps]}\Norm{\xM^1(\nu)-\xM^2(\nu)}$.
    To that end, define the continuous function $\zeta:[\hat{t},\hat{t}+\eps]\to\R$, $t\mapsto\sup_{\nu\in[\hat{t},t]}\Norm{\xM^1(\nu)-\xM^2(\nu)}$.
    There exists a compact set $K$ with $\xM^{i}(t)\in K$ for all $t\in[\hat{t},\hat{t}+\eps]$ and $i=1,2$.
    Moreover, due to the BIBO property~\ref{Item:OperatorPropBIBO} of operator $\oTM$, 
    there exists a compact set $\tilde{K}$ with $\oTM(\xM^{i})(t)\in K$  
    for all $t$ in $[\hat{t},\hat{t}+\eps]$ and $i=1,2$.
    As the function $\fM$ is an element of $\Lip_{\loc}(\R^q,\R^m)$, $\FM$ is
    Lipschitz continuous with constant $L_{\FM}\geq0$ on the set
    $K\times\tilde{K}$. Similarly, there exists a Lipschitz constant
    $L_{\GM}\geq0$ for the function $\GM$ on the set $\tilde{K}$. 
    Using the above considerations and the solution representation~\eqref{eq:ModIntegral},
    $\zeta$ satisfies the following estimate for all $t$ in $[\hat{t},\hat{t}+\eps]$.
    \begin{align*}
        \zeta(t)
        &=\sup_{\nu\in[\hat{t},t]}\Norm{\xM^1(\nu)-\xM^2(\nu)}\\
        &=\sup_{\nu\in[\hat{t},t]}\left\|\xMh(\hat{t})+\int_{\hat{t}}^\nu\FM(\xM^1(s),\oTM(\xM^1)(s))+\GM(\oTM(\xM^1)(s))u(s)\d{s}\right.\\
                               &\hphantom{=\sup_{\nu\in[\hat{t},t]}}\hspace{5pt-\widthof{$-$}}
                               -\left.\xMh(\hat{t})-\int_{\hat{t}}^\nu\FM(\xM^2(s),\oTM(\xM^2)(s))+\GM(\oTM(\xM^2)(s))u(s)\d{s}\right\|\\
        &\leq\sup_{\nu\in[\hat{t},t]}\left(\int_{\hat{t}}^{\nu}\Norm{\FM(\xM^1(s),\oTM(\xM^1)(s))-\FM(\xM^2(s),\oTM(\xM^2)(s))}\d{s}\right.\\
                               &\hphantom{=\sup_{\nu\in[0,t]}}\hspace{0pt}\left.+\int_{\hat{t}}^\nu\Norm{\GM(\oTM(\xM^1)(s))u(s)-\GM(\oTM(\xM^2)(s))u(s)}\d{s}\right)\\
        &\leq\int_{\hat{t}}^t\Norm{\FM(\xM^1(s),\oTM(\xM^1)(s))-\FM(\xM^2(s),\oTM(\xM^2)(s))}\d{s}\\
                               &\hphantom{=}+\int_{\hat{t}}^t\Norm{\GM(\oTM(\xM^1)(s))u(s)-\GM(\oTM(\xM^2)(s))u(s)}\d{s}\\
        &\leq\int_{\hat{t}}^tL_{\FM}\rbl\Norm{\xM^1(s)-\xM^2(s)}+\Norm{\oTM(\xM^1)(s))-\oTM(\xM^2)(s))}\rbr\d{s}\\
                               &\hphantom{=}+\int_{\hat{t}}^tL_{\GM}\Norm{\oTM(\xM^1)(s))-\oTM(\xM^2)(s))}\SNorm{u|_{[\hat{t},t]}}\d{s}\\
        &\leq \rbl L_{\FM}(1+c)+L_{\GM}c\SNorm{u|_{[\hat{t},\hat{t}+\eps]}}\rbr\int_{\hat{t}}^t\sup_{s\in[\hat{t},s]}\Norm{\xM^1(s)-\xM^2(s)}\d{s}\\
        &\leq \rbl L_{\FM}(1+c)+L_{\GM}c\SNorm{u|_{[\hat{t},\hat{t}+\eps]}}\rbr\int_{\hat{t}}^t\zeta(s)\d{s}.
    \end{align*}
    Now, Grönwall's inequality yields $\zeta(t)=0$ for all $t\in[\hat{t},\hat{t}+\eps]$. This shows $\xM^1(t)=\xM^2(t)$ on the interval $t\in[\hat{t},\hat{t}+\eps]$.

    \noindent
    \emph{Step 2}: We show that $\xM^1$ and $\xM^2$ coincide on the interval $[\hat{t},\min\{\omega_1,\omega_2\})$.
    Suppose  $\xM^1|_{[\hat{t},\min\{\omega_1,\omega_2\})}\neq \xM^2|_{[\hat{t},\min\{\omega_1,\omega_2\})}$. 
    Then, there exists a time instant $t\in[\hat{t},\min\{\omega_1,\omega_2\})$ with ${\xM^1(t)\neq\xM^2(t)}$.
    Let 
    \[
        \tilde{t}\coloneqq \inf \setdef{t\in [\hat{t},\min\{\omega_1,\omega_2\})}{\xM^1(t)\neq\xM^2(t)}.
    \]
    In view of Step~1, we have $\tilde{t}>\hat{t}$. 
    Since $\xM^1$ and $\xM^2$ are solutions of the initial value problem~\eqref{eq:Model_r} with initial value $\xMh$ and $\oTMh$ at initial time $\hat{t}$, 
    we have 
    \[
    \xM^1|_{[\hat{t}-\tau,\hat{t}]\cap [0,\hat{t}]}= \xMh=\xM^2|_{[\hat{t}-\tau,\hat{t}]\cap [0,\hat{t}]}
    \text{ and } \oTM(\xM^1)|_{[\hat{t}-\tau,\hat{t}]\cap [t_0,\hat{t}]}= \oTMh=\oTM(\xM^2)|_{[\hat{t}-\tau,\hat{t}]\cap [t_0,\hat{t}]}.
    \]
    Thus, $\xM^1(t)=\xM^2(t)$ for all $t\in[\hat{t}-\tau,\tilde{t}]\cap [0,\tilde{t}]$.
    This implies $\oTM(\xM^1)(t)=\oTM(\xM^2)(t)$ for all $t\in[\tilde{t}-\tau,\tilde{t}]\cap[t_0,\tilde{t}]$
    because $\tau$ is greater than or equal to the memory limit  of $\oTM$, see property~\ref{Item:OperatorPropLimitMemory}.
    Define $\xMt\coloneqq \xM^1|_{[\tilde{t}-\tau,\tilde{t}]\cap[0,\tilde{t}]}$
    and $\oTMt\coloneqq \oTM(\xM^1)|_{[\tilde{t}-\tau,\tilde{t}]\cap[t_0,\tilde{t}]}$.
    These functions are elements of
    $\cR([\tilde{t}-\tau,\tilde{t}]\cap[0,\tilde{t}],\R^{rm})$ and
    $L^\infty_{\loc} ([\tilde{t}-\tau,\tilde{t}]\cap[t_0,\tilde{t}], \R^{q})$,
    respectively.
    Both $\xM^1$ and $\xM^2$ are solutions to the initial value
    problem~\eqref{eq:Model_r} with initial value $\xMt$ and $\oTMt$ at initial
    time $\tilde{t}$ in the sense of~\Cref{Def:ModSolution}. According to {Step
    1}, there exists $\eps>0$ such that $\xM^1(t)=\xM^2(t)$ for all
    $t\in[\tilde{t},\tilde{t}+\eps]$~--~a contradiction to the definition of
    $\tilde{t}$. This completes the proof.
\end{proof}

\begin{remark}
    \Cref{Prop:SolutionUnique} shows that if $\tau$ is greater than or equal to the memory limit of $\oTM$,
    then the maximal solution $\xM:[0,\omega)\to \R^{rm}$ 
    of the initial value problem~\eqref{eq:Model_r} with initial value $\xMh$ and $\oTMh$ at time $\hat{t}$
    is uniquely determined on the interval $[\hat{t},\omega)$.
    Since we will consider solutions of the initial value problem~\eqref{eq:ModDiff} mostly for $t\geq \hat{t}$, 
    we will speak in this case also of \emph{the} maximal solution and \emph{the} response associated with $u$
    when referring to~$\xM(\cdot;\hat{t},\xMh,\oTMh,u)$.
\end{remark}

In the previous~\Cref{Sec:FunnelStageCostFunctions}, we introduced the concept of funnel penalty function to be used in MPC.
Basis for our considerations was the assumption of a Lipschitz continuous solution trajectory of the model~\eqref{eq:Intro:ModelEquation}.
The following~\Cref{Prop:SolutionIsLPath} shows that this assumption is justified and fulfilled for our model class~$\cM^{m,r}_{t_0}$.
\begin{prop}\label{Prop:SolutionIsLPath}
    Consider model~\eqref{eq:Model_r} with $(\fM,\gM,\oTM)\in\cM^{m,r}_{t_0}$.
    Let $\hat{t}\geq t_0$, $\tau\in\Rp$, and initial data $\xMh\in\cR([\hat{t}-\tau,\hat{t}]\cap[0,\hat{t}],\R^{rm})$ and
    $\oTMh\in L^\infty_{\loc} ([\hat{t}-\tau,\hat{t}]\cap[t_0,\hat{t}], \R^{q})$ such that for a control
    $u\in L^\infty_{\loc}([\hat{t},\infty), \R^m)$ the  
    initial value problem~\eqref{eq:Model_r} has a solution $\xM:[0,\omega)\to \R^{rm}$ with $\omega>\hat{t}$ in the sense of~\Cref{Def:ModSolution}.
    Then, for every interval length $T\in(0,\omega-\hat{t})$, the restriction $\xM|_{[\hat{t},\hat{t}+T]}:[\hat{t},\hat{t}+T]\to\R^{rm}$
    is a Lipschitz path.
\end{prop}
\begin{proof}
    We prove the assertion by showing that every component $x_{\mathrm{M}, i}:[0,\omega)\to\R^m$  
    of $\xM=(x_{{\mathrm{M}},1},\ldots,x_{{\mathrm{M}},r})$ for $i=1,\ldots, r$ 
    is a Lipschitz continuous function on the interval $[\hat{t},\hat{t}+T]$.
    By \Cref{Def:ModSolution}, the function $x_{\mathrm{M}, i}$ is continuous on the interval $[\hat{t},\hat{t}+T]$.
    It therefore is sufficient to show that $\dot{x}_{\mathrm{M}, i}$ is an essentially bounded function on the interval $[\hat{t},\hat{t}+T]$.
    For $i=1,\ldots, r-1$,  we have $\dot{x}_{\mathrm{M}, i} =x_{\mathrm{M}, i+1}$ since $\xM$ fulfils the ordinary differential equation~\eqref{eq:ModDiff} 
    on the interval $[\hat{t},\hat{t}+T]$. Due to the compactness of $[\hat{t},\hat{t}+T]$, the continuous function $x_{\mathrm{M}, i+1}$ is bounded.
    For $i=r$, we have 
    \[
        \dot{x}_{\mathrm{M},r}(t) =
        \fM\big(\oTM(\xM)(t) \big) + \gM \big(\oTM(\xM)(t)\big) u(t)
    \]
    for almost all $t\in[\hat{t},\hat{t}+T]$.
    The control $u$ as an element of $L^\infty_{\loc}([\hat{t},\infty), \R^m)$ is bounded. 
    Moreover, due to the compactness of the considered interval and the continuity of the involved functions, 
    $\fM\big(\oTM(\xM)(t) \big)$ and $\gM \big(\oTM(\xM)(t)\big)$ are bounded.
    Thus, $\dot{x}_{\mathrm{M},r}$ is essentially bounded and the proof is complete.
\end{proof}

\section{MPC with funnel stage costs}\label{Sec:MPCWithFunnelStageCost}
In this section, we will analyse how funnel stage cost functions (introduced in~\Cref{Sec:FunnelStageCostFunctions})
can be integrated into model predictive control (MPC) to solve the reference tracking problem formulated in~\Cref{Sec:ControlObjective}.
We employ a (functional) differential equation of the form~\eqref{eq:Model_r},
belonging to the model class~$\cM^{m,r}_{t_0}$ introduced in~\Cref{Sec:ModelClass}, 
as the predictive model within the MPC framework.
To establish \emph{initial and recursive feasibility}  of the resulting MPC scheme and its compliance 
with the control objective, we develop the theoretical groundwork in three key steps:
\begin{enumerate}
    \item \textbf{Relative degree analysis} (\Cref{Sec:HighRelativeDegree}): 
        We investigate the role of the model's relative degree $r\in\N$
        in~\eqref{eq:Model_r} and introduce auxiliary error signals to reduce the complexity of the considered control objective.
    \item \textbf{Feasibility guarantees} (\Cref{Sec:InitialRecFeasibilty}):
        We address the existence of control signals that solve the tracking problem at every iteration of the MPC~\Cref{Algo:MPC}.
    \item \textbf{Optimal control problem solvability} (\Cref{Sec:OptimalControlProblem}):
        We prove that the optimal control problem using funnel stage costs admits a solution, which inherently satisfies the tracking objective. 
\end{enumerate}
While the existence of a solution may appear purely technical, 
it is non-trivial due to the inherent challenges of funnel stage costs: these functions are highly non-linear and generally discontinuous.
Finally, in \Cref{Sec:FunnelMPC}, we synthesise these results into the funnel MPC~\Cref{Algo:FunnelMPC}, 
ensuring adherence to funnel constraints.

\subsection{The higher relative degree}\label{Sec:HighRelativeDegree}
We now develop a control framework to address the reference tracking problem outlined in~\Cref{Sec:ControlObjective},
accounting for the relative degree~$r\in\N$ of the model~\eqref{eq:Model_r}.
While the relative degree~$r$ might initially appear to be a minor technicality 
-- and extending control strategies from $r=1$ to $r>1$ seemingly straightforward -- 
the structural complexity introduced by higher relative degrees poses significant analytical and design challenges.
This difficulty is well-documented: in adaptive control, these challenges were highlighted in~\cite{Mors96}.
Concerning funnel control, the progress was incremental. First proposed for relative degree $r=1$ systems in 
2002 in~\cite{IlchRyan02b}, it took eleven years to extend the framework to $r=2$ in~\cite{HackHopf13} 
and further five years to achieve generalisation for arbitrary $r\in\N$ in~\cite{BergLe18a}.

To meet the control objective in \Cref{Sec:ControlObjective}, 
we introduce auxiliary error variables, 
circumventing the structural limitations imposed by higher relative degrees. 
This approach simplifies the design while ensuring compatibility with the funnel 
stage cost functions discussed in \Cref{Sec:FunnelStageCostFunctions}.
Define, for $(z_1,\ldots,z_r)\in\R^{rm}$ with $z_i\in\R^m$ and for parameters 
$k_1,\ldots, k_{r-1}\in\Rp$, the functions $\eM_i:\R^{rm}\to\R^m$ recursively by 
\begin{equation} \label{eq:ErrorVar}
    \begin{aligned}
        \eM_1 (z_1,\ldots,z_r)&\coloneqq  z_1,\\
        \eM_{i+1} (z_1,\ldots,z_r)&\coloneqq \eM_{i}(\LShift(z_1,\ldots,z_r))+k_i \eM_{i}(z_1,\ldots,z_r),
    \end{aligned}
\end{equation}%
for $i=1,\ldots,r-1$, where 
\begin{equation}\label{eq:LeftShift}
    \LShift:\R^{rm}\to\R^{rm},\ \LShift(z_1,\ldots,z_r)\coloneqq  (z_2,\ldots,z_r,0)
\end{equation}
is the left shift operator.

\begin{remark}
Using the shorthand notation 
\begin{equation}\label{eq:DefOperatorChi}
    \OpChi(\zeta)(t)\coloneqq (\zeta(t),\dot{\zeta}(t),\ldots,\zeta^{(r-1)}(t))\in\R^{rm}    
\end{equation}
for a function  $\zeta\in W^{r,\infty}(I,\R^m)$ on an interval $I\subset\Rp$ and $t\in I$, we get 
\begin{equation} \label{eq:ErrorVarDyn}
    \begin{aligned}
    \eM_1    (\OpChi(\zeta)(t))&=\zeta(t),\\
    \eM_{i+1}(\OpChi(\zeta)(t))&=\dd{t}\eM_{i}(\OpChi(\zeta)(t))+k_i \eM_{i}(\OpChi(\zeta)(t))
    \end{aligned}
\end{equation}
for $i= 1,\ldots,r-1$.
Furthermore, using the polynomials $p_i(s)=\prod_{j=1}^i(s+k_j)\in\R[s]$, 
the function $\eM_{i+1}(\OpChi(\zeta)(t))$ can be represented as
\[
    \eM_{i+1}(\OpChi(\zeta)(t))=p_i(\dd{t})\zeta(t)
\]
for $i=1,\ldots, r-1$.
\end{remark}
Observe that for a solution $\xM$ of the model differential equation~\eqref{eq:Model_r},
the auxiliary error variable $\eM_1(\xM(t)-\OpChi(y_{\rf})(t))$ coincides with the tracking error $\eMTrack(t)=\yM(t)-y_{\rf}(t)$.
Leveraging this equivalence, we solve the tracking problem outlined in~\Cref{Sec:ControlObjective} 
by ensuring $\xM(t)-\OpChi(y_{\rf})(t)$ remains within the set
\begin{equation} \label{eq:DefSetD}
    \cD^{\Psi}_{t}\coloneqq \setdef
    {z\in\R^{rm}}
    {
        \Norm{\eM_{i}(z)}<\Funnel_i(t),\ i=1,\ldots, r
    } 
\end{equation}
for all $t\geq t_0$, where $\Psi\coloneqq (\Funnel_1,\ldots,\Funnel_r)\in\cG^r$
is a vector of suitable funnel functions.
While this approach initially appears to compound the original problem 
-- replacing a single constraint with $r$ time-variant inequalities --
it simplifies the task when the funnel functions are strategically designed. 
Crucially, if a control $u\in L^\infty_{\loc}([\hat{t},\infty), \R^m)$ applied
to the model~\eqref{eq:Model_r} ensures that $\xM(t)-\OpChi(y_{\rf})(t)$ is an
element of the set $\cD^{\Psi}_{t}$ for all $t\geq t_0$, then $\xM(t)$ remains
bounded. By \Cref{Prop:SolutionExists}, this guarantees that a maximal solution
is indeed a global solution over $[t_0,\infty)$, i.e. it has no finite escape time.
In contrast, control strategies that merely confine the tracking error $\yM(t)-y_{\rf}(t)$ to $\cF_{\Funnel}$ lack 
this inherent boundedness guarantee, as illustrated by the following example.
\begin{example}
    The scalar  differential equation 
    \[
        \ddot{y}(t) = 2 \dot{y}^3+u(t),\qquad y(0)=0,\quad \dot{y}(0)=\tfrac{1}{2}
    \]
    of order two belongs to the model class $\cM^{1,2}_{0}$.
    If the constant control $u\equiv0$ is applied to the differential equation,
    then the initial value problem has the unique maximal solution
    $y:[0,1)\to\R$, $t\mapsto1-\sqrt{1-t}$.
    This constant control allows tracking of the constant reference $y_{\rf}\equiv1$
    within a funnel~$\cF_{\Funnel}$ given a constant funnel function $\Funnel \equiv 2$
    because we have, for all $t\in[0,1)$,
    \[
        \Norm{y(t)-y_{\rf}(t)}=\Norm{1-\sqrt{1-t}-1}=\sqrt{1-t}<2=\Funnel(t).
    \]
    However, the solution has finite escape time and cannot be extended to a global solution as the derivative
    $\dot{y}(t)=\tfrac{1}{2\sqrt{1-t}}$ is unbounded and has a pole at $t=1$.
\end{example}
To ensure that the tracking error $\eMTrack(t)=\yM(t)-y_{\rf}(t)$ evolves within the funnel~$\cF_{\Funnel}$ 
(defined by a function $\Funnel\in\cG$ as outlined in~\Cref{Sec:ControlObjective}), we construct
auxiliary funnel functions $\Psi\coloneqq (\Funnel_1,\ldots,\Funnel_r)\in\cG^r$ to simplify the tracking problem.
A fundamental prerequisite is that the initial error -- the mismatch between the
model's initial trajectory $\yM^0$ from~\eqref{eq:Model_r} and the reference
trajectory~$y_{\rf}$ -- lies within the funnel, i.e. satisfies
$\Norm{\yM^0(t_0)-y_{\rf}(t_0)}<\Funnel(t_0)$.
Under this condition,  there exists $\gamma\in(0,1)$ such that 
\begin{equation}\label{eq:DefParameterGamma}
    \Norm{\yM^0(t_0)-y_{\rf}(t_0)}\leq\gamma^{r}\Funnel(t_0).
\end{equation}
Furthermore, the funnel function $\Funnel\in\cG$ satisfies $\dot \Funnel(t) \ge -\FunDeriv\Funnel(t) + \FunDiam$ for all $t\geq0$, 
where $\FunDeriv,\FunDiam>0$ are constants with $\Funnel(t_0)\geq \tfrac{\FunDiam}{\FunDeriv}$, see also~\eqref{eq:DefSetOfFunnelFunctions} for the definition of the set $\cG$.
For completeness, we briefly prove this existence result.
\begin{lemma}
Let $\Funnel\in\cG$, then there exists $\FunDeriv,\FunDiam>0$ such that    
\begin{equation}\label{eq:DefinitionAlphaBeta}
    \Funnel(t_0)\geq \frac{\FunDiam}{\FunDeriv}\quad \text{ and}\quad \dot \Funnel(t) \ge -\FunDeriv\Funnel(t) + \FunDiam\quad \fa t\geq 0. 
\end{equation}
\end{lemma}
\begin{proof}
    We have $\inf_{s\geq 0}\dot{\Funnel}(s) \leq 0$ due to the boundedness of $\Funnel$.
    In the case $\inf_{s\geq 0}\dot{\Funnel}(s)=0$, set $\FunDeriv\coloneqq 1$ and $\FunDiam\coloneqq \inf_{s\geq 0}{\Funnel}(s)>0$.
    Then, $\Funnel(t_0)\geq \tfrac{\FunDiam}{\FunDeriv}$ and, for all $t\geq 0$, 
    \[
        -\FunDeriv\Funnel(t) + \FunDiam= -\Funnel(t) + \inf_{s\geq 0}{\Funnel}(s)\leq 0\leq\inf_{s\geq 0}\dot{\Funnel}(s)\leq \dot{\Funnel}(t).
    \]
    If $\inf_{s\geq 0}\dot{\Funnel}(s)<0$, then set $\FunDeriv\coloneqq \tfrac{-\inf_{s\geq 0}\dot{\Funnel}(s)}{\tfrac{1}{2}\inf_{s\geq 0}{\Funnel}(s)}>0$ and $\FunDiam\coloneqq -\inf_{s\geq 0}{\dot{\Funnel}}(s)>0$.
    Then,  
    \[
        -\FunDeriv\Funnel(t) + \FunDiam= 
        -\FunDeriv\Funnel(t) + \tfrac{\FunDeriv}{2}\inf_{s\geq 0}{\Funnel}(s)
        \leq -\tfrac{\FunDeriv}{2}\inf_{s\geq 0}{\Funnel}(s) 
        \leq\inf_{s\geq 0}\dot{\Funnel}(s)\leq \dot{\Funnel}(t)
    \]
    for all $t\geq 0$. Moreover, $\Funnel(t_0)\geq \tfrac{\FunDiam}{\FunDeriv}$.
    This completes the proof.
\end{proof}
Given $\gamma\in(0,1)$ as in~\eqref{eq:DefParameterGamma}, constants $\FunDeriv,\FunDiam>0$ satisfying~\eqref{eq:DefinitionAlphaBeta}, and the initial time  $t_0\geq0$,
recursively select parameters $k_1,\ldots,k_{r-1}$ such that 
\begin{equation}\label{eq:cond-k_i}
    \begin{aligned}
        k_1 &\ge \frac{2\Norm{(\yMd^0-\dot y_{\rf})(t_0)}}{\gamma^{r-1}(1-\gamma)\Funnel(t_0)} +\frac{2\left(\FunDeriv+\frac{1}{\gamma^{r-1}}\right)}{1-\gamma},\\
        k_i &\ge \frac{2\gamma\Norm{\dd{t} \eM_{i}(\OpChi(\yM^0-y_{\rf})(t_0))}}{(1-\gamma)\left(\Norm{\eM_{i}(\OpChi(\yM^0-y_{\rf})(t_0))}+\frac{\FunDiam}{\FunDeriv \gamma^{i-2}}\right)} + \frac{2(1+\FunDeriv)}{1-\gamma}
    \end{aligned}
\end{equation}
for all $i=2,\ldots,r-1$, where $\eM_i$ are defined as in~\eqref{eq:ErrorVar}.
Using the shorthand notation  ${\eM_i^0\coloneqq \eM_{i}(\OpChi(\yM^0-y_{\rf})(t_0))}$, 
define the vector of auxiliary funnel functions $\Psi$ as $\Funnel_1\coloneqq \Funnel$, and $\Funnel_2,\ldots,\Funnel_r$ as follows:
\begin{equation}\label{eq:psi_i}
    \Funnel_{i+1}(t)\coloneqq  \frac{1}{\gamma^{r-i}}\rbl\Norm{\dot \eM_i^0} + k_i \Norm{\eM_i^0\vphantom{\dot \eM_{i}^0}}\rbr 
    \me^{-\FunDeriv (t-t_0)} + \frac{\FunDiam}{\FunDeriv\gamma^{r-1}}
\end{equation}
for $t\ge 0$ and $i=1,\ldots,r-1$.
Critically, the parameters~$k_i$ and functions~$\Funnel_{i}$ 
do \emph{only} depend on~$\OpChi(\yM^0-y_{\rf})(t_0)$, i.e. the value of $\yM^0-y_{\rf}$ and its derivatives at the initial time~$t_0$,
rather than the entire trajectories $\yM^0$ and $y_{\rf}$.
By construction of~$\Funnel_i$ and by observation~\eqref{eq:ErrorVarDyn}, we have
\[
    \Norm{\eM_i^0\vphantom{\dot \eM_{i}^0}} \le \Norm{\dot \eM_{i-1}^0} + k_{i-1} \Norm{\eM_{i-1}^0\vphantom{\dot \eM_{i}^0}} < \Funnel_i(t_0)
\]
for all $i=2,\ldots,r$, and, by assumption~\eqref{eq:DefParameterGamma},
\[
    \Norm{\eM_1^0} \le \gamma^r \Funnel(t_0) < \Funnel_1(t_0).
\]
Therefore, $\OpChi(\yM^0-y_{\rf})(t_0)$ is an element of $\cD_{t_0}^{\Psi}$ as defined in \eqref{eq:DefSetD}.
Utilising the so constructed parameters~$k_i$ from~\eqref{eq:cond-k_i} and auxiliary funnel functions $\Funnel_i$ from~\eqref{eq:psi_i},
\Cref{Prop:OnlyLastFunnelVary} establishes the following:
If, at time $\hat{t}\geq t_0$, all auxiliary error variables lie within their respective funnels for 
a function $\zeta\in \cC^{r-1}([\hat t,\infty),\R^m)$, i.e. ${\OpChi(\zeta)(\hat t)\in\cD_{\hat t}^\Psi}$, 
and thereafter the last error variable $\eM_r(\OpChi(\zeta)(t))$ evolves within its funnel given by~$\Funnel_{r}$, 
then \emph{all} auxiliary error variables~$\eM_i(\OpChi(\zeta)(t))$ remain within their respective funnels given by~$\Funnel_{i}$ for all $t\geq \hat{t}$.
This has remarkable implications for the reference tracking problem from~\Cref{Sec:ControlObjective}.
If a control function~$u\in L_{\loc}^{\infty}([t_0,\infty),\R^m)$ is applied to the model~\eqref{eq:Model_r} and achieves that 
$\Norm{\eM_r(\xM(t)-\OpChi(y_{\rf})(t))}< \Funnel_{r}(t)$ for $t\geq\hat{t}$, then ${\xM(t)-\OpChi(y_{\rf})(t)\in\cD_{t}^\Psi}$
for $t\geq\hat{t}$, assuming initially  $\xM(\hat{t})-\OpChi(y_{\rf})(\hat{t})\in\cD_{\hat{t}}^\Psi$.
This implies, in particular, ${\Norm{\yM(t)-y_{\rf}(t)}<\Funnel(t)}$ for $t\geq\hat{t}$ 
because of the definitions of the set $\cD_{t}^\Psi$ in~\eqref{eq:DefSetD}, 
the error variables $\eM_{i}$ in~\eqref{eq:ErrorVar}, and the function $\Funnel_{1}=\Funnel$.
In summary, a control~$u$ ensuring that the last auxiliary error variable 
$\eM_r(\xM-\OpChi(y_{\rf}))$ evolves within the funnel~$\cF_{\Funnel_{r}}$
defined by $\Funnel_{r}$ solves the reference tracking problem in~\Cref{Sec:ControlObjective}.

\begin{prop}\label{Prop:OnlyLastFunnelVary}
For $\Funnel\in\cG$ and $t_0\geq0$,
let the parameters $k_i\geq 0$ be given for $i=1\ldots, r-1$ as in~\eqref{eq:cond-k_i} and 
$\Psi=(\Funnel_1,\ldots,\Funnel_r)\in\cG^r$ be given as in~\eqref{eq:psi_i}.
Further, let $\hat t\ge t_0$ and $\zeta\in \cC^{r-1}([\hat t,\infty),\R^m)$ be such that $\OpChi(\zeta)(\hat t)\in\cD_{\hat t}^\Psi$. If  
$\Norm{\eM_r(\OpChi(\zeta)(t))}< \Funnel_r(t)$ for all $t\in[\hat t,s)$ for some $s>\hat t$, then 
$\OpChi(\zeta)(t)\in\cD_{t}^\Psi$ for all $t\in[\hat t,s)$.
\end{prop}
\begin{proof}
    Seeking a contradiction, we assume that, for at least one ${i \in \{1,\ldots,r-1\}}$, there exists $t \in (\hat t,s)$
    such that $\Norm{\eM_i(\OpChi(\zeta)(t))} \geq \Funnel_i(t)$.
    W.l.o.g.\ let~$i$ be the largest index with this property.
    In the following, we use the shorthand notation $\eM_i(t)\coloneqq \eM_i(\OpChi(\zeta)(t))$ and
    ${\eM_i^0\coloneqq \eM_{i}(\OpChi(\yM^0-y_{\rf})(t_0))}$, as before in~\eqref{eq:psi_i}.
    However, we like to emphasise that $\eM_i(t_0) \neq \eM_i^0$ (if $\eM_i(\cdot)$ is defined at~$t_0$) in general,
    since $\OpChi(\zeta)(t_0)\neq \OpChi(y^0)(t_0)$ is possible.
    Invoking $\Norm{\eM_i(\hat t)}<\Funnel_i(\hat t)$ and the continuity of the involved functions, define 
    ${t^\star \coloneqq  \min \setdef{ t \in [\hat t,s) }{ \Norm{\eM_i(t) } = \Funnel_i(t)}}$.
    Set  $\eps\coloneqq\max\cbl\sqrt{\tfrac12 (1+\FunContr)},\Norm{\tfrac{\eM_i(\hat{t})}{\Funnel_{i}(\hat{t})}}\cbr\in(0,1)$.
    Due to continuity of the involved functions,
    there exists $t_\star \coloneqq  \max \setdef{ t \in [\hat t,t^\star) }{ \Norm{\tfrac{\eM_i(t)}{\Funnel_{i}(t)} } = \eps}$.
    We have $\eps\le\Norm{\tfrac{\eM_i(t)}{\Funnel_{i}(t)} } \le 1$ for all~$t\in[t_\star,t^\star]$.
    Utilising~\eqref{eq:ErrorVarDyn} and omitting the dependency on $t$, we calculate for~$t \in [t_{\star},t^\star]$:
    \begin{align*}
        \tfrac{1}{2}\dd{t}\Norm{\frac{\eM_i}{\Funnel_i}}^2
        &= \al\frac{\eM_i}{\Funnel_i},\frac{\dot{\eM}_i\Funnel_i-\eM_i\dot{\Funnel}_i}{\Funnel_i^2}\ar
        = \al\frac{\eM_i}{\Funnel_i},-\rbl k_i+\frac{\dot{\Funnel}_i}{\Funnel_i}\rbr\frac{\eM_{i}}{\Funnel_i}+\frac{\eM_{i+1}}{\Funnel_i}\ar\\
        &\le -\rbl k_i+\frac{\dot{\Funnel}_i}{\Funnel_i}\rbr\Norm{\frac{\eM_i}{\Funnel_i}}^2 + \Norm{\frac{\eM_i}{\Funnel_i}} \frac{\Norm{\eM_{i+1}}}{\Funnel_i}
        \leq -\rbl k_i+\frac{\dot{\Funnel}_i}{\Funnel_i}\rbr\eps^2 +\frac{\Funnel_{i+1}}{\Funnel_i},
    \end{align*}
    where we used $\Norm{\eM_{i+1}(t)}\leq \Funnel_{i+1}(t)$ due to the maximality of $i$.
    Now, we distinguish the two cases $i=1$ and $i>1$. For $i=1$, note that $\Funnel_1 = \Funnel$ and by properties of~$\cG$ it follows
    \[
        -\frac{\dot{\Funnel}(t)}{\Funnel(t)}\le \frac{\FunDeriv \Funnel(t) - \FunDiam}{\Funnel(t)} \le \FunDeriv.
    \]
    Furthermore, we have that $\Funnel(t) \ge \rbl\Funnel(t_0) -\frac{\FunDiam}{\FunDeriv}\rbr\me^{-\FunDeriv (t-t_0)} + \frac{\FunDiam}{\FunDeriv}$ for all $t\ge t_0$. Therefore,
    \begin{align*}
        \frac{\Funnel_2(t)}{\Funnel(t)}
        &\le  \frac{1}{\FunContr^{r-1}} \frac{\rbl\Norm{\dot \eM_1^0} + k_1 \Norm{\eM_1^0\vphantom{\dot \eM_{i}^0}}\rbr \me^{-\FunDeriv (t-t_0)}}{\rbl\Funnel(t_0) -\frac{\FunDiam}{\FunDeriv}\rbr \me^{-\FunDeriv (t-t_0)} + \frac{\FunDiam}{\FunDeriv}}
        + \frac{\FunDiam}{\FunDeriv\FunContr^{r-1}\rbl\rbl\Funnel(t_0) -\frac{\FunDiam}{\FunDeriv}\rbr \me^{-\FunDeriv (t-t_0)} + \frac{\FunDiam}{\FunDeriv}\rbr} \\
        &\le \frac{1}{\FunContr^{r-1}} \frac{\Norm{\dot \eM_1^0} + k_1 \Norm{\eM_1^0\vphantom{\dot \eM_{i}^0}}}{\Funnel(t_0)} + \frac{1}{\FunContr^{r-1}}
        \le \FunContr k_1 + \frac{\Norm{\dot \eM_1^0}}{\FunContr^{r-1} \Funnel(t_0)} + \frac{1}{\FunContr^{r-1}}
    \end{align*}
    for all $t\ge t_0$, where the estimate $\Norm{\eM_1^0}\le \FunContr^r \Funnel(t_0)$ was used, see~\eqref{eq:DefParameterGamma}. 
    Hence, we obtain 
    \begin{align*}
         \tfrac{1}{2}\dd{t}\Norm{\frac{\eM_1}{\Funnel}}^2 
         &\le -\tfrac12 (k_1 - \FunDeriv) (1+\FunContr) + \FunContr k_1 + \frac{\Norm{\dot \eM_1^0}}{\FunContr^{r-1} \Funnel(t_0)} + \frac{1}{\FunContr^{r-1}}\\
         &\le -\tfrac12 (1-\FunContr) k_1 + \FunDeriv + \frac{\Norm{\dot \eM_1^0}}{\FunContr^{r-1} \Funnel(t_0)} + \frac{1}{\FunContr^{r-1}} \le 0
    \end{align*}
    for all $t \in [t_{\star},t^\star]$, where the last inequality follows from~\eqref{eq:cond-k_i}. Now, consider the case~${i>1}$. Then, we have
    $-\tfrac{\dot{\Funnel}_i(t)}{\Funnel_i(t)}\le  \FunDeriv$ for all $t\ge t_0$ and, invoking that by~\eqref{eq:ErrorVarDyn}
    \[
        \Norm{\eM_i^0\vphantom{\dot \eM_{i}^0}} \le \Norm{\dot \eM_{i-1}^0} + k_{i-1} \Norm{\eM_{i-1}^0\vphantom{\dot \eM_{i}^0}},
    \]
    we find that
    \begin{align*}
        \frac{\Funnel_{i+1}(t)}{\Funnel_i(t)}&=  
        \frac{\frac{1}{\FunContr^{r-i}}\rbl\Norm{\dot \eM_i^0} + k_i \Norm{\eM_i^0\vphantom{\dot \eM_{i}^0}}\rbr \me^{-\FunDeriv (t-t_0)} + \frac{\FunDiam}{\FunDeriv\FunContr^{r-1}}}
        {\frac{1}{\FunContr^{r-i+1}}\rbl\Norm{\dot \eM_{i-1}^0} + k_{i-1} \Norm{\eM_{i-1}^0\vphantom{\dot \eM_{i}^0}}\rbr\me^{-\FunDeriv (t-t_0)} + \frac{\FunDiam}{\FunDeriv\FunContr^{r-1}}}\\
        &\le \FunContr \frac{\Norm{\dot \eM_i^0} + k_i \Norm{\eM_i^0\vphantom{\dot \eM_{i}^0}}}{\Norm{\dot \eM_{i-1}^0} + k_{i-1} \Norm{\eM_{i-1}^0\vphantom{\dot \eM_{i}^0}} + \frac{\FunDiam}{\FunDeriv \FunContr^{i-2}}} + 1
        \le \FunContr k_i + \FunContr \frac{\Norm{\dot \eM_i^0}}{\Norm{\eM_i^0}+ \frac{\FunDiam}{\FunDeriv \FunContr^{i-2}}} + 1
    \end{align*}
    for all $t\ge t_0$. Hence, we obtain that
    \begin{align*}
         \tfrac{1}{2}\dd{t}\Norm{\frac{\eM_i}{\Funnel_i}}^2 
         &\le -\tfrac12 (k_i - \FunDeriv) (1+\FunContr) + \FunContr k_i + \FunContr \frac{\Norm{\dot \eM_i^0}}{\Norm{\eM_i^0}+ \frac{\FunDiam}{\FunDeriv \FunContr^{i-2}}} + 1\\
         &\le -\tfrac12 (1-\FunContr) k_i + \FunDeriv + \FunContr \frac{\Norm{\dot \eM_i^0}}{\Norm{\eM_i^0}+ \frac{\FunDiam}{\FunDeriv \FunContr^{i-2}}} + 1 \le 0
    \end{align*}
    for all $t \in [t_{\star},t^\star]$, where the last inequality follows from~\eqref{eq:cond-k_i}. Summarising, in each case the contradiction 
    \[
        1 \leq \Norm{ \frac{\eM_i(t^*)}{\Funnel_i(t^*)}}^2 \leq \Norm{\frac{\eM_i(t_*)}{\Funnel_i(t_*)}}^2 = \eps^2<1
    \]
    arises, which completes the proof. 
\end{proof}

The proof of \Cref{Prop:OnlyLastFunnelVary} not only shows that all 
auxiliary error variables $\eM_i(\OpChi(\zeta)(t))$ stay within their respective funnels given by~$\Funnel_{i}$ for $i=1,\ldots, r-1$,
if the last error variable $\eM_r(\OpChi(\zeta)(t))$ evolves within its funnel given by~$\Funnel_{r}$, 
but it moreover shows that the auxiliary error variables~$\eM_i(\OpChi(\zeta)(t))$ 
always uphold an $\eps$-distance to the funnel boundaries~$\Funnel_{i}$. 
For systems with order $r>1$, this means that the tracking error $\yM(t)-y_{\rf}(t)$ fulfils 
\[
    \Norm{\yM(t)-y_{\rf}(t)}<\eps\Funnel(t)
\]
on every interval $[\hat{t},s]$ with $s>\hat{t}$, assuming $\Norm{\yM(\hat{t})-y_{\rf}(\hat{t})}<\eps\Funnel(\hat{t})$.
We sum up this observation in the following.

\begin{corollary}\label{Cor:ErrorEpsDistance}
    For $\Funnel\in\cG$ with $r>1$, 
    let the parameters $k_i\geq 0$ be given  as in~\eqref{eq:cond-k_i} for ${i=1\ldots, r-1}$ and 
    $\Psi=(\Funnel_1,\ldots,\Funnel_r)\in\cG^r$ be given as in~\eqref{eq:psi_i}.
    Further, let $s>\hat t\ge t_0$  and $\zeta\in \cC^{r-1}([\hat t,\infty),\R^m)$ be given such that $\OpChi(\zeta)(t)\in\cD_{t}^\Psi$
    for all $t\in[\hat{t},s]$.   
    There exists $\eps\in(0,1)$, independent of $t$, $s$, and $\zeta$, such that 
    if $\Norm{\eM_1(\OpChi(\zeta)(\hat{t}))}<\eps\Funnel_1(\hat{t})$, then
    \[
        \Norm{\eM_1(\OpChi(\zeta)(t))}<\eps\Funnel_1(t)
    \]
    for all $t\in[\hat{t},s]$.
\end{corollary}
\begin{proof}
    Setting $\eps\coloneqq \sqrt{\tfrac12 (1+\FunContr)}$, this is can be directly seen following the argument for $\eM_1$ of the proof of \Cref{Prop:OnlyLastFunnelVary}.
\end{proof}

The construction of the parameters~$k_i$ in~\eqref{eq:cond-k_i} 
and the auxiliary funnel functions~$\Funnel_i$ in~\eqref{eq:psi_i} assumed given and fixed  
initial trajectory $\yM^0$ for the model~\eqref{eq:Model_r}.
These parameters and functions were tailored to enable the analysis in~\Cref{Prop:OnlyLastFunnelVary},
requiring  the funnel functions to be large enough to accommodate the initial errors
values ${\eM_i^0\coloneqq \eM_{i}(\OpChi(\yM^0-y_{\rf})(t_0))}$ for $i=1,\ldots,r$.
However, this approach imposes intricate constraints on the parameters~$k_i$
and time-varying functions~$\Funnel_i$, complicating their selection.
Crucially, the initial trajectory $\yM^0$ -- a modelling parameter for the model~\eqref{eq:Model_r} 
of the original system~\eqref{eq:Sys} -- often admits flexibility.
By strategically choosing the function $\yM^0$, we simplify the design of
$k_i$ and $\Funnel_i$.
We therefore present a simplified parameter design in the following.
As before, we assume $\Funnel\in\cG$ to be given with associated constants
$\FunDeriv,\FunDiam>0$ fulfilling~\eqref{eq:DefinitionAlphaBeta}.
Define $\Funnel_1\coloneqq \Funnel$ and  
\begin{equation}\label{eq:cond-k_i_psi_i_simple}
    \begin{array}{rcccl}
        k_1=        & \ldots & = & k_{r-1}      & \geq  \FunDeriv+2,\\
        \Funnel_{2}(t)=& \ldots & = & \Funnel_{r}(t)  & \coloneqq \frac{\FunDiam}{\FunDeriv}.
    \end{array}
\end{equation}
This yields the simplified constraints:
\[
    \Norm{\eM_{i}(\OpChi(\yM^0-y_{\rf})(t_0))}<\Funnel_i(t_0),\quad i=1,\ldots,r
\]
on the initial parameter $\yM^0$.
The following \Cref{Lemma:NewPsiFunktions} adapts \Cref{Prop:OnlyLastFunnelVary} to this streamlined framework.
\begin{prop}\label{Lemma:NewPsiFunktions}
    For $\Funnel\in\cG$, let the parameters $k_i\geq \FunDeriv+2$ be given for $i=1,\ldots, r-1$ and 
    $\Psi=(\Funnel_1,\ldots,\Funnel_r)\in\cG^r$ be given as in~\eqref{eq:cond-k_i_psi_i_simple}.
    Let $\hat t\ge t_0$ and $\zeta\in \cC^{r-1}([\hat t,\infty),\R^m)$ such that $\OpChi(\zeta)(\hat t)\in\cD_{\hat t}^{\Psi}$. 
    If $\Norm{\xi_r(\OpChi(\zeta)(t))}< \tfrac{\FunDiam}{\FunDeriv}$ for all $t\in[\hat t,s)$, $s>\hat t$, then 
    $\OpChi(\zeta)(t)\in\cD_{t}^\Psi$ for all $t\in[\hat t,s)$.
\end{prop}
\begin{proof}
We modify the proof of~\Cref{Prop:OnlyLastFunnelVary} to the changed setting.
Seeking a contradiction, we assume that there exists $t \in (\hat t, s)$ such that 
$ \| \eM_i(\OpChi(\zeta)(t)) \| \geq \Funnel_i(t)$ for at least one ${i \in \{1,\ldots,r-1\}}$.
W.l.o.g.\ let~$i$ be the largest index with this property.
We use the shorthand notation $\eM_i(t)\coloneqq \eM_i(\OpChi(\zeta)(t))$.    
Define  $\eps\coloneqq\max\cbl\sqrt{\tfrac12},\Norm{\tfrac{\eM_i(\hat{t})}{\Funnel_{i}(\hat{t})}}\cbr\in(0,1)$.
Invoking continuity of~$\eM_i$, there exist time instants
${t^\star \coloneqq \min \setdef{ t \in [\hat t, s] }{ \| \eM_i(t) \| = \Funnel_i(t) }}$
and 
$t_\star \coloneqq  \max \setdef{ t \in [\hat t, t^\star) }{ \fa s\in[t,t^\star]: \| \eM_i(s) \|= \eps  \Funnel_i(s) }$.
We separately consider the two cases $i=1$ and $i>1$.
First, we suppose $i=1$.
Then, note that $\Funnel(t)\geq \Norm{\eM_1(t)}\geq \eps\Funnel(t)$ for all $t\in[t_\star,t^\star]$.
By properties of~$\Funnel \in \cG$, we have
\[
    -\frac{\dot{\Funnel}(t)}{\Funnel(t)}\le \frac{\FunDeriv\Funnel(t) - \FunDiam}{\Funnel(t)} \le \FunDeriv,
\]
and
$\Funnel(t) \ge \rbl\Funnel(t_0)-\frac{\FunDiam}{\FunDeriv}\rbr \me^{-\FunDeriv(t-t_0)} + \frac{\FunDiam}{\FunDeriv}\geq\frac{\FunDiam}{\FunDeriv}$ for all $t\ge t_0$.
Omitting the dependency on $t$, we calculate for $t \in [t_{\star},t^\star]$
\begin{align*}
    &\tfrac{1}{2}\dd{t}\Norm{\tfrac{\eM_1}{\Funnel}}^2
    = \al\tfrac{\eM_1}{\Funnel},\tfrac{\dot{\eM}_1\Funnel-\eM_1\dot{\Funnel}}{\Funnel^2}\ar
    = \al\tfrac{\eM_1}{\Funnel},- \rbl k_1+\tfrac{\dot{\Funnel}}{\Funnel}\rbr  \tfrac{\eM_{1}}{\Funnel}+\tfrac{\eM_{2}}{\Funnel}\ar\\
    &\le  -\rbl k_1+\tfrac{\dot{\Funnel}}{\Funnel}\rbr\Norm{\tfrac{\eM_1}{\Funnel}}^2  + \tfrac{ \| \eM_1 \| \|\eM_2\|}{\Funnel^2} 
    \!\leq\! -\rbl k_1+\tfrac{\dot{\Funnel}}{\Funnel}\rbr \tfrac{1}{2} + \tfrac{\FunDeriv}{\FunDiam}\|\eM_2\| \\
    &\leq\! -\rbl k_1-\FunDeriv\rbr \tfrac{1}{2} +1\leq 0,
\end{align*}
where we used $k_1\geq\FunDeriv+2$ and $\Norm{\eM_{2}(t)}\leq \frac{\FunDiam}{\FunDeriv}$ for all~$t \in [t_\star,t^\star]$.
Thus, upon integration, the contradiction
\begin{equation*}
    1 = \Norm{\tfrac{\eM_1(t^\star)}{\Funnel(t^\star)} }^2 \le \Norm{\tfrac{\eM_1(t_\star)}{\Funnel(t_\star)}}^2 =\eps^2<1
\end{equation*}
arises.
Now, we consider the case $\| \eM_i(t) \|\ge \Funnel_i(t)$ for $i>1$.
By the choice of $\eps$, we have $\FunDiam/\FunDeriv=\Funnel_i(t)\geq \Norm{\eM_i(t)}\geq \FunDiam/(\FunDeriv\sqrt{2})$ for all $t\in[t_\star,t^\star]$.
Thus,  we calculate
\begin{align*}
    \tfrac{1}{2 }\dd{t} \| \eM_i \|^2
    =  \langle \eM_i , \dot{\eM}_i \rangle
    =  \langle \eM_i, -k_i \eM_i + \eM_{i+1} \rangle
     \le - k_i \| \eM_i \|^2 + \| \eM_i \| \| \eM_{i+1} \|
    \leq \tfrac{\FunDiam^2}{ \FunDeriv^2} \left( - \tfrac{k_i}{2} + 1 \right)  \le  0
\end{align*}
for $t \in [t_{\star},t^\star]$, where we used~$k_i\geq 2$ and that 
$\| \eM_{i+1}(t)\| \leq \tfrac{\FunDiam}{\FunDeriv}$ for all~$t \in [t_\star,t^\star]$  by maximality of $i$.
Hence, the contradiction
\[
  \tfrac{\FunDiam^2}{\FunDeriv^2} \leq \| \eM_i(t^*)\|^2 \leq \| \eM_i(t_*)\|^2< \tfrac{\FunDiam^2}{\FunDeriv^2}
\]
arises, which completes the proof.
\end{proof}

As in~\Cref{Prop:OnlyLastFunnelVary}, the proof of \Cref{Lemma:NewPsiFunktions} shows that 
the auxiliary error variables~$\eM_i(\OpChi(\zeta)(t))$ 
always uphold an $\eps$-distance to the funnel boundaries~$\Funnel_{i}$. 
Similarly to~\Cref{Cor:ErrorEpsDistance}, we get the following result for systems with order $r>1$.

\begin{corollary}\label{Cor:ErrorEpsDistanceAlt}
    For $\Funnel\in\cG$ with $r>1$, let the parameters $k_i\geq \FunDeriv+2$ be given for ${i=1\ldots, r-1}$ and 
    $\Psi=(\Funnel_1,\ldots,\Funnel_r)\in\cG^r$ be given as in~\eqref{eq:cond-k_i_psi_i_simple}.
    Further, let $s>\hat t\ge t_0$  and $\zeta\in \cC^{r-1}([\hat t,\infty),\R^m)$ be given such that $\OpChi(\zeta)(t)\in\cD_{t}^\Psi$
    for all $t\in[\hat{t},s]$.   
    There exists $\eps\in(0,1)$, independent of $t$, $s$, and $\zeta$, such that 
    if $\Norm{\eM_1(\OpChi(\zeta)(\hat{t}))}<\eps\Funnel_1(\hat{t})$, then
    \[
        \Norm{\eM_1(\OpChi(\zeta)(t))}<\eps\Funnel_1(t)
    \]
    for all $t\in[\hat{t},s]$.
\end{corollary}
\begin{proof}
    Setting $\eps\coloneqq \sqrt{\tfrac12}$, this is can be directly seen following the argument for $\eM_1$ in the proof of \Cref{Lemma:NewPsiFunktions}.
\end{proof}

\begin{remark}
    Building on \Cref{Lemma:NewPsiFunktions}, a low-complexity funnel controller
    for systems of higher-order was proposed in~\cite{Dennst25}. Contrary to
    prior works, this control approach eliminates the use of time-varying
    reciprocal penalty terms, replacing them with constant gains. 
    The simpler controller design has the potential to
    mitigate numerical issues and enhance its practicality for real-world
    applications. 
\end{remark}

The two presented parameters designs show that there exists a delicate interplay between 
the choice of initial trajectory $\yM^0$, the parameters $k_i\geq0$, 
the associated error variables~$\eM_i$ from~\eqref{eq:ErrorVar},
and the corresponding auxiliary funnel functions~$\Funnel_i$.
Even though other parameter designs are conceivable, 
in the remaining part of this presented thesis,  the
error variables~$\eM_i$ are always defined as in~\eqref{eq:ErrorVar},
the vector $\Psi\coloneqq  (\Funnel_1,\ldots,\Funnel_r)\in\cG^r$ of funnel functions 
and the corresponding parameters $k_i$ for $i=1,\ldots, r$ are 
always chosen either according to~\eqref{eq:cond-k_i} and~\eqref{eq:psi_i} 
or according to~\eqref{eq:cond-k_i_psi_i_simple}.
We will use the abbreviated notation 
\begin{equation}\label{eq:DefFunnelBoundaryFunctions}
    \Psi\in\FunnelBoundaryFuncs
\end{equation}
to refer to one of the presented cases for the parameter design.
\begin{assn}
    We will implicitly always assume that the 
    initial auxiliary errors are within their respective funnels,
    i.e. $\Norm{\eM_{i}(\OpChi(\yM^0-y_{\rf})(t_0))}<\Funnel_i(t_0)$ for $i=1,\ldots,r$.
\end{assn}
\begin{remark}\label{Rem:InitTrajectoryInFunnel}
    Since the control problem is formulated merely for $t\geq t_0$, it is
    possible that $\Norm{\eM_{i}(\OpChi(\yM^0-y_{\rf})(t))}\geq \Funnel_i(t)$ for some
    $t\in[0,t_0)$ and some $i=1,\ldots,r$. To avoid treating this interval as a
    special case, we assume without loss of generality that
    ${\OpChi(\yM^0-y_{\rf})(t)\in\cD_{t}^{\Psi}}$ for all $t\in[0,t_0]$ in the
    remaining part of this presented thesis.
    It is clear that this is no restriction on the construction of the parameters $k_i$ and auxiliary funnel 
    functions $\Funnel_i$ in~\eqref{eq:cond-k_i_psi_i_simple} since the initial trajectory $\yM^0$ is chosen in order to fit to these parameters.
    In the parameter design setting~\eqref{eq:cond-k_i} and~\eqref{eq:psi_i}, the initial error values~$\eM_i^0$ and~$\dot{\eM}_i^0$ can be
    replaced by their respective suprema on the compact interval~$[0,t_0]$ due to the continuity of the involved functions.
\end{remark}
The following~\Cref{Prop:OnlyLastFunnel} summarises the main observation made about the parameter construction, 
namely that the reference tracking problem from~\Cref{Sec:ControlObjective} can be solved by
a control~$u$ that achieves that last auxiliary error variable 
$\eM_r(\xM-\OpChi(y_{\rf}))$ evolves within the funnel~$\cF_{\Funnel_{r}}$
defined by $\Funnel_{r}$.
\begin{prop}\label{Prop:OnlyLastFunnel}
Let $\Psi\in\FunnelBoundaryFuncs$,
$\hat t\ge t_0$ and $\zeta\in \cC^{r-1}([\hat t,\infty),\R^m)$ be given such that $\OpChi(\zeta)(\hat t)\in\cD_{\hat t}^\Psi$. If  
$\Norm{\eM_r(\OpChi(\zeta)(t))}< \Funnel_r(t)$ for all $t\in[\hat t,s)$ and for some $s>\hat t$, then 
${\OpChi(\zeta)(t)\in\cD_{t}^\Psi}$ for all $t\in[\hat t,s)$.
\end{prop}
\begin{proof}
    This is an immediate consequence of~\Cref{Prop:OnlyLastFunnelVary} and~\Cref{Lemma:NewPsiFunktions}.
\end{proof}

\begin{remark}
    Note that the results of \Cref{Prop:OnlyLastFunnelVary}, \Cref{Lemma:NewPsiFunktions}, and \Cref{Prop:OnlyLastFunnel} also hold true 
    if one allows for $\Norm{\eM_r(\OpChi(\zeta)(t))}\leq \Funnel_r(t)$.
    To be more precise:
    Let $\Psi\in\FunnelBoundaryFuncs$,
    $\hat t\ge t_0$ and $\zeta\in \cC^{r-1}([\hat t,\infty),\R^m)$ be such that $\Norm{\eM_i(\OpChi(\zeta)(t))}< \Funnel_i(t)$ for 
    all $i=1,\ldots, r-1$.
    If $\Norm{\eM_r(\OpChi(\zeta)(t))}\leq \Funnel_r(t)$ for all $t\in[\hat t,s)$ for some $s>\hat t$, then 
    $\Norm{\eM_i(\OpChi(\zeta)(t))}< \Funnel_i(t)$ for $t\in[\hat t,s)$ and all $i=1,\ldots, r-1$.
\end{remark}

We will sum up the observations made in~\Cref{Cor:ErrorEpsDistance} and~\Cref{Cor:ErrorEpsDistanceAlt} in the following.

\begin{corollary}\label{Cor:FMPCExistenceLambda}
    Let $r>1$ and $\Psi=(\Funnel_1,\ldots,\Funnel_r)\in\FunnelBoundaryFuncs$.
    There exists $\eps\in(0,1)$ with the following property.
    If, for $s>\hat t\ge t_0$,  a function $\zeta\in \cC^{r-1}([\hat t,\infty),\R^m)$ 
    fulfils $\OpChi(\zeta)(t)\in\cD_{t}^\Psi$ for all $t\in[\hat{t},s]$   
    and $\Norm{\eM_1(\OpChi(\zeta)(\hat{t}))}<\eps\Funnel_1(\hat{t})$, then
    \[
        \Norm{\eM_1(\OpChi(\zeta)(t))}<\eps\Funnel_1(t)
    \]
    for all $t\in[\hat{t},s]$.
\end{corollary}

\subsection{Feasible control signals}\label{Sec:InitialRecFeasibilty}
Prior to formulating the optimal control problem (OCP) with funnel stage cost 
functions (to be solved in the MPC~\Cref{Algo:MPC}), we first address 
the issue of ensuring \emph{initial and recursive feasibility}.
In the previous~\Cref{Sec:HighRelativeDegree}, we derived sufficient
conditions for a control $u\in L_{\loc}^{\infty}([\hat{t},\infty),\R^m)$ to solve the
reference tracking problem outlined in~\Cref{Sec:ControlObjective}
when applied to the model~\eqref{eq:Model_r} at time
$\hat{t}\geq t_0$. However, the existence of such a control function 
-- an essential prerequisite for the successful application of the MPC~\Cref{Algo:MPC}~--  
remains unverified. We now establish sufficient conditions to guarantee this existence.

Suppose the MPC~\Cref{Algo:MPC} solves the reference tracking problem up to time
$\hat{t}\in[t_0,\infty]$ ensuring $\xM(t)-\OpChi(y_{\rf})(t)$ is an element of the
set~$\cD^\Psi_t$ (from~\eqref{eq:DefSetD}) for all $t\leq\hat{t}$, where $\xM$ solves  
the model differential equation~\eqref{eq:Model_r}. 
The concatenated solution trajectory then belongs to the set
\begin{equation}\label{eq:Def:FunnelTrajectories}
    \FunnelTrajectories_{\hat{t}}\coloneqq \setdef
        {\zeta\in \cR(\Rp,\R^{rm})}
        {
             \fa t\in[0,\hat{t}]: \zeta(t)-\OpChi(y_{\rf})(t)\in\cD_t^{\Psi}
        }.
\end{equation}
The concatenated solution trajectory is a solution of the differential
equation~\eqref{eq:Model_r} in the sense of~\Cref{Def:ModSolution}
on the intervals of the form~$[t_k,t_{k+1}]$ for $k\in\N_0$ and differentiable on these interval.
However, it is, in general, in its entirety not a continuous function but
merely a regulated function because of the potential non-continuous
re-initialisation of the model in Step~\ref{agostep:MPCFirst}
of~\Cref{Algo:MPC}. Note that, if the concatenated solution trajectory is only defined
on a finite interval, then  we implicitly assume a right extension of the solution
as in~\Cref{Rem:PropertiesOperator}~\ref{Item:Rem:PropertiesOperator:Causality} 
in the definition of $\FunnelTrajectories_{\hat{t}}$ in~\eqref{eq:Def:FunnelTrajectories}.
During the execution of the MPC~\Cref{Algo:MPC}, the initial
values $\xMh$ and $\oTMh$  for the model~\eqref{eq:Model_r} at time~$\hat{t}$ as
in~\Cref{Def:ModSolution} are determined by the hitherto existing concatenated solution trajectory.
Building on these considerations, we define in the following the set of \emph{feasible initial values} for the model.

\begin{definition}[Feasible initial values $\InitValues(\hat{t})$]\label{Def:SetInitialValues}
Let $y_{\rf}\in W^{r,\infty}(\Rp,\R^{m})$, $\Psi\in\FunnelBoundaryFuncs$, and $\tau\geq0$.
Using the notation $I_{t_0}^{\hat{t},\tau}\coloneqq [\hat{t}-\tau,\hat{t}]\cap[t_0,\hat{t}] $, we
define the set of \emph{feasible initial values} for the
model~\eqref{eq:Model_r} at time $\hat{t}\geq t_0$ as
\begin{equation}\label{eq:Def:InitValues}
    \InitValues(\hat{t}) \coloneqq  
    \setdef
        {
           (\xMh,\oTMh) \in
           \cR(I_0^{\hat{t},\tau},\R^{rm})\times L^\infty_{\loc}( I_{t_0}^{\hat{t},\tau},\R^q)
        }
        {
             \ex \zeta\in\FunnelTrajectories_{\hat{t}}:
             \begin{array}{rr}
              \zeta|_{I_{0}^{\hat{t},\tau}}=\xMh,\\
              \oTM(\zeta)|_{I_{t_0}^{\hat{t},\tau}}=\oTMh
             \end{array}
        }.
\end{equation}
\end{definition}

\begin{remark}\label{Rem:InitialValuesNotEmpty}
    Note that $\InitValues(\hat{t})$ is, for all $\hat{t}\geq t_0$ and $\tau\geq0$,
    never empty  since $y_{\rf}\in \FunnelTrajectories_{\hat{t}}$ and, therefore, the
    pair $(\OpChi(y_{\rf})|_{I_{0}^{\hat{t},\tau}},\oTM(\OpChi(y_{\rf}))|_{I_{t_0}^{\hat{t},\tau}})$
    is an element of $\InitValues(\hat{t})$. 
    Moreover, $(\OpChi(\yM^0)|_{I_{0}^{\hat{t},\tau}},\oTM(\OpChi(\yM^0))|_{I_{t_0}^{\hat{t},\tau}})\in\InitValues(\hat{t})$
    for $\hat{t}=t_0$ according to our assumption in~\Cref{Rem:InitTrajectoryInFunnel}.
    In the remainder of this thesis, we use the notation ${\InitState\coloneqq(\xMh,\oTMh)\in\InitValues(\hat{t})}$
    to refer to  the initial values $\xMh$ and $\oTMh$ as a pair, since we will
    only rarely consider them independently of each other.
\end{remark}

\begin{remark}\label{Rem:InitialValueInFunnel}
    We want to highlight that choosing a feasible initial value $\InitState\in\InitValues(\hat{t})$ at time $\hat{t}$ 
    for the model~\eqref{eq:Model_r} implies 
    \[
        \xM(\hat{t}) - \OpChi(y_{\rf})(\hat{t})\in\cD_{\hat{t}}^{\Psi}.
    \]
    This means that the tracking error~$\eMTrack$ and all auxiliary error variables~$\xi_i$ as in~\eqref{eq:ErrorVar}
    are within their respective funnels at time $\hat{t}$.
\end{remark}

Although the general assumption within this chapter is that the model-plant mismatch $\eS$ as in~\eqref{eq:Intro:ModelPlantMismatch}
is always identical to zero, we already want to lay the fundamentals to allow for the initialisation of the model based on measurement data
from the system~\eqref{eq:Sys}, as in Step~\ref{agostep:MPCFirst} of~\Cref{Algo:MPC}.
Therefore, we define an \emph{initialisation strategy} as a function selecting a feasible initial value based on 
measurements $\hat{x}$ at time $\hat{t}\geq t_0$.
In application, one will replace $\hat{x}$ with $\OpChi(y)(\hat{t})$ where $y$ is the output of the system~\eqref{eq:Sys}.
\begin{definition}[Initialisation strategy]\label{Def:InitialisationStrategy}
    Let $y_{\rf}\in W^{r,\infty}(\Rp,\R^{m})$, $\Psi\in\FunnelBoundaryFuncs$, and $\tau\geq0$.
    Using the notation $I_{t_0}^{\hat{t},\tau}\coloneqq [\hat{t}-\tau,\hat{t}]\cap[t_0,\hat{t}]$, 
    we call a function 
    \[
        \InitStrategy:\bigcup_{\hat{t}\geq t_0}\cR(I_0^{\hat{t},\tau},\R^{rm})\to
        \bigcup_{\hat{t}\geq t_0}\cR(I_0^{\hat{t},\tau},\R^{rm})\times L^\infty_{\loc}( I_{t_0}^{\hat{t},\tau},\R^q)
    \]
    with $\InitStrategy(\hat{x})\in\InitValues(\hat{t}) $ for $\hat{x} \in
    \cR(I_0^{\hat{t},\tau},\R^{rm})$ 
    and $\hat{t}\geq t_0$ a \emph{$\tau$-initialisation strategy} for the model~$\eqref{eq:Model_r}$.
\end{definition}
In view of the limited memory property~\ref{Item:OperatorPropLimitMemory} of operator $\oTM$, 
we always utilise a $\tau$-initialisation strategy with $\tau$ being greater than or equal to the memory limit of $\oTM$.
In~\Cref{Def:InitialisationStrategy}, the domain of the initialisation strategy $\InitStrategy$ 
was chosen to be consistent with its codomain. However, this choice is entirely arbitrary. 
If deemed beneficial in a given setting, the domain can be adapted such that $\InitStrategy$
acts on signals defined on different time intervals of a different length, i.e. 
$\InitStrategy$ can be defined on $\bigcup_{\hat{t}\geq t_0}\cR(I_0^{\hat{t},\tilde{\tau}},\R^{rm})$ with $\tilde{\tau}\neq\tau$.
Changing the domain of $\InitStrategy$ in such a way does not change the validity of the results presented.

Let a feasible initial value $\InitState\in\InitValues(\hat{t})$ be given 
for the model~\eqref{eq:Model_r} at time $\hat{t}$. 
If a control ${u\in L^{\infty}([\hat{t},\hat{t}+T],\R^m)}$,
bounded by some constant $\umax\geq 0$, ensures that $\xM(t) - \OpChi(y_{\rf})(t)$ evolves within 
$\cD_{t}^{\Psi}$ for all $t$ over the next time interval of length $T>0$, then it is an element of the set
\begin{align}\label{eq:Def-U}
        \Controls(\umax,\InitState) \coloneqq 
        \setdef{
        u\in L^{\infty}([\hat{t},\hat{t}+T],\R^m)
        }
        {
        \begin{array}{l}
            \xM(t;\hat{t},\InitState,u) - \OpChi(y_{\rf})(t)\in\cD_{t}^{\Psi}\\
            \text{for all }t\in [\hat{t},\hat{t}+T],\SNorm{u}\le \umax 
        \end{array} 
        }.
\end{align}
We want to point out that, in the definition of the set $\Controls(\umax,\InitState)$ in~\eqref{eq:Def-U},
we implicitly assume that the solution~$\xM(t;\hat{t},\InitState,u)$ exists on the whole interval $[\hat{t},\hat{t}+T]$.
Due to the construction of funnel functions $\Psi\in\FunnelBoundaryFuncs$ in the previous~\Cref{Sec:HighRelativeDegree},
all functions $u\in\Controls(\umax,\InitState)$ solve the outlined tracking problem from~\Cref{Sec:ControlObjective}.
To guarantee the initial and recursive feasibility of the MPC~\Cref{Algo:MPC},
we therefore want to ensure that the set
$\Controls(\umax,\InitState)$ is always non-empty and that, if a
control $u\in\Controls(\umax,\InitState)$ is applied to the
model~\eqref{eq:Model_r}, then the model's state after applying this control can
again be used as a feasible initial value for the model.
We first address the latter question.
\begin{theorem}\label{Thm:RecursiveInitialValues}
    Consider model~\eqref{eq:Model_r} with $(\fM,\gM,\oTM)\in\cM^{m,r}_{t_0}$.
    Let $\tau\geq 0$ be greater than or equal to the memory limit of $\oTM$,
    $y_{\rf}\in W^{r,\infty}(\Rp,\R^{m})$, $\Psi\in\FunnelBoundaryFuncs$, 
    and $\InitState\in \InitValues(\hat{t})$ for $\hat{t}\geq t_0$.
    Further, let $\umax\geq 0$ and $T> 0$ such that $\Controls(\umax,\InitState)\neq\emptyset$.
    If a control ${u\in\Controls(\umax,\InitState)}$ is applied to the model~\eqref{eq:Model_r}, then
    there exists a solution of the initial value problem~\eqref{eq:Model_r} in the sense of~\Cref{Def:ModSolution}
    on the interval~$[0,\hat{t}+T]$ and every solution~$\xM:[0,\hat{t}+T]\to\R^m$ 
    fulfils 
    \[
        \fa \delta\!\in\![0,T]\!: (\xM(\cdot;\hat{t},\InitState,u)|_{[\hat{t}+\delta-\tau,\hat{t}+\delta]\cap[0,\hat{t}+\delta]},
        \oT(\xM(\cdot;\hat{t},\InitState,u))|_{[\hat{t}+\delta-\tau,\hat{t}+\delta]\cap[t_0,\hat{t}+\delta]})
        \!\in\!\InitValues(\hat{t}+\delta).
    \]
\end{theorem}
\begin{proof}
    Let $(\xMh,\oTMh)=\InitState\in\InitValues(\hat{t})$ be arbitrary but fixed.
    By~\Cref{Def:SetInitialValues} there exists 
    a function $\zeta\in\FunnelTrajectories_{\hat{t}}$ such that
    $\zeta|_{[\hat{t}-\tau,\hat{t}]\cap[0,\hat{t}]}=\xMh$
    and $\oTM(\zeta)|_{[\hat{t}-\tau,\hat{t}]\cap[t_0,\hat{t}]}=\oTMh$.
    If a control $u\in\Controls(\umax,\InitState)$
    is applied to the model~\eqref{eq:Model_r}, then there
    exists a solution of the initial value problem in the sense of~\Cref{Def:ModSolution} on the interval $[\hat{t},\hat{t}+T]$
    due to the definition of the set $\Controls(\umax,\InitState)$ as in~\eqref{eq:Def-U}.
    Let $\xM:[0,\hat{t}+T]\to\R^m$ be a solution of the initial value problem~\eqref{eq:ModDiff}.
    Since $\xM$ fulfils the initial conditions~\eqref{eq:Model_InitialValue} with $(\xMh,\oTMh)$,
    we have
    $\zeta|_{[\hat{t}-\tau,\hat{t}]\cap[0,\hat{t}]}=\xMh=\xM|_{[\hat{t}-\tau,\hat{t}]\cap[0,\hat{t}]}$
    and 
    $ 
        \oTM(\zeta)|_{[\hat{t}-\tau,\hat{t}]\cap[t_0,\hat{t}]}
        =\oTMh
        = \oTM(\xM)|_{[\hat{t}-\tau,\hat{t}]\cap[t_0,\hat{t}]}
    $.
    Define the function $\tilde{\zeta}\in\cR(\Rp,\R^{rm})$ by
    \begin{equation}\label{eq:RecursiveInitialValuesExtension}
        \tilde{\zeta}(t)=\begin{cases}
            \xM(t),& t\in [\hat{t},\hat{t}+T]\\
            \zeta(t),& t\in \Rp\backslash[\hat{t},\hat{t}+T],
        \end{cases}
    \end{equation}
    which fulfils 
    $\tilde{\zeta}|_{[\hat{t}-\tau,\hat{t}+\delta]\cap[0,\hat{t}+\delta]}
    =\xM|_{[\hat{t}-\tau,\hat{t}+\delta]\cap[0,\hat{t}+\delta]}$ 
    for all $\delta \in [0,T]$ and additionally
    $ 
        {\oTM(\tilde{\zeta})|_{[\hat{t}-\tau,\hat{t}]\cap[t_0,\hat{t}]} = \oTM(\xM)|_{[\hat{t}-\tau,\hat{t}]\cap[t_0,\hat{t}]}}
    $.
    Since $\tau\geq 0$ is greater than or equal to the memory limit of operator $\oTM$, 
    see property~\ref{Item:OperatorPropLimitMemory} in~\Cref{Def:OperatorClass},
    it follows that
    \[
        \oTM(\tilde{\zeta})|_{[\hat{t}+\delta-\tau,\hat{t}+\delta]\cap[t_0,\hat{t}+\delta]}
        =
        \oTM(\xM)|_{[\hat{t}+\delta-\tau,\hat{t}+\delta]\cap[t_0,\hat{t}+\delta]}
    \]
    for all $\delta \in [0,T]$.
    By choice of $u$, we have $\xM(t)-\OpChi(y_{\rf})(t)\in\cD_{t}^\Psi$ for all
    $t\in[\hat{t},\hat{t}+T]$.
    Hence, $\tilde{\zeta}\in\FunnelTrajectories_{\hat{t}+T}$. 
    We therefore have
    \[
    (\xM|_{[\hat{t}+\delta-\tau,\hat{t}+\delta]\cap[0,\hat{t}+\delta]},
    \oT(\xM)|_{[\hat{t}+\delta-\tau,\hat{t}+\delta]\cap[t_0,\hat{t}+\delta]})
    \in\InitValues(\hat{t}+\delta)
    \]
    for all $\delta \in [0,T]$.
    This completes the proof.
\end{proof}

\begin{remark}
    Although~\Cref{Def:SetInitialValues} was formulated for arbitrary $\tau\geq 0$, we utilised that $\tau$ was 
    greater than or equal to the memory limit of $\oTM$ in~\Cref{Thm:RecursiveInitialValues} 
    in order to ensure that the image of the operator $\oTM(\xM)$ is independent of the chosen left extension of $\xM$.
    The function~$\tilde{\zeta}$ in~\eqref{eq:RecursiveInitialValuesExtension} is, in general, only one of many possible extensions of $\xM$.
    \null
\end{remark}

\Cref{Thm:RecursiveInitialValues} shows that the model's state,
i.e. the state $\xM$ in combination with value of the operator $\oTM(\xM)$,
is at any point during application of a
control function ${u\in\Controls(\umax,\InitState)}$ a feasible initial value
(assuming $\InitState$ is a feasible initial value to begin with). 
This will be essential for the re-initialising of the model during the application of the MPC~\Cref{Algo:MPC}.
We now show that if $\umax\geq0$ is chosen large enough, then   
$\Controls(\umax,\InitState)$ is non-empty.
We first prove \Cref{Lemma:DynamicBounded} showing that the functions 
$(\fM,\gM,\oTM)\in\cM^{m,r}_{t_0}$ are bounded for functions $\zeta$ evolving within $\cD_{t}^\Psi$.
As a consequence, the dynamics of the model~\eqref{eq:Model_r} are bounded if a control is applied  
that ensures ${\xM(t) - \OpChi(y_{\rf})(t)}$ evolves within  $\cD_{t}^{\Psi}$.
\begin{lemma}\label{Lemma:DynamicBounded}
    Consider the model~$\eqref{eq:Model_r}$ with 
    $(\fM,\gM,\oTM)\in\cM^{m,r}_{t_0}$.
    Further, let ${y_{\rf}\in W^{r,\infty}(\Rp,\R^{m})}$ and $\Psi\in\FunnelBoundaryFuncs$.
    Then, there exist constants~$\fMmax$, $\gMmax,\gMInvmax\geq0$ such that for all $\hat{t}\in[t_0,\infty]$
    and $\zeta\in\FunnelTrajectories_{\hat{t}}$:
    \[
        \fMmax\geq\SNorm{\fM(\oTM(\zeta)|_{[0,\hat{t})})},\qquad
        \gMmax\geq\SNorm{\gM(\oTM(\zeta)|_{[0,\hat{t})})},
    \]
    and 
    \[
        \gMInvmax\geq\SNorm{\gM(\oTM(\zeta)|_{[0,\hat{t})})^{-1}}.
    \]
    If the function $\gM$ is, in addition, positive definite, i.e. $\al z, \gM({x}) z \ar > 0$ for all $x \in \R^{q}$ and  all $z \in \R^m \setminus \{0\}$,
    then there exists $\gMmin>0$ such that for all $z\in\R^m\backslash\cbl0\cbr$, all $\hat{t}\in [t_0,\infty]$,
    and $\zeta\in\FunnelTrajectories_{\hat{t}}$: 
    \[
        \gMmin\leq\frac{\al z, \gM(\oTM(\OpChi(\zeta))|_{[0,\hat{t})}(t))z\ar}{\Norm{z}^2}.
    \] 
\end{lemma}
\begin{proof}
    To prove the assertion, we invoke the continuity of the functions $\fM,\gM$ and the resulting boundedness on compact sets.
    By definition of~$\FunnelTrajectories_{\infty}$ and $\cD_{t}^\Psi$ in \eqref{eq:Def:FunnelTrajectories} and \eqref{eq:DefSetD}, we have for all $i=1,\ldots,r$
    \[
        \fa\zeta\in\FunnelTrajectories_{\infty}\fa t\geq 0:\quad\Norm{\eM_i(\zeta(t)-\OpChi(y_{\rf})(t))} < \Funnel_i(t).
    \]
    Due to the definition of the error variables~$\eM_i$ in \eqref{eq:ErrorVar} there exists an invertible matrix $S\in\R^{rm\times rm}$ such that
    \begin{align}\label{eq:reps_ei_chi}
        \begin{pmatrix} \eM_1(\zeta-\OpChi(y_{\rf})) \\ \vdots \\  \eM_r(\zeta-\OpChi(y_{\rf}))\end{pmatrix} 
        = S (\zeta-\OpChi(y_{\rf})).
    \end{align}
    Hence, by boundedness of $\Funnel_i$ and $y_{\rf}^{(i)}$ for all~$i=1,\ldots, r$,
    there exists a compact set $K\subset\R^{rm}$ with 
    \begin{equation}\label{eq:ExistenceCompactsetAllFunnelTrajectories}
        \fa \zeta\in\FunnelTrajectories_{\infty}\fa t\geq 0:\quad \zeta(t)\in K.
    \end{equation}
    Invoking the BIBO property of the operator~$\oTM$, there exists a compact set $K_q\subset\R^q$ with  
    $\oTM(z)(\Rp)\subset K_q$ for all $z\in\cR(\Rp,\R^{rm})$ with $z(\Rp)\subset K$.
    For arbitrary ${\hat{t}\in(0,\infty)}$ and $\zeta\in\FunnelTrajectories_{\hat{t}}$,
    we have $\zeta(t)\in K$ for all $t\in[0,\hat{t})$.
    For every element ${\zeta\in\FunnelTrajectories_{\hat{t}}}$, the restriction $\zeta|_{[0,\hat{t})}$ 
    can be extended to a function~$\tilde{\zeta}\in\cR(\Rp,\R^{rm})$ 
    with ${\tilde{\zeta}(t)\in K}$ for all $t\in\Rp$.
    We have $\oTM(\tilde{\zeta})(t)\in K_q$ for all $t\in\Rp$ because of the BIBO property of the operator~$\oTM$. 
    This implies $\oTM(\zeta)|_{[0,\hat{t})}(t)\in K_q$ for all $t\in [0,\hat{t})$ and $\zeta\in\FunnelTrajectories_{\hat{t}}$ since $\oTM$ is causal.
    Since~$\fM(\cdot)$, $\gM(\cdot)$, and $\gM^{-1}(\cdot)$ are continuous, the constants ${\fMmax=\max_{z\in K_q}\Norm{\fM(z)}}$,
    $\gMmax=\max_{z\in K_q}\Norm{\gM(z)}$
    and $\gMInvmax=\max_{z\in K_q}\Norm{\gM(z)^{-1}}$ are well-defined.
    For all $\hat{t}\in[t_0,\infty]$ and $\zeta\in\FunnelTrajectories_{\hat{t}}$, we have 
    \[
        \fa t\in [0,\hat{t}):\ \oTM(\zeta)(t)\in K_q.
    \]
    Furthermore, if $\gM(x)$ is positive definite for every $x\in K_q$, then there exists $\gMmin>0$ such that 
    $ \gMmin \leq \frac{\al z, \gM(\oTM(\OpChi(\zeta)))|_{[0,\delta)}(t))z\ar }{\Norm{z}^2}$ for all $z\in\R^m\backslash\cbl0\cbr$, 
    which proves the assertion.
\end{proof}

To prove the existence of $\gMmin>0$ in~\Cref{Lemma:DynamicBounded}, it is assumed that $\gM$ is positive definite.
In general, this is not the case when considering a model $(\fM,\gM,\oTM)\in\cM^{m,r}_{t_0}$.
However, we will restrict the model class $\cM^{m,r}_{t_0}$ and utilise this result in~\Cref{Chapter:DiscretFMPC}.

We are now in the position to prove the existence of a sufficiently large
$\umax\geq0$ such that the set of controls $\Controls(\umax,\InitState)$ is
non-empty.
\begin{theorem}\label{Th:ExUmax}
    Consider model~\eqref{eq:Model_r} with $(\fM,\gM,\oTM)\in\cM^{m,r}_{t_0}$.
    Let $\tau\geq 0$ be greater than or equal to the memory limit of operator $\oTM$.
    Further, let $y_{\rf}\in W^{r,\infty}(\Rp,\R^{m})$ and $\Psi\in\FunnelBoundaryFuncs$. 
    Then, there exists $\umax\geq0$ such that, 
    for $\hat{t}\geq t_0$, $\InitState\in\InitValues(\hat{t})$, and  all $T>0$, we have
    \begin{equation}\label{eq:ExistanceControl}
        \Controls(\umax,\InitState)\neq\emptyset.
    \end{equation}
\end{theorem}
\begin{proof}
    \emph{Step 1}: We define a candidate value of $\umax\geq0$. 
    To that end, define, for ${i=1,\ldots, r-1}$ and $j=0,\ldots,r-i-1$
    \[
        \mu_i^0 \coloneqq  \SNorm{\Funnel_i},\quad \mu_{i}^{j+1}\coloneqq  \mu_{i+1}^{j}+k_{i}\mu_{i}^{j},
    \]
    where $k_i\geq0$, for $i=1,\ldots r-1$, are the to $\Psi$ associated constants,
    which are also used to define the error variables $\eM_i$ as in~\eqref{eq:ErrorVar}.
    Using the constants $\fMmax$ and $\gMInvmax$ from \Cref{Lemma:DynamicBounded}, define
    \[
        \umax\coloneqq \gMInvmax
        \rbl
            \fMmax+\SNorm{y_{\rf}^{(r)}}+\sum_{j=1}^{r-1}k_j\mu_{j}^{r-j}+\SNorm{\dot{\Funnel}_r}
        \rbr.
    \]

    \noindent
    \emph{Step 2}: Let $T>0$, $\hat{t}\geq t_0$, and $(\xMh,\oTMh)=\InitState\in\InitValues(\hat{t})$ be arbitrary but fixed.
    We construct a control function $u$ and show that $u\in\Controls(\umax,\InitState)$.
    To this end, for some $u\in L^{\infty}([\hat{t},\hat{t}+T],\R^m)$, we use the shorthand notation
    $\xM(t)\coloneqq  \xM(t;\hat{t},\InitState,u)$ and $\eM_i(t)\coloneqq  \eM_{i}(\xM(t)-\OpChi(y_{\rf})(t))$
    for $i=1,\ldots,r$.
    The application of the feedback control 
    \[
        u(t)\coloneqq \gM(\oTM(\xM)(t))^{-1}
       \rbl 
            -\fM(\oTM(\xM)(t))+y^{(r)}_{\rf}(t)
            -\sum_{j=1}^{r-1} k_{j}\eM^{(r-j)}_j(t)+\eM_r(t)\tfrac{\dot{\Funnel}_r(t)}{\Funnel_r(t)}
        \rbr
    \]
    to the system~\eqref{eq:Model_r} leads to a closed-loop system.
    If this initial value problem is considered on the interval~$[\hat{t},\hat{t}+T]$ with initial conditions $(\hat{t},\InitState)$ as in \eqref{eq:Model_InitialValue},
    then an application of~\Cref{Prop:SolutionExists}
    yields the existence of a maximal solution~$\xM:[0,\omega)\to\R^{rm}$ with ${\omega>\hat{t}}$ in the sense of~\Cref{Def:ModSolution}.
    If~$\xM$ is bounded, then $\omega=\infty$, see~\Cref{Prop:SolutionExists}~\ref{Item:Theorem:BoundedSolution}.
    In this case, the solution exists on $[0, \hat t + T]$.
    Utilising~\eqref{eq:ErrorVarDyn}, one can show by induction that
    \[
        \eM_r(t)=\eM_1^{(r-1)}(t)+\sum_{j=1}^{r-1}k_j\eM_j^{(r-j-1)}(t).
    \]
    Omitting the dependency on $t$, we calculate for $t \in [\hat{t},\omega)$:
    \begin{align*}
    \frac{\dot{\eM}_r\Funnel_r-\eM_r\dot{\Funnel}_r}{\Funnel_r}
    &=\eM_1^{(r)}+\sum_{j=1}^{r-1}k_j\eM_j^{(r-j)}-\eM_r\frac{\dot{\Funnel}_r}{\Funnel_r}\\
    &=\fM(\oTM(\xM))+\gM(\oTM(\xM))u-y_{\rf}^{(r)}+\sum_{j=1}^{r-1}k_j\eM_j^{(r-j)}-\eM_r\frac{\dot{\Funnel}_r}{\Funnel_r}
    =0.
    \end{align*}
    Therefore,
    \[
    \dd{t}\tfrac{1}{2}\Norm{\frac{\eM_r}{\Funnel_r}}^2
    = \al\frac{\eM_r}{\Funnel_r},\frac{\dot{\eM}_r\Funnel_r-\eM_r\dot{\Funnel}_r}{\Funnel_r^2}\ar
    =0.
    \]
    We have $\Norm{\tfrac{\eM_r(\hat{t})}{\Funnel_r(\hat{t})}}<1$ by the assumption $\xM(\hat{t}) - \OpChi(y_{\rf})(\hat{t})\in\cD_{\hat{t}}^\Psi$, see also~\Cref{Rem:InitialValueInFunnel}.
    This yields $\Norm{\tfrac{\eM_r(t)}{\Funnel_r(t)}}<1$ for all $t\in[\hat{t},\omega)$.
    This implies, according to \Cref{Prop:OnlyLastFunnel},  $\xM(t)-\OpChi(y_{\rf})(t)\in\cD_{t}^\Psi$ for all $t\in[\hat{t},\omega)$, 
    i.e. $\Norm{\eM_i(t)}<\Funnel_i(t)$ for all $i=1,\ldots, r$.
    Thus, $\Norm{\eM_i(t)}\leq\mu_i^0$ for all $i=1,\ldots, r$.
    Invoking boundedness of $y_{\rf}^{(i)}$, $i=0,\ldots,r$, and the relation in~\eqref{eq:reps_ei_chi}, we may infer that $\xM$ is bounded on $[\hat{t},\omega)$.
    Hence, $\omega=\infty$. 
    Since ${(\xMh,\oTMh)=\InitState\in\InitValues(\hat{t})}$, there exists 
    a function $\zeta\in\FunnelTrajectories_{\hat{t}}$ such that
    ${\zeta|_{[\hat{t}-\tau,\hat{t}]\cap[0,\hat{t}]}=\xMh}$
    and $\oTM(\zeta)|_{[\hat{t}-\tau,\hat{t}]\cap[t_0,\hat{t}]}=\oTMh$.
    Moreover, 
    as $\xM$ fulfils the initial conditions~\eqref{eq:Model_InitialValue},
    we have 
    $\xM(t)|_{[\hat{t}-\tau,\hat{t}]\cap[0,\hat{t}]}=\xMh$
    and $\oTM(\xM)|_{[\hat{t}-\tau,\hat{t}]\cap[t_0,\hat{t}]}=\oTMh$.
    Define the regulated function ${\tilde{\zeta}\in\cR(\Rp,\R^{rm})}$ by 
    \[
        \tilde{\zeta}(t)=\begin{cases}
            \xM(t),& t\in [\hat{t},\hat{t}+T]\\
            \zeta(t),& t\in \Rp\backslash[\hat{t},\hat{t}+T].
        \end{cases}
    \]
    The function $\tilde{\zeta}$ is an element of $\FunnelTrajectories_{\hat{t}+T}$ because 
    $\zeta\in\FunnelTrajectories_{\hat{t}}$ and 
    $\xM(t)-\OpChi(y_{\rf})(t)\in\cD_{t}^\Psi$ for all $t\in[\hat{t},\hat{t}+T]$.
    Therefore,
    $\Norm{\fM(\oTM(\tilde{\zeta})(t))}\leq \fMmax$ and 
    $\Norm{\gM(\oTM(\tilde{\zeta})(t))^{-1}}\leq \gMInvmax$ for all $t\in[\hat{t},\hat{t}+T]$ 
    according to~\Cref{Lemma:DynamicBounded}.
    Since $\tau\geq 0$ is greater than or equal to the memory limit of operator $\oTM$, we have
    \[
        \oTM(\xM)(t)=\oTM(\tilde{\zeta})(t)
    \]
    for all $t\in[\hat{t},\hat{t}+T]$. Thus, 
    $\Norm{\fM(\oTM(\xM)(t))}\leq \fMmax$ and 
    $\Norm{\gM(\oTM(\xM)(t))^{-1}}\leq \gMInvmax$ for all $t\in[\hat{t},\hat{t}+T]$. 
    Finally, using~\eqref{eq:ErrorVarDyn} and the definition of~$\mu_i^j$, it follows that
    \[
        \Norm{\eM^{(j+1)}_i(t)}
        =
        \Norm{\eM^{(j)}_{i+1}(t)-k_i\eM^{(j)}_{i}(t)}
        \leq
        \mu_{i+1}^j+k_i\mu^{j}_{i}
        =\mu^{j+1}_{i}
    \]
    inductively for all $i=1,\ldots, r-1$ and $j=0,\ldots,r-i-1$. Thus, by definition of~$u$ and~$\umax$, we have $\SNorm{u}\leq \umax$ and hence $u\in\Controls(\umax,\InitState)$.
\end{proof}

\begin{remark}
    Note that  $\umax\geq0$ in~\Cref{Th:ExUmax} is independent of the time
    $\hat{t}\geq t_0$, the initial value $\InitState\in\InitValues(\hat{t})$,
    and the considered time horizon $T>0$.
    It solely depends on the system dynamics, i.e. the functions
    $(\fM,\gM,\oTM)\in\cM^{m,r}_{t_0}$, the reference $y_{\rf}$, the funnel
    functions $\Psi=(\Funnel_1,\ldots,\Funnel_r)\in\FunnelBoundaryFuncs$, and the
    associated parameters $k_i$, $i=1,\ldots, r-1$.
    \null
\end{remark}

We have seen in this section that at there exists a control function $u\in\Controls(\umax,\InitState)$
at every time $\hat{t}\geq t_0$, assuming that the initial value $\InitState$ is feasible for the model~\eqref{eq:Model_r}
and that $\umax\geq0$ is large enough. Such a control function solves the tracking problem formulated in~\Cref{Sec:ControlObjective}
for the model~\eqref{eq:Model_r}.  Moreover, it was shown that the state of the model at every time during application of the control $u$ 
is a feasible initial state for the model. 
Since the initial value given by $\yM^0$ is feasible for the model at time $\hat{t}=t_0$, 
this means that an iterative application of controls $u\in\Controls(\umax,\InitState)$ to the model~\eqref{eq:Model_r}
solves the tracking problem for the concatenated solution.

\subsection{Optimal control problem}\label{Sec:OptimalControlProblem}
We are now in the position to address the problem of solving the tracking
problem from~\Cref{Sec:ControlObjective} by means of an optimal control problem
utilising the concept of funnel stage cost functions from~\Cref{Sec:FunnelStageCostFunctions}.
Given auxiliary funnel
functions ${\Psi=(\Funnel_1,\ldots,\Funnel_r)\in\FunnelBoundaryFuncs}$ 
around the reference trajectory $y_{\rf}\in W^{r,\infty}(\Rp,\R^{m})$
and the corresponding error variables~$\eM_i$ for $i=1,\ldots, r$ as
in~\eqref{eq:ErrorVar}, we saw in~\Cref{Sec:FunnelStageCostFunctions} that it is
sufficient to ensure that the last  auxiliary error $\eM_r$ evolves within its
funnel given by~$\Funnel_r$ in order to guarantee that all~$\eM_i$ evolve within
their respective funnels defined by~$\Funnel_i$ for $i=1,\ldots, r-1$.
Therefore, choose a funnel stage cost $\FunnelStageCost$ for the last auxiliary funnel function $\Funnel_r$.
Let $\tau\geq0$ be greater than or equal to the memory limit of operator~$\oTM$.
Then, define for $T>0$, $\hat{t}\geq t_0$, and  $\InitState\in\InitValues(\hat{t})$, the \emph{cost functional} 
$J^{\Psi}_T(\cdot;\hat{t},\InitState):L^\infty([\hat{t},\hat{t}+T],\R^{m})\to\R\cup\{\infty\}$ by
\begin{align}\label{eq:DefCostFunctionJ}
    J^{\Psi}_T(u;\hat{t},\InitState)\coloneqq 
        \int_{\hat{t}}^{\hat{t} + T}\FunnelStageCost(s,\eM_{r}(\xM(s;\hat{t},\InitState,u)-\OpChi(y_{\rf})(s)),u(s))\d{s}.
\end{align}
Although it is known that, for every $u\in L^\infty([\hat{t},\hat{t}+T],\R^{m})$, there exists a maximal solution~$\xM(t;\hat{t},\InitState,u)$,
according to~\Cref{Prop:SolutionExists}, this solution might have finite escape time, i.e. $\omega<\hat{t}+T$. 
In this case, and whenever the Lebesgue integral in~\eqref{eq:DefCostFunctionJ} does not exist, 
(i.e. both the Lebesgue integrals of the positive and negative part of
$\FunnelStageCost(s,\eM_{r}(\xM(s;\hat{t},\InitState,u)-\OpChi(y_{\rf})(s)),u(s))$ are infinite), then its value
is treated as infinity.
Further, note that the solution~$\xM(t;\hat{t},\InitState,u)$ is unique on the interval $[\hat{t},\hat{t}+T]$,
according to~\Cref{Prop:SolutionUnique}, rendering $J^{\Psi}_T(u;\hat{t},\InitState)$ well-defined.
In the following, we will study properties of $J^{\Psi}_T(\cdot;\hat{t},\InitState)$ a bit more closely
and analyse the associated
\emph{Optimal Control Problem (OCP)}
\begin{equation}\label{eq:OptimalControlProblem}
            \mathop{\operatorname{minimise}}_{
                \substack
                {
                    u\in L^{\infty}([\hat{t}, \hat{t}+T],\R^{m}),\\
                    \SNorm{u}  \leq \umax 
                }
            }
                J^{\Psi}_T(u;\hat{t},\InitState),
\end{equation}
where $\umax\geq 0$ is a bound on the maximal control input. 
We will prove that the OCP \eqref{eq:OptimalControlProblem} has a solution 
and that this solution is an element of $\Controls(\umax,\InitState)$. 
Thus, it solves the tracking problem from~\Cref{Sec:ControlObjective},
according to our considerations in~\Cref{Sec:InitialRecFeasibilty}.

The cost function $J^{\Psi}_T(u;\hat{t},\InitState)$ is defined
in~\eqref{eq:DefCostFunctionJ} as the integral of funnel stage cost
$\FunnelStageCost$ evaluated over the auxiliary error~$\eM_r$. 
The concept of funnel stage cost functions
from~\Cref{Sec:FunnelStageCostFunctions} was based on the usage of Lipschitz
paths. It has already been proven in~\Cref{Prop:SolutionIsLPath} that the
solution trajectories of the model are Lipschitz continuous. It is evident that
this is also the case for the auxiliary errors $\eM_i$. Nevertheless, we will
briefly formalise this.
\begin{prop}\label{Prop:ErrorIsLPath}
    Consider model~\eqref{eq:Model_r} with $(\fM,\gM,\oTM)\in\cM^{m,r}_{t_0}$ with reference trajectory $y_{\rf}\in W^{r,\infty}(\Rp,\R^m)$.
    Let $\Psi\in\FunnelBoundaryFuncs$, $\hat{t}\geq t_0$, $\tau\in\Rp$,  and $\InitState\in \InitValues(\hat{t})$.
    Moreover, let $u\in L^\infty_{\loc}([\hat{t},\infty), \R^m)$ be a control
    such that initial value problem~\eqref{eq:Model_r} has a solution
    $\xM:[0,\omega)\to \R^{rm}$ with $\omega>\hat{t}$ in the sense
    of~\Cref{Def:ModSolution}. Then, for every $T\in(0,\omega-\hat{t})$, the
    restriction
    $\eM_i(\xM-\OpChi(y_{\rf}))|_{[\hat{t},\hat{t}+T]}:[\hat{t},\hat{t}+T]\to\R^{m}$
    is a Lipschitz path for all $i=1,\ldots, r$. 
\end{prop}
\begin{proof}
    For $y_{\rf}\in W^{1,\infty}$, the function $\OpChi(y_{\rf})(\cdot)$ is a Lipschitz continuous on every compact interval.
    Therefore, due to the definition of the error variables~$\eM_i$ for $i=1,\ldots, r$ in~\eqref{eq:ErrorVar},
    the statement of~\Cref{Prop:ErrorIsLPath} is an immediate consequence of~\Cref{Prop:SolutionIsLPath}.
\end{proof}
    The following~\Cref{Th:JFiniteCost} not only shows that $J^{\Psi}_T(u;\hat{t},\InitState)$ has a finite value 
    for all $u\in\Controls(\umax,\InitState)$,
    which is to be expected since this is the set of controls ensuring the evolution  
    of the auxiliary errors~$\eM_i$ within their respective funnels,
    see~\eqref{eq:Def-U}, but that $\Controls(\umax,\InitState)$ 
    is in fact the set of controls for which $J^{\Psi}_T(u;\hat{t},\InitState)$ is finite.
    
\begin{theorem}\label{Th:JFiniteCost}
    Consider model~\eqref{eq:Model_r} with $(\fM,\gM,\oTM)\in\cM^{m,r}_{t_0}$ with reference trajectory $y_{\rf}\in W^{r,\infty}(\Rp,\R^m)$.
    Let $\Psi\in\FunnelBoundaryFuncs$ and $\tau\geq0$ be greater than or equal to the memory limit of operator~$\oTM$.
    Further, let  $\hat{t}\geq t_0$, $\InitState\in \InitValues(\hat{t})$, $T> 0$,
    and $\umax\geq0$ such that $\Controls(\umax,\InitState)\neq\emptyset$.
    Then, the following identity holds:
    \begin{align*}
        \Controls(\umax,\InitState)&=
        \setdef
            { u\in  L^{\infty}([\hat{t},\hat{t}+T],\R^m)}
            {
                J^{\Psi}_T(u;\hat{t},\InitState)<\infty,\ \SNorm{u} \le \umax 
            }.
    \end{align*}
\end{theorem}
\begin{proof}
    Since $\FunnelStageCost$ is a funnel stage cost, it has the form 
    \begin{align*}
        \FunnelStageCost:\Rp\times\R^{m}\times\R^{m}\to\R\cup\cbl\infty\cbr,\quad
        (t,z,u)\mapsto {\FunnelPenaltyFuncR}(t,z)+\lambda_u\Norm{u}^2,
    \end{align*}
    where $\FunnelPenaltyFuncR$ is a funnel penalty function  and $\lambda_{u}\geq 0$, see~\Cref{Def:FunnelStageCostFunc}.
    Note that, for $u\in L^{\infty}([\hat{t},\hat{t}+T],\R^m)$  such that
    $\xM(s;\hat{t},\InitState,u)$ satisfies~\eqref{eq:ModDiff} for all
    $s\in[\hat{t},\hat{t}+T]$, the function
    $\eM_{r}(s)\coloneqq \eM_{r}(\xM(s;\hat{t},\InitState,u)-\OpChi(y_{\rf})(s)$ is,
    according to~\Cref{Prop:ErrorIsLPath}, a Lipschitz path with
    $(\hat{t},\eM_{r}(\hat{t}))\in\cF_{\UltMFunnel}$ since $\InitState\in
    \InitValues(\hat{t})$, see~\Cref{Rem:InitialValueInFunnel}.
    
    For $u\in \Controls(\umax,\InitState)$,
    we have $\SNorm{u}\leq\umax$ and
    $\xM(s;\hat{t},\InitState,u)-\OpChi(y_{\rf})(s)\in\cD_{t}^\Psi$ for all
    $s\in[\hat{t},\hat{t}+T]$.
    In particular, this implies $(s,\eM_{r}(s))\in\cF_{\UltMFunnel}$ for all $s\in[\hat{t},\hat{t}+T]$.
    Thus, 
    \[ 
        J^{\Psi}_T(u;\hat{t},\InitState)=
        \int_{\hat{t}}^{\hat{t} + T}\FunnelStageCost(s,\eM_{r}(s),u(s))\d{s}=
        \int_{\hat{t}}^{\hat{t} + T}\FunnelPenaltyFuncR(s,\eM_{r}(s))+\lambda_{u}\Norm{u(s)}^2\d{s}<\infty
    \]
    due to the boundedness of $u$ and the property~\ref{Item:FunnelPenaltyIntegral} 
    of the funnel penalty function $\FunnelPenaltyFuncR$, see~\Cref{Def:FunnelPenaltyFunction}.
    
    Conversely, let $u\in L^{\infty}([\hat{t},\hat{t}+T],\R^m)$ with $\SNorm{u}\leq\umax$ such that the cost functional 
    ${J^{\Psi}_T(u;\hat{t},\InitState)=\int_{\hat{t}}^{\hat{t} + T}\FunnelStageCost(s,\eM_{r}(s),u(s))\d{s}}$ is finite.
    Since $u$ is bounded and 
    \[
    {\FunnelStageCost(s,\eM_{r}(s),u(s))=\FunnelPenaltyFuncR(s,\eM_{r}(s))+\lambda_u\Norm{u(s)}}
    \]
    for all $s\in[\hat{t},\hat{t}+T]$, we have
    $\int_{\hat{t}}^{\hat{t} + T}\FunnelPenaltyFuncR(s,\eM_{r}(s))\d{s}<\infty$.
    Thus, $\Norm{\eM_{r}(s)}<\UltMFunnel(s)$ for all $s\in[\hat{t},\hat{t}+T]$
    because the property~\ref{Item:FunnelPenaltyIntegral} 
    of the funnel penalty function $\FunnelPenaltyFuncR$, see~\Cref{Def:FunnelPenaltyFunction}.
    This implies $\xM(s;\hat{t},\InitState,u)-\OpChi(y_{\rf})(s)\in\cD_{t}^\Psi$ for all $s\in[\hat{t},\hat{t}+T]$
    according to \Cref{Prop:OnlyLastFunnel}.
    Hence, $u\in\Controls(\umax,\InitState)$.
\end{proof}

\begin{remark}\label{Rem:PropertiesJ}
    The following statements hold under the assumptions of~\Cref{Th:JFiniteCost}:
    \begin{enumerate}[(a)]
        \item $0\leq J^{\Psi}_T(u;\hat{t},\InitState)<\infty$ for all $u\in \Controls(\umax,\InitState)$ because the
            funnel stage cost~$\FunnelStageCost$ is non-negative while the error~$\eM_{r}$ evolves within its funnel given by $\Funnel_{r}$,
            see~\Cref{Rem:FunnelStageCostPositive}.
        \item
        The optimal control problem~\eqref{eq:OptimalControlProblem} can be reformulated as
        \[
            \mathop{\operatorname{minimise}}_{
                u \in \Controls(\umax,\InitState)
            }
            J^{\Psi}_T(u;\hat{t},\InitState).
        \]
    \end{enumerate}
\end{remark}
If the initial value $\InitState$ is feasible for the model~\eqref{eq:Model_r}, then any control function $u$ with
${J^{\Psi}_T(u;\hat{t},\InitState)<\infty}$ guarantees that, if applied to the model~\eqref{eq:Model_r}, 
all errors $\eM_{i}$  remain (strictly) within their respective funnels~$\Funnel_i$, $i=1,\ldots, r$.
Since $J^{\Psi}_T(u;\hat{t},\InitState)$ is non-negative 
for all control functions $u \in \Controls(\umax,\InitState)$, this raises the question as to whether there exists an
optimal  $u^{\star}$ which minimises $J^{\Psi}_T(u;\hat{t},\InitState)$ and is a solution to the optimal
control problem~\eqref{eq:OptimalControlProblem}. The answer is affirmative and shown in the next
theorem.

\begin{theorem}\label{Th:SolutionExists}
    Consider model~\eqref{eq:Model_r} with $(\fM,\gM,\oTM)\in\cM^{m,r}_{t_0}$ with reference trajectory $y_{\rf}\in W^{r,\infty}(\Rp,\R^m)$.
    Let $\Psi\in\FunnelBoundaryFuncs$ and $\tau\geq0$ be greater than or equal to the memory limit of operator~$\oTM$.
    Further, let  $\hat{t}\geq t_0$, $(\xMh,\oTMh)=\InitState\in \InitValues(\hat{t})$, $T> 0$,
    and $\umax\geq0$ such that $\Controls(\umax,\InitState)\neq\emptyset$.
    Then, there exists a function  $u^{\star} \in \Controls(\umax,\InitState)$ such that
    \[
            J^{\Psi}_T(u^\star;\hat{t},\InitState) =
            \mathop{\min}_{
                u \in \Controls(\umax,\InitState)
            }
            J^{\Psi}_T(u;\hat{t},\InitState) =
            \mathop{\min}_{
                \substack
                {
                    u\in L^{\infty}([\hat{t}, \hat{t}+T],\R^{m}),\\
                    \SNorm{u}  \leq \umax 
                }
            }
                J^{\Psi}_T(u;\hat{t},\InitState).
    \]
\end{theorem}
\begin{proof}
    The proof essentially follows the lines of~\cite[Prop.~2.2]{sakamoto2023}.\\
    It follows from \Cref{Rem:PropertiesJ} that $J^{\Psi}_T(u;\hat{t},\InitState)\geq 0$ for all $ u \in
    \Controls(\umax,\InitState)$.
    Hence, the infimum ${J^{\star}: = \inf_{u \in \Controls(\umax,\InitState)} J^{\Psi}_T(u;\hat{t},\InitState)}$ exists.
    Let $(u_k)\in\rbl\Controls(\umax,\InitState)\rbr^\N$ be a minimising sequence,
    meaning $J^{\Psi}_T(u_k;\hat{t},\InitState)\to J^{\star}$.
    By definition of $\Controls(\umax,\InitState)$, we have ${\SNorm{ u_k }\leq \umax}$ for all~$k\in\N$.
    Since $L^\infty([\hat{t},\hat{t}+T],\R^{m})\sub L^2([\hat{t},\hat{t}+T],\R^{m})$,
    we conclude  that $(u_k)$ is a bounded sequence in the Hilbert space~$L^2$. Thus, there
    exists a function $u^{\star}\in L^2([\hat{t},\hat{t}+T],\R^{m})$ and a weakly convergent subsequence
    $u_k\rightharpoonup u^{\star}$ (which we do not relabel). More precisely, $u_k|_{[\hat{t},t]}
    \rightharpoonup u^{\star}|_{[\hat{t},t]}$ weakly in $L^2([\hat{t},t],\R^m)$ for all $t\in [\hat{t},\hat{t}+T]$ as a
    straightforward argument shows.
    We define $(x_k)\coloneqq \rbl\xM(\cdot;\hat{t},\InitState,u_k)\rbr\in \cR([0,\hat{t}+T],\R^n)^\N$
    as the sequence of associated responses. Note that,
    although we are only considering the optimal control problem~\eqref{eq:DefCostFunctionJ} on the interval $[\hat{t},\hat{t}+T]$,
    the functions $\xM(\cdot;\hat{t},\InitState,u_k)$ 
    are defined on the entire interval $[0,\hat{t}+T]$ for all $k\in\N$ as they are solutions 
    of the differential equation~\eqref{eq:Model_r} in the sense of~\Cref{Def:ModSolution}.
    This allows us to formally evaluate $\oTM(x_k)$.
    We show the assertion of \Cref{Th:SolutionExists} in the following seven steps.
    
    \noindent
    \emph{Step 1}:
    We construct a uniformly bounded sequence of solutions 
    of the differential equation~\eqref{eq:Model_r} in the sense of~\Cref{Def:ModSolution}.
    Since $(\xMh,\oTMh)=\InitState\in\InitValues(\hat{t})$, there exists 
    a function $\zeta\in\FunnelTrajectories_{\hat{t}}$  such that
    $\zeta|_{[\hat{t}-\tau,\hat{t}]\cap[0,\hat{t}]}=\xMh$
    and $\oTM(\zeta)|_{[\hat{t}-\tau,\hat{t}]\cap[t_0,\hat{t}]}=\oTMh$.
    Moreover, we have 
    $x_{k}|_{[\hat{t}-\tau,\hat{t}]\cap[0,\hat{t}]}=\xMh$
    and $\oTM(x_{k})|_{[\hat{t}-\tau,\hat{t}]\cap[t_0,\hat{t}]}=\oTMh$ 
    because $x_{k}$ fulfils the initial conditions~\eqref{eq:Model_InitialValue}.
    Define the function $\tilde{x}_k\in\cR([0,\hat{t}+T],\R^{rm})$ by 
    \begin{equation}\label{eq:SolutionConstructionXK}
        \tilde{x}_k(t)=\begin{cases}
            \zeta(t),& t\in [0,\hat{t})\\
            x_k(t),& t\in   [\hat{t},\hat{t}+T].
        \end{cases}
    \end{equation}
    The function $\tilde{x}_k$ is, by construction, a solution of the differential equation~\eqref{eq:ModDiff} 
    with initial values $\InitState$ in the sense of~\Cref{Def:ModSolution}.
    We have $\tilde{x}_k|_{[\hat{t},\hat{t}+T]}=x_k$ and additionally
    ${\oTM(x_k)|_{[\hat{t},\hat{t}+T]}=\oTM(\tilde{x}_k)|_{[\hat{t},\hat{t}+T]}}$
    because $\tau\geq 0$ is greater than or equal to the memory limit of operator $\oTM$.
    Without loss of generality, we therefore assume in the following that $x_k$ has the form \eqref{eq:SolutionConstructionXK}, i.e.
    we relabel $\tilde{x}_k$ as $x_k$.
    By $u_k\in \Controls(\umax,\InitState)$, we have $x_k(t)-\OpChi(y_{\rf})(t)\in\cD_{t}^{\Psi}$ for all $t\in[\hat{t},\hat{t}+T]$.
    Thus, $x_k\in\FunnelTrajectories_{\hat{t}+T}$ because $\zeta\in\FunnelTrajectories_{\hat{t}}$.
    This implies $x_k(t)-\OpChi(y_{\rf})(t)\in\cD_{t}^{\Psi}$ for all $t\in[0,\hat{t}+T]$.
    Invoking boundedness of $y_{\rf}^{(i)}$, for $i=0,\ldots,r$,
    and the relation in~\eqref{eq:reps_ei_chi}, we may infer that $x_k$ is uniformly bounded 
    on the entire interval $[0,\hat{t}+T]$.

\noindent
    \emph{Step 2}:
    We show that the sequence of restrictions $(x_k|_{[\hat{t},\hat{t}+T]})$ 
    is uniformly equicontinuous on the interval $[\hat{t},\hat{t}+T]$.
    As $x_k$ is a solution of~\eqref{eq:Model_r} in the sense of~\Cref{Def:ModSolution}, we have
    \begin{equation}\label{eq:SolutionRepresentation}
        x_k(t)=\xMh(\hat{t})+\int_{\hat{t}}^{t}\FM(x_k(s),\oTM(x_k)(s))+\GM(\oTM(x_k)(s))u_k(s)\d s,
    \end{equation}
    for all $k\in\N$ and $t\in [\hat{t},\hat{t}+T]$, where $\FM$ and $\GM$ are defined as in~\Cref{Def:ModSolution}.
    Since the sequence $(u_k)$ is bounded,
    $\bar{u}\coloneqq \sup_{k\in\N} \LNorm{ u_k}$ exists. 
    Furthermore, using the considerations in from~\Cref{Lemma:DynamicBounded}, there exists 
    constants $\FM^{\max},\GM^{\max}\geq 0$ such that 
    $\FM^{\max}\geq\SNorm{\FM(\tilde{\zeta},\oTM(\tilde{\zeta}))}$ and 
    $\GM^{\max}\geq\SNorm{\GM(\oTM(\tilde{\zeta}))}$ 
    for all $\tilde{\zeta}\in\FunnelTrajectories_{\hat{t}+T}$.
    Now, let $\eps>0$ and define $\varsigma\coloneqq \min\cbl1,\frac{1}{\eps} (\FM^{\max}+  \GM^{\max} \bar{u})\cbr$.
    Let $k\in \N$ and $t_1,t_2\in[\hat{t},\hat{t}+T]$ such that $|t_2-t_1|<\varsigma^2$. 
    Then, using $x_k\in\FunnelTrajectories_{\hat{t}+T}$ and H\"older's inequality in the third estimate,
    \begin{align*}
        \Norm{x_k(t_2)-x_k(t_1)}
        &\leq \int_{t_1}^{t_2}\Norm{ \FM(x_k(s),\oTM(x_k)(s))}+ \Norm{\GM(\oTM(x_k)(s))} \Norm{u_k(s)} \d{s}\\
        &\leq \FM^{\max}\Abs{t_2 - t_1}        + \GM^{\max}\int_{t_1}^{t_2}\Norm{u_k(s)} \d{s}\\
        &\leq \FM^{\max}\sqrt{\Abs{t_2 - t_1}} + \GM^{\max}\sqrt{\Abs{t_2 - t_1}}\LNorm{u_k}\\
        &\leq \sqrt{\Abs{t_2 - t_1}}(\FM^{\max}+\GM^{\max}\bar{u})
        < \eps,
    \end{align*}
    which shows that $(x_k|_{[\hat{t},\hat{t}+T]})$ is uniformly equicontinuous.

\noindent
    \emph{Step 3}:
    By the Arzel\`{a}-Ascoli theorem,
    there exists a function $x^{\star}\in \cC([\hat{t},\hat{t}+T],\R^{rm})$ and a 
    subsequence (which we do not relabel) such that the restriction $x_{k}|_{[\hat{t},\hat{t}+T]}$ 
    to the interval $[\hat{t},\hat{t}+T]$ is uniformly convergent, i.e.  $x_{k}|_{[\hat{t},\hat{t}+T]}\to x^{\star}$ .
    As in~\eqref{eq:SolutionConstructionXK}, we extend $x^{\star}$ by $\zeta$ on the interval $[0,\hat{t}]$ (we do not relabel $x^\star$). 
    By construction~\eqref{eq:SolutionConstructionXK}, we have $x_k(t)=\zeta(t)=x^\star(t)$ for all $t\in[t_0,\hat{t})$.
    Thus, $x_{k}$ converges uniformly to  $x^{\star}$ on the whole interval $[t_0,\hat{t}+T]$.
    Now we prove that $x^{\star}=\xM(\cdot;\hat{t},\InitState,u^{\star})$, which means to show that
    ${x^{\star}(t)=\xMh+\int_{\hat{t}}^{t}\FM(x^{\star}(s),\oTM(x^{\star})(s))+\GM(\oTM(x^{\star}(s))) u^{\star}(s)\d s}$ for all $t\in
    [\hat{t},\hat{t}+T]$. 
    On the interval $[\hat{t},\hat{t}+T]$, the values of $\oTM(x^{\star})(s)$  are completely determined by 
    $\InitState$ and $x^{\star}|_{[\hat{t},\hat{t}+T]}$ 
    since $\tau\geq 0$ is greater than or equal to the memory limit of operator $\oTM$.
    The same is true for $\oTM(x_k)(s)$ on the $[\hat{t},\hat{t}+T]$ for all $k\in\N$. 
    We, therefore, will in the following abuse the notation slightly by only writing $\FM(x^\star(s))$ and $\FM(x_k(s))$ instead of
    $\FM(x^{\star}(s),\oTM(x^{\star})(s))$ and  $\FM(x_k(s),\oTM(x_k(s)))$, respectively.
    We will use the same shorthand notation for $\GM$.
    Due to the representation $x_k$ as in~\eqref{eq:SolutionRepresentation}
    and since $x_k$ in particular converges pointwise to~$x^{\star}$ and the sequence $(\FM(x_k))$ 
    is uniformly bounded as $(x_k)$ is uniformly bounded and $\FM$ is continuous, the bounded convergence theorem gives that
    \[
        \fa t\in[\hat{t},\hat{t}+T]:\ \int_{\hat{t}}^{t}\FM(x_k(s)) \d s \To \int_{\hat{t}}^{t}\FM(x^{\star}(s)) \d s.
    \]
    Therefore, it remains to show
    \[
       \fa t\in[\hat{t},\hat{t}+T]:\  \int_{\hat{t}}^{t}\GM(x_k(s))u_k(s)\d s\To \int_{\hat{t}}^{t}\GM(x^{\star}(s))u^{\star}(s)\d s.
    \]
    The argument~$s$ is omitted in the following.
    Since $\GM(x^{\star})$ is bounded on $[\hat{t},\hat{t}+T]$,
    it is an element of $L^2([\hat{t},\hat{t}+T],\R^{n\times m})$, thus the weak
    convergence of~$(u_k)$ implies 
    ${\int_{\hat{t}}^{t}\GM(x^{\star})u_k\d s \to \int_{\hat{t}}^{t}\GM(x^{\star})u^{\star}\d s}$
    for all~$t\in[\hat{t},\hat{t}+T]$.
    Therefore, using H\"older's inequality  in the second estimate, we obtain, for all $t\in [\hat{t},\hat{t}+T]$,
    \begin{align*}
        \Norm{\int_{\hat{t}}^{t}\GM(x_k)u_k -\GM(x^{\star})u^{\star}\d s}
        &= \Norm{\int_{\hat{t}}^{t}\GM(x_k)u_k \!+\!\GM(x^{\star})u_k \!-\! \GM(x^{\star}) u_k \!-\! \GM(x^{\star})u^{\star}\d s } \\
        &\hspace{-4.85cm}\leq \int_{\hat{t}}^{t}\Norm{ \GM(x_k)-\GM(x^{\star})} \Norm{u_k}\d s
        +\Norm{\int_{\hat{t}}^{t}  \GM(x^{\star})u_k - \GM(x^{\star})u^{\star}\d s}\\
        &\hspace{-4.85cm}\leq \rbl\int_{\hat{t}}^{t} \Norm{ \GM(x_k)-\GM(x^{\star})}^2\d s
            \rbr^{\frac{1}{2}}\rbl\int_{\hat{t}}^{t}\Norm{u_k}^2\d s\rbr^{\frac{1}{2}}
            +\Norm{\int_{\hat{t}}^{t}\GM(x^{\star})u_k -\GM(x^{\star})u^{\star}\d s}\\
        &\hspace{-4.85cm}\leq \sup_{m\in\N}\LNorm{u_m}
              \underbrace{\rbl\int_{\hat{t}}^{t} \Norm{ \GM(x_k)-\GM(x^{\star})}^2\d
              s\rbr^{\frac{1}{2}}}_{\to 0}
        +\underbrace{\Norm{\int_{\hat{t}}^{t} \GM(x^{\star})u_k -\GM(x^{\star})u^{\star}\d s}}_{\to 0} \to 0.
    \end{align*}

\noindent
    \emph{Step 4}:
    We show $\SNorm{ u^{\star}}\leq \umax$. To this end, define the sets
    \[
        A_m\coloneqq \setdef{t\in [\hat{t},\hat{t}+T]}{\Norm{u^{\star} (t)}^2\ge \umax^2+\tfrac{1}{m}},
        \quad m\in\N.
    \]
    Let $\Indic_{A_m}$ denote the indicator function of the set $A_m$, then, since $u_k\rightharpoonup
    u^{\star}$, we have that
    $
        \langle u_k, \Indic_{A_m} u^{\star}\rangle_{L^2} \to \langle u^{\star}, \Indic_{A_m}
        u^{\star}\rangle_{L^2} = \|\Indic_{A_m} u^{\star}\|_{L^2}^2.
    $
    On the other hand, by the Cauchy-Schwarz inequality we have that
    $
        \langle u_k, \Indic_{A_m} u^{\star}\rangle_{L^2} \le \|\Indic_{A_m} u_k\|_{L^2} \|\Indic_{A_m}
        u^{\star}\|_{L^2},
    $
    thus
    \[
        \|\Indic_{A_m} u^{\star}\|_{L^2} 
         =  \|\Indic_{A_m} u^{\star}\|_{L^2}^{-1} \liminf_{k\to\infty} \langle u_k, \Indic_{A_m} u^{\star}\rangle_{L^2} 
        \le \liminf_{k\to\infty} \|\Indic_{A_m} u_k\|_{L^2}
    \]
    and hence
    $
        \int_{A_m} \|u^{\star}(s)\|^2 \d s \le \liminf_{k\to\infty}  \int_{A_m} \|u_k(s)\|^2 \d s.
    $
    Since $\|u_k\|_\infty\leq \umax$, we then find the following for all $m\in\N$ and $k\in\N$:
    \[
        \Lebesgue\rbl A_m\rbr
        =    \int_{A_m}1\d s
        \leq m\int_{A_m}\Norm{u^{\star}(s)}^2-\umax^2 \d s
        \leq m\int_{A_m}\Norm{u^{\star}(s)}^2-\Norm{u_k(s)}^2  \d s,
    \]
    where $\Lebesgue$ denotes the Lebesgue measure, thus
    \[
        {0 \le \Lebesgue\rbl A_m\rbr \le \liminf\limits_{k\to\infty} m\int_{A_m}\Norm{u^{\star}(s)}^2-\Norm{u_k(s)}^2 \d s \le 0}.
    \]
    Due to the $\sigma$-continuity of $\Lebesgue$ we get
    \[
        \Lebesgue\rbl\setdef
            {t\in [\hat{t},\hat{t}+T]}
            {\Norm{u^{\star} (t)}>\umax}\rbr
        = \Lebesgue\rbl\bigcup_{m\in\N}A_m\rbr
        = \lim_{m\to\infty}\Lebesgue(A_m)=0.
    \]
    This implies $\SNorm{ u^{\star}}\leq \umax$. 

    \noindent
    \emph{Step 5}:
    We prove $u^{\star} \in \Controls(\hat{t},\InitState)$, which means to show~$x^\star(t)-\OpChi(y_{\rf})(t)\in\cD_{t}^\Psi$
    for all $t\in[\hat{t},\hat{t}+T]$. 
    Since $\FunnelStageCost$ is a funnel stage cost, it has the form 
    \begin{align*}
        \FunnelStageCost:\Rp\times\R^{m}\times\R^{m}\to\R\cup\cbl\infty\cbr,\quad
        (t,z,u)\mapsto {\FunnelPenaltyFuncR}(t,z)+\lambda_u\Norm{u}^2,
    \end{align*}
    where $\FunnelPenaltyFuncR$ is a funnel penalty function  and $\lambda_{u}\geq 0$, see~\Cref{Def:FunnelStageCostFunc}.
    The supremum $\sup_{k\in\N}J^{\Psi}_T(u_k;\hat{t},\InitState)<\infty$ exists because $J^{\Psi}_T(u_k;\hat{t},\InitState)\to J^{\star}$.
    Thus, due to the uniform boundedness of $\Norm{u_k}$ and the definition of the function $J^{\Psi}_T$, see~\eqref{eq:DefCostFunctionJ},
    there exists $M\geq0$ such that 
    \[
        \int_{\hat{t}}^{\hat{t}+T}\FunnelPenaltyFuncR(s,\eM_{r}(x_k(s)-\OpChi(y_{\rf})(s)))\d{s}\leq M
    \]
    for all $k\in\N$.
    In the following, we use the shorthand notation $\eM_{r}^k(\cdot)\coloneqq \eM_{r}(x_k(\cdot)-\OpChi(y_{\rf})(\cdot))$
    and $\eM_{r}^\star(\cdot)\coloneqq \eM_{r}(x^\star(\cdot)-\OpChi(y_{\rf})(\cdot))$.
    The functions $\eM_{r}^k(\cdot)$ are, according to~\Cref{Prop:ErrorIsLPath}, Lipschitz paths 
    with $(\hat{t},\eM_{r}^k(\hat{t}))\in\cF_{\UltMFunnel}$ 
    since $\InitState\in \InitValues(\hat{t})$, see~\Cref{Rem:InitialValueInFunnel}.
    The uniform convergence of $x_k$ to $x^{\star}$ implies the uniform convergence of $\eM_{r}^k$
    to $\eM_{r}^\star$. Therefore, ${(\hat{t},\eM_{r}^\star(\hat{t}))\in\cF_{\UltMFunnel}}$
    and $\int_{\hat{t}}^{\hat{t}+T}\FunnelPenaltyFuncR(s,\eM_{r}^\star(s))\d{s}<\infty$
    because of the property~\ref{Item:FunnelPenaltyConvergence} 
    of the funnel penalty function $\FunnelPenaltyFuncR$, see~\Cref{Def:FunnelPenaltyFunction}.
    In particular, this implies $J^{\Psi}_T(u^{\star};\hat{t},\InitState)<\infty$.
    Furthermore, we have $\Norm{\eM_{r}^\star(s)}<\UltMFunnel(s)$ for all $s\in[\hat{t},\hat{t}+T]$
    due to the property~\ref{Item:FunnelPenaltyIntegral} of the function $\FunnelPenaltyFuncR$, see~\Cref{Def:FunnelPenaltyFunction}.
    This implies $x^\star(s)-\OpChi(y_{\rf})(s)\in\cD_{t}^\Psi$ for all $s\in[\hat{t},\hat{t}+T]$
    according to \Cref{Prop:OnlyLastFunnel}.
    Hence, $u^{\star} \in \Controls(\umax,\InitState)$.

    \noindent
    \emph{Step 6}:
    We show $J^{\Psi}_T(u^{\star};\hat{t},\InitState) = J^{\star}$. 
    According to Step~5, we have $\Norm{\eM_{r}^\star(s)}<\UltMFunnel(s)$ for all $s\in[\hat{t},\hat{t}+T]$.
    By the continuity of the involved functions and the compactness of the interval, 
    there exists $\eps>0$ such that $\Norm{\eM_{r}^{\star}(s)}\leq\UltMFunnel(s)-\eps$ for all $s\in[\hat{t},\hat{t}+T]$.
    Moreover, due to the uniform convergence of $\eM_r^k$ to $\eM_{r}^\star$, there exists $N\in\N$ such that 
    $\SNorm{\eM_{r}^k-\eM_{r}^{\star}}<\tfrac{\eps}{2}$ for $k\geq N$. Thus,
    \[
        \fa k\geq N\ \fa s\in[\hat{t},\hat{t}+T]:\quad
        \Norm{\eM_{r}^k(s)}\leq
        \Norm{\eM_{r}^k(s)-\eM_{r}^{\star}(s)}+\Norm{\vphantom{\eM_{r}^{k}}\eM_{r}^{\star}(s)}<\UltMFunnel(s)-\frac{\eps}{2}.
    \]
    Since the function $\FunnelPenaltyFuncR$ restricted to ${\cF_{\UltMFunnel}}$ is continuous, see~\Cref{Def:FunnelPenaltyFunction},
    the sequence $\rbl\FunnelPenaltyFuncR(\cdot,\eM_{r}^k(\cdot))^{i}\rbr$ therefore is uniformly bounded with $i=1,2$.
    Hence, the bounded convergence theorem gives that
    \[
        \FunnelPenaltyFuncR(\cdot,\eM_{r}^k(\cdot))^{i}\to\FunnelPenaltyFuncR(\cdot,\eM_{r}^\star(\cdot))^{i}
    \]
    strongly and, thus, also weakly in $L^2([\hat{t},\hat{t}+T],\R)$ for $i=1,2$.
    Since the $L^2$-norm is weakly lower semi-continuous and since 
    $J^{\Psi}_T(u_k;\hat{t},\InitState)\to J^{\star} = \inf_{u \in \Controls(\hat{t},\InitState)} J^{\Psi}_T(u;\hat{t},\InitState)$,
    the following holds.
    \begin{align*}
        J^{\Psi}_T(u^{\star};\hat{t},\InitState)
                &=\int_{\hat{t}}^{\hat{t}+T}\FunnelStageCost(s,\eM_{r}^{\star}(s))\d s
                = \LNorm{ \FunnelPenaltyFuncR(\cdot,\eM_{r}^{\star}(\cdot))^\frac12}^2
                 +\lambda_u\LNorm{u^{\star}}^2\\
                &\leq \liminf_{k\to\infty} \LNorm{ \FunnelPenaltyFuncR(\cdot,\eM_{r}^{k}(\cdot))^\frac12}^2+\liminf_{k\to\infty}
                     \lambda_u\LNorm{u_k}^2
                \leq \liminf_{k\to\infty} J^{\Psi}_T(u_k;\hat{t},\InitState) = J^{\star}.
    \end{align*}
    Therefore $J^{\Psi}_T(u^{\star};\hat{t},\InitState) = \min_{u\in \Controls(\umax,\InitState)}J^{\Psi}_T(u;\hat{t},\InitState)$.

\noindent
    \emph{Step 7}:
    We show that $J^{\Psi}_T(u^{\star};\hat{t},\InitState) 
    = \min_{\substack
                {
                    u\in L^{\infty}([\hat{t}, \hat{t}+T],\R^{m}),\\
                    \SNorm{u}  \leq \umax 
                }
            }
     J^{\Psi}_T(u;\hat{t},\InitState)$.\\ 
     Since $\Controls(\umax,\InitState) \neq \emptyset$ by assumption this follows from Remark~\ref{Rem:PropertiesJ}\,(ii) and completes the proof. 
\end{proof}

\begin{remark}
    \Cref{Th:SolutionExists} shows that there exists a solution $u^\star$ to the optimal control problem~\eqref{eq:OptimalControlProblem} and that 
    the application of this solution to the model~\eqref{eq:Model_r} solves the tracking problem from~\Cref{Sec:ControlObjective}, i.e.
    it ensures that the output tracking error ${\eMTrack(t)=\yM(t)-y_{\rf}(t)}$ evolves within the funnel $\cF_{\Funnel}$ given by $\Funnel$.
    In application, however, it is often not possible to compute the solution of an OCP.
    One has to utilise numerical approximations instead. 
    As a consequence of~\Cref{Prop:ErrorIsLPath}, every approximation $\tilde{u}$ of $u^\star$, for which the cost function $J^{\Psi}_T$ is finite,
    still guarantees that the tracking error evolves within the prescribed performance funnel.
\end{remark}

\section{The funnel MPC algorithm}\label{Sec:FunnelMPC}
We now want to summarise our findings in the funnel MPC~\Cref{Algo:FunnelMPC}.
With the definitions, concepts, and results  so far at hand, we will prove that it is initial and recursive feasible 
and that the application of this control scheme to the model~\eqref{eq:Model_r} solves the tracking problem laid out in~\Cref{Sec:ControlObjective},
i.e. it  ensures that the distance between the model's output~$\yM$ and a given reference signal~$y_{\rf}\in W^{r,\infty}(\Rp,\R^{m})$ 
evolves within the funnel~$\cF_{\Funnel}$ given by a function~$\Funnel\in\cG$.

\begin{algo}[Funnel MPC]\label{Algo:FunnelMPC}\ \\
    \textbf{Given:}\\[-4ex] %
    \begin{itemize}
        \item 
    Model~\eqref{eq:Model_r} with initial time $t_0\in\Rp$ and initial value ${\yM^0\in\cC^{r-1}([0,t_0],\R^m)}$,
    reference signal $y_{\rf}\in W^{r,\infty}(\Rp,\R^{m})$ , signal memory length~$\tau\geq0$,
    \item 
    a set of funnel boundary function $\Psi=(\Funnel_1,\ldots,\Funnel_r)\in\FunnelBoundaryFuncs$ with corresponding parameters~$k_i$ for $i=1,\ldots, r$,
    input saturation level $\umax\geq0$, funnel stage cost function~$\FunnelStageCost$, and 
    a $\tau$-initialisation strategy $\InitStrategy$ as in~\Cref{Def:InitialisationStrategy}.
    \end{itemize} 
    \textbf{Set} the time shift $\delta >0$, 
                 the prediction horizon $T\geq\delta$, index $k\coloneqq 0$, 
                 and $\xMh^0\coloneqq \OpChi(\yM^0)$.\\
    \textbf{Define} the time sequence~$(t_k)_{k\in\N_0} $ by $t_k \coloneqq  t_0+k\delta$.\\ 
    \textbf{Steps:}
    \begin{enumerate}[(a)]
    \item\label{agostep:FunnelMPCFirst}
        Select initial model state $\InitStateK_{k}\coloneqq \InitStrategy(\xMh^k)\in\InitValues(t_k)$ at current time $t_k$
        based on $\xMh^k$.
    \item
        Compute a solution $\uFMPCk\in L^\infty([t_k,t_k +T],\R^{m})$ of the Optimal
        Control Problem~(OCP)
    \begin{equation}\label{eq:FunnelMpcOCP}
        \mathop
                {\operatorname{minimise}}_{\substack
                {
                    u\in L^{\infty}([t_k,t_k+T],\R^{m}),\\
                    \SNorm{u}  \leq \umax 
                }
            }\      \int_{t_k}^{t_k + T}\FunnelStageCost(s,\eM_{r}(\xM(s;t_k,\InitStateK_{k},u)-\OpChi(y_{\rf})(s)),u(s))\d{s}.
    \end{equation}
    \item\label{agostep:FunnelMPCLast} Apply the control law
        \begin{equation}\label{eq:FMPC-fb}
            \mu:[t_k,t_{k+1})\times\InitValues(t_k)\to\R^m, \quad \mu(t,\InitStateK_{k}) =\uFMPCk(t)
        \end{equation}
        to model~\eqref{eq:Model_r} with initial time and data
        $(t_k,\InitStateK_{k})$ and obtain, on the interval
        ${I_{0}^{t_{k+1},\tau}\coloneqq [t_{k+1}-\tau,t_{k+1}]\cap[0,t_{k+1}]}$
        a measurement of the model's output and its derivatives
        ${\xMh^{k+1}\coloneqq
        \xM(\cdot;t_k,\InitStateK_{k},\uFMPCk)|_{I_{t_0}^{t_{k+1},\tau}}}$.
        Increment $k$ by $1$ and go to Step~\ref{agostep:FunnelMPCFirst}.
    \end{enumerate}
\end{algo}

Before proving the correct functioning of the funnel MPC~\Cref{Algo:FunnelMPC} a comment on the solution concept seems in order.
At every iteration of~\Cref{Algo:FunnelMPC} the model~\eqref{eq:Model_r} 
is re-initialised with new initial values $(t_k,\InitStateK_{k})$ and a solution $\xM$ on the interval $[t_k,t_k+T]$ 
while solving the OCP~\eqref{eq:FunnelMpcOCP}.
Note that this solution is in fact defined on the whole interval $[0,t_k+T]$
according to our understanding of a solution of the initial value problem, see~\Cref{Def:ModSolution}.
However, since the in~Step~\ref{agostep:FunnelMPCFirst} selected initial values do not necessarily coincide
with the solution of the initial value problem from the previous iteration, applying~\Cref{Algo:FunnelMPC} to the model~\eqref{eq:Model_r} 
does not result in a closed-loop system with a global solution in the classical sense. 
In the following~\Cref{Def:SolutionClosedLoop}, we therefore define a notion of a concatenated solution
which takes the re-initialisation into account and is, in a certain sense, a solution of all the considered initial value problems
on the subintervals~${[t_k,t_{k+1})}$.

\begin{definition}[Concatenated model solution]\label{Def:SolutionClosedLoop}
    Let  ${(\fM,\gM,\oTM)\in\cM^{m,r}_{t_0}}$ and consider the model \eqref{eq:Model_r}.
    Let $y_{\rf}\in W^{r,\infty}(\Rp,\R^{m})$, $\Psi=(\Funnel_1,\ldots,\Funnel_r)\in\FunnelBoundaryFuncs$, $t_0\in\Rp$, and $\delta>0$ be given. 
    Define the sequences $(t_k)_{k\in N_0}$ and $(\InitStateK_k)_{k\in N_0}$
    by $t_k\coloneqq t_0+k\delta$ and $\InitStateK_k\in \InitValues(t_k)$.
    Further, suppose $u\in L^\infty_{\loc}([t_0,\infty),\R^m)$ is a control such that initial value problem~\eqref{eq:Model_r}
    with initial data~$\InitStateK_k$ at time $t_k$ has a solution $\xM^k:[0,t_{k+1}]\to\R^{rm}$ in the sense of~\Cref{Def:ModSolution}
    for every $k\in\N$.
    We call the function $\xM: \Rp \to \R^{rm}$ that is piecewise defined by 
    \begin{equation*}
        \xM (t)=\begin{cases}
            \xM^0(t), &t<t_1,\\
            \xM^k(t), &t\in[t_k,t_{k+1})
        \end{cases}
    \end{equation*}    
    a \emph{concatenated solution} of the initial value problem~\eqref{eq:Model_r} with sequence of initial values $(t_k,\InitStateK_k)_{k\in\N_0}$.
    Its first $m-$dimensional component $x_{\mathrm{M},1}:\Rp \to \R^{m}$ is  denoted by~$\yM$.
\end{definition}

We will now prove the initial and recursive feasibility of the funnel MPC~\Cref{Algo:FunnelMPC} and 
that applying this algorithm to model~\eqref{eq:Model_r} results in a system that has a concatenated solution in the sense of~\Cref{Def:SolutionClosedLoop}.

\begin{theorem}\label{Thm:FMPC}
    Consider model~\eqref{eq:Model_r} with $(\fM,\gM,\oTM)\in\cM^{m,r}_{t_0}$ with initial value $\yM^0\in\cC^{r-1}([0,t_0],\R^m)$.
    Let $y_{\rf}\in W^{r,\infty}(\Rp,\R^{m})$ and $\Psi=(\Funnel_1,\ldots,\Funnel_r)\in\FunnelBoundaryFuncs$ be given.  
    Further, let $\tau\geq0$ be greater than or equal to the memory limit of operator~$\oTM$ and 
    choose a $\tau$-initialisation strategy
    ${\InitStrategy:\bigcup_{\hat{t}\geq t_0}\cR(I_0^{\hat{t},\tau},\R^{rm})\to
        \bigcup_{\hat{t}\geq t_0}\cR(I_0^{\hat{t},\tau},\R^{rm})\times L^\infty_{\loc}( I_{t_0}^{\hat{t},\tau},\R^q)
    }$ as in~\Cref{Def:InitialisationStrategy}. 
    Then, there exists $\umax\geq0$ such that the funnel MPC~\Cref{Algo:FunnelMPC} with $\delta>0$ and $T\ge\delta$ is initially and recursively feasible, i.e. 
    \begin{itemize}
        \item at every time instant $t_k \coloneqq  t_0+k\delta $ for $k\in\N_0$ the OCP~\eqref{eq:FunnelMpcOCP} has a solution $\uFMPCk\in L^\infty([t_k,t_k+T],\R^m)$, and
        \item the model~\eqref{eq:Model_r} with applied funnel MPC feedback~\eqref{eq:FMPC-fb} 
        has a concatenated solution $\xM : [0,\infty) \to \R^{rm}$ in the sense of \Cref{Def:SolutionClosedLoop}. 
    \end{itemize}
    The corresponding input is given by 
    \[
        \uFMPC(t) = \uFMPCk(t),
    \]
    for $t\in[t_k,t_{k+1})$ and $k\in\N_0$.
    Each global solution $\xM$ with corresponding output $\yM$ and input $u_{\mathrm{FMPC}}$ satisfies: 
    \begin{enumerate}[label = (\roman{enumi}), ref=(\roman{enumi})]
        \item\label{th:item:BoundedInput} the control input is bounded by $\umax$, i.e.
        \[
            \fa t \ge t_0 :\quad \Norm{u_{\mathrm{FMPC}}(t)}\leq\umax,
        \]
        \item\label{th:item:ErrorInFunnel}
        the tracking error between the model output and the reference evolves within prescribed boundaries, i.e.
        \[
            \fa t \ge t_0 :\quad \Norm{\yM(t)  - y_{\rf}(t)} < \Funnel_1(t) .
        \]
    \end{enumerate}
\end{theorem}

\begin{proof}
    \emph{Step 1}:
    According to~\Cref{Rem:InitialValuesNotEmpty}, the set $\InitValues(\hat{t})$ is non-empty for all $\hat{t}\geq t_0$.
    In particular, $\InitValues(t_0)\neq\emptyset$. 
    Using the notation $I_{0}^{t_0,\tau}\coloneqq [t_k-\tau,t_k]\cap[0,t_0]$, let $\xMh^0\coloneqq \OpChi(\yM^0)|_{I_{0}^{t_0,\tau}}$   
    be the measurement of the initial model output $\yM^0$ and its derivatives.
    Note that we identify $\cC^r([0,t_0],\R^m)$ with the vector space $\R^{rm}$ if $t_0=0$, see two cases $t_0>0$ and $t_0=0$ in~\eqref{eq:Model_r}.
    Since $\InitValues(t_0)\neq\emptyset$, it is possible to select the initial model state $\InitStateK_{0}=\InitStrategy(\xMh^0)$.
    
    \noindent 
    \emph{Step 2}:
    There exists $\umax\geq0$ such that, $\Controls(\umax,\InitState)\neq\emptyset$
    for all $\hat{t}\geq t_0$, $\InitState\in\InitValues(\hat{t})$, and $T>0$, according to~\Cref{Th:ExUmax}.
    Assume that we have given $\xMh^k\in\cR(I_0^{t_{k},\tau},\R^{rm})$ for $k\in\N_0$,
    where $I_{0}^{t_{k},\tau}\coloneqq [t_{k}-\tau,t_{k}]\cap[0,t_{k}]$.
    $\xMh^k$ is the model's output and its derivatives on the interval $I_{0}^{t_k,\tau}$ from the previous iteration 
    of the funnel MPC~\Cref{Algo:FunnelMPC}.  
    Since the set $\InitValues(t_k)$ is non-empty according to~\Cref{Rem:InitialValuesNotEmpty}, 
    there exists a model state $\InitStateK_{k}\coloneqq \InitStrategy(\xMh^k)\in\InitValues(t_k)$.
    \Cref{Th:SolutionExists} yields the existence of some function ${\uFMPCk\in\cU_{[t_k,t_k+T]}(\umax,\InitStateK_k)}$
    such that the functional $J^{\Psi}_T$ as in~\eqref{eq:DefCostFunctionJ} has a minimum, that~is 
    \[
        J^{\Psi}_T(\uFMPCk;t_k,\InitStateK_k)   =
        \mathop{\min}_{
            \substack
            {
                u\in L^{\infty}([t_k, t_k+T],\R^{m}),\\
                \SNorm{u}  \leq \umax 
            }
        }
            J^{\Psi}_T(u;t_k,\InitStateK_k). 
    \]
    Thus, $\uFMPCk$ is a solution of the OCP~\eqref{eq:FunnelMpcOCP}.
    If the control $\uFMPCk$ is applied to the model~\eqref{eq:Model_r} at initial time $t_k$ initial value $\InitStateK_k$
    in Step~\ref{eq:FMPC-fb}, then the initial value problem~\eqref{eq:Model_InitialValue} has a solution 
    $\xM^k:[0,t_{k+1}]\to\R^{rm}$ in the sense of~\Cref{Def:ModSolution} as a consequence of the definition of 
    $\cU_{[t_k,t_k+T]}(\umax,\InitStateK_k)$, see~\eqref{eq:Def-U}. 
    As $\xM^k$ is defined on the whole interval $[0,t_{k+1}]$, the function 
    $\xMh^{k+1}\coloneqq \xM(\cdot;t_k,\InitStateK_{k},\uFMPCk)|_{I_{t_0}^{t_{k+1},\tau}}\in\cR(I_0^{t_{k+1},\tau},\R^{rm})$
    is well-defined in Step~\ref{agostep:FunnelMPCLast} of~\Cref{Algo:FunnelMPC}.

    \noindent 
    \emph{Step 3}:
    The recursive application of Step~2 in combination with Step~1 of this proof yield 
    the existence of a sequence $(\InitStateK_{k})_{k\in\N_0}$ 
    of initial values $\InitStateK_{k}\in\InitValues(t_k)$,
    control signals $\uFMPCk\in\cU_{[t_k,t_k+T]}(\umax,\InitStateK_k)$ solving the OCP~\eqref{eq:FunnelMpcOCP},
    and corresponding solutions ${\xM^k:[0,t_{k+1}]\to\R^{rm}}$ 
    of the initial value problem~\eqref{eq:Model_InitialValue}
    in the sense of~\Cref{Def:ModSolution}.
    Hence, the funnel MPC~\Cref{Algo:FunnelMPC} is initially and recursively feasible.
    Define ${u_{\mathrm{FMPC}}\in L^\infty_{\loc}([t_0,\infty),\R^m)}$  by $\uFMPC(t) = \uFMPCk(t)$ for $k\in \N_0$ and 
    $\xM: \Rp \to \R^{rm}$ by 
    \begin{equation*}
        \xM (t)=\begin{cases}
            \xM^0(t), &t<t_1,\\
            \xM^k(t), &t\in[t_k,t_{k+1}).
        \end{cases}
    \end{equation*}     
    Then, $\xM$ is a concatenated solution of the initial value problem~\eqref{eq:Model_r}
    with sequence of initial values $(t_k,\InitStateK_k)_{k\in\N_0}$ in the sense of~\Cref{Def:SolutionClosedLoop}. 
    
    \noindent 
    \emph{Step 4}: As $\uFMPCk\in\cU_{[t_k,t_k+T]}(\umax,\InitStateK_k)$, we have $\SNorm{\uFMPCk}\leq\umax$ for all $k\in\N_0$.
    Thus, $\Norm{\uFMPC(t)}\leq \umax$ for all $t\geq t_0$. This shows \ref{th:item:BoundedInput}.
    By definition of the set $\cU_{[t_k,t_k+T]}(\umax,\InitStateK_k)$, 
    all error variables $\eM_{i}(\xM(t;t_k,\InitStateK_k,u^\star_k) - \OpChi(y_{\rf}))$ from \eqref{eq:ErrorVar} evolve within their respective 
    funnels $\cF_{\Funnel_{i}}$ given by $\Funnel_i$ for $i=1,\ldots, r$, i.e.
    \[
        \Norm{\eM_{i}(\xM(t;t_k,\InitStateK_k,u^\star_k) - \OpChi(y_{\rf}))}<\Funnel_{i}(t)
    \]
    for all $t\in[t_k,t_{k+1})$.
    The output $\yM$ is the first $m-$dimensional component of $\xM$, see~\Cref{Def:SolutionClosedLoop}.
    According to the definition of $\eM_{1}$ as in \eqref{eq:ErrorVar}, we thus have 
    \[
        \Norm{\yM(t)  - y_{\rf}(t)}=\Norm{\eM_{1}(\xM(t;t_k,\InitStateK_k,u^\star_k) - \OpChi(y_{\rf}))}< \Funnel_1(t)
    \]
    for all $t\in[t_k,t_{k+1})$ and $k\in\N_0$.
    This shows \ref{th:item:ErrorInFunnel} and completes the proof.
\end{proof}

\begin{remark}\label{Rem:FMPCAlternativeFirstStep}
    As the inclined reader has probably already noticed, 
    using the term \emph{recursive feasibility} with respect to the funnel MPC~\Cref{Algo:FunnelMPC} 
    is a bit of a stretch since applying the solution $\uFMPCk\in L^\infty([t_k,t_k +T],\R^{m})$
    of the optimal control problems~\eqref{eq:FunnelMpcOCP} 
    to model~\eqref{eq:Model_r} with initial data $(t_k,\InitStateK_{k})$ is an open-loop control problem for every $k\in\N_{0}$.
    These individual problems are only loosely coupled via the initialisation strategy~$\kappa$.
    It would be more precise to say that the funnel MPC~\Cref{Algo:FunnelMPC} solves an infinite sequence of open-loop problems.
    However, \Cref{Thm:RecursiveInitialValues} in combination with~\Cref{Th:SolutionExists} and~\Cref{Th:JFiniteCost}
    shows that the funnel MPC~\Cref{Algo:FunnelMPC} 
    ensures that the state of the model at the end of each iteration
    is  again a feasible initial value for the next iteration, i.e.
    \[
        (\xM(\cdot;t_k,\InitStateK_k,u)|_{[t_k-\tau,t_{k+1}]\cap[0,t_{k+1}]},
        \oT(\xM(\cdot;t_k,\InitStateK_k,u))|_{[t_k-\tau,t_{k+1}]\cap[t_0,\hat{t}+\delta]})
        \in\InitValues(t_{k+1}).
    \]
    If this model state is always selected in Step~\ref{agostep:FunnelMPCFirst} of~\Cref{Algo:FunnelMPC}, 
    then applying~\Cref{Algo:FunnelMPC} to the model~\eqref{eq:Model_r} results 
    in a closed-loop system with a global solution $\xM:[0,\infty)\to\R^{rm}$ of the initial value problem~\eqref{eq:Model_r}
    in a more classical sense, meaning the differential equation is fulfilled on the whole interval $[t_0,\infty)$.
    Using the notation ${I_{t_0}^{t_k,\tau}\coloneqq [t_k-\tau,t_k]\cap[t_0,t_k]}$
    and setting 
    $\InitStateK_0\coloneqq (\OpChi(\yM^0)|_{I_{0}^{t_0,\tau}},\oTM(\OpChi(\yM^0))|_{I_{t_0}^{t_0,\tau}})$,
    it is therefore possible to 
    replace Step~\ref{agostep:FunnelMPCFirst} of~\Cref{Algo:FunnelMPC}, for $k\geq 1$, by
    \begin{enumerate}[(a')]
    \item\label{agostep:FunnelMPCFirstAlternative}
    Obtain a measurement of the model state $\xM$ and~$\oTM(\xM)$ of~\eqref{eq:Model_r} on the interval 
    $I_{t_0}^{t_k,\tau}$ 
    and set $\InitStateK_k\coloneqq (\xM|_{I_{0}^{t_k,\tau}},\oT(\xM)|_{I_{t_0}^{t_k,\tau}})$.
    \end{enumerate}
    However, utilising an initialisation strategy $\InitStrategy$ in Step~\ref{agostep:FunnelMPCFirst} of~\Cref{Algo:FunnelMPC}
    opens up the possibility of directly applying the funnel MPC algorithm to a system which does not coincide with the model~\eqref{eq:Model_r}.
    In this case, the initial model state is selected based on measurement data from the system's output $y$ and its derivatives.
    \Cref{Algo:FunnelMPC} still remains feasible, meaning the optimal control problem~\eqref{eq:FunnelMpcOCP} has a solution and
    the model output $\yM$ is ensured to evolve within the funnel $\cF_{\Funnel}$. 
    While such guarantees can in such cases not be given for the actual system to be controlled,  
    the controller's performance might still be adequate if system and model only slightly diverge due to measurement errors and small disturbances.
    In~\Cref{Chapter:RobustFunnelMPC}, we will examine in more detail how the funnel MPC \Cref{Algo:FunnelMPC} can be adapted 
    in order to give guarantees on the tracking error for the actual system in the presence of a model-plant mismatch~$\eSTrack$ 
    as in~\eqref{eq:Intro:ModelPlantMismatch}.
\end{remark}

\begin{remark}
    \begin{enumerate}[(a)]
        \item
            The OCP~\eqref{eq:FunnelMpcOCP} has neither state nor terminal constraints.
            Nevertheless, application of the funnel MPC~\Cref{Algo:FunnelMPC} to the
            model~\eqref{eq:Model_r} ensures the existence of a global solution of the initial value problem
            in the sense of~\Cref{Def:SolutionClosedLoop} or the solution of the closed-loop system if 
            Step~\ref{agostep:FunnelMPCFirst} is replaced by Step~\ref{agostep:FunnelMPCFirstAlternative} from~\Cref{Rem:FMPCAlternativeFirstStep}.
            However, note that in neither case this solution is unique in general.
            One of the reasons is that the solution of the OCP~\eqref{eq:FunnelMpcOCP} found in
            each step may not be unique. The MPC algorithm has to select a particular optimal
            control. In particular, \Cref{Thm:FMPC} shows that the
            properties~\ref{th:item:BoundedInput} and~\ref{th:item:ErrorInFunnel} are independent
            of the particular choice made within the MPC algorithm, since they hold for every such solution.
            However, $\xM|_{[t_k,t_{k-1})}$ is uniquely determined by the choice of $\InitStateK_k$ and $u_k^\star$ 
            for every $k\in\N_{0}$ as~\Cref{Prop:SolutionUnique} shows.
        \item
            Funnel MPC is initially and recursively feasible for every choice of $T>0$. Usually, recursive
            feasibility for model predictive control can only be guaranteed when the prediction
            horizon is sufficiently long, see e.g.~\cite{boccia2014stability}, or when additional
            terminal constraints are added to the OCP, see e.g.~\cite{rawlings2017model}.
            For funnel MPC merely the boundary on the control input $\umax\geq0$ must be sufficiently large.
    \end{enumerate}
\end{remark}

\begin{remark}
While the primary funnel function $\Funnel$ is user-defined based on application-specific
tracking error constraints, the funnel MPC~\Cref{Algo:FunnelMPC}
introduces additional parameters -- notably the auxiliary funnel functions
$\Psi=(\Funnel_1,\ldots,\Funnel_r)\in\FunnelBoundaryFuncs$ and
their associated error gains $ k_i$ for $i=1,\ldots, r$ -- whose impact on the controller
performance warrants discussion.

The functions~$\Funnel_i$ for $i=1,\ldots, r-1$ do not directly influence the controller, 
as the optimal control problem \eqref{eq:FunnelMpcOCP} only optimises over a funnel
stage cost function~$\FunnelStageCost$ linked to the last auxiliary funnel
function~$\Funnel_r$.
These auxiliary functions $\Funnel_i$ are determined by the gains $k_i$ and the
constants $\FunDeriv$, $\FunDiam$ and $\gamma$ (determined by the function
$\Funnel$), as defined in \eqref{eq:cond-k_i},~\eqref{eq:psi_i}
or~\eqref{eq:cond-k_i_psi_i_simple}. Consequently, we focus on the effects of $k_i$ and $\Funnel_r$:
\begin{itemize}
    \item Tighter boundary $\Funnel_r$ improves track precision but restricts the optimiser's flexibility 
    to accommodate secondary objectives, e.g. minimising the control effort.
    \item Larger $k_i$ values intensify penalisation of error variables~$\eM_i$ (see \eqref{eq:ErrorVar}),  
    enhancing accuracy at the cost of higher control inputs.
    \item The minimum bound $\umax\geq0$ for admissible control inputs depends on $k_i$ and $\sup\Funnel_i$, as shown in
    \Cref{Th:ExUmax}. However, these derived bounds are typically rather conservative; refining them requires a problem-specific analysis.
\end{itemize}
Moreover, as noted in~\Cref{Sec:HighRelativeDegree}, the construction of the
parameters $k_i$ and $\Funnel_r$ (via~\eqref{eq:cond-k_i} and~\eqref{eq:psi_i}) can
be simplified using the design~\eqref{eq:cond-k_i_psi_i_simple}. 
This simplification replaces the time-varying funnel penalties with a constant
cost funnel penalty function  in the MPC stage cost $\FunnelStageCost$, reducing
computational complexity. However, the initial model trajectory $\yM^0$ for the
model \eqref{eq:Model_r} must be freely selectable to satisfy the constraints
imposed by this simplified design.
\end{remark}

\begin{remark}\label{Rem:MPCQuadraticWithConstraints}
    The proof of the funnel MPC’s recursive feasibility hinges primarily on
    \Cref{Th:ExUmax}. Crucially, this result does not depend on the use of
    \emph{funnel penalty functions},  implying that recursive feasibility can
    also be guaranteed for the MPC scheme in \Cref{Algo:MPC} with
    alternative stage cost functions $\ell:\Rp\times\R^m\times\R^{m}$
    -- such as the classical quadratic cost function in \eqref{eq:stageCostClassicalMPC} --
    provided the optimal control problem \eqref{eq:MpcOCP} incorporates 
    the additional constraint 
    \[
       \fa s\in [t_k,t_k+T]:\quad\Norm{\eM_{r}(\xM(s;t_k,\InitStateK_{k},u)-\OpChi(y_{\rf})(s))}< \Funnel_{r}(s).
    \]
    However, without funnel penalty functions, alternative methods are required
    to ensure the optimal control problem always admits a solution, as the
    proof of \Cref{Th:SolutionExists} relies explicitly on their use.
    Assuming this is achievable, \Cref{Algo:MPC} can also fulfil the
    control objective outlined in \Cref{Sec:ControlObjective}.
    While the funnel MPC (Algorithm~\Cref{Algo:FunnelMPC}) is hypothesised to
    offer superior performance -- due to its cost function dynamically penalising
    proximity to the funnel boundary -- a numerical case study~\cite{Oppeneiger24}
    found comparable results between the two approaches. A comprehensive
    comparative analysis, however, remains an open research question.
\end{remark}

\section{Simulation}\label{Sec:FMPC:Simulation}
This section illustrates the application of the funnel MPC algorithm
(\Cref{Algo:FunnelMPC}) using two numerical examples. In this section, we do
not distinguish between the actual system and its model. This distinction will
be revisited in subsequent chapters. Consequently, we omit the subscript
$\mathrm{M}$ used to denote model equations in this section.

The \textsc{Matlab} source code for the simulations performed in this thesis
can be found on \textsc{GitHub} under the link \url{https://github.com/ddennstaedt/FMPC_Simulation}.

\toclesslab\subsection{Exothermic chemical reaction}{Sec:FMPC:Sim:Reactor}
To demonstrate the funnel MPC~\Cref{Algo:FunnelMPC}, we consider a model of a
chemical reactor where an exothermic reaction Substance-1 $\to$ Substance-2 takes
place.
This example was also used in~\cite{IlchTren04} to study funnel control with
input saturation and in~\cite{LibeTren10} to demonstrate the feasibility of the
bang-bang funnel controller. According to \cite{VielJado97}, this type of reactor can
be modelled by the following system of equations of order one:
\begin{equation}\label{eq:ExampleExothermicReaction}
    \begin{aligned}
        \dot{x}_1(t) &= c_1\, p(x_1(t), x_2(t), y(t)) +d(x_1^{\mathrm{in}}-x_1(t)),\\
        \dot{x}_2(t) &= c_2\, p(x_1(t), x_2(t), y(t)) +d(x_2^{\mathrm{in}}-x_2(t)),\\
        \dot{y}(t)   &= b\,   p(x_1(t), x_2(t), y(t)) -q\,y(t) + u(t),
    \end{aligned}
\end{equation}
where $x_1$ is the concentration of the reactant Substance-1, $x_2$ the
concentration of the product Substance-2 and $y$ describes the reactor
temperature; $u$ is the feed temperature/coolant control input.
Further, the constant $b > 0$ describes the exothermicity of the reaction, $d>0$
is associated with the dilution rate and $q > 0$ is a constant consisting of the
combination of the dilution rate and the heat transfer rate.
Further, $c_1<0$ and $c_2 \in \R$ are the stoichiometric coefficients and $p : \Rp
\times \Rp \times \Rp \to \Rp$ is a locally Lipschitz continuous function with
$p(0,0,t)=0$ for all $t>0$ that models the reaction heat.
As in~\cite{IlchTren04}, we consider for the function $p$ the Arrhenius law
\begin{equation}\label{eq:Ex:ArrheniusLaw}
    p(x_1,x_2,y) = k_0 e^{- \frac{k_1}{y}} x_1,
\end{equation}
where $k_0, k_1$ are positive parameters. Since $c_1<0$, it is easy
to see that the subsystem
\begin{align*}
    \dot{x}_1(t) &= c_1 p(x_1(t), x_2(t), y(t)) +d(x_1^{\text{in}}-x_1(t)),\\
    \dot{x}_2(t) &= c_2 p(x_1(t), x_2(t), y(t)) +d(x_2^{\text{in}}-x_2(t)),
\end{align*}
satisfies the~BIBS condition~\eqref{eq:BIBO-ID} from \Cref{Ex:ControlAffineMod}, 
when~$y$ is restricted to the set $\setdef{y\in W^{1,\infty}(\Rp,\R)}{\forall\, t\ge 0:\ y(t)>0}$.
We like to emphasise that the control must thus guarantee that $y$ is always positive, 
which is also from a practical point of view a reasonable objective.
The control objective is to steer the reactor's temperature to a certain given reference
value $y_{\rf}(t)$ within boundaries given by a function~$\Funnel(t)$. The
reactor's temperature should follow a given heating profile specified as
\begin{equation}\label{eq:Ex:Reactor:RefTemp}
    y_{\rf}(t) = 
    \begin{cases}
        y_{\rf,\mathrm{start}} + \frac{y_{\rf ,\mathrm{final}} - y_{\rf,\mathrm{start}}}{t_{\mathrm{final}}} t, & t \in [0,t_{\mathrm{final}}), \\
        y_{\rf ,\mathrm{final}}, & t \ge t_{\mathrm{final}} . 
    \end{cases}
\end{equation}
Note that this heating profile has a kink at~$t = t_{\mathrm{final}}$.
Starting at $y_{\rf ,\mathrm{start}} = 270 \, \mathrm{K}$, the reactor is heated up to
$y_{\rf,\mathrm{final}} = 337.1 \, \mathrm{K}$ within the prescribed time~$[0,t_{\mathrm{final}}]$, here we choose~$t_{\mathrm{final}} = 2$.
The maximal control value is limited to $\umax=600$.
During the heating phase, the tolerated temperature deviation from the heating
profile decreases from~$\pm 24 \, \mathrm{K}$ to~$\pm 4.4 \, \mathrm{K}$ (time-varying
output constraints). After reaching the desired level, the temperature in the
reactor is kept constant with deviation of no more than $\pm 4.4 \, \mathrm{K}$ after
four units of time after beginning of the heating process.
We therefore choose the funnel function $\Funnel\in\cG$ given by 
\[
    \Funnel(t)\coloneqq20\me^{-2t}+4.
\]
To achieve the control objective with funnel MPC~\Cref{Algo:FunnelMPC}, 
we use the strict funnel stage cost function 
$\ell_{\Funnel}:\Rp\times\R\times\R\to\R\cup\{\infty\}$ given by 
\begin{equation}\label{eq:Ex:ReactorFMPCCostFunction}
\begin{aligned}
    \ell_{\Funnel}(t,y,u) =
    \begin{dcases}
        \frac {\Norm{y-y_{\rf}(t)}}{\Funnel(t)^2 - \Norm{y-y_{\rf}(t)}^2} + \lambda_u \Norm{u - 360}^2,
            & \Norm{y-y_{\rf}(t)} \neq \Funnel(t)\\
        \infty,&\text{else},
    \end{dcases}
\end{aligned}
\end{equation}
with design parameter ${\lambda_u\in\Rp}$.
This is a slightly modified variant of stage cost function from~\eqref{eq:stageCostFunnelMPC}.
As in~\cite{VielJado97,IlchTren04}, the initial data is chosen as $[x_1^0,
x_2^0, y^0] = [0.02,0.9,270]$ describing the initial concentration of the two substances 
and the initial reactor temperature.
The parameters of the system are
\begin{equation}\label{eq:Ex:Reactor:Sys:Params}
\begin{aligned}
     c_1&=-1,&
     k_0&=\me^{25},&
     x_{1}^{\mathrm{in}}&=1,&
     d&= 1.1,\\
     c_2&=1,&
     k_1&=8700,&
     x_{2}^{\mathrm{in}}&=0,&
     q&=1.25,&
     b&=209.2.
\end{aligned}
\end{equation}

\begin{figure}[h]
    \begin{subfigure}[b]{0.49\textwidth}
        \centering
        \includegraphics[width=.95\linewidth]{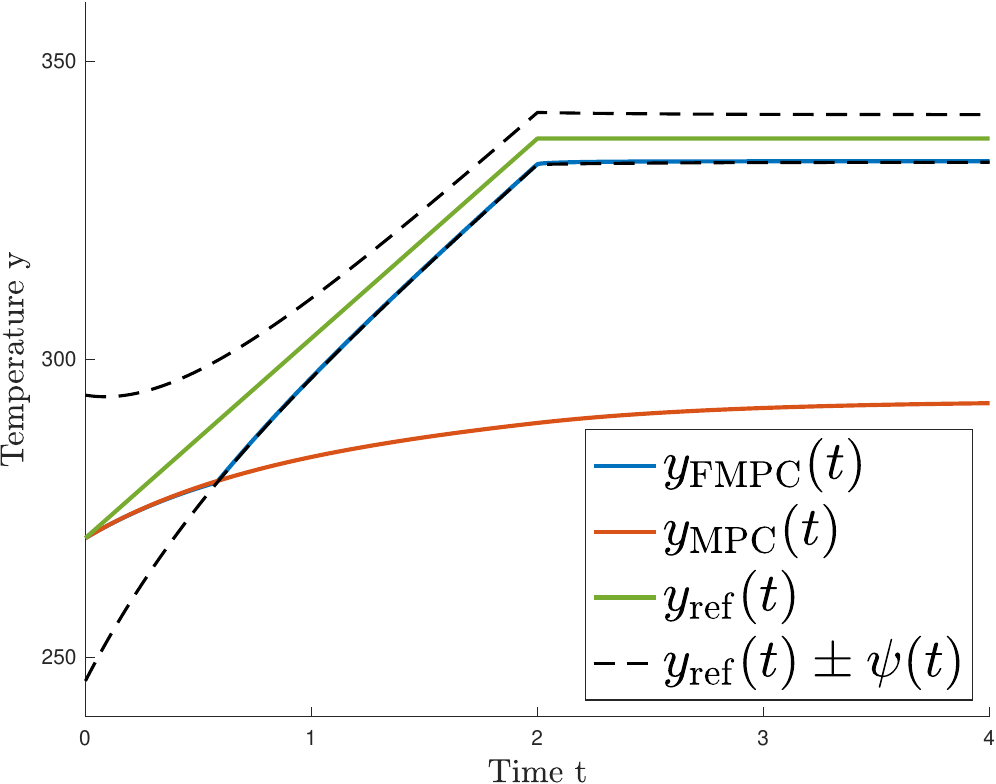}
        \caption{Outputs and reference, with boundary~$\Funnel$.}
        \label{fig:sim:fmpc_vs_mpc_short:output}
    \end{subfigure}
      \begin{subfigure}[b]{0.49\textwidth}
        \centering
        \includegraphics[width=\linewidth]{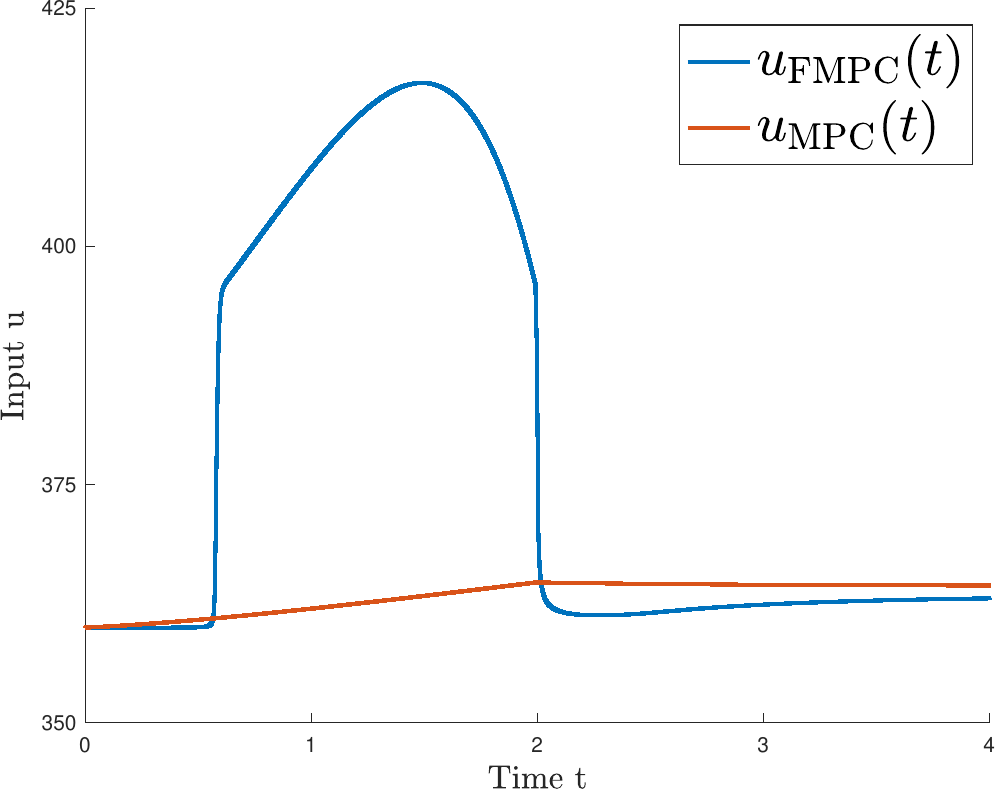}
        \caption{Control inputs.}
        \label{fig:sim:fmpc_vs_mpc_short:input}
    \end{subfigure}
    \caption{Simulation of system~\eqref{eq:ExampleExothermicReaction} under the
    control generated by the funnel MPC~\Cref{Algo:FunnelMPC} and the
    MPC~\Cref{Algo:MPC} with stage cost~\eqref{eq:Ex:QuadraticStageCost} with
    parameters $T= 10^{-2}$, $\delta=5\cdot 10^{-4}$, and
    $\lambda_u=0.1$.}
    \label{fig:sim:fmpc_vs_mpc_short}
\end{figure}

To demonstrate that the funnel MPC~\Cref{Algo:FunnelMPC} is initially and recursively feasible
even for a short prediction horizon, we choose the time shift $\delta=5\cdot 10^{-4}$ and ${T= 20\cdot\delta= 10^{-2}}$.
Due to discretisation, only step functions with constant
step length~$\SampleTime\coloneqq\delta$ are considered\footnote{
By a step function on an interval~$[a,b]$ with constant step length~$\SampleTime>0$, we mean a
mapping ${f:[a,b]\to\R}$ which is constant on every interval $[a+k\SampleTime,a+(k+1)\SampleTime)\cap[a,b]$ for
$k=0,\ldots,\lceil\frac{b-a}{\SampleTime}\rceil-1$, see also~\Cref{Def:PartitionAndStepFunction}.
} for the optimal control problem~\eqref{eq:FunnelMpcOCP} of the funnel MPC~\Cref{Algo:FunnelMPC}.
As perfect system knowledge is assumed here, the model is initialised, at every
iteration of the algorithm, with the model's state from the previous iteration.
We compare this control approach with the MPC~\Cref{Algo:MPC} using a standard quadratic cost function 
\begin{equation}\label{eq:Ex:QuadraticStageCost}
    \ell(t,y,u)= \Norm{y-y_{\rf}(t)}^2+\lambda_u \Norm{u-360}^2
\end{equation}
as in \eqref{eq:stageCostClassicalMPC}.
For both control schemes, the parameter $\lambda_{u}$ is chosen as $\lambda_{u}=0.1$.
The simulations are performed  with \textsc{Matlab} and the
toolkit \textsc{CasADi}\footnote{http://casadi.org}~\cite{Andersson2019}
over the time interval $[0,4]$ and depicted in~\Cref{fig:sim:fmpc_vs_mpc_short}.
\Cref{fig:sim:fmpc_vs_mpc_short:output} shows that the output of the system
evolves within the funnel boundaries when the control signal is generated by the
funnel MPC~\Cref{Algo:FunnelMPC} (labelled with $y_{\mathrm{FMPC}}$).
The standard MPC~\Cref{Algo:MPC} with the quadratic stage cost
function~\eqref{eq:Ex:QuadraticStageCost} does however not achieve the control objective.
The corresponding system output (labelled with $y_{\mathrm{MPC}}$) evolves outside of 
the prescribed boundaries. 
This observation is not surprising since no information about the
funnel~$\cF_{\Funnel}$ is included in stage cost
function~\eqref{eq:Ex:QuadraticStageCost} of the \Cref{Algo:MPC}. 
The incorporation of output constraints of the form
\begin{equation}\label{eq:Ex:ChemicalReactor:Constraints}
    \fa t\in[\hat{t},\hat{t}+T]:\quad \Norm{\yM(t)-y_{\rf}(t)} <\Funnel(t) 
\end{equation}
as in \eqref{eq:ConstraintClassicalMPC} in the corresponding
OCP~\eqref{eq:MpcOCP} is necessary in order to ensure that MPC with stage cost
\eqref{eq:Ex:QuadraticStageCost} is feasible with the selected prediction
horizon $T$, time shift $\delta$,  and design parameter $\lambda_u$.
In this case, the standard MPC scheme also achieves the control objective in accordance to
\Cref{Rem:MPCQuadraticWithConstraints}.

\begin{figure}[h]
    \begin{subfigure}[b]{0.49\textwidth}
        \centering
        \includegraphics[width=\linewidth]{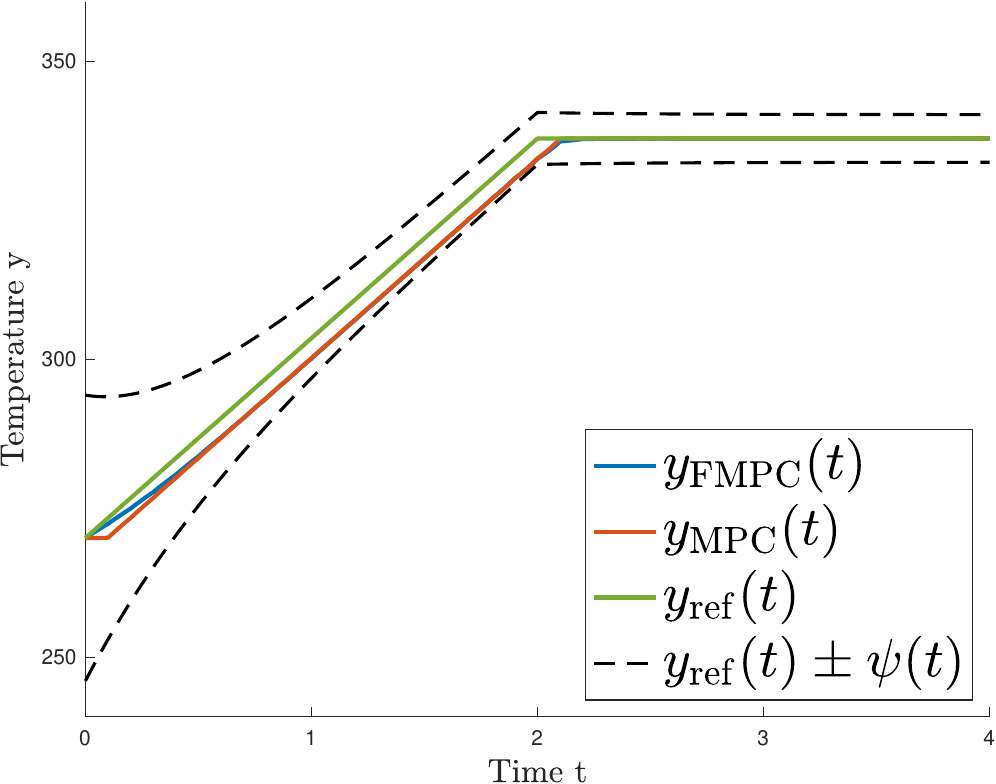}
        \caption{Outputs and reference, with boundary~$\Funnel$.}
        \label{fig:sim:fmpc_vs_mpc_long:output}
    \end{subfigure}
      \begin{subfigure}[b]{0.49\textwidth}
        \centering
        \includegraphics[width=\linewidth]{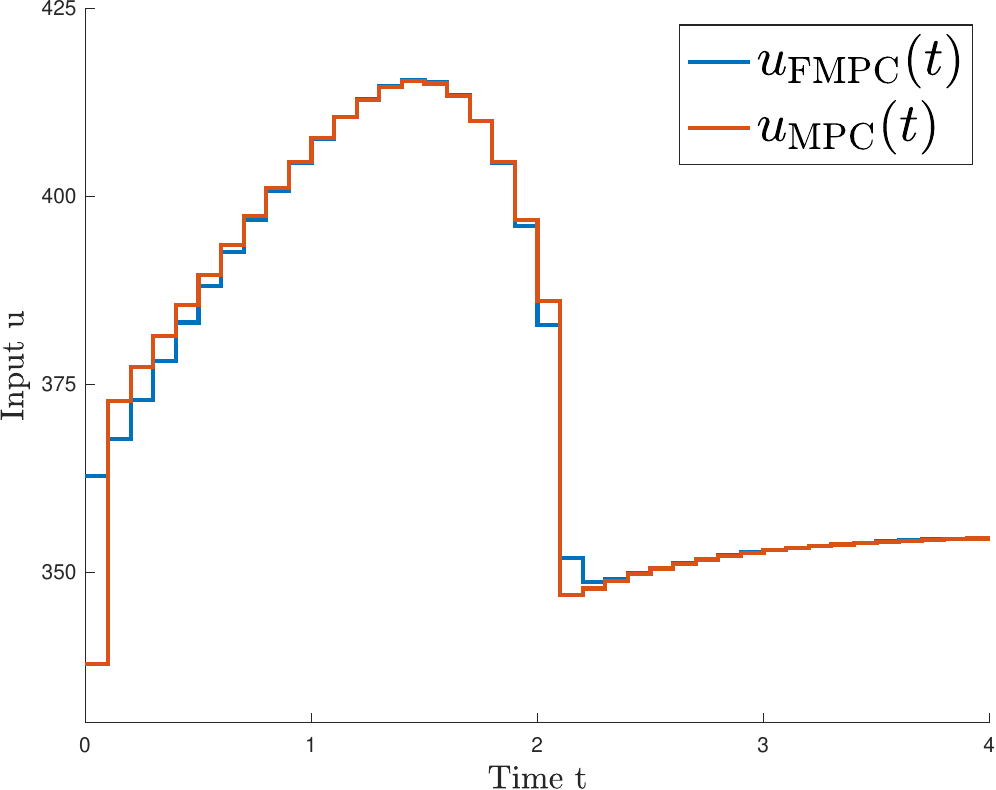}
        \caption{Control inputs.}
        \label{fig:sim:fmpc_vs_mpc_long:input}
    \end{subfigure}
    \caption{Simulation of system~\eqref{eq:ExampleExothermicReaction} under the
    control generated by the funnel MPC~\Cref{Algo:FunnelMPC} and the
    MPC~\Cref{Algo:MPC} with stage cost~\eqref{eq:Ex:QuadraticStageCost} with
    parameters $T= 1$, $\delta=0.1$, and
    $\lambda_u=10^{-4}$.}
    \label{fig:sim:fmpc_vs_mpc_long}
\end{figure}

Alternatively, if an appropriately long prediction horizon $T=1$ is chosen 
and the penalisation of the control signal is reduced by choosing~$\lambda_u=10^{-4}$,
then standard MPC~\Cref{Algo:MPC} is also able to achieve the control objective by chance.
Additionally, the time shift  is increased to $\delta =0.1$.
The performance of both control schemes with the longer horizon and time
shift and the adapted design parameter $\lambda_u=10$ is depicted
in~\Cref{fig:sim:fmpc_vs_mpc_long}.
While \Cref{fig:sim:fmpc_vs_mpc_long:output} shows the output of the system
evolving within the funnel boundaries under the both control schemes,
\Cref{fig:sim:fmpc_vs_mpc_long:input} shows the corresponding input signals.
It is evident that both control techniques generate very similar control signals
and achieve the control objective. 
However, while it is, to a certain extend, incidental that the tracking error of the
system under the control generated by the MPC~\Cref{Algo:MPC} evolves within 
$\cF_{\Funnel}$, the funnel stage cost function \eqref{eq:Ex:ReactorFMPCCostFunction}
provably ensures the adherence of the system to the funnel boundaries for 
the funnel MPC~\Cref{Algo:FunnelMPC}, see~\Cref{Th:JFiniteCost}. 
However, note that the discontinuous funnel penalties function can lead to  
compatibility issues with standard optimisation frameworks causing the utilised numerical 
solvers to fail. In this case, the incorporation of additional output constraints 
like~\eqref{eq:Ex:ChemicalReactor:Constraints} in the funnel MPC~\Cref{Algo:FunnelMPC} 
can mitigate these issues although they are, from a theoretical point of view, redundant.

In the following, we compare the funnel MPC~\Cref{Algo:FunnelMPC} 
to the funnel controller which was the inspiration for the development and 
usage of funnel penalty functions of the form~\eqref{eq:stageCostFunnelMPC}.
The original funnel controller proposed in~\cite{IlchRyan02b} takes the form
\begin{equation}\label{eq:Ex:ExothermicReaction:FC}
    \begin{aligned}
        \uFC(t)&=-\frac{1}{\Funnel(t)^2-\Norm{e(t)}^2}e(t).
    \end{aligned}
\end{equation}
For the comparison of the two controllers , we choose, as before, the strict
funnel stage cost~$\ell_{\Funnel}$ as in~\eqref{eq:Ex:ReactorFMPCCostFunction}
and the parameters $T= 1$, $\delta=0.1$, and $\lambda_u=10^{-4}$ for the MPC scheme
and restrict the set of control functions considered in the OCP~\eqref{eq:FunnelMpcOCP}
to step functions with constant step length $\SampleTime\coloneqq\delta=0.1$. 
To numerically compute the solution of the closed-loop system under the both control laws, 
the explicit four stage Runge-Kutta method (RK4) with a constant step size~$h>0$ is used.
\begin{figure}[h]
    \begin{subfigure}[b]{0.49\textwidth}
        \centering
        \includegraphics[width=\linewidth]{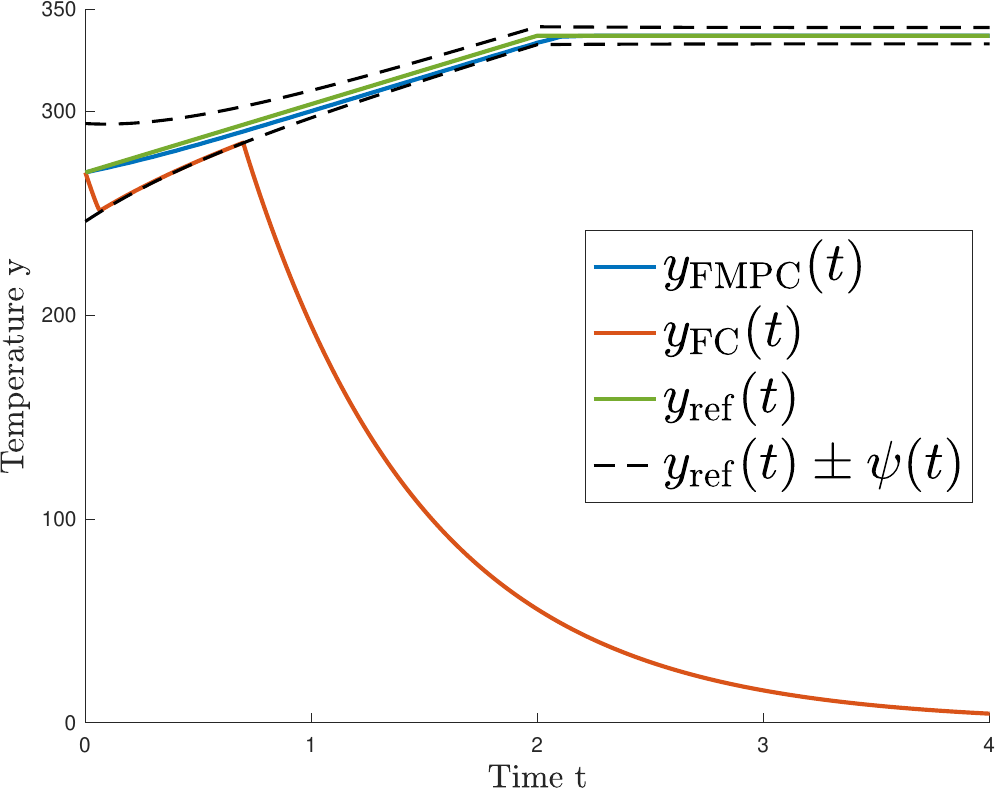}
        \caption{Outputs and reference, with boundary~$\Funnel$.}
        \label{fig:sim:fmpc_vs_fc_fail:output}
    \end{subfigure}
      \begin{subfigure}[b]{0.49\textwidth}
        \centering
        \includegraphics[width=\linewidth]{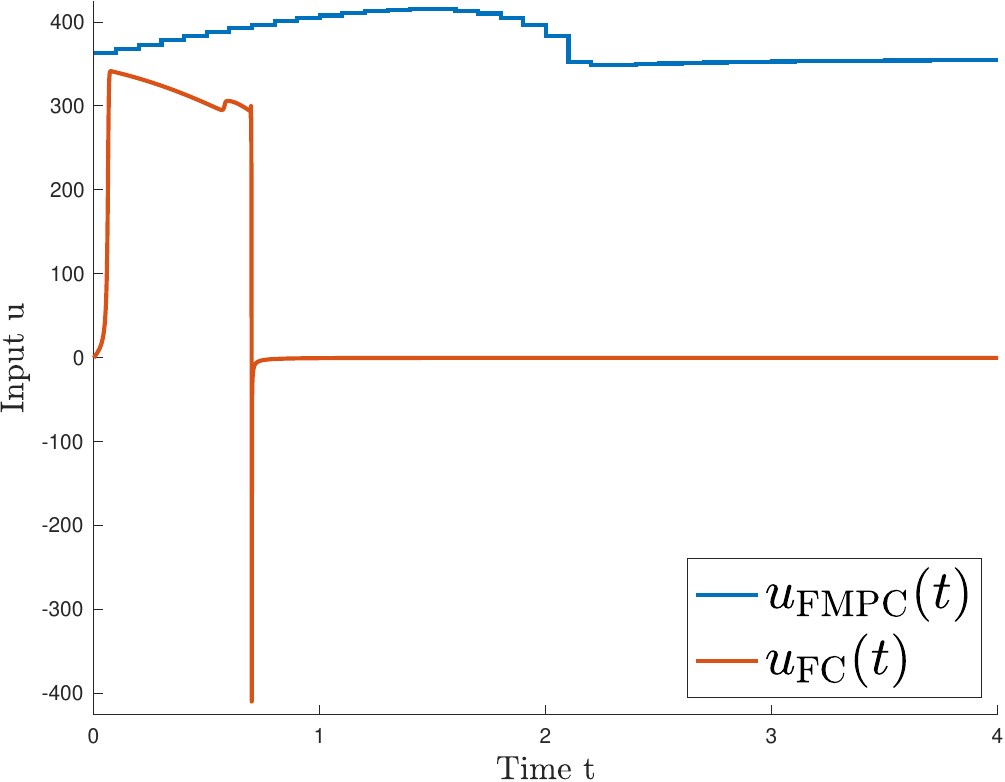}
        \caption{Control inputs.}
        \label{fig:sim:fmpc_vs_fc_fail:input}
    \end{subfigure}
    \caption{Simulation of system~\eqref{eq:ExampleExothermicReaction} under the
    control generated by the funnel MPC~\Cref{Algo:FunnelMPC} 
    with parameters $T= 1$, $\delta=0.1$, and $\lambda_u=10^{-4}$
    and the funnel control law~\eqref{eq:Ex:ExothermicReaction:FC} with a constant step size~$h=10^{-3}$.
    } 
    \label{fig:sim:fmpc_vs_fc_fail}
\end{figure}

\Cref{fig:sim:fmpc_vs_fc_fail} depicts the performance both of the funnel
controller~\eqref{eq:Ex:ExothermicReaction:FC} and the funnel
MPC~\Cref{Algo:FunnelMPC} when a constant step size~$h=10^{-3}$ is
used in the RK4 method. 
\Cref{fig:sim:fmpc_vs_fc_fail:output} shows that, while initially both control 
schemes are feasible, the system's output when controlled by the funnel controller 
(labelled with~$y_{\mathrm{FC}}$) breaches the funnel boundary at $t\approx0.7$
and  evolves from then onward outside the prescribed boundaries. 
The controller reacts at this time instant with a large peak in its control signal, see~\Cref{fig:sim:fmpc_vs_fc_fail:input}.
It is however not able to achieve the control objective on the entire considered time interval.
Although the funnel MPC~\Cref{Algo:FunnelMPC} is restricted to step functions as
its control signals a relatively wide step length of $\SampleTime=0.1$ and
therefore adapts its control signal significantly less often than the funnel
controller, the MPC algorithm is feasible and the system output
$y_{\mathrm{FMPC}}$ evolves within the performance funnel. Funnel MPC actually
still achieves the control objective if an even larger step size
of~$h=10^{-2}$ is used to solve the ordinary differential equation.

To ensure that the simulation of the funnel
controller~\eqref{eq:Ex:ExothermicReaction:FC} also achieves the control
objective the usage of a smaller step size in the RK4 method is required.
\begin{figure}[h]
    \begin{subfigure}[b]{0.49\textwidth}
        \centering
        \includegraphics[width=\linewidth]{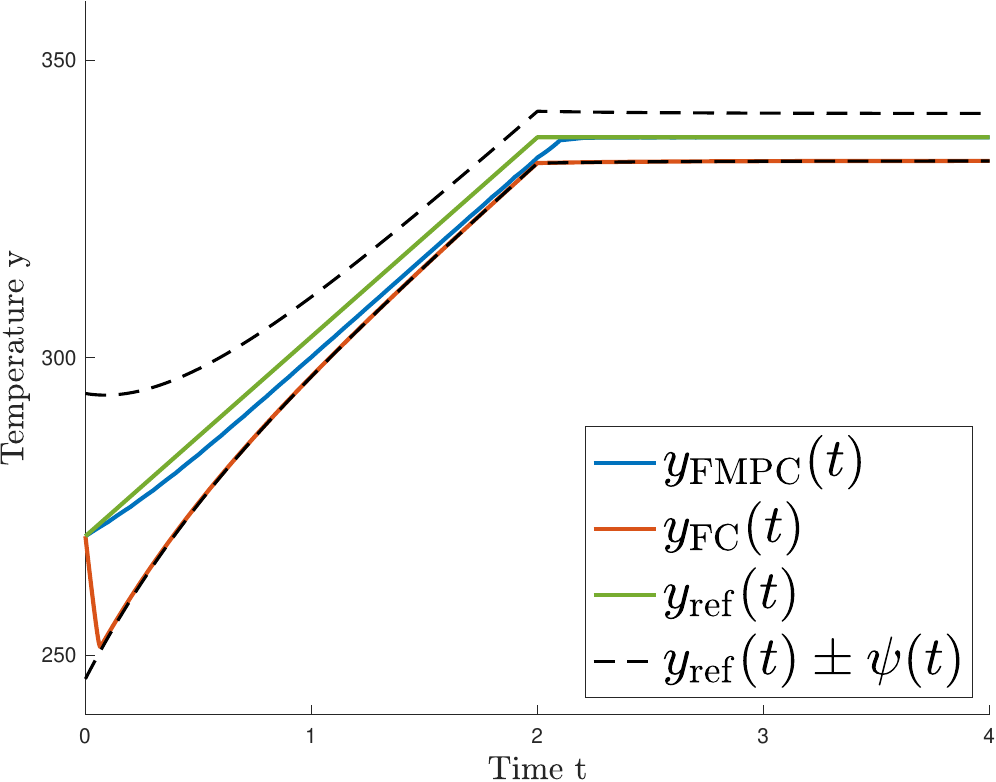}
        \caption{Outputs and reference, with boundary~$\Funnel$.}
        \label{fig:sim:fmpc_vs_fc:output}
    \end{subfigure}
      \begin{subfigure}[b]{0.49\textwidth}
        \centering
        \includegraphics[width=\linewidth]{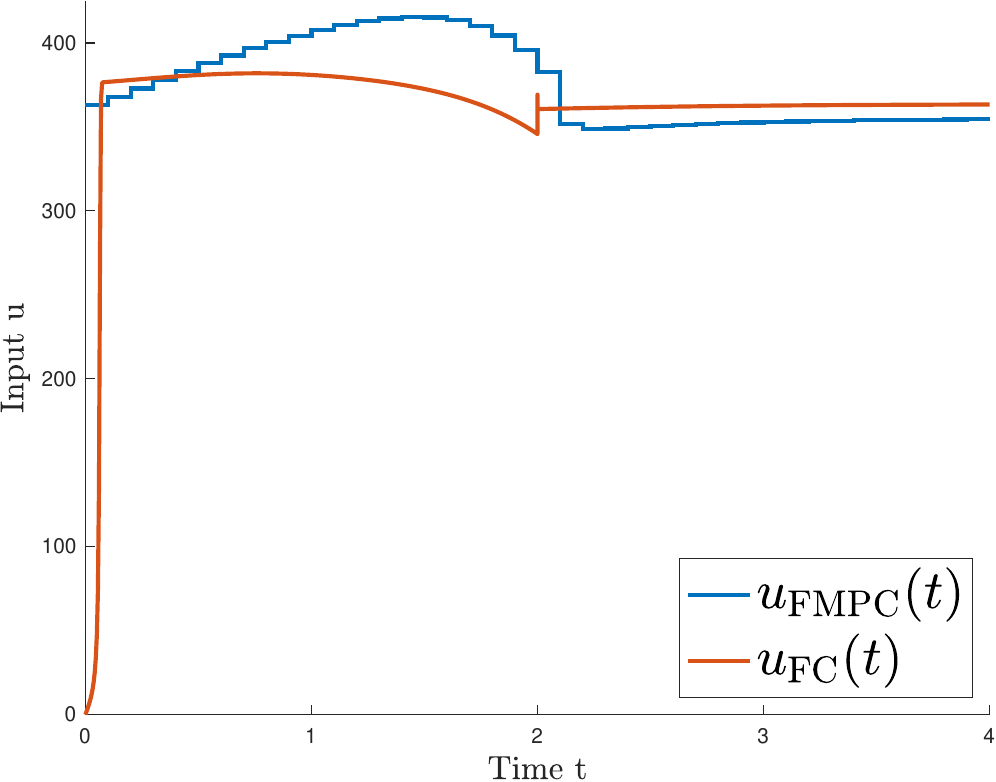}
        \caption{Control inputs.}
        \label{fig:sim:fmpc_vs_fc:input}
    \end{subfigure}
    \caption{Simulation of system~\eqref{eq:ExampleExothermicReaction} under the
    control generated by the funnel MPC~\Cref{Algo:FunnelMPC} 
    with parameters $T= 1$, $\delta=0.1$, and $\lambda_u=10^{-4}$
    and the funnel control law~\eqref{eq:Ex:ExothermicReaction:FC} with a constant step size~$h=10^{-4}$.
    } 
    \label{fig:sim:fmpc_vs_fc}
\end{figure}
The system output and the corresponding control signals are depicted
in~\Cref{fig:sim:fmpc_vs_fc} when a constant step size~$h=10^{-4}$ is
used for the simulation. \Cref{fig:sim:fmpc_vs_fc:output} shows the system
output under the control of the two approaches evolving within the funnel
boundaries and \Cref{fig:sim:fmpc_vs_fc:input} depicts the corresponding input
signals. It is evident that both control techniques are feasible and achieve the
control objective in this case. 
As the initial tracking error is zero, the funnel
controller~\eqref{eq:Ex:ExothermicReaction:FC} does not act in the beginning.
It's input signal $\uFC$ is zero and then increases when the system output
$y_{\mathrm{FC}}$ approaches the funnel boundary. 
Afterwards, $y_{\mathrm{FC}}$ evolves close to the funnel boundary over 
the entire time interval.
The system when controlled by the funnel MPC~\Cref{Algo:FunnelMPC}
exhibits a more accurate tracking. 
The system's output $y_{\mathrm{FMPC}}$ evolves closer to the reference
signal~$y_{\rf}$ than $y_{\mathrm{FC}}$.
Thanks to its predictive capabilities, the funnel MPC~\Cref{Algo:FunnelMPC}
applies already at the beginning a larger control signal $\uFMPC$ and does not
wait until the system's output is close to the boundary until it reacts.
After the system reached the desired temperature $y_{\rf,\mathrm{final}}$ 
at $t_{\mathrm{final}} = 2$, the output $y_{\mathrm{FC}}$ tracks 
the reference signal~$y_{\rf}$ almost perfectly. 
The system output $y_{\mathrm{FC}}$ under the control of $\uFC$ 
has a constant offset to the reference.
It is worth noting that funnel MPC~\Cref{Algo:FunnelMPC} achieves this better
performance while applying less control input than the funnel
controller~\eqref{eq:Ex:ExothermicReaction:FC}.

\toclesslab\subsection{Mass-on-car system}
\begin{figure}[h]
    \begin{center}
    \includegraphics[trim=2cm 4cm 5cm 15cm,clip=true,width=6.5cm]{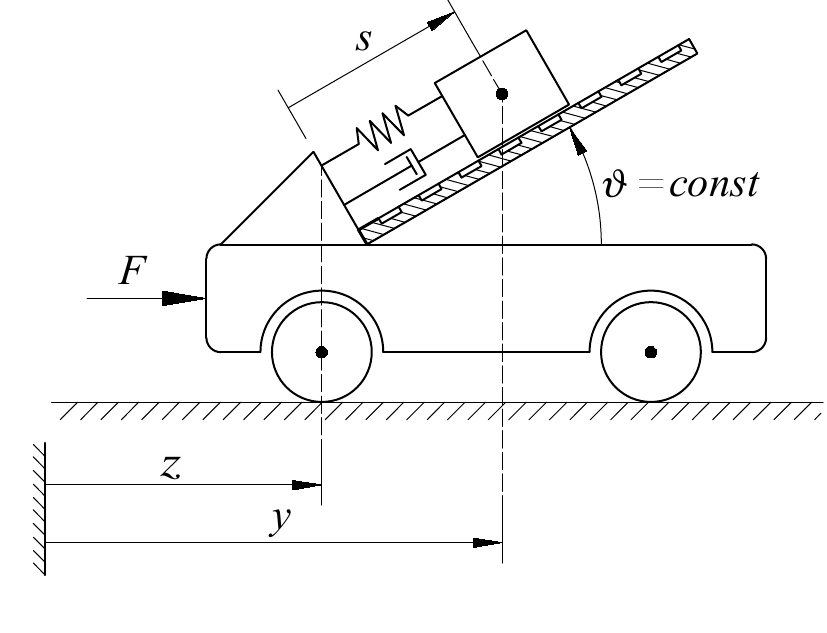}
    \end{center}
    \caption{Mass-on-car system. The figure is based on the respective figures in~\cite{BergIlch21}, and~\cite{SeifBlaj13}.}
    \label{Mass.on.car}
\end{figure}

{Sec:FMPC:Sim:MassOnCar}
For purposes of illustration that funnel MPC~\Cref{Algo:FunnelMPC} can also
successfully applied to systems with fixed higher relative degree, i.e. $r>1$,
we consider the example of a mass-on-car system from~\cite{SeifBlaj13}.
This example was also examined in~\cite{BergIlch21}
and~\cite{berger2019learningbased} to compare different versions of funnel control.
On a car with mass~$m_1$, to which a force~$F=u$ can be applied, a ramp is
mounted on which a second mass~$m_2$ moves passively, see \Cref{Mass.on.car}.
The second mass is coupled to the car by a spring-damper combination, and the
ramp is inclined by a fixed angle~$\vartheta \in (0,\pi/2)$.
The equations of motion are given by
\begin{align}\label{eq:ExampleMassOnCarSystem}
    \begin{bmatrix}
        m_1 + m_2& m_2\cos(\vartheta)\\
        m_2 \cos(\vartheta) & m_2
    \end{bmatrix}
    \begin{bmatrix}
        \ddot{z}(t)\\
        \ddot{s}(t)
    \end{bmatrix}
    +
    \begin{bmatrix}
        0\\
        k s(t) +d\dot{s}(t)
    \end{bmatrix}
    =
    \begin{bmatrix}
        u(t)\\
        0
    \end{bmatrix},
\end{align}
where $z(t)$ is the horizontal position of the car and $s(t)$ is the relative position of the mass on
the ramp at time $t$. The physical constants $k>0$ and $d>0$ are the coefficients of the spring and
damper, respectively. The horizontal position of the mass on the ramp is the output $y$ of the system, i.e.
\[
    y(t)=z(t)+s(t)\cos(\vartheta).
\]
The objective is tracking the reference signal $y_{\rf}:t\mapsto \cos(t)$, such
that for $\Funnel\in\cG$ the error function ${t\mapsto e(t)\coloneqq
y(t)-y_{\rf}(t)}$ evolves within the prescribed performance funnel
$\cF_{\Funnel}$, i.e. ${\Norm{e(t)} < \Funnel(t)}$ for all $t\geq 0$.
For this example, we choose the funnel boundary function
\[
    \Funnel(t)\coloneqq5\me^{-2t}+0.1,
\]
which fulfils \eqref{eq:DefinitionAlphaBeta} for $\FunDeriv=2$ and $\FunDiam=0.2$.
By setting $\mu\coloneqq m_2 ( m_1 + m_2\sin^2(\vartheta))$, $\mu_1\coloneqq\frac{m_1}{\mu}$, and
$\mu_2\coloneqq\frac{m_2}{\mu}$, the system takes the form~\eqref{eq:LTISystem}, with
\begin{small}
\[
    x(t)\coloneqq
    \begin{bmatrix}
        z(t)\\
        \dot{z}(t)\\
        s(t)\\
        \dot{s}(t)\\
    \end{bmatrix}
    \!,\,
    A\coloneqq
    \begin{bmatrix}
        0& 1& 0& 0\\
        0& 0& \mu_2k\cos(\vartheta)& \mu_2d\cos(\vartheta)\\
        0& 0& 0& 1\\
        0& 0&-(\mu_1+\mu_2)k& -(\mu_1+\mu_2)d\\
    \end{bmatrix}
    \!,\,
    B\coloneqq
    \begin{bmatrix}
        0\\
        \mu_2\\
        0\\
        -\mu_2\cos(\vartheta)\\
    \end{bmatrix}
    \!,\,
    C\coloneqq
    \begin{bmatrix}
        1\\
        0\\
        \cos(\vartheta)\\
        0\\
    \end{bmatrix}^\top\!\!.
\]
\end{small}
As outlined in~\cite[Sec.~3]{BergIlch21}, 
the system has global relative degree
\[
    r =
        \begin{dcases}
            2, & \vartheta\in\rbl0,\tfrac{\pi}{2}\rbr,\\
            3, & \vartheta=0,
        \end{dcases}
\]
bounded-input bounded-output internal dynamics, and the positive scalar high-frequency
gain $\Gamma=CA^{r-1}B$, see also \Cref{Ex:LTISystem}.
For the simulation, we choose the same system parameters
\begin{equation}\label{eq:Ex:MassOnCar:Params}
    m_1=4,\quad m_2=1,\quad k=2,\quad d=1,\quad \vartheta=\frac{\pi}{4}
\end{equation}
and initial values $z(0)=s(0)=\dot{z}(0)=\dot{s}(0) = 0$ as in~\cite{BergIlch21}. 
Given $\vartheta=\frac{\pi}{4}$, the system has relative degree $r=2$.
Following \cite{IlchWirt13}, the system can equivalently be written in the form 
\begin{equation}\label{eq:MassOnCarInputOutputDeg2}
    \begin{aligned}
        \ddot y(t)    &= R_1y(t)+ R_2\dot{y}(t) + S\eta(t) +\Gamma u(t)\\
        \dot \eta(t)  &= Q \eta(t) + P y(t), 
    \end{aligned}
\end{equation}
with initial conditions $[y(0),\dot{y}(0)]=[y_0^0,y_1^0]\in\R^2$ and $\eta(0)=\eta^0\in\R^2$.
For the given parameters~\eqref{eq:Ex:MassOnCar:Params}, the matrices are
\begin{equation}\label{eq:ParamsMassOnCarInputOutputDeg2}
    R_1=0,\enspace
    R_2=\frac{8}{9},\enspace
    S=\frac{-4\sqrt{2}}{9}
    \begin{bmatrix}
    2 & 1 
    \end{bmatrix},\enspace
    \Gamma=\frac{1}{9},\enspace
    Q=\begin{bmatrix}
    0&1\\
    -4&-2
    \end{bmatrix},\enspace
    P=2\sqrt{2}
    \begin{bmatrix}
       1\\
       0
    \end{bmatrix}.
\end{equation}
To apply the funnel MPC~\Cref{Algo:FunnelMPC} to the system
\eqref{eq:MassOnCarInputOutputDeg2} with relative degree $r=2$, 
we choose the
construction \eqref{eq:psi_i} and \eqref{eq:cond-k_i} for the auxiliary funnel
function $\Funnel_{2}$ and the associated parameter $k_1$. Straightforward
calculations show that  $\Funnel_{2}$ takes the form 
\begin{equation}\label{ex:Sim:MassOnCar:AuxFunnel}
    \Funnel_{2}(t)\coloneqq  \frac{1}{\gamma} k_1  \me^{-\FunDeriv (t-t_0)} + \frac{\FunDiam}{\FunDeriv\gamma}
\end{equation}
with  $k_1=14$ and $\gamma =0.2$ satisfying \eqref{eq:DefParameterGamma}. 
To achieve the control objective, we use the strict funnel stage cost function 
$\ell_{\Funnel_{2}}:\Rp\times\R\times\R\to\R\cup\{\infty\}$ given by 
\begin{equation}\label{eq:Ex:MassOnCar:CostFunction}
\begin{aligned}
    \ell_{\Funnel_2}(t,\zeta,u) =
    \begin{dcases}
        \frac {\Norm{\zeta}}{\Funnel_{2}(t)^2 - \Norm{\zeta}^2} + \lambda_u \Norm{u}^2,
            & \Norm{\zeta} \neq \Funnel_{2}(t)\\
        \infty,&\text{else},
    \end{dcases}
\end{aligned}
\end{equation}
with design parameter ${\lambda_u\in\Rp}$. Utilising the error variables 
$\eM_1(z_1,z_2)\coloneqq z_1$, and $\eM_{2} (z_1,z_2)\coloneqq
\eM_{1}(z_2,0)+k_1 \eM_{1}(z_1,z_2)$ as in \eqref{eq:ErrorVar}, the variable
$\zeta$ in \eqref{eq:Ex:MassOnCar:CostFunction} is replaced, in the optimal
control problem~\eqref{eq:OptimalControlProblem}, by
$\eM_2(\chi(y-y_{\rf})(t)) = \dot{e}(t) + k_1 e(t)$, where
$e(t)  \coloneqq y(t)-y_{\rf}(t)$.

\begin{figure}[h]
    \begin{subfigure}[b]{0.49\textwidth}
        \centering
        \includegraphics[width=\linewidth]{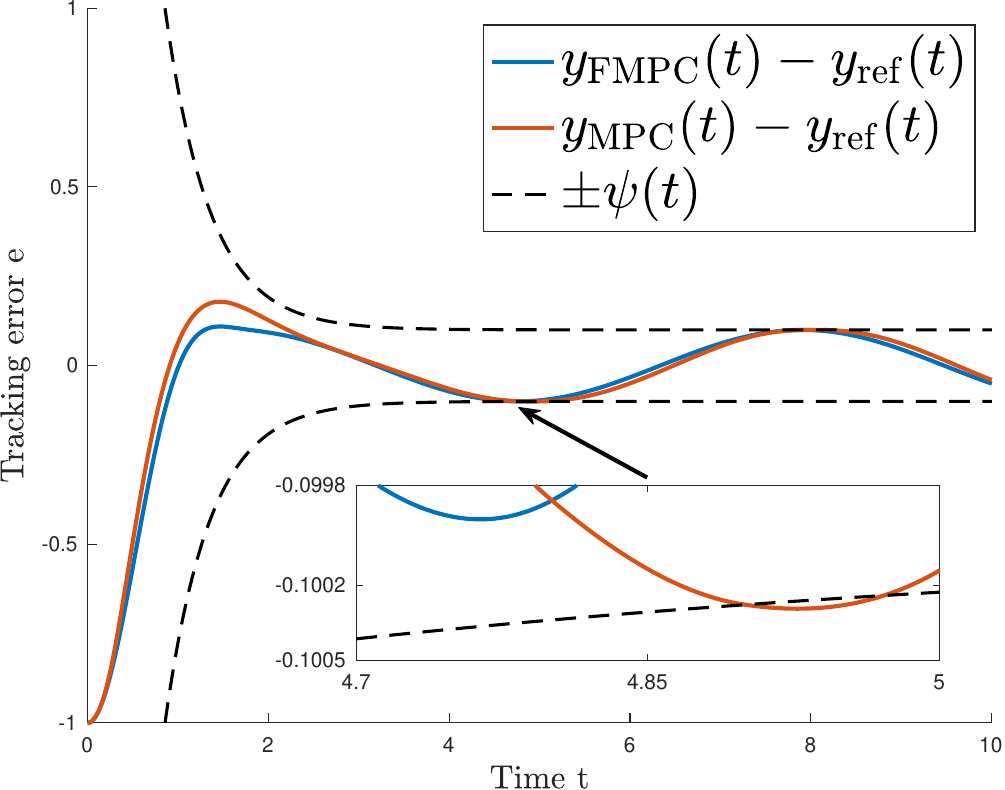}
        \caption{Tracking error $e=y-y_{\rf}$ with boundary~$\Funnel$.}
        \label{fig:sim:massoncar:fmpc_vs_mpc:output}
    \end{subfigure}
      \begin{subfigure}[b]{0.49\textwidth}
        \centering
        \includegraphics[width=\linewidth]{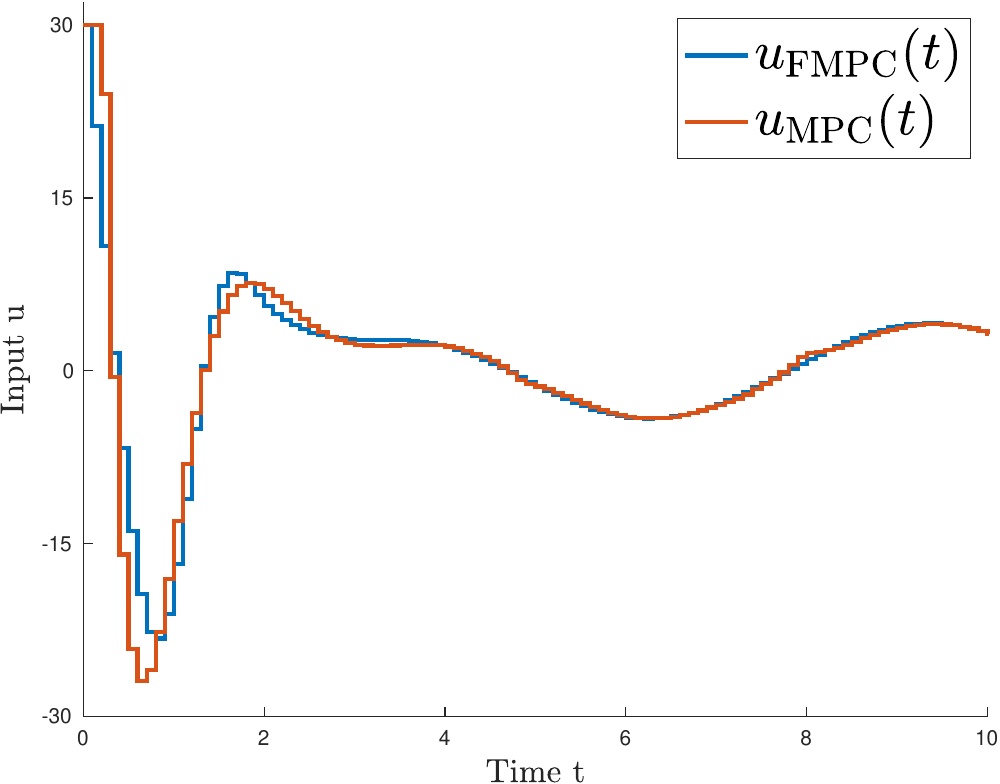}
        \caption{Control inputs.}
        \label{fig:sim:massoncar:fmpc_vs_mpc:input}
    \end{subfigure}
    \caption{Simulation of system~\eqref{eq:MassOnCarInputOutputDeg2} under the
            control generated by funnel MPC~\Cref{Algo:FunnelMPC} 
            and the MPC~\Cref{Algo:MPC} with cost function \eqref{eq:stageCostClassicalMPC}
            and \cref{eq:ConstraintClassicalMPC}. 
            The  parameters are $T= 1$, $\delta=0.1$, and $\lambda_u=10^{-4}$.
    } 
    \label{fig:sim:massoncar:fmpc_vs_mpc}
\end{figure}

We compare the funnel MPC~\Cref{Algo:FunnelMPC} with the standard
MPC~\Cref{Algo:MPC} using the quadratic cost function 
\begin{equation}\label{eq:Ex:QuadraticStageCostMassOnCar}
    \ell(t,y,u)= \Norm{y-y_{\rf}(t)}^2+\lambda_u \Norm{u}^2
\end{equation}
as in \eqref{eq:stageCostClassicalMPC} and additional output constraints
\[
    \fa t\in[t_k,t_k+T]:\quad \Norm{y(t)-y_{\rf}(t)} <\Funnel(t) 
\]
for $t_k\in\delta\N_0$ as in \eqref{eq:ConstraintClassicalMPC} in the OCP~\eqref{eq:MpcOCP}.
For both MPC schemes, we choose the prediction horizon $T=1$, the time shift
$\delta=0.1$, the parameter $\lambda_u=10^{-4}$, and allow for a maximal control
value of $\umax=30$. Due to discretisation, only step functions with constant
step length~$\SampleTime\coloneqq\delta$ are considered when solving the respective
optimal control problems.
The simulations are performed on the time interval $[0,10]$ with \textsc{Matlab}
and the toolkit \textsc{CasADi} and displayed in \Cref{fig:sim:massoncar:fmpc_vs_mpc}.
While \Cref{fig:sim:massoncar:fmpc_vs_mpc:input} shows that the control signals generated by the two 
MPC schemes are relatively similar, \Cref{fig:sim:massoncar:fmpc_vs_mpc:output} displays that 
the MPC~\Cref{Algo:MPC} with cost function \eqref{eq:Ex:QuadraticStageCostMassOnCar} does, 
contrary to the funnel MPC~\Cref{Algo:FunnelMPC}, not achieve the control objective.
Despite the output constraints~\eqref{eq:ConstraintClassicalMPC}, the tracking error 
$y_{\mathrm{MPC}}(t)-y_{\rf}(t)$ leaves the performance funnel. 
The reason for this behaviour could be that the solver implements these barriers
internally with a certain tolerance level using barrier functions.
To ensure an adherence to the funnel boundaries,  an adaptation of the parameter
$\lambda_u$, a smaller step length $\delta$, or a longer prediction horizon $T$ are
sufficient as demonstrated in the previous example. Even though the tracking
error $y_{\mathrm{FMPC}}-y_{\rf}$ evolves at times close to the funnel boundary,
the cost function \eqref{eq:Ex:MassOnCar:CostFunction} ensures that the control
objective is achieved when the system is controlled by the funnel
MPC~\Cref{Algo:FunnelMPC}.

Now, we compare the funnel MPC~\Cref{Algo:FunnelMPC} to the funnel controller from~\cite{BergIlch21}.
For the system \eqref{eq:MassOnCarInputOutputDeg2}, the funnel control law takes the form
\begin{equation}\label{eq:Ex:MassOnCarRelDeg2:FC}
    \begin{aligned}
        w(t)&=\frac{\dot{e}(t)}{\Funnel(t)}+\FCBijec\rbl\frac{e(t)^2}{\Funnel(t)^2}\rbr\frac{e(t)}{\Funnel(t)},\qquad e(t)=y(t)-y_{\rf}(t),\\
        \uFC(t)&=-\FCBijec\rbl w(t)^2\rbr w(t),
    \end{aligned}
\end{equation}
with $\FCBijec(s)=\frac{1}{1-s}$ for $s\in[0,1)$. 
For the funnel MPC scheme, we choose the prediction horizon $T=1$, the time
shift $\delta=0.1$, the parameter $\lambda_u=10^{-3}$, and allow for a maximal
control value of $\umax=30$.

\begin{figure}[h]
    \begin{subfigure}[b]{0.49\textwidth}
        \centering
        \includegraphics[width=\linewidth]{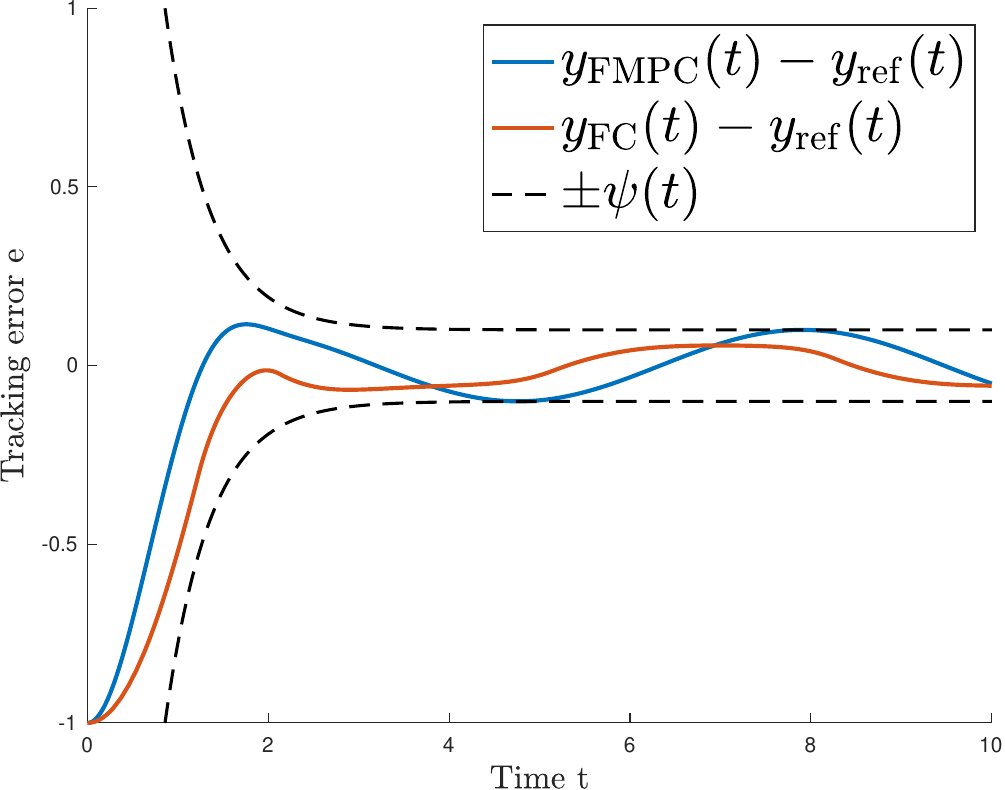}
        \caption{Tracking error $e=y-y_{\rf}$ within boundary~$\Funnel$.}
        \label{fig:sim:massoncar:fmpc_vs_FC:output}
    \end{subfigure}
      \begin{subfigure}[b]{0.49\textwidth}
        \centering
        \includegraphics[width=\linewidth]{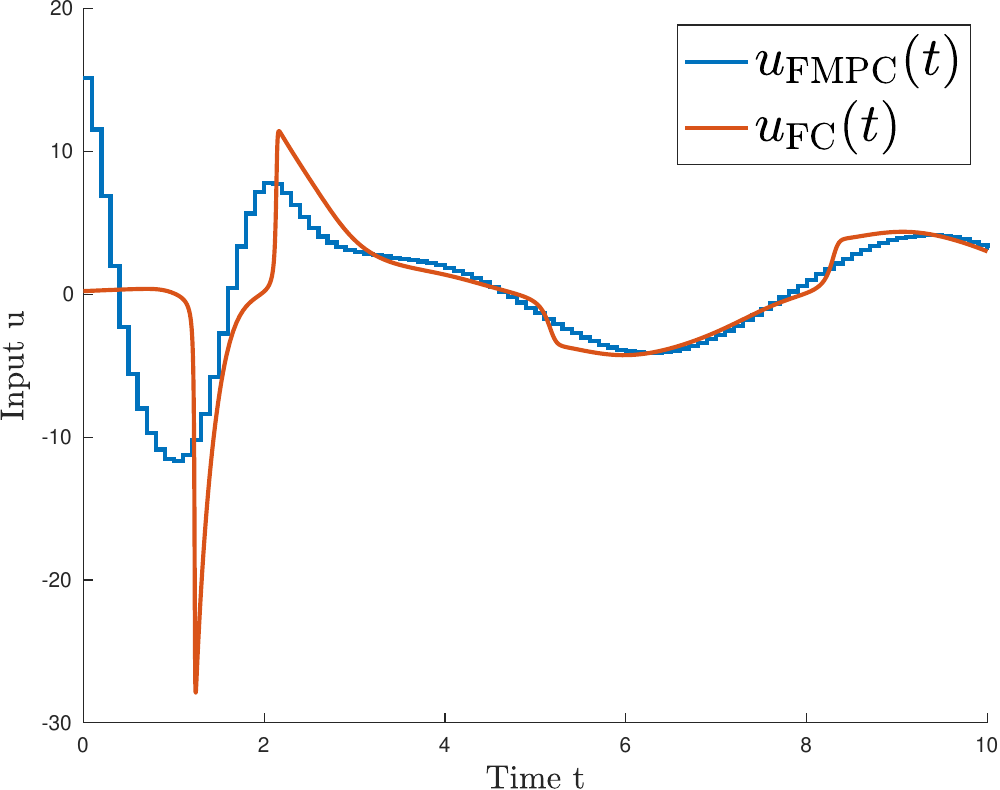}
        \caption{Control inputs.}
        \label{fig:sim:massoncar:fmpc_vs_FC:input}
    \end{subfigure}
    \caption{Simulation of system~\eqref{eq:MassOnCarInputOutputDeg2} under the
            control generated by funnel MPC~\Cref{Algo:FunnelMPC} 
            with parameters $T= 1$, $\delta=0.1$, and $\lambda_u=10^{-3}$
            and the funnel control law~\eqref{eq:Ex:MassOnCarRelDeg2:FC}.
    } 
    \label{fig:sim:massoncar:fmpc_vs_FC}
\end{figure}
The performance of the funnel controller~\eqref{eq:Ex:MassOnCarRelDeg2:FC} and the funnel MPC 
\Cref{Algo:FunnelMPC} is depicted in \Cref{fig:sim:massoncar:fmpc_vs_FC}. While 
\Cref{fig:sim:massoncar:fmpc_vs_FC:output} shows the tracking error of the two controllers evolving within
the funnel boundaries, \Cref{fig:sim:massoncar:fmpc_vs_FC:input} displays the respective input signals.
It is evident that both control techniques are feasible and achieve the control objective.
The funnel controller generates a smooth input signal, while the OCP
\eqref{eq:OptimalControlProblem} of the funnel MPC \Cref{Algo:FunnelMPC} is
solved over step functions with constant step length~$\SampleTime:=0.1$.
The funnel MPC seemingly takes more advantage of the available error tolerance
boundaries resulting in a smaller range of employed control values.
Funnel control tends to change the control values
very quickly and the control signal shows peaks. The MPC scheme avoids this 
undesirable behaviour thanks to its predictive capabilities.
\chapter{Robust funnel MPC}\label{Chapter:RobustFunnelMPC}
Optimisation-based control techniques, such as model predictive control, achieve  
high-performance control while rigorously adhering to state and input constraints. 
These schemes~-- including the funnel MPC~\Cref{Algo:FunnelMPC} developed
in~\Cref{Chapter:FunnelMPC} -- fundamentally depend on accurate system models.
Without such models, essential closed-loop properties -- stability,
performance, and constraint satisfaction --  are generally not preserved.
Significant challenges arise from model uncertainties and external
disturbances, as even high-fidelity models deviate from real-world systems, while 
disturbances are omnipresent. Moreover, to mitigate computational complexity, 
practitioners often opt for simplified, lower-dimensional approximations -- such
as discretised representations of partial differential equations -- over
intricate models. For a comprehensive treatment of model order reduction
techniques, see for example the textbook~\cite{schilders2008model}.

The development of robust MPC methods to address structural model-plant
mismatches and external disturbances therefore remains an active research area,
see e.g.~\cite{Buja21,KohlSolo20,SunDai19,RakoDai22} and the references therein.
Key approaches include:
\begin{itemize}
    \item \emph{Scenario-based optimisation}: Handles uncertainties via sampling a
    suitable number of randomly selected disturbance realisations in a receding
    horizon fashion \cite{calafiore2012robust}.
    \item \emph{Barrier-augmented MPC}: Ensure states/outputs to remain within safe
    regions as (relaxed) barrier functions penalise proximity to constraint
    boundaries. Safety and constraint satisfaction is enforced through dynamic penalty adjustment 
    and inherently accounting for deviations \cite{Petsag2021,Feller2016,Yin2023}.
    \item \emph{Feedback MPC}: Solves for an optimal and stabilising feedback policy rather than an open-loop input signal \cite{Scokaert1998,goulart2006optimization}.
    The applied (robust) feedback controller counteracts occurring disturbances between two iterations of the MPC algorithm.
    \item \emph{Adaptive MPC}: Dynamically updates model parameters online using
    techniques like moving horizon estimation (MHE)~\cite{Haseltine2005},
    (non)-linear state observers \cite{Kalman61,Besancon2007,Korder2022}, or system
    identification methods \cite{Rao2006} bridging model-system gaps, see e.g.
    \cite{Adetola2011,Sasfi2023}. For a comprehensive overview on adaptive MPC,
    see also the survey paper~\cite{KIM2010}. 
    \item \emph{Stochastic MPC}: Employs chance constraints or risk-aware formulations
    for quantifiable probabilistic uncertainties. It offers probabilistic
    guarantees for systems with measurable noise distributions \cite{Kouvaritakis2015,mesbah2016,Singh2019}.
    \item \emph{Learning-augmented MPC}: Integrates data-driven models, such as
    Gaussian processes or neural networks, to refine predications and quantify
    uncertainties~\cite{Aswa13}. We explore this integration in more detail in
    \Cref{Chapter:LearningRobustFMPC}; see also \cite{HewingWaber20} for a survey.
\end{itemize}
Central to robust MPC are  constraint tightening
techniques~\cite{chisci2001systems}, particularly \emph{tube-based}
MPC~\cite{langson2004robust}. 
To robustly achieve output tracking, these
methods construct tubes around reference trajectories to guarantee the actual
system output remains within prescribed bounds. For linear systems, 
foundational work in \cite{MaynSero05} demonstrates this approach, while non-linear extensions
in \cite{FaluMayn14,KohlSolo20,RakoDai22} address geometric and
dynamic complexities. Notably, \cite{Lopez2019} introduces
co-optimisation of tubes and reference trajectories, adapting tube geometry
based on proximity to boundaries. 

To enforce tube invariance, terminal
conditions are embedded within the optimisation problem, ensuring recursive
feasibility. For linear systems, \cite{CairBorr16} 
achieves reference tracking within constant bounds via \emph{robust control invariant
(RCI)} sets, which satisfy state, input, and performance constraints.
\cite{yuan2019bounded} extends this framework to external disturbances, though
RCI computation remains non-trivial, with algorithms potentially failing to
terminate finitely \cite{CairBorr16}. For non-linear systems,
\cite{yu2013tube,singh2017robust} employ incremental Lyapunov functions and
precomputed stabilising feedback laws to ensure control objectives. While
effective, these methods face challenges in balancing conservatism and
computational tractability, as tube design must inherently account for system
uncertainty magnitude. 

Despite advancements in robustification methods for MPC, critical challenges persist:
\begin{itemize}
    \item Computational complexity: Scaling methods for high-dimensional systems \cite{KohlSolo20,Gesser2018}.
    \item Conservatism vs performance: Balancing conservatism and performance, in particular in tube-based approaches.
    \item Safety certification: Ensuring reliability in learning-augmented components \cite{tambon2022certify}.
\end{itemize}

\subsubsection{Robust funnel MPC: Bridging prediction and adaptation}
To address the challenge of output tracking within prescribed performance
boundaries while retaining the predictive power of MPC and the disturbance
rejection capabilities of adaptive control, this chapter proposes \emph{robust funnel MPC}.
This method relaxes the assumption from
\Cref{Chapter:FunnelMPC} that the system \eqref{eq:Sys} and surrogate model
\eqref{eq:Intro:ModelEquation} coincide, explicitly accounting for external
disturbances and (structural) model-plant mismatches.
The controller synergises two complementary strategies:
\begin{enumerate}
    \item \textbf{Funnel MPC}: Leverages model-based predictions to compute feed-forward control signals.
    \item \textbf{Funnel control}: A model-free, high-gain adaptive feedback
    loop (introduced in \Cref{Sec:FunnelControl}) that refines the control
    signal using real-time measurements to reject disturbances and compensate
    mismatches.
\end{enumerate}
The synergy of these techniques ensures arbitrary output constraint
satisfaction: the predictive component (funnel MPC) plans trajectories using the
surrogate model, while the model-free adaptive component (funnel control)
instantaneously compensates for unmodelled dynamics or disturbances. This two
component approach marries the predictive power of MPC with the robustness of
adaptive feedback, addressing key limitations of stand-alone methods in
uncertain environments.

\section{System class}\label{Sec:SystemClass}

In this section, we concretise the structural properties of the
system~\eqref{eq:Sys} and formally introduce the system class under
consideration. 
To briefly recapitulate, we consider non-linear multi-input multi-output control
systems of order $r\in\N$ of the form 
\begin{equation}\tag{\ref{eq:Sys} revisited}
    \begin{aligned}
    &\textover[r]{$y^{(r)}(t)$}{$\big(y(t_0),\ldots,y^{(r-1)}(t_0)\big)\ $}   = F(\oT(y,\dot{y},\ldots, y^{(r-1)})(t),u(t)),\\
    &\left.
        \begin{aligned}
       y|_{[0,t_0]}&= y^0  \in \cC^{r-1}([0,t_0],\R^m), && \mbox{if } t_0 >0,\\
       \big(y(t_0),\ldots,y^{(r-1)}(t_0)\big)&= y^0  \in\R^{rm}, && \mbox{if } t_0 = 0,
        \end{aligned} \right\}
    \end{aligned}
\end{equation}
with $t_0\geq 0$,  initial trajectory $y^0$, input $u\in L_{\loc}^{\infty}([t_0,\infty),\R^m)$,
and output $y(t)\in\R^m$ at time $t\geq t_0$. 
The following definition formalises the properties of the function $F$ and the operator~$\oT$.
\begin{definition}[System class $\cN^{m,r}_{t_0}$]\label{Def:SystemClass}
    We say that the system~\eqref{eq:Sys} belongs to the \emph{system class}~$\cN^{m,r}_{t_0}$ for $m,r\in\N$, and $t_0\in\Rp$,
    written~$(F,\oT)\in\cN^{m,r}_{t_0}$, if, for some $q\in\N$, the following holds:
    \begin{enumerate}[(a)]
    \item\label{Item:Prop:OperatorSys}
        $\oT:\cR(\Rp,\R^n) \to L^\infty_{\loc} ([t_0,\infty), \R^{q})$ has the 
        \emph{causality}~\ref{Item:OperatorPropCasuality},
        \emph{local Lipschitz}~\ref{Item:OperatorPropLipschitz},
        and the \emph{bounded-input bounded-output (BIBO)}~\ref{Item:OperatorPropBIBO}
        property as defined in~\Cref{Def:OperatorClass}.
        \item\label{Item:PerturbationHighGain}$F\in\cC(\R^q \times \R^m,\R^m)$ has the \emph{perturbation high-gain property}, i.e. 
        for every compact set $K_m \subset \R^m$ there exists~$\nu\in(0,1)$ such that
        for every compact set  $K_q\subset\R^q$ the function
        \begin{equation}\label{eq:Def:HighGainFunction}
            \HighGainFunc\colon\R\to\R, \
            s \mapsto \min
            \setdef{\langle v, F(z, d-s v)\rangle}
            {
                 d\in K_m, z \in  K_q, v\in\R^m,~\nu \leq \|v\| \leq 1
            }
        \end{equation}
        satisfies $\sup_{s\in\R} \HighGainFunc(s)=\infty$. 
    \end{enumerate}
\end{definition}

We already discussed examples for operators $\oT$ satisfying the properties
\ref{Item:Prop:OperatorSys} from \Cref{Def:SystemClass} in
\Cref{Ex:LTISystem,Ex:ControlAffineMod}. To also gain a better understanding for
the high-gain property of the function $F$ in~\eqref{eq:Sys}, we briefly discuss
a simple example of a differential equation belonging to the considered system
class.

\begin{example}\label{Ex:ControlAffineSignDefinite}
    Let $p:\R^n\to \R^n$ and $\Gamma:\R^n\to \R^{n\times m}$ 
    be continuous non-linear functions. 
    Assume $\Gamma(x)\in\GL_{m}(\R)$ for all $x\in\R^n$.
    We show that the function $F:\R^n\times\R^m\to\R^n$ defined by 
    \begin{equation}\label{eq:ExHighGain}
        F(x,u)\coloneqq p(x)+\Gamma(x)u
    \end{equation}
    has the perturbation high-gain property~\eqref{eq:Def:HighGainFunction} if, and only if,
    $\Gamma(x)$ is sign-definite for all $x\in\R^n$, i.e. the scalar product
    $\al v,\Gamma(x)v\ar$ is positive (negative) for all $v\in\R^m\backslash\{0\}$.
    We show this equivalence by adapting~\cite[Sec. 2.1.3]{BergIlch21} to the given context. 

    Assume that the function $F$ has the perturbation high-gain property~\eqref{eq:Def:HighGainFunction} 
    and suppose that $\Gamma$ is not sign-definite.
    Then, there exists $z\in\R^n$ and $v\in\R^m\backslash\{0\}$ with $\al v,\Gamma(z)v\ar=0$.
    Define $K_{m}\coloneqq \{0\}$ and $K_{n}=\{z\}$. 
    For~$\nu\in(0,1)$, set $K_{\nu}\coloneqq \setdef{v\in\R^m}{\nu\leq\Norm{v}\leq 1}$.
    As $\al v,\Gamma(z)v\ar=0$, there exists $\hat{v}\in K_{\nu}$ with $\al \hat{v},\Gamma(z)\hat{v}\ar=0$.
    For $s\in\R$, we have
    \begin{align*}
        \HighGainFunc(s)
        &=\min_{v\in K_{\nu}}\al v,F(z,-sv)\ar
        =\min_{v\in K_{\nu}}\al v,p(z)-\Gamma(z)sv\ar\\
        &\leq\nu\Norm{p(z)}+\min_{v\in K_{\nu}}-s\al v,\Gamma(z)v\ar
        \leq\nu\Norm{p(z)}-s\al \hat{v},\Gamma(z)\hat{v}\ar
        = \nu\Norm{p(z)}. 
    \end{align*}
    This is a contradiction to the perturbation high-gain property~\eqref{eq:Def:HighGainFunction}.

    Assume $\Gamma$ is sign-definite. Due to the continuity of $\Gamma(\cdot)$,
    there exists $\sigma\in\{-1,1\}$ such that $\sigma\Gamma(z)$ is positive definite for all $z\in\R^n$.
    We show that the function $F$ has the perturbation high-gain property~\eqref{eq:Def:HighGainFunction}.
    Let $K_{m}\in\R^m$, $K_{n}\in\R^n$  be compact sets and set $\nu=\tfrac{1}{2}$. 
    Define $K_{\nu}\coloneqq \setdef{v\in\R^m}{\nu\leq\Norm{v}\leq 1}$.
    Set $G(z)\coloneqq \tfrac{\sigma}{2}(\Gamma(z)+\Gamma(z)^\top)$ and let $\lambda_{\min}$ be the smallest eigenvalue of
    $G(z)$ for all $z\in K_n$ which exists because of the compactness of $K_n$.
    Moreover, due to the continuity of the involved functions and the compactness of the considered sets, there exists
    \[
        c\coloneqq \min \setdef{\langle v, p(z)+\Gamma(z)d\rangle}
            {
                 d\in K_m, z \in  K_n, v\in K_\nu
            }\in\R.
    \]
    Let $(s_j)_{j\in\N}\in\R^\N$ be a sequence with $s_j\sigma <0$ for all $j\in\N$ and $s_j\sigma\to-\infty$ for $j\to\infty$.
    It follows that
    \begin{align*}
    \HighGainFunc(s_j)&=\min
            \setdef{\langle v, F(z, d-s_j v)\rangle}
            {
                \ d\in K_m, z \in  K_n, v\in K_\nu
            }\\
            &\geq
            \min\setdef{\langle v, p(z)+\Gamma(z)d\rangle}
            {
                 \ d\in K_m, z \in  K_n, v\in K_\nu
            }\\
           &\hphantom{= } + 
            \min\setdef{-\langle v, \Gamma(z)s_j v\rangle}
            {
               \   z \in  K_n, v\in K_\nu
            }\\
            &
            =c+ 
            \min\setdef{-s_n\sigma\langle v, G(z)v\rangle}
            {
                \  z \in  K_n, v\in K_\nu
            }\\
            &\geq c+ 
            \min\setdef{-s_j\sigma \lambda_{\min}\Norm{v}^2}
            {
                \  v\in K_\nu
            }\\
            &\geq  c-\frac{s_j\sigma \lambda_{\min}}{4}.
    \end{align*}
    Thus, $\HighGainFunc(s_j)\to\infty$ for $j\to\infty$ proving that the function $F$ has the perturbation high-gain property~\eqref{eq:Def:HighGainFunction}.
    
    We saw in~\Cref{Ex:ControlAffineMod} that non-linear differential equations of the form 
    \begin{equation} \tag{\ref{eq:NonLinSysStates} revisited}
    \begin{aligned}
        \dot{x}(t)  & = f(x(t)) + g(x(t)) u(t),\quad x(t_0)=x^0,\\
        y(t)        & = h(x(t)),
    \end{aligned}
    \end{equation}
    with~$t_0\in\R_{\ge 0}$, $x^0\in\R^n$, and non-linear functions~$f:\R^n\to \R^n$, $g:\R^n\to \R^{n\times m}$ and $h : \R^n \to \R^m$,
    are admissible candidates for a model by transforming it into the {Byrnes-Isidori form}~\eqref{eq:BIF}.
    This was achieved by, among other things, assuming that 
    \eqref{eq:NonLinSysStates} has a strict (global) relative degree~$r \in \mathbb{N}$, i.e.
    \begin{align*}
        \fa k \in \{1,\ldots,r-1\}\ \fa x \in \R^n:\
            (L_g L_f^{k-1} h)(x)  &= 0  \\
          \text{and}\quad
            (L_g L_f^{r-1} h)(x)  &\in \GL_{m}(\R).
    \end{align*}
    A consequence of our considerations regarding function $F$ in~\eqref{eq:ExHighGain} is that 
    \eqref{eq:NonLinSysStates} is also an admissible system if, in addition to the strict relative degree,
    one assumes $(L_g L_f^{r-1} h)(x)$ to be sign-definite. 
    Similarly, linear time invariant systems
    given by matrices $A \in \R^{n \times n}$ and $C^\top,B\in\R^{n\times m}$, as discussed in~\Cref{Ex:LTISystem}, 
    are admissible systems in the sense of~\Cref{Def:SystemClass} if, 
    in addition to the assumptions from ~\Cref{Ex:LTISystem}, 
    the matrix $CA^{r-1}B$ is sign-definite, where $r>0$ is the relative degree of the linear system~\eqref{eq:LTISystem}.
\end{example}

\begin{remark}\label{Rem:SystemClass}
    We want to comment on a few aspects of the system class $\cN^{m,r}_{t_0}$.
    \begin{enumerate}[(a)]
    \item\label{Item:OperatorContainsDisturbences}
    For $t_0\geq0$ and $m,r\in\N$, let $(F,\oT)\in\cN^{m,r}_{t_0}$.
    For $d\in L^\infty([t_0,\infty),\R^p)$,  
    define the operator $\tilde{\oT}$ by
    \[
        \tilde{\oT}(\zeta)(t)\coloneqq (d(t),\oT(\zeta)(t))
    \]
    for $\zeta\in\cR(\Rp,\R^n)$.
    Straightforward calculations show that $\tilde{\oT}$ also fulfils 
    the properties~\ref{Item:OperatorPropCasuality},~\ref{Item:OperatorPropLipschitz},
    and~\ref{Item:OperatorPropBIBO} of~\Cref{Def:OperatorClass}.
    The system class~$\cN^{m,r}_{t_0}$ therefore implicitly contains differential equations of the form
    \[
        y^{(r)}(t)= F(d(t),\oT(y,\dot{y},\ldots, y^{(r-1)})(t),u(t)),
    \]
    with unknown disturbance $d\in L^\infty([t_0,\infty),\R^p)$.
    \item 
    The system class~$\cN^{m,r}_{t_0}$ allows for 
    the usage of more general operators $\oT$ than the model class~$\cM^{m,r}_{t_0}$ because
    $\oT$ is not required to have the \emph{limited memory} property~\ref{Item:OperatorPropLimitMemory} of~\Cref{Def:OperatorClass}.
    Many physical phenomena such as \emph{backlash} and \emph{relay hysteresis},
    and \emph{non-linear time delays} can be modelled by means of a general operator~$\oT$, cf.~\cite[Sec.~1.2]{BergIlch21}.
    Moreover, the operator~$\oT$ can even be the solution operator of an infinite-dimensional dynamical system,
    e.g. a partial differential equation.
    Thus, systems with such internal dynamics can be represented by~\eqref{eq:Sys}, see~\cite{BergPuch20}.
    For a practically relevant example of infinite-dimensional internal dynamics (modelled by an operator~$\oT$),
    we refer to~\cite{BergPuch22}, where a moving water tank was
    subject to funnel control, and the water in the tank was modelled by the linearised Saint-Venant equations. 
    While we deem the limited memory property~\ref{Item:OperatorPropLimitMemory} not to be a major restriction posed 
    on the operator $\oTM$ used in the model,  
    it still remains to be verified whether the mentioned examples can also be modelled by an operator with 
    property~\ref{Item:OperatorPropLimitMemory}.
    \item 
    The perturbation high-gain property of the function $F$ in~\ref{Item:PerturbationHighGain} of~\Cref{Def:SystemClass}
    is a modification of the so-called \emph{high-gain} property, see e.g. \cite[Def.~1.2]{BergIlch21},
    and,  at first glance, a stronger assumption.
    The high-gain property is essential in high-gain adaptive control and,
    roughly speaking, guarantees that, if a large enough input is applied,
    the system reacts sufficiently fast.
    For linear systems, as in \Cref{Ex:LTISystem}, having the high-gain property implies that the system can
    be stabilised via high-gain output feedback, cf.~\cite[Rem.~1.3]{BergIlch21}. In
    order to account for possible bounded perturbations of the input, we require the
    modified property from~\ref{Item:PerturbationHighGain}. It is an open question whether the
    perturbation high-gain property and the high-gain property are equivalent. 
    \item Although there are many systems belonging to both the model class~$\cM^{m,r}_{t_0}$ from~\Cref{Def:ModelClass}
    and the system class~$\cN^{m,r}_{t_0}$ from~\Cref{Def:SystemClass},
    neither the set of admissible models~$\cM^{m,r}_{t_0}$ is a subset of all considered systems~$\cN^{m,r}_{t_0}$ nor the opposite is true.  
    Every system $(F,\oT)\in\cN^{m,r}_{t_0}$ which does not have a control affine representation of the form~\eqref{eq:Model_r}
    cannot belong to $\cM^{m,r}_{t_0}$. On the other hand, \Cref{Ex:ControlAffineSignDefinite}
    shows that differential equations of the form~\eqref{eq:NonLinSysStates} are admissible models if $(L_g L_f^{r-1} h)(x)$
    is invertible but only admissible systems if $(L_g L_f^{r-1} h)(x)$ is in addition sign definite.
    \item Throughout this thesis, we always assume that the parameters $m$ and
    $r$ for the system class $\cN^{m,r}_{t_0}$ and the model class $\cM^{m,r}_{t_0}$ coincide.
    This means that the system~\eqref{eq:Sys} and the model~\eqref{eq:Model_r}
    are of the same order $r\geq1$ and have the identical output/input dimension $m\geq1$.
   \end{enumerate}\unskip
\end{remark}

In \Cref{Def:ModSolution}, we introduced a solution concept for the initial
value problem~\eqref{eq:Model_r}, which is used as model for the funnel
MPC~\Cref{Algo:FunnelMPC}. Mainly due to the inherent conflict between the
domain of the operator and the re-initialisation of the model, it
had certain peculiarities distinguishing it from more traditional solution
concept. As the system's differential equation~\eqref{eq:Sys} is not
re-initialised during operation of any controller, we utilise conventional solutions 
in sense of \textit{Carath\'{e}odory}.
For the sake of completeness, we recall this solution concept.
\begin{definition}[System solution]\label{Def:FCSolutionConcept}
For initial trajectory $y^0\in\cC^{(r-1)}([0,t_0],\R^m)$ for $t_0>0$ or $y_0\in\R^{rm}$ in the case $t_0=0$
and a control function $u\in L^\infty_{\loc}([t_0,\infty), \R^m)$, 
an absolutely continuous function $x=(x_1,\ldots,x_r):[0,\omega)\to\R^{rm}$  with $\omega\in(t_0,\infty]$ 
is called a solution of~\eqref{eq:Sys} (in the sense of \textit{Carath\'{e}odory})
if 
\begin{align*}
    \dot{x}_{i}(t)&=x_{i+1}(t),\hspace{2cm} i=1,\ldots, r-1,\\
    \dot{x}_r(t)&= F(\oT(x)(t),u(t)),
\end{align*}
for almost all $t\in[t_0,\omega)$ and $x|_{[0,t_0]}= \OpChi(y^0)$ if $t_0>0$ or $x(t_0)= y^0$ in the case $t_0=0$.
A solution is maximal if it has no proper right extension that is also a solution.
A maximal solution is also called a \emph{response} of the system associated with $u$
and denoted by $x(\cdot;t_0,y^0,u)$. We denote its first component $x_1$ by $y(\cdot,t_0,y^0,u)$.
\end{definition}

In the \nameref{Chapter:Appendix}, we show that \eqref{eq:Sys} has  a solution
$x:[0,\omega)\to\R^{rm}$ in the sense of \Cref{Def:FCSolutionConcept} for every
$u\in L^\infty_{\loc}([t_0,\infty), \R^m)$ and that every solution can be
extended to a maximal solution, see \Cref{Appendix:Cor:SystemSolutionExists}.

\section{Controller structure}\label{Sec:RobustFMPCControllerStructure}

\begin{figure}[ht]
    \centering
    \scalebox{.85}{
        \begin{tikzpicture}[very thick,%
        scale=0.64,%
        node distance = 9ex,
        box/.style={fill=white,rectangle, draw=black},
        blackdot/.style={inner sep = 0, minimum size=3pt,shape=circle,fill,draw=black},%
        blackdotsmall/.style={inner sep = 0, minimum size=0.1pt,shape=circle,fill,draw=black},%
        plus/.style={fill=white,circle,inner sep = 0,very thick,draw},%
        metabox/.style={inner sep = 3ex,rectangle,draw,dotted,fill=gray!20!white}]
        \begin{scope}[scale=0.5]
            \node (sys) [box,minimum size=9ex,xshift=-1ex, fill=orange!60]  {System \eqref{eq:Sys}};
            \node(FC) [box, below of = sys,yshift=-8ex,minimum size=9ex] {Funnel controller};
            \node(fork1) [plus, right of = FC, xshift=18ex] {$+$};
            \node(fork9) [blackdot, inner sep = 0pt, right of = fork1, xshift=-2ex ] {};
            \node(fork2) [plus, left of = FC, xshift=-15ex] {$+$};
            \node(fork3) [blackdot, left of = fork2, xshift=-0ex] {};
            \node(MPC) [box, left of = fork3,xshift=-8ex,minimum size=9ex] {FMPC};
            \node(MPCin1) [minimum size=0pt, inner sep = 0pt, below of = MPC, yshift=4.5ex, xshift=2ex] {};
            \node(MPCin1Desc) [minimum size=0pt, inner sep = 0pt, below of = MPCin1, yshift=5ex, xshift=2.5ex] {$y$};
            \node(MPCin2) [minimum size=0pt, inner sep = 0pt, below of = MPC, yshift=4.5ex, xshift=-2ex]{};
            \node(MPCin2Desc) [minimum size=0pt, inner sep = 0pt, below of = MPCin2, yshift=5ex, xshift=-2.5ex] {$y_{\rf}$};
            \node(refin) [minimum size=0pt, inner sep = 0pt, below of = MPC, yshift=-2ex, xshift=-2ex] {};
            \node(Mod) [box, above of = MPC,yshift=8ex,minimum size=9ex] {Model \eqref{eq:Model_r}};
            \node(fork4) [blackdot, left of = MPC, xshift=-5ex] {};
            \node(fork5) [minimum size=0pt, inner sep = 0pt, below of = fork4, yshift=-5ex] {};
            \node(fork6) [minimum size=0pt, inner sep = 0pt, below of = fork9, yshift=-7ex] {};
            \node(fork7) [minimum size=0pt, inner sep = 0pt, below of = MPCin1, yshift=-2.5ex] {};
            \draw[->] (refin) -- (MPCin2) node[pos=0.4,left] {};
            \draw[->] (MPC) -- (fork2) node[pos=0.2,above] {$u_{\mathrm{FMPC}}$};
            \draw[->] (fork3) |- (Mod);
            \draw (Mod) -| (fork4) node[pos=0.3,above] {$y_{\mathrm{M}}$};
            \draw[->] (fork4) -- (MPC);
            \draw[-] (sys) -| (fork9) node[pos=0.05,above] {$y$} ;
            \draw[->] (fork9) -- (fork1) node[pos=0.6,right, above] {$+$};
            \draw[-] (fork9) -- (fork6.south);
            \draw[-] (fork6.east) -- (fork7.west);

            \path[->,name path=line1] (fork7.south) -- (MPCin1){};
            \draw[->,name path=line2] (fork5.west) -| (fork1) node[pos=0.9,left] {$-$};
            \path [name intersections={of = line1 and line2}];
            \coordinate (S)  at (intersection-1);          
            \path[name path=circle] (S) circle(5.mm);
            \path [name intersections={of = circle and line1}];
            \coordinate (I1)  at (intersection-1);
            \coordinate (I2)  at (intersection-2);
            \tkzDrawArc[color=black, very thick](S,I1)(I2);
            \draw[-] (fork7.south) |- (I2);
            \draw[->] (I1) -- (MPCin1);
            \draw[->] (fork1) -- (FC) node[midway,above] {$e_{\mathrm{S}}=y - y_{\mathrm{M}}$};
            \draw[->] (FC) -- (fork2) node[pos=0.4,above]{$u_{\mathrm{FC}}$};
            \draw[->] (fork2) |- (sys) node[pos=0.69,above] {$u=u_{\mathrm{FMPC}} + u_{\mathrm{FC}}$};
            \draw (fork4) -- (fork5.south);
        \end{scope}
        \begin{pgfonlayer}{background}
            \fill[red!20] (-18,-7.5) rectangle (-8,2.);
            \fill[blue!20] (-5.,-7.5) rectangle (8,-2.);
            \node at (-12.9,-9.5) {\color{red}{\large Model-based controller component}};
            \node at (1.5,-9.5) {\color{blue}{\large Model-free controller component}};
        \end{pgfonlayer}
\end{tikzpicture} 
    } 
    \caption{Structure of the robust funnel MPC scheme}
    \label{Fig:Robust_FMPC}
\end{figure}
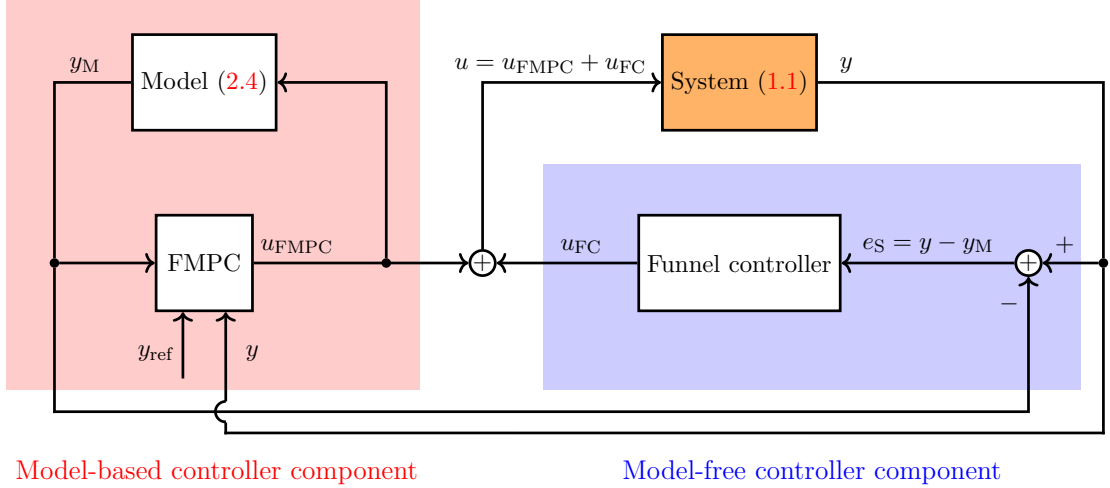

We propose \emph{robust funnel MPC}, a two component control architecture that
synergises the model-based funnel MPC~\Cref{Algo:FunnelMPC} with the
\emph{model-free} funnel controller to achieve reference tracking within
prescribed boundaries despite a mismatch between the true system \eqref{eq:Sys}
and the nominal model \eqref{eq:Model_r}. 
The overall structure is depicted in \Cref{Fig:Robust_FMPC}.
This framework addresses the inherent tension between optimality and robustness
by combining the predictive capabilities of MPC with the disturbance rejection
of adaptive feedback.

The left (red) block  of~\Cref{Fig:Robust_FMPC} comprises the 
surrogate model~\eqref{eq:Model_r}, the funnel MPC~\Cref{Algo:FunnelMPC}, 
and a given reference trajectory $y_{\rf}$.
By Theorem \ref{Thm:FMPC}, for any given funnel function $\Funnel\in \cG$,
the funnel‐MPC controller produces an input $\uFMPC$ that minimises the stage cost 
\eqref{eq:stageCostFunnelMPC} while guaranteeing the model’s output $\yM$ tracks $y_{\rf}$
within the prescribed funnel $\Funnel$, i.e.
\[
    \Norm{\eMTrack(t)} = \Norm{\yM(t) - y_{\rf}(t) } < \Funnel(t)\text{ for all }t \ge t_0.
\]
This controller component relies on the model’s accuracy but delivers optimality
by design. Since the model is chosen by the designer, it is however known
exactly. 

In contrast, the right block contains the actual system \eqref{eq:Sys} and a model-free funnel control loop 
(blue box in~\Cref{Fig:Robust_FMPC}). 
Given an arbitrary reference signal $\rho \in W^{1,\infty}(\Rp,\R^m)$ and a funnel function $\phi\in \cG$, 
the funnel control (cf. \Cref{Sec:FunnelControl}) ensures the system's output~$y$ satisfies 
$\phi(t)\Norm{y(t) - \rho(t)} < 1$ for all $t \ge t_0$, as shown in
\cite{IlchRyan02b,BergLe18a,BergIlch21}. This component requires no model
knowledge to compute the control signal $\uFC$ and track the reference with
predefined accuracy, operating purely on real-time measurements, and is
inherently robust to disturbances and system uncertainties, provided the
initial error lies within the funnel boundary.

While funnel MPC prioritises optimality through minimisation of a
designer‐specified cost, funnel control ensures robustness by adaptively
rejecting disturbances. Merging these two approaches \emph{robustifies} the funnel MPC
scheme against model uncertainties and disturbances. The combined control signal
$u = \uFMPC + \uFC$ sacrifices strict optimality (due to the corrective $\uFC$)
and model independence (due to reliance on~\eqref{eq:Model_r}) but achieves a
critical balance: the funnel controller remains dormant unless the
model-predicted error $\eMTrack$ approaches the funnel boundary $\psi$. In such
critical states -- where the model inaccuracies or disturbances threaten
constraint violation -- $\uFC$ activates to realign the system with the model's
prediction.
By keeping the funnel controller’s activation intentionally sparse, it
intervenes only as much as necessary to reject disturbances. This minimises 
deviations from the optimal control signal~$\uFMPC$. For instance, if the model
inaccurately predicts a disturbance’s impact, the funnel controller adjusts
$\uFC$ instantaneously using high-gain feedback. This ensures the system’s
output $y(t)$ adheres to constraints even when the model’s predictions $\yM$
diverge from reality. This minimal intervention strategy preserves near-optimal
performance whenever the model is accurate, while enforcing robustness in the
presence of mismatches.

Before detailing the precise interconnection and proving that the overall scheme
meets the control objective stated in \Cref{Sec:ControlObjective}, we first
outline in more detail the operating principles of the model‐free funnel
controller.
To this end, we will utilise the funnel controller
from~\cite{BergIlch21}. This controller uses error variables structurally
similar to~$\eM_i$ that we have defined in~\eqref{eq:ErrorVar} to be used by the
funnel MPC~\Cref{Algo:FunnelMPC}.
For $\phi>0$, a bijection $\FCBijec\in\cC^{1}([0,1),[1,\infty))$, $\eps\in(0,1]$, and $z=(z_1,\ldots,z_r)\in\R^{rm}$ with $z_i\in\R^m$,
we formally introduce auxiliary error variables $\eS_i$ for $i=1,\ldots, r$ in the following. Define 
\[
    \eS_1(\phi,z)\coloneqq \phi z_1,\quad
    \cE^\eps_1(\phi)\coloneqq \setdef{z\in\R^{rm}}{\Norm{\eS_1(\phi,z)}<\eps},
\]
and recursively for $z\in \cE^\eps_{i}(\phi)$ define
\begin{equation}\label{eq:ek_FC} 
    \begin{aligned}
        \eS_{i+1}(\phi,z)&\coloneqq \phi z_{i+1}+\FCBijec\rbl\Norm{\eS_i(\phi,z)}^2\rbr\eS_{i}(\phi,z),\\
        \cE^\eps_{i+1}(\phi)&\coloneqq \setdef{z\in\R^{rm}}{\Norm{\eS_i(\phi,z)}<\eps,\ j=1,\ldots, i+1},
    \end{aligned}
\end{equation}
for $i=1,\ldots,r-1$. A suitable choice for the bijection is for example $\FCBijec(s)\coloneqq 1\slash(1-s)$.
Note that in the definition of $\eS_{i}$ and $\cE^\eps_{i}(\phi)$  the value $\phi$ can be replaced  
with a time-varying function $\phi(\cdot)$ with $\phi(t)>0$ for all $t$. We will make use of this observation.

In \Cref{Sec:HighRelativeDegree}, we saw that the auxiliary error variables $\eM_{i}$ introduced in~\eqref{eq:ErrorVar}
have the property that, for $\hat{t}\geq t_0$ and a function $\zeta\in \cC^{r-1}([\hat t,\infty),\R^m)$, 
all error signals~$\eM_i(\OpChi(\zeta)(t))$ for $i=1,\ldots, r-1$ evolve within their respective funnels given by~$\Funnel_{i}$ 
if the last error variable $\eM_r(\OpChi(\zeta)(t))$ evolves within its funnel given by~$\Funnel_{r}$, see \Cref{Prop:OnlyLastFunnel}.
In the following, we show that the error variables $\eS_{i}$ in~\eqref{eq:ek_FC} exhibit a similar property.
To that end, we define, for a function $\phi\in\cG$, the set
\begin{equation}\label{eq:DefFCTrajectories}
    \FCTrajectories_{\hat{t}}\coloneqq \setdef
        {\zeta\in \cC^{r-1}(\Rp,\R^{m})}
        {
            \zeta|_{[0,t_0]}=y^0,
            \fa t\in [t_0,\hat{t}):\OpChi(\zeta)(t)\in\cEFC{1}(\phi(t))
        }.
\end{equation}
This is the set of all functions $\zeta\in \cC^{r-1}(\Rp,\R^{m})$ coinciding with 
$y^0$ and for which $\OpChi(\zeta)$ evolves within $\cEFC{1}$ on the interval $[t_0, \hat{t})$ for $\hat{t}> t_0$.
We show that all error signals $\eS_{i}(\phi(t),\OpChi(\zeta))$ for $i=1,\ldots, r-1$ evolve within $\cE^\eps_{r}(\phi(t))$ 
if the norm of the last auxiliary error $\eS_{r}(\phi(t),\OpChi(\zeta))$ remains lower than one for all $t\in[t_0,\hat{t})$ and 
if all error values $\eS_{i}$ at initial time $t_0$ are an element of $\cE^\eps_{r}(\phi(t_0))$.

\begin{lemma}\label{Lem:FCExistenceEpsMu}
    Let $\phi \in \cG$, $\FCBijec\in\cC^{1}([0,1),[1,\infty))$ be a bijection,
    and $y^0\in\cC^{r-1}([0,t_0],\R^m)$ with $\OpChi(y^0)(t_0)\in\cEFC{1}(\phi(t_0))$ be given. 
    Then, there exist constants $\eps_i,\mu_i>0$ such that for all $\hat{t}\in(t_0,\infty]$ and  all $\zeta\in \FCTrajectories_{\hat{t}}$ 
    the functions $e_i$ defined in~\eqref{eq:ek_FC} satisfy
    \begin{enumerate}[i)]
        \item $\| \eS_i(\phi(t),\OpChi(\zeta)(t))\|  \ {\leq}\ \eps_i < 1$,
        \item $\Norm{\dd{t} \eS_i(\phi(t),\OpChi(\zeta)(t)) } \ \leq\ \mu_i$,
    \end{enumerate}
    for all $t\in [t_0,\hat{t}]$ and for all $i=1,\ldots,r-1$.
\end{lemma}

\begin{proof}
    We introduce the constants $\eps_i,\mu_i$.
    Let $\eps_0=0$ and $\bar\eta_0 \coloneqq 0$.
    Utilising the bijectivity of $\FCBijec$, define successively   
    \begin{equation}
    \begin{aligned}
    \hat \eps_i &\in (0,1)  \text{ s.t. } 
        \FCBijec(\hat \eps_i^2) \hat \eps_i \geq  \SNorm{\frac{\dot \phi}{\phi} } ( 1 + \FCBijec(\eps_{i-1}^2) \eps_{i-1}) + 1 + \bar \eta_{i-1}, \\
        \eps_i &\coloneqq  \max \{ \| \eS_i(\phi(t_0),\OpChi(y_0)(t_0))\|,  \hat \eps_i\} < 1, \label{eq:ve_mu_gam} \\
        \mu_i & \coloneqq  \SNorm{\frac{\dot \phi}{\phi} } ( 1  +  \FCBijec(\eps_{i-1}^2) \eps_{i-1} )  +  1+ \FCBijec(\eps_i^2) \eps_i   +  \bar \eta_{i-1} ,\\
        \bar \eta_i & \coloneqq  2 \dot{\FCBijec}(\eps_i^2) \eps_i^2 \mu_i + \FCBijec(\eps_i^2) \mu_i, 
    \end{aligned}
    \end{equation}
    for $i={1},\ldots,r-1$.
    To improve legibility, we use the notation $\eS_i(t)\coloneqq \eS_i(\phi(t),\OpChi(\zeta)(t))$ for $\zeta \in \FCTrajectories_{\hat{t}}$.
    Let $\hat{t}\in(t_0,\infty]$ and  $\zeta\in \FCTrajectories_{\hat{t}}$ be arbitrary but fixed.
    We define the auxiliary functions $\eta_i(t) \coloneqq  \FCBijec(\|\eS_i(t)\|^2) \eS_i(t)$, and set $\eta_0(\cdot) = \dot \eta_0(\cdot) = 0$.
    To further increase readability, we omit the dependency of these functions on $t$ in the following.
    Note that, for $i=1,\ldots,r-1$, each of the error signals defined in~\eqref{eq:ek_FC} satisfies
    \[
        \dot \eS_i 
        = \dot{\phi}\zeta^{(i)}+\phi\zeta^{(i+1)}+ \dot \eta_{i-1}
        = \frac{\dot \phi}{\phi} ( \eS_i -\eta_{i-1} ) + \eS_{i+1} 
        - \FCBijec(\|\eS_i\|^2) \eS_i 
        + \dot \eta_{i-1}
    \]
    for $t\in[t_0,\hat{t})$. We observe 
    \begin{align*}
        \dot \eta_i &= 2 \dot{\FCBijec}(\| \eS_i\|^2) \al \eS_i, \dot \eS_i \ar \eS_i + \FCBijec(\| \eS_i\|^2) \dot \eS_i.
    \end{align*}
    Seeking a contradiction, we assume that, for at least one ${j \in \{1,\ldots,r-1\}}$, there exists $t^\star \in (t_0,\hat{t})$
    such that $\|\eS_j(t^\star)\|^2 > \eps_{j}$. W.l.o.g. we assume that this is the smallest possible~$j$.
    Invoking the assumption $\OpChi(y^0)\in\cEFC{1}(\phi(t_0))$ and the continuity of the involved functions, we may define
    ${t_\star \coloneqq  \max \setdef{ t \in [t_0,t^\star) }{ \| \eS_j(t) \|^2 = \eps_j}}$.
    Then,  we calculate
    \begin{align*}
        \dd{t} \tfrac{1}{2} \Norm{\eS_j}^2 
         &= \al \eS_j, \tfrac{\dot \phi}{\phi} ( \eS_j - \eta_{j-1} ) + \eS_{j+1}  + \dot \eta_{j-1} - \FCBijec(\|e_j\|^2) e_j\ar \\
         & \le \| \eS_j \| \left( \SNorm{\frac{\dot \phi}{\phi}} ( 1 + \FCBijec(\eps_{j-1}^2) \eps_{j-1}) 
         +  1   +  \bar \eta_{j-1}  - \FCBijec(\eps_j^2) \eps_j \right) 
        \le 0
    \end{align*}
    for $t \in [t_\star,t^\star]$.
    In this estimation, we used the monotonicity of $\FCBijec(\cdot)$,
    the definition of $\eps_j$, and the fact that $\dot \eta_{j - 1}$ is bounded due to the minimality of~$j$. 
    Hence, the contradiction 
    ${\eps_j < \| \eS_j(t^\star)\|^2 \leq \| \eS_j(t_\star)\|^2 = \eps_j}$
    arises after integration.
    This yields boundedness of $\eS_j, \eta_j$. 
    Using the derived bounds, we estimate
    \[
        \Norm{\dot \eS_j}  
        \le  \SNorm{\frac{\dot \phi}{\phi} } ( 1 + \FCBijec(\eps_{j-1}^2) \eps_{j-1} ) 
        + 1  + \FCBijec(\eps_j^2) \eps_j 
        + \bar \eta_{j-1}
        = \mu_j. 
    \]
    We conclude $\| \eS_i(t) \| \le \eps_i < 1$ and $\| \dot \eS_i(t)\| \le
    \mu_i$ for all $i = 1,\ldots,r-2$ and all $t \in [t_0,\hat{t}$). 
    For $i=r-1$, the same arguments are valid invoking $\eS_r : [t_0,\hat{t}) \to \cB_1$.
\end{proof}

Comparable to result in~\Cref{Cor:FMPCExistenceLambda} about the error signal $\eM_1$,
\Cref{Lem:FCExistenceEpsMu} shows that the auxiliary error signals $\eS_{i}$ for $i=1,\ldots, r-1$ maintain a 
uniform $\eps$ distance to the boundary of $\cE^\eps_{r}(\phi(t))$. 
If the initial errors are small enough, then the $\eps_i$ in \eqref{eq:ve_mu_gam}
can be chosen independent of the concrete values of $\Norm{\eS_i(\phi(t_0),\OpChi(y_0)(t_0))}$ for $i=1,\ldots, r-1$.
We summarise this in the following.

\begin{corollary}
   Let $\phi \in \cG$, $\FCBijec\in\cC^{1}([0,1),[1,\infty))$ be a bijection, and $r>1$. 
   Then, there exists $\eps\in (0,1)$ such that 
   for all $y^0\in\cC^{r-1}([0,t_0],\R^m)$ with $\OpChi(y^0)(t_0)\in\cE^{\eps}_{r-1}(\phi(t_0))$ and 
   for all $\hat{t}\in(t_0,\infty]$ every $\zeta\in \FCTrajectories_{\hat{t}}$ satisfies
   \[
        \fa t\in[t_0,\hat{t}):\quad \OpChi(\zeta)(t)\in\cE_{r-1}^\eps(\phi(t)).
   \]
\end{corollary}
\begin{proof}
    The claim immediately follows from the proof of~\Cref{Lem:FCExistenceEpsMu} by choosing $\eps$ as the minimum of 
    all $\hat{\eps}_i$ in~\eqref{eq:ve_mu_gam}.
\end{proof}

Building on \Cref{Lem:FCExistenceEpsMu}, we demonstrate that the funnel control
law from~\cite{BergIlch21} guarantees the system~\eqref{eq:Sys} tracks  a given
reference signal $\rho\in W^{r,\infty}(\Rp,\R^m)$ within predefined boundaries
governed by a function~$\phi\in\cG$. This result is generalised to accommodate
bounded disturbances $d \in  L^\infty([t_0,\infty),\R^m)$ in the input channel.
To achieve this, we leverage the \emph{perturbation high-gain property} defined
in \Cref{Def:SystemClass}~\ref{Item:PerturbationHighGain}, ensuring robustness
to such disturbances while maintaining tracking performance.

\begin{prop}\label{Prop:FC_dist}
Consider a system~\eqref{eq:Sys} with $(F,\oT) \in \cN^{m,r}_{t_0}$ as in~\Cref{Def:SystemClass}.
Let $\FCSurjec \in \cC(\Rp,\R)$ be a surjection, $\FCBijec\in\cC^{1}([0,1),[1,\infty))$ be a bijection.
Further, let the functions ${y^0\in \cC^{r-1}([0,t_0],\R^m)}$, $\rho\in W^{r,\infty}(\Rp,\R^m)$,  and 
$\phi \in \cG$ be given such that ${\OpChi(y^0 - \rho)(t_0)\in\cEFC{1}(\phi(t_0))}$ and let 
$d \in  L^\infty([t_0,\infty),\R^m)$ be an arbitrary disturbance. 
Then, the application of
\begin{equation}\label{eq:DefFC}
    u(t) = (\FCSurjec\circ\FCBijec) \rbl\Norm{\eS_{r}(\phi(t),\OpChi(e)(t))}^2\rbr\eS_{r}(\phi(t),\OpChi(e)(t)), \quad e(t) \coloneqq  y(t) - \rho(t),
\end{equation}
to the system 
\begin{equation}\label{eq:SysDisturbedControl}
    y^{(r)}(t) =  F\big(\oT(\OpChi(y)(t), d(t)+u(t)\big), \quad y|_{[0,t_0]}=y^0,
\end{equation}
yields a closed-loop initial value problem, which has a solution, every solution can be maximally extended, and every maximal solution $y: [0,\omega) \to \R^m$ has the following properties
\begin{enumerate}[label = (\roman{enumi})]
    \item\label{Item:FCGlobalSolution} the solution is global, i.e. $\omega = \infty$,
    \item\label{Item:FCControlBounded} all signals are bounded, in particular, $u\in L^\infty([t_0,\infty),\R^m)$ and $y \in W^{r,\infty}(\Rp,\R^m)$, 
    \item\label{Item:FCErrorsBounded} there exists $\eps\in (0,1)$ such that the error signals given by~$\eS_{i}$ for $i=1,\ldots, r$ as in~\eqref{eq:ek_FC} are uniformly bounded by $\eps$, i.e.
    \[
        \fa t \ge t_0  :\OpChi(y-\rho)(t)\in\cEFC{\eps}(\phi(t)).
    \]
    This implies, in particular, that the tracking error evolves within prescribed error bounds, i.e. 
    \begin{equation*}
        \fa t \ge t_0  :  \Norm{ \phi(t)(y(t) - \rho(t)) } <1.
    \end{equation*}
\end{enumerate}
\end{prop}
\begin{proof}
    We modify the proof of \cite[Thm.~1.9]{BergIlch21} to the current setting.
   
    \noindent 
    \emph{Step 1}: We show the existence of a solution of the feedback-controlled initial value problem~\eqref{eq:SysDisturbedControl} with funnel control~\eqref{eq:DefFC}.
    To this end, define the set 
    \[
        \cE\coloneqq \setdef{(t,z)\in\Rp\times\R^{rm}}{z-\OpChi(\rho)(t)\in\cEFC{1}(\phi(t))}
    \]
    where $\cEFC{1}$ is defined as in~\eqref{eq:ek_FC}.
    Moreover, formally define the function $\tilde{F}:\cE\times\R^q\to\R^{rm}$ mapping $(t,z,\eta)=(t,z_1,\ldots,z_r,\eta)$ to
    \[
        \tilde{F}(t,z,\eta)\coloneqq 
        \begin{bmatrix}
        z_2\\
        \vdots\\
        z_r\\
        F(\eta,d(t)+(\FCSurjec\circ\FCBijec) \rbl\Norm{\eS_{r}(\phi(t),z-\OpChi(\rho)(t)}^2\rbr\eS_{r}(\phi(t),z-\OpChi(\rho)(t))
        \end{bmatrix}.
    \]
    Using the notation $x(t)=\OpChi(y)(t)$,
    the initial value problem~\eqref{eq:SysDisturbedControl} with feedback control~\eqref{eq:DefFC} takes the form
    \begin{equation}\label{eq:SysWithFC}
        \dot{x}=\tilde{F}(t,x(t),\oT(x)(t)),\quad x|_{[0,t_0]}=\OpChi(y^0)\in\cC([0,t_0],\R^{rm}).
    \end{equation}
    By assumption, we have $(t_0,x(t_0))\in\cE$.
    Application of \Cref{Appendix:Th:SolutionExists} yields the existence of a maximal solution $x:[0,\omega)\to\R^{rm}$, $\omega\in(t_0,\infty]$ of \eqref{eq:SysWithFC}
    with
    \[
        \graph\rbl x|_{[t_0,\omega)}\rbr\subset \cE.
    \]
    Moreover, the closure of $\graph\rbl x|_{[t_0,\omega)}\rbr$ is not a compact subset of $\cE$.

    \noindent
    \emph{Step 2}: 
    We define several constants for later use. 
    To improve legibility, we use the notation $\eS_k(t)\coloneqq \eS_k(\phi(t),x(t)-\OpChi(\rho)(t))$ for $k=1,\ldots,r$ and $t\in[t_0,\omega)$ 
    where $\eS(\cdot,\cdot)$ is defined as in~\eqref{eq:ek_FC}.
    Further, denote with $\eS(t)\coloneqq y(t)-\rho(t)$ the tracking error between $y$ (the first $m$-dimensional component of $x$) and $\rho$.
    For the auxiliary function $\eta_k(t) \coloneqq  \FCBijec(\|\eS_k(t)\|^2) \eS_k(t)$ with $k=1,\ldots, r-1$, we observe
    \[
        \dot \eta_k = 2 \dot{\FCBijec}(\| \eS_k\|^2) \al \eS_k, \dot \eS_k \ar \eS_k + \FCBijec(\| \eS_k\|^2) \dot \eS_k,
    \]
    omitting the dependency on $t$.
    \Cref{Lem:FCExistenceEpsMu} yields the existence of  $\eps_k,\mu_k>0$ such
    that $\| \eS_k(t)\|  \ \leq\ \eps_k < 1$ and $ \Norm{\dd{t} \eS_k(t) } \
    \leq\ \mu_k$ for all $t\in[t_0,\omega)$ and all $k=1,\ldots, r-1$. Thus,
    there exists $\bar{\eta}_{r-1}\geq0$ such that
    $\Norm{\dot{\eta}_{r-1}(t)}\leq\bar{\eta}_{r-1}$ for all  $t\in[t_0,\omega)$
    (in the case $r=1$ set $\eta_0(\cdot) = \dot \eta_0(\cdot) = 0$).
    Moreover, $\| \eS_k(t)\|  \ \leq\ \eps_k$ for $k=1,\ldots r-1$ and $\| \eS_r(t)\| <1 $ for all $t\in [t_0,\omega)$
    implies the boundedness of $x(\cdot)$ in $\R^{rm}$ on the interval $[t_0,\omega)$ 
    because $\OpChi(\rho)(\cdot)$ is bounded by assumption and $\inf_{t\geq 0}\phi(t)>0$,
    see definition of $\eS_k$ in~\eqref{eq:ek_FC}.
    Thus, there exists a compact set $K_q\subset\R^q$ with $\oT(x)(t)\in K_q$
    for all $t\geq t_0$ according to the bounded-input bounded-output 
    property~\ref{Item:OperatorPropBIBO} of operator~$\oT$.
    Choose a compact set $K_m\subset\R^m$ with $d(t)\in K_m$ for all $t\geq t_0$.
    As $F$ has the perturbation high-gain property, let $\nu\in (0,1)$ such that the function 
    \[
        \HighGainFunc(s)\coloneqq 
        \min
        \setdef{\langle v, F(z, d-s v)\rangle}
        {
             d\in K_m, z \in  K_q, v\in\R^m,~\nu \leq \|v\| \leq 1
        }
    \]
    is unbounded from above, see~\Cref{Def:SystemClass}~\ref{Item:PerturbationHighGain}.
    Due to the unboundedness of the function~$\HighGainFunc$ and the surjectivity of~$\FCSurjec\circ\FCBijec$,
    it is possible to choose $\eps_r\in(0,1)$ such that 
    $\eps_r>\max\cbl{\nu,\Norm{ e_r(t_0)}}\cbr$ and  
    \begin{equation}\label{eq:DefFCEpsr}
        \tfrac{1}{2}\HighGainFunc(\FCSurjec\circ\FCBijec(\eps_r^2))\geq\theta:
        =\SNorm{\frac{\dot{\phi}}{\phi^2}}(1+\FCBijec(\eps_{r-1}^2)\eps_{r-1})+\SNorm{\frac{\bar{\eta}_{r-1}}{\phi}}+ \SNorm{\rho^{(r)}\vphantom{\frac{\dot{\phi}}{\phi}}},
    \end{equation}
    with $\eps_0=0$.
    
    \noindent
    \emph{Step 3}: We show $\Norm{\eS_r(t)}\leq\eps_r$ for all $t\in[t_0,\omega)$.
    Seeking a contradiction, assume there exists $t^\star\in[t_0,\omega)$ with $\Norm{\eS_r(t^\star)}>\eps_r$.
    Due to the continuity of $e_r$ on $[t_0,t^\star]$, there exists
    \[
        t_{\star}\coloneqq \sup\setdef{t\in[t_0,t^\star)}{\Norm{{\eS_r(t)}}=\eps_r}<t^{\star}.
    \]
    Then, we have $\Norm{\eS_r(t)}\geq\eps_r\geq\nu$ for all $t\in[t_\star,t^\star]$
    and $\HighGainFunc(\FCSurjec\circ\FCBijec(\Norm{e_r(t_{\star})}^2))\geq2\theta$.
    Thus, there exists $\tilde{t}\in [t_\star,t^\star]$ such that 
    $\HighGainFunc(\FCSurjec\circ\FCBijec(\eS_r(t)))\geq\theta$ for all $t\in[t_\star,\tilde{t}]$.
    Utilising the definition of $\eS_{r}$ in~\eqref{eq:ek_FC}, we have 
    \[
        \Norm{\eS^{(r-1)}(t)}=\Norm{\frac{1}{\phi(t)}(\eS_r(t)-\FCBijec(\eS_{r-1}^2)\eS_{r-1})}<\frac{1}{\Abs{\phi(t)}}\rbl 1+\FCBijec(\eps_{r-1}^2)\eps_{r-1}\rbr
    \]
    for all $t\in[t_\star,\tilde{t}]$ and with $\eS_0(t)=\dot{\eS}_0(t)=0$ in the case of $r=1$.
    Omitting the dependency on $t$,
    we calculate that, for almost all $t \in [t_\star,\tilde{t}]$,
    \begin{align*}
    \dd{t} \tfrac{1}{2} \Norm{ \eS_r}^2
    &= \al \eS_r, \dot{e}_r\ar\\
    &= \al \eS_r, \dot{\phi}e^{(r-1)}+\phi e^{(r)} + \dot{\eta}_{r-1}\ar\\
    &= \dot{\phi}\al \eS_r, e^{(r-1)}\ar +\phi\al\eS_r , F(\oT(x),d+u)-\rho^{(r)} \ar +\al\eS_r,\vphantom{e^{(r)}} \dot{\eta}_{r-1}\ar\\
    &\leq \Abs{\dot{\phi}}\Norm{e_r\vphantom{e^{(r)}}}\Norm{e^{(r-1)}}+\Norm{e_r\vphantom{e^{(r)}}}\Norm{\dot{\eta}_{r-1}}+ \phi\Norm{e_r\vphantom{e^{(r)}}}\Norm{\rho^{(r)}} + \phi\al\eS_r , F(\oT(x),d+u)\ar\\
    &\leq \Abs{\dot{\phi}}{\Abs{\phi}}(1+\FCBijec(\eps_{r-1}^2)\eps_{r-1})+\bar{\eta}_{r-1}+ \phi\Norm{\rho^{(r)}} + \phi\al\eS_r , F(\oT(x),d+u)\ar\\
    & \leq \phi\cdot\rbl \SNorm{\frac{\dot{\phi}}{\phi^2}}(1+\FCBijec(\eps_{r-1}^2)\eps_{r-1})+\SNorm{\frac{\bar{\eta}_{r-1}}{\phi}}+ \SNorm{\rho^{(r)}\vphantom{\frac{\dot{\phi}}{\phi}}}+ \al\eS_r , F(\oT(x),d+u)\ar\rbr\\
    & =   \phi\cdot\rbl \theta+ \al\eS_r , F(\oT(x),d+u)\ar\rbr\\
    & =   \phi\cdot\rbl \theta+ \al\eS_r , F(\oT(x),d+(\FCSurjec\circ\FCBijec)(\Norm{e_r}^2)e_r)\ar\rbr\\
    & \le \phi\cdot\rbl \theta
    - \min \setdef{ \al v, F\rbl z, d- (\FCSurjec  \circ  \FCBijec)\rbl\Norm{\eS_r}^2\rbr v\rbr \ar }
    {   \begin{array}{l}
         d \in K_{m},  \\
         z \in K_q, \\
         \nu \leq \|v\| \leq 1
    \end{array}}\rbr \\
    & \le \phi\cdot \rbl \theta-  \HighGainFunc\rbl(\FCSurjec \circ \FCBijec)\rbl\Norm{\eS_r}^2\rbr\rbr\rbr
     \le 0.
    \end{align*}
    Integration yields $\eps<\Norm{\eS_r(\tilde{t})}\leq \Norm{\eS_r(t_\star)}=\eps$, a contradiction.
    Therefore, we have $\Norm{\eS_r(t)}\leq \eps_r$ for all $t\in[t_0,\omega)$.
    
    \noindent
    \emph{Step 4}: As a consequence of \Cref{Lem:FCExistenceEpsMu} and Step~3 $\Norm{\eS_k(t)}\leq \eps_k$ for all $t\in [t_0,\omega)$ and all $k=1,\ldots, r$. 
    Choosing $\eps\in (0, 1)$ with $\eps > \eps_i$ for all $i=1,\ldots, r$ shows~\ref{Item:FCErrorsBounded}.
    By the definition of $\eS_k$ in~\eqref{eq:ek_FC} and the boundedness of the function $\OpChi(\rho)$, the solution $x$ is a bounded function, too.
    Since the closure of $\graph\rbl x|_{[t_0,\omega)}\rbr$ is not a compact subset of $\cE$,
    this implies $\omega=\infty$ and thereby shows \ref{Item:FCGlobalSolution}.
    Further, $\Norm{\eS_r(t)}\leq \eps_r<1$ implies the boundedness of $u$ in~\eqref{eq:DefFC}.
    Together with the definition of $y$ as the first $m$-dimensional component of $x$, see~\Cref{Def:FCSolutionConcept}, shows~\ref{Item:FCControlBounded}
    and completes the proof.
\end{proof}

\begin{remark}\label{Rem:FCSurjectionKnownDirection}
    The perturbation high-gain property~\ref{Item:PerturbationHighGain} holds for $F\in\cC(\R^q \times \R^m,\R^m)$
    if, and only if, for every compact set $K_m \subset \R^m$ there exists $\nu\in(0,1)$ such that, for every 
    compact set  $K_q\subset\R^q$, the function $\HighGainFunc$ defined in~\eqref{eq:Def:HighGainFunction} fulfils
    \[
        \sup_{s>0}\HighGainFunc(s)=\infty \quad\text{ or }\quad\sup_{s<0}\HighGainFunc(s)=\infty.
    \]
    If $\sup_{s>0}\HighGainFunc(s)=\infty$ for such $K_m$, $\nu$ and $K_q$, then we say that $F$ has the 
    \emph{negative-definite perturbation high-gain property} 
    (respectively, \emph{positive-definite perturbation high-gain property} if $\sup_{s<0}\HighGainFunc(s)=\infty$). 
    If it is a priori known that the negative-definite perturbation high-gain property holds for $F$,
    then the surjection $\FCSurjec$   in~\eqref{eq:DefFC} can be replaced by any surjection $\Rp\to [0,\infty)$. 
    The simplest example is the identity map $s\mapsto s$. 
    The feedback law~\eqref{eq:DefFC} then takes the form $u(t)=\FCBijec(\Norm{\eS_r(t)}^2)\eS_r(t)$, 
    where $\eS_r(t)=\eS_{r}(\phi(t),\OpChi(e)(t)$.
    Similarly, if $F$ has the positive-definite perturbation high-gain property, then the surjection $\FCSurjec$  
    in~\eqref{eq:DefFC} can be replaced by an arbitrary surjection $\Rp\to (-\infty,0]$.
\end{remark}

\Cref{Prop:FC_dist} demonstrates that applying the funnel controller~$\uFC$
(as defined in~\eqref{eq:DefFC}) to the system~\eqref{eq:Sys} forces 
the system's output~$y$ to track any given reference signal ${\rho\in W^{r,\infty}(\Rp,\R^m)}$
within given accuracy bounds governed by a function~$\phi\in\cG$.
The funnel controller generates its control signal~$\uFC$
solely from instantaneous measurements of the error
signal~$\eS_r(t)=\eS_r(\phi(t),\OpChi(y-\rho)(t))$ 
and requires no model information or look-ahead.
However, since the controller lacks predictive capacities, it may yield
suboptimal tracking performance or excessive control effort over extended
horizons. 
Crucially, naively deploying the \emph{same} reference signal~$y_{\rf}$
and funnel function~$\Funnel$ for both the model-based (MPC) and the model-free component (funnel control)
risks rendering the MPC signal~$\uFMPC$ a disruptive disturbance to the funnel controller.

To leverage model‐based prediction while retaining the funnel’s robustness, we
propose a refined integration of the funnel MPC~\Cref{Algo:FunnelMPC} and the
funnel controller~\eqref{eq:DefFC}. Instead of sharing $y_{\rf}$ and~$\Funnel$,
the MPC's predicted model output~$\yM$ serves as a reference signal for the
funnel controller. As depicted in \Cref{Fig:Robust_FMPC}, the combined
controller structure operates as follows:
\begin{itemize}
    \item \textbf{Funnel MPC (red box)}: Computes the control signal $\uFMPC(t)$  and the
   corresponding model output $\yM(t)$ over the intervals $[t_k,t_{k+1})$ with
   $t_{k}\in t_0+\delta\N_{0}$ and~$\delta>0$. 
   \item \textbf{Funnel controller (blue box)}: Receives~$\yM$ as its reference, ensuring the system output~$y$ tracks $\yM$
   with prescribed accuracy:
   \[
        \|\eSTrack(t)\| = \phi(t)\| y(t) - \yM(t) \| < 1.
   \]
\end{itemize}
The control signal applied to the system then is $u=\uFMPC+\uFC$.
The combined controller leverages the strengths of both components in a complementary framework:
\begin{enumerate}
    \item \textbf{Model accuracy}: When the model output $\yM$ aligns perfectly with the system output $y$,
    the funnel controller remains inactive ($\uFC=0$), as the MPC-generated control signal $\uFMPC$ alone achieves tracking within prescribed boundaries:
    \[
        \| y(t) - y_{\rf}(t) \| = \| \yM(t) - y_{\rf}(t) \| < \Funnel(t).
    \]
    Here, the MPC’s predictive planning dominates, optimising performance over the horizon without requiring corrective intervention.
    \item \textbf{Model uncertainty}: Under discrepancies between the model and system, the funnel controller dynamically compensates.
    The tracking error~$\eSTrack\geq0$ activates~$\uFC$, ensuring robustness by enforcing $\phi(t)\| y(t) - \yM(t) \| < 1$.
    The magnitude of~$\uFC$ scales intuitively with the model mismatch -- greater
    deviations demand stronger corrective action, while closer alignment shifts
    dominance to~$\uFMPC$.
\end{enumerate}
This dynamic interaction between the components creates a synergetic
self-regulating control hierarchy: The MPC component provides optimal foresight,
minimising control effort and improving long-term tracking and the funnel
controller acts as a safety layer, guaranteeing transient performance and
stability despite uncertainties. Utilising different reference signals ($\yM$
for the funnel controller vs. $y_{\rf}$ for the funnel MPC component), the
design avoids conflict, ensuring~$\uFMPC$ enhances -- rather than disrupts --
the funnel controller’s corrective role.

\subsection{Funnel boundary and proper initialisation}\label{Sec:RobustFMPCInit}
The funnel controller~\eqref{eq:DefFC} permits the utilisation of quite general boundary functions $\phi: \Rp \to \Rpp$.
We design $\phi$ to ensure that the feedback controller not only compensates for
model-plant mismatch $\eSTrack=y-\yM$ but also guarantees that the system output~$y$  
tracks the given reference signal~$\yM$ within the predefined error bound $\psi$ imposed 
on the MPC component. To achieve this, we propose
\begin{equation}\label{eq:DefPhi}
    \phi(t) \coloneqq  \frac{1}{\Funnel(t) - \| \yM(t) - y_{\rf}(t)\|}
\end{equation}
motivated by the following rationale:
\begin{itemize}
    \item If the MPC component ensures accurate reference tracking (i.e. $\yM \approx y_{\rf}$), then
    the boundary function for the funnel controller is $\phi \approx 1/\psi$.
    This corresponds to a ``safe'' scenario where larger deviations between
    the system $y$  and model $\yM$ are permissible.
    \item In safety-critical situations ($\yM$ deviates significantly
    from~$y_{\rf}$), $\phi$ adaptively tightens the funnel for the model-free  
    controller component, forcing the system to mimic the model and $y$ to closely follow
    $\yM$. This ensures the MPC’s optimal control input affects both dynamics
    comparably, preventing $\uFMPC$ from acting as a disturbance to the funnel
    controller.
\end{itemize}
Crucially, deviations between $\yM$ and $y_{\rf}$ are evaluated relative to the current funnel width~$\Funnel$:
Smaller $\Funnel$ tolerates less absolute deviation between the system and the
model and heightens sensitivity to mismatches, while larger $\Funnel$ permits
greater flexibility.
The function $\phi$ in \eqref{eq:DefPhi} inherently scales the allowable deviation in relationship to $\Funnel$.
This proposed design ensures that the total tracking error~$e=y-y_{\rf}$ satisfies
\[
    \Norm{e}=\Norm{y-\yM+\yM-y_{\rf}}\leq\underbrace{\Norm{\eSTrack}}_{<1/\phi}+\Norm{\eMTrack}<\Funnel-\Norm{\eMTrack}+\Norm{\eMTrack} =\Funnel,
\]
where time arguments are omitted for clarity.

In the following, we discuss mathematical difficulties arising from this
particular choice of funnel $\phi$  and reference $\yM$. A notable initial
concern is the potential discontinuity of $\yM$ (and consequently of the
function $\phi$) due to the model's re-initialisation in
Step~\ref{agostep:FunnelMPCFirst} of~\Cref{Algo:FunnelMPC}. At each time
$t_k=t_0+\delta\N_0$, the model's initial state in~\eqref{eq:Model_r} is set to
$\InitStateK_{k}$, which may introduce jumps in the concatenated trajectory $\yM$.
While~\Cref{Prop:FC_dist} assumes continuity of $\phi$ and $\rho$, this
discontinuity is largely a technical nuance.
However, careful initialisation of the combined controller is critical to ensure
compatibility between the funnel MPC and funnel controller component.

To preserve the feasibility of Algorithm~\ref{Algo:FunnelMPC} (as established in
Section~\ref{Sec:InitialRecFeasibilty}), the initial model state  $(\xM^k,\oTM^k)=\InitStateK_{k}$ 
at time $t_k$ must be an element of $\InitValues(t_k)$.
In particular, this implies 
\[
    \xM^k(t_k) - \OpChi(y_{\rf})(t_k)\in\cD_{t_k}^{\Psi},
\]
as per~\Cref{Rem:InitialValueInFunnel}.
Beyond this constraint, the MPC component permits considerable freedom in
selecting the initial state $\InitStateK_{k}$.

To maximise the effectiveness of the funnel MPC component, we want to achieve the
control objective of tracking the reference signal $y_{\rf}$ 
primarily through the (piecewise) optimal MPC control signal $\uFMPC$, with 
ideally minimal funnel controller interventions to correct deviations between
the system output $y$ and the model output $\yM$.
The model’s re-initialisation by $\InitStateK_{k}$ at each time $t_k=t_0+\delta\N_{0}$ is pivotal 
for maintaining a small model-system mismatch. 
A sophisticated initialisation strategy, leveraging system output measurements
is therefore advisable. Let $y(t_k)$ denote the system output
(from~\eqref{eq:Sys}) and $\yM^k$ the model output at time
$t_k=t_0+\delta\N_{0}$ after initialisation with $\InitStateK_k$, i.e. the first
$m$-dimensional component of $\xM^k(t_k)$. For the funnel
controller~\eqref{eq:DefFC} to function correctly when applied to
system~\eqref{eq:Sys} and tracking a given reference~$\rho$ within
boundaries~$\phi$, \Cref{Prop:FC_dist} requires 
\[
    \OpChi(y-\rho)(t_k)\in \cEFC{1}(\phi(t_k)).
\]
This restricts potential choices for the initialisation of the model.
The primary mathematical difficulty however lies in ensuring that the funnel
controller component remains uniformly bounded on the entire interval $[t_0,\infty)$.
Crucially, the maximal control input of~\eqref{eq:DefFC} depends on the maximal
value of the error variables~$\eS_i$ as defined in~\eqref{eq:ek_FC} for
$i=1,\ldots, r$. As we choose $\rho$ to be the model's output $\yM$ and $\phi$
according to \eqref{eq:DefPhi}, these error variables are in a sense
``re-initialised'' with every initialisation of the model~\eqref{eq:Model_r}.
While the funnel controller guarantees the boundedness of these error signals
between every iteration of the MPC loop, the initialisation of the model with
value~$\InitStateK_k$ has to ensure that the values $\eS_i(t_k)$ remain
uniformly bounded over all time instants $t_k= t_0+\delta\N_0$.
For the combined controller, this poses the condition  
\[
    \OpChi(y-\yM^k)(t_k)\in\cEFC{\eps}\rbl\frac{1}{\psi(t_k) - \| \yM^k - y_{\rf}(t_k)\|}\rbr
\]
for some $\eps\in(0,1)$.
The maximal control input  moreover depends on $\SNorm{\yM^{(r)}}$,
$\SNorm{\tfrac{1}{\phi}}$, and $\SNorm{\vphantom{\tfrac{1}{\phi}}\dot{\phi}}$,
see proof of~\Cref{Prop:FC_dist}.
For systems of order $r=1$, boundedness of $\tfrac{\dot{\phi}}{\phi^2}$ 
instead of~$\dot{\phi}$ suffices, see definition of $\eps_r$ in the aforementioned proof. 
The boundedness of $\yM^{(r)}$ directly follows from $\SNorm{\uFMPC}\leq\umax$ and~\Cref{Lemma:DynamicBounded}.
Moreover,
\[
    \SNorm{1\slash\phi}=\SNorm{\psi-\Norm{\yM-y_{\rf}}}\leq \SNorm{\psi}+\SNorm{\yM-y_{\rf}}\leq 2\SNorm{\psi}<\infty
\]
since $\yM$, $\psi$ and $y_{\rf}$ are bounded.
For systems of order $r>1$, we additionally have to ensure the existence of some
$\lambda\in(0,1)$ with
\begin{equation}\label{eq:FMPCLambdaDistance}
    \fa t\in [t_0,\infty):\quad \Norm{\yM(t)-y_{\rf}(t)}<\lambda\Funnel(t) <\Funnel(t)
\end{equation}
in order to guarantee the uniform boundedness of $\dot{\phi}$. 
While Theorem~\ref{Thm:FMPC} only mandates $\Norm{\yM(t)-y_{\rf}(t)}<\Funnel(t)$ for all $t\in[t_0,\infty)$,
\Cref{Cor:FMPCExistenceLambda} confirms that \eqref{eq:FMPCLambdaDistance} holds  provided 
the model~\eqref{eq:Model_r} is initialised with a sufficient distance from the funnel boundary, i.e. 
$\yM^k$  fulfils~\ref{eq:FMPCLambdaDistance} at each time instant $t_k\in t_0+\delta\N_{0}$.

The following definition formalises the requirements for initialising the model~\eqref{eq:Model_r} in the combined controller (see Figure~\ref{Fig:Robust_FMPC}).

\begin{definition}[Proper initial values $\PropInitValues(\hat{t},\hat{x})$]
    \label{Def:ProperInitValues}
    Let $y_{\rf}\in W^{r,\infty}(\Rp,\R^{m})$, $\tau\geq0$,  ${\eps,\lambda\in(0,1)}$,
    and ${\Psi=(\Funnel_1,\ldots,\Funnel_r)\in\FunnelBoundaryFuncs}$.
    Given the system data $\hat{x}\in\R^{rm}$, 
    we define the set of \emph{proper ($\eps$, $\lambda$)-initial values} for the model~\eqref{eq:Model_r}
    at time $\hat{t}\geq t_0$ as 
    \[
        \PropInitValues(\hat{t},\hat{x})\! \coloneqq \!
        \setdef
            {
               (\xMh,\oTMh) \in\InitValues(\hat{t})
            }
            {
            \!\!
            \begin{array}{ll}
                \Norm{\hat{x}_{{\mathrm{M}},1}(\hat{t})-y_{\rf}(\hat{t})}<\lambda\cdot\Funnel_1(\hat{t}),\\
                \hat{x}-\xMh(\hat{t})\in\cEFC{\eps}\rbl 1\slash\rbl\Funnel_1(\hat{t})-\Norm{\hat{x}_{{\mathrm{M}},1}(\hat{t})-y_{\rf}(\hat{t})}\rbr\rbr
            \end{array}\!\!\!
            }.
    \]
    We call $\InitState\in\PropInitValues(\hat{t},\hat{x})$ a \emph{proper ($\eps$, $\lambda$)-initialisation} at time $\hat{t}$ given 
    system data $\hat{x}\in\R^{rm}$.
\end{definition}

    By system data $\hat{x}\in\R^{rm}$ in~\Cref{Def:ProperInitValues}, we mean the measurement of the system output and its derivatives at time~$\hat{t}$, i.e.
    we will replace $\hat{x}$ later with $\OpChi(y)(\hat{t})$ where $y$ is the output of the system~\eqref{eq:Sys}.
    Further note that we implicitly allow $\lambda=1$ for systems with order $r=1$ according to our considerations regarding the boundedness of~$\tfrac{\dot{\phi}}{\phi^2}$.
    
\begin{remark}\label{Rem:PropInitialValuesNonEmpty}
    For $\hat{x}\in\R^{rm}$ with 
    $\hat{x}-\OpChi(y_{\rf})(\hat{t})\in\cEFC{\eps}(1/\Funnel_1(\hat{t}))$,
    the set  $\PropInitValues(\hat{t},\hat{x})$ is non-empty since  the pair 
    $(\OpChi(y_{\rf})|_{I_{0}^{\hat{t},\tau}},\oTM(\OpChi(y_{\rf}))|_{I_{t_0}^{\hat{t},\tau}})$
    is an element of $\InitValues(\hat{t})$, see~\Cref{Rem:InitialValuesNotEmpty}.
\end{remark}

According to~\Cref{Thm:RecursiveInitialValues}, the state of the 
model~\eqref{eq:Model_r} from the previous iteration of the funnel MPC loop
can be used to re-initialise the model at every time instant $t_k\in t_0+\delta\N$, see also~\Cref{Rem:FMPCAlternativeFirstStep}.
We will see in the proof of~\Cref{Thm:RobustFMPC} that it is possible to
operate the MPC component of the  combined controller as depicted
in~\Cref{Fig:Robust_FMPC} also in such an ``open-loop fashion'', meaning that no data
from the system is handed over to the MPC.
To be a bit more precise, we will recursively prove that, during the operation of the combined controller,  
the state of the model, when initialised with $\InitStateK_k\in\PropInitValues(t_k,\hat{x})$ at time $t_k$, 
is an element of the set $\PropInitValues(t_{k+1},\hat{x}_{k+1})$ at the next time instant $t_{k+1}$, where 
$\hat{x}_{k}\coloneqq \OpChi(y)(t_k)$ and $\hat{x}_{k+1}\coloneqq \OpChi(y)(t_{k+1})$ are the measurements of the output $y$ of the system~\eqref{eq:Sys}
at the respective time instants.
In short: 
\[ (\xM^k|_{[t_{k+1}-\tau,t_{k+1}]\cap[0,t_k]}, \oTM(\xM^k))|_{[t_{k+1}-\tau,t_{k+1}]\cap[t_0,t_{k+1}]})\in\PropInitValues(t_{k+1},\hat{x}_{k+1}),
\]
where $\xM^k\coloneqq \xM(\cdot;t_k,\InitStateK, u_k)$ is the solution of the model differential equation~\eqref{eq:Model_r}
with initial data $\InitStateK$ on the time interval $[t_k,t_k+T]$
when control $u_k\in\cU_{[t_k,t_k+T]}(\umax,\InitStateK_k)$ is applied to it.
When the computing capacity are limited, applying the combined controller  with
the model predictive control component operating in an open-loop fashion is a
simple way of potentially improving the performance of the funnel
controller~\eqref{eq:DefFC} without sacrificing speed and ease of implementation
as it is possible to pre-compute the MPC's control signal $\uFMPC$ in this case.

However, initialising the model predictive controller component with system measurement 
data sets the control algorithm on a foundation that reflects the current state of the real system~\eqref{eq:Sys}.
Such initialisation is therefore crucial to reduce prediction errors made by the model predictive controller component, 
to minimise the impact of the model-plant mismatch, and to improve the performance of the combined controller.
If the system and the model are of order $r=1$, then it is always possible to find an initialisation $\InitStateK$ at time $t_k$ 
such that the model output coincides with the system output. 
To see this, we assume for now that the combined controller as depicted in~\Cref{Fig:Robust_FMPC} achieves the  
control objective as laid out in~\Cref{Sec:ControlObjective} (we will prove in \Cref{Thm:RobustFMPC} that this is actually the case).
Let $y$ be the output of the system~\eqref{eq:Sys}. Then, 
\[
    \Norm{y(t)-y_{\rf}(t)}<\Funnel(t)
\]
for all $t\in [t_0,t_k]$. Thus, $y$ can be extended to a function $\tilde{y}\in\FunnelTrajectories_{t_k}$ with $\tilde{y}|_{[t_0,t_k]}=y$. 
This implies  $(\tilde{y}|_{I_{0}^{\hat{t},\tau}},\oTM(\tilde{y}))|_{I_{t_0}^{\hat{t},\tau}})\in\InitValues(t_k)$ where 
$I_{t_0}^{t_k,\tau}\coloneqq [\hat{t}-\tau,t_k]\cap[t_0,t_k]$. The function $\tilde{y}$ fulfils 
both 
\begin{align*}
    \Norm{\tilde{y}(t_k)-y_{\rf}(t_k)}&=\Norm{y(t_k)-y_{\rf}(t_k)}<\Funnel(t_k)
\shortintertext{and}
    \Norm{y(t_k)-\tilde{y}(t_k)}&=0<\frac{\eps}{\Funnel(t_k)-\Norm{y(t_k)-y_{\rf}(t_k)}} 
\end{align*}
for all $\eps\in(0,1)$.
It is therefore an element of the set $\PropInitValues(\hat{t},y)$ for $\lambda=1$ and all $\eps\in(0,1)$ 
(note that we allow $\lambda=1$ in the case $r=1$), see \Cref{Def:ProperInitValues}.
It is therefore possible to initialise the model with
$(\tilde{y}|_{I_{0}^{\hat{t},\tau}},\oTM(\tilde{y}))|_{I_{t_0}^{\hat{t},\tau}})$
at time $t_k$ and the model output then coincides with the system output.

For systems of higher order, it is in general not possible to initialise the model 
such that $\OpChi(\xM)(t_k)$ coincides with the system's measurement data $\OpChi(y)(t_k)$ at time of initialisation $t_k\in t_0+\delta\N_0$. 
We illustrate this in the following example.
\begin{example}
Consider a scalar system of order $r=2$. The control objective is to track the
constant reference trajectory $y_{\rf}(t) \equiv 0$ within constant boundaries
given by the funnel function~$\Funnel \equiv 1$. 
With the bijection $\FCBijec(s)\coloneqq 1\slash(1-s)$ for the funnel controller component, 
the combined controller utilises the error variables given in~\eqref{eq:ErrorVar},~\eqref{eq:ek_FC}
\begin{equation*}
    \begin{aligned}
        \xi_1(\OpChi(\yM-y_{\rf})) &= \yM ,  &&\xi_2(\OpChi(\yM-y_{\rf})) = \yMd + k \yM  , \\
        e_1(\phi,\OpChi(y-\yM)) &= \phi\cdot(y - \yM), && e_2(\phi,\OpChi(y-\yM)) = \phi \cdot\left( (\dot y -  \yMd) + \tfrac{y - \yM}{1-\|e_1\|^2} \right) 
    \end{aligned}
\end{equation*}
with  $\phi(t) =  \frac{1}{\Funnel(t) - \| \yM(t) - y_{\rf}(t)\|}$.
As parameters for the funnel MPC algorithm, we choose the constants $\FunDeriv= 1$, $\FunDiam = 1/6$,
and $k_1= 2+ \FunDeriv= 3$, and the auxiliary funnel $\Funnel_2=\tfrac{\FunDiam}{\FunDeriv}$.
Further, assume the system measurement $\OpChi(y) = (y(\hat{t}),\dot y(\hat{t})) = (2/3,0)$ at time $\hat{t}\geq t_0$.
When initialising the model with this measurement, i.e. $\OpChi(\yM)(\hat{t})=(\yM(\hat{t}),\yMd(\hat{t}))\coloneqq (2/3,0)$,
we have $\phi(\hat{t})=1/(\Funnel(\hat{t})-\Abs{\yM(\hat{t})-y_{\rf}(\hat{t})})=3$. Moreover, 
$e_1(\phi(\hat{t}),\OpChi(y-\yM)(\hat{t}))=0$ and $e_2(\phi(\hat{t}),\OpChi(y-\yM)(\hat{t}))=0$.
Thus, $\OpChi(y-\yM)(\hat{t})\in \cEFC{\eps}(\phi(\hat{t}))$ for all $\eps\in(0,1)$.
For the auxiliary variable~$\xi_1$, we have $\Abs{\xi_1(\OpChi(\yM-y_{\rf}))} = \Abs{y(\hat{t})}=2/3<1=\psi(\hat{t})$.
However, 
\[
    \xi_2(\OpChi(\yM-y_{\rf})) = \yMd(\hat{t})+ k_1\yM(\hat{t})= k_1y(\hat{t})=2>1/6=\Funnel_2(\hat{t}).
\]
This means $(2/3,0)\notin\cD_{t}^{\Psi}$.
Therefore, there exists no element of
$\xMh\in\PropInitValues(\hat{t},\OpChi(y))$ coinciding with $\OpChi(y)$ at time
$\hat{t}$, i.e. $\xMh(\hat{t})=\OpChi(y)$.
\end{example}

Just as there exist a multitude of possibilities to initialise the funnel MPC~\Cref{Algo:FunnelMPC}
via an initialisation strategy as defined in \Cref{Def:InitialisationStrategy}, 
there are also many conceivable methods to select a proper ($\eps$, $\lambda$)-initialisation
$\InitState\in\PropInitValues(\hat{t},\hat{x})$ given measurements 
$\hat{x}\coloneqq\OpChi(y)(\hat{t})$ at time $\hat{t}$.
A versatile strategy is solving an optimisation problem of the form 
\begin{equation}\label{eq:MinimizePropInit}
    \mathop {\operatorname{minimise}}_{\substack{(\xMh,\oTMh)\in\PropInitValues(\hat{t},\hat{x})}}  J(\hat{t},\hat{x},\xMh,\oTMh), 
\end{equation}
where $J$ is a cost function  that takes the desired aspects into account. For
example,  as it is in general not possible to find initialisation $\InitStateK$
at time $t_k$ such that the model output $\OpChi(\xM)(t_k)$ coincides with the
system output $\OpChi(y)(t_k)$, one could instead minimise the euclidean
distance between the two vectors. Another possibility would be to give more
weight to the lower derivatives, as these are presumably less affected by
disturbances. A large number of potential approaches are conceivable, which can
be described by such an optimisation problem.

While many MPC schemes assume access to the full
system state, we consider scenarios where only output measurements $\OpChi(y)$ are
available. To address the challenge of state estimation in uncertain or
disturbed linear discrete-time systems, a Luenberger observer was employed to
reconstruct the system state in the works~\cite{mayne2009robust,kogel2017robust}. By
integrating this observer with a tube-based MPC framework, the control scheme
ensures robust constraint satisfaction and preserves recursive feasibility. This
approach demonstrates how observer-based strategies can effectively compensate
for state unavailability whilst maintaining closed-loop performance.
Clearly, the employment of methods beyond the Luenberger observer like 
moving horizon estimation (MHE)~\cite{Haseltine2005} or non-linear state
observers \cite{Besancon2007,Korder2022} is also conceivable.
Similarly, observers can be leveraged to estimate the internal state of the system 
and thus find more suitable initial states~$\oTMh$ of the model, i.e. initial values 
for the operator~$\oTM$. 
While our analysis is indifferent with regard to the selected initial value, it
is clear that the performance of the model predictive component may
significantly be improved by accurate estimates of~$\oTM$.
The deployment of state observers is particularly well-suited to our problem
setting when the structure of the model in~\eqref{eq:Model_r} aligns with the
dynamics of the physical system described in~\eqref{eq:Sys}. 

\subsubsection{Activation function}
Minor deviations between the system output $y(t)$ and the predicted model output
$\yM(t)$ are often negligible in practice, posing no risk of violating  the
funnel boundaries $\Funnel$.
This is inherently addressed by the design of the  function~$\phi$
in~\eqref{eq:DefPhi}, as $\phi\approx 1/\Funnel$ when $\yM \approx y_{\rf}$.
From an application standpoint, it may seem advantageous to fully ``deactivate'' the
funnel feedback controller during nominal operation and only engage it in
safety-critical scenarios.
To this end, we highlight the option of incorporating an \emph{activation function}, i.e. 
a continuous function $\ActivFunc: [0,1] \to [0,\ActivFunc^+]$, $\ActivFunc^+ > 0$,
with $\ActivFunc(1)=\ActivFunc^+$
into the funnel controller. This continuous function modulates the control
signal $\uFC$ based on the magnitude of the error $\eS_r$, effectively scaling the gain
term  $(\FCSurjec\circ\FCBijec)$ in the control law~\eqref{eq:DefFC}. Crucially,
while $\ActivFunc(\cdot)$ adjusts the gain magnitude, the adaptive gain
mechanism remains unaffected -- ensuring it retains the necessary magnitude to
enforce error bounds. The use of such an activation function is rigorously
justified by the following theoretical result.

\begin{lemma}\label{lemma:ActivationSurjective}
    Let $\FCSurjec\in\cC(\Rp,\R)$ be a surjection, $\FCBijec\in\cC([0,1),[1,\infty))$ be a bijection,
    and $\ActivFunc\in\cC([0,1],[0,\ActivFunc^+])$ be an activation function with $\ActivFunc^+>0$
    and $\ActivFunc(1)=\ActivFunc^+$.
    Then, the function ${\tilde{\FCSurjec}\coloneqq (\ActivFunc\circ\sqrt{\FCBijec^{-1}})\cdot \FCSurjec\in\cC(\Rp,\R)}$ is surjective.
\end{lemma}
\begin{proof}
$\FCSurjec\in\cC(\Rp,\R)$ being a surjection is equivalent to 
$\limsup_{s\to \infty} \FCSurjec(s)=\infty$ and $\liminf_{s\to \infty} \FCSurjec(s)=-\infty$.
Since $\lim_{s\to \infty} (\ActivFunc\circ\sqrt{\FCBijec^{-1}})(s)=\ActivFunc^+>0$, we have
\[
    \limsup_{s\to \infty} \tilde{\FCSurjec}(s)=\infty\quad \text{ and }\quad \liminf_{s\to \infty} \tilde{\FCSurjec}(s)=-\infty.
\]
This implies that $\tilde{\FCSurjec}=(\ActivFunc\circ\sqrt{\FCBijec^{-1}})\cdot \FCSurjec\in\cC(\Rp,\R)$ is surjective as well.
\end{proof}
A reasonable and simple choice for an activation function can be
    \begin{equation*}
        \ActivFunc(s) = \begin{dcases}
            0, & 0\leq s \le S_{\mathrm{crit}}, \\
            s-S_{\mathrm{crit}}, &  S_{\mathrm{crit}} \leq s\leq 1,
        \end{dcases}
    \end{equation*}
for $S_{\mathrm{crit}} \in (0,1)$. In this particular case we may set $\ActivFunc^+ = 1-S_{\mathrm{crit}}$.
In the context of machine learning, in particular, artificial neural networks,
this type of functions is known as \emph{rectified linear unit} (ReLU), see e.g.~\cite{ramachandran2017searching} and references therein.
Note that $\ActivFunc$ defined above satisfies $\ActivFunc(S_{\mathrm{crit}}) = 0$,
whereby it is a continuous function and thus the funnel controller contributes continuously to the overall control signal. 

\Cref{lemma:ActivationSurjective} shows that, instead of control law~\eqref{eq:DefFC}, it is possible to use the funnel controller $\uFC$
with an activation function $\ActivFunc$ in~\Cref{Prop:FC_dist}, i.e. the control law 
\[
    u(t) = \ActivFunc(\Norm{\eS_{r}(t)})\cdot(\FCSurjec\circ\FCBijec) \rbl\Norm{\eS_{r}(t)}^2\rbr\eS_{r}(t),
\]
where $\eS_{r}(t)\coloneqq \eS_{r}(\phi(t),\OpChi(y-\rho)(t))$.
In fact, a such scaled funnel controller has already been a potential controller candidate since its development in~\cite{IlchRyan02b}.
However, most examples in the literature utilise the functions $\FCBijec(s)=1/(1-s)$ and $\FCSurjec(s)=s\sin(s)$ for the control law
($\FCSurjec(s)=\pm s$ in case of a known control direction, see~\Cref{Rem:FCSurjectionKnownDirection}). 
To our knowledge,~\cite{BergDenn24b} was the first work to explicitly mention the possibility to ``deactivate'' the funnel controller for small error signals.

\subsection{The robust funnel MPC algorithm}
We now consolidate our findings into the  robust funnel
MPC~\Cref{Algo:RobustFMPC}, formally defining the controller structure
illustrated in \Cref{Fig:Robust_FMPC}. Building on the definitions,
concepts, and results established thus far, we prove that this scheme is
initially and recursively feasible and that its application to the
model~\eqref{eq:Model_r} solves the tracking problem formulated in
\Cref{Sec:ControlObjective}. In particular, the scheme guarantees that
the deviation between the system output $y$ and a given reference
signal~$y_{\rf}\in W^{r,\infty}(\Rp,\R^{m})$ evolves within the
funnel~$\cF_{\Funnel}$ defined by a function~$\Funnel\in\cG$.

\begin{algo}[Robust funnel MPC]\label{Algo:RobustFMPC}\ \\
    \textbf{Given:}\\[-4ex] %
    \begin{itemize}%
        \item instantaneous measurements of the output $y$ and its derivatives of system~\eqref{eq:Sys},
            initial time $t^0\in\Rp$, initial trajectory $y^0\in\cC^{(r-1)}([0,t_0],\R^m)$, 
            reference signal $y_{\rf}\in W^{r,\infty}(\Rp,\R^{m})$, 
            funnel function $\psi\in\cG$.
        \item model~\eqref{eq:Model_r}, signal memory length~$\tau\geq0$,
        auxiliary funnel boundary function ${\Psi=(\psi_1,\ldots,\psi_r)\in\FunnelBoundaryFuncs}$ with corresponding parameters $k_i$ for $i=1,\ldots, r$,
        input saturation level $\umax\geq0$, and funnel stage cost function~$\FunnelStageCost$,
        \item initialisation parameters $\eps,\lambda\in(0,1)$,
        \item a surjection $\FCSurjec \in \cC(\Rp,\R)$ and  a bijection $\FCBijec\in\cC([0,1), [1,\infty))$.
    \end{itemize} 
    \textbf{Set} the time shift $\delta >0$, 
                 the prediction horizon $T\geq\delta$, and index $k\coloneqq 0$.\\
    \textbf{Define} the time sequence~$(t_k)_{k\in\N_0} $ by $t_k \coloneqq  t_0+k\delta$.\\ 
    \textbf{Steps:}
    \begin{enumerate}[(a)]
    \item \label{agostep:RobustFMPCFirst} 
    Obtain a measurement $\hat{x}_k\coloneqq \OpChi(y)(t_k)$ of the system output~$y$ and its derivatives at the current time $t_k$ 
    and  choose a  \emph{proper ($\eps$,$\lambda$)-initialisation} $\InitStateK_{k}\in\PropInitValues(t_k,\hat{x}_k)$ for the model.
    \item \textbf{\textsc{Funnel MPC}}\\ Compute a solution $\uFMPCk\in L^\infty([t_k,t_k +T],\R^{m})$ of the optimal control problem
    \begin{equation}\label{eq:RobustFMPCOCP}
        \mathop
                {\operatorname{minimise}}_{\substack
                {
                    u\in L^{\infty}([t_k,t_k+T],\R^{m}),\\
                    \SNorm{u}  \leq \umax 
                }
            }\      \int_{t_k}^{t_k + T}\FunnelStageCost(s,\eM_{r}(\xM(s;t_k,\InitStateK_{k},u)-\OpChi(y_{\rf})(s)),u(s))\d{s}.
    \end{equation}
    Predict the output~$\yM^k(t;t_k,\InitStateK_k,\uFMPCk)$ of the model on the 
    interval~{$[t_k,t_{k+1}]$}, and define the adaptive funnel $\phi_k: [t_k,t_{k+1}]\to\Rpp $ by 
    \begin{equation} \label{alg:eq:vp}
        \phi_k(t)\coloneqq \frac{1}{\Funnel_1(t)-\Norm{\eMTrack^k(t)}},
    \end{equation}
    where $\eMTrack^k(t) = \yM^k(t) - y_{\rf}(t)$.
    \item\label{alg:step:FC} \textbf{\textsc{Funnel control}}\\
    Using the error variables $\eS_i$ for $i=1,\ldots, r$ as in~\eqref{eq:ek_FC}, define the funnel control law~$\uFC$
    with reference $\yM^k$ and funnel function~$\phi_k$ as in~\eqref{alg:eq:vp} by
    \begin{equation}\label{eq:uFCRobustFMPC}
        \uFCk(t) \coloneqq   (\FCSurjec\circ\FCBijec) (\Norm{\eS_{r}(\phi_{k}(t),\eSTrack(t))}^2)\eS_{r}(\phi_{k}(t),\eSTrack(t)),
    \end{equation}
    with $\eSTrack(t) =y(t)-\yM^k(t)$.
    Apply the control law
    \begin{equation}\label{eq:uRobustFMPC}
        u_k:[t_k,t_{k+1})\to\R^m, \ u_k(t) 
        = \uFMPCk(t)+ \uFCk(t)
    \end{equation}
    to system~\eqref{eq:Sys}.
    Increment $k$ by $1$ and go to Step~\ref{agostep:RobustFMPCFirst}.
    \end{enumerate}
\end{algo}

\begin{remark}
    \Cref{Algo:RobustFMPC} integrates the funnel MPC~\Cref{Algo:FunnelMPC} (from
    \Cref{Chapter:FunnelMPC}) with the model-free funnel controller
    of~\cite{BergIlch21} via Step~\ref{alg:step:FC}. By employing the model
    output $\yM$ as the reference signal for the funnel controller, the combined
    scheme leverages the MPC’s predictive capabilities even in safety-critical
    scenarios, while ensuring the MPC’s optimal control input $\uFMPC$ enhances
    -- rather than disrupts -- the funnel controller’s operation. Coupled with
    the funnel function~$\phi$ (computed using MPC predictions), this guarantees
    the tracking error remains within the prescribed performance funnel
    $\Funnel$, as formalised in~\Cref{Thm:RobustFMPC}.
    The principal mathematical challenges involve
    ensuring that the funnel MPC algorithm remains feasible under 
    ($\eps$,$\lambda$)-initialisation of the model based on system output measurements.
    To this end, we adapt the results from~\cite{BergIlch21} (resp. \Cref{Prop:FC_dist}) to the current setting.
    However, the findings in~\cite{BergIlch21} cannot be directly applied since the
    reference signal for the funnel controller is assumed to be a priori given
    and to be continuous -- conditions violated in~\Cref{Algo:RobustFMPC} due to the MPC-generated
    reference $\yM$.
\end{remark}

\begin{theorem}\label{Thm:RobustFMPC}
    Consider a system~\eqref{eq:Sys} with $(F,\oT) \in \cN^{m,r}_{t_0}$ as in~\Cref{Def:SystemClass} 
    and choose a model~\eqref{eq:Model_r} with $(\fM,\gM,\oTM) \in \cM^{m,r}_{t_0}$ as in~\Cref{Def:ModelClass}.
    Let ${t_0\geq 0}$ be the initial time and let $y_{\rf} \in W^{r,\infty}(\Rp,\R^m)$ 
    and $\Psi=(\Funnel_1,\ldots,\Funnel_r)\in \FunnelBoundaryFuncs$ be given and 
    let $\tau\geq0$ be greater than or equal to the memory limit of operator~$\oTM$.
    Further, let $y^0\in\cC^{(r-1)}([0,t_0],\R^m)$ with ${\OpChi(y_0-y_{\rf})(t_0)\in\cEFC{1}(1/\Funnel(t_0))}$ be the initial trajectory for the system~\eqref{eq:Sys}.
    Then, there exist $\eps,\lambda\in(0,1)$ ($\lambda =1$ in the case $r=1$), and  $\umax\geq0$ such that the robust funnel MPC \Cref{Algo:RobustFMPC} 
    with $\delta>0$ and $T\ge\delta$ is initially and recursively feasible,
    i.e. at every time instant $t_k \coloneqq  t_0+k\delta $ for $k\in\N_0$
    \begin{itemize}
        \item there exists a proper initialisation $\InitStateK_{k}\in\PropInitValues(t_k,\hat{x}_k)$ and
        \item the OCP~\eqref{eq:RobustFMPCOCP} has a solution $\uFMPCk\in L^\infty([t_k,t_k+T],\R^m)$.
    \end{itemize}
    Moreover,  the closed-loop system consisting of the system~\eqref{eq:Sys} and the feedback law~\eqref{eq:uRobustFMPC} 
    has a global solution $y : [0,\infty) \to \R^m$. 
    Each global solution $y$ satisfies that
    \begin{enumerate}[label = (\roman{enumi}), ref=(\roman{enumi})]
        \item \label{Assertion:learning_y_u_bounded}
        all signals are bounded, in particular, $u\in L^\infty([t_0,\infty),\R^m)$ and $y \in W^{r,\infty}(\Rp,\R^m)$, 
        \item \label{Assertion:learning_tracking_error}
        the tracking error between the system's output and the reference evolves within prescribed boundaries, i.e.
        \begin{equation*}
            \fa t \ge t_0 : \Norm{y(t)  - y_{\rf}(t)}< \Funnel_1(t) .
        \end{equation*}
    \end{enumerate}
\end{theorem}
\begin{proof}
    \emph{Step 1}: We define the constants $\lambda$ and $\umax$.
    To that end, set $\lambda =1$ if the order of the model~\eqref{eq:Model_r} is $r =1$.
    Otherwise, choose $\lambda\in(0,1)$ such that 
    for all $s>\hat t\ge t_0$  every function $\zeta\in \cC^{r-1}([\hat t,\infty),\R^m)$  with $\OpChi(\zeta)(t)\in\cD_{t}^\Psi$
    and $\Norm{\eM_1(\OpChi(\zeta)(\hat{t}))}<\lambda \cdot\Funnel_1(\hat{t})$ fulfils 
    \[
        \Norm{\eM_1(\OpChi(\zeta)(t))}<\lambda\cdot\Funnel_1(t)
    \]
    for all $t\in[\hat{t},s]$.
    Here, $\eM_1$ is the first auxiliary error variable used in the funnel MPC~\Cref{Algo:FunnelMPC} as introduced in~\Cref{Sec:HighRelativeDegree}.
    A constant $\lambda$ with this properties exists according to~\Cref{Cor:FMPCExistenceLambda}.
    Further, choose  $\umax\geq0$ such that $\Controls(\umax,\InitState)\neq\emptyset$
    for all $\hat{t}\geq t_0$, $\InitState\in\InitValues(\hat{t})$, and $T>0$.
    Such bound $\umax\geq 0$ exists according to~\Cref{Th:ExUmax}.

    \noindent
    \emph{Step 2}:
    Similarly to~\Cref{Lem:FCExistenceEpsMu}, we define several constants for later use.
    By assumption, we have $\OpChi(y_0-y_{\rf})(t_0)\in\cEFC{1}(1/\Funnel_1(t_0))$. 
    Thus, there exists $\bar{\eps}\in(0,1)$ with ${\OpChi(y_0-y_{\rf})(t_0)\in\cEFC{\bar{\eps}}(1/\Funnel_1(t_0))}$.
    In the case of $r>1$, define  
    \[
        \bar{\phi} \coloneqq  2\SNorm{\Funnel_1} \text{ and }\hat{\phi}\coloneqq \tfrac{\SNorm{\dot{\Funnel}_1}
            +\SNorm{\vphantom{\dot{\Funnel}}\Funnel_2}
            +k_1\SNorm{\vphantom{\dot{\Funnel}_1}\Funnel_1} }{\rbl(1-\lambda)\inf_{s\geq0}\Funnel_1(s)\rbr^2},
    \]
    where $k_1\geq 0$ is the first parameter corresponding to the auxiliary funnel function~$\Psi\in\FunnelBoundaryFuncs$.
    Let $\eps_0=0$ and $\bar\eta_0 \coloneqq 0$.
    Utilising the bijectivity of $\FCBijec$, define successively   
    \begin{equation}
    \begin{aligned}
    \hat \eps_i &\in (0,1)  \text{ s.t. } 
       \FCBijec(\hat \eps_i^2) \hat \eps_i =  \frac{\hat{\phi}}{\bar{\phi}}  ( 1 + \FCBijec(\eps_{i-1}^2) \eps_{i-1}) + 1 + \bar \eta_{i-1}, \\
       \eps_i &\coloneqq  \max \{ \bar{\eps},  \hat \eps_i\} < 1, \label{eq:Robust_FMPC_ve_mu_gam} \\
       \mu_i & \coloneqq  \frac{\hat{\phi}}{\bar{\phi} } ( 1  +  \FCBijec(\eps_{i-1}^2) \eps_{i-1} )  +  1+ \FCBijec(\eps_i^2) \eps_i   +  \bar \eta_{i-1} ,\\
       \bar \eta_i & \coloneqq  2 \dot{\FCBijec}(\eps_i^2) \eps_i^2 \mu_i + \FCBijec(\eps_i^2) \mu_i, 
    \end{aligned}
    \end{equation}
    for $i={1},\ldots,r-1$. 

    \noindent
    \emph{Step 3}: We define $\eps\in(0,1)$. 
    To that end, define the set 
    \[
        \cE\coloneqq \setdef{(t,z)\in\Rp\times\R^{rm}}{z-\zeta(t)\in\cEFC{1}(1/\Funnel_1(t)),\zeta\in\FunnelTrajectories_{\infty} }, 
    \] 
    with $\FunnelTrajectories_{\infty}$ as in~\eqref{eq:Def:FunnelTrajectories}.
    According to the proof of \Cref{Lemma:DynamicBounded}, 
    there exists a compact set $\hat{K}\subset\R^{rm}$ with 
    \[
        \fa \zeta\in\FunnelTrajectories_{\infty}\fa t\geq 0:\quad \zeta(t)\in \hat{K},
    \]
    see~\eqref{eq:ExistenceCompactsetAllFunnelTrajectories}.
    Thus, the set $\cE$ is bounded.
    Due to the bounded-input bounded-output property~\ref{Item:OperatorPropBIBO} in~\Cref{Def:OperatorClass},
    the operator~$\oT$ is bounded for all functions $\zeta\in\cR(\Rp, \R^{rm})$ evolving within $\cE$, 
    see also the definition of the set $\cEFC{1}(1/\Funnel_1(t))$~in~\eqref{eq:ek_FC}.
    Hence, there exists a compact set $K$ with $\oT(\zeta)\subset K$ for all $\zeta\in\cR(\Rp, \R^{rm})$ evolving within $\cE$.
    As $F$ has the perturbation high-gain property, let $\nu\in (0,1)$ such that the function 
    \[
        \HighGainFunc(s)\coloneqq 
        \min
        \setdef{\langle v, F(z, d-s v)\rangle}
        {
             d\in \bar{\cB}_{\umax}, z \in  K_q, v\in\R^m,~\nu \leq \|v\| \leq 1
        }
    \]
    is unbounded from above, see~\Cref{Def:SystemClass}~\ref{Item:PerturbationHighGain}.
    Due to the unboundedness of the function~$\HighGainFunc$ and the surjectivity of~$\FCSurjec\circ\FCBijec$ it is possible to choose $\eps_r\in(0,1)$ such that 
    $\eps_r>\max\cbl{\nu,\Norm{ e_r(t_0)}}\cbr$ and  
    \begin{equation}\label{eq:DefEpsr}
        \tfrac{1}{2}\HighGainFunc(\FCSurjec\circ\FCBijec(\eps_r^2))\geq\theta:
        =\frac{\hat{\phi}}{\bar{\phi}^2}(1+\FCBijec(\eps_{r-1}^2)\eps_{r-1})+\frac{\bar{\eta}_{r-1}}{\bar{\phi}}+ \fMmax+\gMmax\umax,
    \end{equation}
    where $\fMmax$ and $\gM^{\max}$ are the constants from~\Cref{Lemma:DynamicBounded}.
    In the case $r=1$, replace $\frac{\hat{\phi}}{\bar{\phi}^2}$ with $\|\dot{\Funnel}_1\|+\SNorm{\dot{y}_{\rf}}+\fMmax+\gMmax\umax$ in \eqref{eq:DefEpsr}.
    Choose $\eps\in(0,1)$ with $\eps>\eps_i$ for all~$i=1,\ldots,r$.
    
    \noindent
    \emph{Step 4}:
    Let $\delta>0$ and $T\ge\delta$ be arbitrary but fixed.
    When applying the robust funnel MPC~\Cref{Algo:RobustFMPC} to the system~\eqref{eq:Sys}, the system's dynamics 
    on each interval $[t_k,t_{k+1}]$ for $t_k= t_0+k\delta$ and $k\in\N_0$
    are given by
    \begin{equation}\label{eq:SysPiecewise}
       y^{(r)}_k(t)=F(\oT(\OpChi(y_k))(t),u_k(t)),\quad y_k|_{[0,t_k]}= y_{k-1}|_{[0,t_k]}
    \end{equation}
    where $y_{-1}\coloneqq y^0$ and $u_k$ is the control given by \eqref{eq:uRobustFMPC}.
    Note that ${\OpChi(y_k)(t_k)=\OpChi(y_{k-1})(t_k)}$.
    In the following, we show via induction that the robust funnel MPC~\Cref{Algo:RobustFMPC} is initially and recursively feasible.
    This means, in particular, that there exists a proper initialisation $\InitStateK_{k}\in\PropInitValues(t_k,\OpChi(y_{k-1})(t_k))$
    at every time instant $t_k$, that $u_k$ as in~\eqref{eq:uRobustFMPC} is well defined on every interval $[t_k,t_{k+1}]$,
    and that~\eqref{eq:SysPiecewise} has a maximal solution $y_k$ defined on the entire interval $[t_k,t_{k+1}]$.
    
    \noindent
    \emph{Step 4.1}:
    When obtaining the measurement of the system's output and its derivatives at the initial time $t_0$ in Step~\ref{agostep:RobustFMPCFirst}
    of the robust funnel MPC~\Cref{Algo:RobustFMPC}, we have ${\hat{x}_0=\OpChi(y)(t_0)=\OpChi(y^0)(t_0)}$.
    The construction of $\eps$, which is larger or equal to $\bar{\eps}$, yields
    ${\hat{x}_0-\OpChi(y_{\rf})(t_0)\in\cEFC{\eps}(1/\Funnel_1(t_0))}$.
    Thus, the set $\PropInitValues(t_0,\OpChi(y^0)(t_k))$ of proper $(\eps,\lambda)$ initial values is non-empty according to~\Cref{Rem:PropInitialValuesNonEmpty}.

    \noindent
    \emph{Step 4.2}:
    Let $y_{k-1}$ be a solution of~\eqref{eq:SysPiecewise} defined on the interval $[0,t_{k-1}]$ for some $k\in\N_0$. 
    Note that $y_{-1}=y^0$ for $k=0$.
    Let $\hat{x}_k\coloneqq \OpChi(y_{k-1})(t_k)$ be the system's output~$y_{k-1}$ and its derivatives at time instant $t_k$.
    Further, assume that there exists a proper initialisation $(\xMh^k,\oTMh^k)\coloneqq \InitStateK_k\in\PropInitValues(t_k,\hat{x}_k)$.
    We show that the control signal~$u_k$ as in~\eqref{eq:uRobustFMPC} is well-defined and that 
    when applying $u_k$ to the system~\eqref{eq:Sys} the initial value problem~\eqref{eq:SysPiecewise} has a solution $y_k:[t_k,t_{k+1}]\to\R^m$.
    As $\PropInitValues(t_k,\hat{x}_k)\subset\InitValues(t_k)$, the choice of $\umax\geq0$ ensures the non-emptiness  
    of the set $\cU_{[t_k,t_k+T]}(\umax,\InitStateK_k)$.
    Therefore, \Cref{Th:SolutionExists} yields the existence of a solution $\uFMPCk\in\cU_{[t_k,t_k+T]}(\umax,\InitStateK_k)$
    of the OCP~\eqref{eq:RobustFMPCOCP}.
    Let $\yM^k(\cdot;t_k,\InitStateK_k,\uFMPCk):[t_k,t_{k+1}]\to\R^m$ be the corresponding output of the model~\eqref{eq:Model_r} 
    when applying the control $\uFMPCk$ with initial time $t_k$ and initial value $\InitStateK_k$ over the time interval~$[t_k,t_{k+1}]$.
    Note that $\yM^k$ is, in fact, defined on the whole interval  $[0,t_k+T]$ according to the solution 
    concept for the model differential equation~\eqref{eq:Model_r}, see~\Cref{Def:ModSolution}.
    Moreover, $\yM^k$ restricted to the interval $[t_k,t_{k+1}]$ is an element of~$W^{r,\infty}([t_k,t_{k+1}],\R^m)$.
    By ${\uFMPCk\in\cU_{[t_k,t_k+T]}(\umax,\InitStateK_k)}$,
    we have $\Norm{\yM^k(t)-y_{\rf}(t)}<\Funnel_1(t)$ for all $t\in[t_k,t_{k+1}]$.
    Thus, the function~$\phi_k:[t_k,t_{k+1}]\to\Rpp$, $\phi_k(t)=1/(\Funnel_1(t)-\Norm{\yM^k(t)-y_{\rf}(t)})$ 
    in~\eqref{alg:eq:vp} is well defined.
    $\phi_k$ is bounded with a bounded derivative due to the compactness of the interval~$[t_k,t_{k+1}]$.
    Note that $\hat{x}_k-\xMh^k\in\cEFC{\eps}(\phi_0(t_k))$ because $\InitStateK_k$ is a proper initial value, see~\Cref{Def:ProperInitValues}.
    Applying the control signal $u_k$ as in~\eqref{eq:uRobustFMPC} consisting of sum of $\uFMPCk$ 
    and the funnel control signal $\uFCk$ as in~\eqref{eq:uFCRobustFMPC} with reference $\yM^k$ and funnel function~$\phi_k$
    to the system~\eqref{eq:Sys} with initial value $y_k|_{[0,t_k]}= y_{k-1}|_{[0,t_k]}$ to the loop system~\eqref{eq:SysPiecewise}.
    This initial value problem has a solution $y_k:[t_k,t_{k+1}]\to\R^m$, see~\Cref{Prop:FC_dist}.
    
    \noindent
    \emph{Step 4.3}:
    Assuming the existence of a proper initialisation $\InitStateK_k\in\PropInitValues(t_k,\hat{x}_k)$,
    we show certain bounds for $\yM^k$ and $\phi_k$ on the interval $[t_k,t_{k+1}]$.
    $\OpChi(\yM^k-y_{\rf})(t)\in\cD_{t}^\Psi$ for all $t\in[t_k,t_{k+1}]$ because $\uFMPCk\in\cU_{[t_k,t_k+T]}(\umax,\InitStateK_k)$, 
    see definition of $\cU_{[t_k,t_k+T]}(\umax,\InitStateK_k)$ in~\eqref{eq:Def-U}.
    Since
    \[
        \yM^{(r)}(t)=\fM(\OpChi(\yM)(t))+\gM(\OpChi(\yM)(t))\uFMPCk(t)
    \]
    for $t\in[t_k,t_{k+1}]$, the function $\yM^{(r)}$ is bounded on the interval $[t_k,t_{k+1}]$ by $\fMmax+\gMmax\umax$,
    see~\Cref{Lemma:DynamicBounded}.
    We observe 
    \[
        \Norm{1\slash\phi_k(t)}=\Norm{\Funnel_1(t)-\Norm{\yM^k(t)-y_{\rf}(t)}}\leq \SNorm{\Funnel_1}+\SNorm{\yM^k-y_{\rf}}\leq 2\SNorm{\Funnel_1}=\bar{\phi}
    \]
    on the interval $[t_k,t_{k+1}]$.
    Moreover, if the order of the system is $r=1$, then 
    \[
        \Norm{\tfrac{\dot{\phi_k}(t)}{\phi_k(t)^2}}
        \leq\Norm{\dot{\psi}(t)}+\Norm{\vphantom{\dot{\psi}}\yMdk(t)}+\Norm{\vphantom{\dot{\psi}}\dot{y}_{\rf}(t)}
        \leq\SNorm{\vphantom{\dot{\psi}}\psi}+\fMmax+\gMmax\umax+\SNorm{\vphantom{\dot{\psi}}\dot{y}_{\rf}}
    \]
    for almost all $t\in[t_k,t_{k+1}]$.
    As $(\xMh^k,\oTMh^k)\coloneqq \InitStateK_k\in\PropInitValues(t_k,\hat{x}_k)$, we have
    \[ 
        \Norm{\eM_1(\xMh^k(t_k)-\OpChi(y_{\rf})(t_k))}=\Norm{\hat{x}_{{\mathrm{M}},1}^k(t_k)-y_{\rf}(t_k)}=\Norm{\yM^k(t_k)-y_{\rf}(t_k)}<\lambda\cdot\Funnel_1(t_k).
    \]
    If the order of the system is $r>1$, this yields $\Norm{\yM^k(t)-y_{\rf}(t)}<\lambda\cdot\Funnel_1(t)$ for all $t\in[t_k,t_{k+1}]$
    due to the choice of $\lambda$, see~\Cref{Cor:FMPCExistenceLambda}.
    Since $\OpChi(\yM^k-y_{\rf})(t)\in\cD_{t}^\Psi$ for all $t\in[t_k,t_{k+1}]$, 
    \[
        \Norm{\yMdk(t)-\dot{y}_{\rf}(t)}=\Norm{\eM_2(\OpChi(\yM^k-y_{\rf})(t))-k_1\eM_1(\OpChi(\yM^k-y_{\rf})(t))}\leq \SNorm{\Funnel_2}+k_1\SNorm{\Funnel_1},
    \]
    where $k_1\geq 0$ is the parameter corresponding to the auxiliary error variable $\eM_2$,
    see definition of $\eM_i$ in~\eqref{eq:ErrorVar}.
    Therefore, 
    \[
        \Norm{\dot{\phi}_k(t)}\leq\frac{\Norm{\dot{\Funnel}_1(t)}+\Norm{\vphantom{\dot{\Funnel_1}}\yMdk(t)-\dot{y}_{\rf}(t)}}{\rbl\Funnel_1(t)-\Norm{\yM^k(t)-y_{\rf}(t)}\rbr^2}
        \leq\frac{\SNorm{\dot{\Funnel}_1}+\SNorm{\vphantom{\dot{\Funnel_1}}\Funnel_2}+k_1\SNorm{\vphantom{\dot{\Funnel_1}}\Funnel_1}}{\rbl(1-\lambda)\inf_{s\geq0}\Funnel(s)\rbr^2},
    \]
    for almost all $t\in[t_k,t_{k+1}]$.
    Note that the derived boundaries for $\yM^{(r)}$, $1/\phi_k$,  $\dot{\phi}_k$, and $\tfrac{\dot{\phi}_k}{\phi_k^2}$ 
    are independent of the time instant $t_k$ and the particular choice of $\InitStateK_k\in\PropInitValues(t_k,\hat{x}_k)$.
    
    \noindent
    \emph{Step 4.4}:
    We show that if $\PropInitValues(t_k,\hat{x}_k)$ is non-empty, then 
    $\PropInitValues(t_{k+1},\hat{x}_{k+1})$ is non-empty after applying a control~$u_k$ as in~\eqref{eq:uRobustFMPC} to the system~\eqref{eq:Sys},
    where $\hat{x}_k\coloneqq \OpChi(y_{k-1})(t_k)$ and $\hat{x}_{k+1}\coloneqq \OpChi(y_{k})(t_k)$.
    Let $(\xMh^k,\oTMh^k)\coloneqq \InitStateK_k\in\PropInitValues(t_k,\hat{x}_k)$ be an arbitrary but fixed proper initialisation.
    We have $\hat{x}_k-\xMh^k(t_k)\in\cEFC{\eps}(\phi_k(t_k))$, see \Cref{Def:ProperInitValues}.
    According to~\Cref{Prop:FC_dist}, there exists $\tilde{\eps}\in(0,1)$ with 
    \[
        \OpChi(y_k-\yM^k)(t)\in\cEFC{\tilde{\eps}}(\phi_k(t))
    \]
    for all $t\in[t_k,t_{k+1}]$.
    In the proof of~\Cref{Prop:FC_dist}, $\tilde{\eps}$ is constructed as 
    the maximum of $\eps_i$, $i=1,\ldots,r-1$ as defined in~\eqref{eq:ve_mu_gam} and $\eps_r$ in~\eqref{eq:DefFCEpsr}.
    Due to the boundaries derived in {Step 4.3}, $\eps$ as defined in Step~3 fulfils the estimates
    for~$\tilde{\eps}$ in~\eqref{eq:ve_mu_gam}.
    Regarding $\eps_r$, note the following.
    As $\yM^k$ can be extended to an element of $\FunnelTrajectories_{\infty}$, the function $y_{k}$ 
    can be extended to a function evolving within the set $\cE$.
    Thus, the bound~\eqref{eq:DefFCEpsr} for $\eps_r$ can be proven with the same calculations as in the proof of~\Cref{Prop:FC_dist}. 
    Therefore, $\eps$ as defined in Step~3 fulfils the estimates
    for~$\tilde{\eps}$ in~\eqref{eq:ve_mu_gam} and in~\eqref{eq:DefFCEpsr}.
    Or in other words, $\tilde{\eps}$ in~\Cref{Prop:FC_dist} can be chosen smaller or equal $\eps$~from Step~3 of the current proof.
    Thus, $ \OpChi(y_k-\yM^k)(t)\in\cEFC{\eps}(\phi_k(t))$ for all $t\in[t_k,t_{k+1}]$.
    In particular, $\OpChi(y_k-\yM^k)(t_{k+1})\in \cEFC{\eps}\rbl 1\slash\rbl\Funnel_1(t_{k+1})-\Norm{\yM^k(t_{k+1})-y_{\rf}(t_{k+1})}\rbr\rbr$.
    Further note that $\Norm{\yM^k(t)-y_{\rf}(t)}<\lambda\cdot\Funnel_1(t)$ for all $t\in[t_k,t_{k+1}]$
    due to the choice of $\lambda$, see~\Cref{Cor:FMPCExistenceLambda}.
    According to~\Cref{Thm:RecursiveInitialValues}, we have
    \[
        (\OpChi(\yM^k)|_{[t_{k+1}-\tau,t_{k+1}]\cap[0,t_k]},
        \oTM(\OpChi(\yM^k))|_{[t_{k+1}-\tau,t_{k+1}]\cap[t_0,t_{k+1}]})
        \in\InitValues(t_{k+1}).
    \]
    Thus, 
    \[
    (\OpChi(\yM^k)|_{[t_{k+1}-\tau,t_{k+1}]\cap[0,t_k]}, \oTM(\OpChi(\yM^k))|_{[t_{k+1}-\tau,t_{k+1}]\cap[t_0,t_{k+1}]})\in\PropInitValues(t_{k+1},\hat{x}_{k+1}).
    \]

    \noindent
    \emph{Step 4.5}:
    We sum up Step~4.
    Under the assumption that the set of proper initial values $\PropInitValues(t_k,\hat{x}_k)$ at time instant $t_k$ is non-empty, 
    we showed in Step~4.2 that one iteration of the robust funnel MPC~\eqref{Algo:RobustFMPC} can be executed.
    This means, in particular, that the optimisation problem~\eqref{eq:RobustFMPCOCP} has a solution $\uFMPCk\in\cU_{[t_k,t_k+T]}(\umax,\InitStateK_k)$,
    that the output~$\yM^k(t;t_k,\InitStateK_k,\uFMPCk)$ of the model~\eqref{eq:Model_r} exists on the 
    entire interval~{$[t_k,t_{k+1}]$}, and that the adaptive funnel $\phi_k: [t_k,t_{k+1}]\to\Rpp $ given by~\eqref{alg:eq:vp} is well-defined.
    Furthermore, applying the control $u_k$ as defined in~\eqref{eq:uRobustFMPC}
    to the system~\eqref{eq:Sys} with initial value ${y_k|_{[0,t_k]}= y_{k-1}|_{[0,t_k]}}$
    leads to the loop system which has  a maximal solution $y_k:[t_k,t_{k+1}]\to\R^m$.
    Utilising the bounds derived in Step~4.3, it was shown in Step~4.4 that
    $\PropInitValues(t_{k+1},\hat{x}_{k+1})$ is non-empty after applying a control~$u_k$ 
    as in~\eqref{eq:uRobustFMPC} to the system~\eqref{eq:Sys}
    if $\PropInitValues(t_k,\hat{x}_k)$ is non-empty.
    Step~4.1 shows that initially the set $\PropInitValues(t_0,\hat{x}_0)$ is non-empty.
    Therefore, it follows inductively that the robust funnel MPC~\eqref{Algo:RobustFMPC} can recursively be applied to the system~\eqref{eq:Sys} 
    and that the closed-loop system consisting of the system~\eqref{eq:Sys} and the control law~\eqref{eq:uRobustFMPC} has a global solution 
    $y : [0,\infty) \to \R^m$.
    
    \noindent
    \emph{Step 5}:
    Let $y : [0,\infty) \to \R^m$ be a global solution of the closed-loop system consisting of the system~\eqref{eq:Sys} and the control law~\eqref{eq:uRobustFMPC}.
    We show~\ref{Assertion:y_u_bounded} and~\ref{Assertion:tracking_error}.
    Let $\yM : [0,\infty) \to \R^m$ be the associated concatenated solution of the model differential equation~\eqref{eq:Model_r} 
    with the sequence of initial values $(t_k,\InitStateK_k)_{k\in\N_{0}}$ and control signals $\uFMPCk$. 
    Further let $\phi:[t_0,\infty)\to\Rp$ with $\phi(t)\coloneqq 1/(\Funnel_1(t)-\Norm{\yM(t)-y_{\rf}(t)})$.
    Then, $\yM|_{[t_k,t_{k+1})}=\yM^k(\cdot;t_k,\InitStateK_k,\uFMPCk)|_{[t_k,t_{k+1})}$ 
    and $\phi|_{[t_k,t_{k+1})} = \phi_k|_{[t_k,t_{k+1})}$ for all~$k\in\N_{0}$.
    We have 
    \[
        \OpChi(y-\yM)(t)\in\cEFC{\tilde{\eps}}(\phi(t))
    \]
    for all $t\in[t_0,\infty)$.
    Since $\yM$ and $\phi$ are bounded functions, $y \in W^{r,\infty}(\Rp,\R^m)$, see definition of $\cEFC{\tilde{\eps}}$ in~\eqref{eq:ek_FC}.
    The funnel MPC signal $\uFMPCk$ is bounded by $\umax$ for all $k\in\N_{0}$.
    The funnel control signal $\uFCk$ is bounded by $\HighGainFunc(\FCSurjec\circ\FCBijec(\eps^2))$ for all $k\in\N_{0}$, see definition of $\eps_r$ in~\eqref{eq:DefEpsr} 
    and the calculations in Step~3 of the proof of~\Cref{Prop:FC_dist}.
    This shows~\ref{Assertion:learning_y_u_bounded}.
    Moreover, we have   
    \[
        \Norm{y(t)-y_{\rf}(t)}\leq \Norm{y(t)-\yM(t)}+\Norm{\yM(t)-y_{\rf}(t)}<\phi(t)+\Norm{\yM(t)-y_{\rf}(t)}=\Funnel_1(t)
    \]
    for all $t\geq t_0$.
    This shows~\ref{Assertion:learning_tracking_error} and completes the proof.
\end{proof}

\begin{remark}\label{Rem:OpenLoopInit}
We comment on the difference between the proposed control scheme and a
straightforward combination of a MPC scheme with a feedback control law.
\begin{enumerate}[(a)]
    \item\label{Item:Rem:OpenLoopInit} 
    The integration of feed-forward and feedback control is a widely adopted strategy.  
    Prior work in~\cite{BergOtto19,BergDrue21} explores combining funnel control with feed-forward methods.
    Similarly, model predictive control -- specifically funnel MPC -- can be augmented with a feedback controller. 
    This approach can be implemented in the robust funnel MPC~\Cref{Algo:RobustFMPC}
    by omitting the feedback loop between the funnel MPC and the system.
    Instead, at each MPC cycle, the model is re-initialised using only the prior prediction of the model state:
    \[
    \InitStateK_{k+1}\coloneqq (\xM^k|_{[t_{k+1}-\tau,t_{k+1}]\cap[0,t_k]}, \oTM(\xM^k))|_{[t_{k+1}-\tau,t_{k+1}]\cap[t_0,t_{k+1}]}).
    \]
    This is an element of  $\PropInitValues(t_{k+1},\hat{x}_{k+1})$
    independently of $\eps$ and $\lambda$, making it a special case of a
    \emph{proper initialisation}. Here, the funnel MPC signal $\uFMPC$ can be computed
    offline via the model and applied as an open-loop control to the system.
    Concurrently, the feedback controller compensates for errors arising from
    discrepancies between the model and the physical system.
    \item An alternative to the open-loop operation of
    \Cref{Algo:RobustFMPC} involves feedback based on system output
    measurements $\hat{x}_k\coloneqq \OpChi(y)(t_k)$. By properly initialising
    the model with $\InitStateK_{k}\in\PropInitValues(t_k,\hat{x}_k)$, two objectives are achieved:
    recursive feasibility of the MPC scheme is preserved and 
    the model state $\InitStateK_{k}$ mirrors the system’s actual state $\hat{x}_k$.
    This re-initialisation at each MPC cycle incorporates the impact of the
    control signal $\uFMPC$ on the system. Furthermore, it may enhance the
    efficacy of the optimal control signal in improving the system’s tracking
    performance.
\end{enumerate}
\end{remark}

\begin{remark}
    \Cref{Thm:RobustFMPC} demonstrates that the robust funnel MPC~\Cref{Algo:RobustFMPC} is model-agnostic.
    For any system~\eqref{eq:Sys} with $(F,\oT) \in \cN^{m,r}_{t_0}$,
    the algorithm remains functional regardless of the chosen
    model~\eqref{eq:Model_r}, provided $(\fM,\gM,\oTM) \in \cM^{m,r}_{t_0}$.
    Crucially, the system and model need not share structural similarity. For instance:
    \begin{itemize}
        \item The model may be a lower-dimensional approximation of a higher-dimensional system.
        \item The model could represent a linearised version of a non-linear system.
        \item The model might omit time delay effects.
    \end{itemize}
    This flexibility ensures applicability across diverse modelling paradigms.
\end{remark}

\section{Simulation}
In this section, we revisit the numerical examples
from~\Cref{Sec:FMPC:Simulation} to illustrate the robust funnel
MPC~\Cref{Algo:RobustFMPC}.
The \textsc{Matlab} source code for the performed simulations 
can be found on \textsc{GitHub} under the link \url{https://github.com/ddennstaedt/FMPC_Simulation}.

\subsection*{Exothermic chemical reaction}
To demonstrate the application of the robust funnel MPC~\Cref{Algo:RobustFMPC} by
a numerical simulation, we consider again a continuous-time chemical reactor and
concentrate on the control goal of steering the reactor's temperature to a predefined
reference value $y_{\rf}(t)$ within boundaries given by a function~$\Funnel(t)$.
As in \Cref{Sec:FMPC:Sim:Reactor}, we consider a reactor described by  
the following non-linear system of order one:
\begin{equation}\tag{\ref{eq:ExampleExothermicReaction} revisited}
    \begin{aligned}
        \dot{x}_1(t) &= c_1\, p(x_1(t), x_2(t), y(t)) +d(x_1^{\mathrm{in}}-x_1(t)),\\
        \dot{x}_2(t) &= c_2\, p(x_1(t), x_2(t), y(t)) +d(x_2^{\mathrm{in}}-x_2(t)),\\
        \dot{y}(t)   &= b\,   p(x_1(t), x_2(t), y(t)) -q\,y(t) + u(t).
    \end{aligned}
\end{equation}
The reactor's temperature should follow a given heating profile specified in \eqref{eq:Ex:Reactor:RefTemp} within
tolerance limits defined by the funnel function $\Funnel(t)\coloneqq 20\me^{-2t}+4$.
To achieve the control objective with robust funnel MPC~\Cref{Algo:RobustFMPC}, 
we again use the strict funnel stage cost function 
$\ell_{\Funnel}:\Rp\times\R\times\R\to\R\cup\{\infty\}$ given by 
\begin{equation}\tag{\ref{eq:Ex:ReactorFMPCCostFunction} revisited}
\begin{aligned}
    \ell_{\Funnel}(t,y,u) =
    \begin{dcases}
        \frac {\Norm{y-y_{\rf}(t)}}{\Funnel(t)^2 - \Norm{y-y_{\rf}(t)}^2} + \lambda_u \Norm{u - 360}^2,
            & \Norm{y-y_{\rf}(t)} \neq \Funnel(t)\\
        \infty,&\text{else},
    \end{dcases}
\end{aligned}
\end{equation}
with design parameter ${\lambda_u\in\Rp}$.
We restrict the MPC control signal to ${\|\uFMPC \|_\infty \le 600}$.
Further, we choose the design parameters $\lambda_u =
10^{-4}$, prediction horizon $T = 1$, and time shift $\delta = 0.1$.
However, unlike before, we do not utilise the actual differential equations
describing the system~\eqref{eq:ExampleExothermicReaction} as a model for the MPC algorithm.
Instead, we consider a linearisation of this non-linear
reaction process  obtained by linearising the Arrhenius function $p(x_1,x_2,y) =
k_0 \me^{- \frac{k_1}{y}} x_1$ around the desired final temperature $\bar y =
337.1K$ and $x_1 = \tfrac12 x_1^{\mathrm{in}}$.
This results in
\begin{equation*}
    p_{\mathrm{lin}}(x_1,x_2,y) = k_0 \me^{-\frac{k_1}{\bar y}} x_1 + \frac{k_0 k_1 e^{-\frac{k_1}{\bar y}}}{\bar y^2}  \frac{x_1^{\mathrm{in}}}{2}  (y-\bar y).
\end{equation*}
Set $a_1 \coloneqq \tfrac{k_0 k_1 e^{-\frac{k_1}{\bar y}}}{\bar y^2} \tfrac{x_1^{\mathrm{in}}}{2}$, $a_2 \coloneqq k_0 e^{-\frac{k_1}{\bar y}}$ and define the expressions
\begin{equation*}
    A = \begin{bmatrix}
     c_1 a_2 - d & 0   & c_1 a_1   \\
     c_2 a_2     & -d  & c_2 a_1   \\
     b a_2       & 0   & b a_1 - q 
    \end{bmatrix}\in\R^{3\times3}, \quad
    D = \begin{bmatrix}
        -c_1 a_1 \bar y + d x_{\mathrm{M},1}^{\mathrm{in}} \\
        -c_2 a_1 \bar y + d x_{\mathrm{M},2}^{\mathrm{in}}\\
        - b a_1 \bar y 
    \end{bmatrix} \in \R^{3}.
\end{equation*}
Then, with $\xM \coloneqq [x_{\mathrm{M},1}, x_{\mathrm{M},2},\yM]^\top \in \R^3$, the model to be used in the funnel MPC controller component is given by
\begin{equation}\label{Eq:Ex:Reactor:LinearModel}
\begin{aligned}
    \xM(t) &= A \xM(t) + B \uFMPC(t) + D, \\
    \yM(t) &= C \xM(t),
\end{aligned}
\end{equation}
where  $C = B^\top = [0,0,1] \in \R^{1\times 3}$.
We choose the same parameters as in~\eqref{eq:Ex:Reactor:Sys:Params} and 
assume initial values of the system and the model to coincide, i.e.
\[
    [x_1(0),x_2(0),y(0)]^\top=[x_{{\mathrm{M}},1}(0), x_{{\mathrm{M}},2}(0),\yM(0)]^\top \coloneqq [0.02,0.9,270]^\top.
\]
Due to discretisation, we consider only step functions with a constant step length of 
$\SampleTime\coloneqq\delta= 0.1$ to solve the OCP~\eqref{eq:RobustFMPCOCP}.

\begin{figure}[h]
    \begin{subfigure}[b]{0.49\textwidth}
        \centering
        \includegraphics[width=\textwidth]{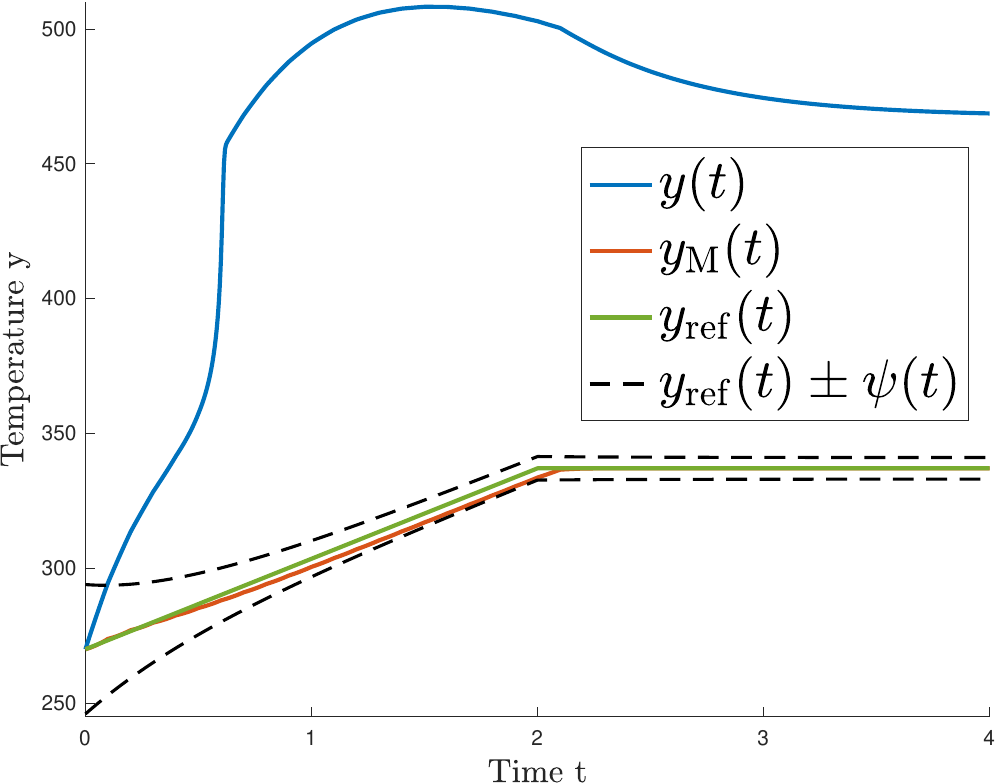}
        \caption{Outputs and reference, with boundary~$\Funnel$.}
        \label{fig:sim:robust_fmpc_nonrobust:output}
    \end{subfigure}
      \begin{subfigure}[b]{0.49\textwidth}
        \centering
        \includegraphics[width=\textwidth]{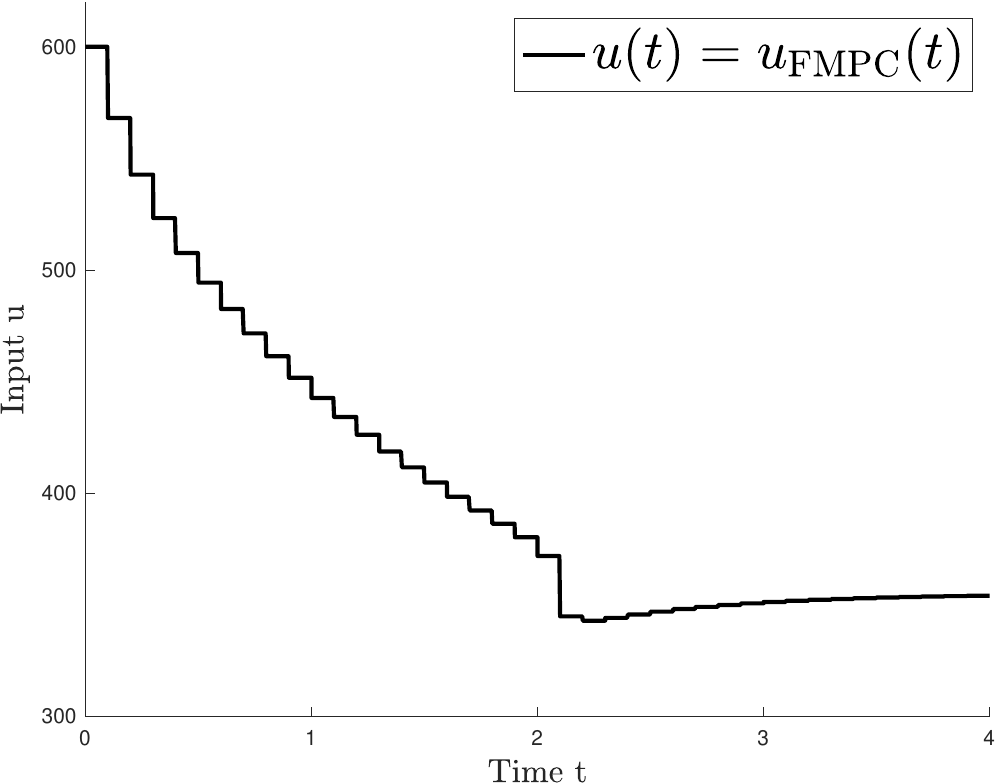}
        \caption{Control inputs.}
        \label{fig:sim:robust_fmpc_nonrobust:input}
    \end{subfigure}
    \caption{Simulation of system~\eqref{eq:ExampleExothermicReaction} under the
    control generated by the funnel MPC~\Cref{Algo:FunnelMPC} without additional
    funnel control feedback loop.}
    \label{fig:sim:robust_fmpc_nonrobust}
\end{figure}

\begin{figure}[h]
    \begin{subfigure}[b]{0.49\textwidth}
        \centering
        \includegraphics[width=\textwidth]{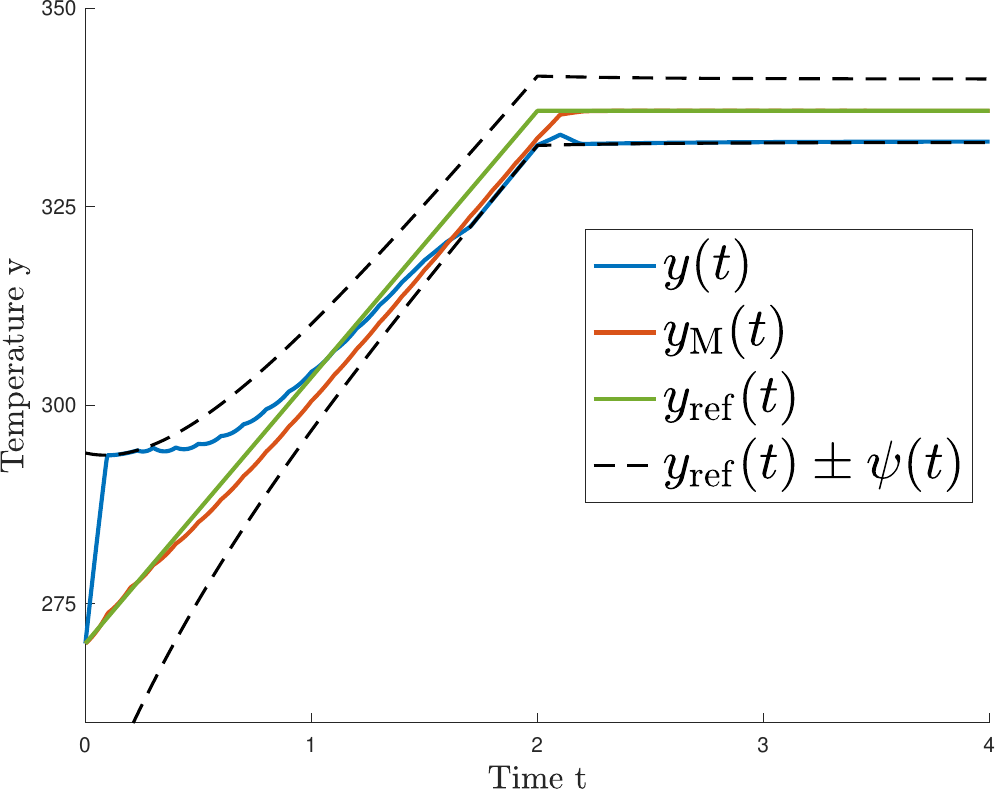}
        \caption{Outputs and reference, with boundary~$\Funnel$.}
        \label{fig:sim:robust_fmpc_trivial:output}
    \end{subfigure}
      \begin{subfigure}[b]{0.49\textwidth}
        \centering
        \includegraphics[width=\textwidth]{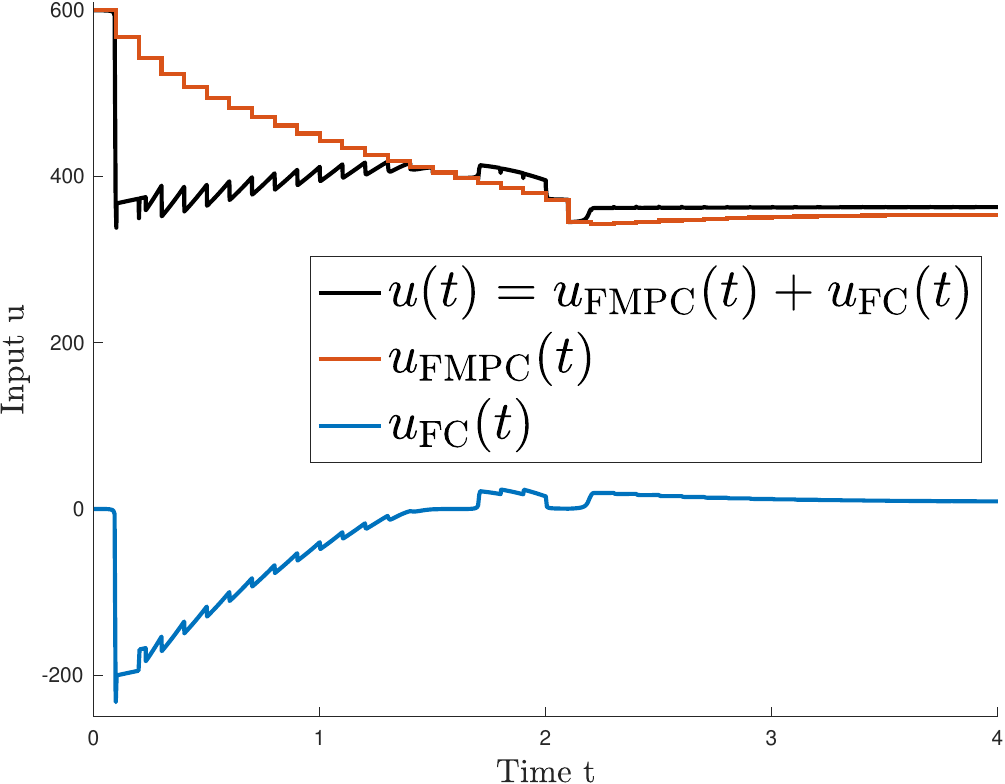}
        \caption{Control inputs.}
        \label{fig:sim:robust_fmpc_trivial:input}
    \end{subfigure}
    \caption{Simulation of system~\eqref{eq:ExampleExothermicReaction} under the
    control generated by the robust funnel MPC~\Cref{Algo:RobustFMPC} with 
    funnel control feedback, without model re-initialisation.}
    \label{fig:sim:robust_fmpc_trivial}
\end{figure}

For the control law of funnel control component, we choose the bijection $\FCBijec(s)=1/(1-s)$
and the function $\FCSurjec(s)=- s$. 
This choice for $\FCSurjec$ is justified since we assume the
control direction to be known, see~\Cref{Rem:FCSurjectionKnownDirection}. 
This assumption is also realistic from a practical point of view.
To additionally demonstrate that the funnel controller can be combined with an activation
function $\ActivFunc$, as discussed in \Cref{Sec:RobustFMPCInit}, we interconnect the
controller with a ReLU-like map
\begin{equation*}
\ActivFunc(s) = 
    \begin{cases}
        0, & s\le S_{\mathrm{crit}}, \\
        s- S_{\mathrm{crit}}, & s \ge S_{\mathrm{crit}},
    \end{cases}
\end{equation*}
where we choose $S_{\mathrm{crit}} = 0.4$. 
The funnel controller therefore is only active, if the error $e=y-\yM$
exceeds $40 \%$ of the maximal distance to its funnel boundary.
We run the simulation on an interval of $[0, 4]$ 
and consider  the following
scenarios:
\begin{itemize}
    \item \textit{Case 1:} Funnel MPC without robustification, i.e. $\uFMPC$ is
    computed via the funnel MPC~\Cref{Algo:FunnelMPC} and applied to the system
    without an additional funnel control loop. The model is initialised, at
    every iteration of the algorithm, with the model's state from the previous
    iteration. The results are shown in \Cref{fig:sim:robust_fmpc_nonrobust}.
    \item \textit{Case 2:} Robust funnel MPC with a trivial proper re-initialisation,
    i.e. model is initialised, at
    every iteration of the algorithm, with the model's state from the previous
    iteration ($\InitStateK_{k} = \xM^{k-1}(t_k)$ in Step~\ref{agostep:RobustFMPCFirst} of \Cref{Algo:RobustFMPC}).
    The results are depicted in~\Cref{fig:sim:robust_fmpc_trivial}.
    \item \textit{Case 3:} Robust funnel MPC with a proper initialisation according to the system's output, 
    i.e. $\InitStateK_{k}\in \PropInitValues(t_k,y(t_k))$ is selected such
    that $y(t_k) = \yM(t_k)$ in Step~\ref{agostep:RobustFMPCFirst} of
    \Cref{Algo:RobustFMPC}. To this end, we initialise the model \eqref{Eq:Ex:Reactor:LinearModel} 
    with  ${\xM(t_k) \coloneqq  [x_{\mathrm{M},1}^{k-1}(t_k), x_{\mathrm{M},2}^{k-1}(t_k),y(t_k)]^\top}$
    at every time instant $t_k\in\delta\N$, i.e. the states $x_{\mathrm{M},1}$ and $x_{\mathrm{M},2}$
    remain unchanged during initialisation and $\yM(t_k)$ is set to $y(t_k)$.
    The results are displayed in~\Cref{fig:sim:robust_fmpc_reinit}.
\end{itemize}
In the following figures, the control signal generated via the funnel MPC component is labelled
with the subscript~\textsc{FMCP} ($\uFMPC$); the signal generated by the additional
funnel controller component is labelled with the subscript~\textsc{FC} ($\uFC$).
\Cref{fig:sim:robust_fmpc_nonrobust} shows the application of the control
signal computed with funnel MPC~\Cref{Algo:FunnelMPC} in Case~1 to the system without
an additional funnel control feedback loop. The error $\eMTrack(t) = \yM(t) - y_{\rf}(t)$ 
between the model's output $\yM(t)$ and the
reference $y_{\rf}(t)$ evolves within the funnel boundaries $\Funnel(t)$.
However, the control signal computed with funnel MPC using the linear model is
not sufficient to achieve that the tracking error $e(t) = y(t) - y_{\rf}(t)$ of 
the non-linear system evolves within the funnel boundaries $\Funnel(t)$.
Obviously, the deviation is induced during the initial phase. After about $t=2$, 
the linearised model is a good approximation of the system. In this
region, the control $\uFMPC$ has a comparable effect on both dynamics; however, the
error $y(t)- y_{\rf}(t)$ already evolves outside the funnel boundaries $\Funnel(t)$.
\Cref{fig:sim:robust_fmpc_trivial} shows the application of the control signal
computed with robust funnel MPC~\Cref{Algo:RobustFMPC} in Case~2, i.e. besides
the funnel MPC control signal the additional funnel controller is applied in
order to guarantee that the error $y(t) - y_{\rf}(t)$ evolves within the
boundaries~$\Funnel(t)$. Since the model and the system do not coincide, the system
evolves differently from the model and hence the funnel controller has to compensate the
model-plant mismatch. 
However, after the system has reached the desired temperature
$y_{\rf,\mathrm{final}}$ at $t_{\mathrm{final}} = 2$, the system's states evolve
close to the linearisation point of the model \eqref{Eq:Ex:Reactor:LinearModel}.
Hence, the linear model closely approximates the non-linear system
\eqref{eq:ExampleExothermicReaction}. Consequently, the control signal $\uFMPC$
generated by the MPC controller component is nearly sufficient to maintain the
system output $y$ at $y_{\rf,\mathrm{final}}$ within the desired temperature
range and the funnel controller intervenes only slightly  with a small control
signal~$\uFC$.

Note that, in both Cases 1 and 2, it is possible to pre-compute the control
signal $\uFMPC$ as no system measurement data is fed back to the funnel MPC
component, see \Cref{Rem:OpenLoopInit}~\ref{Item:Rem:OpenLoopInit}.

\begin{figure}[h]
    \begin{subfigure}[b]{0.49\textwidth}
        \centering
        \includegraphics[width=\textwidth]{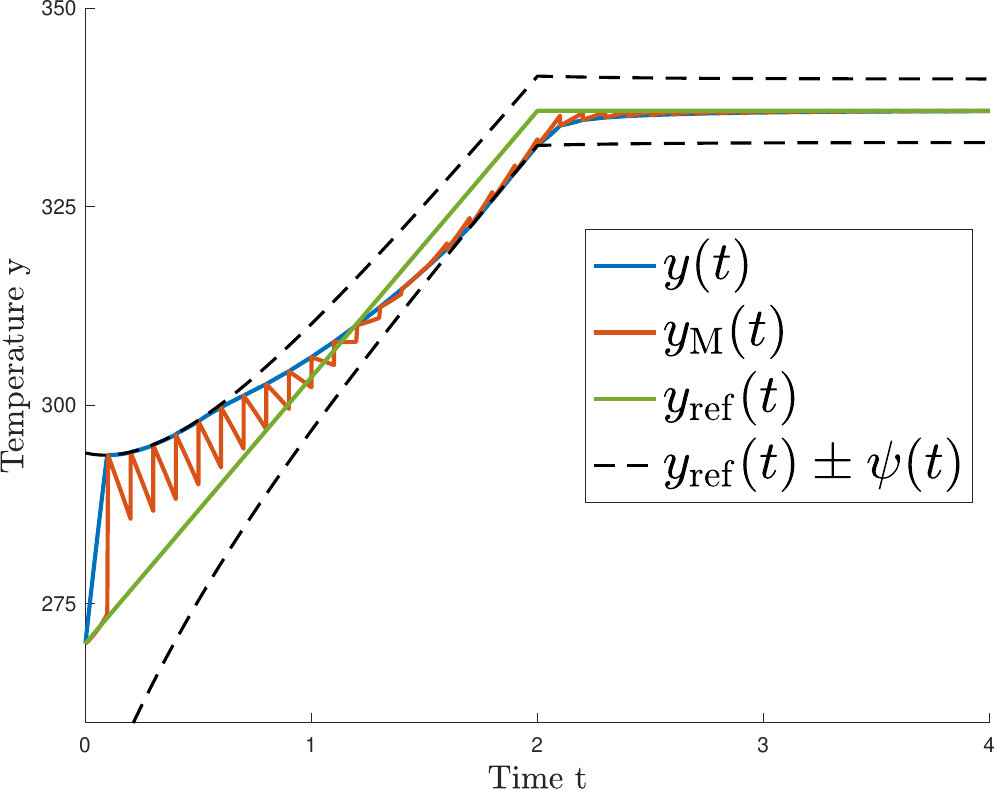}
        \caption{Outputs and reference, with boundary~$\Funnel$.}
        \label{fig:sim:robust_fmpc_reinit:output}
    \end{subfigure}
      \begin{subfigure}[b]{0.49\textwidth}
        \centering
        \includegraphics[width=\textwidth]{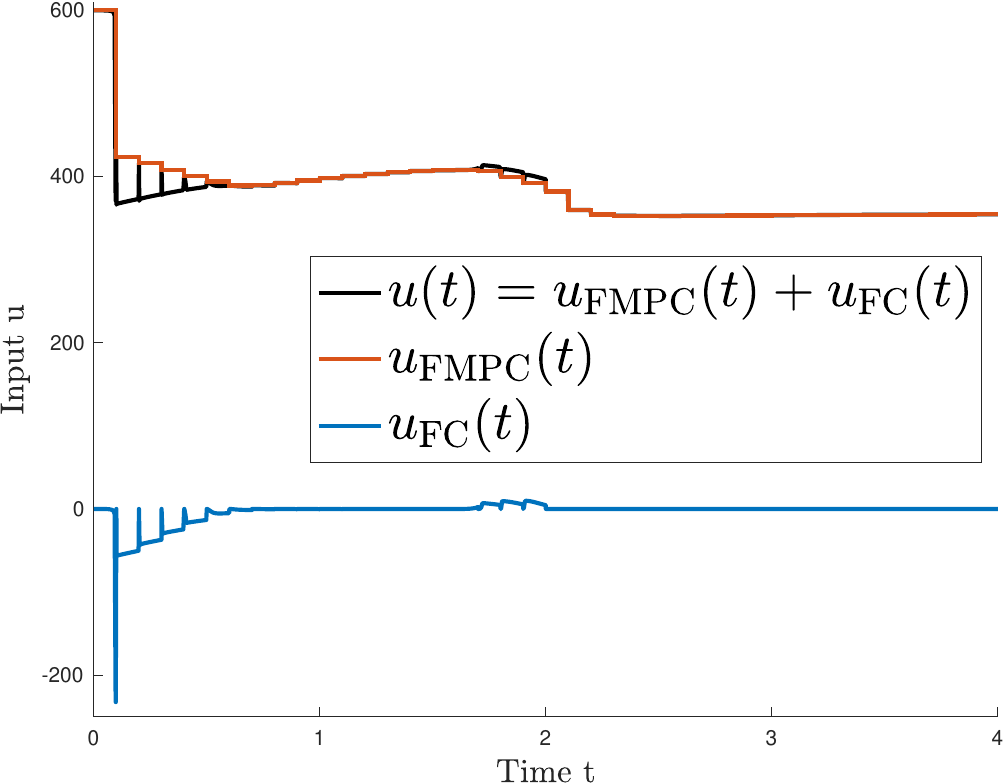}
        \caption{Control inputs.}
        \label{fig:sim:robust_fmpc_reinit:input}
    \end{subfigure}
    \caption{Simulation of system~\eqref{eq:ExampleExothermicReaction} under the
    control generated by the robust funnel MPC~\Cref{Algo:RobustFMPC} with 
    funnel control feedback and model re-initialisation.}
    \label{fig:sim:robust_fmpc_reinit}
\end{figure}

\Cref{fig:sim:robust_fmpc_reinit} shows the application of
\Cref{Algo:RobustFMPC} in Case~3. Besides the additional application of the
funnel controller, the model's state is updated with 
\[
    \xM(t_k) \coloneqq  [x_{\mathrm{M},1}^{k-1}(t_k), x_{\mathrm{M},2}^{k-1}(t_k),y(t_k)]^\top
\]
at the beginning of every MPC cycle. This results in ${\yM(t_k) =y(t_k)}$ at every 
time instant $t_k\in\delta\N_{0}$.
The internal states $x_{\mathrm{M},1}$ and $x_{\mathrm{M},2}$ of the model
remain, however, unchanged during initialisation, i.e. they are initialised with 
their values from the end of previous iteration.
Note that, the proportion of the control signal generated by the MPC component is larger
than in the previous case \eqref{fig:sim:robust_fmpc_trivial} and the funnel controller 
component does overall intervene less. 
Moreover, after $t\approx 0.5$, the funnel controller is inactive most of the
time in~\Cref{fig:sim:robust_fmpc_reinit}, i.e. the applied control signal can
be viewed to be close to \emph{optimal} with respect to the cost
function~\eqref{eq:stageCostFunnelMPC}, since it is computed via the
OCP~\eqref{eq:RobustFMPCOCP}.
In the beginning, the funnel controller however has to compensate for the model inaccuracies 
in order to ensure that the system's output evolves within the funnel boundaries $\Funnel$. 

\subsection*{Mass-on-car system}
We revisit the example of the mass-on-car system from~\Cref{Sec:FMPC:Sim:MassOnCar}. 
The relative degree two system is described by the differential equation
\begin{equation}\tag{\ref{eq:MassOnCarInputOutputDeg2} revisited}
    \begin{aligned}
        \ddot y(t)    &= R_1y(t)+ R_2\dot{y}(t) + S\eta(t) +\Gamma u(t)\\
        \dot \eta(t)  &= Q \eta(t) + P y(t). 
    \end{aligned}
\end{equation}
Assuming the mass $m_2=2$, on the ramp inclined by the angle $\vartheta=\frac{\pi}{4}$,
is connected to the car with mass $m_1=4$ 
via spring and damper system with spring constant $k=2$ and damper constant $d=1$, 
the matrices $R_1$, $R_2$, $S$, $\Gamma$ $Q$, and $R$ have the values as in \eqref{eq:ParamsMassOnCarInputOutputDeg2}.
The objective is tracking of the reference signal $y_{\rf}(t) = \cos(t)$ such
that the tracking error ${y(t)-y_{\rf}(t)}$ evolves within the prescribed
performance funnel given by the function $\Funnel\in\cG$ with $\Funnel(t) =
5\me^{-2t}+0.1$.
To achieve the control objective with robust funnel MPC~\Cref{Algo:RobustFMPC}, 
we use the strict funnel stage cost function 
${\ell_{\Funnel_{2}}:\Rp\times\R\times\R\to\R\cup\{\infty\}}$ given by 
\begin{equation}\tag{\ref{eq:Ex:MassOnCar:CostFunction} revisited}
\begin{aligned}
    \ell_{\Funnel_{2}}(t,\zeta,u) =
    \begin{dcases}
        \frac {\Norm{\zeta}}{\Funnel_{2}(t)^2 - \Norm{\zeta}^2} + \lambda_u \Norm{u}^2,
            & \Norm{\zeta} \neq \Funnel_{2}(t)\\
        \infty,&\text{else},
    \end{dcases}
\end{aligned}
\end{equation}
with design parameter ${\lambda_u\in\Rp}$ and auxiliary funnel function 
$\Funnel_{2}(t)\coloneqq  \frac{1}{\gamma} k_1  \me^{-\FunDeriv (t-t_0)} + \frac{\FunDiam}{\FunDeriv\gamma}$
with $k_1=14$ and $\gamma =0.2$ as in \eqref{ex:Sim:MassOnCar:AuxFunnel}.
For the simulation, the  MPC control signal is restricted to ${\|\uFMPC
\|_\infty \le \umax=30}$ and we choose the design parameters $\lambda_u =
10^{-4}$, prediction horizon $T = 1$, and time shift $\delta = \tfrac{T}{12}\approx 0.083$.
We assume that the MPC component uses a model with incorrect parameters 
\[
    m_1=6,\quad m_2=2,\quad k=3,\quad d=0.75,
\]
for the mass of the car, the mass, the spring constant, and the damper constant. 
This results in a differential equation comparable to \eqref{eq:MassOnCarInputOutputDeg2}.
When referring to this model equation, we use the subscript $\textsc{M}$ in the following.
Moreover, the OCP~\eqref{eq:RobustFMPCOCP} is restricted to step functions with a constant step length of 
$\SampleTime \coloneqq\delta\approx 0.083$ due to discretisation.
For the control law of funnel control component, we choose the bijection $\FCBijec(s)=1/(1-s)$
and the function $\FCSurjec(s)=- s$. The funnel feedback law takes the form 
\begin{equation}\label{eq:Sim:MassOnCar:RobustFC}
    \begin{aligned}
        w(t)&=\phi(t)\eSTrackdot(t)+\FCBijec\rbl\phi(t)^2\eSTrack(t)^2\rbr\phi(t)\eSTrack(t),& \eSTrack(t)&=y(t)-\yM(t),\\
        \uFC(t)&=-\FCBijec\rbl w(t)^2\rbr w(t),& \phi(t)&=\frac{1}{\Funnel(t)-\Norm{\yM(t)-y_{\rf}(t)}},
    \end{aligned}
\end{equation}
where $\yM$ is the prediction for the system output computed by the MPC component. 

We run the simulation on the interval $[0, 10]$ and the system and the
model both use the origin as initial value, i.e.
${[y(0),\dot{y}(0),\eta^1(0),\eta^2(0)]
=[\yM(0),\yMd(0),\eta_{\mathrm{M}}^1(0),\eta_{\mathrm{M}}^2(0)] = [0, 0, 0, 0]}$.
In the following figures, the control signal generated via funnel MPC component is labelled
with the subscript~\textsc{FMCP} ($\uFMPC$); the signal generated by the additional
funnel controller component is labelled with the subscript~\textsc{FC} ($\uFC$).
\begin{figure}[h]
    \begin{subfigure}[b]{0.49\textwidth}
        \centering
        \includegraphics[width=\linewidth]{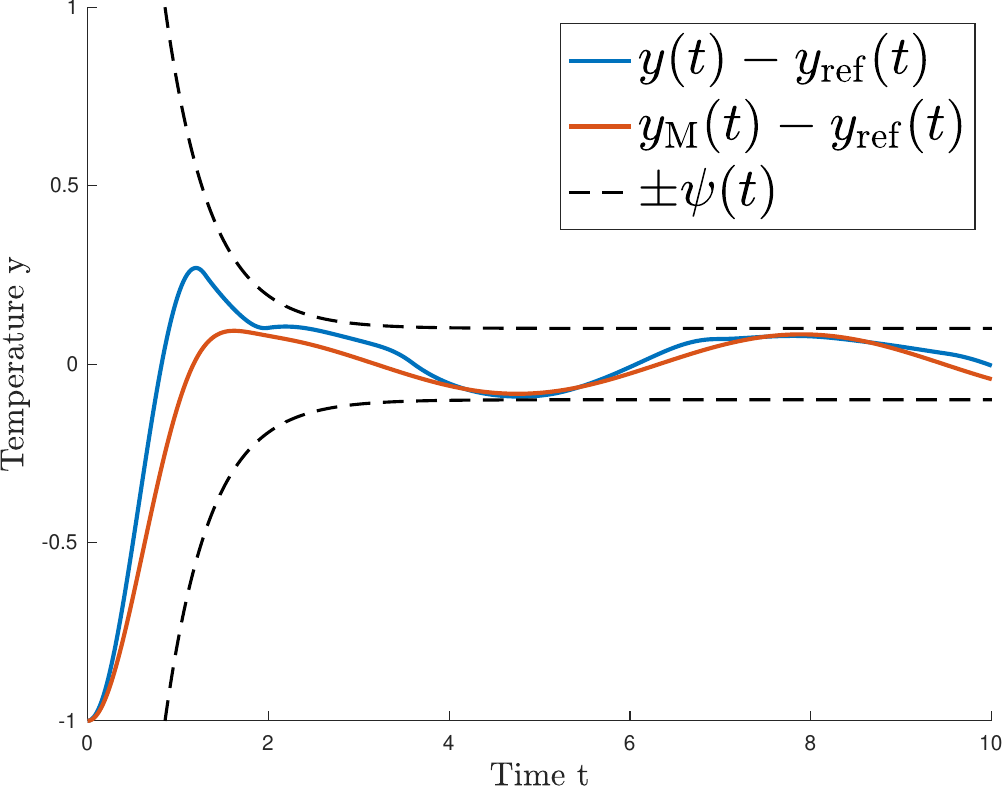}
        \caption{Tracking error $e=y-y_{\rf}$ within boundary~$\Funnel$.}
        \label{fig:sim:massoncar:robust:trivial:output}
    \end{subfigure}
      \begin{subfigure}[b]{0.49\textwidth}
        \centering
        \includegraphics[width=\linewidth]{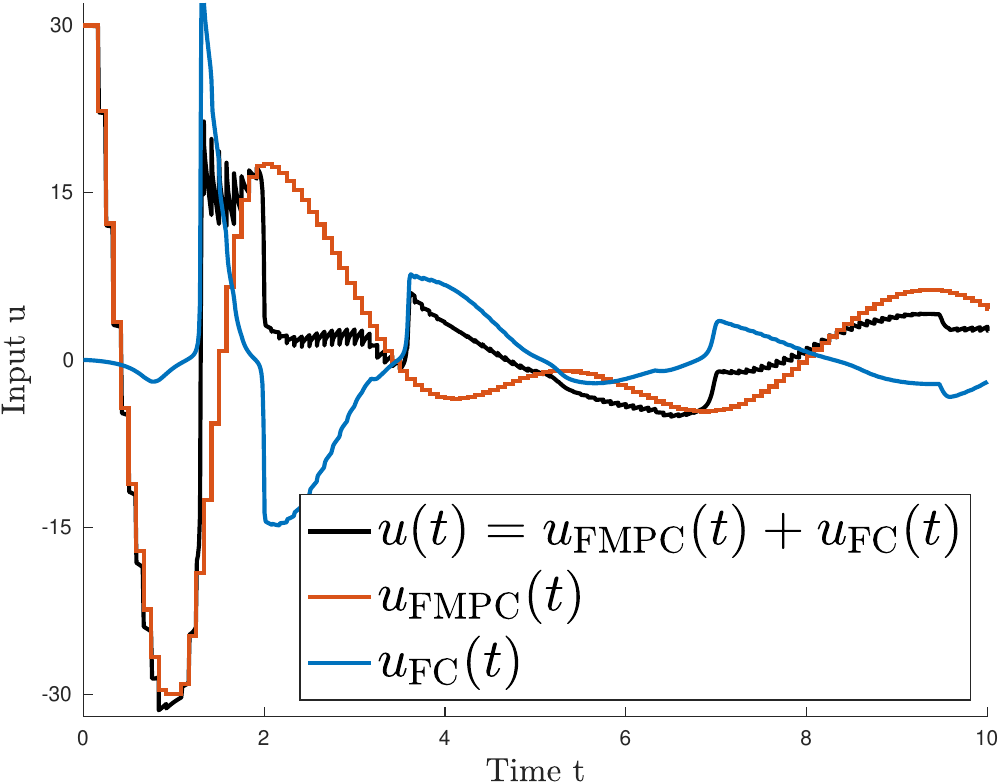}
        \caption{Control inputs.}
        \label{fig:sim:massoncar:robust:trivial:input}
    \end{subfigure}
    \caption{Simulation of system~\eqref{eq:MassOnCarInputOutputDeg2} under the
    control generated by the robust funnel MPC~\Cref{Algo:RobustFMPC} without model re-initialisation.}
    \label{fig:sim:massoncar:robust:trivial}
\end{figure}
\Cref{fig:sim:massoncar:robust:trivial} shows the application of the control signal computed with the
robust funnel MPC~\Cref{Algo:RobustFMPC} to the system~\eqref{eq:MassOnCarInputOutputDeg2} when 
the model is not re-initialised with data from the system. 
The model's state from the previous iteration ($\InitStateK_{k} =
\xM^{k-1}(t_k)$ in Step~\ref{agostep:RobustFMPCFirst} of
\Cref{Algo:RobustFMPC}). The funnel MPC control signal is applied to the system in an open-loop fashion
and it is hence possible to pre-compute the control signal $\uFMPC$.
As the control signal $\uFC$ mitigates the discrepancies between the model's predictions and the system's output, 
both the model's tracking error $\eMTrack(t)=\yM(t)-y_{\rf}(t)$ and the system's tracking error $\eSTrack(t)=y(t)-y_{\rf}(t)$ 
evolve within the boundaries given by $\Funnel$, see \Cref{fig:sim:massoncar:robust:trivial:output}.
Thus, the controller achieves the control objective. 
However, the funnel controller is active over the whole considered time interval to compensate for 
the deviation between the model and system. 
The resulting control signal $u(t)=\uFMPC(t)+\uFC(t)$ shows large fluctuations with peaks, see \Cref{fig:sim:massoncar:robust:trivial:input}.

\begin{figure}[h]
    \begin{subfigure}[b]{0.49\textwidth}
        \centering
        \includegraphics[width=\linewidth]{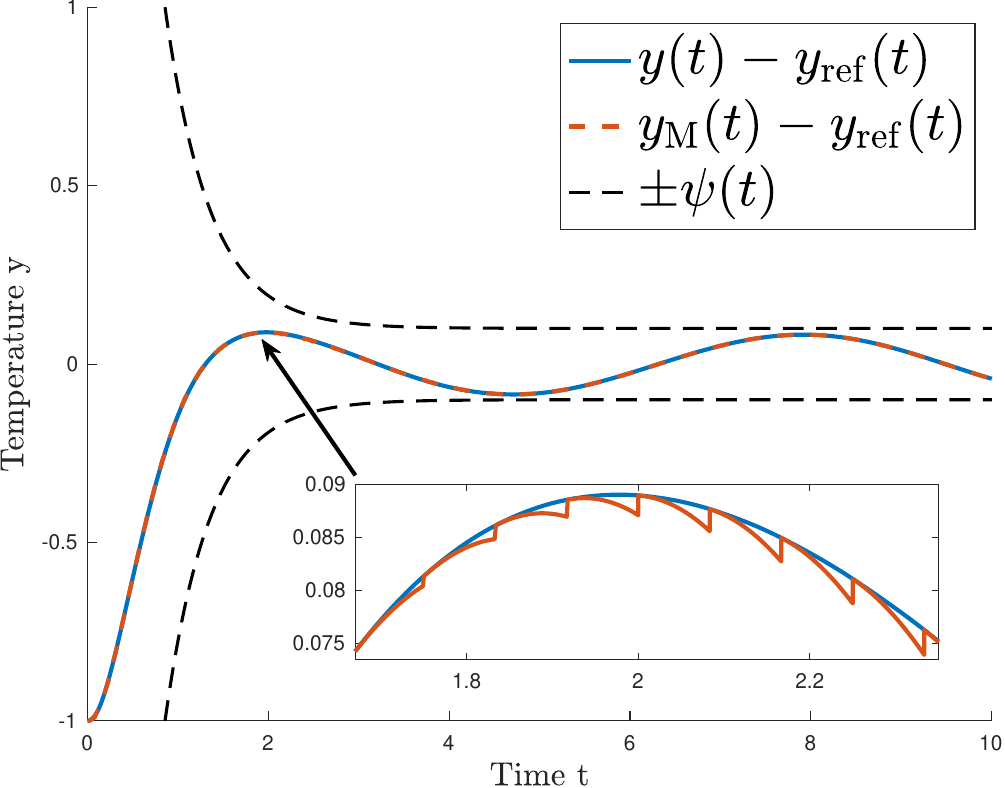}
        \caption{Tracking error $e=y-y_{\rf}$ within boundary~$\Funnel$.}
        \label{fig:sim:massoncar:robust:reinit:output}
    \end{subfigure}
      \begin{subfigure}[b]{0.49\textwidth}
        \centering
        \includegraphics[width=\linewidth]{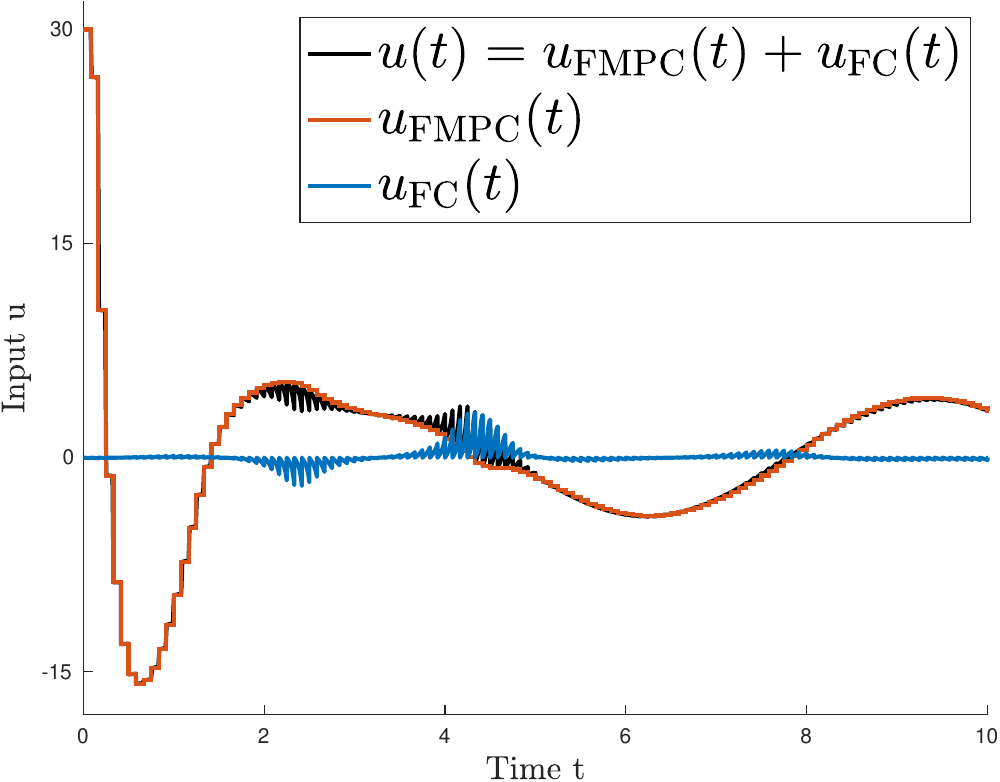}
        \caption{Control inputs.}
        \label{fig:sim:massoncar:robust:reinit:input}
    \end{subfigure}
    \caption{Simulation of system~\eqref{eq:MassOnCarInputOutputDeg2} under the
    control generated by the robust funnel MPC~\Cref{Algo:RobustFMPC} with 
    proper model re-initialisation via the optimisation problem \eqref{eq:Ex:MassOnCar:Sim:Initial}.}
    \label{fig:sim:massoncar:robust:reinit}
\end{figure}

In a second simulation, the model is re-initialised with data from the system
but we leave the rest of the setup unchanged.
To properly initialise the model in accordance with \Cref{Def:ProperInitValues},
we solve, given measurements $\dot{y}(t_k)$ and $\dot{y}(t_k)$ at time $t_k\in\delta\N$,
the following optimisation problem at every iteration of \Cref{Algo:RobustFMPC}
following the ideas from \eqref{eq:MinimizePropInit}.
\begin{equation}\label{eq:Ex:MassOnCar:Sim:Initial}
    \begin{alignedat}{2}
             &\hspace{-0.5cm}\mathop {\operatorname{minimise}}_{\yM^0,\yMd^0\in\R}  \
             \Norm{\yM^0-y(t_k)}^2+\Norm{\yMd^0-\dot{y}(t_k)}^2\\
            &\begin{matrix}
            \text{s.t.}\\
                \null \\
                \null \\
                \null \\
                \null \\
                \null \\
                \null 
            \end{matrix}\quad 
            \begin{aligned}
                      \lambda\Funnel(t_k)&>\Norm{\yM^0-y_{\rf}(t_k)},\\
                      \Funnel_{2}(t_k)&>\Norm{\yMd^0-\dot{y}_{\rf}(t_k)+k_1(\yM^0-y_{\rf}(t_k))},\\
                      \eps&>\hat{\phi}\Norm{\yM^0-y(t_k)},\\
                      \eps&>\Norm{\hat{\phi}(\yMd^0-\dot{y}(t_k))+\FCBijec\rbl\hat{\phi}^2(\yM^0-y(t_k))^2\rbr\hat{\phi}(\yM^0-y(t_k))},\\
                      \hat{\phi}&=\frac{1}{\Funnel(t_k)-\Norm{\yM^0-y_{\rf}(t_k)}}.
            \end{aligned}
    \end{alignedat}
\end{equation}
Afterwards, the solution $\yM^0$ and $\yMd^0$ serves as an initial value for the
model's differential equation \eqref{eq:MassOnCarInputOutputDeg2} at time
$t_k\in\delta\N$. The state $\eta_{\mathrm{M}}$ remains unchanged, i.e. 
the second differential equation in \eqref{eq:MassOnCarInputOutputDeg2} is initialised
with the value of $\eta_{\mathrm{M}}$ from the previous iteration.
The results are displayed in \Cref{fig:sim:massoncar:robust:reinit}.
It is evident that the control scheme is feasible and achieves the control objective.
Both errors~$\yM-y_{\rf}$ and~$y-y_{\rf}$ evolve within the funnel boundaries given by~$\Funnel$.
The model's output $\yM$ diverges from the system's output~$y$ due to the
modelling error. However, it is set back to the system's trajectory at the
beginning of every iteration of the robust funnel MPC~\Cref{Algo:FunnelMPC} as
\Cref{fig:sim:massoncar:robust:reinit:output} shows.
This results in a control signal in which the predominant portion is contributed
by the MPC component. The funnel controller remains mainly inactive only 
compensating for the modelling errors when the system is in a critical state, i.e. 
close the to funnel boundary (i.e. for $t\in [4,5]$), see \Cref{fig:sim:massoncar:robust:reinit:input}.
Its contribution is relatively small but suffices to ensure the adherence of the system's output to the prescribed boundaries.
The overall control signal is less fluctuating and demonstrates a smaller range
of applied control values compared to the previous case.
\chapter{Learning-based robust funnel MPC}\label{Chapter:LearningRobustFMPC}
MPC fundamentally depends on the availability and accuracy of the model 
for the underlying dynamical system.
However, model-plant mismatches and external disturbances pose significant
challenges, driving research into robustification and adaptive strategies.
Building on the previous \Cref{Chapter:RobustFunnelMPC}, which 
introduced \emph{robust funnel MPC} by synergising the funnel MPC~\Cref{Algo:FunnelMPC}
with model-free funnel control, this chapter extends the architecture through
integrated online learning. The original hybrid approach dynamically compensates
for model discrepancies through combined predictive optimisation and reactive
feedback, enabling robust tracking even under severe model-plant mismatches.

Complementing direct robustification efforts, an alternative research direction
focuses on adapting the underlying model to achieve robust constraint satisfaction.
Examples include:
\begin{itemize}
    \item \emph{Data-Driven model refinement}: Techniques like those in
    \cite{BerbKoeh20} leverage persistently exciting data (cf. \cite{WRMDM05,faulwasser2023behavioral})
    for iterative model updates,
    while ensuring initial and recursive feasibility. \emph{Set-membership identification}
    \cite{MILANESE2004957} extend this paradigm by bounding model uncertainties
    using online data, enabling adaptive MPC with guaranteed robust feasibility
    under bounded disturbances \cite{LoreCann19}.
    \item \emph{Iterative Learning Control} (ILC): Leveraging historical trial data,
    ILC refines control inputs cycle-to-cycle for repetitive tasks \cite{Bristow2006}.
    Combined with MPC, modern variants improve controller performance 
    in presence of model mismatch and periodic disturbances \cite{HosseinNia2015,Ma2021}.
    \item \emph{Gaussian process (GP) integration}:  Frameworks, such as
    \cite{hewing2019cautious,maiworm2021online}, combine MPC with Gaussian
    process regression for probabilistic safety guarantees. The latter
    incorporates a non-linear autoregressive exogenous model (NARX) model, while
    the former validates its approach via an autonomous racing case study with
    chance constraints. Similar methods enable safe learning-based control in
    robotics \cite{matschek2023safe}. Hybrid \emph{physics-informed machine
    learning} architectures \cite{Raissi2019} enhance these approaches by
    embedding domain knowledge into learned models, reducing data requirements
    while preserving interpretability \cite{Sanyal2023}. 
    \item \emph{Constrained neural networks}: Utilising tubes containing all possible
    state trajectories \cite{zieger2022non} restricts neural networks to remain
    near predefined nominal models. This ensures safe operation despite
    potential learning failures.
\end{itemize}
In addition, due to the recent advancements in the field of machine learning,
there have been also attempts to utilise such techniques, especially Reinforcement
Learning (RL), to learn an optimal control policy and mimic the behaviour of
(robust) MPC algorithms~\cite{amos2018differentiable,Cao2020,Tagliabue2024}.
Practical applications include chemical reactor control via industrial MPC 
implementations~\cite{hassanpour2024practically}.
\emph{Transfer learning} can further extend this concept by transferring
(safety-critical) control policies across different but related domains,
reducing dependence on large number of system-specific data
needs~\cite{Zhuang2021}.
\emph{Predictive safety filters} \cite{Wabersich18,Wabersich21,Wabersich23}
bridge learning-based control and robust MPC. These filters validate inputs
proposed by learning algorithms (e.g. reinforcement learning) against a
safety-critical model. If unsafe, inputs are modified as little as necessary to
ensure constraint compliance, enabling safe operation while leveraging the
benefits of learning-based control. Similar in idea, hybrid frameworks pair
data-driven controllers with reactive feedback to safeguard the transient
behaviour: 
\begin{itemize}
    \item  Policy iteration~\cite{GottschalkLanza24} and  Q-learning~\cite{LanzaDenn24} paired with safeguards,
    \item  Koopman operator-based MPC for non-linear systems~\cite{BoldLanzWoth2024_Koopman},
    \item  Data-enabled predictive control (DeePC~\cite{coulson2019data}) leveraging the
fundamental lemma by Willems and co-authors~\cite{WRMDM05} for LTI
systems~\cite{Schmitz23}.
\end{itemize}

Surveys \cite{HewingWaber20,wabersich2023data,brunke2022safe} document the
progress of application of various safe learning methods in MPC,
yet ensuring (runtime) safety in complex non-linear systems remains
challenging -- particularly when balancing performance and robustness in
uncertain environments.

Building upon these foundations, this chapter extends the robust funnel MPC
approach presented in~\Cref{Chapter:RobustFunnelMPC} with a general online
learning architecture. This framework continuously improves the surrogate model
using historical data -- system outputs, model predictions, and applied control
signals -- from both the model-based funnel MPC and the model-free feedback
component. Robust tracking within predefined boundaries is achieved while allowing for:
\begin{itemize}
    \item \textbf{Varying model complexity}:  The framework accommodates both fine-tuning
    of detailed models and learning of entirely unknown dynamics. It handles
    low-order linear approximations to high-dimensional non-linear models and
    allows for changes in model dimensionality. 
    \item \textbf{Continual improvement}: By refining the model’s predictive capability
    the controller progressively enhances its performance.
    \item \textbf{Methodological agnosticism}: Rather than prescribing a specific
    learning architecture diverse paradigms and methodologies are supported.
\end{itemize}
By combining learning techniques with both model-based prediction and adaptive
control, this architecture bridges the gap between robustness and adaptability
in uncertain environments.

\section{Controller structure}\label{Sec:StructureLearningFMPC}
\begin{figure}[ht]
\begin{center}
\resizebox{0.95\textwidth}{!}{
    \begin{tikzpicture}[very thick,%
    scale=0.64,%
    node distance = 9ex,
    box/.style={fill=white,rectangle, draw=black},
    blackdot/.style={inner sep = 0, minimum size=3pt,shape=circle,fill,draw=black},%
    blackdotsmall/.style={inner sep = 0, minimum size=0.1pt,shape=circle,fill,draw=black},%
    plus/.style={fill=white,circle,inner sep = 0,very thick,draw},%
    metabox/.style={inner sep = 3ex,rectangle,draw,dotted,fill=gray!20!white}]
    \begin{scope}[scale=0.5]
        \node (sys) [box,minimum size=9ex,xshift=-1ex, fill=orange!60]  {System \eqref{eq:Sys}};
        \node(FC) [box, below of = sys,yshift=-6ex,minimum size=9ex] {Funnel Controller};
        \node(fork1) [plus, right of = FC, xshift=18ex] {$+$};
        \node(fork2) [plus, left of = FC, xshift=-15ex] {$+$};
        \node(fork3) [blackdot, left of = fork2, xshift=-0ex] {};
        \node(MPC) [box, left of = fork3,xshift=-8ex,minimum size=9ex] {Funnel MPC};
        \node(MPCin1) [minimum size=0pt, inner sep = 0pt, below of = MPC, yshift=4.5ex, xshift=2ex] {};
        \node(MPCin1Desc) [minimum size=0pt, inner sep = 0pt, below of = MPCin1, yshift=5.5ex, xshift=2.5ex] {$y$};
        \node(MPCin2) [minimum size=0pt, inner sep = 0pt, below of = MPC, yshift=4.5ex, xshift=-2ex]{};
        \node(MPCin2Desc) [minimum size=0pt, inner sep = 0pt, below of = MPCin2, yshift=5.5ex, xshift=-2.5ex] {$y_{\rf}$};
        \node(refin) [minimum size=0pt, inner sep = 0pt, below of = MPC, yshift=-2ex, xshift=-2ex] {};
        \node(Mod) [box, above of = MPC,yshift=5ex,minimum size=9ex] {Model \eqref{eq:Model_r}};
        \node(fork4) [blackdot, left of = MPC, xshift=-5ex] {};
        \node(fork5) [minimum size=0pt, inner sep = 0pt, below of = fork4, yshift=-3.5ex] {};
        \node(fork9) [blackdot, inner sep = 0pt, right of = fork1, xshift=-2ex ] {};
        \node(fork6) [minimum size=0pt, inner sep = 0pt, below of = fork9, yshift=-2ex] {};
        \node(fork11)[blackdot, inner sep = 0pt, above of = fork9, yshift=6ex ] {};
        \node(fork12)[blackdot, inner sep = 0pt, above of = fork3, yshift=5ex ] {};
        \node(fork13)[blackdot, inner sep = 0pt, above of = fork4, yshift=5ex ] {};
        \node(ML) [box,above of = Mod,minimum size=9ex,xshift=15ex,yshift=8ex] %
            {$\begin{array}{c} \text{Machine}\\ \hspace*{0.55cm}\text{learning}\hspace*{0.55cm} \end{array}$};
        \node(MLu) [minimum size=0pt, inner sep = 0pt, right of = ML, yshift=2ex] {};
        \node(MLo) [minimum size=0pt, inner sep = 0pt, right of = ML, yshift=-2ex] {};
        \node(fork8) [minimum size=0pt, inner sep = 0pt, above of = fork1, yshift=-1ex] {};

        \node(fork7) [minimum size=0pt, inner sep = 0pt, below of = MPCin1, yshift=2.5ex] {};
        \draw[-] (MPC) -- (fork3) node[pos=0.2,above, xshift=2ex] {$u_{\mathrm{FMPC}}$};
        \draw[->] (refin) -- (MPCin2) node[pos=0.4,left] {};
        \draw[-] (fork3) -- (fork12);
        \draw[->] (fork12) -- (Mod);
        \node(ML_uFMPC)[below of = ML, yshift=5.4ex, xshift=2ex] {};
        \node(ML_yM1)[above of = ML, xshift=2ex] {};
        \coordinate[above of = ML_uFMPC, yshift=3ex] (ML_yM1) ;
        \node(ML_yM2)[above of = ML,yshift=-25, xshift=2ex] {};
        \draw[->] (fork12) -- (ML_uFMPC) node[pos=0.85,right] {$u_{\mathrm{FMPC}}$};
        \draw[-] (fork13) |- (ML_yM1.east);
        \draw[-] (fork13) -- (fork4) ;
        \draw[->] (ML_yM1.north) -| (ML_yM2) node[pos=0.7,right] {$y_{\mathrm{M}}$};
        \draw (Mod) -- (fork13) node[pos=0.4,above] {$y_{\mathrm{M}}$}; 
        \draw[->] (fork4) -- (MPC);
        \draw[->] (fork3) -- (fork2.west);
        \draw[-, name path=line4] (sys) -- (fork11) node[pos=0.07,above] {$y$} ;
        \draw[->] (fork9) -- (fork1) node[pos=0.6,right, above] {$+$};
        \draw[-] (fork9) -- (fork6.south);
        \draw[->,name path=line2] (fork5.west) -| (fork1) node[pos=0.9,left] {$-$};
        \draw[->] (fork7.south) -- (MPCin1);
        \path[->,name path=line1] (fork7.west) -- (fork6.east){};
        \path [name intersections={of = line1 and line2}];
        \path [name intersections={of = line1 and line2}];
        \coordinate (S)  at (intersection-1);          
        \path[name path=circle] (S) circle(5.mm);
        \path [name intersections={of = circle and line1}];
        \coordinate (I1)  at (intersection-1);
        \coordinate (I2)  at (intersection-2);
        \tkzDrawArc[color=black, very thick](S,I2)(I1);
        \draw[-] (fork6.east) -- (I2);
        \draw[-] (fork7) -- (I1);
        \draw[-] (fork11) -- (fork9);
        \draw[->] (fork11) |- (MLu) node[pos=0.975,above] {$y$};
        \draw[->] (fork1) -- (FC) node[midway,above] {$y - y_{\mathrm{M}}$};
        \draw[->] (FC) -- (fork2) node[pos=0.4,above, xshift=2ex] {$u_{\mathrm{FC}}$};
        \draw[->] (fork2) |- (sys) node[pos=0.69,above] {$u=u_{\mathrm{FMPC}} + u_{\mathrm{FC}}$};
        \draw (fork4) -- (fork5.south);
        \node(fork10) [blackdot, right of = fork2, xshift=-1ex] {};
        \draw (fork10) |- (fork8.east);
        \draw[->] (ML) -| (Mod) node[pos=0.6,left] {$\begin{array}{l} \text{Model}\\ \text{update}\end{array}$};
        \path[->,name path=line5] (fork8.north) |- (MLo)  node[pos=0.96,below] {$u_{\mathrm{FC}}$};
        \path [name intersections={of = line5 and line4}];
        \coordinate (S3)  at (intersection-1);         
        \path[name path=circle3] (S3) circle(5.mm);
        \path[name path=circle3] (S3) circle(5.mm);
        \path [name intersections={of = circle3 and line5}];
        \coordinate (I5)  at (intersection-1);
        \coordinate (I6)  at (intersection-2);
        \tkzDrawArc[color=black, very thick](S3,I5)(I6);
        \draw[-] (fork8.north) -- (I6);
        \draw[->] (I5) |- (MLo);
    \end{scope}
    \begin{pgfonlayer}{background}
        \fill[gray!25] (-18,-7.4) rectangle (9,1.8);
        \fill[red!20] (-17.5,-6.4) rectangle (-8,1.5);
        \fill[blue!20] (-5,-6.4) rectangle (8,-2.);
        \fill[green!20] (-17.5,2) rectangle (-5,7);
        \node at (-12.5,-8) {\color{red!75}{Model-based controller component}};
        \node at (1.5,-8) {\color{blue!75}{Model-free controller component}};
        \node at (-2,5.35) {\color{green}{Learning component}};
    \end{pgfonlayer}
\end{tikzpicture}
    }
    \caption{Structure of the learning-based robust FMPC scheme.
   The grey box (containing both the red (funnel MPC) and the blue (funnel control) structures) 
   represents the two-component controller \emph{robust funnel MPC} as discussed in~\Cref{Chapter:RobustFunnelMPC}.
   The green box represents the learning component, which receives the four signals: system output~$y$, model output~$\yM$,
   funnel MPC control signal~$\uFMPC$, and funnel control signal~$\uFC$.} 
    \label{Fig:Learning_RFMPC}
    \end{center}
\end{figure}
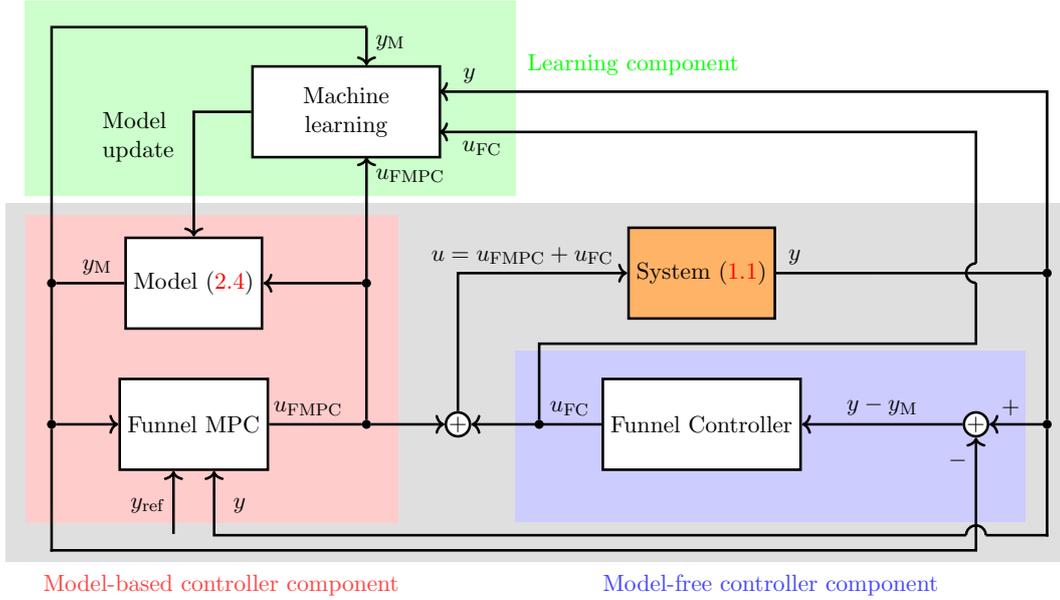

To achieve the overall control task of output reference tracking within
prescribed bounds on the tracking error, we developed a model predictive
controller in~\Cref{Chapter:FunnelMPC}, which ensures superior controller
performance while rigorously maintaining input and output constraints.
However, given the inevitability of model-plant mismatches in practice,
\Cref{Chapter:RobustFunnelMPC} augmented this framework with 
the model-free funnel controller. This addition safeguards the funnel MPC scheme
by guaranteeing satisfaction of the output-tracking criterion even under
disturbances and model uncertainties.
We now introduce a third component -- a data-based learning module -- 
integrated alongside funnel MPC and funnel control, see~\Cref{Fig:Learning_RFMPC}.
This learning component iteratively updates the system model to reduce
model-plant mismatch, thereby progressively enhancing overall control
performance.
A critical challenge lies in ensuring proper functioning of the interplay of these
three components, which necessitates the introduction of additional consistency conditions
(see~\Cref{Def:RestrictedModelClass} and~\Cref{Def:LearningScheme}) for the
model updates -- the key novelty of this approach compared to the robust
funnel MPC~\Cref{Algo:RobustFMPC} (which combines the first two components) proposed
in~\Cref{Chapter:RobustFunnelMPC}.

For the sake of readability and completeness, we recall the robust funnel
MPC~\Cref{Algo:RobustFMPC} and explain the general ideas. In the following, we
simplify the explanation and leave out some details in order to improve
comprehensibility. We refer the reader to~\Cref{Chapter:FunnelMPC,Chapter:RobustFunnelMPC}
for the technical details.

\noindent \textbf{Robust funnel MPC} (grey box in \Cref{Fig:Learning_RFMPC})
is a two-component controller that achieves the control objective 
of tracking a given reference signal $y_{\rf}\in W^{r,\infty}(\Rp,\R^{m})$
within a prescribed performance funnel $\cF_{\Funnel}$ given by $\Funnel\in\cG$,
as laid out in~\Cref{Sec:ControlObjective}.
The controller combines the continuous-time funnel MPC scheme with the adaptive funnel controller.
The model-based funnel MPC component (red box in \Cref{Fig:Learning_RFMPC}) uses a model of the form
\[
   \yM^{(r)}(t) = \fM \big(\oTM(\OpChi(\yM))(t) \big) + \gM \big(\oTM(\OpChi(\yM))(t) \big) u(t) 
\]
as an approximation of the system~\eqref{eq:Sys}, where $(\fM,\gM,\oTM)$ is an element of the model class $\cM^{m,r}_{t_0}$.
At time instants $\hat{t}\in t_0+\delta\N_0$ with $\delta>0$, the current output $y(\hat{t})$ of the system~\eqref{eq:Sys} 
is measured and used to initialise the model, i.e. to select an initial value $\InitState\in\InitValues(\hat{t})$. 
The model is used to predict the future behaviour of the system over the
next time interval of length $T>0$. A control signal $\uFMPC\in
L^{\infty}([\hat{t},\hat{t}+T],\R^m)$ 
satisfying  a given bound $\umax\geq 0$ on the maximal control value  
is computed as a solution of a finite horizon optimal control
problem. The computed model output $\yM$ when applying $\uFMPC$ 
serves as a prediction for the system behaviour.
Utilising a time-varying \emph{funnel penalty function} $\FunnelStageCost$
ensures that the control signal $\uFMPC$ achieves the control objective for the
model, i.e. the model tracking error $\eMTrack(t)\coloneqq \yM(t)-y_{\rf}(t)$ evolves
within the performance funnel $\cF_{\Funnel}$.
Formally, this means that $\uFMPC$ is an element of the set $\Controls(\umax,\InitState)$.

The model-free funnel control component (blue box in \Cref{Fig:Robust_FMPC} and \Cref{Fig:Learning_RFMPC}) computes an 
instantaneous control signal~$\uFC$ based on the deviation between the output~$y$ of system~\eqref{eq:Sys} and the funnel MPC-based predicted~$\yM$. 
The combined control ${u(t) = \uFMPC(t) + \uFC(t)}$ is then applied to the actual system~\eqref{eq:Sys} {at time~$t$}. 
The signal~$\uFC$ from the funnel controller compensates for occurring disturbances,
uncertainties in the model~\eqref{eq:Model_r} and unmodelled dynamics.
In other words, the funnel controller ensures that the model-plant mismatch
$\eSTrack(t)  \coloneqq  y(t)-\yM(t)$ remains small. 
By doing so, not only the model output $\yM$ tracks the reference signal $y_{\rf}$ 
within prescribed boundaries but also the system output $y$, i.e. 
the combined controller achieves the control objective as laid out in~\Cref{Sec:ControlObjective}.
Note that the control signal $\uFC$ is solely determined by the instantaneous
values of the system output~$y$, the funnel function~$\psi$, and the
prediction~$\yM$ made by the model. Therefore, the model-free component cannot
\textit{plan ahead}. This may result in large control values and a rapidly
changing control signal if the actual output significantly deviates from its
predicted counterpart, where the term \emph{significant} is to be understood in
comparison to the current funnel size. 

\noindent \textbf{Learning and improving the model} is the objective of 
the third component that we now incorporate in the overall control scheme (green box in~\Cref{Fig:Learning_RFMPC}).
Since funnel MPC exhibits better controller performance but the robust funnel MPC is able to compensate for model-plant mismatches,
it is desirable to improve the model {so that, preferably,} the control $\uFMPC$ is sufficient to achieve the tracking task 
{with prescribed performance} for the unknown system while satisfying the input constraints~-- in other words,
it is desirable that the funnel controller component is inactive most of the time.
In the following, we identify and establish properties of the learning component such that learning and updating the 
model preserves the structure necessary for robust funnel MPC~\Cref{Algo:RobustFMPC}.
We emphasise that, in the present work, we do not focus on a particular learning scheme 
but develop an abstract learning framework suitable to be combined with robust funnel MPC.
In \Cref{Sec:LearningLinearModel}, we discuss a variant of parameter
identification as one possible instance of a learning scheme; however, we
emphasise that the presented methodology is not restricted to this scheme. As a
result, the particular robustness with respect to model-plant mismatches of
robust funnel MPC even allows to start with ``no model'', e.g. only an
integrator chain, and then learn the remaining drift-dynamics.

The idea of the learning component is to use measurement data from the system output~$y$,
the model output~$\yM$ and its derivatives, i.e. the model state $\xM=\OpChi(\yM)$,
the funnel MPC signal~$\uFMPC$ and the funnel control signal~$\uFC$
to improve the model  used for computation of~$\uFMPC$ in the next iteration of the MPC algorithm (cf.~\Cref{Fig:Learning_RFMPC}).
The data $(y,\xM,\uFMPC,\uFC)$ collected up to the time $\hat{t}\geq t_0$ over the interval $[t_0,\hat{t}]$ in order to be used to 
update the model $(\fM,\gM,\oTM)\in\cM^{m,r}_{t_0}$ is an element of the set
\begin{equation}\label{eq:DefSetOfSignals}
    \SetOfSignals_{\hat{t}}\coloneqq \cC^{r-1}([t_0,\hat{t}],\R^m)\times \cR([t_0,\hat{t}],\R^m)^r\times L^\infty([t_0,\hat{t}],\R^m)\times L^\infty([t_0,\hat{t}],\R^m).
\end{equation}
Note that the image spaces of all signals have the same dimension~$m\in\N$ 
(here we consider $\xM$ to be an element of $\cR([t_0,\hat{t}],\R^m)^r$ instead of $\cR([t_0,\hat{t}],\R^{rm})$).
In order to incorporate an abstract learning scheme~$\cL$ into the funnel MPC algorithm, 
it is imperative that, after updating the model, both other controller components -- 
the model-based funnel MPC and the model-free funnel control -- maintain functionality.
For the functioning of the funnel MPC component, it is necessary to ensure that at every iteration of the MPC scheme that
there exists a control signal feasible for the model. Meaning: given $\InitState\in\InitValues(\hat{t})$
at time $\hat{t}\in t_0+\delta\N_{0}$,
there exists a control $u\in L^{\infty}([\hat{t},\hat{t}+T],\R^m)$ bounded by the constant 
$\umax\geq 0$ that, if applied to the model~\eqref{eq:Model_r}, 
ensures that $\xM(t) - \OpChi(y_{\rf})(t)$ evolves within 
$\cD_{t}^{\Psi}$ for all $t$ over the next time interval of length $T>0$.
In short, 
the set $\Controls(\umax,\InitState)$ as in~\eqref{eq:Def-U} has to be non-empty given $\umax\geq 0$.
\Cref{Th:ExUmax} shows that for every model $(\fM,\gM,\oTM)\in\cM^{m,r}_{t_0}$ there 
exists $\umax\geq 0$ such that $\Controls(\umax,\InitState)\neq\emptyset$ for all $\hat{t}\geq t_0$, $T> 0$,
and $\InitState\in\InitValues(\hat{t})$.
However, the difficulty now lies in ensuring that the input saturation level
$\umax\geq0$ does not increase over time.
We want to a priori choose a uniform $\umax\geq0$ for all models generated 
by the learning component during the operation of the control algorithm.
Moreover, incorporating a learning scheme~$\cL$  must not lead to a globally unbounded 
control signal~$\uFC$ of the model-free funnel controller component.
Establishing the existence of such bound already has been the main challenge in
proving the functioning of the robust funnel MPC~\Cref{Algo:RobustFMPC}
in~\Cref{Thm:RobustFMPC}.
The bound on the maximal control effort of the funnel controller component~$\uFC$
derived in the proof of~\Cref{Thm:RobustFMPC}
depends among other terms on 
\[
    \fMmax+\gMmax\umax, 
\]
where $\fMmax\geq\SNorm{\fM(\oTM(\zeta)|_{[0,\hat{t})})}$ and $\gMmax\geq\SNorm{\gM(\oTM(\zeta)|_{[0,\hat{t})})}$
for all $\hat{t}\in[t_0,\infty]$ and $\zeta\in\FunnelTrajectories_{\hat{t}}$,
see~\eqref{eq:DefEpsr} in the proof of~\Cref{Thm:RobustFMPC} and also~\Cref{Lemma:DynamicBounded}.
To impose uniform maximal control values on both the model-based and the
model-free controller component, we restrict the considered model class
$\cM^{m,r}_{t_0}$ in the following definition.

\begin{definition}[Restricted model class $\cM^{m,r}_{t_0, \umax,\bar{\rho}}$]
    \label{Def:RestrictedModelClass}
    Let $y_{\rf}\in W^{r,\infty}(\Rp,\R^{m})$, ${\bar{\rho}\geq0}$, ${\umax\geq 0}$,
    and  ${\Psi=(\Funnel_1,\ldots,\Funnel_r)\in\FunnelBoundaryFuncs}$. 
    We say that the model~\eqref{eq:Model_r} belongs to the restricted model class~$\cM^{m,r}_{t_0, \umax,\bar{\rho}}$ for $m,r\in\N$, and $t_0\in\Rp$,
    written~$(\fM,\gM,\oTM)\in\cM^{m,r}_{t_0, \umax,\bar{\rho}}$, if 
    \begin{enumerate}
    \item[\Itemlabel{Item:LearningModel}{(L.1)}]
    $(\fM,\gM,\oTM)\in\cM^{m,r}_{t_0}$,
    \item[\Itemlabel{Item:LearningExControl}{(L.2)}]
    $\Controls(\umax,\InitState)\neq\emptyset$ for all $\hat{t}\geq t_0$, $T>0$, and $\InitState\in\mathfrak{I}_{t_0,\hat{t}}^{\Psi}(\hat{t})$,
    \item[\Itemlabel{Item:LearningFCBounded}{(L.3)}]
    $\bar{\rho}\geq\SNorm{\fM(\oTM(\zeta)|_{[0,s)})}+\SNorm{\gM(\oTM(\zeta)|_{[0,s)})}\umax$ for all $s\in[t_0,\infty]$ and $\zeta\in\FunnelTrajectories_{s}$.
    \end{enumerate}
\end{definition}

As~\Cref{Def:RestrictedModelClass} restricts the model class $\cM^{m,r}_{t_0}$ by the properties \ref{Item:LearningExControl}
and \ref{Item:LearningFCBounded}, the question arises for which parameters $\umax,\bar{\rho}\geq0$ the 
restricted model class~$\cM^{m,r}_{t_0, \umax,\bar{\rho}}$ is non empty. The following lemma gives an answer to this question.
\begin{lemma}\label{Lem:SetOfModelsForLearningNonEmpty}
    Let $y_{\rf}\in W^{r,\infty}(\Rp,\R^{m})$, $\Psi=(\Funnel_1,\ldots,\Funnel_r)\in\FunnelBoundaryFuncs$.
    For every $\umax>0$, there exists $\bar{\rho}>0$ such that $\cM^{m,r}_{t_0, \umax,\bar{\rho}}\neq\emptyset$. 
\end{lemma}
\begin{proof}
    Given $\umax>0$, let $\eps>0$ such that 
    \[
        \umax\geq\eps
        \rbl
            \SNorm{y_{\rf}^{(r)}}+\sum_{j=1}^{r-1}k_j\mu_{j}^{r-j}+\SNorm{\dot{\Funnel}_r}
        \rbr, 
    \]
    where $k_j$ with $j=1,\ldots,r-1$ are the parameters associated to the auxiliary
    funnel functions~$(\Funnel_1,\ldots,\Funnel_r)$ and the constants 
    $\mu_{i}^{j}$  are recursively defined via 
    $\mu_i^0 \coloneqq  \SNorm{\Funnel_i}$  and
    $\mu_{i}^{j+1}\coloneqq  \mu_{i+1}^{j}+k_{i}\mu_{i}^{j} $ 
    for $i=1,\ldots, r-1$ and $j=0,\ldots,r-i-1$.
    Set $\fM\equiv0$, $\oTM\equiv0$, and $\gM\equiv\tfrac{1}{\eps}I_{m}$, where $I_{m}$ denotes the identity matrix in $\R^{m\times m}$.
    Utilising these functions, it is easy to see that $(\fM,\gM,\oTM)\in\cM^{m,r}_{t_0}$.
    According to~\Cref{Th:ExUmax}, we have
    $\Controls(\umax,\InitState)\neq\emptyset$ for all $\hat{t}\geq t_0$, $T>0$,
    and $\InitState\in\mathfrak{I}_{t_0,\hat{t}}^{\Psi}(\hat{t})$.
    We choose $\bar{\rho}\geq \frac{1}{\eps}\umax$, then  
    \[
        \SNorm{\fM(\oTM(\zeta)|_{[0,s)})}+\SNorm{\gM(\oTM(\zeta)|_{[0,s)})}\umax =\SNorm{\tfrac{1}{\eps}I_{m}}\umax= \tfrac{1}{\eps}\umax\leq \bar{\rho},
    \]
    for $s\in[t_0,\infty]$ and $\zeta\in\FunnelTrajectories_{s}$.
    Therefore, $(\fM,\gM,\oTM)\in\cM^{m,r}_{t_0, \umax,\bar{\rho}}$.
\end{proof}

\begin{remark}
    For order $r=1$, the set $\cM^{m,r}_{t_0, \umax,\bar{\rho}}$ is non-empty for
    $\umax=\bar{\rho}=0$ if $\dot{y}_{\rf}\equiv 0 $ and $\dot{\psi}\equiv0$.
    Utilising~\Cref{Th:ExUmax}, this can be easily proven by showing
    $(0,I_{m},0)\in\cM^{m,r}_{t_0, \umax,\bar{\rho}}$, where $I_{m}$ denotes the
    identity matrix in $\R^{m\times m}$.
\end{remark}

With~\Cref{Def:RestrictedModelClass} at hand, we define a learning scheme~$\cL$
mapping the signals $(y,\xM,\uFMPC,\uFC)$ collected up to the time $\hat{t}\geq t_0$ to a model 
$(\fM,\gM,\oTM)\in\cM^{m,r}_{t_0,\umax,\bar{\rho}}$.

\begin{definition}(Feasible learning scheme $\cL$)\label{Def:LearningScheme}
    Let $y_{\rf}\in W^{r,\infty}(\Rp,\R^{m})$, $\umax,\bar{\rho}\geq0$, and
    ${\Psi=(\Funnel_1,\ldots,\Funnel_r)\in\FunnelBoundaryFuncs}$
    such that $\cM^{m,r}_{t_0, \umax,\bar{\rho}}\neq\emptyset$.
    We call a function  
    \[
        \cL : \bigcup_{t\geq t_0}\SetOfSignals_{t}\to \cM^{m,r}_{t_0,\umax,\bar{\rho}}
    \]
    a \emph{($\umax$,$\bar{\rho}$)-feasible learning scheme for robust funnel MPC}.
\end{definition}

\begin{remark}\label{Rem:LearningScheme}
    The function $\cL$ maps the hitherto available data at time $\hat{t}$,
    i.e. the signals $(\hat{y},\xMh,\uFMPCh,\uFCh)\in\SetOfSignals_{\hat{t}}$, 
    to a suitable model $(\fM,\gM,\oTM)\in\cM^{m,r}_{t_0, \umax,\bar{\rho}}$.
    Due to the quite abstract nature of~\Cref{Def:RestrictedModelClass,Def:LearningScheme}, a few comments are in order. 
    \begin{enumerate}[(a)]
    \item\label{Item:Rem:LearningControlNonEmpty} Condition \ref{Item:LearningExControl} in \Cref{Def:RestrictedModelClass} can be ensured by prescribing 
    two constants $\fMmax$,$\gMInvmax\geq0$ fulfilling  
    \[
    \fMmax   \geq\SNorm{\fM(\oTM(\zeta)|_{[0,s)})} \text{ and }
    \gMInvmax\geq\SNorm{\gM(\oTM(\zeta)|_{[0,s)})^{-1}}
    \]
    for all
    $s\in[t_0,\infty]$ and $\zeta\in\FunnelTrajectories_{s}$.
    Then, the set $\Controls(\umax,\InitState)$ is non-empty for all $T>0$, $\InitState\in\InitValues(\hat{t})$,  and 
    \[
        \umax\geq\gMInvmax
        \rbl
            \fMmax+\SNorm{y_{\rf}^{(r)}}+\sum_{j=1}^{r-1}k_j\mu_{j}^{r-j}+\SNorm{\dot{\Funnel}_r}
        \rbr,
    \]
    where $k_j$ with $j=1,\ldots,r-1$ are the parameters associated to the auxiliary
    funnel functions~$(\Funnel_1,\ldots,\Funnel_r)$ and the constants 
    $\mu_{i}^{j}$  are recursively defined via 
    $ \mu_i^0 \coloneqq  \SNorm{\Funnel_i}$, 
    $\mu_{i}^{j+1}\coloneqq  \mu_{i+1}^{j}+k_{i}\mu_{i}^{j} $ 
    for $i=1,\ldots, r-1$ and $j=0,\ldots,r-i-1$,
    see~\Cref{Th:ExUmax}. 
    If one additionally prescribes a constant $\gMmax\geq 0$ with ${\gMmax\geq\SNorm{\gM(\oTM(\zeta)|_{[0,s)})}}$
    for all $s\in[t_0,\infty]$ and $\zeta\in\FunnelTrajectories_{s}$, then 
    condition \ref{Item:LearningFCBounded} is fulfilled for $\bar{\rho}\geq \fMmax+\gMmax\umax$.
    \item Condition~\ref{Item:LearningFCBounded} in \Cref{Def:RestrictedModelClass} guarantees that 
    $\yM^{(r)}$ is uniformly bounded by
    \[
        \Norm{\yM^{(r)}(t)}\leq \fMmax+\gMmax\umax
    \]
    independent of the chosen model. In \Cref{Thm:learningFRMPC}, we use this
    estimate to prove the uniform boundedness of the control signal $\uFC$
    generated by the funnel control component.
    \item\label{Rem:Item:ModDomain} The function $\cL$ in \Cref{Def:RestrictedModelClass}
          is defined on the set $\bigcup_{t\geq t_0}\SetOfSignals_{t}$.
          However, it is clear that the domain of $\cL$ can be modified to take additional 
          aspects relevant to the control problem into account.
          We want to comment on certain possibilities.
        \begin{enumerate}[(i)]
        \item\label{Rem:Item:ModDomainSlidingWindow}
        The learning scheme $\cL$ utilises the entire measured data up to the current time instant $\hat{t}$,
        meaning the signals $(\hat{y},\xMh,\uFMPCh,\uFCh)$ are collected over the whole interval $[t_0,\hat{t}]$.
        With increasing time, this results in ever growing memory requirements for the measurements.
        Obviously, this is not suitable in practice. 
        Thus, it is beneficial to use a sliding window approach and use measurements over a time window of length $\hat{\tau}\geq 0$, 
        i.e. the measurements $(\hat{y},\xMh,\uFMPCh,\uFCh)$ are only defined on the interval $[\hat{t}-\hat{\tau},\hat{t}]\cap [t_0,\hat{t}]$. 
        However, to avoid introducing another parameter $\hat{\tau}$ and further complicating \Cref{Def:LearningScheme},
        we assume in this work that signals are indeed available for the whole time interval $[t_0,\hat{t}]$.
        \item\label{Rem:Item:ModDomainModel} 
        In many applications, sufficiently accurate models are often already
        available. Typically, only specific parameters remain unknown,
        inaccurately estimated, or require refinement. Furthermore, as most
        optimisation algorithms inherently require an initialisation point, the
        current model $(\fM,\gM,\oTM)$ can serve as a natural additional input to the 
        learning module~$\cL$.
        This approach achieves dual benefits: reducing computational effort by
        leveraging prior knowledge, while simultaneously mitigating the risk of
        algorithmic instability -- avoiding abrupt, destabilising changes to the
        model structure during successive executions of the function~$\cL$.
        \item The function $\cL$ need not operate solely as a learning algorithm -- it
        can also be utilised to dynamically switch between distinct models
        within the model-based funnel MPC component of the control framework.
        For instance, in systems that operate at different setpoints for
        extended periods, it may be advantageous to employ separate models
        tailored to each operating regime. Here, $\cL$ triggers model switching after
        setpoint transitions, enabling the use of simpler, locally accurate
        models rather than relying on a single complex global model. This
        approach can result in overall improved accuracy while reducing computational
        overhead.
    \end{enumerate}
    \item Since~\Cref{Def:LearningScheme} is rather general, the set of potential learning functions~$\cL$ 
    can be fairly large and difficult to grasp, including with the restrictions~\ref{Item:LearningExControl}
    and \ref{Item:LearningFCBounded} in \Cref{Def:RestrictedModelClass} on the set $\cM^{m,r}_{t_0}$.
    Depending on the specific application, it can therefore be advisable to
    restrict oneself to a subset of potential models in order to simplify the selection of a suitable 
    function~$\cL$ and to be able to compare different learning algorithms more easily.
    In~\Cref{Sec:LearningLinearModel}, we will derive conditions for a 
    learning scheme restricted to linear models to be ($\umax$,$\bar{\rho}$)-feasible.
    \item Given a system~\eqref{eq:Sys} with $(F,\oT) \in \cN^{m,r}_{t_0}$, then $\oT$ is an operator mapping from 
    $\cR(\Rp,\R^n)$ to $L^\infty_{\loc} ([t_0,\infty), \R^{\kappa})$ for some $\kappa\geq0$, see~\Cref{Def:SystemClass}.
    For systems with state representation, see~\Cref{Ex:LTISystem,Ex:ControlAffineMod}, 
    this dimension can be interpreted as the dimension of the internal dynamics of the system. 
    The dimension~$\kappa$ is unknown but fixed. 
    In contrast, for the operator $\oTM\in\cT^{rm,\nu}_{t_0}$ of the model,
    the dimension $\nu \in \N_0$ of the model's internal dynamics can be considered as a parameter in the learning step.
    This means, in order to improve the model such that it ``explains'' the system measurements,
    the dimension of the internal state can be varied. 
    Note that $\nu = 0$ (no internal dynamics) is explicitly allowed for the model.
    \end{enumerate}
\end{remark}

Now, we summarise the reasoning so far in the following algorithm,
which achieves the tracking control objective formulated in~\Cref{Sec:ControlObjective}.
It is a modification of the robust funnel MPC~\Cref{Algo:RobustFMPC}.
Here, the proper re-initialisation {of the model at every iteration} done in~\Cref{Algo:RobustFMPC}
is substituted by the learning component~$\cL$.

\begin{algo}[Learning-based robust funnel MPC]\label{Algo:LearningRFMPC}\ \\
\textbf{Given:}\\[-4ex] %
    \begin{itemize}%
        \item instantaneous measurements of the output $y$ and its derivatives of system~\eqref{eq:Sys},
            initial time $t^0\in\Rp$, initial trajectory $y^0\in\cC^{(r-1)}([0,t_0],\R^m)$, 
            reference signal $y_{\rf}\in W^{r,\infty}(\Rp,\R^{m})$, 
            funnel function $\psi\in\cG$,
        \item  auxiliary funnel boundary function $\Psi=(\psi_1,\ldots,\psi_r)\in\FunnelBoundaryFuncs$ 
        with corresponding parameters $k_i$ for $i=1,\ldots, r$,
        input saturation level $\umax\geq0$,
        parameter $\bar{\rho}$ such that $\cM^{m,r}_{t_0, \umax,\bar{\rho}}\neq \emptyset$,  
        initial model $(\fM^{0},\gM^0,\oTM^0)\in\cM^{m,r}_{t_0, \umax,\bar{\rho}}$, 
        and funnel stage cost function~$\FunnelStageCost$,
        \item initialisation parameters $\eps,\lambda\in(0,1)$ and a $(\umax,\bar{\rho})$-feasible learning scheme $\cL$ as in~\Cref{Def:LearningScheme},
        \item a surjection $\FCSurjec \in \cC(\Rp,\R)$ and  a bijection $\FCBijec\in\cC([0,1), [1,\infty))$.
    \end{itemize} 
    \textbf{Set} the time shift $\delta >0$, 
                 the prediction horizon $T\geq\delta$, and index $k\coloneqq 0$.\\
    \textbf{Define} the time sequence~$(t_k)_{k\in\N_0} $ by $t_k \coloneqq  t_0+k\delta$.\\ 
    \textbf{Steps:}
    \begin{enumerate}[(a)]
    \item\label{agostep:LearningRobustFMPCInit}
    Obtain a measurement $\hat{x}_k\coloneqq \OpChi(y)(t_k)$ of the system output~$y$ and its derivatives at the current time $t_k$ 
    and  choose a  \emph{proper ($\eps$, $\lambda$)-initialisation} $\InitStateK_{k}\in\mathfrak{PI}_{t_0,t_k}^{\Psi,\eps,\lambda}(t_k,\hat{x}_k)$.
    \item \textbf{\textsc{Funnel MPC}}\\ Compute a solution $\uFMPCk\in L^\infty([t_k,t_k +T],\R^{m})$ of the optimal control problem
    \begin{equation}\label{eq:LearningRobustFMPCOCP}
        \mathop
                {\operatorname{minimise}}_{\substack
                {
                    u\in L^{\infty}([t_k,t_k+T],\R^{m}),\\
                    \SNorm{u}  \leq \umax 
                }
            }\      \int_{t_k}^{t_k + T}\FunnelStageCost(s,\eM_{r}(\xM^k(s;t_k,\InitStateK_{k},u)-\OpChi(y_{\rf})(s)),u(s))\d{s}
    \end{equation}
    utilising the model $(\fM^{k},\gM^k,\oTM^k)$.
    Predict the state $\xM^k(t;t_k,\InitStateK_k,\uFMPCk)$ and output~$\yM^k(t;t_k,\InitStateK_k,\uFMPCk)$ of the model on the 
    interval~{$[t_k,t_{k+1}]$}, and define the adaptive funnel $\phi_k: [t_k,t_{k+1}]\to\Rpp $ by 
    \begin{equation} \label{alg:Learning_FMPC:vp}
        \phi_k(t)\coloneqq \frac{1}{\Funnel_1(t)-\Norm{\eMTrack^k(t)}},
    \end{equation}
    where $\eMTrack^k(t) = \yM^k(t) - y_{\rf}(t)$.
    \item\label{alg:step:LearningFC} \textbf{\textsc{Funnel control}}\\
    Using the error variables $\eS_i$ for $i=1,\ldots, r$ as in~\eqref{eq:ek_FC}, define the funnel control law~$\uFC$
    with reference $\yM^k$ and funnel function~$\phi_k$ as in~\eqref{alg:Learning_FMPC:vp} by
    \begin{equation}\label{eq:uLearningFCRobustFMPC}
        \uFCk(t) \coloneqq   (\FCSurjec\circ\FCBijec) (\Norm{\eS_{r}(\phi_{k}(t),\eSTrack(t))}^2)\eS_{r}(\phi_{k}(t),\eSTrack(t)),
    \end{equation}
    with $\eSTrack(t) =y(t)-\yM^k(t)$.
    Apply the control law
    \begin{equation}\label{eq:uLearningRobustFMPC}
        u_k:[t_k,t_{k+1})\to\R^m, \ u_k(t) 
        = \uFMPCk(t)+ \uFCk(t)
    \end{equation}
    to system~\eqref{eq:Sys}. 
    \item \textbf{\textsc{Continual learning}}\label{agostep:LearningRobustFMPCLearning} \\
    Increment~$k$ by $1$, find a feasible model 
        \begin{equation*}
          \cL\big( (y,\xM,\uFMPC,\uFC)|_{[t_0,t_k]} \big) = (\fM^{k},\gM^k,\oTM^k)
        \end{equation*}
    based on the measurement of the signals on the interval $[t_0,t_k]$. 
    Then, go to Step~\ref{agostep:LearningRobustFMPCInit}.
    \end{enumerate}
\end{algo}

\begin{remark} We comment on some aspects of the learning-based robust funnel MPC~\Cref{Algo:LearningRFMPC}.
    \begin{enumerate}[(a)]
    \item The signals $(y,\xM,\uFMPC,\uFC)$ used for the learning scheme~$\cL$ during Step~\ref{agostep:LearningRobustFMPCLearning}
    are the whole trajectories of the individual functions up to the current time $t_{k+1}$.
    This means that $y$ is the solution of the system differential equation~\eqref{eq:Sys} up to the current time,
    the control signals $\uFMPC$ and $\uFC$ are the concatenation of the 
    signals $\uFCk$ and $\uFMPCk$ applied at every interval $[t_i,t_{i+1}]$ for $i=0,\ldots, k$, 
    and $\xM$ is the concatenation of the solutions $\xM^i$ of the model differential equation~\eqref{eq:Model_r} 
    with model $(\fM^{i},\gM^i,\oTM^i)$ and initial value $\InitStateK_{i}$ on the interval $[t_i,t_{i+1}]$ for $i=0,\ldots, k-1$.
    To be more precise: 
    \[
        \uFMPC(t) = \uFMPCk[i](t),\quad  \uFC(t) = \uFCk[i](t),\quad \xM(t)=\xM^i(t;t_{i},\InitStateK_{i},\uFMPCk[i])
    \]
    for $t\in[t_i,t_{i+1})$ and $i=0,\ldots, k$.
    Note that $\xM$ is not a concatenated solution in the sense of~\Cref{Def:SolutionClosedLoop} 
    as the model $(\fM^{k},\gM^k,\oTM^k)$ changes at every iteration of the~\Cref{Algo:LearningRFMPC}.
    \item Let
        $
        \umax\geq
        \rbl
            \SNorm{y_{\rf}^{(r)}}+\sum_{j=1}^{r-1}k_j\mu_{j}^{r-j}+\SNorm{\dot{\Funnel}_r}
        \rbr
        $
    and $\bar{\rho}=\umax$
    where $k_j$ are the parameters associated to the 
    funnel functions $(\Funnel_1,\ldots,\Funnel_r)$ for ${j=1,\ldots,r-1}$. The constants 
    $\mu_{i}^{j}$  are recursively defined via 
    $ \mu_i^0 \coloneqq  \SNorm{\Funnel_i}$,
    $\mu_{i}^{j+1}\coloneqq  \mu_{i+1}^{j}+k_{i}\mu_{i}^{j}$
    for ${i=1,\ldots, r-1}$ and $j=0,\ldots,r-i-1$.
    Then, the integrator chain 
    \[
        \yM^{(r)}(t)=u(t)
    \]
    is a model in $\cM^{m,r}_{t_0, \umax,\bar{\rho}}$, 
    see proof of~\Cref{Lem:SetOfModelsForLearningNonEmpty}.
    The model-based MPC component of the control scheme can operate in this sense ``without'' a model.
    It therefore is possible to apply the learning-based robust funnel MPC~\Cref{Algo:LearningRFMPC} without an initial model
    or an offline learning phase.
    \item In practice, it may often not be desirable to update the model $(\fM^{k},\gM^k,\oTM^k)$ 
    at every iteration of the~\Cref{Algo:LearningRFMPC}. 
    Especially, if the execution of the learning procedure is very time-consuming, 
    it may be advantageous to evaluate~$\cL$ only every $i$-th iteration for $i>1$.
    \item Note that, the initialisation in Step~\ref{agostep:LearningRobustFMPCInit} of~\Cref{Algo:LearningRFMPC} 
    at time $t_k\in t_0+\delta\N_{0}$ is independent of the current model~$(\fM^k,\gM^k,\oTM^k)$.
    It only depends on $y_{\rf}$, $\Psi$, $\eps,\lambda$, and $\hat{x}_k=\OpChi(y)(t_k)$, see~\Cref{Def:ProperInitValues}.
    Instead of $\InitStateK_{k}\in\PropInitValues(t_k,\hat{x}_k)$ for a fixed $\tau\geq0$
    as in~Step~\ref{agostep:RobustFMPCFirst} of the robust funnel MPC~\Cref{Algo:RobustFMPC},
    we require $\InitStateK_{k}\in\mathfrak{PI}_{t_0,t_k}^{\Psi,\eps,\lambda}(t_k,\hat{x}_k)$ 
    in \eqref{eq:LearningRobustFMPCOCP} of~\Cref{Algo:LearningRFMPC}, i.e. 
    both components of $\InitStateK_{k}$ are defined on their entire maximal time intervals up to $t_k$.
    By doing so, we avoid having to deal with changing memory limits for the operators $\oTM^k$.
    This set is non empty for $\hat{x}_k$ with ${\hat{x}_k-\OpChi(y_{\rf})(t_k)\in\cEFC{\eps}(1/\Funnel_1(t_k))}$, see~\Cref{Rem:PropInitialValuesNonEmpty}.
    In case all operators generated by the learning scheme~$\cL$ during Step~\ref{agostep:LearningRobustFMPCLearning} have a memory limit 
    lower or equal than a pre-specified bound $\bar{\tau}\geq0$, 
    the initialisation  can alternatively be chosen as  
    $\InitStateK_{k}\in\mathfrak{PI}_{t_0,\bar{\tau}}^{\Psi,\eps,\lambda}(t_k,\hat{x}_k)$
    in Step~\ref{agostep:LearningRobustFMPCInit} of~\Cref{Algo:LearningRFMPC}.
    \end{enumerate}
\end{remark}

We are now in the position to formulate the main result of this chapter,
which extends~\Cref{Algo:RobustFMPC} and the corresponding~\Cref{Thm:RobustFMPC} by the learning component.

\begin{theorem}\label{Thm:learningFRMPC}
    Consider a system~\eqref{eq:Sys} with $(F,\oT) \in \cN^{m,r}_{t_0}$ as in~\Cref{Def:SystemClass}. 
    Let $t_0\geq 0$ be the initial time, $y_{\rf} \in W^{r,\infty}(\Rp,\R^m)$,
    and $\Psi=(\Funnel_1,\ldots,\Funnel_r)\in \FunnelBoundaryFuncs$ be given, and
    let $y^0\in\cC^{(r-1)}([0,t_0],\R^m)$ be the initial trajectory for the system~\eqref{eq:Sys}
    with $\OpChi(y_0-y_{\rf})(t_0)\in\cEFC{1}(1/\Funnel(t_0))$.
    Further, let $\umax,\bar{\rho}>0$ such that the set of models $\cM^{m,r}_{t_0, \umax,\bar{\rho}}$ is non-empty.
    There exist $\eps,\lambda\in(0,1)$ such that,
    for every initial model ${(\fM^{0},\gM^0,\oTM^0)\in\cM^{m,r}_{t_0, \umax,\bar{\rho}}}$ and 
    for every $(\umax,\bar{\rho})$-feasible learning scheme
    ${\cL : \bigcup_{t\geq t_0}\SetOfSignals_{t}\to \cM^{m,r}_{t_0,\umax,\bar{\rho}}}$,
    the robust funnel MPC \Cref{Algo:RobustFMPC} with $\delta>0$ and $T\ge\delta$ is initially and recursively feasible,
    i.e. at every time instant $t_k \coloneqq  t_0+k\delta $ for $k\in\N_0$
    \begin{itemize}
        \item there exists a proper initialisation $\InitStateK_{k}\in\mathfrak{PI}_{t_0,t_k}^{\Psi,\eps,\lambda}(t_k,\hat{x}_k)$ and
        \item the OCP~\eqref{eq:RobustFMPCOCP} has a solution $\uFMPCk\in L^\infty([t_k,t_k+T],\R^m)$.
    \end{itemize}
    Moreover,  the closed-loop system consisting of the system~\eqref{eq:Sys} and the feedback law~\eqref{eq:uLearningRobustFMPC}
    has a global solution $y : [0,\infty) \to \R^m$. 
    Each global solution $y$ satisfies that
    \begin{enumerate}[label = (\roman{enumi}), ref=(\roman{enumi})]
        \item \label{Assertion:y_u_bounded}
        all signals are bounded, in particular, $u\in L^\infty([t_0,\infty),\R^m)$ and $y \in W^{r,\infty}(\Rp,\R^m)$, 
        \item \label{Assertion:tracking_error}
        the tracking error between the system's output and the reference evolves within prescribed boundaries, i.e.
        \begin{equation*}
            \fa t \ge t_0 : \Norm{y(t)  - y_{\rf}(t)}< \Funnel_1(t) .
        \end{equation*}
    \end{enumerate}
\end{theorem}
\begin{proof}
    The learning-based robust funnel MPC~\Cref{Algo:LearningRFMPC} differs from the robust funnel MPC~\Cref{Algo:RobustFMPC}
    only in two aspects, the utilisation of the learning scheme~$\cL$ in Step~\ref{agostep:LearningRobustFMPCLearning} of the algorithm
    and the usage of changing models $(\fM^k,\gM^k,\oTM^k)$ in Step~\ref{agostep:LearningRobustFMPCInit}.
    We will show how the proof of~\Cref{Thm:RobustFMPC} can be adapted to the current setting.
    However, as the proof of~\Cref{Thm:RobustFMPC} does, in large parts, not depend on the used model,
    we will not repeat all technical details. 

    \noindent
    \emph{Step 1}:
    Let $\delta>0$ and $T\ge\delta$ be arbitrary but fixed. 
    Note that, the set of controls~$\Controls(\hat{t},\InitStateK)$ is non-empty 
    for all  $(\fM,\gM,\oTM)\in\cM^{m,r}_{t_0,\umax,\bar{\rho}}$, 
    all $\hat{t}\geq t_0$, and all $\InitState\in\mathfrak{I}_{t_0,\hat{t}}^{\Psi}(\hat{t})$
    due to property \ref{Item:LearningExControl} of $\cM^{m,r}_{t_0,\umax,\bar{\rho}}$, see \Cref{Def:RestrictedModelClass}.
    Define $\lambda$ and $\eps$ as in Step~1--3 in the proof of~\Cref{Thm:RobustFMPC}.
    Then, $\eps\in(0,1)$ is constructed in a way such that  
    we have $\hat{x}_0-\OpChi(y_{\rf})(t_0)\in\cEFC{\eps}(1/\Funnel_1(t_0))$ for $\hat{x}_0\coloneqq \OpChi(y)(t_0)=\OpChi(y^0)(t_0)$, 
    see definition of $\eps$ in the proof of~\Cref{Thm:RobustFMPC}.
    Thus, $\mathfrak{PI}_{t_0,\hat{t}}^{\Psi}(\hat{t})(t_0,\OpChi(y)(t_0))\neq\emptyset$, see~\Cref{Rem:PropInitialValuesNonEmpty}.
     
    \noindent
    \emph{Step 2}:
    Let $\cL : \bigcup_{t\geq t_0}\SetOfSignals_{t}\to \cM^{m,r}_{t_0,\umax,\bar{\rho}}$ be a  $(\umax,\bar{\rho})$-feasible learning scheme.
    When applying the learning-based robust funnel MPC~\Cref{Algo:LearningRFMPC} to the system~\eqref{eq:Sys}, the system's dynamics 
    on each interval $[t_k,t_{k+1}]$ are given by
    \begin{equation}\label{eq:SysPiecewiseLearning}
       y^{(r)}_k(t)=F(\oT(\OpChi(y_k))(t),u_k(t)),\quad y_k|_{[0,t_k]}= y_{k-1}|_{[0,t_k]}
    \end{equation}
    where $y_{-1}\coloneqq y^0$ and $u_k$ is the control law given by~\eqref{eq:uLearningRobustFMPC}.
    In Step~4 of the proof of~\Cref{Thm:RobustFMPC}, it was inductively shown that the
    robust funnel MPC~\Cref{Algo:RobustFMPC} is initially and recursively feasible.
    This means, in particular, that there exists a proper initialisation $\InitStateK_{k}\in\mathfrak{PI}_{t_0,t_k}^{\Psi}(t_k,\OpChi(y_{k-1})(t_k))$
    at every time instant $t_k$, that $u_k$ as in~\eqref{eq:uRobustFMPC} is well defined on every interval $[t_k,t_{k+1}]$,
    and that~\eqref{eq:SysPiecewiseLearning} has a maximal solution $y_k$ defined on the entire interval $[t_k,t_{k+1}]$.
    Step~4 of the proof of~\Cref{Thm:RobustFMPC} does not depend on the concrete choice 
    $(\fM,\gM,\oTM)\in\cM^{m,r}_{t_0,\umax,\bar{\rho}}$ of the model used on the time interval $[t_k,t_{k+1}]$.
    Only two of the model's aspects are used within the proof: the non-emptiness of the set $\cU_{[t_k,t_k+T]}(\umax,\InitStateK_k)$ 
    and the uniform boundedness of
    \[ 
        \Norm{\yM^{(r)}(t)}=\Norm{\fM(\OpChi(\yM)(t))+\gM(\OpChi(\yM)(t))\uFMPCk(t)}.
    \]
    The former one is directly fulfilled by property~\ref{Item:LearningExControl} of $\cM^{m,r}_{t_0,\umax,\bar{\rho}}$,
    see \Cref{Def:RestrictedModelClass}.
    The latter one is also satisfied since property~\ref{Item:LearningFCBounded} ensures 
    \[
        \SNorm{\fM(\oTM(\zeta)|_{[0,s)})}+\SNorm{\gM(\oTM(\zeta)|_{[0,s)})}\umax\leq \bar{\rho}
    \]
    for all $s\in[t_0,\infty]$ and $\zeta\in\FunnelTrajectories_{s}$.
    One therefore can adapt Step~4 of the proof of~\Cref{Thm:RobustFMPC} to the current setting 
    in order to show that there exists a proper initialisation ${\InitStateK_{k}\in\mathfrak{PI}_{t_0,t_k}^{\Psi}(t_k,\OpChi(y_{k-1})(t_k))}$
    at every time instant $t_k$, that $u_k$ as in~\eqref{eq:uLearningRobustFMPC} is well defined on every interval $[t_k,t_{k+1}]$,
    and that~\eqref{eq:SysPiecewiseLearning} has a maximal solution $y_k$ defined on the entire interval $[t_k,t_{k+1}]$
    if the learning-based robust funnel MPC~\Cref{Algo:LearningRFMPC} is applied to the system~\eqref{eq:Sys}.
    The existence of a solution
    $\uFMPCk\in\cU_{[t_k,t_k+T]}(\umax,\InitStateK_k)$ of the
    OCP~\eqref{eq:LearningRobustFMPCOCP} at time instant $t_k$ for $k\in\N_0$ is a
    direct consequence of \Cref{Th:SolutionExists} and the non-emptiness of
    $\cU_{[t_k,t_k+T]}(\umax,\InitStateK_k)$.

    \noindent
    \emph{Step 3}:
    The signal $\uFMPCk$ as an element of $\cU_{[t_k,t_k+T]}(\umax,\InitStateK_k)$ is bounded by $\umax\geq 0$ 
    for all $k\in\N_{0}$. 
    The funnel control signal $\uFCk$ is bounded by $\HighGainFunc(\FCSurjec\circ\FCBijec(\eps^2))$ for all $k\in\N_0$ 
    because of the construction of $\eps\in(0,1)$, c.f.~Step~5 of the proof of~\Cref{Thm:RobustFMPC}.
    Moreover, we have 
    \[
        \OpChi(y-\yM)(t)\in\cEFC{\tilde{\eps}}(\phi(t))
    \]
    for all $t\in[t_0,\infty)$ for some $\tilde{\eps}\in(0,1)$, c.f.~Step~5 of the proof of~\Cref{Thm:RobustFMPC}.
    Since $\yM$ and $\phi$ are bounded functions, $y \in W^{r,\infty}(\Rp,\R^m)$, see definition of $\cEFC{\tilde{\eps}}$ in~\eqref{eq:ek_FC}.
    Finally, 
    \[
        \Norm{y(t)-y_{\rf}(t)}\leq \Norm{y(t)-\yM(t)}+\Norm{\yM(t)-y_{\rf}(t)}<\phi(t)+\Norm{\yM(t)-y_{\rf}(t)}=\Funnel_1(t)
    \]
    for all $t\geq t_0$.
    This shows~\ref{Assertion:learning_tracking_error} and completes the proof.
\end{proof}

\begin{remark}
    Conditions \ref{Item:LearningExControl} and \ref{Item:LearningFCBounded} in
    \Cref{Def:RestrictedModelClass} ensure
    $\cU_{[t,t+T]}(\umax,\InitStateK_k)\neq\emptyset$ for all $t\geq t_0$ and
    $T>0$ and that the funnel control signal $\uFCk$ is uniformly bounded for
    all $k\in\N$.
    As the attentive reader might have noticed, it is possible to relax
    these conditions during operation of the learning-based robust
    funnel MPC~\Cref{Algo:LearningRFMPC}.
    Firstly, it is possible to fix the prediction horizon $T>0$.
    Moreover,  it is sufficient that the model~$(\fM^k,\gM^k,\oTM^k)$ chosen at time instant $t_k\in t_0+\delta\N_{0}$ 
    can ensure these properties for all future time $t\geq t_k$.
    It is not required to ensure them for the past, i.e. $t\leq t_k$. 
    To be precise, one can replace \ref{Item:LearningExControl} at time $t_k$ and with given $T>0$ by 
    \begin{enumerate} 
    \item[\Itemlabel{Item:LearningExControlAlt}{(L.2')}]
    $\Controls(\umax,\InitState)\neq\emptyset$ for all $\hat{t}\geq t_k$ and $\InitState\in\mathfrak{I}_{t_0,\hat{t}}^{\Psi}(\hat{t})$,
    \end{enumerate}
    and the condition \ref{Item:LearningFCBounded} can be relaxed by 
    \begin{enumerate} 
    \item[\Itemlabel{Item:LearningFCBoundedAlt}{(L.3')}]
    $\bar{\rho}\geq\SNorm{\fM(\oTM(\zeta)|_{[t_k,s)})}+\SNorm{\gM(\oTM(\zeta)|_{[t_k,s)})}\umax$ for all $s\in[t_k,\infty]$ and $\zeta\in\FunnelTrajectories_{s}$.
    \end{enumerate}
    In order to avoid introducing a time dependency and thus an additional parameter which introduces even more  
    technicalities, we refrained from formulating \Cref{Def:RestrictedModelClass} in this more general way.
\end{remark}

\section{On learning schemes}\label{Sec:LearningLinearModel}
In recent years, data-driven control has attracted significant attention, with a
proliferation of research contributions in the field. These results can broadly
be categorised into control schemes for linear systems and techniques developed
for non-linear systems. 
Bolstered by successful applications~\cite{Mezi13}, powerful numerical methods
such as extended dynamic mode
decomposition~\cite{williams:kevrekidis:rowley:2015}, and theoretical advances
-- including convergence guarantees in the infinite-data
limit~\cite{korda:mezic:2018b}, finite-data error bounds~\cite{KohnPhil24},
and extensions to stochastic control systems~\cite{nuske2023finite} -- the Koopman
formalism~\cite{BrunKutz22}, originally proposed  
in~\cite{koopman1931hamiltonian}, has emerged as a cornerstone for
data-driven controller design \cite{GoswPale21,OttoRowl21,StraScha24:generator}. 
Recent work has further extended this framework to model predictive control, 
establishing rigorous closed-loop guarantees \cite{KordMezi18:MPC,BoldGrun25,BoldScha25}. 
For linear time-invariant systems, Subspace Predictive Control~\cite{favoreel1999spc} 
has gained prominence, while the so-called fundamental lemma by Willems and
co-authors~\cite{WRMDM05} enables direct data-driven methods such as DeePC~\cite{coulson2019data}.
Complementary approaches include Reinforcement Learning (RL)~\cite{Sutton2018},
Gaussian processes for uncertainty-aware
designs~\cite{Kocijan2004,maiworm2021online,hewing2019cautious}, and SINDY for sparse identification of non-linear dynamics
\cite{brunton2016discovering}. Deep neural
networks (DNNs) have further expanded the scope of data-driven control, enabling
approximation of complex dynamics and control policies for high-dimensional
systems~\cite{Pillonetto2025,Schussler2019,Cao2020}. Recent advances also address
safety-critical scenarios through Hamilton-Jacobi reachability
analysis~\cite{bansal2017hamilton}.

The structural conditions provided in \Cref{Def:LearningScheme} can be used to
define suitable learning algorithms based on the previously discussed
techniques -- for linear as well as for non-linear systems.

In this section, we derive sufficient conditions on the parameters of models to be learned
in order to make them eligible for a learning scheme~$\cL$ as defined in~\Cref{Def:LearningScheme}.
Since in many applications a linear model may serve as a good prediction model, 
we derive sufficient conditions on the parameters of linear systems of the form
\begin{equation}\label{eq:LearningLinearExample}
\begin{aligned}
    \yM^{(r)}(t) &= \sum_{j=1}^{r}R_j\yM^{(j-1)}(t)+S\eta+D_1+\Gamma u(t),&\OpChi(\yM)(t_0)&=\yM^0,\\
    \dot \eta(t)  &= Q \eta(t)+ P \yM(t)+D_2, &\eta(t_0)&=\eta^0
\end{aligned}
\end{equation}
where $R_j \in \R^{m \times m}$ for all $j=1,\ldots,r$, $S, P^\top\in\R^{m\times \nu}$, $D_1\in\R^m$, $D_2\in\R^{\nu}$,
$Q\in\SGroup_{\nu}^{--}$, and $\Gamma\in\GL_m(\R)$. We use in the following the notation $R\coloneqq (R_1,\ldots,R_r)$
and denote the largest eigenvalue of the symmetric negative definite matrix $Q$ by $\lambda_{\max}(Q)<0$.
Define the functions 
\begin{equation} \label{eq:f_ParameterAndg_Parameter}
    \begin{aligned}
        f_{D_1} : \R^{\nu}&\to\R^m, &\eta &\mapsto \eta + D_1, \\
        g_{\Gamma}:  \R^{\nu}&\to\R^{m\times m}, &\eta &\mapsto \Gamma, \\
    \end{aligned}
\end{equation}
and the linear integral operator $\oT_{R,S,Q,P,D_2,\eta_0}: \cR(\Rp,\R^{m})^r \to L^\infty_{\loc} ([t_0,\infty), \R^{\nu})$ by
\begin{equation}\label{eq:oTM_Parameter}
    \oT_{R,S,Q,P,D_2,\eta_0}(z_1,\ldots,z_r)(t)\coloneqq \sum_{j=1}^{r}R_jz_j(t)
    +S\rbl\me^{Q(t-t_0)}\eta_0+\int_{t_0}^t\me^{Q(t-s)}\rbl Pz_1(s) + D_2 \rbr\d{s}\rbr\!\!.
\end{equation}
Using these functions, the model~\eqref{eq:LearningLinearExample} can be written in the form~\eqref{eq:Model_r}, i.e. 
\[
    \yM^{(r)}(t)=f_{D_1}(\oT_{R,S,Q,P,D_2,\eta_0}(\OpChi(\yM))(t))+ g_{\Gamma}(\oT_{R,S,Q,P,D_2,\eta_0}(\OpChi(\yM))(t))u(t)
\]
with initial value $\OpChi(\yM)(t_0)=\yM^0$.
Let $y_{\rf} \in W^{k,\infty}(\Rp,\R^m)$ and $\Psi=(\Funnel_1,\ldots,\Funnel_r)\in \FunnelBoundaryFuncs$, let 
\begin{equation}\label{eq:FunnelFrak}
    \mathfrak{F}\coloneqq \setdef{x\in\R^{rm}}{t\geq t_0, x-\OpChi(y_{\rf})(t)\in\cDSet}.
\end{equation}
Due to the boundedness of the involved functions, the set $\mathfrak{F}$ is bounded as well and $\sup_{x\in\mathfrak{F}}\Norm{x}$ is finite.
Thus, for  $\bar{y}\geq \sup_{x\in\mathfrak{F}}\Norm{x}$ and given numbers~$\bar \eta, \bar{r},\bar{s},\bar \gamma,\bar{p}, \bar{d} \ge 0$, 
we define the following set of matrices, where we do not indicate the dependence on the parameters. Let 
\begin{equation*}
    \bar \cK  \coloneqq  
    \rbl\R^{m \times m}\rbr^r  \times \R^{m \times \nu}  \times \GL_m(\R)  \times \R^{m}  \times \SGroup_{\nu}^{--}  \times \R^{\nu \times m}  \times \R^{\nu}  \times {\R^\nu} ,
\end{equation*}
and define
\begin{equation}\label{eq:SetOfRetrictions}
    \mathcal{K} \coloneqq  \setdef{(R_1,\ldots,R_r,S,\Gamma,D_1,Q,P,D_2,\eta^0) \in \bar \cK}{\eqref{eq:ConditionsMatrices}} ,
\end{equation}
where
\begin{equation} \label{eq:ConditionsMatrices}
    \begin{aligned}
    \| S \| &\le \bar{s}, 
    &\| \Gamma\|, \| \Gamma^{-1}\|&\leq \bar{\gamma},
    &\| D_1\|, \| D_2\| &\le \bar d,\qquad
    \| P \| \le \bar{p},\\
    \| \eta^0 \| & \le \bar \eta,
    &\lambda_{\max}(Q) & \le - \frac{\bar{p} \bar y + \bar{d}}{\bar \eta},
    &\| R_i \| &\le \bar{r}\ \text{ for all } i=1,\ldots, r.\\
    \end{aligned}
\end{equation}
Then, we may derive the following statement. 
    
\begin{prop}\label{Prop:RestrictedMod}
    Let $y_{\rf} \in W^{k,\infty}(\Rp,\R^m)$ and $\Psi=(\Funnel_1,\ldots,\Funnel_r)\in \FunnelBoundaryFuncs$ 
    with associated parameters $k_i\geq 0$ for $i=1,\ldots,r-1$, and $\bar y \geq  \sup_{x\in\mathfrak{F}}\Norm{x}$.
    Further, let~$\bar \eta, \bar{r},\bar{s},\bar \gamma,\bar{p}, \bar{d} \ge 0$ be given
    and define recursively 
    $ \mu_i^0 \coloneqq  \SNorm{\Funnel_i}$,
    $\mu_{i}^{j+1}\coloneqq  \mu_{i+1}^{j}+k_{i}\mu_{i}^{j}$
    for ${i=1,\ldots, r-1}$ and $j=0,\ldots,r-i-1$.
    Choose 
    \[
        \umax\geq\bar{\gamma}
        \rbl
            r\bar{r}\bar{y} +\bar{s}\bar{\eta} +\bar{d}+\SNorm{y_{\rf}^{(r)}}+\sum_{j=1}^{r-1}k_j\mu_{j}^{r-j}+\SNorm{\dot{\Funnel}_r}
        \rbr,
    \] 
    and 
    \[
        \bar{\rho}\geq
            r\bar{r}\bar{y} +\bar{s}\bar{\eta} +\bar{d}+\bar{\gamma}\umax.
    \] 
    Then, the set~$\mathcal{K}$ defined in~\eqref{eq:SetOfRetrictions} satisfies the implication
\begin{multline*}
    (R,S,\Gamma,D_1,Q,P,D_2,\eta^0) \in \mathcal{K}
    \implies 
    (f_{R,S,D_1}, g_{\Gamma}, \oT_{R,S,Q,P,D_2,\eta_0}) \in \cM^{m,r}_{t_0, \umax,\bar{\rho}}, 
\end{multline*}
    where $f_{R,S,D_1}$, $g_{\Gamma}$, and $\oT_{R,S,Q,P,D_2,\eta_0}$ 
    are defined as in~\eqref{eq:f_ParameterAndg_Parameter} and~\eqref{eq:oTM_Parameter}.
\end{prop}
\begin{proof}
    Let $ (R,S,\Gamma,D_1,Q,P,D_2,\eta^0) \in \mathcal{K}$ be arbitrary but fixed.
    
    \noindent 
    \emph{Step 1}: Repeating the arguments from \Cref{Ex:LTISystem},
    one can easily see that 
    \[
        (f_{R,S,D_1}, g_{\Gamma}, \oT_{R,S,Q,P,D_2,\eta_0}) \in \cM^{m,r}_{t_0}.
    \]
    \noindent 
    \emph{Step 2}: We show properties~\ref{Item:LearningExControl} and~\ref{Item:LearningFCBounded} from~\Cref{Def:RestrictedModelClass}.
    Following the reasoning from \Cref{Rem:LearningScheme}~\ref{Item:Rem:LearningControlNonEmpty}, it is  sufficient to show
    that 
    \[
     r\bar{r}\bar{y} +\bar{s}\bar{\eta} +\bar{d}\geq\SNorm{f_{D_1}(\oT_{R,S,Q,P,D_2,\eta_0}(\zeta)|_{[0,s)})}, 
    \]
    and
    \[
    \bar{\gamma}\geq\SNorm{g_{\Gamma}(\oT_{R,S,Q,P,D_2,\eta_0}(\zeta)|_{[0,s)})},\quad
    \bar{\gamma}\geq\SNorm{g_{\Gamma}(\oT_{R,S,Q,P,D_2,\eta_0}(\zeta)|_{[0,s)})^{-1}} 
    \]
    for all $s\in[t_0,\infty]$ and $\zeta\in\FunnelTrajectories_{s}$.
    The last two inequalities are trivially fulfilled due to the definition of $g_{\Gamma}$ and $\bar{\gamma}$.
    We show that the first inequality is also satisfied. 
    To this end,  let $s\in[t_0,\infty]$ and $\zeta=(\zeta_1,\ldots,\zeta_r)\in\FunnelTrajectories_{s}$ 
    with $\zeta_i\in\cR(\Rp,\R^m)$ for $i=1,\ldots, r$ be arbitrary but fixed.
    By construction of $\bar{y}$, we have $\Norm{\zeta(t)}\leq \bar{y}$ for all $t\in[t_0,s]$. 
    Let $\eta(\cdot;t_0,\eta^0,\zeta_1)$ be the maximal solution of the initial value problem  
    \[ 
        \dot{\eta}(t)=Q \eta(t)+D_2 + P \zeta_1(t), \quad \eta(t_0)=\eta^0.
    \]
    For $t\in[t_0,s]$, we calculate  
    \begin{align*}
        \dd{t}\tfrac{1}{2}\Norm{\eta(t;t_0,\eta^0,\zeta_1)}^2
        &=\eta(t;t_0,\eta^0,\zeta_1)\rbl Q\eta(t;t_0,\eta^0,\zeta_1) +P\zeta_1(t)+D_2\rbr\\
        &\leq \Norm{\eta(t;t_0,\eta^0,\zeta_1)}\rbl\lambda^+(Q)\Norm{\eta(t;t_0,\eta^0,\zeta_1)}+\bar{p}\bar{y}+\Norm{D_2}\rbr,
    \end{align*}
    which is non-positive for $\Norm{\eta(t;t_0,\eta^0,\zeta_1)}\geq(\bar{p}\bar{y}+\Norm{D_2})/\Abs{\lambda^+(Q)}$ as  $\lambda^+(Q)<0$.
    Therefore, \cite[Thm.~4.3]{Lanz21} yields
    \[
        \Norm{\eta(t;t_1,\eta^0,\zeta_1)}\leq \max\cbl(\bar{p}\bar{y}+\bar{d})/\Abs{\lambda^+(Q)},\Norm{\eta^0}\cbr
    \]
    for all $t\in [t_0,s]$.
    By assumption~\eqref{eq:ConditionsMatrices}, we have $\|\eta^0\|\leq \bar{\eta}$ and $\Abs{\lambda^+(Q)}\geq (\bar{p}\bar{y}+\bar{d})/\bar{\eta}$.
    Hence, $\Norm{\eta(t;t_1,\eta^0,\zeta_1)}\leq \bar{\eta}$ for all $t\in [t_0,s]$.
    As $\oT_{R,S,Q,P,D_2,\eta_0}(\zeta)= \sum_{j=1}^{r}R_j\zeta_j+S\eta(\cdot;t_1,\eta^0,\zeta_1)$, we estimate
    \[
    \Norm{\oT_{R,S,Q,P,D_2,\eta_0}(\zeta)(t)} \leq r\bar{r}\bar{y} +\bar{s}\bar{\eta}.
    \]
    for all $t\in [t_0,s]$. Thus, $\Norm{f_{R,S,D_1}(\oT_{R,S,Q,P,D_2,\eta_0}(\zeta)(t))}\leq r\bar{r}\bar{y} +\bar{s}\bar{\eta}+\bar{d}$ for all $t\in[t_0,s]$. 
    As~$s\in[t_0,\infty]$ and $\zeta\in\FunnelTrajectories_{s}$ are arbitrarily chosen,  
    this shows 
    \[
        {(f_{R,S,D_1},g_{\Gamma}, \oT_{R,S,Q,P,D_2,\eta_0}) \in \cM^{m,r}_{t_0, \umax,\bar{\rho}}}
    \] 
    and completes the proof.
\end{proof}

With the set of parameters~$\cK$, the functions $f_{D_1}, g_{\Gamma},\oT_{R,S,Q,P,D_2,\eta_0}$, and 
defined in~\eqref{eq:f_ParameterAndg_Parameter} and~\eqref{eq:oTM_Parameter},
and \Cref{Prop:RestrictedMod}, we may define a learning scheme $\cL$ mapping from $\bigcup_{\hat{t}\geq t_0}\SetOfSignals_{\hat{t}}$
to the subset
\[
\setdef{ (f_{D_1}, g_{\Gamma},\oT_{R,S,Q,P,D_2,\eta_0}) }
{(R,S,\Gamma,D_1,Q,P,D_2,\eta^0)  \in  \cK }  
\]
of $\cM^{m,r}_{t_0, \umax,\bar{\rho}}$, defined by
\begin{align*}
    \cL:\ & ((y,\xM, \uFMPC, \uFC)|_{[t_0,\hat{t}]}) 
    \mapsto (f_{D_1}, g_{\Gamma},\oT_{R,S,Q,P,D_2,\eta_0})
\end{align*}
for some $\hat{t}\ge t_0$, where $(f_{D_1}, g_{\Gamma},\oT_{R,S,Q,P,D_2,\eta_0})$ 
is determined by the solution of an optimisation problem involving measurements of the system data $y$
and the applied control signals $\uFMPC$ and $\uFC$ over the time interval $[t_0,\hat{t}]$ of the form 
\begin{equation}\label{eq:LearningOP}
    \begin{alignedat}{2}
             &\hspace{-1cm}\mathop {\operatorname{minimise}}_{\substack{(R,S,\Gamma,D_1,Q,P,D_2,\eta^0) \in \cK}}  \
             J((y,z)|_{[t_0,\hat{t}]})\\
            \text{s.t.}\ 
                      &\OpChi(z)(t_0) =  z^0,\\
                      &{\makebox[\widthof{$\OpChi(z)(t_0)$}]{$z^{(r)}(t)$}} =  f_{D_1}(\oT_{R,S,Q,P,D_2,\eta_0}(\OpChi(z))(t))\\
                      &\hspace{1.75cm}+ g_{\Gamma}(\oT_{R,S,Q,P,D_2,\eta_0}(\OpChi(z))(t))(\uFMPC+\uFC)(t),
    \end{alignedat}
\end{equation}
where $J(\cdot)$ is a suitable cost function.
Here, $\hat{t}$ refers to time of the execution of the learning algorithm, 
i.e. the current time instant $t_k$ during operation of the learning-based robust funnel MPC~\Cref{Algo:LearningRFMPC}.
Note that solving the differential equation  
$z^{(r)}(t) =  f_{D_1}(\oT_{R,S,Q,P,D_2,\eta_0}(\OpChi(z))(t))+ g_{\Gamma}(\oT_{R,S,Q,P,D_2,\eta_0}(\OpChi(z))(t))(\uFMPC+\uFC)(t)$
is equivalent to solving the linear differential equation with the state~\eqref{eq:LearningLinearExample}.
\begin{remark}
In application, measurement of the system data $y$ (and its derivatives) and the applied control signals $\uFMPC$ and $\uFC$ is only 
available at discrete time instants~${t_0+i\SampleTime}$ with $\SampleTime>0$ and $i\in\N$.
In this case, it is reasonable to replace the control used in constraints~\eqref{eq:LearningOP} by 
the piecewise constant $\tilde{u}$ defined as 
\[
     \tilde{u}(t)=(\uFMPC+\uFC)(t_0+(i-1)\SampleTime),
\]
for $t\in[t_0+(i-1)\SampleTime,t_0+i\SampleTime )$ and all $i\leq (\hat{t}-t_0)/\SampleTime$,
and use a cost function~$J(\cdot)$ which evaluates $y$ and $z$ only at time instants $t_0+i\SampleTime$. 
In the following \Cref{Chapter:DiscretFMPC}, we will discuss the matter of using piecewise constant control signals 
for the overall control problem in more detail. 
However, we want to discuss, in the following, 
some possible choices for the cost function~$J(\cdot)$ when only discrete measurements are available.
\begin{enumerate}[(a)]
    \item $J((y,z)|_{[t_0,\hat{t}]}) \coloneqq  \sum_{i=0}^{\lfloor (\hat{t}-t_0)/\SampleTime\rfloor} a_i \| 
    \OpChi(z)(t_0+i\SampleTime) - \OpChi(y)(t_0+i\SampleTime)\|^2$ with weights $a_i \ge 0$. 
    The idea is to find a model in the set~$\cK$ which minimises the weighted squared measured output errors.
    The weights $a_i$ reflect the relative importance of the measurements $\OpChi(y)(t_0+i\SampleTime)$. 
    In certain cases, it might be beneficial to weight data points that are far in the past lower than current data points. 
    By choosing $a_i>0$ for all $i>0$, all measured past data is taken into account.
    With increasing runtime of the algorithm, this results in a growing complexity of the optimisation problem,
    computation time,  and required memory space for the measurements.
    Therefore, this is not suitable in practice.
    Thus, it is beneficial to use a moving horizon estimation approach and 
    only take the last $N$ measurements into account and set $a_{i}=0$ for $i<\lfloor (\hat{t}-t_0)/\SampleTime\rfloor-N$.
    In application, one has to find a good balance between considering many data points (large $N$),
    thus having a probably more accurate model, and low computation time and memory requirements (small $N$).
    This is comparable to the sliding window approach discussed in 
    \Cref{Rem:LearningScheme}\ref{Rem:Item:ModDomain}\ref{Rem:Item:ModDomainSlidingWindow}.
    
    \item If the computation of the solution of the optimisation problem has to be done very quickly,
    it is also possible to  only consider the last measurement $y(t_0+\lfloor (\hat{t}-t_0)/\SampleTime\rfloor\SampleTime)$.
    Thus, one might choose the cost function   
    \[
    J((y,z)|_{[t_0,\hat{t}]}) \coloneqq  \| \OpChi(z)(t_0+\lfloor (\hat{t}-t_0)/\SampleTime\rfloor\SampleTime) - \OpChi(y)(t_0+\lfloor (\hat{t}-t_0)/\SampleTime\rfloor\SampleTime)\|^2.
    \]
    The idea is to find a model, which best explains the last MPC period in terms of output error,
    i.e., a model on the prediction interval $[t_k,t_{k+1}]$ so that, with $\SampleTime = \delta$,
    the error $\| \OpChi(z)(t_{k+1}) - \OpChi(y)(t_{k+1})\|$ at the end of the interval is minimal. 
    
    \item In addition, it is worth considering to include the used model in the cost function as discussed in
    \Cref{Rem:LearningScheme}\ref{Rem:Item:ModDomain}\ref{Rem:Item:ModDomainModel}.
    This can be done by adding regularisation terms for the model parameters in the cost function. 
    For the parameter vector ${\cK_i=(R_i,S_i,\Gamma_i,D_{1,i},Q_i,P_i,D_{2,i},\eta^0_i) \in  \cK}$,
    one could either penalise the weighted distance of $\cK_i$ to a priori known parameters 
    $\cK^\star=(R^\star,S^\star,\Gamma^\star,D^\star,Q^\star,P^\star,D_2^\star,{\eta^0}^\star)$ 
    and thus allow only small adaptions of the a priori known model or 
    penalise the change of parameters $\cK_i$ such that the model does only change slightly between two learning steps. 
    This results in a cost function of the form 
    \[
    J((y,z)|_{[t_0,t]}) \coloneqq  \sum_{i=0}^{\lfloor (\hat{t}-t_0)/\SampleTime\rfloor} \big (a_i \|  \OpChi(z)(t_0+i\SampleTime) - \OpChi(y)(t_0+i\SampleTime)\|^2 + \sum_{j=1}^8b_i^j\|(\cK_i^j-\tilde{\cK}^j)\|\big ),
    \]
    where $\tilde{\cK}=\cK^\star$ or $\tilde{\cK}=\cK_{i-1}$ and
    with weights $a_i ,b_i^j \ge 0 $.
    Here the expressions $\cK_i^j, \tilde \cK_i^j$ with $j=1,\ldots,8$ refer to the
    $j^{\mathrm{th}}$ entry of the tuple $\cK_i$, $\tilde{\cK}_i$, respectively; for instance, $\cK_i^2 = S_i$.
\end{enumerate}
\end{remark}
    
\begin{remark}
    The bounds for~$\umax$ and $\bar{\rho}$ derived~\Cref{Prop:RestrictedMod} are rather conservative and can clearly be improved.
    However, \Cref{Prop:RestrictedMod} exemplifies how to construct a subset of models belonging to $\cM^{m,r}_{t_0, \umax,\rho}$
    by prescribing bounds  ${\fMmax, \gMmax, \gMInvmax\geq0}$ on the dynamics.
    \Cref{Prop:RestrictedMod} relies on the following abstract idea.
    For a compact set $K\subset \R^\nu$, choose a set of operators $\cT\sub\cT^{rm,\nu}_{t_0}$ with 
    \begin{equation}\label{eq:ConditionCompactSet}
        \fa \oTM \in\cT\ \fa\zeta\in\FunnelTrajectories_{\infty}:\quad \oTM(\zeta)(\Rp)\subset K.
    \end{equation}
    Moreover, consider only functions $\fM \in (\R^\nu ,\R^m)$ and $\gM \in (\R^\nu ,\R^{m\times m})$ satisfying ${\gM(z)\in\GL_{m}(\R)}$ for all $z\in\R^\nu$ and
    \[ 
          \Norm{\vphantom{\gM(x)^{-1}}\fM(x)}\leq\fMmax,\quad\Norm{\vphantom{\gM(x)^{-1}}\gM(x)}\leq\gMmax,\quad\Norm{\gM(x)^{-1}} \leq\gMInvmax 
    \]
    for all $x\in K$. Using this approach, one can construct a set of models of the form 
    \begin{subequations}
        \begin{align}
        \yM^{(r)}(t) &= p\big(\OpChi(\yM)(t),\eta(t)\big) + \Gamma\big(\OpChi(\yM)(t),\eta(t)\big)\,u(t),\\
        \dot \eta(t) &= q\big(\OpChi(\yM)(t),\eta(t)\big),\label{eq:CompactSetEtaEquation}
        \end{align}
    \end{subequations}
    with $p:\R^{rm\times \nu}\to \R^{m}$, $q: \R^{rm\times\nu}\to \R^{\nu}$, $\Gamma :\R^{rm\times\nu}\to\R^{m\times m}$,
    belonging to $\cM^{m,r}_{t_0, \umax,\bar{\rho}}$ where $\umax$ and $\bar{\rho}$ are given as in
    \Cref{Rem:LearningScheme}~\ref{Item:Rem:LearningControlNonEmpty}.
    We already saw in \Cref{Ex:ControlAffineMod} that these models belong to $\cM^{m,r}_{t_0}$.
    The main difficulty lies in constructing a compact set  $K\subset \R^\nu$ and ensuring \eqref{eq:ConditionCompactSet}.
    The matter ultimately comes down to finding a uniform bound $\bar{\eta}\geq 0$ of 
    \[
        \Norm{\eta(t;t_0,\eta^0,\zeta)}\leq \bar{\eta}
    \]
    for all $t\geq t_0$ and all $\zeta\in\FunnelTrajectories_{\infty}$ where $\eta(t;t_0,\eta^0,\zeta)$
    is the global solution of the equation~\eqref{eq:CompactSetEtaEquation} 
    where $\OpChi(\yM)(t)$ is replaced by $\zeta$.
    One way to verify the satisfaction of such a uniform bound is 
    to apply \cite[Thm.~4.3]{Lanz21}, which states the following.
    Assume there exists ${V \in \cC^1( \R^\nu , \Rp)}$ 
    with $V(\eta) \to \infty$ as $\|\eta\| \to \infty$ and, for $q \in \cC(\R^{rm}\times \R^\nu , \R^\nu)$, 
    ${V'(\eta) \cdot q(z,\eta) \le 0}$ for all $z\in\mathfrak{F}$ as in~\eqref{eq:FunnelFrak} 
    and $\eta\in\R^\nu$ with $\Norm{\eta}>\tilde\eta$ for a predefined value $\tilde{\eta}\geq 0$.
    Then, ${\Norm{\eta(t;t_0,\eta^0,\zeta)}\leq \max\cbl\eta^0,\tilde{\eta}\cbr}$ for all $t\geq t_0$, all $\eta^0\in\R^\nu$,
    and all $\zeta\in\FunnelTrajectories_{\infty}$.
    Hence, fixing $V(\cdot)$ and $\tilde{\eta}>0$ in advance can be used to restrict choices of~$q(\cdot)$ satisfying
    $\SNorm{\eta(\cdot;t_0,\eta^0,\zeta)}\leq \bar{\eta}$.
    We made use of this fact in the proof of~\Cref{Prop:RestrictedMod}.
\end{remark}

\section{Simulation}
In this section, we  illustrate the application of the learning-based robust
funnel MPC~\Cref{Algo:LearningRFMPC} to the numerical examples
from~\Cref{Sec:FMPC:Simulation}.
The \textsc{Matlab} source code for the performed simulations 
can be found on \textsc{GitHub} under the link \url{https://github.com/ddennstaedt/FMPC_Simulation}.

\subsection*{Exothermic chemical reaction}
To demonstrate the functioning of the robust funnel MPC~\Cref{Algo:FunnelMPC}, 
we revisit the example of a continuous chemical reactor from \Cref{Sec:FMPC:Sim:Reactor}.
The system is described by the following non-linear differential equation:
\begin{equation}\tag{\ref{eq:ExampleExothermicReaction} revisited}
    \begin{aligned}
        \dot{x}_1(t) &= c_1\, p(x_1(t), x_2(t), y(t)) +d(x_1^{\mathrm{in}}-x_1(t)),\\
        \dot{x}_2(t) &= c_2\, p(x_1(t), x_2(t), y(t)) +d(x_2^{\mathrm{in}}-x_2(t)),\\
        \dot{y}(t)   &= b\,   p(x_1(t), x_2(t), y(t)) -q\,y(t) + u(t),
    \end{aligned}
\end{equation}
where the function $p$ is the Arrhenius law \eqref{eq:Ex:ArrheniusLaw},
the parameters are given in \eqref{eq:Ex:Reactor:Sys:Params}, and
the initial data is $[x_1^0, x_2^0, y^0] = [0.02,0.9,270]$.
Following the given heating profile $y_{\rf}(t)$ given in~\eqref{eq:Ex:Reactor:RefTemp}
within boundaries defined by the funnel function~$\Funnel(t)\coloneqq20\me^{-2t}+4$, 
the control objective is to steer the reactor's temperature $y$ to a certain
desired constant value $y_{\rf ,\mathrm{final}}$. 
To achieve the control objective with the learning-based robust funnel MPC~\Cref{Algo:LearningRFMPC},
we consider linear models of order $r=1$ of the form~\eqref{eq:LearningLinearExample}
with $R,  D_1\in\R$, $S,D_2^\top, P^\top\in\R^{1\times 2}$, and $Q\in\R^{2\times 2}$.
To learn the model from the measured data, we use linear regression subject to
the constraints introduced in \Cref{Def:RestrictedModelClass,Def:LearningScheme}. Hence, feasibility
of the data-based models is guaranteed by \Cref{Prop:RestrictedMod}.
We assume $\Gamma= 1$ and, as initial model, we choose ${R=D_1=0\in\R}$,
${S=D_2^\top= P^\top=0\in\R^{1\times 2}}$, ${Q=0\in\R^{2\times 2}}$, and
${\eta^0=[x_{1}^0,x_{1}^0]=[0.02,0.9]}$, which represents an integrator chain with decoupled
internal dynamics.
To improve this (deliberately poorly chosen) model over time, we adapt the matrices over a compact set~$\cK$
as in~\eqref{eq:SetOfRetrictions} at every fifth time step $t_k$
by minimising the model-plant mismatch based on the data of the last system output $y(t_{k-1})$, i.e.
we solve the optimisation problem
\begin{equation*}
    \begin{aligned}
       &\mathop{\operatorname{minimize}}_{\mathclap{\substack{(R,S,1,D_1,Q,P,D_2,\eta(0))\in \cK}}}
            \qquad\quad\|{\yM(t_k)-y(t_k)}\|^2 \\
        &\begin{matrix}\textrm{s.t.} \\
        \phantom{x}\\
        \phantom{x}\\
        \phantom{x}\\
        \end{matrix}\qquad
        \begin{aligned}
       \frac{\textrm{d}}{\textrm{d} t} \begin{bmatrix}
        \eta(t)\\
        \yM(t)
        \end{bmatrix}  &=  
       \begin{bmatrix}
        S & R\\
        Q & P
        \end{bmatrix}
       \begin{bmatrix}
            \eta(t) \\
            \yM(t)
        \end{bmatrix} 
         + 
       \begin{bmatrix}
        0\\
        1
        \end{bmatrix}
         u(t) +
       \begin{bmatrix}
        D_1 \\
        D_2 
        \end{bmatrix},\\ 
        \begin{bmatrix}
                \eta(t_{k-1})\\
                \yM(t_{k-1})
            \end{bmatrix}  &= 
         \begin{bmatrix}
             \eta(0)\\
             y(t_{k-1})
         \end{bmatrix},
         \end{aligned}
    \end{aligned}
\end{equation*}
where $u(t)= \uFMPC(t_{k-1})+\uFC(t_{k-1})$ which was applied to the model at the last time step~$t_{k-1}$
and $\eta(0)\coloneqq [x_{1}^0,x_{1}^0]$ is the vector of initial concentrations of the substances~$x_1$ and~$x_2$.
As before, we choose the strict funnel stage cost 
${\ell_{\Funnel}:\Rp\times\R^m\times\R^{m}\to\R\cup\{\infty\}}$ defined in  
\eqref{eq:Ex:ReactorFMPCCostFunction} with $\lambda_{u}=10^{-4}$, the prediction
horizon $T=1$, and time shift $\delta=0.1$ for the funnel MPC component of the control algorithm 
and restrict the OCP \eqref{eq:LearningRobustFMPCOCP} to step functions with a constant step length of 
$\delta= 0.1$.
We choose for the set~$\cK$ as in~\eqref{eq:SetOfRetrictions} 
the parameters in \eqref{eq:ConditionsMatrices} as 
$\bar{r}=1.3$,
$\bar{s}=1.4$,
$\bar{\eta} =0.91$,
$\bar{\gamma}=1$,
$\bar{p}=1/400$,
$\bar{d} = 2.5$,
and $\bar{y}=341.4$.
We have $\SNorm{\dot{y}_{\rf}}=33.55$ given by the heating profile and
$\|\dot{\Funnel}\|_{\infty}=40$ by choice of the funnel function.
Thus, we  restrict the funnel MPC control signal to $\SNorm{\uFMPC}\leq \umax\coloneqq 600$ to satisfy  
the requirements of  \Cref{Prop:RestrictedMod} for~$\bar{\rho}=1125$. 
The learning scheme is therefore ($\umax$,$\bar{\rho}$)-feasible.
For the control law of funnel control component, we choose the bijection $\FCBijec(s)=1/(1-s)$
and the function $\FCSurjec(s)=-10s$. 
\begin{figure}[h]
    \begin{subfigure}[b]{0.49\textwidth}
        \centering
        \includegraphics[width=\textwidth]{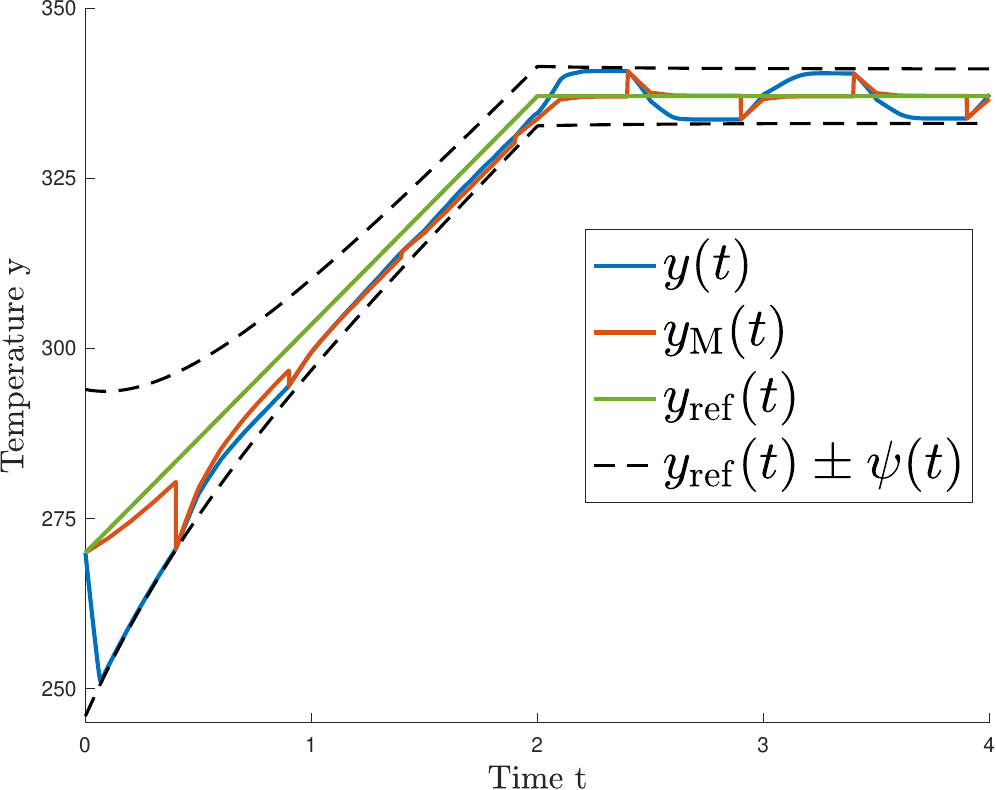}
        \caption{Outputs and reference, with boundary~$\Funnel$.}
        \label{fig:sim:learning_fmpc:output}
    \end{subfigure}
      \begin{subfigure}[b]{0.49\textwidth}
        \centering
        \includegraphics[width=\textwidth]{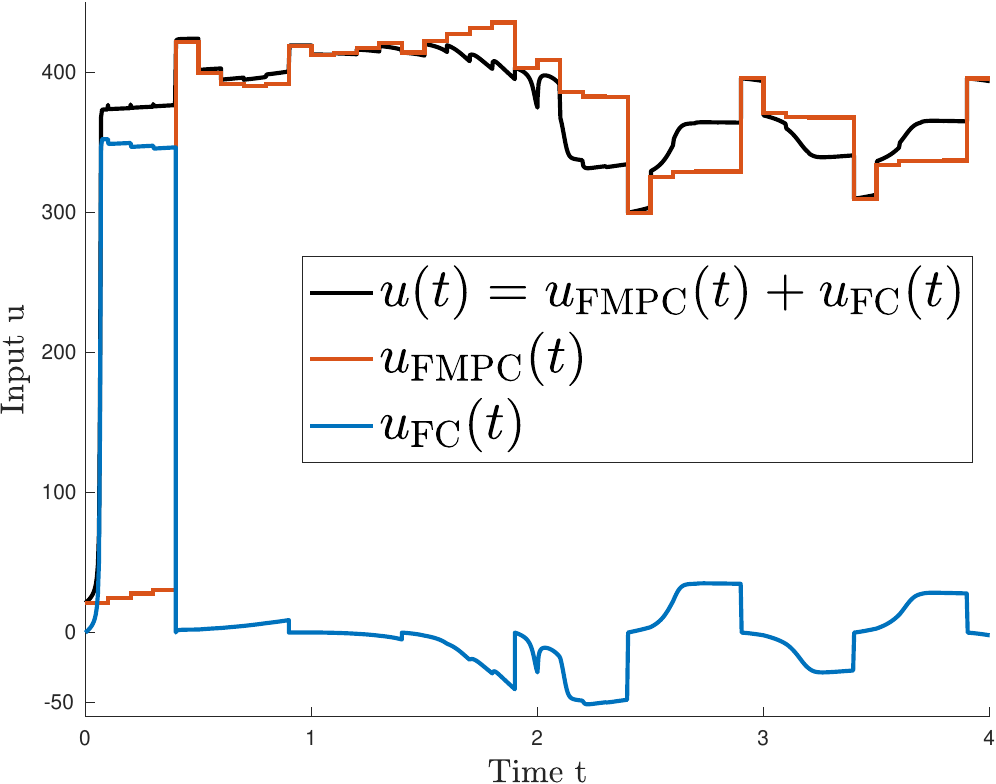}
        \caption{Control inputs.}
        \label{fig:sim:learning_fmpc:input}
    \end{subfigure}
    \caption{Simulation of system~\eqref{eq:ExampleExothermicReaction} under the
    control generated by the learning-based robust funnel MPC~\Cref{Algo:LearningRFMPC}
    with model update every five iterations.}
    \label{fig:sim:learning_fmpc}
\end{figure}

\Cref{fig:sim:learning_fmpc} shows the control signals and the system and model output errors, respectively.
It is evident that both $\yM-y_{\rf}$
and $y-y_{\rf}$ remain within the predefined funnel boundaries~$\Funnel$.
Before the first learning step for $t \in [0, 0.5)$, the tracking error $y-y_{\rf}$ and the
predicted error~$\yM-y_{\rf}$  diverge due to the poor quality of the initial model.
However, since the tracking error is not close to the funnel boundary,
the funnel controller remains inactive in the beginning and 
only reacts when the tracking error is close to the boundary.
After the first learning step, the general direction of the predicted tracking error is consistent 
with the actual tracking error. The funnel controller still has to slightly
compensate for the model inaccuracies in order to guarantee that the tracking
error remains within the boundaries, but with a significantly smaller
contribution to the control signal.
After each learning step, the model output jumps $\yM$ to the system output $y$
due to the newly updated model. 
The control signal $\uFC$ is zero after each learning step 
since the system and model output  coincide, 
and it becomes larger afterwards to compensate for the model inaccuracy.
After the heating phase, the model, only being updated every five iterations of the MPC algorithm,
does not adequately describe the system dynamics. 
The funnel controller therefore has to compensate these inaccuracies during the whole  
operation of the algorithm, but with a significantly smaller control signal than before the first 
learning step.
\begin{figure}[h]
    \begin{subfigure}[b]{0.49\textwidth}
        \centering
        \includegraphics[width=\textwidth]{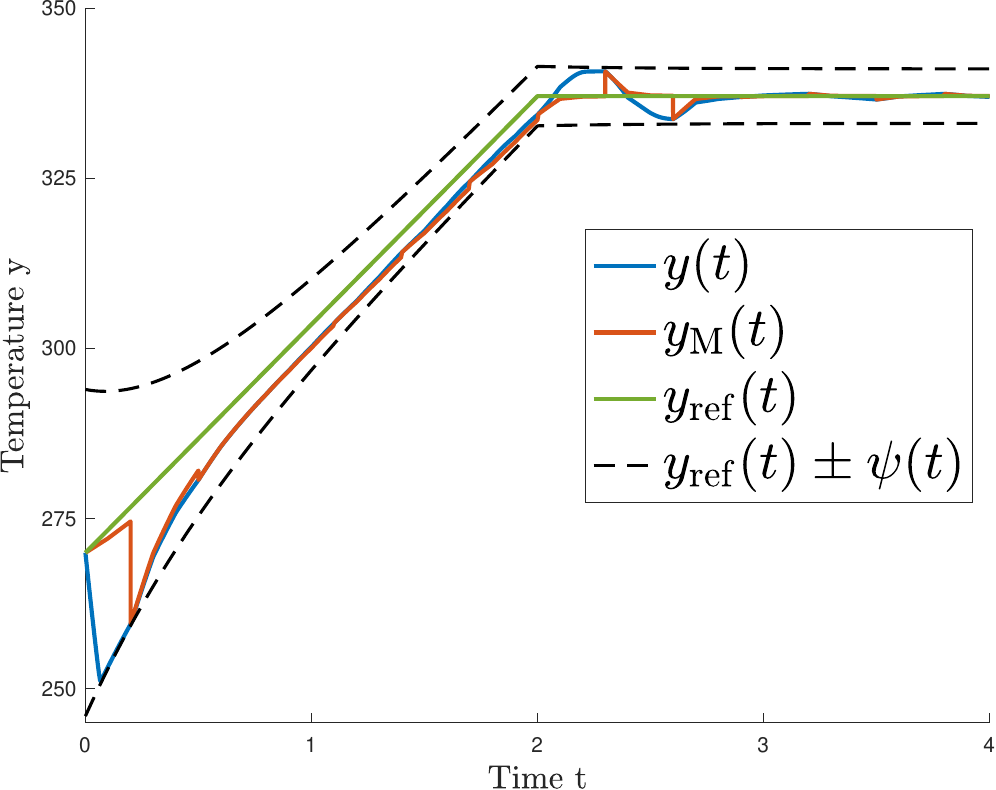}
        \caption{Outputs and reference, with boundary~$\Funnel$.}
        \label{fig:sim:learning_fmpc_often:output}
    \end{subfigure}
      \begin{subfigure}[b]{0.49\textwidth}
        \centering
        \includegraphics[width=\textwidth]{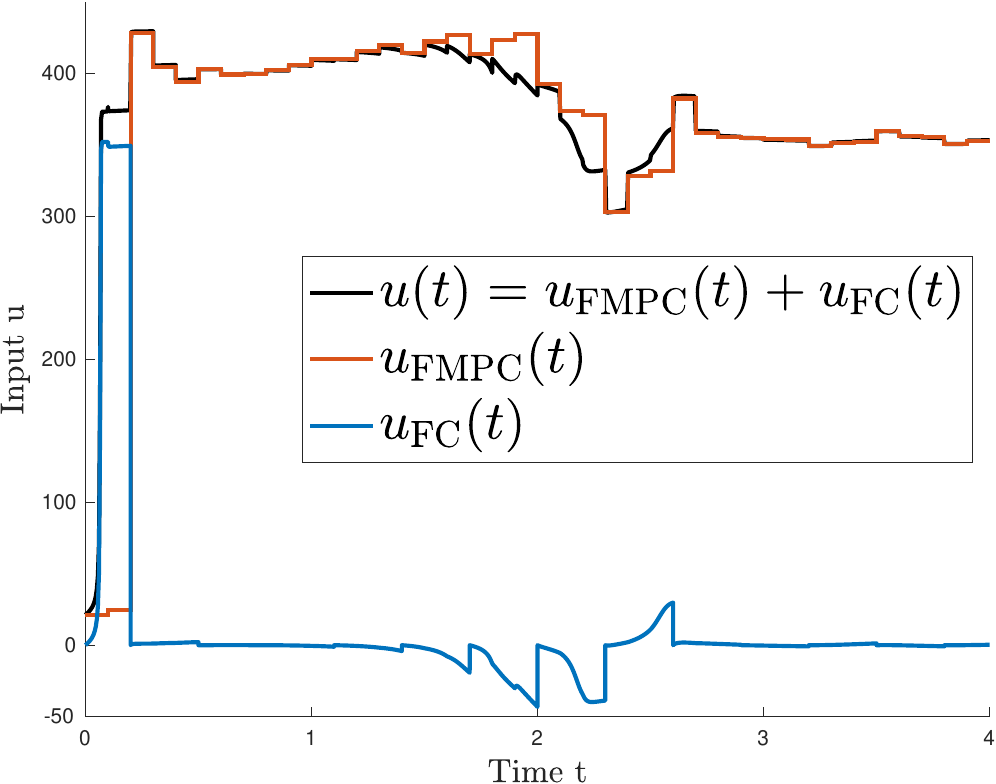}
        \caption{Control inputs.}
        \label{fig:sim:learning_fmpc_often:input}
    \end{subfigure}
    \caption{Simulation of system~\eqref{eq:ExampleExothermicReaction} under the
    control generated by the learning-based robust funnel MPC~\Cref{Algo:LearningRFMPC}
    with model update every three iterations.}
    \label{fig:sim:learning_fmpc_often}
\end{figure}
In a second simulation, we update the model every third time step $t_k$ instead of every fifth 
but leave rest of the controller configuration unchanged. The results are depicted in \Cref{fig:sim:learning_fmpc_often}.
As one can see, the combined controller is able to achieve the control objective.
Before the initial learning step, the funnel controller has to compensate for the inaccuracies 
of the model with a large control signal comparable to the setting before. 
Already after the first update of the model,
the principal portion of the control signal is generated by the MPC component.
The funnel controller only has to intervene during the transition heating phase
(before $t_{\mathrm{final}} = 2$) to the constant temperature phase of the
system (after $t_{\mathrm{final}} = 2$). Thenceforth, the linear model is
adequate to predict the system behaviour and the control signal computed by
funnel MPC is sufficient to achieve the tracking objective.
In contrast to the case before, the funnel controller remains mainly inactive 
after $t\approx2.6$. This shows that the ``quality'' of the learning scheme 
and the update frequency of the model can have a significant impact on the controller
behaviour and its performance.
The more accurate the model is, the less control is required  by the funnel
controller to mitigate the model-system mismatch. However, updating the model
more frequently can lead to increased computation costs.

We note that this example merely serves as an illustration that the learning-based
robust funnel MPC~\Cref{Algo:LearningRFMPC} can be combined with any
($\umax$,$\bar{\rho}$)-feasible learning scheme~$\cL$. We do not claim that the
learning algorithm used is superior to other methods. 

\subsection*{Mass-on-car system}
To illustrate that the learning-based funnel MPC~\Cref{Algo:LearningRFMPC} can be successfully applied to systems with relative degree $r>1$, 
we revisit the example of the mass-on-car system from~\Cref{Sec:FMPC:Sim:MassOnCar}.
Assuming the mass $m_2=2$, on the ramp inclined by the angle $\vartheta=\frac{\pi}{4}$, is connected to the car with mass $m_1=4$ 
via a spring and damper system with spring constant $k=2$ and damper constant $d=1$,
the system can be described by the differential equation
\begin{equation}\tag{\ref{eq:MassOnCarInputOutputDeg2} revisited}
    \begin{aligned}
        \ddot y(t)    &= R_1y(t)+ R_2\dot{y}(t) + S\eta(t) +\Gamma u(t)\\
        \dot \eta(t)  &= Q \eta(t) + P y(t). 
    \end{aligned}
\end{equation}
with matrices given in \eqref{eq:ParamsMassOnCarInputOutputDeg2}.
The objective is to track the reference signal $y_{\rf}(t) = \cos(t)$ such
that the tracking error ${y(t)-y_{\rf}(t)}$ evolves within the prescribed
performance funnel given by the function $\Funnel\in\cG$ with $\Funnel(t) =
5\me^{-2t}+0.1$.  To achieve this  control objective with the learning-based robust
funnel MPC~\Cref{Algo:LearningRFMPC}, we use the strict funnel stage cost
function 
${\ell_{\Funnel_{2}}:\Rp\times\R\times\R\to\R\cup\{\infty\}}$ as defined in \eqref{eq:Ex:MassOnCar:CostFunction}.
For the simulation, the  MPC control signal is further restricted to ${\|\uFMPC
\|_\infty \le \umax=30}$ and we choose the design parameters $\lambda_u =
10^{-4}$, prediction horizon $T = 0.5$, and time shift $\delta = \tfrac{T}{20}= 0.025$.
For the model-free component of the controller, we use a slightly modified form of
the control law~\eqref{eq:Sim:MassOnCar:RobustFC}: 
\begin{equation*}
    \begin{aligned}
        w(t)&=\phi(t)\eSTrackdot(t)+\FCBijec\rbl\phi(t)^2\eSTrack(t)^2\rbr\phi(t)\eSTrack(t),& \eSTrack(t)&=y(t)-\yM(t),\\
        \uFC(t)&=-2\FCBijec\rbl w(t)^2\rbr w(t),& \phi(t)&=\frac{1}{\Funnel(t)-\Norm{\yM(t)-y_{\rf}(t)}},
    \end{aligned}
\end{equation*}
where $\yM$ is the prediction for the system output computed by the MPC component. 

\begin{figure}[h]
    \begin{subfigure}[b]{0.49\textwidth}
        \centering
        \includegraphics[width=\linewidth]{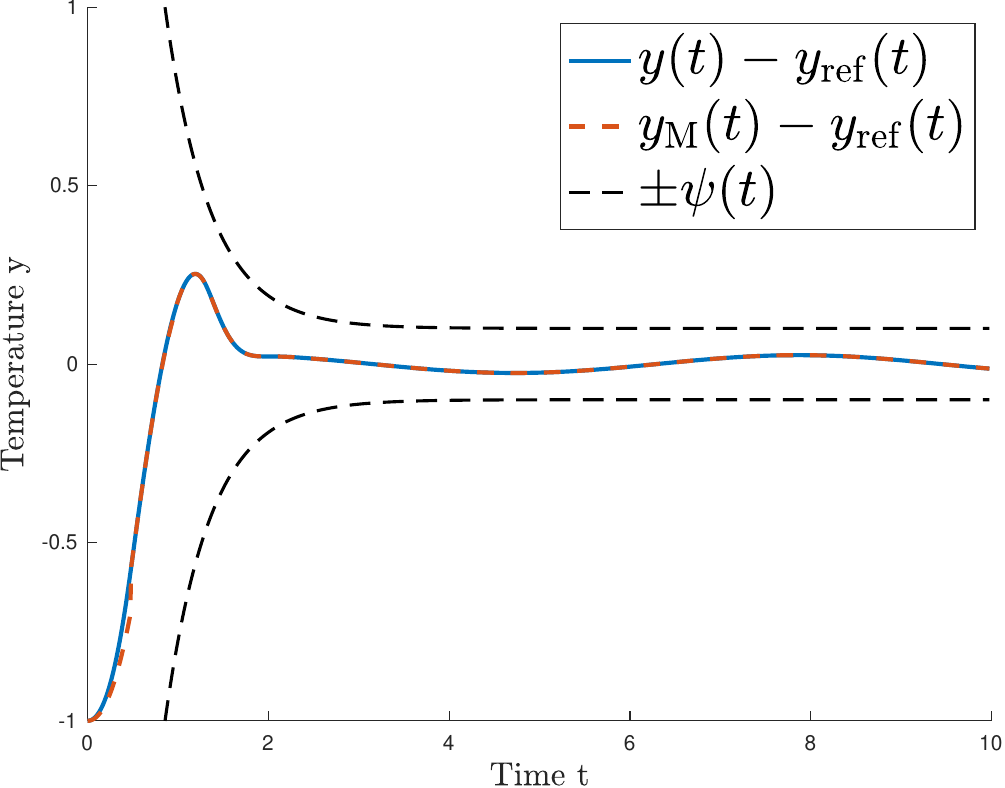}
        \caption{Tracking error $e=y-y_{\rf}$ within boundary~$\Funnel$.}
        \label{fig:sim:massoncar:learning:output}
    \end{subfigure}
      \begin{subfigure}[b]{0.49\textwidth}
        \centering
        \includegraphics[width=\linewidth]{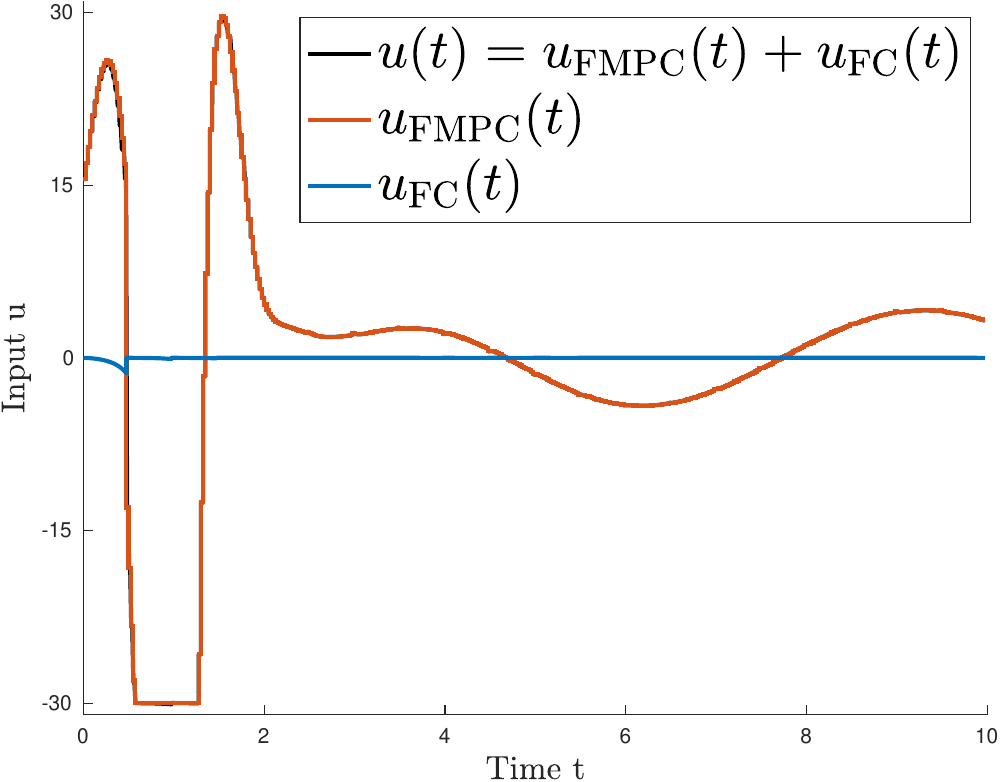}
        \caption{Control inputs.}
        \label{fig:sim:massoncar:learning:input}
    \end{subfigure}
    \caption{Simulation of system~\eqref{eq:MassOnCarInputOutputDeg2} under the
    control generated by the learning-based robust funnel MPC~\Cref{Algo:LearningRFMPC}
    } 
    \label{fig:sim:massoncar:learning}
\end{figure}

Similar to~\cite{berger2019learningbased}, where this problem was studied in the context of model identification 
for the learning component during runtime, 
we assume knowledge about the structure of the system, but only limited information about its parameters.
We assume to know $m_1,m_2\in [0.5,10]$ and $k,d\in[0.5,5]$.
As an initial model, we choose the parameters $m_1=6$, $m_2=2$, $k=3$,and $d=0.75$.
To learn or update the model parameters, we take measurements of the system's input-output data
$((\uFMPC+\uFC)(ih),y(ih))$ for $h =2.5\cdot10^{-4}$ and $i\in\N_{0}$
and update the model every twentieth iteration of the MPC algorithm,
i.e. at $t\in T\N$, by solving the optimisation problem
\begin{align*}
    \mathop {\operatorname{minimise}}_{\substack{m_1,m_2\in[0.5,10],\\k,d\ \in[0.5,5] }}  \quad
            & \sum_{i=0}^{2000j}\Norm{\tilde{y}_{\mathrm{M}}(i h) - y(i h)}^2\\
            \text{s.t.} \ \ 
                       z(0) &=  0  \text{ and for all }i=1,\ldots, 2000j:\\
                       z(i h )  &=  z(h; z((i-1)h ),(\uFMPC+\uFC)((i-1)h)),\\
                       \tilde{y}_{\mathrm{M}}(ih)&=[1, 0, 0, 0 ] z(ih),
\end{align*}
at every time $t =2000jh$ for $j\in\N$, where $z=[y,\dot{y},\eta_1,\eta_{2}]^\top$
denotes the state of the mass-on-car system~\eqref{eq:MassOnCarInputOutputDeg2} and
$z(\cdot; z((i-1)h ),(\uFMPC+\uFC)((i-1)h))$ denotes its solution under the 
initial condition $z(0) = z((i-1)h )$ and with constant control~$u(\cdot)\equiv(\uFMPC+\uFC)((i-1)h)$.
Since only the interval $[0,10]$ is considered for the simulation,
the entire history of input-output data is considered in the optimisation problem instead of a moving horizon approach.
After every execution of this learning scheme, the model is properly initialised by solving the
optimisation problem~\eqref{eq:Ex:MassOnCar:Sim:Initial}. 
Between two updates of the model, the MPC component's control signal $\uFMPC$ is applied to 
the system in open-loop fashion, i.e. the model is initialised with its state from the
previous iteration as initial value.
All simulations are depicted in~\Cref{fig:sim:massoncar:learning}.
It is evident that the control scheme is feasible and achieves the control
objective. Both errors~$\yM-y_{\rf}$ and~$y-y_{\rf}$ evolve within the
funnel boundaries given by~$\Funnel$, see \Cref{fig:sim:massoncar:learning:output}. 
While the model output $\yM$ and the system $y$ initially diverge, both trajectories
evolve almost identically already following the first model update at $t=0.5$.
Note that already after the first learning step, the quality of the model is apparently 
good enough such that the funnel controller remains henceforth inactive and does not 
have to compensate for model errors.
The control signal primarily consists of the control~$\uFMPC$ generated by the model-based controller 
component, see \Cref{fig:sim:massoncar:learning:input}.

\begin{figure}[H]
    \begin{subfigure}[b]{0.49\textwidth}
        \centering
        \includegraphics[width=\linewidth]{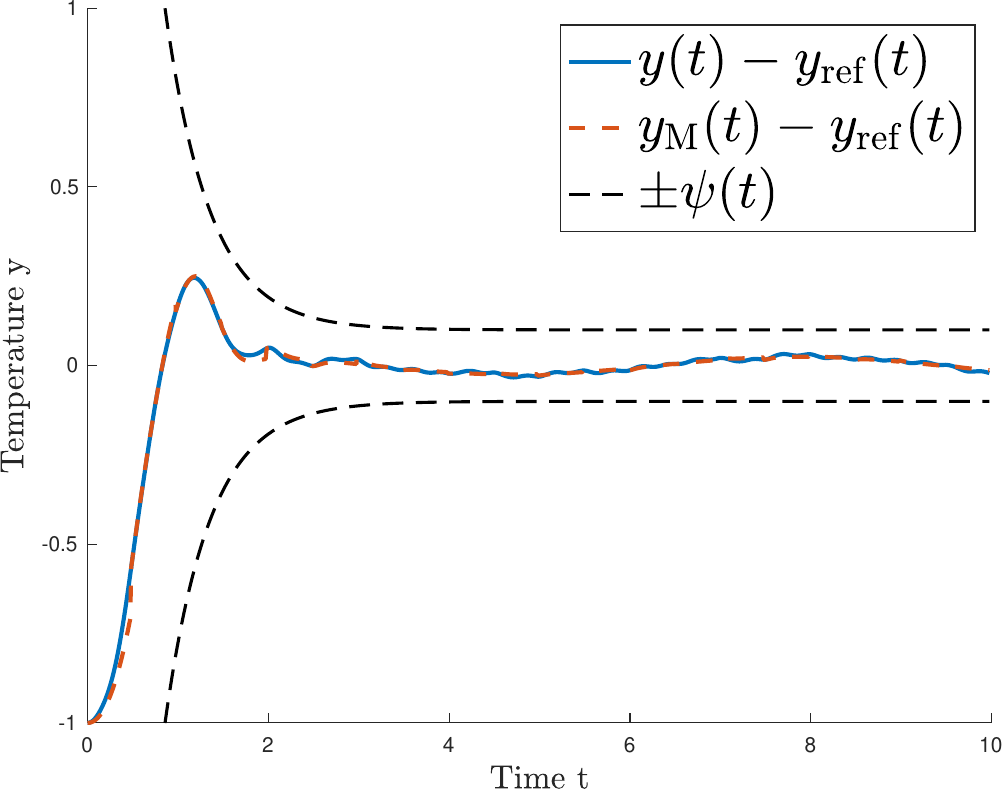}
        \caption{Tracking error $e=y-y_{\rf}$ within boundary~$\Funnel$.}
        \label{fig:sim:massoncar:learning_dist:output}
    \end{subfigure}
      \begin{subfigure}[b]{0.49\textwidth}
        \centering
        \includegraphics[width=\linewidth]{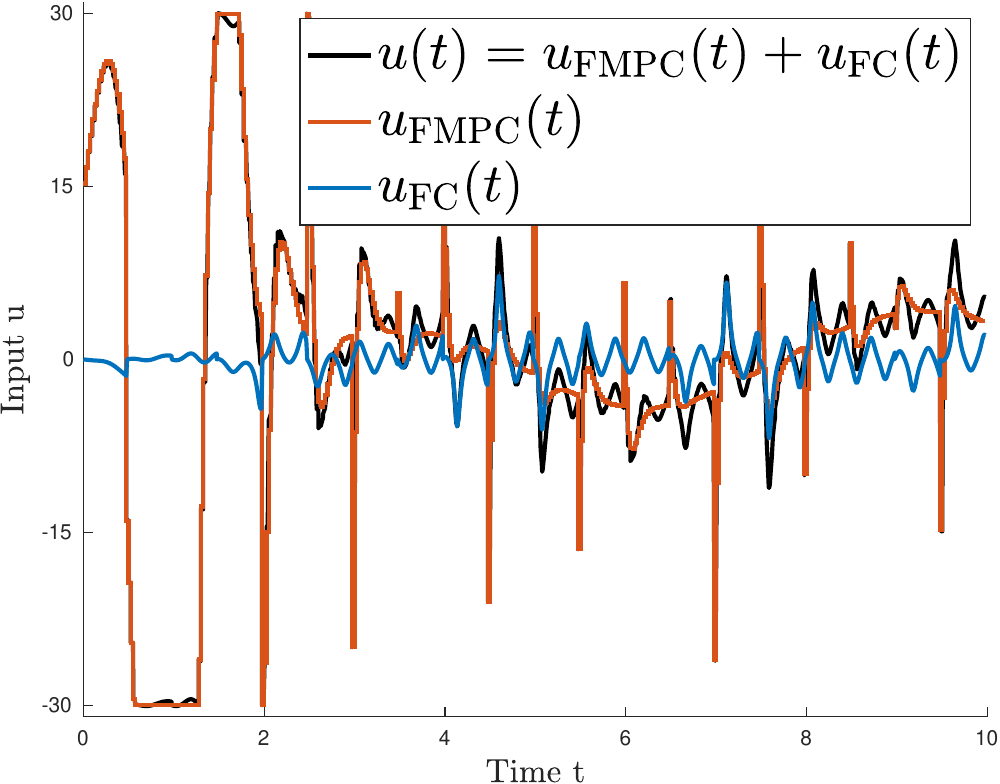}
        \caption{Control inputs.}
        \label{fig:sim:massoncar:learning_dist:input}
    \end{subfigure}
    \caption{Simulation of disturbed system~\eqref{eq:MassOnCarInputOutputDeg2dist} under the
    control generated by the learning-based robust funnel MPC~\Cref{Algo:LearningRFMPC}
    } 
    \label{fig:sim:massoncar:learning_dist}
\end{figure}

In a second simulation, we add an artificial
additive disturbance $d$ to the differential equation, i.e. the system takes the form
\begin{equation}\label{eq:MassOnCarInputOutputDeg2dist}
    \begin{aligned}
        \ddot y(t)    &= R_1y(t)+ R_2\dot{y}(t) + S\eta(t) +\Gamma u(t) + d(t)\\
        \dot \eta(t)  &= Q \eta(t) + P y(t). 
    \end{aligned}
\end{equation}
The disturbance is unknown to the controller and, for the simulation, we choose
the periodic disturbance $d(t)= \cos(20\cdot t)$ and 
leave the controller as it is. The results are depicted in \Cref{fig:sim:massoncar:learning_dist}.
The controller evidently still achieves the control objective.
The two tracking errors $\yM-y_{\rf}$ and~$y-y_{\rf}$ evolve within the
funnel boundaries given by~$\Funnel$, see \Cref{fig:sim:massoncar:learning_dist:output}.
The system output $y$ closely tracks the model $\yM$, which in turn tracks the
given reference $y_{\rf}$ within the prescribed funnel boundaries, despite the
added disturbance. 
Contrary to the prior case, the funnel controller remains active during the
whole operation of the controller. It has to compensate the high-frequency
additive disturbance. However, its contribution remains relatively modest.
The predominant portion of the control signal consists of the  MPC component's control action
suggesting that the learning component still successfully identifies the underlying system dynamics.

\chapter{Sampled-data robust funnel MPC}\label{Chapter:DiscretFMPC}
When applying control strategies to real-world systems, both model predictive and
adaptive control algorithms are nowadays commonly implemented on digital devices.
Unlike the idealised, continuously measured signals assumed in classical control theory, 
practical digital controllers only measure system outputs at discrete sampling
intervals. Consequently, the controller observes the plant at discrete time
points, computes a new input, and then holds that input constant until the next
sample -- introducing two fundamental challenges. 
First, dynamics or disturbances occurring between samples may go undetected, and
high-frequency components can alias as lower-frequency behaviour if the sampling
rate violates the Nyquist–Shannon criterion~\cite{Shannon1949}.
Second, because most digital hardware can usually only generate piecewise‐constant inputs,
the controller cannot apply an arbitrarily varying {(dis-)continuous} actuation
signal, potentially degrading performance relative to a continuous design.
As a result, the controller must be implemented as a
sampled-data controller, specifically designed to operate under these
discrete-time conditions.
In its simplest form, a sampled-data controller  samples the system output at
regular intervals and uses this information to compute a control action. This
control action is then held constant over the entire sampling period, only
updating at the next sampling time.
Although conceptually straightforward, this arrangement requires careful
attention to preserve stability and performance. 
Potential challenges include:
\begin{itemize}
    \item \emph{Stability criteria shift}: Stability of (linear) discrete systems
    require poles inside the unit circle (vs. left half-plane in continuous-time). 
    Discretisation can alter pole locations, destabilising an otherwise
    stable design, see~\Cref{Ex:UnstableDiscr}.
    \item \emph{Model discretisation errors}: 
    Converting a continuous system to a discrete‐time model -- for example via
    Zero-order-Hold (ZoH) approximations -- introduces approximation errors that can
    degrade accuracy~\cite{NESIC1999259,yuz2005sampled}.
    \item \emph{Performance loss}: Sampled-data controllers can reduce performance
    of the closed-loop systems \cite{Leung1991} and exhibit slower responsiveness,
    increased overshoot~\cite{Mita1980}, and steady-state errors \cite{Chen2008}.
    \item \emph{Inter-sample constraint violation}: When safety or performance
    constraints must hold continuously, a controller updated only at discrete
    time instants can inadvertently violate them due to insufficiently fast
    sampling~\cite{Yang2020,Breeden2022}.
\end{itemize}
These issues have motivated a rich body of research in digital control, see \cite{astrom_wittenmark_1997,Laila2006}.
To mitigate discretisation effects and balance trade-offs between sampling
frequency, computational load, and performance in digital implementations,
several mitigation techniques have been developed:
\begin{itemize}
    \item \emph{Sampled-data redesign}: Explicitly account for discrete-time dynamics during controller synthesis,
    rather than simply discretising a continuous design \cite{Grune2008,Grune2008b}.
    \item \emph{Multi-rate sampling}: Use varying sampling frequencies for subsystems with different time scales \cite{Monaco2001,Giovanni2015}.
    \item \emph{Event-triggered and self-triggered Control}: Update control actions
    only when certain conditions are met (e.g. when errors exceed thresholds)
    rather than at fixed intervals, reducing computational load \cite{Heemels2021}. 
\end{itemize}
By accounting for digital implementation from the outset, these approaches help
bridge the gap between continuous‐time theory and real‐world sampled‐data
systems.

In this chapter, we show that it is possible to modify the robust funnel
MPC~\Cref{Algo:RobustFMPC} from~\Cref{Chapter:RobustFunnelMPC} such that it
achieves the output tracking problem with prescribed performance as outlined
in~\Cref{Sec:ControlObjective} with sampled-data control.
In contrast to the robust funnel MPC from~\Cref{Chapter:RobustFunnelMPC}, the
space of admissible controls is restricted to step functions, i.e. the control
signal can only change finitely often between two sampling instants. 
Thus, the control signal applied to the system has the form
\[
    u(t) \equiv u_i \qquad \forall\, t \in [t_i, t_{i+1}),
    \ i \in \N_{0},
\]
where the data to compute the control signal $u_i$ is collected at sample times $(t_i)_{i\in\N_0}$.
To introduce the control scheme properly, we formally define step functions in the following definition.
\begin{definition}[Step function] \label{Def:PartitionAndStepFunction}
    Let $I\subset\R$ be an interval of the form $I=[a,b]$ with $b>a$ or ${I=[a,\infty)}$.
    We call a strictly increasing sequence $\cP=(t_i)_{i\in\N_0}$ with $\lim_{i\to\infty}t_i=\infty$ 
    and ${t_0 = a}$ a~\textit{partition} of $I$.
    The norm of~$\cP$ is defined as $\Abs{\cP}\coloneqq  \sup\setdef{t_{i+1}-t_{i}}{i\in\N_0}$.
    A function~$f:I\to\R^m$ is called~\textit{step function with partition}~$\cP$ if $f$ is constant on every
    interval~$[t_i,t_{i+1})\cap I$ for all~$i\in\N_0$. 
    We denote the space of all step functions on~$I$ with partition~$\cP$ by~$\cT_{\cP}(I,\R^m)$.
\end{definition}
Note that in the case of finite intervals $I=[a,b]$ with $b>a$,
\Cref{Def:PartitionAndStepFunction} can also be formulated 
using finite sequences  $\cP=(t_i)_{i=0}^N$ with $N\in\N$ and $t_N=b$.
However, using infinite sequences every partition~$\cP$ of $[a,\infty)$
is also a partition of $[a,b]$ for all $b>a$.
Using this observation simplifies formulating our results.
Further, note that~\Cref{Def:PartitionAndStepFunction} allows for the usage of a
non-uniform step length, i.e. for $\SampleTime_i \coloneqq |t_{i-1} - t_i|$ we allow $\SampleTime_i
\neq \SampleTime_j$ for $i \neq j$, where $i,j \in \N$. However, in practice, a
uniform step length will be used often.

The robust funnel MPC~\Cref{Algo:RobustFMPC}
from~\Cref{Chapter:RobustFunnelMPC} consists of two components, the model-free
funnel controller~\eqref{eq:DefFC} and the model based funnel
MPC~\Cref{Algo:FunnelMPC}, see also~\Cref{Fig:Robust_FMPC}.
In the  following~\Cref{Sec:FCZoH,Sec:FMCPZoH},
we restrict ourselves to showing that both components individually can be designed 
to work with the restricted space of step functions as control signals.
However, we refrain from  integrating both controllers in one single 
control scheme like done for the robust funnel MPC~\Cref{Algo:RobustFMPC}
and proven in~\Cref{Thm:RobustFMPC}.
The arguments and considerations for such an integration are the same as
in~\Cref{Chapter:RobustFunnelMPC} and do not provide any new insights into the underlying issue.
The restriction to step functions merely adds another level of technicalities.

For the controller design in this chapter, we restrict both the class of
potential systems $\cN^{m,r}_{t_0}$ and associated models $\cM^{m,r}_{t_0}$.
For both the system and the model, we consider non-linear multi-input
multi-output differential equations of order $r\in\N$ of the form 
\begin{equation}\label{eq:SysMod_Discrete}
    \begin{aligned}
    y^{(r)}(t) &= f \big(\oT(\OpChi(y))(t) \big)  + g \big(\oT(\OpChi(y))(t) \big) u(t), \\
    y|_{[0,t_0]} &= y^0  \in \cC^{r-1}([0,t_0],\R^m),
    \end{aligned}
\end{equation}
where $f\in\Lip_{\loc}(\R^q,\R^m)$, $g\in\Lip_{\loc}(\R^q,\R^{m\times m})$, and $\oT\in\cT_{t_0}^{rm,q}$.
In addition, we assume that the matrix valued function $g$ is strictly positive definite, 
that is 
\[
    \forall {x} \in \R^{q} \ \forall z \in \R^m \setminus \{0\}: \quad  \al z, g({x}) z \ar > 0.
\]
Note that by replacing $u$ in~\eqref{eq:SysMod_Discrete} by $-u$
all results presented in this chapter remain valid if $g$ is strictly negative definite.
Note further that, while some authors only use the term \emph{strictly positive definite}
for symmetric matrices, we do not assume $g(x)$ to be symmetric.

We use the notation $(f,g,\oT)\in\ModSysClassDiscr$ to refer to a system,
respectively a model, of the form~\eqref{eq:SysMod_Discrete} with the
aforementioned properties. 
When it is necessary to distinguish between the system and the model,
we will use an index $\mathrm{M}$, i.e. $\fM$, $\gM$, $\oTM$, to refer to the model's functions as done in~\Cref{Chapter:FunnelMPC}.
However, we want to emphasise that the system and the model are not assumed to be identical even though 
we use the same class $\ModSysClassDiscr$ of functions for the system and the model.
In order to avoid having to differentiate between the model and system class in this chapter, 
we assume for the sake of simplicity that $\oT\in\cT_{t_0}^{rm,q}$
However, all presented results hold true if the operator~$\oT$ of the system 
merely fulfils the \emph{causality}~\ref{Item:OperatorPropCasuality},
\emph{local Lipschitz}~\ref{Item:OperatorPropLipschitz}, and the
\emph{bounded-input bounded-output (BIBO)}~\ref{Item:OperatorPropBIBO} property
as defined in~\Cref{Def:OperatorClass}. It is not required to fulfil the
\emph{limited memory} property \ref{Item:OperatorPropLimitMemory}.
\begin{remark}
    As previously pointed out in~\Cref{Rem:SystemClass}~\ref{Item:OperatorContainsDisturbences}, 
    an unknown disturbance $d\in L^\infty([t_0,\infty),\R^p)$ in the system~\eqref{eq:SysMod_Discrete}
    can be modelled in terms of the operator $\oT\in\cT_{t_0}^{rm,q}$.
     Systems of the form
    \[
        y^{(r)}(t) = f \big(d(t),\oT(\OpChi(y))(t) \big)  + g \big(d(t),\oT(\OpChi(y))(t) \big) u(t) 
    \]
    are therefore implicitly contained in the system class $\ModSysClassDiscr$.
\end{remark}

\section{Funnel control with zero-order-hold}\label{Sec:FCZoH}
Funnel control is an adaptive high-gain control methodology guaranteeing 
satisfaction of a priori fixed, possibly time-varying output constraints
while only imposing structural assumptions but not requiring knowledge about the system dynamics,
see e.g. \cite{BergIlch21} and the survey paper \cite{BergIlch23}. 
However, the availability of the system's output as a continuous-time signal and
the ability to continuously adapt the input signal is pivotal for its
functioning, cf. \Cref{Intro:Prop:FC,Prop:FC_dist}.

Although funnel control has been successfully implemented in a sampled-data
system with Zero-order-Hold (ZoH) for a sufficiently small sampling
time~in~\cite{berger2019learningbased}, we are not aware of any results prior to~\cite{LanzaDenn24} 
rigorously showing that the output signal stays within the prescribed boundaries for 
ZoH funnel control. 
In this section, we present the in~\cite{LanzaDenn24} proposed sampled-data
feedback controller with ZoH. We show that the controller ensures
output tracking of a given reference signal within prescribed, possibly
time-varying performance bounds -- at every time instant meaning that also the
intersampling behaviour is fully taken into account. To balance the need for a
sufficiently large feedback gain for output tracking and avoidance of
overshooting (which could violate error bounds within one sampling period), we
use results from the previous chapters to infer uniform bounds on sampling
rates and control inputs. This allows us to ensure that  the imposed
output constraints are satisfied along the closed loop leveraging coarse bounds
on the system dynamics. To the best of our knowledge, in funnel control uniform
bounds on the input signal are only known if the region of feasible initial
values is further restricted \textit{and} the dynamics are
known~\cite{BergIlch21}. While there have been several attempts to deal with the
closely related issue of input saturation~\cite{Berg24,HuTren22,IlchTren04} and
bang-bang controller designs~\cite{LibeTren10,LibeTren13b} exhibiting similarities to our
approach, an analysis of combining a ZoH with funnel control has not been
conducted prior to the work~\cite{LanzaDenn24}.

Before presenting the results from~\cite{LanzaDenn24}, we want to motivate why
applying a controller in a sampled-data fashion to a system poses additional
challenges. When applied to system~\eqref{eq:SysMod_Discrete}, a high-gain
feedback controller, e.g. the funnel controller, achieves the control
objective as laid out in~\Cref{Sec:ControlObjective} if the gain is large
enough. When applied in a sample-and-hold form, however, such approaches can
fail if the gain or the sampling time is too large, respectively. To see this
consider the following example.
\begin{example}\label{Ex:UnstableDiscr}
    Consider the scalar linear system
    \[
        \dot{x}(t)=ax(t) +u (t),
    \]
    with $a\in\R$. As is well known, 
    every linear feedback $u(t)=-kx(t)$ with $k>|a|$ stabilises the system. If
    $u$ is applied in a sample-and-hold form with sampling rate $\SampleTime>0$, then
    the solution at the time instants $i\SampleTime$ with $i\in\N$ has the form
    \[
        x((i+1)\SampleTime) = \me^{a\SampleTime}x(i\SampleTime) -\tfrac{1}{a}\big( \me^{a\SampleTime}-1\big)kx(i\SampleTime)=
        \big(e^{a\SampleTime} -\tfrac{1}{a}\big( \me^{a\SampleTime}-1\big)k\big)x(i\SampleTime).
    \]
    Thus, for $k>\left|a\tfrac{\me^{a\SampleTime}+1}{\me^{a\SampleTime}-1}\right|$,
    we have $|x((i+1)\SampleTime)|>|x(i\SampleTime)|$.
    Therefore, the system is unstable, even if the
    initial uncontrolled system is stable, i.e. $a<0$.
\end{example}
To design a zero-order-hold control strategy able to achieve the control objective using 
data only collected at discrete time instants given a partition $\cP=(t_i)_{i\in\N_0}$ 
of the interval $[t_0,\infty)$ we utilise the auxiliary error variables 
$\eS_{k}$ for $k=1,\ldots, r$  as in~\eqref{eq:ek_FC}.
As in~\Cref{Chapter:RobustFunnelMPC}, they are recursively given  
for $\phi>0$, a bijection $\FCBijec\in\cC^{1}([0,1),[1,\infty))$, and ${z=(z_1,\ldots,z_r)\in\R^{rm}}$ with $z_k\in\R^m$
by
\begin{equation}\tag{\ref{eq:ek_FC} revisited}
        \eS_1(\phi,z)\coloneqq \phi z_1,\quad
        \eS_{k+1}(\phi,z)\coloneqq \phi z_{k+1}+\FCBijec\rbl\Norm{\eS_k(\phi,z)}^2\rbr\eS_{k}(\phi,z),
\end{equation}
for $k=1,\ldots, r-1$.
For details we refer to~\Cref{Chapter:RobustFunnelMPC}.
Given a funnel function $\phi\in\cG$ and a reference trajectory ${y_{\rf}\in W^{r,\infty}(\Rp,\R^m)}$,
we use in the following the short notation $\eS_r(t)\coloneqq \eS_r(\phi(t),\OpChi(y-y_{\rf})(t))$, where
$y$ is the output of the system~\eqref{eq:SysMod_Discrete}.
We propose the following controller structure for~${i \in \N_{0}}$
   \begin{equation}\label{eq:controller_recursive}
   \forall  t \in  [t_i, t_i + \SampleTime)  : \, \uZoH(t) =
    \begin{cases}
        0,                                                          &  \|  \eS_r(t_i)\| < \FCDiscreteThresh, \\
         - \FCDiscreteGain \tfrac{\eS_r(t_i)}{\|\eS_r(t_i)\|^2}  ,  &  \|  \eS_r(t_i)\| \ge \FCDiscreteThresh,
    \end{cases}
    \end{equation}
where $\FCDiscreteThresh\in (0,1)$ is an \emph{activation threshold},
and~$\FCDiscreteGain > 0$ is the \emph{input gain}.
In~\Cref{Thm:DiscreteFC}, we derive lower bounds on the input gain $\FCDiscreteGain$ and 
upper bounds on the maximal sampling time, i.e. $\SampleTime=\Abs{\cP}$, which 
ensure that the control objective is achieved when applying the
controller~\eqref{eq:controller_recursive} to the
system~\eqref{eq:SysMod_Discrete}.
We do this by showing $\eS_r(t)\in\cB_1$ for all $t\geq t_0$.
Thus, the control signal $\uZoH$ is then uniformly bounded since
\[
    \forall t \geq t_0:\quad \|\uZoH(t) \| \leq \frac{\FCDiscreteGain}{\FCDiscreteThresh}.
\]
The controller design can be considered to be similar to funnel control, see
\cite{BergIlch21,BergLe18a,IlchRyan02b}, in terms of its ability to achieve
output reference tracking within predefined error boundaries, as well as
concerning the used intermediate error variables~\eqref{eq:ek_FC}. On the other
hand, contrary to the standard funnel controller, the feedback
law~\eqref{eq:controller_recursive} is a normalised linear sample-and-hold
output feedback with uniformly bounded sampling rate.
A further essential difference to continuous funnel control is that in the
present approach the control objective is achieved by using estimates about the
system dynamics, while in continuous-time funnel control no such information is used to the
price that the maximal control effort cannot be estimated a priori.

In order to formulate and prove the main result of~\cite{LanzaDenn24} about 
feasibility of the proposed ZoH controller~\eqref{eq:controller_recursive},
we recall some results from the previous~\Cref{Chapter:FunnelMPC,Chapter:RobustFunnelMPC}. 
To ensure that the controller achieves the control objective, namely that the system output~$y$ 
tracks a given reference signal $y_{\rf}$ with prescribed performance in terms of a function $\phi\in\cG$,
we show that the norm of the axillary error variables $\eS_{k}$ for $k=1,\ldots, r$ as in~\eqref{eq:ek_FC} evaluated along
$\OpChi(y-y_{\rf})$ is always below one, i.e.
\begin{equation}\label{eq:DiscrAuxillaryErrorsLowerOne}
    \fa k=1,\ldots,r\ \fa t\geq t_0:\quad \Norm{\eS_{k}(\phi(t),\OpChi(y-y_{\rf})(t))}<1.
\end{equation}
Assuming that \eqref{eq:DiscrAuxillaryErrorsLowerOne} is fulfilled at the initial time $t=t_0$,  
\Cref{Lem:FCExistenceEpsMu} states that all error signals $\eS_{k}(\phi(t),\OpChi(y-y_{\rf})(t))$ for $k=1,\ldots, r$
satisfy \eqref{eq:DiscrAuxillaryErrorsLowerOne} for all $t\geq t_0$ given that 
the norm of the last auxiliary error $\eS_{r}(\phi(t),\OpChi(y-y_{\rf})(t))$ remains below one for all $t\geq t_0$.
In the proof of~\Cref{Prop:FC_dist}, we used this result to show that the funnel controller~\eqref{eq:DefFC}
achieves the control objective.
In a similar fashion, we derive bounds on the input gain~$\FCDiscreteGain$ and upper bounds on the sampling rate, i.e. $\SampleTime=\Abs{\cP}$, 
that ensure the feasibility of the proposed ZoH controller~\eqref{eq:controller_recursive}
by guaranteeing that the norm of the last auxiliary error $\eS_{r}(\phi(t),\OpChi(y-y_{\rf})(t))$ remains bounded by one.
In addition to the mentioned statement, \Cref{Lem:FCExistenceEpsMu} 
states that $\Norm{\eS_{k}(\phi,\OpChi(\zeta))}$ for $k=1,\ldots, r-1$ remain bounded away from one
by some $\eps\in(0,1)$ for all signals $\zeta\in\FCTrajectories_{\hat{t}}$ and all $\hat{t}\geq t_0$, 
where $\FCTrajectories_{\hat{t}}$ is the set of all functions coinciding with $y_0$ on
the interval $[0,t_0]$ and fulfilling~\eqref{eq:DiscrAuxillaryErrorsLowerOne}
where $y-y_{\rf}$ is replaced by $\zeta$, see~\eqref{eq:DefFCTrajectories} 
This yields the existence of a compact set in which all functions $\zeta\in\FCTrajectories_{\hat{t}}$ evolve until~$\hat{t}$.
The existence of such a compact set allows us to adapt~\Cref{Lemma:DynamicBounded} to the current setting
stating that the  system~\eqref{eq:SysMod_Discrete} with
$(f,g,\oT)\in\ModSysClassDiscr$ is uniformly bounded for every
$\zeta\in\FCTrajectories_{\hat{t}}$ and all $\hat{t}\geq t_0$.
\begin{lemma}\label{Lemma:DynamicBoundedDiscreteSys}
    Consider the system~$\eqref{eq:SysMod_Discrete}$ with 
    $(f,g,\oT)\in\ModSysClassDiscr$ and 
    reference trajectory $y_{\rf}\in W^{r,\infty}(\Rp,\R^{m})$.  
    Let $\phi\in\cG$.
    Then, there exist constants~$\fmax$, $\gmax\geq0$ such that for all $\hat{t}\in(t_0,\infty]$
    and $\zeta\in\FCTrajectories_{\hat{t}}$:
    \[
        \fmax\geq\SNorm{f(\oTM(\OpChi(\zeta))|_{[0,\hat{t})})},\qquad
        \gmax\geq\SNorm{g(\oTM(\OpChi(\zeta))|_{[0,\hat{t})})}.
    \]
    Moreover, there exists $\gmin>0$ such that for all $z\in\R^m\backslash\cbl0\cbr$ and all $\hat{t}\in (t_0,\infty]$:
    \[
        \gmin\leq\frac{\al z, g(\oT(\OpChi(\zeta))|_{[0,\hat{t})}(t))z\ar}{\Norm{z}^2}.
    \]
\end{lemma}
\begin{proof}
    To prove the assertion, we adapt the proof of~\Cref{Lemma:DynamicBounded} to the current setting. 
    According to \Cref{Lem:FCExistenceEpsMu}, there exist 
    constants $\eps_k>0$ such that all functions $\zeta\in
    \FCTrajectories_{\infty}$  fulfil $\|\eS_k(\phi(t),\OpChi(\zeta)(t))\|  \
    {\leq}\ \eps_k < 1$ for all $t\in[t_0,\infty)$ and all $k=1,\ldots, r-1$.
    Hence, by boundedness of $\phi$ and $y_{\rf}^{(i)}$ for all~$i=1,\ldots, r$,
    there exists a compact set $K\subset\R^{rm}$ with 
    \[
        \fa \zeta\in\FCTrajectories_{\infty}\fa t\geq 0:\quad \OpChi(\zeta)(t)\in K.
    \]
    Invoking the BIBO property of the operator~$\oT$, there exists a compact set $K_q\subset\R^q$ with  
    $\oT(z)(\Rp)\subset K_q$ for all $z\in\cC(\Rp,\R^{rm})$ with $z(\Rp)\subset K$.
    For arbitrary~$\hat{t}\in(t_0,\infty)$ and $\zeta\in\FCTrajectories_{\hat{t}}$,
    we have $\OpChi(\zeta)(t)\in K$ for all $t\in[0,\hat{t}]$.
    For every element $\zeta\in\FCTrajectories_{\hat{t}}$ the restriction $\OpChi(\zeta)|_{[0,\hat{t})}$ 
    can be extended to a function~$\tilde{\zeta}\in\cR(\Rp,\R^{rm})$ 
    with $\tilde{\zeta}(t)\in K$ for all $t\in\Rp$.
    We have $\oT(\tilde{\zeta})(t)\in K_q$ for all $t\in\Rp$ because of the BIBO property of the operator~$\oT$. 
    This implies $\oT(\OpChi(\zeta))|_{[0,\hat{t})}(t)\in K_q$ for all $t\in [0,\hat{t})$ and $\zeta\in\FCTrajectories_{\hat{t}}$ since $\oT$ is causal.
    Since~$f(\cdot)$ and $g(\cdot)$ are continuous, the constants $\fMmax=\max_{z\in K_q}\Norm{\fM(z)}$
    and $\gMmax=\max_{z\in K_q}\Norm{\gM(z)}$ are well-defined.
    For all $\hat{t}\in(0,\infty]$ and $\zeta\in\FCTrajectories_{\hat{t}}$ we have 
    \[
        \fa t\in [0,\hat{t}):\ \oT(\OpChi(\zeta))(t)\in K_q.
    \]
    Furthermore, since $g(x)$ is positive definite for every $x\in K_q$, there exists $\gmin>0$ such that 
    $ \gmin \leq \frac{\al z, g(\oT(\OpChi(\zeta)))|_{[0,\hat{t})}(t))z\ar }{\Norm{z}^2}$ for all $z\in\R^m\backslash\cbl0\cbr$. 
    This completes the proof. 
\end{proof}

A consequence of~\Cref{Lemma:DynamicBoundedDiscreteSys} is that the dynamics of
system~\eqref{eq:SysMod_Discrete} are bounded if a control is applied that
ensures that all error signals $\eS_{k}(\phi(t),\OpChi(y-y_{\rf})(t))$ for
$k=1,\ldots, r$ satisfy~\eqref{eq:DiscrAuxillaryErrorsLowerOne}.
In the following~\Cref{Thm:DiscreteFC}, 
we use these bounds to derive an input gain~$\FCDiscreteGain > 0$ 
large enough to counteract the system dynamics.
When applying the ZoH controller~\eqref{eq:controller_recursive} to the system
the large enough gain guarantee that the norm of the auxiliary error signal 
$\eS_{r}(\phi(t),\OpChi(y-y_{\rf})(t))$ decreases at the sampling instants~$t_i$,
if the error signal is greater than or equal to the activation threshold $\FCDiscreteThresh\in (0,1)$,
see Step~2.b in the proof of \Cref{Thm:DiscreteFC}.
Based on bound on the system dynamics and the maximal control value applied to the system, 
we compute a uniform bound on the sampling time $\SampleTime=\Abs{\cP}$ required 
to avoid overshooting of the error signal between two sampling instants.

\begin{theorem}\label{Thm:DiscreteFC}
Given a reference {$y_{\rf} \in W^{r,\infty}(\Rp,\R^m)$} and a funnel function ${\phi \in \cG}$,
consider the system \eqref{eq:SysMod_Discrete} with $(f,g,\oT) \in \ModSysClassDiscr$.
Assume that the initial trajectory ${y^0\in \cC^{r-1}([0,t_0],\R^m)}$ satisfies 
${\OpChi(y^0 - y_{\rf})(t_0)\in\cEFC{1}(\phi(t_0))}$, i.e.
the error variables in~\eqref{eq:ek_FC}  satisfy
${\| \eS_k(\phi(t_0),\OpChi(y^0-y_{\rf})(t_0))\| < 1}$
for all $k=1,\ldots,r$.
With the constants given in~\eqref{eq:ve_mu_gam} and in~\Cref{Lemma:DynamicBoundedDiscreteSys}, set
\[
        \kappa_0 \coloneqq  \SNorm{\frac{\dot \phi}{\phi}} ( 1 + \FCBijec(\eps_{r-1}^2) \eps_{r-1} )
                    + \SNorm{\phi} ( \fmax + \|y_{\rf}^{(r)}\|_\infty )+ \bar \eta_{r-1},
 \]
and choose the input gain
    \begin{equation*}
        \FCDiscreteGain > \frac{2 \kappa_0}{ \gmin \inf_{s \ge 0} \phi(s) }.
    \end{equation*}
Further,  for an activation threshold ${\FCDiscreteThresh \in (0,1)}$, define the constant
$
   \kappa_1 \coloneqq  \kappa_0 +  \SNorm{\phi} \tfrac{\FCDiscreteGain}{\FCDiscreteThresh} \gmax
$
and  let $\cP$ be a partition of the interval $[t_0,\infty)$
for which the maximal sampling time $\SampleTime\coloneqq \Abs{\cP}$ fulfils  
    \begin{equation} \label{eq:tau_r}
        0 <      \SampleTime \le \min \left\{ \frac{ \kappa_0 }{\kappa_1^2}, \frac{1-\FCDiscreteThresh}{\kappa_0} \right\}.
    \end{equation}
Then, the ZoH controller~\eqref{eq:controller_recursive}
applied to a system~\eqref{eq:SysMod_Discrete} yields 
\[
    \Norm{\eS_k(\phi(t),\OpChi(y-y_{\rf})(t))} < 1
\]
for all $k=1,\ldots,r-1$ and $\|\eS_r(t)\| \le 1$ for all $t \ge t_0$.
This is initial and recursive feasibility of the ZoH control law~\eqref{eq:controller_recursive}.
In particular, the tracking error $e\coloneqq y-y_{\rf}$ satisfies $\|e(t)\| < 1/\phi(t)$ for all $t \geq t_0$.
\end{theorem}

\begin{proof}
The proof consists of two main steps. 
In the first step, we establish the existence of a solution of the initial value problem~\eqref{eq:SysMod_Discrete},~\eqref{eq:controller_recursive}.
In the second step, we show feasibility of the proposed control law, 
i.e. all error variables $\eS_k(\phi(t),\OpChi(y-y_{\rf})(t))$, $k=1,\ldots, r$ are bounded by one. 
Thus, the tracking error evolves within the funnel boundaries given by~$\phi$.
In the following, we use the shorthand notation  $\eS_k(t)\coloneqq \eS_k(\phi(t),\OpChi(y-y_{\rf})(t))$.

\noindent
\emph{Step 1}:
The application of the control signal~\eqref{eq:controller_recursive} to system~\eqref{eq:SysMod_Discrete} leads to an initial value problem. 
If this problem is considered on the interval~$[t_0,\SampleTime]$, then there exists a unique maximal solution on $[t_0,\omega)$ with $\omega\in(t_0,\SampleTime]$.
If all error variables $\eS_k(t)$ evolve within the set $\cB_1$ for all $t\in [t_0,\omega)$,
then $\| \OpChi(y)(\cdot) \|$ is bounded on the interval $[t_0,\omega)$ and,
as a consequence of the BIBO condition of the operator, $\oT(\cdot)$ is bounded as well. 
Then $\omega =\SampleTime$, cf. \cite[\S~10, Thm.~XX]{Walt98} and there is nothing else to show.
Seeking a contradiction, assume the existence of $t \in [t_0,\omega)$ such that $\| \eS_k(t)\| \ge 1$ for at least one $k = 1,\ldots,r$.
Invoking \Cref{Lem:FCExistenceEpsMu}, it remains only to show that the last error variable $\eS_r$ satisfies $\| \eS_r(t) \| \le 1$ for all $t \in [t_0,\omega)$.
Before doing so, record the following observation.
For $\eta_{r-1}(t) \coloneqq  \FCBijec(\|\eS_{r-1}(t)\|^2) \eS_{r-1}(t)$ and ${z}(\cdot) \coloneqq  \oT(\chi(y))(\cdot)$, we calculate
\begin{equation} \label{eq:J}
\begin{aligned}
     \dot \eS_r(t) - \phi(t) g({z}(t)) u(t) &= \dot \phi(t) e^{(r-1)}(t) + \phi(t) e^{(r)}(t) + \dot \eta_{r-1}(t) - \phi(t) g({z}(t)) u \\
     &= \frac{\dot \phi(t)}{\phi(t)} (\eS_r(t) - \eta_{r-1}(t)) + \dot \eta_{r-1}(t)  + \phi(t) ( f({z}(t)) - y_{\rf}^{(r)}(t)  ) =: J(t).
\end{aligned}
\end{equation}
\noindent
\emph{Step 2}: We show $\|\eS_r(t)\| \le 1$ for all $t \in [t_0,\omega)$.
We separately investigate the two cases $\| \eS_r(t_0)\| < \FCDiscreteThresh$ and $\| \eS_r(t_0)\| \ge \FCDiscreteThresh$. \\
\emph{Step 2.a}: Consider $\| \eS_r(t_0)\| < \FCDiscreteThresh$. In this case, the constant control signal ${u (t)=\uZoH(t) = 0}$ is applied to the system.
Seeking a contradiction, we suppose that there exists a time instant ${t^* \coloneqq  \inf \setdef{ t \in (t_0,\omega)}{ \| \eS_r(t) \|> 1}}$. For the function $J(\cdot)$ introduced in~\eqref{eq:J}, we observe ${\| J|_{[t_0,t^*)} \|_\infty \le \kappa_0}$ according to \Cref{Lem:FCExistenceEpsMu,Lemma:DynamicBoundedDiscreteSys}.
Then, we calculate 
\begin{align*}
  1 &=  \| \eS_r(t^*) \| \le \|\eS_r(t_0)\|  + \int_{t_0}^{t^*}  \| \dot \eS_r(s) \|  \d s \\
    &= \|\eS_r(t_0)\| +  \int_{t_0}^{t^*}   \| J(s) \| \d s \\
    &\le  \|\eS_r(t_0)\| +   \int_{t_0}^{t^*} \kappa_0  \d s 
    < \FCDiscreteThresh + \kappa_0 \omega < 1,
\end{align*}
where $t^* < \omega \le \SampleTime < (1-\FCDiscreteThresh)/\kappa_0$ was used.
This contradicts the definition of $t^*$. \\
\emph{Step 2.b}:
Consider $\| \eS_r(t_0)\| \ge \FCDiscreteThresh$. In this case, ${u(t) = \uZoH(t)= - \FCDiscreteGain \eS_r(t_0)/\|\eS_r(t_0)\|^2}$
is applied to the system.
We show again $\|\eS_r(t)\| \le 1$ for all $t \in [t_0,\omega)$.
To this end, seeking a contradiction, we suppose the existence of $t^* = \inf \setdef{(t_0,\omega)}{ \| \eS_r(t) \| > 1 }$.
For the function~$J(\cdot)$, we observe $\| J|_{[t_0,t^*)} \|_\infty \le \kappa_0$ according to
\Cref{Lem:FCExistenceEpsMu,Lemma:DynamicBoundedDiscreteSys}. 
Moreover, $\| \dot \eS_r|_{[t_0,t^*]} \| \le \kappa_1$ due to equation \eqref{eq:J} 
and the bound $\SNorm{\uZoH}\leq \tfrac{\FCDiscreteGain}{\FCDiscreteThresh}$.
Invoking the initial conditions and continuity of the involved functions, {and~\eqref{eq:J}}, we calculate for~$t \in [t_0,t^*]$:
\begin{align*}
    \dd{t} \tfrac{1}{2} \| \eS_r(t)\|^2 &= \al \eS_r(t), \dot \eS_r(t) \ar 
     = \al \eS_r(t_0) +   \int_{t_0}^t \dot \eS_r(s) \d s, \dot \eS_r(t) \ar \\
    & \le \|\eS_r(t_0)\| \| J(t) \| + \omega \| \dot \eS_r|_{[t_0,t^*]}\|^2_\infty + \phi(t) \al \eS_r(t_0), g({z}(t)) \uZoH(t) \ar \\
    & = \|\eS_r(t_0)\| \| J(t) \| + \omega \| \dot \eS_r|_{[t_0,t^*]}\|^2_\infty -\phi(t) \FCDiscreteGain \tfrac{\al \eS_r(t_0), g({z}(t)) \eS_r(t_0)\ar}{\|\eS_r(t_0)\|^2} \\
    & \le \|\eS_r(t_0)\| \kappa_0 + \omega \| \dot \eS_r|_{[t_0,t^*]}\|^2_\infty - \inf_{s \ge t_0} \phi(s) g_{\min} \FCDiscreteGain \\
    & \le \kappa_0 + \omega \kappa_1^2 - \inf_{s \ge t_0} \phi(s) g_{\min} \FCDiscreteGain  \\
    &\le 2 \kappa_0 - \inf_{s \ge t_0} \phi(s) g_{\min} \FCDiscreteGain  < 0.
\end{align*}
Here, the second line holds true due to $t^* < \omega \le \SampleTime$, 
the penultimate line via the definition of $\SampleTime$,
and the last line by definition of~$\FCDiscreteGain$. 
In particular, this yields $\tfrac{1}{2} \dd{t} \| \eS_r(t)|_{t=t_0}\|^2 < 0$, by which $t^* > t_0$.
Therefore, we find the contradiction $1 = \| \eS_r(t^*)\|^2 < \| \eS_r(t_0)\|^2 \le 1$.
Repeated application of the arguments in \emph{Steps 1} and \emph{2} on the interval $[t_i,t_i + \SampleTime]$, $i \in \N$, 
yields recursive feasibility.
\end{proof}

The maximal sampling time $\SampleTime$ in~\eqref{eq:tau_r} strongly depends on the {evolution
of the} funnel function and on the reference~$y_{\rf}$. This gives the
possibility of dynamically adapting the sampling time, e.g. in the case of
setpoint transition, where the reference is constant $y_{\rf}^0$ in the first
period and constant $y_{\rf}^1 \neq y_{\rf}^0$ in the last period. At the
setpoints the sampling time can be larger than during the transition.

The parameter $\FCDiscreteThresh\in (0,1)$ in~\eqref{eq:controller_recursive} is
an ``activation threshold'' to set the control input to zero for small tracking errors,
akin to the idea of using funnel control with an activation function as
discussed in~\Cref{Sec:RobustFMPCInit}, the $\lambda$-tracker~\cite{IlchRyan94},
or more broadly event- and self-triggered controller designs, see
e.g.~\cite{Heemels2021} and references therein.
The activation threshold~$\FCDiscreteThresh$ is chosen by the designer and
divides the funnel for the tracking error in a safe and a safety critical
region. A large value of $\FCDiscreteThresh$ implies that the controller will be
inactive for a wide range of values of the last error variable, which, in case
of relative degree one, means inactivity for a wide range of the tracking error,
while still guaranteeing transient accuracy.

Applying a zero-input to the system~\eqref{eq:SysMod_Discrete} while the tracking error 
is within the safe region is mainly done for mathematical reasons as it simplifies the proof of~\Cref{Thm:DiscreteFC}. 
In many situations it might be beneficial to apply different bounded control signal instead.
One potential strategy is to simply \emph{hold the input}, i.e. to apply the control value $u(t_{i-1})$ of the last 
sampling period. As pointed out in~\cite{Schenato09},
neither of these two strategies is consistently superior to the other. 
However, more sophisticated strategies may choose the control value according 
to some data informativity framework~\cite{van2020data} and can outperform the
controller~\eqref{eq:controller_recursive}.
In the following~\Cref{Sec:FCAsSaftyFilter}, we give a short outlook on how
such data-driven approaches can be safeguarded by the proposed
controller~\eqref{eq:controller_recursive}.

An explicit bound on the control input can be computed in advance, since $\| u \|_\infty \le \FCDiscreteGain/\FCDiscreteThresh$. 
This bound depends on the system parameters derived in~\Cref{Lemma:DynamicBoundedDiscreteSys}.
However, precise knowledge about the functions $f$, $g$ and the operator~$\oT$ is not necessary.
Mere (conservative) estimates on the bounds $\fmax$, $\gmax$, and $\gmin$ in~\Cref{Lemma:DynamicBoundedDiscreteSys}
are sufficient to guarantee the functioning of the ZoH controller~\eqref{eq:controller_recursive}.

The controller~\eqref{eq:controller_recursive} only requires for its functioning
measurement data of the system's output and its derivatives at discrete time
instants $t_i$.
It therefore overcomes the funnel controller's requirement of the availability
of continuous output signal.
However, the reliance of the controller~\eqref{eq:controller_recursive} on the
derivatives of the system's output can still be problematic in application as
those signals are very sensitive to noise and might require the usage of
numerical differentiation algorithms. For systems of order $r=2$ the control
approach~\eqref{eq:controller_recursive} was adapted
in~\cite{lanza2024derivative} to overcome this issue and to only rely on the
output signal at discrete time instants but not on its derivatives.
However, a generalisation to higher-order systems is still outstanding.

\subsection{Safeguarded data-based control}\label{Sec:FCAsSaftyFilter}

Dividing the funnel for the tracking error in a safe and a safety critical
region opens up the possibility for the
controller~\eqref{eq:controller_recursive} to act as a safety filter for
data-driven approaches and (online) learning techniques, which have gained a lot
of popularity recently. These techniques, despite their superior performance,
often lack rigorous constraint satisfaction, which is especially important  in
safety-critical applications like medical devices and human-robot interaction,
see e.g.~\cite{brunke2022safe}. We also refer to~\cite{amodei2016concrete}
and~\cite{tambon2022certify} for an overview of the challenges employing
learning-based approaches to safety-critical systems; and for challenges and
recent results in the field of continual learning, we refer to the two comprehensive
surveys~\cite{shaheen2022continual,wang2024comprehensive}.

To address the challenge of ensuring constraint satisfaction while leveraging the benefits of
learning-based control, the field of safe learning has gained prominence and
several safety
frameworks have been proposed~\cite{HewingWaber20, garcia2015comprehensive}, employing various
approaches like control barrier functions~\cite{ames2019control}, Hamilton-Jacobi reachability
analysis~\cite{bansal2017hamilton,chen2018hamilton}, Model Predictive Control (MPC)~\cite{Aswa13}, and Lyapunov
stability~\cite{perkins2002lyapunov}. Predictive safety filters, as exemplified
in~\cite{Wabersich21,Wabersich23}, verify control input signals against a model to ensure compliance
with prescribed constraints. 
Similar ideas are also used in the learning-based robust funnel
MPC~\Cref{Algo:LearningRFMPC} from~\Cref{Chapter:LearningRobustFMPC} as the
funnel controller compensates for the model inaccuracies of the model based
controller component. The model-free controller component serves as 
a safety filter for the learning component which updates (or even replaces) the model at runtime 
while being employed in the funnel MPC algorithm.
In~\cite{GottschalkLanza24} the funnel controller from~\cite{BergIlch21} in
combination with an activation function as presented
in~\Cref{Sec:RobustFMPCInit} was used in a comparable manner as a
safety filter for a model-free Reinforcement Learning (RL) control algorithm, 
namely the Proximal Policy Optimisation~(PPO) algorithm from~\cite{schulman2017proximal}.
In a similar manner, the funnel controller was utilised to ensure safety
guarantees for Koopman operator-based MPC scheme in~\cite{BoldLanzWoth2024_Koopman}.

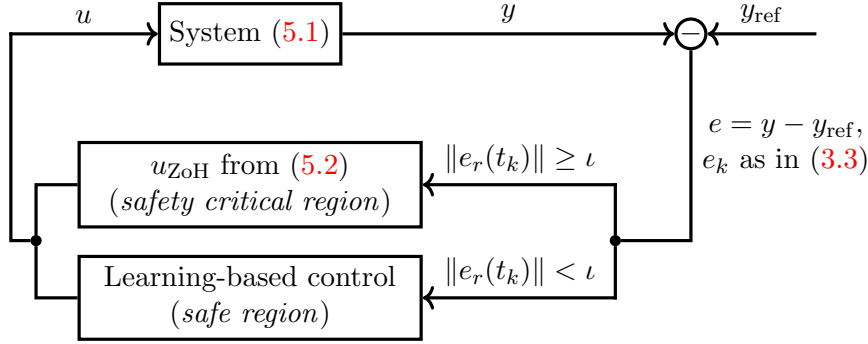
\begin{figure}[ht]
  \begin{center}
    \begin{tikzpicture} [very thick,%
        scale=0.64,%
        node distance = 9ex,
        box/.style={fill=white,rectangle, draw=black},
        blackdot/.style={inner sep = 0, minimum size=3pt,shape=circle,fill,draw=black},%
        blackdotsmall/.style={inner sep = 0, minimum size=0.1pt,shape=circle,fill,draw=black},%
        plus/.style={fill=white,circle,inner sep = 0,very thick,draw},%
        metabox/.style={inner sep = 3ex,rectangle,draw,dotted,fill=gray!20!white},%
        every text node part/.style={align=center}]
\def\hoch{0.8cm};
\def\breit{0cm};
\def\distu{1.3cm};
\def\dista{1.3cm};
\node [block, minimum width = \breit, minimum height = \hoch,] (System) {System~\eqref{eq:SysMod_Discrete}};
\node [block, minimum width = 4.5cm, minimum height = \hoch, below of=System, node distance = 2cm] (ZoH) { $\uZoH$ from~\eqref{eq:controller_recursive} \\ (\emph{safety critical region})};
\node [block, minimum width = 4.5cm, minimum height = \hoch, below of=ZoH, node distance = 1.5cm ] (ZoH+MPC) {Learning-based control \\ (\emph{safe region})};
\node[left of=System, node distance = 3cm, xshift=-2mm] (u) {};
\node(fork1)     [minimum size=0pt, inner sep = 0pt, left of = System, xshift=-10ex] {};
\node(forkZoH)   [minimum size=0pt, inner sep = 0pt, left of = ZoH, xshift=-8ex] {};
\node(forkLearn) [minimum size=0pt, inner sep = 0pt, left of =  ZoH+MPC, xshift=-8ex] {};
\node(joinU)     [blackdot, below of = forkZoH, yshift=0.75cm]{};
\node(ZoHin)   [minimum size=0pt, inner sep = 0pt, right of = ZoH, xshift=20ex] {};
\node(EinFork)   [blackdot, below of = ZoHin, yshift=0.75cm] {};
\node(Learnin) [minimum size=0pt, inner sep = 0pt, right of =  ZoH+MPC, xshift=20ex] {};

\node[circle,,draw=black, fill=white,inner sep=0pt,minimum size=3pt,  right of=System, node distance = 5.8cm] (ref_in) {$-$};
\node[right of=ref_in, node distance = 1.8cm] (ref) {};
\coordinate[below of=ref_in, node distance = 1.9cm] (er);
\draw[->] (fork1.west) -- (System) node[pos=0.5,above] {$u$};
\draw[-] (joinU) |- (ZoH);
\draw[-] (joinU) |- (ZoH+MPC);
\draw[-] (fork1.north) |- (joinU);

\draw[->] (ZoHin.east) -- (ZoH) node[pos=0.5,above] {$\|\eS_r(t_k)\| \ge \FCDiscreteThresh$};
\draw[->] (Learnin.east) -- (ZoH+MPC)node[pos=0.5,above] {$\|\eS_r(t_k)\| < \FCDiscreteThresh$};
\draw[-] (EinFork) -- (ZoHin.north);
\draw[-] (EinFork) -- (Learnin.south);
\draw[-] (ref_in) |- node[pos = 0.25, right]{$e = y - y_{\rf}$,\\ $e_k$ as in~\eqref{eq:ek_FC}} (EinFork) ;

\draw[->] (System) --  node[above]{$y$} (ref_in);
\draw[->] (ref) --  node[above]{$y_{\rf}$} (ref_in);
%\draw[-] (ref_in) -- node[right]{$e = y - y_{\rf}$,\\ $e_k$ as in~\eqref{eq:ek}} (er);
\end{tikzpicture}
    \end{center}
 \caption{Schematic structure of the combined controller~\eqref{eq:combined_controller}.}
 \label{Fig:DiscreteFCController}
\end{figure}

To utilise the controller~\eqref{eq:controller_recursive} as a safety filter,
the idea is to apply a data-driven control algorithm to the system~\eqref{eq:SysMod_Discrete} 
and temporarily interrupt its learning and control process when the activation threshold is
surpassed, resorting to the pure feedback control with ZoH, see~\Cref{Fig:DiscreteFCController}.
The combination of a data-driven control algorithm with the ZoH feedback
control~\eqref{eq:controller_recursive} can be formulated in the following 
switched control strategy.
    \begin{equation} \label{eq:combined_controller}
    \forall  t \in  [t_i, t_{i+1})  :  u(t) =
    \begin{cases}
   \quad  u_{\mathrm{data}}(t),   &  \|  e_r(t_i)\| < \FCDiscreteThresh, \\
         - \FCDiscreteGain \tfrac{e_r(t_i)}{\|e_r(t_i)\|^2}  ,  &\|  e_r(t_i)\| \ge \FCDiscreteThresh.
    \end{cases}
\end{equation}
Since the calculations in the proof of~\Cref{Thm:DiscreteFC} involve worst case
estimates, the application of $u(t) \neq 0$ for $t \in [t_i, t_i + \SampleTime)$, if
$\|e_r(t_i)\| < \FCDiscreteThresh$ requires adaption of the sampling
time~$\SampleTime$. The following~\Cref{Thm:Combined_Controller} formalises this observation.

\begin{theorem}\label{Thm:Combined_Controller}
Given a reference {$y_{\rf} \in W^{r,\infty}(\Rp,\R^m)$} and a function
${\phi \in \cG}$, consider a system \eqref{eq:SysMod_Discrete} with $(f,g,\oT) \in\ModSysClassDiscr$.
Assume the initial trajectory ${y^0\in \cC^{r-1}([0,t_0],\R^m)}$ satisfies 
$\OpChi(y^0 - y_{\rf})(t_0)\in\cEFC{1}(\phi(t_0))$.
Let the constants on the system dynamics be given as in~\Cref{Lemma:DynamicBoundedDiscreteSys}, 
and, for an activation threshold ${\FCDiscreteThresh \in (0,1)}$, $\kappa_0, \kappa_1$ and $\FCDiscreteGain$ be given as in \Cref{Thm:DiscreteFC}.
Further, for $\umax \ge 0$, 
let $\cP$ be a partition of the interval $[t_0,\infty)$
for which the maximal sampling time $\SampleTime\coloneqq \Abs{\cP}$ satisfies  
\[
    0 < \SampleTime \le \min \left\{ \frac{\kappa_0}{\kappa_1^2}, \frac{1-\FCDiscreteThresh}{\kappa_0 + \|\phi\|_\infty g_{\max} \umax} \right\}.
\]
If $\SNorm{u_{\mathrm{data}}} \le \umax$, then the combined
controller~\eqref{eq:combined_controller} applied to a system~\eqref{eq:SysMod_Discrete}
yields 
\[
    \Norm{\eS_k(\phi(t),\OpChi(y-y_{\rf})(t))} < 1
\]
for all $k=1,\ldots,r-1$ and $\|\eS_r(t)\| \le 1$ for all $t \ge t_0$.
This is initial and recursive feasibility of the ZoH control law~\eqref{eq:combined_controller}.
In particular, the tracking error $e\coloneqq y-y_{\rf}$ satisfies $\|e(t)\| < 1/\phi(t)$ for all $t \geq t_0$.
\end{theorem}
\begin{proof}
    By adapting the sampling time~$\SampleTime$ the statement follows with the same proof as for~\Cref{Thm:DiscreteFC}.
\end{proof}

\begin{remark}
    The control schemes applied when $\|e_r(t_i)\| < \FCDiscreteThresh$ is not
    required to achieve any tracking guarantees. The only requirement is that
    the control signal~$u_{\mathrm{data}}$ satisfies
    $\|u_{\mathrm{data}}\|_\infty \le \umax$ for given $\umax\geq 0$.
    In particular, this means that \emph{any} controller (predictive, or
    learning-based, or model inversion-based, or locally stabilising) applied in
    the safe region given it satisfies the input constraints defined by~$\umax$.
    Moreover, a control scheme applied in the safe region is not even supposed
    to be suitable for the system to be controlled. This  means that it is
    possible to apply, for example, controllers designed for discrete-time systems to
    the continuous-time system to be controlled. Maintenance of the tracking
    behaviour is still ensured by \Cref{Thm:Combined_Controller}. 
\end{remark}

The versatility of the proposed framework \eqref{eq:combined_controller} has been demonstrated
in~\cite{LanzaDenn24,Schmitz23} through its application to prominent data-driven predictive control schemes,
specifically data-driven model predictive control and Reinforcement Learning (RL). 
The data-driven MPC scheme presented therein builds on Willems et al.'s so-called fundamental
lemma~\cite{WRMDM05}, which enables a non-parametric description of the system's
input-output behaviour from measurement data, see
also~\cite{MarkDorf21,faulwasser2023behavioral} and the references therein.
This combined control approach elevates standard MPC to a
data-enabled predictive control scheme, cf.~\cite{berberich2022linear,coulson2019data}.
In~\cite{LanzaDenn24}, $Q$-learning -- first developed ~\cite{watkins1989learning}
and now a cornerstone of RL supporting many derivative algorithms~\cite{jang2019q} --
illustrates how the controller~\eqref{eq:controller_recursive} combines with model-free RL
techniques. This integration both safeguards the learning process and 
enhances the control signal via the strategy~\eqref{eq:combined_controller}.
Although \Cref{Thm:Combined_Controller} requires a shorter sampling period 
to ensure compliance with the control objective,
the two-component data-driven controller~\eqref{eq:combined_controller} 
outperformed the pure feedback controller~\eqref{eq:controller_recursive} in both cases.

\subsection{Simulation}
\begin{figure}[h]
    \begin{subfigure}[b]{0.49\textwidth}
        \centering
        \includegraphics[width=\linewidth]{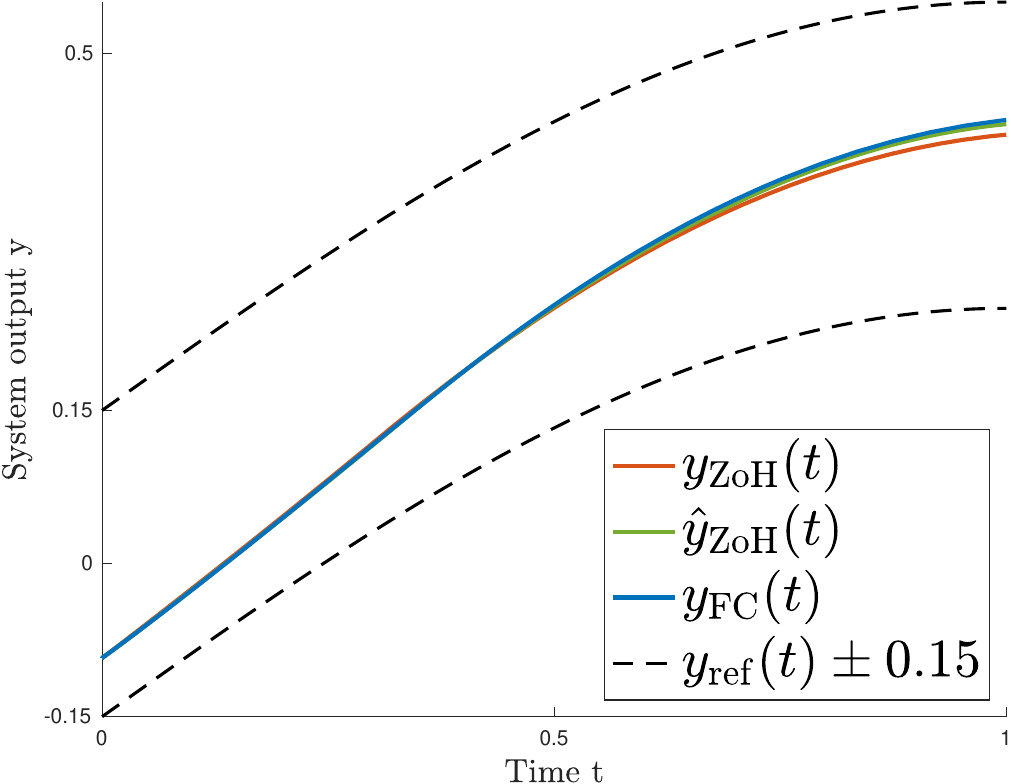}
        \caption{Outputs and reference, with error boundary.}
        \label{fig:sim:fc:zoh:output}
    \end{subfigure}
      \begin{subfigure}[b]{0.49\textwidth}
        \centering
        \includegraphics[width=\linewidth]{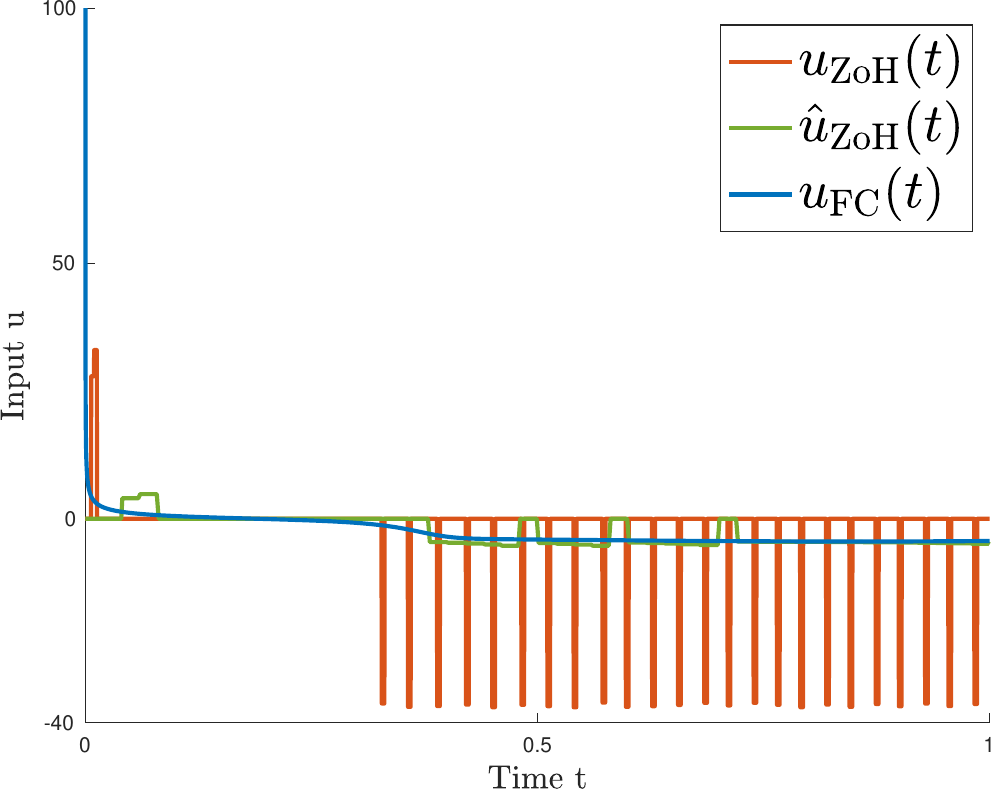}
        \caption{Control inputs.}
        \label{fig:sim:fc:zoh:input}
    \end{subfigure}
    \caption{Simulation of system~\eqref{eq:MassOnCarInputOutputDeg2} under the control of
    the zero-order-hold control law~\eqref{eq:controller_recursive} and the funnel controller \cite{BergIlch21}.}
    \label{fig:sim:fc:zoh}
\end{figure}

For the purpose of illustration, we revisit the mass-on-car
system~\cite{SeifBlaj13} from \Cref{Sec:FMPC:Sim:MassOnCar}, and compare the ZoH
controller~\eqref{eq:controller_recursive} with the funnel controller presented
in~\cite{BergIlch21}. Given the parameters $m_1=1$, $m_2=2$, spring constant $k=1$, damping~$d=1$, 
and angle $\vartheta = \pi/4$,
the system takes the form
\begin{equation}\tag{\ref{eq:MassOnCarInputOutputDeg2} revisited}
    \begin{aligned}
        \ddot y(t)    &= R_1y(t)+ R_2\dot{y}(t) + S\eta(t) +\Gamma u(t)\\
        \dot \eta(t)  &= Q \eta(t) + P y(t), 
    \end{aligned}
\end{equation}
with  initial conditions $[y(0),\dot{y}(0)]=[y_0^0,y_1^0]\in\R^2$ and  $\eta(0)=\eta^0\in\R^2$ for 
\[
    R_1=0,\quad 
    R_2=\frac{1}{4},\quad
    S=\frac{-\sqrt{2}}{8}
    \begin{bmatrix}
       1&1 
    \end{bmatrix},\quad
    \Gamma=\frac{1}{4},\quad
    Q=\begin{bmatrix}
    0&1\\
    -1&-1
    \end{bmatrix},\quad
    P=\sqrt{2}
    \begin{bmatrix}
       1\\
       0
    \end{bmatrix}\!.
\]
We simulate output reference tracking of the signal $y_{\rf}(t) = 0.4 \sin(\tfrac{\pi}{2}t)$ for $t \in [0,1]$,
transporting the mass~$m_2$ on the car from position~$0$ to $0.4$ within chosen error boundaries of $\pm 0.15$.
We choose the activation threshold~${\FCDiscreteThresh = 0.75}$.
With these parameters a brief calculation (using the variation of constants formula for the internal dynamics)
yields $\fmax \le 1.4$, $\gmax= \gmin = 0.25$, and hence,
the sampling time $\SampleTime\le 3.2 \cdot 10^{-3}$, and the input gain $\FCDiscreteGain \ge 27.78$,
which guarantee success of the tracking task according to \Cref{Thm:DiscreteFC}.
Choosing the smallest~$\FCDiscreteGain$, this already results in $\|u_{\mathrm{ZoH}}\|_\infty \le \FCDiscreteGain/\FCDiscreteThresh  \le 37.04$.
We start with a small initial tracking error of $y(0) = -0.0925$, and $\dot y(0) = \dot y_{\rf}(0)$.
The simulation of the controller~\eqref{eq:controller_recursive} in comparison to the continuous-time funnel controller~\cite{BergIlch21}
is displayed in \Cref{fig:sim:fc:zoh}.
The corresponding signals of the continuous-time funnel controller have the subscript $\mathrm{FC}$, i.e. $y_{\mathrm{FC}}$ and $\uFC$.
Since simulating the ZoH controller~\eqref{eq:controller_recursive} is by chance also successful for 
$\SampleTime = 2.0\cdot 10^{-2}$ and $\FCDiscreteGain = 4$ -- beyond the theoretical bounds derived in \Cref{Thm:DiscreteFC} --
the corresponding signals are also displayed and have a circumflex, i.e.~$\hat y_{\mathrm{ZoH}}$ and $\hat u_{\mathrm{ZoH}}$.
\Cref{fig:sim:fc:zoh:input} shows the system's output alongside the reference trajectory within the error tolerance bounds.
Note that although the control input is discontinuous for the control law~\eqref{eq:controller_recursive}, 
the output signal remains continuous due to integration.
The corresponding input signals are shown in \Cref{fig:sim:fc:zoh:output}.
The three considered controllers achieve the tracking task.
The ZoH input consists of separated pulses for two primary reasons. 
First, the control law~\eqref{eq:controller_recursive} uses (undirected) worst-case estimates $\gmin, \gmax$ and~$\fmax$ 
to compute the input signal.
Hence, the control signal is at many time instants unnecessary large; however, it is ensured that the control signal always sufficiently large.
Second, \eqref{eq:controller_recursive} includes the activation threshold~$\FCDiscreteThresh$,
rendering the controller is inactive when the tracking error is small.
If the tracking error exceeds this threshold at a sampling instant,
the applied input is sufficiently large (due to the worst case estimations) to force the error back below the threshold by the next sampling instant.
Thus, at this time instant the input is determined to be zero.
Consequently, the worst-case estimations combined with the ZoH implementation inevitably produce a peaky control signal.
The control signal $\hat u_{\mathrm{ZoH}}$ (green) is also peaky, but exhibits smaller magnitude (due to smaller~$\FCDiscreteGain$)
and larger pulse width (due to larger~$\SampleTime$). Overall, $\hat u_{\mathrm{ ZoH}}$ is comparable to $ \uFC$.
The successful simulation with these parameters suggests potential for finding better estimates 
of sufficient control parameters~$\FCDiscreteGain,\SampleTime$ in future work.
The control performance could also be enhanced using the extension discussed in \Cref{Sec:FCAsSaftyFilter}.
Note that the control signal~$\uFC$ also has a large initial peak, with $\|\uFC\|_\infty \approx 100$.
For the simulation, we used \textsc{Matlab}. The corresponding source code can be
found on \textsc{GitHub} under the link \url{https://github.com/ddennstaedt/FMPC_Simulation}.
For the integration of the dynamics, the
routine \texttt{ode15s} with $\mathrm{AbsTol}=\mathrm{RelTol}=10^{-6}$ and 
adaptive step size was utilised.
To simulate the system behaviour under control of the funnel
controller~\cite{BergIlch21}, \texttt{ode15s} produces a maximal step size
of $\approx 3.99 \cdot 10^{-2}$ and a minimal step size of $\approx 1.21 \cdot
10^{-6}$.
Thus, the largest step is about twelve times larger than $\SampleTime$,
and the smallest time step is about 4000 times smaller than~$\SampleTime$. 
Due to the worst case estimates used in the proof of \Cref{Thm:DiscreteFC} 
to derive the bounds for $\SampleTime$, the proposed framework \eqref{eq:combined_controller}
requires a higher sampling rate than the funnel controller $\uFC$ during most of the time.
However, there are currently no results regarding an upper limit for the sampling rate of $\uFC$.
Especially for unfavourable initial values, it can become arbitrarily large.

\section{Sampled-data funnel MPC}\label{Sec:FMCPZoH}
In this section, we adapt the funnel MPC~\Cref{Algo:FunnelMPC} -- designed to
achieve the control objective outlined in~\Cref{Sec:ControlObjective} -- to
operate under sampled-data constraints. Unlike the prior learning-based and robust
funnel MPC formulations explored in~\Cref{Chapter:FunnelMPC,Chapter:RobustFunnelMPC,Chapter:LearningRobustFMPC},
the space of admissible controls is now restricted to step functions, where the
control signal may only change finitely often between two sampling instants. 

Sampling can have a profound impact on both the stability and performance of both linear and non-linear 
model predictive control schemes, as analysed in \cite{worthmann2014role}.
Consequently, a variety of sampled-data MPC schemes for continuous-time systems have been 
developed~\cite{Geromel2022,Worthmann2015}.
Notably, \cite{yuz2005sampled} derives discrete-time model approximations for
continuous-time systems, whose solution error scales with the sampling time.
Complementary approaches like event-triggered MPC \cite{Brunner2017} 
further optimise digital implementations by updating control actions only when
necessary, reducing computational overhead without sacrificing stability.

In contrast to these existing frameworks, we reformulate the funnel
MPC~\Cref{Algo:FunnelMPC} as a sampled-data scheme building on the ZoH-funnel
controller framework developed in \Cref{Sec:FCZoH}.
By constraining controls to step functions, 
we propose the following modification of the funnel MPC~\Cref{Algo:FunnelMPC}.

\begin{algo}[Sampled-data funnel MPC]\label{Algo:DiscrFunnelMPC}\ \\
    \textbf{Given:} Model~\eqref{eq:Model_r} with initial time $t_0\in\Rp$ and initial value $\yM^0\in\cC^{r-1}([0,t_0],\R^m)$,
    reference signal $y_{\rf}\in W^{r,\infty}(\Rp,\R^{m})$, signal memory length~$\tau\geq0$,
    a set of funnel boundary function $\Psi=(\psi_1,\ldots,\psi_r)\in\FunnelBoundaryFuncs$ with corresponding parameters $k_i$ for $i=1,\ldots, r$,
    input saturation level $\umax\geq0$, a maximal step length~$\SampleTime>0$, funnel stage cost function~$\FunnelStageCost$, and 
    a $\tau$-initialisation strategy $\InitStrategy$ as in~\Cref{Def:InitialisationStrategy}.\\
    \textbf{Set} the time shift $\delta >0$, 
                 the prediction horizon $T\geq\delta$,
                 index $k\coloneqq 0$,
                 and $\xMh^0\coloneqq \OpChi(\yM^0)$.
                 Choose a partition $\cP=(t_i)_{i\in\N_0}$ of the
                 interval~$[t_0,\infty)$ with~$\Abs{\cP}\leq\SampleTime$ and which
                 contains $(t_0 + i\delta)_{i\in\N_0}$ as a subsequence.\\
    \textbf{Define} the time sequence~$(\hat{t}_k)_{k\in\N_0} $ by $\hat{t}_k \coloneqq  t_0+k\delta$.\\ 
    \textbf{Steps}:
    \begin{enumerate}[(a)]
    \item\label{agostep:DiscrFunnelMPCFirst}
        Select initial model state $\InitStateK_k\coloneqq \InitStrategy(\xMh^k)\in\InitValues(\hat{t}_k)$ at current time $\hat{t}_k$
        based on $\xMh$.
    \item
        Compute a solution $\uFMPCk \in \cT_{\cP}([\hat{t}_k,\hat{t}_k +T],\R^{m})$ of 
    \begin{equation}\label{eq:DiscrFunnelMpcOCP}
        \mathop
                {\operatorname{minimise}}_{\substack
                {
                    u\in \cT_{\cP}([\hat{t}_k,\hat{t}_k+T],\R^{m}),\\
                    \SNorm{u}  \leq \umax 
                }
            }\      \int_{\hat{t}_k}^{\hat{t}_k + T}\FunnelStageCost(s,\eM_{r}(\xM(s;\hat{t}_k,\InitStateK_k,u)-\OpChi(y_{\rf})(s)),u(s))\d{s}.
    \end{equation}
    \item Apply the control law
        \begin{equation}\label{eq:DiscrFMPC-fb}
            \mu:[\hat{t}_k,\hat{t}_{k+1})\times\InitValues(\hat{t}_k)\to\R^m, \quad \mu(t,\xMh^k) =\uFMPCk(t)
        \end{equation}
        to model \eqref{eq:SysMod_Discrete} with initial time and data
        $(t_k,\InitStateK_{k})$ and obtain, on the interval
        ${I_{0}^{t_{k+1},\tau}\coloneqq [t_{k+1}-\tau,t_{k+1}]\cap[0,t_{k+1}]}$
        a measurement of the model's output and its derivatives
        ${\xMh^{k+1}\coloneqq
        \xM(\cdot;t_k,\InitStateK_{k},\uFMPCk)|_{I_{t_0}^{t_{k+1},\tau}}}$. 
        Increment $k$ by 1 and go to Step~\ref{agostep:DiscrFunnelMPCFirst}.
    \end{enumerate}
\end{algo}

\begin{remark}
Note that while the time shift $\delta>0$ is an upper bound for the step length~$\SampleTime>0$ of the  control signals,
$\delta$ is allowed to be larger than $\SampleTime$ under the condition that the partition~$\cP$ contains $(t_0+i\delta)_{i\in\N_0}$ as a subsequence.
In this case,  several control signals are applied to the system between two steps of the MPC~\Cref{Algo:DiscrFunnelMPC}.
This can also be interpreted as a multi-step MPC scheme, cf.~\cite{worthmann2014role}.
\end{remark}

\begin{theorem}\label{Thm:DiscreteFMPC}
    Consider model~\eqref{eq:SysMod_Discrete} with $(\fM,\gM,\oTM)\in\ModSysClassDiscr$ with initial trajectory ${\yM^0\in\cC^{r-1}([0,t_0],\R^m)}$.
    Let $y_{\rf}\in W^{r,\infty}(\Rp,\R^{m})$ and $\Psi=(\psi_1,\ldots,\psi_r)\in\FunnelBoundaryFuncs$ be given.  
    Further, let $\tau\geq0$ be greater than or equal to the memory limit of operator~$\oTM$ and 
    ${\InitStrategy:\bigcup_{\hat{t}\geq t_0}\cR(I_0^{\hat{t},\tau},\R^{rm})\to
        \bigcup_{\hat{t}\geq t_0}\cR(I_0^{\hat{t},\tau},\R^{rm})\times L^\infty_{\loc}( I_{t_0}^{\hat{t},\tau},\R^q)
    }$ be an $\tau$-initialisation strategy as in~\Cref{Def:InitialisationStrategy}. 
    Then, there exists $\umax\geq0$ and a maximal step length $\SampleTime>0$
    such that the sampled-data funnel MPC~\Cref{Algo:DiscrFunnelMPC} 
    with $\delta>0$, $T\ge\delta$, and a partition $\cP$ of the interval $[t_0,\infty)$ with $\Abs{\cP}\leq\SampleTime$ 
    is initially and recursively feasible, i.e. 
    \begin{itemize}
        \item the OCP~\eqref{eq:DiscrFunnelMpcOCP} has a solution $\uFMPCk\in \cT_{\cP}([\hat{t}_k,\hat{t}_k+T],\R^m)$ at every time instant $\hat{t}_k\coloneqq  t_0+\delta k$ for $k\in \N_0$, and
        \item the model~\eqref{eq:SysMod_Discrete} with applied funnel MPC feedback~\eqref{eq:DiscrFMPC-fb}
        has a concatenated solution $\xM : [0,\infty) \to \R^{rm}$ in the sense of \Cref{Def:SolutionClosedLoop}. 
    \end{itemize}
    The corresponding input is given by 
    \[
        u_{\mathrm{FMPC}}(t) = \uFMPCk(t),
    \]
    for $t\in[\hat{t}_k,\hat{t}_{k+1})$ and $k\in\N_0$.
    Each global solution $\xM$ with corresponding output $\yM$ and input $u_{\mathrm{FMPC}}$ satisfies: 
    \begin{enumerate}[(i)]
        \item the control input is bounded by $\umax$, i.e.
        \[
            \fa t \ge t_0 :\quad \Norm{u_{\mathrm{FMPC}}(t)}\leq\umax,
        \]
        \item
        the tracking error between the model output and the reference evolves within prescribed boundaries, i.e.
        \[
            \fa t \ge t_0 :\quad \Norm{\yM(t)  - y_{\rf}(t)} < \psi_1(t) .
        \]
    \end{enumerate}
\end{theorem}

To prove \Cref{Thm:DiscreteFMPC}, we reformulate certain results from
\Cref{Chapter:FunnelMPC} adapted to the changed setting. 
Most importantly, one has to show that there exists a step function $u$ that, if
applied to the model~\eqref{eq:Model_r} at time $\hat{t}$, ensures that $\xM(t)
- \OpChi(y_{\rf})(t)$ evolves within $\cD_{t}^{\Psi}$ for all $t$ over the next
time interval of length $T>0$. 
For a step function with partition $\cP$ to achieve this objective it has to be an element of
\begin{align}\label{eq:Def-DiscrU}
        \Controls^{\cP}(\umax,\InitState) \coloneqq \cT_{\cP}([\hat{t},\hat{t}+T],\R^m)\cap\Controls(\umax,\InitState).
\end{align}
\Cref{Th:ExUmax} shows that there exists a bound $\umax\geq0$ on the control
input such that the set $\Controls(\umax,\InitState)$ is non-empty.
To prove that there exists a step function $u\in \cT_{\cP}([\hat{t},\hat{t}+T],\R^m)$
with a uniform minimal step length $\SampleTime>0$ that is an element of $\Controls(\umax,\InitState)$,
we utilise ideas from \Cref{Thm:DiscreteFC}.
The difficulty lies in the usage of different auxiliary error variables.
\Cref{Thm:DiscreteFC} shows that there exists a piece-wise constant control
ensuring the evolution of $\OpChi(y - y_{\rf})(t)$ within the set
$\cEFC{1}(\phi(t))$ for all $t\geq t_0$. 
To be used in the discrete funnel MPC \Cref{Algo:DiscrFunnelMPC},
this result has to be also verified utilising the error signals $\eM_i$ as in~\eqref{eq:ErrorVar}
(the set $\cEFC{1}(\phi(t))$ is defined in terms of the error variables $\eS_i$ in \eqref{eq:ek_FC}). 

\begin{lemma}\label{Lem:DiscrUNonEmpty}
    Consider model~\eqref{eq:SysMod_Discrete} with $(\fM,\gM,\oTM)\in\ModSysClassDiscr$.
    Let $\tau\geq 0$ be greater than or equal to the memory limit of operator $\oTM$.
    Further, let $y_{\rf}\in W^{r,\infty}(\Rp,\R^{m})$ and $\Psi\in\FunnelBoundaryFuncs$. 
    Then, there exists $\umax\geq0$ and $\SampleTime>0$ such that, 
    for $\hat{t}\geq t_0$, $\InitState\in\InitValues(\hat{t})$, $T>0$,
    and every partition $\cP$ of the interval $[\hat{t},\hat{t}+T]$ with $\Abs{\cP}\leq\SampleTime$, we have
    \[
        \Controls^{\cP}(\umax,\InitState)\neq\emptyset.
    \]
\end{lemma}
\begin{proof}
    To prove the existence of a step function achieving the control objective,
    we combine the ideas from~\Cref{Thm:DiscreteFC} and \Cref{Lemma:DynamicBounded} in the following.

    \noindent
    \emph{Step 1}: We define $\umax\geq0$ and $\SampleTime$. 
    As in the proof of \Cref{Th:ExUmax}, define, for $i=1,\ldots, r-1$ and $j=0,\ldots,r-i-1$,
    \[
        \mu_i^0 \coloneqq  \SNorm{\psi_i},\quad \mu_{i}^{j+1}\coloneqq  \mu_{i+1}^{j}+k_{i}\mu_{i}^{j},
    \]
    where $k_i\geq0$ are the to $\Psi$ associated constants,
    which are also used to define the error variables $\eM_i$ as in~\eqref{eq:ErrorVar}.
    Utilising the constants $\fMmax$, $\gMmax$, and $\gMmin$ from~\Cref{Lemma:DynamicBounded}, define 
    \[
        \kappa_0\coloneqq \SNorm{\tfrac{1}{\Funnel_r}}\rbl\fMmax+\SNorm{y_{\rf}^{(r)}}+\sum_{j=1}^{r-1}k_j\mu_j^{r-j}+\SNorm{\dot{\Funnel}}\rbr
    \]
    and choose an input gain
    \[
        \FCDiscreteGain > \frac{2\kappa_0\inf_{s\geq t_0}\Funnel_{r}(s)}{ \gMmin}.
    \]
    With $\kappa_1\coloneqq \kappa_0+2\SNorm{\tfrac{1}{\Funnel_r}}  \gMmax\FCDiscreteGain$, we define the constants
    \[
        \SampleTime \coloneqq  \min\cbl \frac{\kappa_0}{\kappa_1^2},\frac{1}{2\kappa_0}\cbr\quad\text{and}\quad
        \umax\coloneqq 2\FCDiscreteGain.
    \]
    All parameters are chosen in a similar fashion as in~\Cref{Thm:DiscreteFC}.
    
    \noindent
    \emph{Step 2}: Let $T>0$, $\hat{t}\geq t_0$, and $(\xMh,\oTMh)=\InitState\in\InitValues(\hat{t})$ be arbitrary but fixed.
    Further, let  $\cP=(t_i)_{i\in\N_0}$ be a partition of the interval $[\hat{t},\infty)$ with $\Abs{\cP}\leq\SampleTime$.
    Note that by $\cP$ is also a partition of the interval $[\hat{t},\hat{t}+T]$ by being a partition of the interval $[\hat{t},\infty)$,
    see~\Cref{Def:PartitionAndStepFunction}.
    We construct a control step function $u$ and show that $u\in\Controls^{\cP}(\umax,\InitState)$.
    To this end, for some $u\in L^{\infty}([\hat{t},\infty),\R^m)$, we use the shorthand notation
    $\xM(t)\coloneqq  \xM(t;\hat{t},\InitState,u)$ and $\eM_i(t)\coloneqq  \eM_{i}(\xM(t)-\OpChi(y_{\rf})(t))$
    for $i=1,\ldots,r$.
    The application of the ZoH feedback control 
    \begin{equation} \label{eq:ZoHControl}
         u_{\mathrm{ZoH}}(t) =
        \begin{cases}
            0,  &  \Norm{\tfrac{\eM_r(t_i)}{\Funnel_{r}(t_i)}} < \frac{1}{2}\\[2ex]
             - \FCDiscreteGain\tfrac{ \Funnel_r(t_i) \eM_{r}(t_i)}{\|\eM_r(t_i)\|^2}  ,  &\Norm{\tfrac{\eM_r(t_i)}{\Funnel_{r}(t_i)}} \geq \frac{1}{2},
        \end{cases}
    \end{equation}
    to the system~\eqref{eq:SysMod_Discrete} leads to a closed-loop system. 
    If this initial value problem is considered on the interval~$[\hat{t},\hat{t}+T]$ with initial conditions $(\hat{t},\InitState)$ as in \eqref{eq:Model_InitialValue},
    then an application of~\Cref{Prop:SolutionExists}
    yields the existence of a maximal solution~$\xM:[0,\omega)\to\R^{rm}$ in the sense of~\Cref{Def:ModSolution}.
    If~$\xM$ is bounded, then $\omega=\infty$, see~\Cref{Prop:SolutionExists}~\ref{Item:Theorem:BoundedSolution}.
    In this case, the solution exists on $[0, \hat t + T]$.
    
    \noindent
    \emph{Step 3}:
    We show that $\Norm{\tfrac{\eM_{r}(t_i)}{\Funnel_r(t_i)}}< 1$ for $i\in\N_0$ implies 
    $\Norm{\tfrac{\eM_{r}(t)}{\Funnel_r(t)}}<1$ for all $t\in[t_i,t_{i+1}]$.
    Seeking a contradiction, suppose that there exists a maximal $i\in\N_0$ such that 
    we have $\Norm{\tfrac{\eM_{r}(t)}{\Funnel_r(t)}}<1$ for all $t\in[\hat{t},t_{i}]$ and
    $\Norm{\tfrac{\eM_{r}(t)}{\Funnel_r(t)}}\geq 1$ for some $t\in(t_i,t_{i+1})$.
    Then, there exists 
    \[
        t^* \coloneqq  \inf \setdef{ t \in (t_i,t_{i+1})}{ \Norm{\tfrac{\eM_r(t)}{\Funnel_r(t)} }\geq 1} .
    \]
    We have $\Norm{\tfrac{\eM_r(\hat{t})}{\Funnel_r(\hat{t})}}<1$ 
    by the assumption $\xM(\hat{t}) - \OpChi(y_{\rf})(\hat{t})\in\cD_{\hat{t}}^\Psi$, see also~\Cref{Rem:InitialValueInFunnel}.
    This yields $\Norm{\tfrac{\eM_r(t)}{\Funnel_r(t)}}<1$ for all $t\in[\hat{t},t^\star)$.
    This implies, according to \Cref{Prop:OnlyLastFunnel},  $\xM(t)-\OpChi(y_{\rf})(t)\in\cD_{t}^\Psi$ for all $t\in[\hat{t},t^\star)$, 
    i.e. $\Norm{\eM_i(t)}<\Funnel_i(t)$ for all $i=1,\ldots, r$.
    Thus, $\Norm{\eM_i(t)}\leq\mu_i^0$ for all $i=1,\ldots, r$.
    Invoking boundedness of $y_{\rf}^{(i)}$, $i=0,\ldots,r$, and the relation in~\eqref{eq:reps_ei_chi},
    we may infer that $\xM$ is bounded on $[\hat{t},t^\star]$.
    Hence, $\omega> t^\star$.
    Since ${(\xMh,\oTMh)=\InitState\in\InitValues(\hat{t})}$, there exists 
    a function $\zeta\in\FunnelTrajectories_{\hat{t}}$ such that
    ${\zeta|_{[\hat{t}-\tau,\hat{t}]\cap[0,\hat{t}]}=\xMh}$
    and $\oTM(\zeta)|_{[\hat{t}-\tau,\hat{t}]\cap[t_0,\hat{t}]}=\oTMh$.
    Moreover, the function $\xM$ fulfils 
    $\xM(t)|_{[\hat{t}-\tau,\hat{t}]\cap[0,\hat{t}]}=\xMh$
    and ${\oTM(\xM)|_{[\hat{t}-\tau,\hat{t}]\cap[t_0,\hat{t}]}=\oTMh}$ 
    because $\xM$ satisfies the initial conditions~\eqref{eq:Model_InitialValue}.
    Define the function $\tilde{\zeta}\in\cR(\Rp,\R^{rm})$ by 
    \[
        \tilde{\zeta}(t)=\begin{cases}
            \xM(t),& t\in [\hat{t},t^\star)\\
            \zeta(t),& t\in \Rp\backslash[\hat{t},t^\star).
        \end{cases}
    \]
    Then, $\tilde{\zeta}$ is an element of $\FunnelTrajectories_{s}$ for all $s\in[\hat{t},t^\star)$
    because  $\zeta\in\FunnelTrajectories_{\hat{t}}$ and 
    $\xM(t)-\OpChi(y_{\rf})(t)\in\cD_{t}^\Psi$ for all $t\in[\hat{t},t^\star)$.
    Hence, we have 
    $\Norm{\fM(\oTM(\tilde{\zeta})(t))}\leq \fMmax$ and 
    ${\Norm{\gM(\oTM(\tilde{\zeta})(t))^{-1}}\leq \gMInvmax}$ for all $t\in[\hat{t},t^\star)$ 
    according to~\Cref{Lemma:DynamicBounded}.
    Since $\tau\geq 0$ is greater than or equal to the memory limit of operator $\oTM$, we have
    \[
        \oTM(\xM)(t)=\oTM(\tilde{\zeta})(t)
    \]
    for all $t\in[\hat{t},t^\star)$. Thus, 
    $\Norm{\fM(\oTM(\xM)(t))}\leq \fMmax$ and 
    $\Norm{\gM(\oTM(\xM)(t))^{-1}}\leq \gMInvmax$ for all $t\in[\hat{t},t^\star)$. 
    Using~\eqref{eq:ErrorVarDyn} and the definition of~$\mu_i^j$, it follows that
    \[
        \fa t\in [\hat{t},t^\star): \Norm{\eM^{(j+1)}_i(t)}
        =
        \Norm{\eM^{(j)}_{i+1}(t)-k_i\eM^{(j)}_{i}(t)}
        \leq
        \mu_{i+1}^j+k_i\mu^{j}_{i}
        =\mu^{j+1}_{i} 
    \]
    inductively for all $i=1,\ldots, r-1$ and $j=0,\ldots,r-i-1$.
    Utilising again~\eqref{eq:ErrorVarDyn}, it follows by induction that
    \[
        \eM_r(t)=\eM_1^{(r-1)}(t)+\sum_{j=1}^{r-1}k_j\eM_j^{(r-j-1)}(t).
    \]
    Omitting the dependency on $t$, we calculate for $t \in [\hat{t},t^\star)$:
    \begin{equation}\label{eq:DiffEmFunnel}
    \begin{aligned}
    \Dd{t}\frac{\eM_r}{\Funnel_r} &= \frac{\dot{\eM}_r\Funnel_r-\eM_r\dot{\Funnel}_r}{\Funnel_r^2}
    =\frac{1}{\Funnel_r}\rbl\eM_1^{(r)}+\sum_{j=1}^{r-1}k_j\eM_j^{(r-j)}-\eM_r\frac{\dot{\Funnel}_r}{\Funnel_r}\rbr\\
    &=\frac{1}{\Funnel_r}\rbl\fM(\oTM(\xM))+\gM(\oTM(\xM))u-y_{\rf}^{(r)}+\sum_{j=1}^{r-1}k_j\eM_j^{(r-j)}-\eM_r\frac{\dot{\Funnel}_r}{\Funnel_r}\rbr.
    \end{aligned}
    \end{equation}
    We now  consider the two cases $\Norm{\tfrac{\eM_{r}(t_i)}{\Funnel_r(t_i)}}<\tfrac{1}{2}$ 
    and $\Norm{\tfrac{\eM_{r}(t_i)}{\Funnel_r(t_i)}}\geq \tfrac{1}{2}$ separately.
    
    \noindent
    \emph{Step 3.a}: We consider $\Norm{\tfrac{\eM_{r}(t_i)}{\Funnel_r(t_i)}}<\tfrac{1}{2}$. 
    By definition~\eqref{eq:ZoHControl}, we have $u_{\mathrm{ZoH}}(t) = 0$ for all $t\in[t_i,t^\star)$.
    With \eqref{eq:DiffEmFunnel}, we have $\Norm{ \rbl\dd{t}\tfrac{\eM_r}{\Funnel_r}\rbr (t)}\leq\kappa_0$
    for $t\in[\hat{t},t^\star)$.
    Thus, we calculate 
    \begin{align*}
      1 =  \Norm{\tfrac{\eM_r(t^*)}{\Funnel_r(t^\star)} } \le \Norm{\tfrac{\eM_r(t_i)}{\Funnel_r(t_i)}}  + 
            \int_{t_i}^{t^*}  \Norm{ \rbl\dd{t}\tfrac{\eM}{\Funnel_r}\rbr\! (s)}  \d s 
        \le  \Norm{\tfrac{\eM_r(t_i)}{\Funnel_r(t_i)}} +  \int_{t_i}^{t^*}\! \kappa_0  \d s 
        < \tfrac{1}{2} + \kappa_0 \SampleTime\leq 1,
    \end{align*}
    where $t^* < \SampleTime \leq \tfrac{1}{2\kappa_0}$ was used.
    This contradicts the definition of $t^*$.
    
    \noindent
    \emph{Step 3.b}:
    We consider $\Norm{\tfrac{\eM_{r}(t_i)}{\Funnel_r(t_i)}}\geq \tfrac{1}{2}$. 
    Therefore, we have $u_{\mathrm{ZoH}}(t)=- \FCDiscreteGain\tfrac{ \Funnel_r(t_i) \eM_{r}(t_i)}{\|\eM_r(t_i)\|^2}$ for all $t\in[t_i,t^\star)$.
    With \eqref{eq:DiffEmFunnel}, we have $\Norm{ \rbl\dd{t}\tfrac{\eM_r}{\Funnel_r}\rbr (t)}\leq\kappa_1$
    for all $t\in[\hat{t},t^\star)$.
    Moreover, for the expression $J(t)\coloneqq \rbl\Dd{t}\frac{\eM_r}{\Funnel_r}\rbr(t)-\tfrac{1}{\Funnel_{r}(t)}\gM(\oTM(\xM)(t))u(t)$,
    we have $\Norm{J(t)}\leq\kappa_0$ for all $t\in[\hat{t},t^\star)$.
    Thus, we calculate 
    \begin{align*}
        \dd{t} \tfrac{1}{2} \Norm{ \frac{\eM_r(t)}{\Funnel_r(t)}}^2 &= \al \frac{\eM_r(t)}{\Funnel_{r}(t)}, \rbl\Dd{t}\frac{\eM_r}{\Funnel_r}\rbr(t) \ar \\
         &= \al \frac{\eM_r(t_i)}{\Funnel_r(t_i)} +   \int_{t_i}^t \rbl\Dd{t}\frac{\eM_r}{\Funnel_r}\rbr(s) \d s, \rbl\Dd{t}\frac{\eM_r}{\Funnel_r}\rbr (t) \ar \\
         &= \al \frac{\eM_r(t_i)}{\Funnel_r(t_i)} , \rbl\Dd{t}\frac{\eM_r}{\Funnel_r}\rbr (t)\ar 
         +\al \int_{t_i}^t \rbl\Dd{t}\frac{\eM_r}{\Funnel_r}\rbr(s) \d s, \rbl\Dd{t}\frac{\eM_r}{\Funnel_r}\rbr (t)\ar\\
         &\leq \al \frac{\eM_r(t_i)}{\Funnel_r(t_i)} , J(t)+\frac{1}{\Funnel_{r}(t)}\gM(\oTM(\xM)(t))u(t)\ar +(t-t_i)\SNorm{\Dd{t}\frac{\eM_r}{\Funnel_r}}^2 \\
         &\leq \Norm{J(t)}+\al \frac{\eM_r(t_i)}{\Funnel_r(t_i)} , \frac{1}{\Funnel_{r}(t)}\gM(\oTM(\xM)(t))u(t)\ar +\SampleTime\kappa_1^2 \\
         &= \kappa_0-
         \al \frac{\eM_r(t_i)}{\Funnel_r(t_i)} , \frac{1}{\Funnel_{r}(t)}\gM(\oTM(\xM)(t)) \FCDiscreteGain\frac{ \Funnel_r(t_i) \eM_{r}(t_i)}{\|\eM_r(t_i)\|^2}\ar +\SampleTime\kappa_1^2 \\
         &\leq 2\kappa_0- \frac{\FCDiscreteGain}{\Funnel_{r}(t)}\frac{\|\eM_r(t_i)\|^2}{\Funnel_r(t_i)^2}\al \frac{\eM_r(t_i)}{\Funnel_r(t_i)}, 
             \gM(\oTM(\xM)(t)) \frac{ \eM_{r}(t_i)}{\Funnel_r(t_i)}\ar\\
         &\leq 2\kappa_0 - \frac{\FCDiscreteGain}{\inf_{s\geq t_0}\Funnel_{r}(s)}\gMmin< 0. 
    \end{align*}
    In particular, this yields $\rbl\dd{t} \tfrac{1}{2} \Norm{ \frac{\eM_r}{\Funnel_r}}^2\rbr(t_i)< 0$, by which $t^* > t_i$.
    Therefore, we find the contradiction $1 = \Norm{\frac{\eM_r(t^*)}{\Funnel_r(t^*)}}^2 < \Norm{\frac{\eM_r(t_i)}{\Funnel_r(t_i)}}^2 < 1$.

    \noindent
    \emph{Step 4}:
    As $\InitState\in\InitValues(\hat{t})$, we have $\Norm{\tfrac{\eM_{r}(\hat{t})}{\Funnel_r(\hat{t})}}< 1$, see~\Cref{Rem:InitialValueInFunnel}.
    By induction, {Step~3} yields $\Norm{\tfrac{\eM_r(t)}{\Funnel_r(t)}}<1$ for all $t\in[\hat{t},\omega)$.
    This implies, according to \Cref{Prop:OnlyLastFunnel},  ${\xM(t)-\OpChi(y_{\rf})(t)\in\cD_{t}^\Psi}$ for all $t\in[\hat{t},\omega)$, 
    i.e. $\Norm{\eM_i(t)}<\Funnel_i(t)$ for all $i=1,\ldots, r$.
    Invoking boundedness of $y_{\rf}^{(i)}$, $i=0,\ldots,r$, and the relation in~\eqref{eq:reps_ei_chi},
    we may infer that $\xM$ is bounded on $[\hat{t},\omega)$. Thus, $\omega=\infty$.
    Also note that the ZoH feedback control $u_{\mathrm{ZoH}}$ in~\eqref{eq:ZoHControl} fulfils $\SNorm{u_{\mathrm{ZoH}}}\leq\umax$.
    Since the partition  $\cP$ is also a partition of the interval $[\hat{t},\hat{t}+T]$, we have  $u_{\mathrm{ZoH}}\in\Controls^{\cP}(\umax,\InitState)$.
    This completes the proof.
\end{proof}
To prove the functioning of the discrete funnel MPC~\Cref{Algo:DiscrFunnelMPC},
we further have to show that the optimisation
problem~\eqref{eq:DiscrFunnelMpcOCP} has a solution. To this end, we recall, for
$T>0$, $\hat{t}\geq t_0$, and  $\InitState\in\InitValues(\hat{t})$ with
$\tau\geq0$ being greater than or equal to the memory limit of operator~$\oTM$, the
definition of cost functional 
${J^{\Psi}_T(\cdot;\hat{t},\InitState):L^\infty([\hat{t},\hat{t}+T],\R^{m})\to\R\cup\{\infty\}}$ given by 
\begin{align}\tag{\ref{eq:DefCostFunctionJ} revisited}
    J^{\Psi}_T(u;\hat{t},\InitState)\coloneqq 
        \int_{\hat{t}}^{\hat{t} + T}\FunnelStageCost(s,\eM_{r}(\xM(s;\hat{t},\InitState,u)-\OpChi(y_{\rf})(s)),u(s))\d{s}.
\end{align}
We prove that $J^{\Psi}_T(\cdot;\hat{t},\InitState)$, when restricted to the set 
$\Controls^{\cP}(\umax,\InitState)$ as in \eqref{eq:Def-DiscrU} has a minimum
by adapting \Cref{Th:SolutionExists} to the changed setting.
\begin{lemma}\label{Lem:DiscrSolutionExists}
    Consider model~\eqref{eq:SysMod_Discrete} with $(\fM,\gM,\oTM)\in\ModSysClassDiscr$ with reference trajectory $y_{\rf}\in W^{r,\infty}(\Rp,\R^m)$.
    Let $\Psi\in\FunnelBoundaryFuncs$ and $\tau\geq0$ be greater than or equal to the memory limit of operator~$\oTM$.
    Further, let  $\hat{t}\geq t_0$, $(\xMh,\oTMh)=\InitState\in \InitValues(\hat{t})$, $T> 0$, $\umax\geq0$,
    and $\cP$ be a partition of the interval $[\hat{t},\hat{t}+T]$
    such that $\Controls^\cP(\umax,\InitState)\neq\emptyset$.
    Then, there exists a function  $u^{\star} \in \Controls^{\cP}(\umax,\InitState)$ such that
    \[
            J^{\Psi}_T(u^\star;\hat{t},\InitState) =
            \mathop{\min}_{
                u \in \Controls^\cP(\umax,\InitState)
            }
            J^{\Psi}_T(u;\hat{t},\InitState) =
            \mathop{\min}_{
                \substack
                {
                    u\in \cT_{\cP}([\hat{t}, \hat{t}+T],\R^{m}),\\
                    \SNorm{u}  \leq \umax 
                }
            }
                J^{\Psi}_T(u;\hat{t},\InitState).
    \]
\end{lemma}
\begin{proof}
    We adapt the proof of~\Cref{Th:SolutionExists}.
    Since $\Controls^\cP(\umax,\InitState)\subset \Controls(\umax,\InitState)$,
    the set $\Controls(\umax,\InitState)$ is non-empty by assumption.
    By~\Cref{Th:SolutionExists}, there exists a control ${u \in \Controls(\umax,\InitState)}$ 
    minimising the functional $J^{\Psi}_T(u;\hat{t},\InitState)$.
    Thus, the infimum ${J^\star\coloneqq \mathop{\inf}_{ u \in \Controls^\cP(\umax,\InitState) } J^{\Psi}_T(u;\hat{t},\InitState)}$ exists as well.
    Let~$(u_k)_{k\in\N_0}\in\rbl\Controls^\cP(\umax,\InitState) \rbr^{\N_0}$ 
    be a minimising sequence, meaning $J^{\Psi}_T(u;\hat{t},\InitState)\to J^{\star}$. 
    As $\cP=(t_k)_{k\in\N_0}$ is a partition of the interval~$[\hat{t},\hat{t}+T]$,
    we have $t_0=\hat{t}$ and  there exists a minimal $N\in\N_0$  with $t_n>\hat{t}+T$ for all $n>N$. 
    Define $u_{i,k}\coloneqq  u_k(t_i)$ for $i=0,\ldots, N$. 
    For every ${i=0,\ldots,N}$, $(u_{i,k})_{k\in\N_0}$ is a sequence in $\R^m$ 
    with $\Norm{u_{i,k}}\leq u_{\max}$ for all $k\in \N$.
    Thus, it has a limit point~$u_i^{\star}\in\R^m$. 
    The function $u^\star$ defined by $u^\star|_{[t_i,t_{i+1})\cap [\hat{t},\hat{t}+T]}\coloneqq u^\star_i$ 
    is an element of ${\cT_{\cP}([\hat{t},\hat{t}+T],\R^m)}$ with $\Norm{u^\star}\leq u_{\max}$.
    Up to subsequence,~$u_k$ converges uniformly to~$u^\star$.
    We define ${(x_k)\coloneqq \rbl\xM(\cdot;\hat{t},\InitState,u_k)\rbr\in \cR([0,\hat{t}+T],\R^n)^\N}$
    as the sequence of associated responses.
    Repeating Steps~2 and 3 of the proof of~\Cref{Th:SolutionExists},
    the sequence $(x_k)$ has a subsequence (which we do not relabel) that uniformly 
    converges to $x^\star=\xM(\cdot;\hat{t},\InitState,u^\star)$.
    It remains to show that $u^\star\in\Controls^\cP(\umax,\InitState)$, $J^\star=J^{\Psi}_T(u^\star;\hat{t},\InitState)$, and 
    that $J^{\Psi}_T(u^\star;\hat{t},\InitState)= \mathop{\min}_{ u \in
    \Controls^\cP(\umax,\InitState) }$.
    These statements follow along the lines of Steps~{5--7} of the proof of~\Cref{Th:SolutionExists}.
\end{proof}

We are now in the position to summarise our results in the proof of~\Cref{Thm:DiscreteFMPC}.

\begin{proof}[Proof of~\Cref{Thm:DiscreteFMPC}]
    Using the results of~\Cref{Lem:DiscrUNonEmpty,Lem:DiscrSolutionExists} 
    proving~\Cref{Thm:DiscreteFMPC} is a straightforward adaptation of the proof of~\Cref{Thm:FMPC}
    to the changed context.
\end{proof}

\subsection{Simulation}
\begin{figure}[h]
         \centering
          \def\svgwidth{70pt}    
         %% Creator: Inkscape 1.3 (0e150ed, 2023-07-21), www.inkscape.org
%% Accompanies image file 'MassChain_model.pdf' (pdf, eps, ps)
\begingroup%
  \makeatletter%
  \providecommand\color[2][]{%
    \errmessage{(Inkscape) Color is used for the text in Inkscape, but the package 'color.sty' is not loaded}%
    \renewcommand\color[2][]{}%
  }%
  \providecommand\transparent[1]{%
    \errmessage{(Inkscape) Transparency is used (non-zero) for the text in Inkscape, but the package 'transparent.sty' is not loaded}%
    \renewcommand\transparent[1]{}%
  }%
  \providecommand\rotatebox[2]{#2}%
  \newcommand*\fsize{\dimexpr\f@size pt\relax}%
  \newcommand*\lineheight[1]{\fontsize{\fsize}{#1\fsize}\selectfont}%
  \ifx\svgwidth\undefined%
    \setlength{\unitlength}{595.27559055bp}%
    \ifx\svgscale\undefined%
      \relax%
    \else%
      \setlength{\unitlength}{\unitlength * \real{\svgscale}}%
    \fi%
  \else%
    \setlength{\unitlength}{\svgwidth}%
  \fi%
  \global\let\svgwidth\undefined%
  \global\let\svgscale\undefined%
  \makeatother%
  \begin{picture}(1,1.41428571)%
    \lineheight{1}%
    \setlength\tabcolsep{0pt}%
    \put(0,0){\includegraphics[width=\unitlength,page=1]{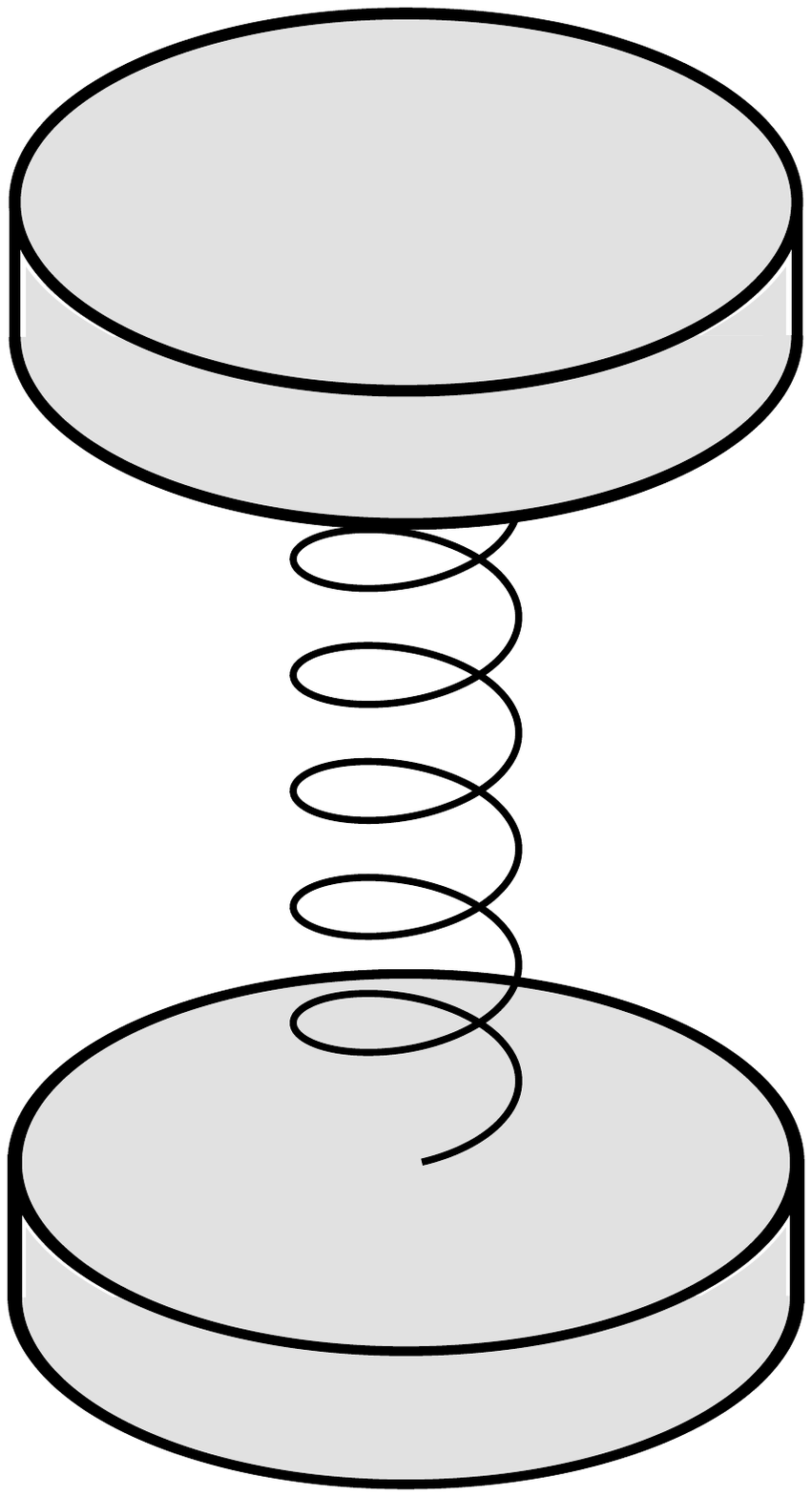}}%
    \put(0.67609275,0.49197402){\color[rgb]{0.25490196,0.46666667,0.24705882}\makebox(0,0)[lt]{\lineheight{0}\smash{\begin{tabular}[t]{l}$z_1$\end{tabular}}}}%
    \put(0.66492404,1.33212705){\color[rgb]{0.25490196,0.46666667,0.24705882}\makebox(0,0)[lt]{\lineheight{0}\smash{\begin{tabular}[t]{l}$z_2$\end{tabular}}}}%
    \put(0,0){\includegraphics[width=\unitlength,page=2]{MassChain_model.pdf}}%
    \put(0.86068449,0.27293541){\color[rgb]{0.63529412,0.10196078,0.10980392}\makebox(0,0)[lt]{\lineheight{0}\smash{\begin{tabular}[t]{l}$u$\end{tabular}}}}%
    \put(0,0){\includegraphics[width=\unitlength,page=3]{MassChain_model.pdf}}%
    \put(0.64236566,0.66745031){\color[rgb]{0,0,0}\makebox(0,0)[lt]{\lineheight{1.25}\smash{\begin{tabular}[t]{l}$k,d$\end{tabular}}}}%
  \end{picture}%
\endgroup%

         \caption{Torsional oscillator. The figure is based on \cite[Fig.~2.7]{Druecker22}, edited to the case of two flywheels for the present purpose.}
         \label{Fig:TorsOscill}
\end{figure}

\begin{figure}[hbtp]
    \begin{subfigure}[b]{0.49\textwidth}
        \centering
        \includegraphics[width=\linewidth]{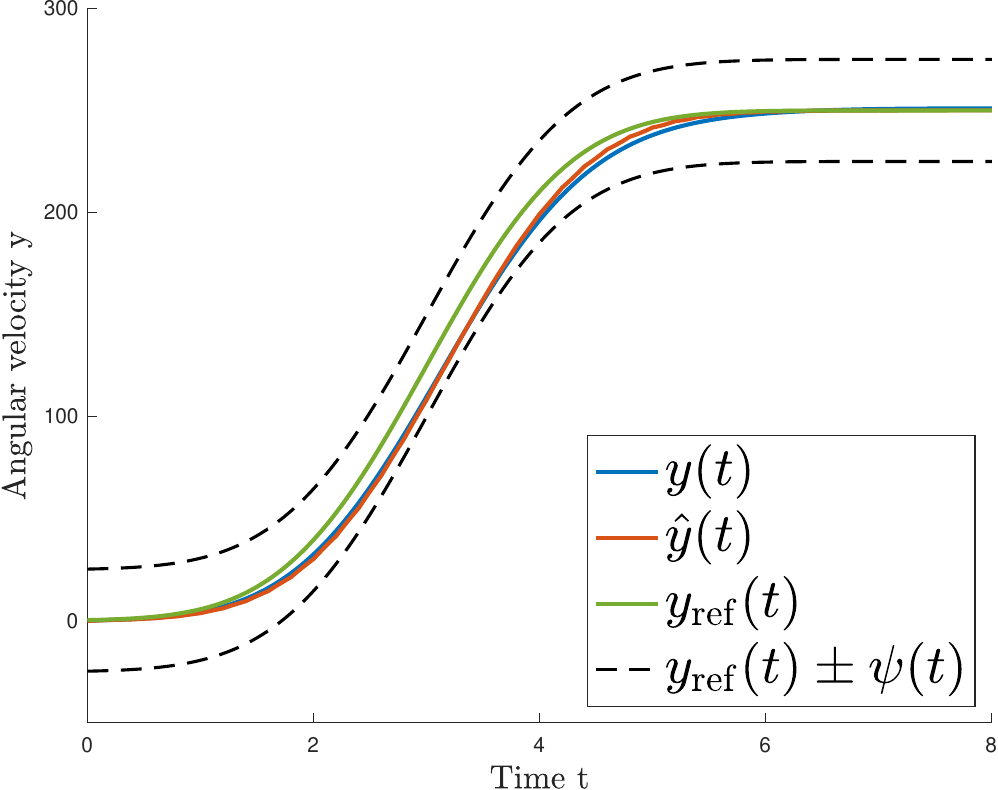}
        \caption{Outputs and reference, with error boundary.}
        \label{fig:sim:discreteFMPC:output}
    \end{subfigure}
      \begin{subfigure}[b]{0.49\textwidth}
        \centering
        \includegraphics[width=\linewidth]{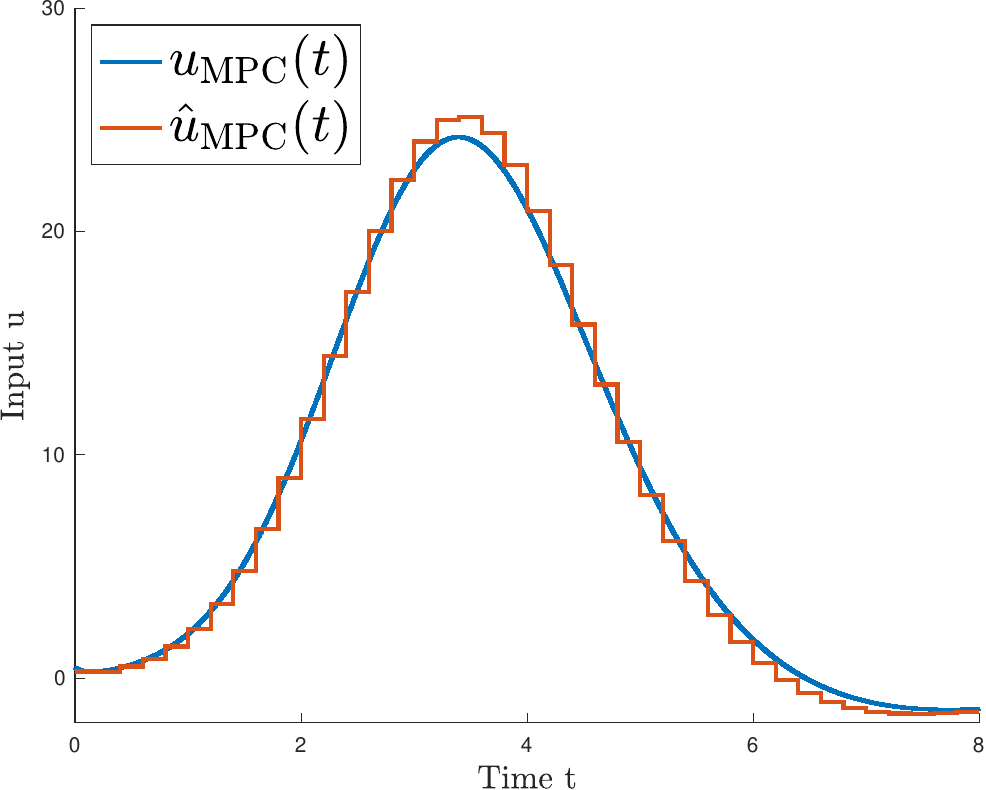}
        \caption{Control inputs.}
        \label{fig:sim:discreteFMPC:input}
    \end{subfigure}
    \caption{Simulation of system~\eqref{eq:TorsOscill-IO} under the control generated by~\Cref{Algo:DiscrFunnelMPC} with $\delta=\SampleTime = 0.002$ and $\delta=\SampleTime = 0.2$.\label{Fig:DiscreteFMPC}}
\end{figure}

To illustrate  the  theoretical results by a numerical example, we consider a
torsional oscillator with two flywheels, which are connected by a rod, see
\Cref{Fig:TorsOscill}. %
Such a system can be interpreted as a simple model of a driving train, cf.~\cite{Druecker22}.
The equations of motion for the torsional oscillator are given by
\begin{equation*} 
\begin{aligned}
    \begin{bmatrix} I_1 & 0 \\ 0 & I_2 \end{bmatrix}
    \begin{bmatrix} 
    \ddot z_1(t) \\ \ddot z_2(t)
    \end{bmatrix} = 
    \begin{bmatrix}
        -d & d \\ d & - d
    \end{bmatrix}
    \begin{bmatrix}
        \dot z_1(t) \\ \dot z_2(t)
            \end{bmatrix}
         + 
        \begin{bmatrix}
            -k & k \\ k & -k
        \end{bmatrix}
        \begin{bmatrix}
            z_1(t) \\ z_2(t)
        \end{bmatrix}
        + \begin{bmatrix}
            1 \\ 0
        \end{bmatrix} u(t),
\end{aligned}
\end{equation*}
where for $i=1,2$ (the index~$1$ refers to the lower flywheel) $z_i$ is the rotational position of the flywheel, $I_i > 0$ is the inertia, $d,k > 0$ are damping and torsional-spring constant, respectively.
We aim to control the oscillator such that the lower flywheel follows a given velocity profile. Hence, we choose 
$
    y(t) = \dot z_1(t)
$
as the output.
To remove the rigid-body motion from the dynamics, we introduce $\hat z:= z_1 - z_2$.
With this new variable, setting $x:=[\hat z, \dot z_1, \dot z_2]$ the dynamics can be written as
\begin{equation*}
    \begin{aligned}
        \dot x(t) &= A x(t) + B u(t),\\ 
        y(t) &= C x(t)  = \dot z_1(t),
    \end{aligned}
\end{equation*}
where
\begin{equation*} %
    \begin{aligned}
        \tilde A &\coloneqq \begin{bmatrix} 0&1&-1\\-k&-d&d\\k&d&-d \end{bmatrix},& 
        \tilde B &\coloneqq \begin{bmatrix}0\\1\\0\end{bmatrix},&
        M &\coloneqq \begin{bmatrix} 1&0&0\\ 0&I_1&0 \\ 0&0&I_2 \end{bmatrix},   \\
        A &\coloneqq M^{-1} \tilde A,&
        B &\coloneqq M^{-1} \tilde B,&  
        C&\coloneqq \begin{bmatrix} 0 & 1 & 0\end{bmatrix}.
    \end{aligned}
\end{equation*}
Using standard techniques, see e.g.~\cite{ilchmann1991non} and also
\Cref{Ex:LTISystem}, the reduced dynamics of the torsional oscillator can then
be written in Byrnes-Isidori form \eqref{eq:LTIByrnesIsidori}
\begin{equation} \label{eq:TorsOscill-IO}
    \begin{aligned}
        \dot y(t) &= R y(t) + S \eta(t) + \Gamma u(t), \\
        \dot \eta(t) &= Q \eta(t) + P y(t),
    \end{aligned}
\end{equation}
where~$\eta$ is the internal state, and 
\begin{equation*}
    \begin{aligned}
        R &= \frac{-d}{I_1},& 
        S &= \frac{1}{I_1} \begin{bmatrix} k & d \end{bmatrix},&
        Q &= \frac{1}{I_2} \begin{bmatrix} 0 & I_2 \\ -k & -d \end{bmatrix},&
        P &= \frac{1}{I_2} \begin{bmatrix} -I_2 \\d           \end{bmatrix}.
    \end{aligned}
\end{equation*}
Note that~$Q$ is a stable matrix, i.e. its eigenvalues are on the left half plane. 
Thus, the internal dynamics are bounded-input bounded-state stable.
The high-gain matrix is given by $\Gamma := CB = 1/I_1 > 0$.
For the purpose of simulation, we choose the reference
\begin{equation*}
    y_{\rf}(t) = \frac{250}{2} \left(1 + \frac{1}{\sqrt{2\pi}} \int_0^t \me^{-\frac{1}{2}(s-3)^2} \text{d}s \right),
\end{equation*}
which is a modified version of the error function ($\textsc{erf}$) and represents a smooth transition from zero rotation to an (approximately) constant angular velocity of~$250$ rotations per unit time. Thus, $\|y_{\rf}\|_\infty \le 250$, 
$\|\dot y_{\rf}\|_\infty = 250/\sqrt{2\pi}$.
Inserting the dimensionless parameters $I_1 = 0.136$, $I_2 = 0.12$, $k=0.1$, and
$d=0.16$, and invoking the reference~$y_{\rf}$ and the constant error
tolerance~$\psi \equiv 25$ (we allow a deviation of $10\%$), we may derive worst case
bounds on the system dynamics by estimating the explicit solution of the linear
equations~\eqref{eq:TorsOscill-IO}. We compute these bounds in order to estimate
a sufficiently large ${\umax\geq0}$ as in the proof of \Cref{Lem:DiscrUNonEmpty}.
For the sake of simplicity, we will assume $\eta(0) = 0$, which does not cause
loss of generality. For $\|y\|_\infty \le \|y_{\rf}\|_\infty + \psi$, we estimate
\begin{equation*}
\begin{aligned}
  \forall \, t \ge 0\, :  \ \   \| \eta(t) \| %
    & \le \frac{M}{\mu}\|P\| ( \|y_{\rf}\|_\infty + \psi) ,%
\end{aligned}
\end{equation*}
where $M := \sqrt{\|K^{-1}\|\|K\|}$ and $\mu := 1/(2\|K\|)$, and $K \in \R^{2
\times 2}$ solves the Lyapunov equation $KQ + Q^\top K + I_2 = 0$. Inserting the
values, we find that the estimates for step length of the control
signal~$\SampleTime$ and maximal control provided in the proof
of~\Cref{Lem:DiscrUNonEmpty} are satisfied with $\SampleTime = 0.002$, and
$\umax = 267$. We choose the time shift $\delta=\SampleTime$, i.e. a constant
control is applied to the system between two iterations of the sampled-data
funnel MPC~\Cref{Algo:DiscrFunnelMPC}. Further, the prediction horizon is set as~$T=10
\delta$. For the purpose of simulation, we use the non-strict funnel penalty function
\begin{equation*} 
\begin{aligned}
    \ell(t,y,u)=
    \begin{cases}
        \| y-y_{\rf}\|^2 + \lambda_u \|u\|^2,
            & \Norm{y-y_{\rf}} \leq  \psi(t)\\
        \infty, &\text{else},
    \end{cases}
\end{aligned}
\end{equation*}
with~$\lambda_u = 10^{-1}.$
The results are depicted in \Cref{Fig:DiscreteFMPC}.
While \Cref{Fig:DiscreteFMPC} displays the system's output evolving within the
funnel boundary, \Cref{fig:sim:discreteFMPC:input} shows the corresponding
control signals.

We stress that the estimates for $\SampleTime$ and $\umax$ in the proof of
\Cref{Lem:DiscrUNonEmpty} are very conservative. To demonstrate this aspect, we
run a second simulation, where we chose $\delta=\SampleTime = 0.2$, $T=1$, and
$\umax=30$. The results of this simulation are labelled as~$\hat y, \hat u$,
respectively. With this much larger uniform step length, the tracking objective
can be satisfied as well, cf.~\Cref{Fig:DiscreteFMPC}. Note that the maximal
applied control value is in both cases much smaller than the (conservative)
estimate~$\umax$ satisfying \Cref{Lem:DiscrUNonEmpty}. These simulations suggest
that the bounds derived \Cref{Lem:DiscrUNonEmpty} leave room for improvement.
As before, all simulations have been performed with \textsc{Matlab} using the
\textsc{CasADi} framework. The corresponding source code can be
found on \textsc{GitHub} under the link \url{https://github.com/ddennstaedt/FMPC_Simulation}.
\chapter{Outlook}
In this thesis, the concept of funnel model predictive control is presented, which integrates ideas 
from the adaptive high-gain control technique funnel control in a model predictive control scheme.
Building upon the framework  outlined in~\Cref{Chapter:FunnelMPC},
three extensions are subsequently introduced in~\Cref{Chapter:RobustFunnelMPC} through~\ref{Chapter:DiscretFMPC}.
The following section summarises the main results and provides a brief outlook on future research directions.\\

\noindent \textbf{Funnel model predictive control}  represents a
novel MPC approach to output tracking for a class of non-linear multi-input
multi-output systems governed by functional differential equations. 
By combining the predictive capabilities of MPC with concepts of the adaptive
funnel control technique, this framework guarantees prescribed transient
performance -- ensuring the tracking error remains within user-defined,
time-varying boundaries for smooth reference signals.
Central to its efficacy are \emph{funnel penalty functions}, which dynamically
penalise the error trajectory's distance to the funnel boundaries eliminating 
the need for conventional mechanisms such as
terminal conditions, artificially extended prediction horizons, or restrictive
output constraints to ensure initial and recursive feasibility.

A critical assumption underpinning funnel MPC is the availability of  sufficiently
large control values, quantified by $\umax\geq0$.
While \Cref{Th:ExUmax} establishes existence of such a bound, its
current formulation is inherently conservative and computationally
intractable -- limiting practical applicability. Addressing this, future research
should prioritise:
\begin{enumerate}
    \item \textbf{Refinement of estimates}: Existing bounds on $\umax\geq 0$, derived as
    worst-case guarantees independent of the prediction horizon $T>0$, likely
    obscure potential synergies between $T$ and the required control effort. A
    rigorous exploration of $T$'s role -- particularly in balancing transient
    performance against input magnitude -- could yield tighter,
    horizon-dependent bounds. 
    \item \textbf{Parametric sensitivity analysis}: A systematic characterisation
    of how auxiliary parameters (e.g. funnel shape, error variables $\eM_i$,
    weighting parameters $k_i$) influence feasibility and performance would
    enhance design flexibility. 
    \item \textbf{Fixed-input feasibility}: Developing mechanisms to ensure
    recursive feasibility under a priori fixed control limits $\umax\geq 0$ 
    remains a pivotal challenge for implementation.
    \item \textbf{Cost function simplification}: Investigating whether the
    weighted sum of the tracking error $\eMTrack=\yM-y_{\rf}$ and its
    derivatives in the funnel penalty function for higher order systems can be
    reduced to the sole error signal $\eMTrack$ -- while ensuring initial and
    recursive feasibility  provided $T>0$ is chosen large enough -- would
    simplify the algorithm's complexity.
    \item\textbf{Generalisation of model class}: 
    The presented results hold for models with a strict global relative degree. Since funnel control has been successfully generalised to systems with vector relative degree \cite{Hoang18,Berger2020Vector}, it is worth investigating a corresponding generalisation of the funnel MPC framework.
    \item \textbf{Numerical implementation}: 
    The incompatibility of discontinuous funnel penalties functions
    with standard optimisation frameworks (e.g.~\textsc{CasADi}) needs to be addressed.
    Future work should explore the development of smooth approximations or
    custom solvers tailored to funnel penalty functions in order to ensure fast
    numerical convergence while adhering to funnel boundaries and maintaining feasibility guarantees.
\end{enumerate}
Beyond these technical refinements, broader questions remain unanswered. A
comprehensive benchmarking study comparing funnel MPC against classical MPC variants
remains an open research question.
Furthermore, extending the developed principles to alternative control objectives, such
as safety-critical set invariance (e.g. confining states to prescribed safe
regions), presents fertile ground for further theoretical and
applied investigations.

\noindent \textbf{Robust funnel MPC} synergises funnel MPC and model-free
adaptive funnel control into a two-component architecture.
This hybrid scheme bridges the often-competing priorities of optimality and
robustness, achieving prescribed tracking performance even under structural
model-plant mismatches and unknown disturbances.
\begin{itemize}
    \item \textbf{Funnel MPC} prioritises optimality by minimising a designer-specified cost functional over receding horizons.
    \item \textbf{Funnel control} ensures robustness through adaptive disturbance rejection, activated only when necessary.
\end{itemize}
Key to their compatibility is the strategic design of the funnel controller’s
reference signal and boundary, derived from the MPC’s predictions. This ensures
the components complement rather than conflict. Further refinement can be achieved
via an activation function, which sparsely engages the funnel controller to
minimally perturb the optimal MPC signal while rejecting disturbances.

The framework periodically updates the model with system measurements via proper
initialisation. While theoretically generalisable, this process remains
cumbersome in practice, prompting the question: Can initialisation be
streamlined without compromising robustness? 
Future research will focus on extracting criteria to find explicit and
beneficial proper initialisation strategies.

Further open challenges and future directions include:
\begin{enumerate}
    \item \textbf{Unified model-system classes}: The model and system currently
    require distinct classes of differential equation. While the model is
    assumed to have a control affine representation, the function $F$ describing
    the system dynamics has the perturbation high-gain property. A unification
    of these two classes would broaden applicability.
    \item \textbf{Explicit combined input bounds}:
    While the MPC component’s control input is bounded by $\umax\geq 0$, the
    model-free funnel controller lacks explicit a-priori bounds. Deriving a composite
    bound for the combined scheme is critical for safety-critical applications.
    \item \textbf{Derivative-free operation}: The funnel controller’s reliance
    on output derivatives poses practical challenges with noisy measurements. 
    Integrating a \emph{funnel pre-compensator}~\cite{BergReis18,lanza2022output} -- 
    to estimate derivatives or bypass their need -- warrants exploration.
    \item \textbf{Order flexibility}: The proposed framework mandates matching
    relative degrees for model and system.  Relaxing this constraint could
    enable simplified models (e.g. lower-order approximations) for complex
    systems.
\end{enumerate}

\noindent \textbf{Learning-based robust funnel MPC} extends the robust funnel
MPC framework by integrating a versatile online learning architecture. This
approach continuously refines the surrogate model using historical data -- system
outputs, model predictions, and applied control signals -- drawn from both the
model-based funnel MPC and the model-free feedback component. It ensures robust
tracking within predefined (time-varying) performance boundaries while
accommodating:
\begin{itemize}
    \item \textbf{Varying model complexity}, from simplified approximations to high-fidelity representations.
    \item \textbf{Continual improvement} via iterative data assimilation.
    \item \textbf{Methodological agnosticism}, allowing integration of diverse learning paradigms.
\end{itemize}
By combining learning techniques with both model-based prediction and adaptive
control, this framework bridges the gap between robustness and adaptability in
uncertain environments.

While the current formulation is abstract and theoretical, future research will address critical open questions:
\begin{enumerate}
    \item\textbf{Learning scheme efficacy}: What defines an effective learning
    scheme? How can controller performance improvement be rigorously verified?
    \item\textbf{Technique compatibility}: Which established methods -- Willems’
    fundamental lemma, Koopman operator theory, or neural networks -- can
    effectively be used to leverage the collected data? 
    \item\textbf{Feasibility guarantees}: How can feasibility be rigorously proven for advanced learning algorithms?
    \item\textbf{Prior knowledge integration}: How should existing system knowledge inform the learning architecture?
\end{enumerate}

\noindent \textbf{Sampled-data robust funnel MPC} demonstrates how output
tracking with prescribed performance can be achieved while restricting
admissible controls to piecewise constant step functions. The key contribution
is \emph{explicit uniform bounds} on sampling rates and maximal control
effort for both the funnel MPC and model-free funnel controller.
This is an important step to bridge the gap between continuous‐time theory and
real‐world sampled‐data implementations. 
For the funnel controller, we further showed that its Zero-order-Hold 
implementation can serve as a safety filter for learning-based control
architectures. However, effective deployment requires addressing the reliance on
noise-sensitive output derivative measurements -- a critical challenge for future
work.
While foundational, the derived bounds remain highly conservative. Relaxing
these estimates is essential for practical applicability. Additionally, the
current system and model classes (tailored for sampled-data control) represent
subsets of those in prior chapters. Generalising these results to broader
classes of systems/models remains an open problem.

To advance digital implementation, three key questions arise:
\begin{itemize}
    \item Can the continuous-time cost function (currently integral-based) used
    in the funnel MPC algorithm be efficiently discretised with uniform error
    bounds? 
    \item  Can the algorithms be redeveloped entirely for discrete-time systems,
    bypassing continuous-time computations?
    \item How might a discrete-time theory for funnel control and funnel MPC be formulated?
\end{itemize}
Presently, all theoretical guarantees assume continuous-time dynamics. A
discrete-time counterpart -- for both components -- remains unexplored.
Furthermore, learning techniques specifically tailored for
sampled-data systems -- such as those leveraging intermittent measurements or
quantised data -- could prove particularly advantageous in enhancing
adaptability while preserving robustness. 
The integration of such methods also promises to be an  interesting direction for future research.

\chapter*{Appendix}\setcurrentname{Appendix}
\label{Chapter:Appendix}
\addcontentsline{toc}{chapter}{Appendix}
\setcounter{chapter}{7}
The existence of solutions of the differential equations  is essential for  both
the system~\eqref{eq:Sys} and model~\eqref{eq:Model_r}. From an
application point of view, the question of the solution's existence is often not
of interest or merely seen as a technical detail. However, it is of utmost importance
mathematically as the foundation of all further investigations and results.
Although several works, see
e.g.~\cite{ryan2001controlled,IlchRyan02a,IlchRyan02b,Ilchmann01102009}, have
already provided answers to this question for systems similar to the ones
considered in this thesis, we would like to provide a rigorous proof in this
work as well for the sake of completeness.

To this end, we consider the initial value problem
\begin{equation}\label{eq:SysAppendix}
\begin{aligned}
    \dot{y}(t)&=F(t,y(t),\oT(y)(t)),\\
    y|_{[0,t_0]}&=y^0\in\cC([0,t_0],\R^n),
\end{aligned}
\end{equation}
and will prove the existence of solutions for this initial value problem in the
following. For the sake of generality, we want to analyse the problem with the
function $F$  being only defined on a domain, i.e. a non-empty connected relatively open
set, but not necessarily on the whole space. Thus, let $\cD_1\subset
\Rp\times\R^n$ be non-empty, connected, relatively open sets with $(t_0,y^0(t_0))\in\cD_1$.
Assume that $F:\cD_1\times \R^q\to\R^n$ is a \emph{Carath\'{e}odory function},
i.e. it has the following properties for every compact interval $I$,
$x^0\in\R^n$ and $\eps>0$ with $I\times\bar{\cB}_{\eps}(x^0)\subset \cD_1$ and
every compact set~$K\subset\R^q$:
\begin{enumerate}
    \item[\Itemlabel{Item:AppendixCaratheodoryProp1}{(C.1)}] $F(t,\cdot,\cdot):\bar{\cB}_{\eps}(x^0)\times K\to\R^n$ is continuous for almost all $t\in I$,
    \item[\Itemlabel{Item:AppendixCaratheodoryProp2}{(C.2)}] $F(\cdot,x,z):I\to\R^n$ is measurable for all fixed $(x,z)\in \bar{\cB}(x^0)\times K$,
    \item[\Itemlabel{Item:AppendixCaratheodoryProp3}{(C.3)}]  there exists an integrable function $\kappa:I\to\Rp$ such that 
    $\Norm{F(t,x,z)}\leq \kappa(t)$ for almost all $t\in I$ and all $(x,z)\in \bar{\cB}_{\eps}(x^0)\times K$.
\end{enumerate}
This notion of a {Carath\'{e}odory function} is based on the definition in~\cite[Appendix B]{Ilchmann01102009}.

We assume that the operator $\oT$ only acts on functions evolving within a domain ${\cD_2\subset\Rp\times\R^n}$.
To formally define the properties of $\oT$, we impose on $\cD_2$ that, for all ${t\in\Rp}$, there exists $x\in\R^n$ with 
$(t,x)\in\cD_2$ and define the set of regulated functions evolving in $\cD_2$ as  
$\cR{\cD_2}\coloneqq \setdef{y\in\cR(\Rp,\R^n)}{\graph(y)\subset\cD_2 }$.
Then, we assume the operator $\oT: \cR{\cD_2} \to L^\infty_{\loc} ([t_0,\infty), \R^{q})$
to fulfil the following properties:
\begin{enumerate}
    \item[\Itemlabel{Item:AppendixOperatorPropCasuality}{(T.1')}]$\fa y_1,y_2\in\cR{\cD_2}$ $\fa t\geq t_0$:
    \[
        y_1\vert_{[0,t]} = y_2\vert_{[0,t]}
        \ \Impl\ 
            \oT(y_1)\vert_{[t_0,t]}=\oT(y_2)\vert_{[t_0,t]}.
    \]
    \item[\Itemlabel{Item:AppendixOperatorPropLipschitz}{(T.2')}]
    $\fa t \ge t_0 $ $\fa y \in \cR([0,t] ; \R^n)$ with $\graph(y)\subset\cD_2$
    $\ex \Delta, \delta, c > 0$ 
    $\fa y_1, y_2 \in \cR{\cD_2}$ with
    $y_1|_{[0,t]} = y_2|_{[0,t]} = y $ 
    and $\Norm{y_1(s) - y(t)} < \delta$,  $\Norm{y_2(s) - y(t)} < \delta $ for all $s \in [t,t+\Delta]$:
    \[
     \esssup_{\mathclap{s \in [t,t+\Delta]}}  \Norm{\oT(y_1)(s) - \oT(y_2)(s) }  
        \le c \ \sup_{\mathclap{s \in [t,t+\Delta]}}\ \Norm{y_1(s)- y_2(s)}.
    \] 
    \item[\Itemlabel{Item:AppendixOperatorPropBIBO}{(T.3')}]
    For every $t\geq t_0$ and every family $(K_s)_{s\in[0,t]}$ of compact sets $K_s\subset\R^n$ such that 
    $\bigcup_{s\in [0,t]} K_{s}$ is a bounded set and $\bigcup_{s\in [0,t]}\cbl s\cbr\times K_s\subset\cD_2$,
    there exists $c_1>0$  such that for all $y \in \cR\cD_2$:
    \[
    \graph(y)\subset \bigcup_{s\in [0,t]}\cbl s\cbr\times K_s
    \Impl \ \sup_{s\in[t_0,t]} \Norm{\oT(y)(s)}  \le c_1,
    \]
\end{enumerate}
The properties \ref{Item:AppendixOperatorPropCasuality},
\ref{Item:AppendixOperatorPropLipschitz}, and \ref{Item:AppendixOperatorPropBIBO}
adapt \ref{Item:OperatorPropCasuality},
\ref{Item:OperatorPropLipschitz}, and \ref{Item:OperatorPropBIBO} from
\Cref{Def:OperatorClass} to accommodate the restriction of $\oT$ to the set of functions
with codomain $\cD_2$.
While defining the operator $\oT$ only on a set of functions restricted to a domain
for a system of the form~\eqref{eq:SysAppendix} was already considered
in~\cite{Hachmeister23}, a proof for the existence of solutions was omitted. 
Since the usage of a such modified operator could be of interest in application
and future research work, we want to provide a proof in the following.
It is clear that modified properties \ref{Item:AppendixOperatorPropCasuality} and 
\ref{Item:AppendixOperatorPropLipschitz}
are equivalent to their original counterparts in the case $\cD_2=\Rp\times \R^n$.
This is not the case for \ref{Item:AppendixOperatorPropBIBO}. 
It is a weaker assumption on~$\oT$ as the following lemma shows. 
\begin{lemma}
    Let $\cD_2=\Rp\times \R^n$ and 
    $\oT: \cR{\cD_2} \to L^\infty_{\loc} ([t_0,\infty), \R^{q})$
    be an operator with property~\ref{Item:AppendixOperatorPropCasuality}.
    If $\oT$ has property \ref{Item:OperatorPropBIBO}, then it also satisfies property 
    \ref{Item:AppendixOperatorPropBIBO}. The opposite is in general not true.
\end{lemma}
\begin{proof}
    Let $\oT: \cR(\Rp,\R^n) \to L^\infty_{\loc} ([t_0,\infty), \R^{q})$ satisfying
    \ref{Item:AppendixOperatorPropCasuality} and \ref{Item:OperatorPropBIBO}. We show $\oT$ has property~\ref{Item:AppendixOperatorPropBIBO}.
    Let $t\geq t_0$ and $(K_s)_{s\in[0,t]}$ be a family of compact sets 
    $K_s\subset\R^n$ with $\bigcup_{s\in [0,t]} K_{s}$ being bounded.
    There exists $c_0>0$ with $\bigcup_{s\in [0,t]} K_{s}\subset\bar{\cB}_{c_0}$. 
    Due to property~\ref{Item:OperatorPropBIBO}, there exists $c_1>0$ with 
    such that
    \[
        \sup_{t \in [t_0,\infty)} \Norm{\oT(y)(t)}  \le c_1
    \]
    for all $y \in \cR(\Rp, \R^n)$ with $\sup_{t \in \Rp} \Norm{y(t)} \le c_0$.
    Let $y \in \cR(\Rp, \R^n)$ be a function with ${\graph(y)\subset \bigcup_{s\in [0,t]}\cbl s\cbr\times K_s}$.
    Define $\tilde{y}\in\cR(\Rp, \R^n)$ by 
    $\tilde{y}(s)=y(s)$ for $s\in [0,t]$ and $\tilde{y}(s)=y(t)$ for $s>t$.
    Due to the causality property~\ref{Item:AppendixOperatorPropCasuality}, we have
    \[
        \sup_{s\in[t_0,t]} \Norm{\oT(y)(s)}=
        \sup_{s\in[t_0,t]} \Norm{\oT(\tilde{y})(s)}  \leq
        \sup_{s\in[t_0,\infty)} \Norm{\oT(\tilde{y})(s)}  \le c_1.
    \]
    This shows that $\oT$ fulfils \ref{Item:AppendixOperatorPropBIBO}.

    To show that the opposite is in general not true, we consider a counter example.
    Define ${\oT: \cR(\Rp,\R) \to L^\infty_{\loc} ([t_0,\infty), \R)}$ by
    $\oT(y)(t)\coloneqq \int_{0}^ty(\tau)\d{\tau}$ for $t\in [t_0,\infty)$.
    It is clear that $\oT$ is a causal, i.e. it fulfils property~\ref{Item:AppendixOperatorPropCasuality}.
    We show that $\oT$ has also property~\ref{Item:AppendixOperatorPropBIBO}.
    To this end, let $t\geq t_0$ and $(K_s)_{s\in[0,t]}$ be a family of compact sets 
    $K_s\subset\R^n$ with $\bigcup_{s\in [0,t]} K_{s}$ being bounded.
    There exists $c_0>0$ with $\bigcup_{s\in [0,t]} K_{s}\subset\bar{\cB}_{c_0}$. 
    Let  $y \in \cR\cD_2$ with ${\graph(y)\subset \bigcup_{s\in [0,t]}\cbl s\cbr\times K_s}$ and 
    set $c_1\coloneqq tc_0$.
    Then, we have
    \[
        \sup_{s\in[t_0,t]} \Norm{\oT(y)(s)}=
        \sup_{s\in[t_0,t]} \Norm{\int_{0}^sy(\tau)\d{\tau}}\leq
        \int_{0}^t\SNorm{y}\d{\tau}\leq tc_0=c_1.
    \]
    Thus, $\oT$ has property~\ref{Item:AppendixOperatorPropBIBO}.
    However, for the constant function $\tilde{y}\equiv1$, we have $\oT(y)(t)\to\infty$ for $t\to\infty$.
    Therefore, $\oT$ does not have the bounded-input bounded-output property~\ref{Item:OperatorPropBIBO}.
\end{proof}

With the assumed properties of $F$ and $\oT$ at hand, we define a solution of the initial value problem~\eqref{eq:SysAppendix}
in the virtue of~\cite[Section 5]{IlchRyan02b} as follows.

\begin{definition}\label{Appendix:Def:Solution}
    For $y^0\in\cC([0,t_0],\R^n)$ with $(t_0,y^0(t_0))\in\cD_1$ and $\graph(y^0)\subset\cD_2$, 
    a function $y:[0,\omega)\to\R^n$ with $\omega\in(t_0,\infty]$ 
    and $[t_0,\omega)\subset I$  is called a \emph{solution} of the 
    initial value problem~\eqref{eq:SysAppendix}, if $y|_{[0,t_0]}=y^0$
    and 
    \[
        \fa t\in[t_0,\omega):\quad y(t)=y^0(t_0)+\int_{t_0}^t F(s,y(s),\oT(y)(s))\d{s}.
    \]
    A solution $y:[0,\omega)\to\R^n$ is said to be \emph{maximal} if it has no proper right extension that is also a solution.
\end{definition}
    Note that we identify $\cC([0,t_0],\R^n)$ with $\R^n$ if $t_0=0$ in \Cref{Appendix:Def:Solution}.  
    Moreover, given a function $y\in\cR([0,\omega),\R^n)$ with $\omega<\infty$ and $\graph(y)\subset\cD_2$,
    $\oT(y)(t)$ is interpreted for $t\in[0,\omega)$ as the evaluation of $\oT(y^e)(t)$ 
    for an arbitrary right extension $y^e\in\cR\cD_2$  of $y$, as elaborated in~\Cref{Rem:PropertiesOperator}~\ref{Item:Rem:PropertiesOperator:Causality}.
    
\begin{remark}
    Although it was not mentioned explicitly in \Cref{Appendix:Def:Solution}, a solution 
    ${y:[0,\omega)\to\R^n}$ of the initial value problem~\eqref{eq:SysAppendix} has the following properties: 
\begin{enumerate}[(i)]
    \item $y_{[t_0,\omega)}$ is absolutely continuous,
    \item $(t,y(t))\in\cD_1\cap\cD_2$ for all $t\in[t_0,\omega)$.
\end{enumerate}
\end{remark}
With the definition of solutions of the initial value problem~\eqref{eq:SysAppendix} established, we now present a key existence theorem.
\begin{theorem}\label{Appendix:Th:SolutionExists}
    Consider the initial value problem~\eqref{eq:SysAppendix} where $F:\cD_1\times\R^q\to\R^n$ 
    is a Carath\'{e}odory function and the operator $\oT: \cR\cD_2 \to L^\infty_{\loc} ([t_0,\infty), \R^{q})$
    has the properties \ref{Item:AppendixOperatorPropCasuality}, \ref{Item:AppendixOperatorPropLipschitz}, and \ref{Item:AppendixOperatorPropBIBO}.
    Let $y^0\in\cC([0,t_0],\R^n)$ with $y^0(t_0)\in \cD_1$ and $\graph(y^0)\subset\cD_2$.
    Then,
    \begin{enumerate}[(i)]
        \item\label{Item:Th:ExSolution} the initial value problem~\eqref{eq:SysAppendix} has a solution $y:[0,\omega)\to\R^n$ with $\omega>t_0$ in the sense of~\Cref{Appendix:Def:Solution},
        \item\label{Item:Th:ExSolutionMaximal} every solution can be extended to a maximal solution,
        \item\label{Item:Th:ExSolutionBounded} if $F$ is locally essentially bounded and $y\in\cC([0,\omega),\R^n)$ is a maximal solution, 
        then the closure of $\graph(y|_{[t_0,\omega)})$ is not a compact subset of $\cD_1\cap\cD_2$.
    \end{enumerate}
\end{theorem}
\begin{proof}
    We adapt the proof of~\cite[Theorem 7.1]{Ilchmann01102009} to the current setting.
    
    \noindent
    \emph{Step 1}:
    We have $y^0\in\cC([0,t_0],\R^n)$ with $\graph(y^0)\subset\cD_2$.
    Thus, using property~\ref{Item:AppendixOperatorPropLipschitz} of operator~$\oT$, there exist 
    $ \Delta, \delta, c > 0$ such that, for all
    $ y_1, y_2 \in \cR\cD_2$ with
    $y_1|_{[0,t_0]} =y^0= y_2|_{[0,t_0]}$ 
    and $\Norm{y_1(s) - y^0(t_0)} < \delta$,  $\Norm{y_2(s) - y^0(t_0)} < \delta $ ,
    we have for all $s \in [t_0,t_0+\Delta]$:
    \[
     \esssup_{\mathclap{s \in [t_0,t_0+\Delta]}}  \Norm{\oT(y_1)(s) - \oT(y_2)(s) }  
        \le c \ \sup_{\mathclap{s \in [t_0,t_0+\Delta]}}\ \Norm{y_1(s)- y_2(s)}.
    \]
    Both $\Delta>0$ and $\delta>0$ can be chosen sufficiently small such that 
    \[
        [t_0,t_0+\Delta]\times \bar{\cB}_{\delta}(y^0(t_0))\subset \cD_1\cap\cD_2.
    \]
    For $t\in[0,t_0+\Delta]$ define the compact set 
    \[
        K_t\coloneqq \begin{cases}
        \cbl y^0(t)\cbr,&t\in[0,t_0),\\
        \bar{\cB}_{\delta}(y^0(t_0)),&t\in[t_0,t_0+\Delta].\\
        \end{cases}
    \]
    Then, $\bigcup_{t\in[0,t_0+\Delta]} K_t$ is a bounded set and $\bigcup_{t\in [0,t_0+\Delta]}\cbl t\cbr\times K_t\subset\cD_2$.
    By property~\ref{Item:AppendixOperatorPropBIBO} of operator $\oT$,  there
    exists $c_1>0$ such that, for $y\in\cR\cD_2$
    with  $\graph(y)\subset \bigcup_{t\in [0,t_0+\Delta]}\cbl t\cbr\times K_t$, we have
    $\Norm{\oT(y)(t)}<c_1$ for all $t\in [t_0,t_0+\Delta]$.  
    Note that, for every right extension $y^e\in\cR([0,t_0+s],\R^n)$ of $y^0$,
    $s\in[0,\Delta]$, with  $y^e([t_0,t_0+s])\subset
    \bar{\cB}_{\delta}(y^0(t_0))$, there exists a function
    $\hat{y}^e\in\cR\cD_2$ with $\hat{y}^e|_{[0,t_0+s]}=y^e$ and
    $\graph(\hat{y}^e|_{[0,t_0+\Delta]})\subset \bigcup_{t\in [0,t_0+\Delta]}\cbl t\cbr\times K_t$.
    Thus,  
    \begin{equation}\label{eq:AppendixObservationBounded}
        \fa t\in [t_0,t_0+s]:\quad \Norm{\oT(y^e)(t)}=\Norm{\oT(\hat{y}^e)(t)}\leq\sup_{t\in[t_0,t_0+\Delta]}\Norm{\oT(\hat{y}^e)(t)}<c_1
    \end{equation}
    because of the causality property~\ref{Item:AppendixOperatorPropCasuality} of operator~$\oT$.
    We will use this observation later.
    As $F$ is a {Carath\'{e}odory function}, 
    property~\ref{Item:AppendixCaratheodoryProp3} yields the existence of an integrable function $\kappa:[t_0,t_0+\Delta]\to\Rp$ with
    \[
        \fa (t,x,z)\in[t_0,t_0+\Delta]\times \bar{\cB}_{\delta}(y^0(t_0))\times\bar{\cB}_{c_1}:\quad \Norm{F(t,x,z)}\leq \kappa(t).
    \]
    Define $\gamma:[0,t_0+\Delta]\to\Rp$ by
    \[
        \gamma(t)\coloneqq \begin{cases}
        0,&t\in[0,t_0),\\
        \int_{t_0}^t\kappa(s)\d{s},&t\in[t_0, t_0+\Delta].
        \end{cases}
    \]
    There exists $\tau>0$ such that $\gamma(t_0+\tau)<\delta$.
    We define a sequence $(y_n)\in\cC([0,t_0+\Delta],\R^n)^\N$ as follows
    \[
        y_n(t)=\begin{cases}
        y^0(t),& t\in[0,t_0],\\
        y^0(t_0),& t\in(t_0,t_0+\tau/n],\\
        y^0(t_0)+\int_{t_0}^{t-\tau/n}F(s,y_n(s),\oT(y_n)(s))\d{s},& t\in(t_0+\tau/n,t_0+\tau].\\
        \end{cases}
    \]
    By construction, $y_n$ is a right extension of $y^0$ with $y_n(t)\in\bar{\cB}_{\delta}(y^0(t_0))$
    for all $t\in [t_0, t_0+\tau/n]$ and $n\in\N$. Thus, $\Norm{\oT(y_n)(t)}<c_1$
    for all $t\in [t_0, t_0+\tau/n]$ and $n\in\N$ because of the observation made in
    \eqref{eq:AppendixObservationBounded}.
    Therefore,
    \begin{equation}\label{eq:AppendixEstimateGamma}
        \Norm{y_n(t)-y^0(t_0)}\leq \int_{t_0}^{t-\tau/n}\Norm{F(s,y_n(s),\oT(y_n)(s))}\d{s}\leq \int_{t_0}^{t-\tau/n}\!\!\kappa(s)\d{s}=\gamma(t-\tau/n)
    \end{equation}
    for $t\in[t_0+\tau/n,t_0+\tau]$ and $n\in\N$. 
    Since $\gamma(t-\tau/n)<\delta$, 
    we infer $y_n(t)\in \bar{\cB}_{\delta}(y^0(t_0))$ and $\Norm{\oT(y_n)(t)}<c_1$ for all $t\in[t_0,t_0+\tau]$.
    Thus, $\Norm{F(t,y_n(t),\oT(y_n)(t)}\leq \kappa(t)$ for all $t\in[t_0,t_0+\tau]$ and all $n\in\N$.
    
    We will prove that the sequence $(y_n)_{n\in\N}$ is equicontinuous.
    To this end, let $\eps>0$ be arbitrary but fixed. 
    The function $\gamma$ is uniformly continuous on the compact interval $[t_0,t_0+\tau]$.
    Thus, there exists $\bar{\delta}>0$ such that 
    \[
        \Abs{\gamma(t)-\gamma(s)}<\eps
    \]
    for all $t,s\in[t_0,t_0+\tau]$ with $\Abs{t-s}<\bar{\delta}$.
    Let $n\in\N$ and $t,s\in[t_0,t_0+\tau]$ with $\Abs{t-s}<\bar{\delta}$. 
    We assume $s\leq t$ without loss of generality and consider three cases.
    First, if ${t_0\leq s\leq t\leq t_0+\tau/n}$, then $y_n(s)=y_n(t)=y^0(t_0)$. 
    Thus, $ \Norm{y_n(s)-y_n(t)}=0$.
    Second, if $t_0\leq s\leq t_0+\tau/n\leq t$, then   
    \[
       \Norm{y_n(t)- y_n(s)}=\Norm{y_n(t)- y^0(t_0)}=\gamma(t-\tau/n)<\eps,
    \]
    where estimate~\eqref{eq:AppendixEstimateGamma} was used.
    Third, if $t_0+\tau/n\leq s\leq t$, then
    \[ 
        \Norm{y_n(t)- y_n(s)}\leq \Abs{\gamma(t-\tau/n)-\gamma(s-\tau/n)}<\eps.
    \]
    As $y_n|_{[0,t_0]}=y^0$ for all $n$, the sequence $(y_n)_{n\in\N}$ is therefore equicontinuous.
    By the Arzel\`{a}-Ascoli theorem,
    there exists a function $y\in \cC([0,t_0+\tau],\R^{n})$ and a 
    subsequence (which we do not relabel) such that $y_n$ is uniformly convergent, 
    i.e.  $y_{n}\to y$. 
    Clearly, $y|_{[0,t_0]}=y^0$ and $y([t_0, t_0+\tau])\subset \bar{\cB}_{\delta}(y^0(t_0))\subset\cD_1\cap\cD_2$ 
    since $\bar{\cB}_{\delta}(y^0(t_0))$ is compact and ${y_n([t_0, t_0+\tau])\subset \bar{\cB}_{\delta}(y^0(t_0))}$ for all $n\in\N$.

    Since the $\oT$ is local Lipschitz continuous, see property~\ref{Item:AppendixOperatorPropLipschitz}, 
    ${\lim\limits_{n\to\infty} \oT(y_n)(t)=\oT(y)(t)}$ for almost all $t\in[0, t_0+\tau]$. 
    Thus, 
    \[
        \lim_{n\to\infty} F(t, y_n(t),\oT(y_n)(t))=F(t, y(t),\oT(y)(t))
    \]
    for almost all $t\in[0, t_0+\tau]$ since $F(t,\cdot,\cdot)$ is continuous, according to property~\ref{Item:AppendixOperatorPropCasuality}. 
    As $\Norm{F(t,y_n(t),\oT(y_n)(t)}<\kappa(t)$ for all $t\in [t_0,t_0+\tau]$, the Lebesgue dominated convergence theorem yields
    \[
        \lim_{n\to\infty}\int_{t_0}^{t}F(s, y_n(s),\oT(y_n)(s))\d{s}=\int_{t_0}^{t}F(s, y(s),\oT(y)(s))\d{s}
    \]
    for all $t\in[t_0,t_0+\tau]$.
    Note that 
    \begin{align*}
        y_n(t)&=y^0(t_0)+\int_{t_0}^{t-\tau/n}F(s, y_n(s),\oT(y_n)(s))\d{s}\\
        &=y^0(t_0)+\int_{t_0}^{t}F(s, y_n(s),\oT(y_n)(s))\d{s}-\int_{t-\tau/n}^{t}F(s, y_n(s),\oT(y_n)(s))\d{s}
    \end{align*}
    for all $t\in (t_0+\tau/n,t_0+\tau]$ and $n\in\N$. 
    For the limit $n\to \infty$, we conclude 
    \[
        y(t)=\begin{cases} 
            y^0(t),& t\in[0,t_0],\\
            y^0(t_0)+\int_{t_0}^{t}F(s, y(s),\oT(y)(s))\d{s},& t\in[t_0,t_0+\tau].
        \end{cases}
    \]
    Therefore, $y$ is a solution of the initial value problem~\eqref{eq:SysAppendix} in the sense of~\Cref{Appendix:Def:Solution}
    proving assertion~\ref{Item:Th:ExSolution}.
    
    \noindent
    \emph{Step 2}:
    We prove assertion~\ref{Item:Th:ExSolutionMaximal}.
    Let $y:[0,\omega)\to\R^n$ be a solution of the initial value problem~\eqref{eq:SysAppendix}.
    Define the set
    \[
        \cE\coloneqq \setdef{(\hat{\omega},\zeta)}{ \hat{\omega}\geq\omega, \zeta\in\cC([0,\hat{\omega}),\R^n)\text{ is a solution of~\eqref{eq:SysAppendix} with } \zeta|_{[0,\omega)}=y}.
    \]
    This is basically the set of all right extensions of $y$ that are also a solution of~\eqref{eq:SysAppendix}.
    As $(\omega,y)\in\cE$, this set is non-empty.
    The relation $\preceq$ given by
    \[
        (\omega_1,y_1)\preceq(\omega_2,y_2)\quad\Longleftrightarrow\quad \omega_1\leq\omega_2 \wedge y_1=y_2|_{[0,\omega_1)}
    \]
    defines a partial order on $\cE$.
    Let $\Omega$ be a chain in $\cE$, i.e. a totally ordered subset of $\cE$.
    Define $\omega^\star\coloneqq \sup\setdef{\hat{\omega}}{(\hat{\omega},\zeta)\in\Omega}$.
    Further, define $y^\star\in\cC([0,\omega^\star),\R^n)$ by $y^\star|_{[0,\hat{\omega})}=\zeta$ for $(\hat{\omega},\zeta)\in\cP$.
    Then, $(\omega^\star,y^\star)\in\Omega$ and $(\hat{\omega},\zeta)\preceq(\omega^\star,y^\star)$ for all $(\hat{\omega},\zeta)\in\cP$, i.e.
    $(\omega^\star,y^\star)$ is an upper bound of $\Omega$. 
    Zorn's lemma yields the existence of an maximal element of $\cE$.  
    By the construction of $\cE$ this is a maximal extension of $y$ that is also a solution.
    This proves~\ref{Item:Th:ExSolutionMaximal}.
    
    \noindent
    \emph{Step 3}:
    We prove assertion~\ref{Item:Th:ExSolutionBounded}.
    Assume that $F$ is locally essentially bounded and let $y\in\cC([0,\omega),\R^n)$ be a maximal solution.
    Seeking a contradiction, suppose that the closure of $\graph(y|_{[t_0,\omega)})$ is a compact subset of $\cD_1\cap\cD_2$.
    This implies, in particular, that $[t_0,\omega)$ is a bounded interval.
    As $y$ is bounded, property~\ref{Item:AppendixOperatorPropBIBO} implies the boundedness of $\oT(y)|_{[t_0,\omega)}$. 
    The local essential boundedness of $F$ yields the existence of $c_2>0$ with
    $
        \Norm{\dot{y}(t)}=\Norm{F(t,y(t),\oT(y)(t)}\leq c_2
    $
    for all $t\in [t_0,\omega)$. 
    Hence, $y$ is uniformly continuous on the interval $[t_0,\omega)$.
    There thus exists a right extension $y^e\in\cC([0,\omega],\R^n)$ 
    of $y$ with $\graph(y^e|_{[t_0,\omega]})\subset\cD_1\cap\cD_2$.
    In particular, $(\omega,y^e(\omega))\in \cD_1\cap\cD_2$.
    Assertion~\ref{Item:Th:ExSolution} yields the existence of a
    solution $\hat{y}:[0,\hat{\omega})\to\R^n$ with $\hat{\omega}>\omega$ of the
    initial value problem
    \begin{align*}
       \dot{y}(t)&=F(t,y(t),\oT(y)(t)),\\
        y|_{[0,\omega]}&=y^e.
    \end{align*}
    As $\hat{y}|_{[0,\omega)}=y$,  the function $\hat{y}$ is a proper extension of $y$ 
    and also a solution of the initial value problem~\eqref{eq:SysAppendix}. 
    This contradicts the maximality of $y$ and completes the proof.
\end{proof}

By reducing a higher order system of the form 
\begin{equation}\label{eq:SysHighOrderAppendix}
\begin{aligned}
    y^{(r)}(t)&=F(t,\OpChi(y)(t),\oT(\OpChi(y))(t)),\\
    y|_{[0,t_0]}&=y^0\in\cC([0,t_0],\R^n),
\end{aligned}
\end{equation}
with $r>1$ to a system of order one,
\Cref{Appendix:Th:SolutionExists} can clearly be also applied to such systems. 
Note that we now identify $\cC([0,t_0],\R^n)$ with $\R^{rn}$ in the case $t_0=0$.
In the spirit of~\Cref{Def:FCSolutionConcept},
a solution of the initial value problem~\eqref{eq:SysHighOrderAppendix} is an
absolutely continuous function $x=(x_1,\ldots,x_r):[0,\omega)\to\R^{rm}$  with $\omega\in(t_0,\infty]$
fulfilling 
\begin{equation}\label{Appendix:eq:DefSolutionSystem}
\begin{aligned}
    \dot{x}_{i}(t)&=x_{i+1}(t),\hspace{3cm} i=1,\ldots, r-1,\\
    \dot{x}_r(t)&= F(t,\OpChi(x)(t),\oT(x)(t)),
\end{aligned}
\end{equation}
for almost all $t\in[t_0,\omega)$ and $x|_{[0,t_0]} =\OpChi(y^0)$ (resp. $x(t_0)= y^0$ in the case $t_0=0$).

\Cref{Appendix:Th:SolutionExists} yields as a straightforward corollary the
existence of solutions of the initial value problem~\eqref{eq:Sys} for the
considered system class if a control $u\in L^\infty_{\loc}([t_0,\infty), \R^m)$
is applied.

\begin{corollary}\label{Appendix:Cor:SystemSolutionExists}
    Consider system~\eqref{eq:Sys} with $(F,\oT)\in\cN^{m,r}_{t_0}$  at initial time $t_0\geq0$.
    Let ${y^0\in \cC^{r-1}([0,t_0],\R^m)}$ be an initial trajectory
    and $u\in L^\infty_{\loc}([t_0,\infty), \R^m)$ be a control function.
    Then, 
\begin{enumerate}[(i)]
    \item the initial value problem~\eqref{eq:Sys} has a solution $x:[0,\omega)\to \R^{rm}$ in the sense of~\Cref{Def:FCSolutionConcept},
    \item every solution can be extended to a maximal solution,
    \item if $x:[0,\omega)\to \R^{rm}$ is a bounded maximal solution, then $\omega=\infty$.
\end{enumerate}
\end{corollary}
\begin{proof}
    The assertions follow directly from~\Cref{Appendix:Th:SolutionExists} since 
    $\oT$ has the properties \ref{Item:AppendixOperatorPropCasuality},
    \ref{Item:AppendixOperatorPropLipschitz},
    and \ref{Item:AppendixOperatorPropBIBO}, and since it is easy to see that
    the function
    \[
        \tilde{F}:\Rp\times\R^q\to\R^m,\quad (t,z)\mapsto \tilde{F}(t,z)\coloneqq F(u(t),z),
    \]
    is {Carath\'{e}odory function}.
\end{proof}
The same holds true for class of models $\cM^{m,r}_{t_0}$ we considered in this thesis. 
\begin{corollary}\label{Appendix:Cor:ModelSolutionExists}
    Consider model~\eqref{eq:Model_r} with $(\fM,\gM,\oTM)\in\cM^{m,r}_{t_0}$  at initial time $t_0\geq0$.
    Let $\yM^0\in \cC^{r-1}([0,t_0],\R^m)$ be an initial trajectory
    and $u\in L^\infty_{\loc}([t_0,\infty), \R^m)$ be a control function.
    Then, 
\begin{enumerate}[(i)]
    \item the initial value problem~\eqref{eq:ModDiff} has a solution $\xM:[0,\omega)\to \R^{rm}$ in the sense of~\eqref{Appendix:eq:DefSolutionSystem},
    \item every solution can be extended to a maximal solution,
    \item if $\xM:[0,\omega)\to \R^{rm}$ is a bounded maximal solution, then $\omega=\infty$.
\end{enumerate}
\end{corollary}

\chapter*{Nomenclature}
The following notation is used in this thesis to refer to commonly known mathematical concepts
\addcontentsline{toc}{chapter}{Nomenclature and Notation}  
\setlength{\extrarowheight}{10pt}
\begin{longtable}{p{2.25cm} p{\linewidth-2.9cm}}
  $\N$                        &{set of positive integers}
\\$\N_{0}$                    &{$\coloneqq \N\cup\cbl0\cbr$, set of non-negative integers}
\\$\Z$                        &ring of integers
\\$X^\N$                      &set of sequences with elements in a set $X$
\\$\R$                        &{field of real numbers}
\\$\Rp$                       &{$\coloneqq [0,\infty)$, set of non-negative real numbers}
\\$\C$                        &{field of complex numbers}
\\$\C_{\geq 0 (>0)}$          &{$\coloneqq \setdef{ z \in \C}{\ \text{Re}(z)\geq 0\ (>0)}$, the complex (open) right half-plane}
\\$\C_{\leq 0 (<0)}$          &{$\coloneqq \setdef{ z \in \C}{\ \text{Re}(z)\leq 0\ (<0)}$, the complex (open) left half-plane}
\\$\R[s]$                     &ring of polynomials with coefficients in $\R$ and indeterminate $s$
\\$\R^n$                      &vector space of real-valued ordered $n$-tuples with $n\in\N$
\\$\al\cdot,\cdot\ar$         &Euclidean scalar product in $\R^n$
\\$\Norm{\cdot}$              &Euclidean norm in $\R^n$
\\$\cB_\eps(x^0)$             &{$\coloneqq \!\setdef{ x \in \R^n}{\Norm{x-x^0} < \eps}$, open ball around $x^0\in\R^n$ with radius $\eps>0$}
\\$\bar{\cB}_\eps(x^0)$       &{$\coloneqq \!\!\setdef{ x \in \R^n}{\Norm{x-x^0} \leq \eps}$, closed ball around $x^0\in\R^n$ with radius $\eps>0$}
\\$\cB_\eps$                  &{$\coloneqq \cB_{\eps}(0)$, open ball around the origin with radius $\eps>0$}
\\$\bar{\cB}_\eps$            &{$\coloneqq \bar{\cB}_{\eps}(0)$, closed ball around the origin with radius $\eps>0$}
\\$\R^{n\times m}$            &vector space of real-valued $n\times m$ matrices with $n,m\in\N$
\\$I_{n}$                     &the identity matrix in $\R^{n\times n}$
\\$\Norm{A}$                  &$\coloneqq \sup_{\Norm{x}=1}\Norm{Ax}$, induced operator norm for $A\in\R^{n\times m}$
\\$\SGroup_{n}(\R)$           &{set of symmetric $\R^{n\times n}$-matrices}
\\$\SGroup_{n}^{--}(\R)$      &{set of symmetric negative definite $\R^{n\times n}$-matrices}
\\$\GL_{n}(\R)$               &{group of invertible $\R^{n\times n}$ matrices}
\\$\det(A)$                   &{determinant of a square matrix $A\in\R^{n\times n}$}
\\$\spec(A)$                  &{spectrum of a square matrix $A\in\R^{n\times n}$ (set of a eigenvalues of $A$)}
\\$\lambda_{\max}(A)$         &{largest eigenvalue of a symmetric matrix $A\in\SGroup_n$}
\\$\lambda_{\min}(A)$         &{smallest eigenvalue of a symmetric matrix $A\in\SGroup_n$}
\\$\lfloor x\rfloor$          &$\coloneqq \max\setdef{n\in\Z}{n\leq x}$, the floor function for real numbers $x\in\R$
\\$\circ$                     &composition of functions
\\$\graph(f)$                 &$\coloneqq \setdef{(x,f(x))\in X\times Y}{x\in X}$, the graph of a function $f:X\to Y$
\\$f|_A$                      &restriction of the function $f:X\to Y$ to the subset $A$ of the set $X$
\\$\Indic_A$                  &indicator function $\Indic_A:X\to \cbl 0,1\cbr$ of a subset $A$ of a set $X$
\\$\Lebesgue(\cdot)$          &Lebesgue measure
\\$f^+$                       &$\coloneqq \max\{f,0\}$, positive part of a function $f:X\to\R$
\\$f^-$                       &$\coloneqq \max\{-f,0\}$, negative part of a function $f:X\to\R$
\\${\cT_{\cP}(I,\R^n)}$       &space of step functions $f:I\to\R^n$ over an interval $I\subset\R$ with partition~$\cP$, see~\Cref{Def:PartitionAndStepFunction}
\\${\cR(I,\R^n)}$             &{space of regulated functions $f:I\to\R^n$ over an interval $I\subset\R$, see~\Cref{Def:RegulatedFunction}}
\\${\Lip(V,\R^n)}$            &{space of Lipschitz continuous functions $f:V\to\R^n$, where $V\subset\R^m$}
\\${\Lip_{\loc}(V,\R^n)}$     &{space of locally Lipschitz continuous functions $f:V\to\R^n$, where $V\subset\R^m$}
\\${\cC^p(V,\R^n)}$           &{linear space of $p$-times continuously differentiable functions $f:V\to\R^n$, where $V\subset\R^m$ and $p\in\N_0\cup \{\infty\}$}
\\${\cC(V,\R^n)}$             &{$\coloneqq \cC^0(V,\R^n)$, linear space of continuous functions $f:V\to\R^n$, $V\subset\R^m$}
\\${L^p(I,\R^n)}$             &{space of measurable $p$-integrable functions $f:I\to\R^n$ over an interval $I\subset \R$ with norm $\LNorm[p]{\cdot}$ and $p\in\N$}
\\$\al \cdot,\cdot\ar_{L^2}$  &scalar product in $L^2(I,\R^n)$
\\$x_k\rightharpoonup x^{\star}$& weak convergence of $(x_k)\in  L^2(I,\R^n)^\N$ to $x^\star\in  L^2(I,\R^n)$
\\${L^\infty(I,\R^n)}$        &{space of measurable essentially bounded functions $f:I\to\R^n$ over an interval $I\subset \R$}
\\$\esssup\limits_{t\in I}\Norm{f(t)}$&essential supremum of a measurable function $f:I\to \R^n$
\\$\SNorm{f}$                 &$\coloneqq \esssup_{t\in I}\Norm{f(t)}$, norm of $f\in L^\infty(I,\R^n)$
\\${L_{\loc}^\infty(I,\R^n)}$ &{space of measurable locally bounded functions $f:I\to\R^n$ over an~interval $I\subset \R$}
\\${W^{k,\infty}(I,\R^n)}$    &{Sobolev space of all $k$-times weakly differentiable functions $f:I\to\R^n$ over an~interval $I\subset \R$ such that $f,\ldots,f^{(k)}\in L^{\infty}(I,\R^n)$}
\end{longtable}

\chapter*{Notation}
Unless explicitly stated otherwise, the following symbols introduced and used in this work always have the meaning given below.
\setlength{\extrarowheight}{10pt}
\begin{longtable} {p{2.7cm} p{\linewidth-3.54cm}}
$t_0$                         & $\geq 0$, initial time
\\$\hat{t}$                   & $\geq t_0$, arbitrary time instant 
\\$\cN^{m,r}_{t_0}$           & system class with $m,r\in\N$, see \Cref{Def:SystemClass}
\\$(F,\oT)$                   &$\in\cN^{m,r}_{t_0}$
\\$y$                         & output of system~\eqref{eq:Sys}
\\$y^0$                       & $\in\cC^{r-1}([0,t_0],\R^m)$, initial system trajectory
\\$\cT_{t_0}^{n,q}$           & operator class as defined in~\Cref{Def:OperatorClass}, with $n,q\in\N$
\\$\tau$                      & $\geq 0$, memory limit of operator $\oTM\in\cT_{t_0}^{n,q}$ 
\\$\cM^{m,r}_{t_0}$           & model class with $m,r\in\N$, see \Cref{Def:ModelClass}
\\$\cM^{m,r}_{t_0, \umax,\bar{\rho}}$& restricted model class with $\bar{\rho},\umax\geq0$, see \Cref{Def:RestrictedModelClass}
\\$(\fM,\gM,\oTM)$            & $\in\cM^{m,r}_{t_0}$, model used in the MPC algorithm
\\$\xM$                       & solution of the model differential equation~\eqref{eq:Model_r}, see \Cref{Def:ModSolution}
\\$\yM$                       & output of model~\eqref{eq:Model_r}
\\$\yM^0$                     & $\in\cC^{r-1}([0,t_0],\R^m)$, initial model trajectory
\\$y_{\rf}$                   & $\in W^{r,\infty}(\Rp,\R^m)$, reference trajectory
\\$e$                         & $\coloneqq y-y_{\rf}$, tracking error 
\\$\eMTrack$                  & $\coloneqq \yM-y_{\rf}$, model's tracking error \eqref{eq:Intro:ModelTrackingerror}
\\$\eSTrack$                  & $\coloneqq y-\yM$, model-system output mismatch \eqref{eq:Intro:ModelPlantMismatch}
\\$\cG$                       & $\coloneqq \setdef{\Funnel\in W^{1,\infty}(\Rp,\R)} { \inf_{t\geq 0}\Funnel(t) > 0 }$, set of admissible funnel functions as defined in~\eqref{eq:DefSetOfFunnelFunctions}
\\$\Funnel$                   & $\in\cG$, funnel boundary function
\\$\FunDeriv,\FunDiam$        & $>0$, constants associated to $\Funnel\in\cG$ satisfying \eqref{eq:DefinitionAlphaBeta} 
\\$\cF_{\Funnel}$             & $\coloneqq \setdef{(t,e)\in \Rp\times\R^{m}}{\Norm{e} < \Funnel(t)}$, funnel for $\psi\in\cG$, defined in~\eqref{eq:DefFunnel} 
\\$\bar{\cF}_{\Funnel}$       & $\coloneqq \setdef{(t,e)\in \Rp\times\R^{m}}{\Norm{e} \leq \Funnel(t)}$,  for $\psi\in\cG$
\\$\delta$                    & $>0$, time-shift of the MPC algorithm
\\$T$                         & $\geq\delta>0$, prediction horizon of the MPC algorithm
\\$(t_k)_{k\in\N_0}$          & time sequence defined by $t_k \coloneqq  t_0+k\delta$ for $k\in\N_0$
\\$\umax$                     & $\geq0$, input saturation level of the MPC algorithm
\\$\uFMPCk$                   &$\in L^\infty([t_k,t_k +T],\R^{m})$, solution of the optimal control problem~\eqref{eq:FunnelMpcOCP}
\\$\xM^k$                     & solution of the model differential equation on the interval $[t_k,t_{k+1}]$
\\$\yM^k$                     & predicted output of model~\eqref{eq:Model_r} on the interval $[t_k,t_{k+1}]$
\\$\uFMPC(t)$                 & $=\uFMPCk(t)$ for $t\in[t_k,t_{k+1})$ and $k\in\N_0$
\\$\FunnelPenaltyFunc$        &funnel penalty function for $\Funnel\in\cG$, see \Cref{Def:FunnelPenaltyFunction}
\\$\OrigFunnelStageCost$      &funnel stage cost function for $\Funnel\in\cG$, see  \Cref{Def:FunnelStageCostFunc}
\\$J^{\Psi}_T(\cdot;\hat{t},\InitState)$ & $L^\infty([\hat{t},\hat{t}+T],\R^{m})\to\R\cup\{\infty\}$, cost functional as defined in \eqref {eq:DefCostFunctionJ}
\\$I_{t_0}^{\hat{t},\tau}$    & $\coloneqq [\hat{t}-\tau,\hat{t}]\cap[t_0,\hat{t}]$ for $\hat{t}\geq t_0\geq0$ and $\tau\geq 0$
\\$\OpChi$                    & shorthand notation for $\OpChi(\zeta)(t)\coloneqq (\zeta(t),\dot{\zeta}(t),\ldots,\zeta^{(r-1)}(t))\in\R^{rm}$ for a function  $\zeta\in W^{r,\infty}(I,\R^m)$ as defined in~\eqref{eq:DefOperatorChi} 
\\$\LShift$                   & left shift operator as defined in \eqref{eq:LeftShift}
\\$\eM_{i}$                   & auxiliary error variable as defined in \eqref{eq:ErrorVar} for $i=1,\ldots, r$
\\$\Psi$                     &$\coloneqq(\Funnel_1,\ldots,\Funnel_r)\in\FunnelBoundaryFuncs$, auxiliary funnel functions, see~\eqref{eq:DefFunnelBoundaryFunctions}
\\$k_i$                       & $\geq 0$, parameter associated to $\Funnel_i$ and $\eM_{i}$ for  $i=1,\ldots, r$
\\$\cD^{\Psi}_{t}$            & $\coloneqq \setdef {z\in\R^{rm}} { \Norm{\eM_{i}(z)}<\psi_i(t),\ i=1,\ldots, r }$ for $t\in\Rp$ and $\Psi\in\cG^r$ as defined in~\eqref{eq:DefSetD}
\\$\FunnelTrajectories_{\hat{t}}$&$\coloneqq \setdef{\zeta\in \cR(\Rp,\R^{rm})} { \fa t\in[0,\hat{t}]: \zeta(t)-\OpChi(y_{\rf})(t)\in\cD_t^{\Psi} }$ for $\hat{t}\geq t_0$ and $\Psi\in\FunnelBoundaryFuncs$, defined in~\eqref{eq:Def:FunnelTrajectories}
\\$\InitValues(\hat{t})$      & set of feasible initial values at time $\hat{t}\geq t_0$ for $\Psi\in\FunnelBoundaryFuncs$ and $\tau\geq 0$, see
\Cref{Def:SetInitialValues}
\\$\PropInitValues(\hat{t},\hat{x})$ & proper initial model state for $\eps,\lambda\in(0,1)$ at time $\hat{t}\geq t_0$ given system data $\hat{x}\in\R^{rm}$, see 
    \Cref{Def:ProperInitValues}
\\$\InitState$                &$\coloneqq(\xMh,\oTMh)\in\InitValues(\hat{t})$, model initial value at time $\hat{t}\geq t_0$
\\$\InitStateK_{k}$           &$\coloneqq(\xM^k,\oTM^k)\in\InitValues(t_k)$, model initial state at time $t_k$
\\$\kappa$                    &initialisation strategy for the model, see \Cref{Def:InitialisationStrategy}
\\$\Controls(\umax,\InitState)$ &set of feasible controls on the interval $[\hat{t},\hat{t}+T]$ as defined in \eqref{eq:Def-U}
\\$\Controls^{\cP}(\umax,\InitState)$ & $\coloneqq \cT_{\cP}([\hat{t},\hat{t}+T],\R^m)\cap\Controls(\umax,\InitState)$, see~\eqref{eq:Def-DiscrU} 
\\ $\begin{array}{l}
     \fMmax, \gMmax,  \\
     \gMInvmax, \gmin 
\end{array}$& $\geq0$, bounds on the dynamics of the model, see \Cref{Lemma:DynamicBounded}
\\ $ \fmax, \gmax, \gmin 
$& $\geq0$, bounds on the dynamics of the system, see \Cref{Lemma:DynamicBoundedDiscreteSys}
\\$\FCSurjec$                 & $\in\cC(\Rp,\R)$, surjection used in the funnel controller
\\$\FCBijec$                  & $\in\cC^1([0,1),[1,\infty))$, bijection used in the funnel controller
\\$\eS_i(\phi,z)$             & auxiliary error variable for $\phi>0$, $z\in\R^{rm}$, $i =1,\ldots, r$, see~\eqref{eq:ek_FC}
\\$\cE^\eps_{i}(\phi)$        & $\coloneqq \setdef{z\in\R^{rm}}{\Norm{\eS_j(\phi,z)}<\eps,\ j=1,\ldots, i}$ for ${i\in \{1,\ldots, r\}}$ and ${\eps\in(0,1]}$, see~\eqref{eq:ek_FC}
\\$\FCTrajectories_{\hat{t}}$ & $\coloneqq\!\! \setdef
        {\!\zeta\in \cC^{r-1}(\Rp,\R^{m})\!}
        {
            \zeta|_{[0,t_0]}=y^0,
            \fa t\in [t_0,\hat{t}):\OpChi(\zeta)(t)\in\cEFC{1}(\phi(t))\!
        }$, see \eqref{eq:DefFCTrajectories}
\\$\phi(t)$ &$\coloneqq \frac{1}{\Funnel(t) - \| \yM(t) - y_{\rf}(t)\|}$, boundary for funnel controller component based on model prediction, see \eqref{eq:DefPhi}
\\$\phi_k$ &$[t_k,t_{k+1}]\to \Rp$, adaptive funnel based on model prediction, see~\eqref{alg:eq:vp}
\\$\uFCk(t)$ &$\coloneqq   (\FCSurjec\circ\FCBijec) (\Norm{\eS_{r}(\phi_{k}(t),\eSTrack(t))}^2)\eS_{r}(\phi_{k}(t),\eSTrack(t))$, control law of the funnel controller component on the interval $[t_k,t_{k+1})$, see \eqref{eq:uFCRobustFMPC}
\\$\uFC(t)$ &$=\uFCk(t)$ for $t\in[t_k, t_{k+1})$ for $k\in\N_0$
\\$\ActivFunc$                & $:[0,1] \to [0,\ActivFunc^+]$ with $\ActivFunc^+ > 0$, activation function, see~\Cref{Sec:RobustFMPCInit}
\\$\mathfrak{S}_{\hat{t}}$    & $\coloneqq \cC^r([t_0,\hat{t}],\R^m)\times \cR([t_0,\hat{t}],\R^m)^r\times L^\infty([t_0,\hat{t}],\R^m)\times L^\infty([t_0,\hat{t}],\R^m)$ for $\hat{t}\geq t_0$, see~\eqref{eq:DefSetOfSignals}
\\ $\cL$&$\bigcup_{t\geq t_0}\SetOfSignals_{t}\to \cM^{m,r}_{t_0,\umax,\bar{\rho}}$, feasible learning scheme, see \Cref{Def:LearningScheme}
\\$\ModSysClassDiscr$ & model/system class considered in \Cref{Chapter:DiscretFMPC} for the sampled-data robust funnel MPC 
\\$\uZoH$ & funnel controller with zero-order-hold, defined in~\eqref{eq:controller_recursive}
\\$\FCDiscreteGain$ &$\geq 0$, input gain of $\uZoH$ 
\\$\FCDiscreteThresh$ &$\geq 0$, activation threshold of $\uZoH$ 
\\$\SampleTime$ &$>0$, sampling time of $\uZoH$ and the sampled-data funnel MPC~\Cref{Algo:DiscrFunnelMPC}
\end{longtable}

%\begin{description}%[labelindent=-0.5cm]

\nocite{BergDenn21}
\nocite{BergDenn22}
\nocite{BergDenn24}
\nocite{BergDenn24b}
\nocite{LanzaDenn24}
\nocite{LanzaDenn24b}
\nocite{Dennst24}
\nocite{Oppeneiger24}
\nocite{Dennst25}
\chapter*{Bibliography}
\addcontentsline{toc}{chapter}{Bibliography}  
All scientific articles written in the context of this dissertation are listed below in the section Publications and Preprints.
All cited publications can be found thereafter in References.
\defbibnote{NoteOwnArticles}{In connection with the work on this dissertation, the following articles were published in international scientific journals.}
\begingroup
\sloppy
\printbibliography[heading=subbibliography,keyword={own},prenote=NoteOwnArticles,title={Publications}]
\printbibliography[heading=subbibliography,notkeyword={own}]

@article{BergDenn21,
    author = {Berger, Thomas and Dennst\"{a}dt, Dario and Ilchmann, Achim and Worthmann, Karl},
    title = {Funnel {M}odel {P}redictive {C}ontrol for {N}onlinear {S}ystems with {R}elative {D}egree {O}ne},
    journal = {SIAM Journal on Control and Optimization},
    volume = {60},
    number = {6},
    pages = {3358-3383},
    year = {2022},
    doi = {10.1137/21M1431655}
}

@article{BergDenn22,
    title={Funnel {MPC} with feasibility constraints for nonlinear systems with arbitrary relative degree},
    author={Berger, Thomas and Dennst{\"a}dt, Dario},
    journal={IEEE Control Systems Letters},
    year={2022},
    publisher={IEEE},
    doi={10.1109/LCSYS.2022.3178478}
}

@article{BergDenn24,
    title = {{Funnel MPC for nonlinear systems with arbitrary relative degree}},
    author = {Thomas Berger and Dario Dennst{\"a}dt},
    journal = {Automatica},
    volume = {167},
    pages = {111759},
    year = {2024},
    issn = {0005-1098},
    doi = {10.1016/j.automatica.2024.111759}
}

@article{BergDenn24b,
    author = {Berger, Thomas and Dennst\"{a}dt, Dario and Lanza, Lukas and Worthmann, Karl},
    title = {{Robust Funnel Model Predictive Control for Output Tracking with Prescribed Performance}},
    journal = {SIAM Journal on Control and Optimization},
    volume = {62},
    number = {4},
    pages = {2071-2097},
    year = {2024},
    doi = {10.1137/23M1551195}
}

@article{LanzaDenn24b,
    author = {Lanza, Lukas and Dennstädt, Dario and Berger, Thomas and Worthmann, Karl},
    title = {{Safe Continual Learning in Model Predictive Control With Prescribed Bounds on the Tracking Error}},
    journal = {International Journal of Robust and Nonlinear Control},
    year={2025},
    doi = {10.1002/rnc.8001}
}

@article{LanzaDenn24,
    title={Sampled-data funnel control and its use for safe continual learning},
    author={Lanza, Lukas and Dennst{\"a}dt, Dario and Worthmann, Karl and Schmitz, Philipp and  {Şen Gökçen Devlet} and  Trenn, Stephan and Schaller, Manuel},
    journal = {Systems \& Control Letters},
    volume = {192},
    pages = {105892},
    year = {2024},
    issn = {0167-6911},
    doi = {10.1016/j.sysconle.2024.105892}
}

@INPROCEEDINGS{Dennst24,
    author={Dennst{\"a}dt, Dario and Lanza, Lukas and  Karl Worthmann},
    booktitle={European Control Conference (ECC)}, 
    title={{On Model Predictive Control with Sampled-Data Input for Output Tracking with Prescribed Performance}}, 
    year={2024},
    volume={},
    number={},
    pages={2978-2984},
    doi={10.23919/ECC64448.2024.10590848}
}

@article{Oppeneiger24,
    title = {Model predictive control of a magnetic levitation system with prescribed output tracking performance},
    author={Oppeneiger, Benedikt and Lanza, Lukas and Schell, Maximilian and  Dennst{\"a}dt, Dario and
          Schaller, Manuel and Zamzow, Bert and Berger, Thomas and Worthmann, Karl},
    journal = {Control Engineering Practice},
    volume = {151},
    pages = {106018},
    year = {2024},
    issn = {0967-0661},
    doi = {10.1016/j.conengprac.2024.106018}
}

@article{Dennst25, 
    title = {A low-complexity funnel control approach for non-linear systems of higher-order},
    journal = {IFAC-PapersOnLine},
    volume = {59},
    number = {14},
    pages = {7-12},
    year = {2025},
    note = {15th IFAC Workshop on Adaptive and Learning Control Systems ALCOS 2025},
    issn = {2405-8963},
    doi = {10.1016/j.ifacol.2025.12.117},
    author = {Dario Dennstädt}
}

@STRING{ adneural ={Advances in Neural Information Processing Systems}}

@STRING{ aiche = {AIChE Journal}}

@STRING{ acc2004 =  {American Control Conference (ACC)}}

@STRING{ acc2006 =  {American Control Conference (ACC)}}

@STRING{ acc2019 =  {American Control Conference (ACC)}}

@STRING{ acc2020 =  {American Control Conference (ACC)}}

@STRING{ acc2021 =  {American Control Conference (ACC)}}

@STRING{ ascasme = {ASCE-ASME Journal of Risk and Uncertainty in Engineering Systems Part B: Mechanical Engineering}}

@STRING{ annurevcon   = {Annual Review of Control, Robotics, and Autonomous Systems}}

@STRING{ annufluid = {Annual Review of Fluid Mechanics}}

@STRING{ arcontrol    = {Annual Reviews in Control}}

@STRING{ atauto       = {at--Automatisierungs\-technik} }

@STRING{ asiacont     = {Asian Journal of Control}}

@STRING{ auto         = {Automatica} }

@STRING{ autosofteng  = {Automated Software Engineering}}

@STRING{ cdc49      = {49th IEEE Annual Conference on Decision and Control (CDC)}}

@STRING{ cdc54 =      {54th IEEE Annual Conference on Decision and Control (CDC)}}

@STRING{ cdc55  =     {55th IEEE Annual Conference on Decision and Control (CDC)}}

@STRING{ cdc56  =     {56th IEEE Annual Conference on Decision and Control (CDC)}}

@STRING{ cdc57  =     {57th IEEE Annual Conference on Decision and Control (CDC)}}

@STRING{ cdc61      = {61st IEEE Annual Conference on Decision and Control (CDC)}}

@STRING{ cep        = {Control Engineering Practice}}

@STRING{ compchem   = {Computers \& Chemical Engineering}}

@STRING{ ecc2014  = {European Control Conference (ECC)}}

@STRING{ ecc2015 =  {European Control Conference (ECC)}}

@STRING{ ecc2016  = {European Control Conference (ECC)}}

@STRING{ ecc2019 =  {European Control Conference (ECC)}}

@STRING{ ecc2022 =  {European Control Conference (ECC)}}

@STRING{ ecc2024 =  {European Control Conference (ECC)}}

@STRING{ ejc        = {European Journal of Control} }

@STRING{ esiamcocv  = {ESAIM: Control, Optimisation and Calculus of Variations}}

@STRING{ ICSTCC23  = {27th International Conference on System Theory, Control and Computing (ICSTCC)}}

@STRING{ iccas2010 = {International Conference on Control, Automation and Systems (ICCAS)}}

@STRING{ IEEProced = {IEE Proceedings - Control Theory and Applications}}

@STRING{ ieeeaccess    = {IEEE Access} }

@STRING{ ieeeproc = {Proceedings of the IEEE}}

@STRING{ ieeetac    = {IEEE Transactions on Automatic Control} }

@STRING{ ieeetcst = {IEEE Transactions on Control Systems Technology}}

@string{ ieeetrob = {IEEE Transactions on Robotics}}

@STRING{ ieeetpami   = {IEEE Transactions on Pattern Analysis and Machine Intelligence} }

@STRING{ ieeetnnls = {IEEE Transactions on Neural Networks and Learning Systems}}

@STRING{ ieeecsl    = {IEEE Control Systems Letters}}

@STRING{ ieeeconsysMag ={IEEE Control Systems Magazine}}

@STRING{ ieeessci = {IEEE Symposium Series on Computational Intelligence (SSCI)}}

@STRING{ ieeeicra17 = {2017 IEEE International Conference on Robotics and Automation (ICRA)}}

@STRING{ ieeeicra23 = {2023 IEEE International Conference on Robotics and Automation (ICRA)}}

@STRING{ industEngChem = {Industrial \& Engineering Chemistry Research}}

@STRING{ cavs2019 = {2019 IEEE 2nd Connected and Automated Vehicles Symposium (CAVS)}}

@STRING {ieeeicit2019 ={2019 IEEE International Conference on Industrial Technology (ICIT)}}

@STRING{ ifaconline = {IFAC-PapersOnLine}}

@STRING{ ifacprec = {IFAC Proceedings Volumes}}

@STRING {facNonPC4 = {4th IFAC Conference on Nonlinear Model Predictive Control}}

@STRING {ifcaworld14 = {14th IFAC World Congress}}

@STRING {ifcaworld17 = {17th IFAC World Congress}}

@STRING {ifcaworld20 = {20th IFAC World Congress}}

@STRING {ifcaworld21 = {21st IFAC World Congress}}

@STRING {ifacnolcos13 ={13th IFAC Symposium on Nonlinear Control Systems NOLCOS}}

@STRING {imacontrol = {IMA Journal of Mathematical Control and Information}}

@STRING { jprocescon ={Journal of Process Control}}

@STRING{ intjacsp  = {International Journal of Adaptive Control and Signal Processing}}

@STRING{ intjoc    = {International Journal of Control} }

@STRING{ intjautomation = {International Journal of Control, Automation and Systems} }

@STRING{ introbust = {International Journal of Robust and Nonlinear Control}}

@STRING{ introbosys = {Journal of Intelligent \& Robotic Systems}}

@STRING{ mcss    = {Mathematics of Control, Signals, and Systems} }

@STRING{ mlr     = {Journal of Machine Learning Research} }

@STRING{ mtns24 = {26th International Symposium on Mathematical Theory of Networks and Systems MTNS 2024}}

@STRING{ jourphysic= {Journal of Computational Physics}}

@STRING{ journonlin = {Journal of Nonlinear Science}}

@STRING{ med2021 = {2021 29th Mediterranean Conference on Control and Automation (MED)}}

@STRING{ rtucon14 = {55th International Scientific Conference on Power and Electrical Engineering of Riga Technical University (RTUCON)}}

@STRING{ scl        = {Systems \& Control Letters} }

@STRING{ siamjco    = {SIAM Journal on Control and Optimization} }

@STRING{ siamads = {SIAM Journal on Applied Dynamical Systems}}

@STRING{ siamrev = {SIAM Review}}

@STRING{ nonlindyn = {Nonlinear Dynamics}}

@STRING{ mathprog = {Mathematical Programming Computation}}

@STRING{ mecsci  = {Mechanical Sciences}}

@article{koopman1931hamiltonian,
    author = {Bernard O. Koopman},
    title = {{Hamiltonian Systems and Transformation in Hilbert Space}},
    journal = {Proceedings of the National Academy of Sciences},
    volume = {17},
    number = {5},
    pages = {315-318},
    year = {1931},
    doi = {10.1073/pnas.17.5.315},
}

@ARTICLE{Shannon1949,
    author={Shannon, Claude E.},
    journal={Proceedings of the IRE}, 
    title={{Communication in the Presence of Noise}}, 
    year={1949},
    volume={37},
    number={1},
    pages={10-21},
    doi={10.1109/JRPROC.1949.232969}
}

@article{Kalman61,
    author = {Kálmán, Rudolf E. and Bucy, Richard S.},
    title = {{New Results in Linear Filtering and Prediction Theory}},
    journal = {Journal of Basic Engineering},
    volume = {83}, 
    number = {1},
    pages = {95-108},
    year = {1961},
    month = {03},
    issn = {0021-9223},
    doi = {10.1115/1.3658902}
}

@book{Rudi76,
    title={Principles of {M}athematical {A}nalysis},
    author={Rudin, Walter},
    volume={3},
    year={1976},
    publisher={McGraw-hill New York},
    isbn = {978-0070542358}
}

@ARTICLE{Mita1980,
    author={Mita, Tsutomu},
    journal=ieeetac, 
    title={A relation between overshoot and sampling period in sampled data feedback control systems}, 
    year={1980},
    volume={25},
    number={3},
    pages={603-604},
    doi={10.1109/TAC.1980.1102334}
}

@phdthesis{watkins1989learning,
    title={{Learning from Delayed Rewards}},
    author={Watkins, Christopher John Cornish Hellaby},
    year={1989},
    school={King's College},
    address = {Cambridge, United Kingdom},
    url = {https://www.cs.rhul.ac.uk/~chrisw/new_thesis.pdf}, 
    urldate  = {2025-06-24},
}

@article{sussmann1990limitations,
    title={Limitations on the stabilizability of globally-minimum-phase systems},
    author={Sussmann, Héctor J.},
    journal= ieeetac,
    volume={35},
    number={1},
    pages={117--119},
    year={1990},
    publisher={IEEE},
    doi={10.1109/9.45159}
}

@article{ilchmann1991non,
    author = {Ilchmann, Achim},
    title = {{Non-Identifier-Based Adaptive Control of Dynamical Systems: A Survey}},
    journal = imacontrol,
    volume = {8},
    number = {4},
    pages = {321-366},
    year = {1991},
    month = {12},
    issn = {0265-0754},
    doi = {10.1093/imamci/8.4.321},
}

@ARTICLE{VielJado97,
    author={Viel, F. and Jadot, Fabrice and Bastin, Georges},
    journal=ieeetac, 
    title={Robust feedback stabilization of chemical reactors}, 
    year={1997},
    volume={42},
    number={4},
    pages={473-481},
    doi={10.1109/9.566657}
}

@book{Walt98,
    AUTHOR    = {Walter, Wolfgang},
    YEAR      = 1998,
    TITLE     = {{Ordinary Differential Equations}},
    PUBLISHER = {Springer},
    Address   = {New York},
    ISBN      = {978-0-387-98459-9},
    doi={10.1007/978-1-4612-0601-9}
}

@article{chen1998quasi,
    title={A quasi-infinite horizon nonlinear model predictive control scheme with guaranteed stability},
    author={Chen, Hong and Allg{\"o}wer, Frank},
    journal=auto,
    volume={34},
    number={10},
    pages={1205--1217},
    year={1998},
    publisher={Elsevier},
    doi={10.1016/S0005-1098(98)00073-9}
}

@ARTICLE{Scokaert1998,
  author={Scokaert, Pierre O.M. and Mayne, David Q.},
  journal=ieeetac, 
  title={Min-max feedback model predictive control for constrained linear systems}, 
  year={1998},
  volume={43},
  number={8},
  pages={1136-1142},
  doi={10.1109/9.704989}
}

@article{sussmann1991peaking,
    title={The peaking phenomenon and the global stabilization of nonlinear systems},
    author={Sussmann, Héctor J. and Kokotovic, Petar V.},
    journal= ieeetac,
    volume={36},
    number={4},
    pages={424--440},
    year={1991},
    publisher={IEEE},
    doi={10.1109/9.75101}
}

@article{ByrnIsid91a,
    AUTHOR    = {Byrnes, Christopher I. and Isidori, Alberto},
    YEAR      = 1991,
    TITLE     = {Asymptotic stabilization of minimum phase nonlinear systems},
    JOURNAL   = ieeetac,
    Volume    = 36,
    Number    = 10,
    Pages     = {1122--1137},
    publisher={IEEE},
    doi={10.1109/9.90226}
}

@article{Leung1991,
    title = {Performance analysis of sampled-data control systems},
    journal = auto,
    volume = {27},
    number = {4},
    pages = {699-704},
    year = {1991},
    issn = {0005-1098},
    doi = {10.1016/0005-1098(91)90060-F},
    author = {Gary M.H. Leung and T.P. Perry and Bruce A. Francis},
}

@article{IlchRyan94,
    AUTHOR    = {Ilchmann, Achim and Ryan, Eugene P.},
    YEAR      = 1994,
    TITLE     = {{Universal $\lambda$-Tracking for Nonlinearly-Perturbed Systems in the Presence of Noise}},
    JOURNAL   = auto,
    Volume    = 30,
    Number    = 2,
    Pages     = {337--346},
    publisher = {Elsevier},
    doi={10.1016/0005-1098(94)90035-3}
}

@book{Isid95,
    title={{Nonlinear Control Systems: An Introduction}},
    author={Isidori, Alberto},
    year={1985},
    publisher={Springer},
    address={Heidelberg},
    isbn={3540506012},
    doi={10.1007/BFb0006368}
}

@article{Mors96,
    AUTHOR    = {Morse, A. Stephen},
    YEAR      = 1996,
    TITLE     = {Overcoming the obstacle of high relative degree},
    JOURNAL   = ejc,
    Volume    = 2,
    Number    = 1,
    Pages     = {29--35},
    doi = {10.1016/S0947-3580(96)70025-0}
}

@book{astrom_wittenmark_1997,
    author    = {Åström, Karl Johan and Wittenmark, Björn},
    title     = {{Computer-controlled Systems: Theory and Design}},
    edition   = {3rd},
    publisher = {Prentice Hall},
    address   = {Upper Saddle River, NJ},
    year      = {1997},
    isbn      = {0-13-314899-8}
}

@InProceedings{bemporad99,
    author={Bemporad, Alberto and Morari, Manfred},
    title={{Robust model predictive control: A survey}},
    booktitle={Robustness in identification and control},
    year={1999},
    publisher={Springer},
    address={London},
    pages={207--226},
    isbn={978-1-84628-538-7},
    doi={doi.org/10.1007/BFb0109870}
}

@article{favoreel1999spc,
    title = {{SPC: Subspace Predictive Control}},
    journal = {IFAC Proceedings Volumes},
    volume = {32},
    number = {2},
    pages = {4004-4009},
    year = {1999},
    note = ifcaworld14,
    issn = {1474-6670},
    doi = {10.1016/S1474-6670(17)56683-5},
    author = {Wouter Favoreel and Bart De Moor and Michel Gevers},
}

@article{NESIC1999259,
    title = {Sufficient conditions for stabilization of sampled-data nonlinear systems via discrete-time approximations},
    journal = scl,
    volume = {38},
    number = {4},
    pages = {259-270},
    year = {1999},
    issn = {0167-6911},
    doi = {10.1016/S0167-6911(99)00073-0},
    author = {Dragan Nešić and Andrew R. Teel and Petar V. Kokotović},
}

@inproceedings{kerrigan2000soft,
    title={Soft constraints and exact penalty functions in model predictive control},
    author={Kerrigan, Eric C. and Maciejowski, Jan M.},
    booktitle={Proceedings of the UKACC International Conference on Control (Control 2000)},
    pages={2319--2327},
    year={2000},
    address= {Cambridge, United Kingdom},
    url={http://www-control.eng.cam.ac.uk/Homepage/papers/cued_control_53.pdf},
    urldate  = {2025-06-24}
}

@article{ryan2001controlled,
    title={Controlled functional differential equations and adaptive stabilization},
    author={Ryan, Eugene P. and Sangwin, C.J.},
    journal=intjoc,
    volume={74},
    number={1},
    pages={77--90},
    year={2001},
    publisher={Taylor \& Francis},
    doi={10.1080/00207170150202698}
}

@article{Monaco2001,
    title = {{Issues on Nonlinear Digital Control}},
    journal = ejc,
    volume = {7},
    number = {2},
    pages = {160-177},
    year = {2001},
    issn = {0947-3580},
    doi = {10.3166/ejc.7.160-177},
    author = {Salvatore Monaco and Dorothée Normand-Cyrot},
}

@article{chisci2001systems,
    title = {Systems with persistent disturbances: predictive control with restricted constraints},
    journal = auto,
    volume = {37},
    number = {7},
    pages = {1019-1028},
    year = {2001},
    issn = {0005-1098},
    doi = {10.1016/S0005-1098(01)00051-6},
    author = {Luigi Chisci and J. Anthony Rossiter and Giovanni Zappa},
}

@article{perkins2002lyapunov,
    title={Lyapunov design for safe reinforcement learning},
    author={Perkins, Theodore J. and Barto, Andrew G.},
    journal=mlr,
    volume={3},
    pages={803--832},
    year={2002},
    url={https://www.jmlr.org/papers/v3/perkins02a.html},
    urldate  = {2025-06-24}
}

@article{IlchRyan02a,
    title={Systems of controlled functional differential equations and adaptive tracking},
    author={Ilchmann, Achim and Ryan, Eugene P. and Sangwin, Christopher J.},
    journal=siamjco,
    volume={40},
    number={6},
    pages={1746--1764},
    year={2002},
    publisher={SIAM},
    doi={10.1137/S0363012900379704}
}

@article{IlchRyan02b,
    AUTHOR    = {Ilchmann, Achim and Ryan, Eugene P. and Sangwin, Christopher J.},
    YEAR      = 2002,
    TITLE     = {Tracking with prescribed transient behaviour},
    JOURNAL   = esiamcocv,
    Volume    = 7,
    Pages     = {471--493},
    doi= {10.1051/cocv:2002064}
}

@article{chen2003terminal,
    title={On the terminal region of model predictive control for non-linear systems with input/state constraints},
    author={Chen, Wen-Hua and O'Reilly, John and Ballance, Donald J.},
    journal= intjacsp,
    volume={17},
    number={3},
    pages={195--207},
    year={2003},
    publisher={Wiley Online Library},
    doi={10.1002/acs.731}
}

@article{QinBadg03,
    title={A survey of industrial model predictive control technology},
    author={Qin, S. Joe and Badgwell, Thomas A.},
    journal=cep,
    volume={11},
    number={7},
    pages={733--764},
    year={2003},
    publisher = {Elsevier},
    doi = {10.1016/S0967-0661(02)00186-7}
}

@article{WILLS20041415,
    title = {Barrier function based model predictive control},
    journal = auto,
    volume = {40},
    number = {8},
    pages = {1415-1422},
    year = {2004},
    issn = {0005-1098},
    doi = {10.1016/j.automatica.2004.03.002},
    author = {Adrian G. Wills and William P. Heath},
}

@article{MILANESE2004957,
    title = {{Set Membership identification of nonlinear systems}},
    journal = auto,
    volume = {40},
    number = {6},
    pages = {957-975},
    year = {2004},
    issn = {0005-1098},
    doi = {10.1016/j.automatica.2004.02.002},
    author = {Mario Milanese and Carlo Novara},
}

@article{langson2004robust,
    title = {Robust model predictive control using tubes},
    journal = auto,
    volume = {40},
    number = {1},
    pages = {125-133},
    year = {2004},
    issn = {0005-1098},
    doi = {10.1016/j.automatica.2003.08.009},
    author = {Wilbur Langson and Ioannis Chryssochoos and Saša V. Raković and David Q. Mayne},
}

@INPROCEEDINGS{Kocijan2004,
    author={Kocijan, Juš and Murray-Smith, Roderick and Rasmussen, Carl E. and Girard, Agathe},
    booktitle=acc2004, 
    title={Gaussian process model based predictive control}, 
    year={2004},
    volume={3},
    number={},
    pages={2214-2219},
    organization={IEEE},
    doi={10.23919/ACC.2004.1383790}
}

@article{IlchTren04,
    AUTHOR    = {Ilchmann, Achim and Trenn, Stephan},
    YEAR      = 2004,
    TITLE     = {Input constrained funnel control with applications to chemical reactor models},
    JOURNAL   = scl,
    Volume    = 53,
    Number    = 5,
    Pages     = {361--375},
    publisher = {Elsevier},
    doi= {10.1016/j.sysconle.2004.05.014}
}

@article{Haseltine2005,
    author = {Haseltine, Eric L. and Rawlings, James B.},
    title = {{Critical Evaluation of Extended Kalman Filtering and Moving-Horizon Estimation}},
    journal = industEngChem,
    volume = {44},
    number = {8},
    pages = {2451-2460},
    year = {2005},
    doi = {10.1021/ie034308l}
}

@article{MaynSero05,
    title = {Robust model predictive control of constrained linear systems with bounded disturbances},
    journal = auto,
    volume = {41},
    number = {2},
    pages = {219-224},
    year = {2005},
    issn = {0005-1098},
    doi = {10.1016/j.automatica.2004.08.019},
    author = {David Q. Mayne and Maria M. Seron and Saša V. Raković},
}

@article{yuz2005sampled,
    title={On sampled-data models for nonlinear systems},
    author={Yuz, Juan I. and Goodwin, Graham C.},
    journal=ieeetac,
    volume={50},
    number={10},
    pages={1477--1489},
    year={2005},
    publisher={IEEE},
    doi={10.1109/TAC.2005.856640}
}

@article{WRMDM05,
    author    = {Jan C. {Willems} and Paolo {Rapisarda} and Ivan {Markovsky} and Bart L.M. {De Moor}},
    title     = {{A note on persistency of excitation}},
    journal   = scl,
    issn      = {0167-6911},
    volume    = {54},
    number    = {4},
    pages     = {325--329},
    year      = {2005},
    address   = {Amsterdam},
    publisher = {Elsevier},
    doi       = {10.1016/j.sysconle.2004.09.003},
    msc2010   = {93B30 93C15},
    zbl       = {1129.93362}
}

@Inbook{Laila2006,
    author="Laila, Dina Shona
    and Ne{\v{s}}i{\'{c}}, Dragan
    and Astolfi, Alessandro",
    editor="Lor{\'i}a, Antonio
    and Lamnabhi-Lagarrigue, Fran{\c{c}}oise
    and Panteley, Elena",
    title={{Sampled-data Control of Nonlinear Systems}},
    bookTitle="Advanced Topics in Control Systems Theory: Lecture Notes from FAP 2005",
    year="2006",
    publisher="Springer",
    address="London",
    pages="91--137",
    isbn={978-1-84628-418-2},
    doi={10.1007/11583592_3},
}

@article{Rao2006,
    author = {Ganti P. Rao  and Heinz Unbehauen},
    title = {Identification of continuous-time systems},
    journal = IEEProced,
    volume = {153},
    issue = {2},
    pages = {185-220},
    year = {2006},
    doi = {10.1049/ip-cta:20045250},
}

@ARTICLE{Bristow2006,
    author={Bristow, Douglas A. and Tharayil, Marina and Alleyne, Andrew G.},
    journal=ieeeconsysMag, 
    title={A survey of iterative learning control}, 
    year={2006},
    volume={26},
    number={3},
    pages={96-114},
    doi={10.1109/MCS.2006.1636313}
}

@article{goulart2006optimization,
    title = {Optimization over state feedback policies for robust control with constraints},
    journal = auto,
    volume = {42},
    number = {4},
    pages = {523-533},
    year = {2006},
    issn = {0005-1098},
    doi = {10.1016/j.automatica.2005.08.023},
    author = {Paul J. Goulart and Eric C. Kerrigan and Jan M. Maciejowski},
}

@INPROCEEDINGS{Tuna2006,
    author={Tuna, S. Emre and Messina, Michael J. and Teel, Andrew R.},
    booktitle=acc2006, 
    title={Shorter horizons for model predictive control}, 
    year={2006},
    volume={},
    number={},
    pages={6 pp.},
    organization={IEEE},
    doi={10.1109/ACC.2006.1655466}
}

@article{IlchRyan07,
    title={Tracking with prescribed transient behavior for nonlinear systems of known relative degree},
    author={Ilchmann, Achim and Ryan, Eugene P. and Townsend, Philip},
    journal=siamjco,
    volume={46},
    number={1},
    pages={210--230},
    year={2007},
    publisher={SIAM},
    doi={10.1137/050641946}
}

@Inbook{Besancon2007,
    author="Besan{\c{c}}on, Gildas",
    title={{An Overview on Observer Tools for Nonlinear Systems}},
    bookTitle={{Nonlinear Observers and Applications}},
    year="2007",
    publisher="Springer",
    address="Berlin, Heidelberg",
    pages="1--33",
    isbn="978-3-540-73503-8",
    doi={10.1007/978-3-540-73503-8_1},
}

@article{BECHLIOULIS2008,
    author={Bechlioulis, Charalampos P. and Rovithakis, George A.},
    journal=ieeetac, 
    title={{Robust Adaptive Control of Feedback Linearizable MIMO Nonlinear Systems With Prescribed Performance}}, 
    year={2008},
    volume={53},
    number={9},
    pages={2090-2099},
    doi={10.1109/TAC.2008.929402}
}

@article{Limon2008,
    title = {{On the design of Robust tube-based MPC for tracking}},
    journal = ifacprec,
    volume = {41},
    number = {2},
    pages = {15333-15338},
    year = {2008},
    note = ifcaworld17,
    issn = {1474-6670},
    doi = {10.3182/20080706-5-KR-1001.02593},
    author = {Daniel Limon and Ignacio Alvarado and Teodoro Alamo and Eduardo F. Camacho},
}

@inbook{Grune2008,
    author = {Lars Grüne and Karl Worthmann},
    title = {Sampled-data redesign for nonlinear multi-input systems},
    booktitle = {Geometric Control and Nonsmooth Analysis: In Honor of the 73rd Birthday of H. Hermes and of the 71st Birthday of RT Rockafellar},
    chapter = {},
    pages = {206-227},
    doi = {10.1142/9789812776075_0011},
    year = {2008},
    publisher={World Scientific}
}

@article{Grune2008b,
    title = {{Continuous-time controller redesign for digital implementation: A trajectory based approach}},
    journal = auto,
    volume = {44},
    number = {1},
    pages = {225-232},
    year = {2008},
    issn = {0005-1098},
    doi = {10.1016/j.automatica.2007.05.003},
    author = {Lars Grüne and Karl Worthmann and Dragan Nešić},
}

@book{schilders2008model,
    title={Model order reduction: theory, research aspects and applications},
    author={Schilders, Wilhelmus H.A. and Van der Vorst, Henk A. and Rommes, Joost},
    volume={13},
    year={2008},
    publisher={Springer},
    address={Heidelberg},
    doi={10.1007/978-3-540-78841-6}
}

@ARTICLE{Chen2008,
    author={Chen, Jie and Hara, Shinji and Qiu, Li and Middleton, Richard H.},
    journal=ieeetac, 
    title={{Best Achievable Tracking Performance in Sampled-Data Systems via LTI Controllers}}, 
    year={2008},
    volume={53},
    number={11},
    pages={2467-2479},
    doi={10.1109/TAC.2008.2006924}
}

@article{Schenato09,
    author={Schenato, Luca},
    journal=ieeetac, 
    title={To {Z}ero or to {H}old {C}ontrol {I}nputs {W}ith {L}ossy {L}inks?}, 
    year={2009},
    volume={54},
    number={5},
    pages={1093-1099},
    doi={10.1109/TAC.2008.2010999}
}

@article{mayne2009robust,
    title = {{Robust output feedback model predictive control of constrained linear systems: Time varying case}},
    journal = auto,
    volume = {45},
    number = {9},
    pages = {2082-2087},
    year = {2009},
    issn = {0005-1098},
    doi = {10.1016/j.automatica.2009.05.009},
    author = {David Q. Mayne and Saša V. Raković and Rolf Findeisen and Frank Allgöwer},
}

@article{gonzalez2009enlarging,
    title={Enlarging the domain of attraction of stable {MPC} controllers, maintaining the output performance},
    author={Gonz{\'a}lez, Alejandro H. and Odloak, Darci},
    journal=auto,
    volume={45},
    number={4},
    pages={1080--1085},
    year={2009},
    publisher={Elsevier},
    doi={10.1016/j.automatica.2008.11.015}
}

@article{Ilchmann01102009,
    author = {Achim Ilchmann and Eugene P. Ryan},
    title = {Performance funnels and tracking control},
    journal = intjoc,
    volume = {82},
    number = {10},
    pages = {1828--1840},
    year = {2009},
    publisher = {Taylor \& Francis},
    doi = {10.1080/00207170902777392}
}

@INPROCEEDINGS{LibeTren10,
    author={Liberzon, Daniel and Trenn, Stephan},
    booktitle=cdc49, 
    title={The bang-bang funnel controller}, 
    year={2010},
    volume={},
    number={},
    pages={690-695},
    organization={IEEE},
    doi={10.1109/CDC.2010.5717742}
}

@article{Limon2010,
    title = {Robust tube-based {MPC} for tracking of constrained linear systems with additive disturbances},
    journal = jprocescon,
    volume = {20},
    number = {3},
    pages = {248-260},
    year = {2010},
    issn = {0959-1524},
    doi = {10.1016/j.jprocont.2009.11.007},
    author = {Daniel Limon and Alvarado Alvarado and Teodoro Alamo and Eduardo F. Camacho},
}

@INPROCEEDINGS{KIM2010,
  author={Kim, Jung-Su},
  booktitle=iccas2010, 
  title={{Recent advances in adaptive MPC}}, 
  year={2010},
  volume={},
  number={},
  pages={218-222},
  doi={10.1109/ICCAS.2010.5669892}
}

@article{Adetola2011,
    author = {Adetola, Veronica and Guay, Martin},
    title = {{Robust adaptive MPC for constrained uncertain nonlinear systems}},
    journal = intjacsp,
    volume = {25},
    number = {2},
    pages = {155-167},
    doi = {10.1002/acs.1193},
    year = {2011}
}

@ARTICLE{bechlioulis2011,
    author={Bechlioulis, Charalampos P. and Rovithakis, George A.},
    journal=ieeetac, 
    title={{Robust Partial-State Feedback Prescribed Performance Control of Cascade Systems With Unknown Nonlinearities}}, 
    year={2011},
    volume={56},
    number={9},
    pages={2224-2230},
    doi={10.1109/TAC.2011.2157399}
}

@phdthesis{worthmann2011stability,
    title={Stability analysis of unconstrained receding horizon control schemes},
    author={Worthmann, Karl},
    year={2011},
    school    = {University of Bayreuth},
    address = {Bayreuth, Germany},
    url={https://epub.uni-bayreuth.de/id/eprint/273},
    urldate  = {2025-06-24},
}

@Article{Lee2011,
    author           = {Lee, Jay H.},
    year             = {2011},
    journal          = intjautomation,
    title            = {Model predictive control: Review of the three decades of development},
    doi              = {10.1007/s12555-011-0300-6},
    issn             = {2005-4092},
    number           = {3},
    pages            = {415--424},
    volume           = {9},
}

@article{IlchWirt13,
    title={On minimum phase},
    author={Ilchmann, Achim and Wirth, Fabian},
    journal=atauto,
    volume={61},
    number={12},
    pages={805--817},
    year={2013},
    publisher={De Gruyter},
    doi={10.1524/auto.2013.1002}
}

@article{Mezi13,
   author = "Mezić, Igor",
   title = {{Analysis of Fluid Flows via Spectral Properties of the Koopman Operator}}, 
   journal= annufluid,
   year = "2013",
   volume = "45",
   pages = "357-378",
   doi = {10.1146/annurev-fluid-011212-140652},
   publisher = "Annual Reviews",
   issn = "1545-4479",
   type = "Journal Article",
}

@Article{SeifBlaj13,
    AUTHOR = {Seifried, Robert and Blajer, Wojciech},
    TITLE = {Analysis of servo-constraint problems for underactuated multibody systems},
    JOURNAL = mecsci,
    VOLUME = {4},
    YEAR = {2013},
    NUMBER = {1},
    PAGES = {113--129},
    DOI = {10.5194/ms-4-113-2013}
}

@ARTICLE{calafiore2012robust,
    author={Calafiore, Giuseppe C. and Fagiano, Lorenzo},
    journal=ieeetac, 
    title={{Robust Model Predictive Control via Scenario Optimization}}, 
    year={2013},
    volume={58},
    number={1},
    pages={219-224},
    doi={10.1109/TAC.2012.2203054}
}

@article{Amrit2009,
    title = {Optimizing process economics online using model predictive control},
    journal = compchem,
    volume = {58},
    pages = {334-343},
    year = {2013},
    issn = {0098-1354},
    doi = {doi.org/10.1016/j.compchemeng.2013.07.015},
    author = {Rishi Amrit and James B. Rawlings and Lorenz T. Biegler},
}

@article{HackHopf13,
    AUTHOR    = {Hackl, Christoph M. and Hopfe, Norman and Ilchmann, Achim and Mueller, Markus and Trenn, Stephan},
    YEAR      = 2013,
    TITLE     = {Funnel control for systems with relative degree two},
    JOURNAL   = siamjco,
    Volume    = 51,
    Number     = 2,
    Pages     = {965--995},
    doi ={10.1137/100799903}
}

@article{LibeTren13b,
    author = {Liberzon, Daniel and Trenn, Stephan},
    title = {The bang-bang funnel controller for uncertain nonlinear systems with arbitrary relative degree},
    journal = ieeetac,
    year = {2013},
    volume = {58},
    number = {12},
    pages = {3126--3141},
    doi={10.1109/TAC.2013.2277631}
}

@article{yu2013tube,
    title = {{Tube MPC scheme based on robust control invariant set with application to Lipschitz nonlinear systems}},
    journal = scl,
    volume = {62},
    number = {2},
    pages = {194-200},
    year = {2013},
    issn = {0167-6911},
    doi = {10.1016/j.sysconle.2012.11.004},
    author = {Shuyou Yu and Christoph Maier and Hong Chen and Frank Allgöwer},
}

@article{Aswa13,
    title = {Provably safe and robust learning-based model predictive control},
    journal = auto,
    volume = {49},
    number = {5},
    pages = {1216-1226},
    year = {2013},
    issn = {0005-1098},
    doi = {10.1016/j.automatica.2013.02.003},
    author = {Anil Aswani and Humberto Gonzalez and S. Shankar Sastry and Claire Tomlin},
}

@article{BECHLIOULIS2014,
    title = {A low-complexity global approximation-free control scheme with prescribed performance for unknown pure feedback systems},
    journal = auto,
    volume = {50},
    number = {4},
    pages = {1217-1226},
    year = {2014},
    issn = {0005-1098},
    doi = {10.1016/j.automatica.2014.02.020},
    author = {Charalampos P. Bechlioulis and George A. Rovithakis}
}

@inproceedings{manrique2014mpc,
    title={{MPC} tracking under time-varying polytopic constraints for real-time applications},
    author={Manrique, Tatiana and Fiacchini, Mirko and Chambrion, Thomas and Mill{\'e}rioux, Gilles},
    booktitle=ecc2014,
    pages={1480--1485},
    year={2014},
    organization={IEEE},
    doi={10.1109/ECC.2014.6862584}
}

@ARTICLE{FaluMayn14,
    author={Falugi, Paola and Mayne, David Q.},
    journal=ieeetac, 
    title={{Getting Robustness Against Unstructured Uncertainty: A Tube-Based MPC Approach}}, 
    year={2014},
    volume={59},
    number={5},
    pages={1290-1295},
    doi={10.1109/TAC.2013.2287727}
}

@article{boccia2014stability,
    title={Stability and feasibility of state constrained {MPC} without stabilizing terminal constraints},
    author={Boccia, Andrea and Gr{\"u}ne, Lars and Worthmann, Karl},
    journal= scl,
    volume={72},
    pages={14--21},
    year={2014},
    publisher={Elsevier},
    doi={10.1016/j.sysconle.2014.08.002}
}

@INPROCEEDINGS{SenfPaug14,
   AUTHOR    = {Senfelds, Armands and Paugurs, Arturs},
   YEAR      = {2014},
   TITLE     = {Electrical drive {DC} link power flow control with adaptive approach},
   BOOKTITLE = rtucon14,
   Pages     = {30--33},
   doi       = {10.1109/RTUCON.2014.6998195}
}

@article{worthmann2014role,
    title={The role of sampling for stability and performance in unconstrained nonlinear model predictive control},
    author={Worthmann, Karl and Reble, Marcus and Grüne, Lars and Allgöwer, Frank},
    journal=siamjco,
    volume={52},
    number={1},
    pages={581--605},
    year={2014},
    publisher={SIAM},
    doi={10.1137/12086652X}
}

@article{limon2014periodic,
    title = {{MPC} for tracking periodic reference signals},
    journal = ifacprec,
    volume = {45},
    number = {17},
    pages = {490-495},
    year = {2012},
    note = facNonPC4,
    issn = {1474-6670},
    doi = {10.3182/20120823-5-NL-3013.00067},
    author = {Daniel Limon and Teodoro Alamo and David Muñoz {de la Peña} and Melanie N. Zeilinger and Colin N. Jones and Mario Pereira},
}

@article{garcia2015comprehensive,
    title={A comprehensive survey on safe reinforcement learning},
    author={Garc{\i}a, Javier and Fern{\'a}ndez, Fernando},
    journal=mlr,
    volume={16},
    number={1},
    pages={1437--1480},
    year={2015},
    urldate  = {2025-06-24},
    url  = {http://jmlr.org/papers/v16/garcia15a.html}
}

@INPROCEEDINGS{Giovanni2015,
    author={Mattei, Giovanni and Monaco, Salvatore and Normand-Cyrot, Dorothée},
    booktitle=ecc2015,
    title={Multi-rate sampled-data stabilization of a class of nonlinear systems}, 
    year={2015},
    volume={},
    number={},
    pages={975-980},
    doi={10.1109/ECC.2015.7330668}
}

@INPROCEEDINGS{HosseinNia2015,
    author={HosseinNia, S. Hassan},
    booktitle=ecc2015, 
    title={{Robust Model Predictive Control using Iterative Learning}}, 
    year={2015},
    volume={},
    number={},
    pages={3514-3519},
    doi={10.1109/ECC.2015.7331078}
}

@article{Kouvaritakis2015,
    author = {Kouvaritakis, Basil and Cannon, Mark},
    title = {{Developments in Robust and Stochastic Predictive Control in the Presence of Uncertainty}},
    journal = ascasme,
    volume = {1},
    number = {2},
    eid = {021003},
    year = {2015},
    month = {04},
    issn = {2332-9017},
    doi = {10.1115/1.4029744},
}

@INPROCEEDINGS{Worthmann2015,
    author={Worthmann, Karl and Reble, Marcus and Grüne, Lars and Allgöwer, Frank},
    booktitle=cdc54, 
    organization={IEEE},
    title={{Unconstrained nonlinear MPC: Performance estimates for sampled-data systems with zero order hold}}, 
    year={2015},
    volume={},
    number={},
    pages={4971-4976},
    doi={10.1109/CDC.2015.7402996}
}

@article{williams:kevrekidis:rowley:2015,
    author = {Williams, Matthew O. and Kevrekidis, Ioannis G. and Rowley, Clarence W.},
    year = {2015},
    pages={1307--1346},
    title = {{A Data-Driven Approximation of the {Koopman} Operator: Extending Dynamic Mode Decomposition}},
    volume = {25},
    number={6},
    journal = journonlin,
    doi ={10.1007/s00332-015-9258-5}
}

@article{amodei2016concrete,
    title={Concrete problems in {AI} safety},
    author={Amodei, Dario and Olah, Chris and Steinhardt, Jacob and Christiano, Paul and Schulman, John and Man{\'e}, Dan},
    journal={arXiv preprint arXiv:1606.06565},
    year={2016},
    doi = {10.48550/arXiv.1606.06565}
}

@INPROCEEDINGS{aydiner2016periodic,
  author={Aydiner, Emre and M{\"u}ller, Matthias A. and Allg{\"o}wer, Frank},
  booktitle=ecc2016, 
  title={Periodic reference tracking for nonlinear systems via model predictive control}, 
  year={2016},
  pages={2602-2607},
  doi={10.1109/ECC.2016.7810682}
}

@ARTICLE{CairBorr16,
    author={Di Cairano, Stefano and Borrelli, Francesco},
    journal=ieeetac, 
    title={{Reference Tracking With Guaranteed Error Bound for Constrained Linear Systems}}, 
    year={2016},
    volume={61},
    number={8},
    pages={2245-2250},
    doi={10.1109/TAC.2015.2491738}
}

@ARTICLE{mesbah2016,
    author={Mesbah, Ali},
    journal=ieeeconsysMag, 
    title={{Stochastic Model Predictive Control: An Overview and Perspectives for Future Research}}, 
    year={2016},
    volume={36},
    number={6},
    pages={30-44},
    doi={10.1109/MCS.2016.2602087}
}

@INPROCEEDINGS{Feller2016,
    author={Feller, Christian and Ouerghi, Meriam and Ebenbauer, Christian},
    booktitle=cdc55, 
    organization={IEEE},
    title={Robust output feedback model predictive control based on relaxed barrier functions}, 
    year={2016},
    volume={},
    number={},
    pages={1477-1483},
    doi={10.1109/CDC.2016.7798475}
}

@article{brunton2016discovering,
    author = {Steven L. Brunton  and Joshua L. Proctor  and J. Nathan Kutz },
    title = {Discovering governing equations from data by sparse identification of nonlinear dynamical systems},
    journal = {Proceedings of the National Academy of Sciences},
    volume = {113},
    number = {15},
    pages = {3932-3937},
    year = {2016},
    doi = {10.1073/pnas.1517384113},
}

@BOOK{rawlings2017model,
    title={Model predictive control: theory, computation, and design},
    author={Rawlings, James Blake and Mayne, David Q. and Diehl, Moritz},
    volume={2},
    year={2017},
    publisher={Nob Hill Publishing Madison, WI},
    isbn={978-0975937785}
}

@article{kogel2017robust,
    title = {{Robust output feedback MPC for uncertain linear systems with reduced conservatism}},
    journal = ifaconline,
    volume = {50},
    number = {1},
    pages = {10685-10690},
    year = {2017},
    note = ifcaworld20,
    issn = {2405-8963},
    doi = {10.1016/j.ifacol.2017.08.2186},
    author = {Markus Kögel and Rolf Findeisen},
    }

@ARTICLE{Ames2017,
    author={Ames, Aaron D. and Xu, Xiangru and Grizzle, Jessy W. and Tabuada, Paulo},
    journal=ieeetac, 
    title={{Control Barrier Function Based Quadratic Programs for Safety Critical Systems}}, 
    year={2017},
    volume={62},
    number={8},
    pages={3861-3876},
    doi={10.1109/TAC.2016.2638961}
}

@book{GrunPann17,
  title = {{Nonlinear Model Predictive Control: Theory and Algorithms}},
  isbn = {978-0-85729-500-2},
  author={Gr{\"u}ne, Lars and Pannek, J{\"u}rgen},
  pages={45--69},
  Address = {London},
  year={2017},
  publisher={Springer},
  doi = {10.1007/978-0-85729-501-9}
}

@INPROCEEDINGS{singh2017robust,
    author={Singh, Sumeet and Majumdar, Anirudha and Slotine, Jean-Jacques and Pavone, Marco},
    booktitle=ieeeicra17, 
    title={Robust online motion planning via contraction theory and convex optimization}, 
    year={2017},
    volume={},
    number={},
    pages={5883-5890},
    doi={10.1109/ICRA.2017.7989693}
}

@book{Hackl17,
    title={{Non-identifier based adaptive control in mechatronics: Theory and Application}},
    author={Hackl, Christoph M.},
    volume={466},
    year={2017},
    publisher={Springer},
    doi={10.1007/978-3-319-55036-7},
    isbn={978-3-319-55034-3}
}

@article{ramachandran2017searching,
    title={Searching for activation functions},
    author={Ramachandran, Prajit and Zoph, Barret and Le, Quoc V.},
    journal={arXiv preprint arXiv:1710.05941},
    year={2017},
    doi = {10.48550/arXiv.1710.05941}
}

@ARTICLE{Brunner2017,
    author={Brunner, Florian David and Heemels, W.P. Maurice H. and Allgöwer, Frank},
    journal=ieeetac, 
    title={{Robust Event-Triggered MPC With Guaranteed Asymptotic Bound and Average Sampling Rate}}, 
    year={2017},
    volume={62},
    number={11},
    pages={5694-5709},
    doi={10.1109/TAC.2017.2702646}
}

@INPROCEEDINGS{bansal2017hamilton,
    author={Bansal, Somil and Chen, Mo and Herbert, Sylvia and Tomlin, Claire J.},
    booktitle=cdc56, 
    organization={IEEE},
    title={{Hamilton-Jacobi reachability: A brief overview and recent advances}}, 
    year={2017},
    volume={},
    number={},
    pages={2242-2253},
    doi={10.1109/CDC.2017.8263977}
}

@article{schulman2017proximal,
    title={Proximal policy optimization algorithms},
    author={Schulman, John and Wolski, Filip and Dhariwal, Prafulla and Radford, Alec and Klimov, Oleg},
    journal={arXiv preprint arXiv:1707.06347},
    year={2017},
    doi = {10.48550/arXiv.1707.06347}
}

@article{limon2018nonlinear,
    title={{Nonlinear MPC for Tracking Piece-Wise Constant Reference Signals}},
    author={Limon, Daniel and Ferramosca, Antonio and Alvarado, Ignacio and Alamo, Teodoro},
    journal=ieeetac,
    volume={63},
    number={11},
    pages={3735--3750},
    year={2018},
    publisher={IEEE},
    doi={10.1109/TAC.2018.2798803}
}

@inproceedings{amos2018differentiable,
    author = {Amos, Brandon and Jimenez, Ivan and Sacks, Jacob and Boots, Byron and Kolter, J. Zico},
    booktitle = adneural,
    pages = {},
    publisher = {Curran Associates, Inc.},
    title = {{Differentiable MPC for End-to-end Planning and Control}},
    url = {https://proceedings.neurips.cc/paper_files/paper/2018/file/ba6d843eb4251a4526ce65d1807a9309-Paper.pdf},
    urldate  = {2025-06-24},
    volume = {31},
    year = {2018}
}

@PhdThesis{Hoang18,
    author    = {L{\^e}, Huy Ho{\`a}ng},
    school    = {Universität Hamburg},
    address   = {Hamburg, Germany},
    title     = {Funnel control for systems with known vector relative degree},
    year      = {2018},
    url       = {https://ediss.sub.uni-hamburg.de/handle/ediss/8142},
    urldate   = {2025-06-24}, 
}

@article{BergLe18a,
    title={Funnel control for nonlinear systems with known strict relative degree},
    author={Berger, Thomas and L{\^e}, Huy Ho{\`a}ng and Reis, Timo},
    journal=auto,
    volume={87},
    pages={345--357},
    year={2018},
    publisher={Elsevier},
    doi ={10.1016/j.automatica.2017.10.017}
}

@article{BergReis18,
    author = {Berger, Thomas and Reis, Timo},
    title = {The funnel pre-compensator},
    journal = introbust,
    volume = {28},
    number = {16},
    pages = {4747-4771},
    doi = {10.1002/rnc.4281},
    year = {2018}
}

@INPROCEEDINGS{Wabersich18,
    author={Wabersich, Kim Peter and Zeilinger, Melanie N.},
    booktitle=cdc57, 
    organization={IEEE},
    title={{Linear Model Predictive Safety Certification for Learning-Based Control}}, 
    year={2018},
    volume={},
    number={},
    pages={7130-7135},
    doi={10.1109/CDC.2018.8619829}
}

@article{korda:mezic:2018b,
    title={On convergence of extended dynamic mode decomposition to the {Koopman} operator},
    author={Korda, Milan and Mezi{\'c}, Igor},
    journal=journonlin,
    volume={28},
    number={2},
    pages={687--710},
    year={2018},
    publisher={Springer},
    doi ={10.1007/s00332-017-9423-0} 
}

@article{chen2018hamilton,
    title={{Hamilton--Jacobi reachability: Some recent theoretical advances and applications in unmanned airspace management}},
    author={Chen, Mo and Tomlin, Claire J.},
    journal=annurevcon,
    volume={1},
    pages={333--358},
    year={2018},
    publisher={Annual Reviews},
    doi ={10.1146/annurev-control-060117-104941}
}

@book{Sutton2018,
    author    = {Sutton, Richard S. and Barto, Andrew G.},
    title     = {{Reinforcement Learning, second edition: An Introduction}},
    publisher = {A Bradford Book},
    address   = {Cambridge, MA, USA},
    year      = {2018},
    isbn      = {9780262039246}
}

@article{Gesser2018,
    title = {{Robust Model Predictive Control: Implementation Issues with Comparative Analysis}},
    journal = ifaconline,
    volume = {51},
    number = {25},
    pages = {478-483},
    year = {2018},
    note = {9th IFAC Symposium on Robust Control Design ROCOND 2018},
    issn = {2405-8963},
    doi = {10.1016/j.ifacol.2018.11.183},
    author = {Rodrigo S. Gesser and Daniel M. Lima and Julio E. Normey-Rico},
}

@article{KordMezi18:MPC,
    title = {Linear predictors for nonlinear dynamical systems: Koopman operator meets model predictive control},
    journal = auto,
    volume = {93},
    pages = {149-160},
    year = {2018},
    issn = {0005-1098},
    doi = {10.1016/j.automatica.2018.03.046},
    author = {Milan Korda and Igor Mezić},
}

@INPROCEEDINGS{Marvi2019,
    author={Marvi, Zahra and Kiumarsi, Bahare},
    booktitle=cavs2019,
    title={{Safety Planning Using Control Barrier Function: A Model Predictive Control Scheme}}, 
    year={2019},
    volume={},
    number={},
    pages={1-5},
    doi={10.1109/CAVS.2019.8887800}
}

@Article{Andersson2019,
    author = {Joel A.E. Andersson and Joris Gillis and Greg Horn
            and James B. Rawlings and Moritz Diehl},
    title = {{CasADi} -- {A} software framework for nonlinear optimization
           and optimal control},
    journal = mathprog,
    volume = {11},
    number = {1},
    pages = {1--36},
    year = {2019},
    publisher = {Springer},
    doi = {10.1007/s12532-018-0139-4}
}

@article{Wu2019,
    title = {{Control Lyapunov-Barrier function-based model predictive control of nonlinear systems}},
    journal = auto,
    volume = {109},
    eid = {108508},
    year = {2019},
    issn = {0005-1098},
    doi = {10.1016/j.automatica.2019.108508},
    author = {Zhe Wu and Fahad Albalawi and Zhihao Zhang and Junfeng Zhang and Helen Durand and Panagiotis D. Christofides},
}

@article{kohler2019nonlinear,
    title={{Nonlinear Reference Tracking: An Economic Model Predictive Control Perspective}},
    author={K{\"o}hler, Johannes and M{\"u}ller, Matthias A. and Allg{\"o}wer, Frank},
    journal=ieeetac,
    volume={64},
    number={1},
    pages={254--269},
    year={2019},
    publisher={IEEE},
    doi={10.1109/TAC.2018.2800789}
}

@inproceedings{ames2019control,
    title={Control barrier functions: {T}heory and applications},
    author={Ames, Aaron D. and Coogan, Samuel and Egerstedt, Magnus and Notomista, Gennaro and Sreenath, Koushil and Tabuada, Paulo},
    booktitle=ecc2019,
    pages={3420--3431},
    year={2019},
    organization={IEEE},
    doi={10.23919/ECC.2019.8796030}
}

@INPROCEEDINGS{Schussler2019,
    author={Schüssler, Max and Münker, Tobias and Nelles, Oliver},
    booktitle=ieeessci, 
    title={{Deep Recurrent Neural Networks for Nonlinear System Identification}}, 
    year={2019},
    volume={},
    number={},
    pages={448-454},
    doi={10.1109/SSCI44817.2019.9003133}
}

@article{berger2019learningbased,
    title={Learning-based {Funnel-MPC} for output-constrained nonlinear systems},
    author={Berger, Thomas and K{\"a}stner, Carolin and Worthmann, Karl},
    journal=ifaconline,
    volume={53},
    number={2},
    pages={5177--5182},
    year={2020},
    publisher={Elsevier},
    doi={10.1016/j.ifacol.2020.12.1186}
}

@ARTICLE{Singh2019,
  author={Singh, Sumeet and Chow, Yinlam and Majumdar, Anirudha and Pavone, Marco},
  journal=ieeetac, 
  title={{A Framework for Time-Consistent, Risk-Sensitive Model Predictive Control: Theory and Algorithms}}, 
  year={2019},
  volume={64},
  number={7},
  pages={2905-2912},
  doi={10.1109/TAC.2018.2874704}
}

@article{jang2019q,
    title={{Q-learning algorithms: A comprehensive classification and applications}},
    author={Jang, Beakcheol and Kim, Myeonghwi and Harerimana, Gaspard and Kim, Jong Wook},
    journal=ieeeaccess,
    volume={7},
    pages={133653--133667},
    year={2019},
    publisher={IEEE},
    doi={10.1109/ACCESS.2019.2941229}
}

@INPROCEEDINGS{Lopez2019,
    author={Lopez, Brett T. and Slotine, Jean-Jacques E. and How, Jonathan P.},
    booktitle=acc2019, 
    title={{Dynamic Tube MPC for Nonlinear Systems}}, 
    year={2019},
    volume={},
    number={},
    pages={1655-1662},
    organization={IEEE},
    doi={10.23919/ACC.2019.8814758}
}

@inproceedings{coulson2019data,
    title={{Data-enabled predictive control: In the shallows of the DeePC}},
    author={Coulson, Jeremy and Lygeros, John and D{\"o}rfler, Florian},
    booktitle=ecc2019,
    pages={307--312},
    year={2019},
    organization={IEEE},
    doi={10.23919/ECC.2019.8795639}
}

@article{LoreCann19,
    title = {{Robust MPC with recursive model update}},
    journal = auto,
    volume = {103},
    pages = {461-471},
    year = {2019},
    issn = {0005-1098},
    doi = {10.1016/j.automatica.2019.02.023},
    author = {Matthias Lorenzen and Mark Cannon and Frank Allgöwer},
}

@article{Raissi2019,
    title = {{Physics-informed neural networks: A deep learning framework for solving forward and inverse problems involving nonlinear partial differential equations}},
    journal = jourphysic,
    volume = {378},
    pages = {686-707},
    year = {2019},
    issn = {0021-9991},
    doi = {10.1016/j.jcp.2018.10.045},
    author = {M. Raissi and P. Perdikaris and G.E. Karniadakis},
}

@ARTICLE{SunDai19,
    author={Sun, Zhongqi and Dai, Li and Liu, Kun and Dimarogonas, Dimos V. and Xia, Yuanqing},
    journal=ieeetac, 
    title={{Robust Self-Triggered MPC With Adaptive Prediction Horizon for Perturbed Nonlinear Systems}}, 
    year={2019},
    volume={64},
    number={11},
    pages={4780-4787},
    doi={10.1109/TAC.2019.2905223}
}

@article{BergOtto19,
    title={Combined open-loop and funnel control for underactuated multibody systems},
    author={Berger, Thomas and Otto, Svenja and Reis, Timo and Seifried, Robert},
    journal=nonlindyn,
    volume={95},
    pages={1977--1998},
    year={2019},
    publisher={Springer},
    doi = {10.1007/s11071-018-4672-5}
}

@INPROCEEDINGS{yuan2019bounded,
    author={Yuan, Meng and Manzie, Chris and Good, Malcolm and Shames, Iman and Keynejad, Farzad and Robinette, Troy},
    booktitle=ieeeicit2019, 
    title={{Bounded Error Tracking Control for Contouring Systems with End Effector Measurements}}, 
    year={2019},
    volume={},
    number={},
    pages={66-71},
    doi={10.1109/ICIT.2019.8755064}
}

@article{HewingWaber20,
    author = {Hewing, Lukas and Wabersich, Kim P. and Menner, Marcel and Zeilinger, Melanie N.},
    title = {{Learning-Based Model Predictive Control: Toward Safe Learning in Control}}, 
    journal = annurevcon,
    volume = {3},
    number = {1},
    pages = {269-296},
    year = {2020},
    publisher = {Annual Reviews},
    issn = {2573-5144},
    doi = {10.1146/annurev-control-090419-075625}
}

@ARTICLE{hewing2019cautious,
    author={Hewing, Lukas and Kabzan, Juraj and Zeilinger, Melanie N.},
    journal=ieeetcst, 
    title={{Cautious Model Predictive Control Using Gaussian Process Regression}}, 
    year={2020},
    volume={28},
    number={6},
    pages={2736-2743},
    doi={10.1109/TCST.2019.2949757}
}

@article{BergPuch20,
    author       = {Thomas Berger and Marc Puche and Felix L. Schwenninger},
    year         = {2020},
    journal      = scl,
    title        = {Funnel control in the presence of infinite-dimensional internal dynamics},
    eid          = {104678},
    volume       = {139},
    doi ={10.1016/j.sysconle.2020.104678}
}

@article{van2020data,
    title={Data informativity: a new perspective on data-driven analysis and control},
    author={Van Waarde, Henk J. and Eising, Jaap and Trentelman, Harry L. and Camlibel, M. Kanat},
    journal=ieeetac,
    volume={65},
    number={11},
    pages={4753--4768},
    year={2020},
    publisher={IEEE},
    doi={10.1109/TAC.2020.2966717}
}

@article{Cao2020,
    title = {{Deep Neural Network Approximation of Nonlinear Model Predictive Control}},
    journal = {IFAC-PapersOnLine},
    volume = {53},
    number = {2},
    pages = {11319-11324},
    year = {2020},
    note = ifcaworld21,
    issn = {2405-8963},
    doi = {10.1016/j.ifacol.2020.12.538},
    author = {Yankai Cao and R. Bhushan Gopaluni},
}

@article{samad2020industry,
    title={{Industry engagement with control research: Perspective and messages}},
    author={Samad, Tariq and Bauer, Margret and Bortoff, Scott and Di Cairano, Stefano and Fagiano, Lorenzo and Odgaard, Peter Fogh and Rhinehart, R. Russell and S{\'a}nchez-Pe{\~n}a, Ricardo and Serbezov, Atanas and Ankersen, Finn and others},
    journal=arcontrol,
    volume={49},
    pages={1--14},
    year={2020},
    publisher={Elsevier},
    doi = {10.1016/j.arcontrol.2020.03.002}
}

@article{CoroGrun20,
    title={Model predictive control, cost controllability, and homogeneity},
    author={Coron, Jean-Michel and Gr{\"u}ne, Lars and Worthmann, Karl},
    journal=siamjco,
    volume={58},
    number={5},
    pages={2979--2996},
    year={2020},
    publisher={SIAM},
    doi={10.1137/19M1265995}
}

@ARTICLE{Dimanidis20,
  author={Dimanidis, Ioannis S. and Bechlioulis, Charalampos P. and Rovithakis, George A.},
  journal=ieeetac, 
  title={{Output Feedback Approximation-Free Prescribed Performance Tracking Control for Uncertain MIMO Nonlinear Systems}}, 
  year={2020},
  volume={65},
  number={12},
  pages={5058-5069},
  doi={10.1109/TAC.2020.2970003}
}

@INPROCEEDINGS{Yang2020,
    author={Yang, Guang and Belta, Calin and Tron, Roberto},
    booktitle=acc2020, 
    title={Continuous-time Signal Temporal Logic Planning with Control Barrier Functions}, 
    year={2020},
    volume={},
    number={},
    pages={4612-4618},
    organization={IEEE},
    doi={10.23919/ACC45564.2020.9147387}
}

@ARTICLE{koehler2020,
    author={Köhler, Johannes and Müller, Matthias A. and Allgöwer, Frank},
    journal=ieeetac, 
    title={{A Nonlinear Model Predictive Control Framework Using Reference Generic Terminal Ingredients}}, 
    year={2020},
    volume={65},
    number={8},
    pages={3576-3583},
    doi={10.1109/TAC.2019.2949350}
}

@article{KohlSolo20,
  author={Köhler, Johannes and Soloperto, Raffaele and Müller, Matthias A. and Allgöwer, Frank},
  journal=ieeetac, 
  title={{A Computationally Efficient Robust Model Predictive Control Framework for Uncertain Nonlinear Systems}}, 
  year={2021},
  volume={66},
  number={2},
  pages={794-801},
  doi={10.1109/TAC.2020.2982585}
}

@InProceedings{Berger2020Vector,
    author={Berger, Thomas
    and L{\^e}, Huy Ho{\`a}ng
    and Reis, Timo},
    title={{Vector Relative Degree and Funnel Control for Differential-Algebraic Systems}},
    booktitle={{Progress in Differential-Algebraic Equations II}},
    year="2020",
    publisher={Springer},
    address="Cham",
    pages="213--255",
    isbn="978-3-030-53905-4",
    doi ={10.1007/978-3-030-53905-4_8},
}

@article{EsteWort21,
    title={{Recursive Feasibility of Continuous-Time Model Predictive Control Without Stabilising Constraints}},
    author={Esterhuizen, Willem and Worthmann, Karl and Streif, Stefan},
    journal=ieeecsl,
    volume={5},
    number={1},
    pages={265--270},
    year={2021},
    publisher={IEEE},
    doi={10.1109/LCSYS.2020.3001514}
}

@Inbook{Heemels2021,
    title={Event-triggered and self-triggered control},
    author={Heemels, W.P. Maurice H. and Johansson, Karl H. and Tabuada, Paulo},
    booktitle={Encyclopedia of Systems and Control},
    pages={724--730},
    year={2021},
    publisher={Springer},
    doi ={10.1007/978-3-030-44184-5}
}

@article{BergIlch21,
    AUTHOR    = {Berger, Thomas and Ilchmann, Achim and Ryan, Eugene P.},
    YEAR      = 2021,
    TITLE     = {Funnel control of nonlinear systems},
    JOURNAL   = mcss,
    volume    = 33,
    pages     = {151--194},
    doi ={10.1007/s00498-021-00277-z}
}

@INPROCEEDINGS{Pfitz2021,
    author={Pfitz, Fabian and Hu, Xuyan and Ebenbauer, Christian},
    booktitle=med2021, 
    title={{Relaxed Barrier MPC for Path Following in Constrained Environments}}, 
    year={2021},
    volume={},
    number={},
    pages={872-879},
    doi={10.1109/MED51440.2021.9480329}
}

@article{maiworm2021online,
    author = {Maiworm, Michael and Limon, Daniel and Findeisen, Rolf},
    title = {Online learning-based model predictive control with {Gaussian} process models and stability guarantees},
    journal = introbust,
    volume = {31},
    number = {18},
    pages = {8785-8812},
    doi = {10.1002/rnc.5361},
    year = {2021}
}

@article{Wabersich21,
    title = {A predictive safety filter for learning-based control of constrained nonlinear dynamical systems},
    journal = auto,
    volume = {129},
    pages = {109597},
    year = {2021},
    issn = {0005-1098},
    author = {Kim P. Wabersich and Melanie N. Zeilinger},
    doi = {10.1016/j.automatica.2021.109597}
}

@article{MarkDorf21,
    title={Behavioral systems theory in data-driven analysis, signal processing, and control},
    author={Markovsky, Ivan and D{\"o}rfler, Florian},
    journal=arcontrol,
    volume={52},
    pages={42--64},
    year={2021},
    publisher={Elsevier},
    doi = {10.1016/j.arcontrol.2021.09.005}
}

@ARTICLE{Zhuang2021,
    author={Zhuang, Fuzhen and Qi, Zhiyuan and Duan, Keyu and Xi, Dongbo and Zhu, Yongchun and Zhu, Hengshu and Xiong, Hui and He, Qing},
    journal=ieeeproc, 
    title={{A Comprehensive Survey on Transfer Learning}}, 
    year={2021},
    volume={109},
    number={1},
    pages={43-76},
    doi={10.1109/JPROC.2020.3004555}
}

@ARTICLE{Ma2021,
    author={Ma, Lele and Liu, Xiangjie and Kong, Xiaobing and Lee, Kwang Y.},
    journal=ieeetnnls, 
    title={{Iterative Learning Model Predictive Control Based on Iterative Data-Driven Modeling}}, 
    year={2021},
    volume={32},
    number={8},
    pages={3377-3390},
    doi={10.1109/TNNLS.2020.3016295}
}

@article{BergDrue21,
    title={Tracking control for underactuated non-minimum phase multibody systems},
    author={Berger, Thomas and Dr{\"u}cker, Svenja and Lanza, Lukas and Reis, Timo and Seifried, Robert},
    journal=nonlindyn,
    volume={104},
    number={4},
    pages={3671--3699},
    year={2021},
    publisher={Springer},
    doi = {10.1007/s11071-021-06458-4}
}

@article{Lanz21,
    title={Internal dynamics of multibody systems},
    author={Lanza, Lukas},
    journal=scl,
    volume={152},
    eid={104931},
    year={2021},
    publisher={Elsevier},
    issn = {0167-6911},
    doi = {10.1016/j.sysconle.2021.104931}
}

@article{OttoRowl21,
   author = "Otto, Samuel E. and Rowley, Clarence W.",
   title = {{Koopman Operators for Estimation and Control of Dynamical Systems}}, 
   journal= annurevcon,
   year = "2021",
   volume = "4",
   pages = "59-87",
   doi = {10.1146/annurev-control-071020-010108},
   publisher = "Annual Reviews",
   issn = "2573-5144",
   type = "Journal Article",
  }

@phdthesis{zanelli2021inexact,
    title={Inexact methods for nonlinear model predictive control},
    author={Zanelli, Andrea},
    year={2021},
    school={University of Freiburg},
    address = {Freiburg, Germany},
    doi = {10.6094/UNIFR/224926}
}

@ARTICLE{BerbKoeh20,
    author={Berberich, Julian and Köhler, Johannes and Müller, Matthias A. and Allgöwer, Frank},
    journal=ieeetac, 
    title={{Data-Driven Model Predictive Control With Stability and Robustness Guarantees}}, 
    year={2021},
    volume={66},
    number={4},
    pages={1702-1717},
    doi={10.1109/TAC.2020.3000182}
}

@inproceedings{Buja21,
    author={Bujarbaruah, Monimoy and Rosolia, Ugo and Stürz, Yvonne R. and Borrelli, Francesco},
    booktitle=acc2021,
    title={{A Simple Robust MPC for Linear Systems with Parametric and Additive Uncertainty}}, 
    year={2021},
    volume={},
    number={},
    pages={2108-2113},
    organization={IEEE},
    doi={10.23919/ACC50511.2021.9482957}
}

@ARTICLE{Petsag2021,
    author={Petsagkourakis, Panagiotis and Heath, William P. and Carrasco, Joaquin and Theodoropoulos, Constantinos},
    journal=ieeetac, 
    title={{Robust Stability of Barrier-Based Model Predictive Control}}, 
    year={2021},
    volume={66},
    number={4},
    pages={1879-1886},
    doi={10.1109/TAC.2020.3010770}
}

@ARTICLE{GoswPale21,
    author={Goswami, Debdipta and Paley, Derek A.},
    journal=ieeetac, 
    title={{Bilinearization, Reachability, and Optimal Control of Control-Affine Nonlinear Systems: A Koopman Spectral Approach}}, 
    year={2022},
    volume={67},
    number={6},
    pages={2715-2728},
    doi={10.1109/TAC.2021.3088802}
}

@ARTICLE{Geromel2022,
  author={Geromel, José C.},
  journal=ieeetac, 
  title={{Sampled-Data Model Predictive Control}}, 
  year={2022},
  volume={67},
  number={5},
  pages={2466-2472},
  doi={10.1109/TAC.2021.3077353}
}

@article{kohler2022constrained,
    author={Köhler, Johannes and Müller, Matthias A. and Allgöwer, Frank},
    journal=ieeetac, 
    title={{Constrained Nonlinear Output Regulation Using Model Predictive Control}}, 
    year={2022},
    volume={67},
    number={5},
    pages={2419-2434},
    doi={10.1109/TAC.2021.3081080}
}

@ARTICLE{Breeden2022,
    author={Breeden, Joseph and Garg, Kunal and Panagou, Dimitra},
    journal=ieeecsl, 
    title={Control Barrier Functions in Sampled-Data Systems}, 
    year={2022},
    volume={6},
    number={},
    pages={367-372},
    doi={10.1109/LCSYS.2021.3076127}
}

@PhdThesis{Druecker22,
    author    = {Svenja Drücker},
    school    = {Technische Universität Hamburg},
    address   = {Hamburg, Germany},
    title     = {Servo-constraints for inversion of underactuated multibody systems},
    year      = {2022},
    series    = {MuM Notes in Mechanics and Dynamics},
    doi       = {10.15480/882.4089} 
}

@article{BergPuch22,
    Author = {Berger, Thomas and Puche, Marc and Schwenninger, Felix L.},
    Title = {Funnel control for a moving water tank},
    journal = auto,
    volume = 135,
    eid = {109999},
    Year = 2022,
    doi ={10.1016/j.automatica.2021.109999}
}

@article{brunke2022safe,
    title={Safe learning in robotics: {F}rom learning-based control to safe reinforcement learning},
    author={Brunke, Lukas and Greeff, Melissa and Hall, Adam W. and Yuan, Zhaocong and Zhou, Siqi and Panerati, Jacopo and Schoellig, Angela P.},
    journal=annurevcon,
    volume={5},
    pages={411--444},
    year={2022},
    publisher={Annual Reviews},
    doi = {10.1146/annurev-control-042920-020211}
}

@article{lanza2022output,
    title={Output feedback control with prescribed performance via funnel pre-compensator},
    author={Lanza, Lukas},
    journal=mcss,
    volume={34},
    number={4},
    pages={715--758},
    year={2022},
    publisher={Springer},
    doi ={10.1007/s00498-022-00322-5}
}

@INPROCEEDINGS{zieger2022non,
    author={Zieger, Tim and Nguyen, Hoang Hai and Schulz, Erik and Oehlschlägel, Thimo and Findeisen, Rolf},
    booktitle=cdc61, 
    organization={IEEE},
    title={Non-diverging neural networks supported tube model predictive control}, 
    year={2022},
    volume={},
    number={},
    pages={3066-3073},
    doi={10.1109/CDC51059.2022.9993089}
}

@article{tambon2022certify,
    title={How to certify machine learning based safety-critical systems? {A} systematic literature review},
    author={Tambon, Florian and Laberge, Gabriel and An, Le and Nikanjam, Amin and Mindom, Paulina Stevia Nouwou and Pequignot, Yann and Khomh, Foutse and Antoniol, Giulio and Merlo, Ettore and Laviolette, Fran{\c{c}}ois},
    journal=autosofteng,
    volume={29},
    number={2},
    pages={38},
    year={2022},
    publisher={Springer},
    doi = {10.1007/s10515-022-00337-x}
}

@article{BrunKutz22,
    author = {Brunton, Steven L. and Budi\v{s}i\'{c}, Marko and Kaiser, Eurika and Kutz, J. Nathan},
    title = {{Modern Koopman Theory for Dynamical Systems}},
    journal = siamrev,
    volume = {64},
    number = {2},
    pages = {229-340},
    year = {2022},
    doi = {10.1137/21M1401243},
}

@ARTICLE{RakoDai22,
    author={Raković, Saša V. and Dai, Li and Xia, Yuanqing},
    journal=ieeetac,
    title={{Homothetic Tube Model Predictive Control for Nonlinear Systems}}, 
    year={2023},
    volume={68},
    number={8},
    pages={4554-4569},
    doi={10.1109/TAC.2022.3207415}
}

@INPROCEEDINGS{Korder2022,
    author={Korder, Kristina and Noack, Matti and Reger, Johann},
    booktitle=cdc61, 
    organization={IEEE},
    title={{Non-asymptotic Observer Design for Nonlinear Systems Based on Linearization}}, 
    year={2022},
    volume={},
    number={},
    pages={615-621},
    doi={10.1109/CDC51059.2022.9992644}
}

@article{shaheen2022continual,
    title={Continual learning for real-world autonomous systems: Algorithms, challenges and frameworks},
    author={Shaheen, Khadija and Hanif, Muhammad Abdullah and Hasan, Osman and Shafique, Muhammad},
    journal=introbosys,
    volume={105},
    number={1},
    pages={9},
    year={2022},
    publisher={Springer},
    doi = {10.1007/s10846-022-01603-6}
}

@inproceedings{HuTren22,
    title = {Funnel control for relative degree one nonlinear systems with input saturation},
    author = {Jiaming Hu and Stephan Trenn and Xiaojin Zhu},
    year = {2022},
    booktitle = ecc2022,
    pages = {227-232},
    address = {London},
    doi = {10.23919/ECC55457.2022.9837979}
}

@article{berberich2022linear,
    title={{Linear tracking {MPC} for nonlinear systems—Part {II}: The data-driven case}},
    author={Berberich, Julian and K{\"o}hler, Johannes and M{\"u}ller, Matthias A. and Allg{\"o}wer, Frank},
    journal=ieeetac,
    volume={67},
    number={9},
    pages={4406--4421},
    year={2022},
    publisher={IEEE},
    doi={10.1109/TAC.2022.3166851}
}

@article{BuPPC23,
    author = {Bu, Xiangwei},
    title = {Prescribed performance control approaches, applications and challenges: {A} comprehensive survey},
    journal = asiacont,
    volume = {25},
    number = {1},
    pages = {241-261},
    doi = {10.1002/asjc.2765},
    year = {2023}
}

@ARTICLE{Yin2023,
    author={Yin, Ji and Dawson, Charles and Fan, Chuchu and Tsiotras, Panagiotis},
    journal={IEEE Robotics and Automation Letters}, 
    title={{Shield Model Predictive Path Integral: A Computationally Efficient Robust MPC Method Using Control Barrier Functions}}, 
    year={2023},
    volume={8},
    number={11},
    pages={7106-7113},
    doi={10.1109/LRA.2023.3315211}
}

@article{sakamoto2023,
    title={When does stabilizability imply the existence of infinite horizon optimal control in nonlinear systems?},
    author={Sakamoto, Noboru},
    journal=auto,
    volume={147},
    eid ={110706},
    year={2023},
    publisher={Elsevier},
    doi = {10.1016/j.automatica.2022.110706}
}

@Article{BergIlch23,
    title = {{Funnel control — A survey}},
    journal = arcontrol,
    volume = {60},
    pages = {101024},
    year = {2025},
    issn = {1367-5788},
    doi = {10.1016/j.arcontrol.2025.101024},
    author = {Thomas Berger and Achim Ilchmann and Eugene P. Ryan},
}

@INPROCEEDINGS{Sanyal2023,
    author={Sanyal, Sourav and Roy, Kaushik},
    booktitle=ieeeicra23, 
    title={{RAMP-Net: A Robust Adaptive MPC for Quadrotors via Physics-informed Neural Network}}, 
    year={2023},
    volume={},
    number={},
    pages={1019-1025},
    doi={10.1109/ICRA48891.2023.10161410}
}

@INPROCEEDINGS{Schmitz23,
    author={Schmitz, Philipp and Lanza, Lukas and Worthmann, Karl},
    booktitle=ICSTCC23, 
    title={Safe data-driven reference tracking with prescribed performance}, 
    year={2023},
    volume={},
    number={},
    pages={454-460},
    doi={10.1109/ICSTCC59206.2023.10308521}
}

@article{Wabersich23,
    author={Wabersich, Kim P. and Zeilinger, Melanie N.},
    journal=ieeetac, 
    title={Predictive {C}ontrol {B}arrier {F}unctions: {E}nhanced {S}afety {M}echanisms for {L}earning-{B}ased {C}ontrol}, 
    year={2023},
    volume={68},
    number={5},
    pages={2638-2651},
    doi={10.1109/TAC.2022.3175628}
}

@ARTICLE{wabersich2023data,
    author={Wabersich, Kim P. and Taylor, Andrew J. and Choi, Jason J. and Sreenath, Koushil and Tomlin, Claire J. and Ames, Aaron D. and Zeilinger, Melanie N.},
    journal=ieeeconsysMag, 
    title={{Data-Driven Safety Filters: Hamilton-Jacobi Reachability, Control Barrier Functions, and Predictive Methods for Uncertain Systems}}, 
    year={2023},
    volume={43},
    number={5},
    pages={137-177},
    doi={10.1109/MCS.2023.3291885}
}

@article{faulwasser2023behavioral,
    title={Behavioral theory for stochastic systems? {A} data-driven journey from {W}illems to {W}iener and back again},
    author={Faulwasser, Timm and Ou, Ruchuan and Pan, Guanru and Schmitz, Philipp and Worthmann, Karl},
    journal= arcontrol,
    year={2023},
    publisher={Elsevier},
    doi= {10.1016/j.arcontrol.2023.03.005}
}

@mastersthesis{Hachmeister23,
    author = {Hachmeister, Jan},
    title = {{Funnel Control of Nonlinear Systems in the Presence of Delays}},
    school = {Paderborn University},
    address = {Paderborn, Germany},
    year = {2023} 
}

@ARTICLE{Facchino23,
  author={Facchino, Matteo and Bemporad, Alberto and Zanon, Mario},
  journal=ieeecsl, 
  title={{Tracking MPC Tuning in Continuous Time: A First-Order Approximation of Economic MPC}}, 
  year={2023},
  volume={7},
  number={},
  pages={2197-2202},
  doi={10.1109/LCSYS.2023.3284763}
}

@ARTICLE{matschek2023safe,
    author={Matschek, Janine and Bethge, Johanna and Findeisen, Rolf},
    journal=ieeetcst, 
    title={{Safe Machine-Learning-Supported Model Predictive Force and Motion Control in Robotics}}, 
    year={2023},
    volume={31},
    number={6},
    pages={2380-2392},
    doi={10.1109/TCST.2023.3262987}
}

@article{Sasfi2023,
    title = {{Robust adaptive MPC using control contraction metrics}},
    journal = auto,
    volume = {155},
    pages = {111169},
    year = {2023},
    issn = {0005-1098},
    doi = {10.1016/j.automatica.2023.111169},
    author = {András Sasfi and Melanie N. Zeilinger and Johannes Köhler},
}

@article{nuske2023finite,
    title={Finite-data error bounds for {K}oopman-based prediction and control},
    author={N{\"u}ske, Feliks and Peitz, Sebastian and Philipp, Friedrich and Schaller, Manuel and Worthmann, Karl},
    journal=journonlin,
    volume={33},
    number={1},
    pages={14},
    year={2023},
    publisher={Springer},
    doi= {10.1007/s00332-022-09862-1}
}

@article{lanza2024derivative,
    title = {On derivative-free sample-and-hold control with prescribed performance},
    journal = ifaconline,
    volume = {58},
    number = {17},
    pages = {121-126},
    year = {2024},
    note = mtns24,
    issn = {2405-8963},
    author = {Lukas Lanza},
    doi = {10.1016/j.ifacol.2024.10.124}
}

@article{wang2024comprehensive,
    title={{A Comprehensive Survey of Continual Learning: Theory, Method and Application}},
    author={Wang, Liyuan and Zhang, Xingxing and Su, Hang and Zhu, Jun},
    journal=ieeetpami,
    year={2024},
    volume={46},
    number={8},
    pages={5362-5383},
    publisher={IEEE},
    doi={10.1109/TPAMI.2024.3367329}
}

@article{GottschalkLanza24,
    title = {{Reinforcement Learning for Docking Maneuvers with Prescribed Performance}},
    journal = ifaconline,
    volume = {58},
    number = {17},
    pages = {196-201},
    year = {2024},
    note = mtns24,
    issn = {2405-8963},
    author = {Simon Gottschalk and Lukas Lanza and Karl Worthmann and Kerstin Lux-Gottschalk},
    doi = {10.1016/j.ifacol.2024.10.163}
}

@article{Berg24,
    author={Berger, Thomas},
    journal=ieeetac, 
    title={{Input-Constrained Funnel Control of Nonlinear Systems}}, 
    year={2024},
    volume={69},
    number={8},
    pages={5368-5382},
    doi={10.1109/TAC.2024.3352362}
}

@INPROCEEDINGS{Berger24Internal,
  author={Berger, Thomas and Hackl, Christoph M. and Trenn, Stephan},
  booktitle=ecc2024, 
  title={{Asymptotic Tracking by Funnel Control with Internal Models}}, 
  year={2024},
  volume={},
  number={},
  pages={1776-1781},
  doi={10.23919/ECC64448.2024.10590783}}

@ARTICLE{Tagliabue2024,
    author={Tagliabue, Andrea and How, Jonathan P.},
    journal=ieeetrob, 
    title={{Efficient Deep Learning of Robust Policies From MPC Using Imitation and Tube-Guided Data Augmentation}}, 
    year={2024},
    volume={40},
    number={},
    pages={4301-4321},
    doi={10.1109/TRO.2024.3431988}
}

@article{hassanpour2024practically,
    author = {Hassanpour, Hesam and Wang, Xiaonian and Corbett, Brandon and Mhaskar, Prashant},
    title = {A practically implementable reinforcement learning-based process controller design},
    journal = aiche,
    volume = {70},
    number = {1},
    eid = {e18245},
    doi = {10.1002/aic.18245},
    year = {2024}
}

@article{KohnPhil24,
    author = {K\"{o}hne, Frederik and Philipp, Friedrich M. and Schaller, Manuel and Schiela, Anton and Worthmann, Karl},
    title = {\(\boldsymbol{L}^{\boldsymbol{\infty }}\)-error {B}ounds for {A}pproximations of the {K}oopman {O}perator by {K}ernel {E}xtended {D}ynamic {M}ode {D}ecomposition},
    journal = siamads,
    volume = {24},
    number = {1},
    pages = {501-529},
    year = {2025},
    doi = {10.1137/24M1650120},
}

@ARTICLE{StraScha24:generator,
    author={Strässer, Robin and Schaller, Manuel and Worthmann, Karl and Berberich, Julian and Allgöwer, Frank},
    journal=ieeetac, 
    title={{Koopman-Based Feedback Design With Stability Guarantees}}, 
    year={2025},
    volume={70},
    number={1},
    pages={355-370},
    doi={10.1109/TAC.2024.3425770}
}

@ARTICLE{BoldGrun25,
    author={Bold, Lea and Grüne, Lars and Schaller, Manuel and Worthmann, Karl},
    journal=ieeetac, 
    title={{Data-Driven MPC With Stability Guarantees Using Extended Dynamic Mode Decomposition}}, 
    year={2025},
    volume={70},
    number={1},
    pages={534-541},
    doi={10.1109/TAC.2024.3431169}
}

@article{BoldScha25,
    title = {Kernel {EDMD} for data-driven nonlinear {K}oopman {MPC} with stability guarantees},
    journal = {IFAC-PapersOnLine},
    volume = {59},
    number = {19},
    pages = {478-483},
    year = {2025},
    note = ifacnolcos13,
    issn = {2405-8963},
    doi = {10.1016/j.ifacol.2025.11.080},
    author = {Lea Bold and Manuel Schaller and Irene Schimperna and Karl Worthmann},
}

@article{BoldLanzWoth2024_Koopman,
    title = {{Two-component controller design to safeguard data-driven predictive control -- A tutorial exemplified with DeePC and Koopman MPC}},
    author = {Lea Bold and Lukas Lanza and Karl Worthmann},
    pages = {383--397},
    volume = {73},
    number = {6},
    journal = atauto,
    doi = {10.1515/auto-2024-0166},
    year = {2025},
}

@article{Pillonetto2025,
    title = {{Deep networks for system identification: A survey}},
    journal = auto,
    volume = {171},
    pages = {111907},
    year = {2025},
    issn = {0005-1098},
    doi = {10.1016/j.automatica.2024.111907},
    author = {Gianluigi Pillonetto and Aleksandr Aravkin and Daniel Gedon and Lennart Ljung and Antônio H. Ribeiro and Thomas B. Schön},
}
\endgroup

\end{document}